\documentclass[11pt,fleqn,openany]{book} 

\usepackage[top=3cm,bottom=3cm,left=3.2cm,right=3.2cm,headsep=10pt,a4paper]{geometry} 

\usepackage{xcolor} 
\definecolor{ocre}{RGB}{243,102,25} 

\usepackage{avant} 
\usepackage{mathptmx} 

\usepackage{microtype} 
\usepackage[utf8]{inputenc} 
\usepackage[T1]{fontenc} 

\usepackage{todonotes}
\usepackage{tikz}
\usepackage{caption}
\usepackage{multicol}
\usepackage[absolute,overlay]{textpos}

\usetikzlibrary{arrows}
\usetikzlibrary{decorations.markings}

\usepackage[style=alphabetic,sorting=nyt,sortcites=true,autopunct=true,babel=hyphen,hyperref=true,abbreviate=false,backref=true,backend=bibtex]{biblatex}
\defbibheading{bibempty}{}

\usepackage{calc} 
\usepackage{makeidx} 
\makeindex 

\newcommand{\mypublic}{1}

\usepackage{titlesec} 

\usepackage{graphicx} 
\graphicspath{{Pictures/}} 

\usepackage{lipsum} 

\usepackage{tikz} 

\usepackage[english]{babel} 

\usepackage{enumitem} 
\setlist{nolistsep} 

\usepackage{booktabs} 

\usepackage{eso-pic} 

\usepackage{titletoc} 

\contentsmargin{0cm} 
\titlecontents{chapter}[1.25cm] 
{\addvspace{15pt}\large\sffamily\bfseries} 
{\color{ocre!60}\contentslabel[\Large\thecontentslabel]{1.25cm}\color{ocre}} 
{}  
{\color{ocre!60}\normalsize\sffamily\bfseries\;\titlerule*[.5pc]{.}\;\thecontentspage} 
\titlecontents{section}[1.25cm] 
{\addvspace{5pt}\sffamily} 
{\contentslabel[\thecontentslabel]{1.25cm}} 
{}
{\sffamily\hfill\color{black}\thecontentspage} 
[]
\titlecontents{subsection}[1.25cm] 
{\addvspace{1pt}\sffamily\small} 
{\contentslabel[\thecontentslabel]{1.25cm}} 
{}
{\sffamily\;\titlerule*[.5pc]{.}\;\thecontentspage} 
[] 

\titlecontents{lsection}[0em] 
{\footnotesize\sffamily} 
{}
{}
{}

\titlecontents{lsubsection}[.5em] 
{\normalfont\footnotesize\sffamily} 
{}
{}
{}
 
\usepackage{fancyhdr} 

\fancypagestyle{plain}{\fancyhead{}} 

\makeatletter
\renewcommand{\cleardoublepage}{
\clearpage\ifodd\c@page\else
\hbox{}
\vspace*{\fill}
\thispagestyle{empty}
\newpage
\fi}

\usepackage{amsmath,amsfonts,amssymb,amsthm} 

\newtheorem{notation}{Notation}[chapter]

\newtheoremstyle{ocrenumbox}
{0pt}
{0pt}
{\normalfont}
{}
{\small\bf\sffamily\color{ocre}}
{\;}
{0.25em}
{\small\sffamily\color{ocre}\thmname{#1}\nobreakspace\thmnumber{\@ifnotempty{#1}{}\@upn{#2}}
\thmnote{\nobreakspace\the\thm@notefont\sffamily\bfseries\color{black}---\nobreakspace#3.}} 

\newtheoremstyle{blacknumex}
{5pt}
{5pt}
{\normalfont}
{} 
{\small\bf\sffamily}
{\;}
{0.25em}
{\small\sffamily{\tiny\ensuremath{\blacksquare}}\nobreakspace\thmname{#1}\nobreakspace\thmnumber{\@ifnotempty{#1}{}\@upn{#2}}
\thmnote{\nobreakspace\the\thm@notefont\sffamily\bfseries---\nobreakspace#3.}}

\newtheoremstyle{blacknumbox} 
{0pt}
{0pt}
{\normalfont}
{}
{\small\bf\sffamily}
{\;}
{0.25em}
{\small\sffamily\thmname{#1}\nobreakspace\thmnumber{\@ifnotempty{#1}{}\@upn{#2}}
\thmnote{\nobreakspace\the\thm@notefont\sffamily\bfseries---\nobreakspace#3.}}

\newtheoremstyle{ocrenum}
{5pt}
{5pt}
{\normalfont}
{}
{\small\bf\sffamily\color{ocre}}
{\;}
{0.25em}
{\small\sffamily\color{ocre}\thmname{#1}\nobreakspace\thmnumber{\@ifnotempty{#1}{}\@upn{#2}}
\thmnote{\nobreakspace\the\thm@notefont\sffamily\bfseries\color{black}---\nobreakspace#3.}} 
\makeatother

\newcounter{dummy} 
\numberwithin{dummy}{section}
\theoremstyle{ocrenumbox}
\newtheorem{theoremeT}[dummy]{Theorem}
\newtheorem{thesisT}[dummy]{Thesis}

\newtheorem{exerciseT}{Exercise}[chapter]
\theoremstyle{blacknumex}
\newtheorem{exampleT}{Example}[chapter]
\theoremstyle{blacknumbox}

\newtheorem{definitionT}{Definition}[section]
\newtheorem{corollaryT}[dummy]{Corollary}
\theoremstyle{ocrenum}
\newtheorem{proposition}[dummy]{Proposition}
\newtheorem{lemma}[dummy]{Lemma}
\newtheorem{claim}[dummy]{Claim}
\newtheorem{widget}[dummy]{Widget}
\newtheorem{question}{Question}[chapter]

\RequirePackage[framemethod=default]{mdframed} 

\newmdenv[skipabove=10pt,
skipbelow=7pt,
backgroundcolor=black!5,
linecolor=ocre,
innerleftmargin=5pt,
innerrightmargin=5pt,
innertopmargin=5pt,
leftmargin=0cm,
rightmargin=0cm,
innerbottommargin=5pt]{tBox}

\newmdenv[skipabove=10pt,
skipbelow=7pt,
rightline=false,
leftline=true,
topline=false,
bottomline=false,
backgroundcolor=ocre!10,
linecolor=ocre,
innerleftmargin=5pt,
innerrightmargin=5pt,
innertopmargin=5pt,
innerbottommargin=5pt,
leftmargin=0cm,
rightmargin=0cm,
linewidth=4pt]{eBox}	

\newmdenv[skipabove=10pt,
skipbelow=7pt,
rightline=false,
leftline=true,
topline=false,
bottomline=false,
linecolor=ocre,
innerleftmargin=5pt,
innerrightmargin=5pt,
innertopmargin=0pt,
leftmargin=0cm,
rightmargin=0cm,
linewidth=4pt,
innerbottommargin=0pt]{dBox}	

\newmdenv[skipabove=10pt,
skipbelow=7pt,
rightline=false,
leftline=true,
topline=false,
bottomline=false,
linecolor=gray,
backgroundcolor=black!5,
innerleftmargin=5pt,
innerrightmargin=5pt,
innertopmargin=5pt,
leftmargin=0cm,
rightmargin=0cm,
linewidth=4pt,
innerbottommargin=5pt]{cBox}

\newenvironment{theorem}{\begin{tBox}\begin{theoremeT}}{\end{theoremeT}\end{tBox}}
\newenvironment{thesis}{\begin{tBox}\begin{thesisT}}{\end{thesisT}\end{tBox}}
				  
\newenvironment{definition}{\begin{dBox}\begin{definitionT}}{\end{definitionT}\end{dBox}}	
\newenvironment{example}{\begin{exampleT}}{\hfill{\tiny\ensuremath{\blacksquare}}\end{exampleT}}		
\newenvironment{corollary}{\begin{cBox}\begin{corollaryT}}{\end{corollaryT}\end{cBox}}	

\newenvironment{proof*}[1][\proofname]{
  
  \begin{proof}[#1]}{\end{proof}}

\newenvironment{answer}[0]{
  
  \begin{proof}[Answer]}{\end{proof}}

\makeatletter
\renewcommand{\@seccntformat}[1]{\llap{\textcolor{ocre}{\csname the#1\endcsname}\hspace{1em}}}                    
\renewcommand{\section}{\@startsection{section}{1}{\z@}
{-4ex \@plus -1ex \@minus -.4ex}
{1ex \@plus.2ex }
{\normalfont\large\sffamily\bfseries}}
\renewcommand{\subsection}{\@startsection {subsection}{2}{\z@}
{-3ex \@plus -0.1ex \@minus -.4ex}
{0.5ex \@plus.2ex }
{\normalfont\sffamily\bfseries}}
\renewcommand{\subsubsection}{\@startsection {subsubsection}{3}{\z@}
{-2ex \@plus -0.1ex \@minus -.2ex}
{.2ex \@plus.2ex }
{\normalfont\small\sffamily\bfseries}}                        
\renewcommand\paragraph{\@startsection{paragraph}{4}{\z@}
{-2ex \@plus-.2ex \@minus .2ex}
{.1ex}
{\normalfont\small\sffamily\bfseries}}

\usepackage{hyperref}
\hypersetup{hidelinks,backref=true,pagebackref=true,hyperindex=true,colorlinks=false,breaklinks=true,urlcolor= ocre,bookmarks=true,bookmarksopen=false,pdftitle={Title},pdfauthor={Author}}

\newcommand{\thechapterimage}{}
\newcommand{\chapterimage}[1]{\renewcommand{\thechapterimage}{#1}}

\def\thechapter{\arabic{chapter}}
\def\@makechapterhead#1{
\thispagestyle{empty}
{\centering \normalfont\sffamily
\ifnum \c@secnumdepth >\m@ne
\if@mainmatter
\startcontents
\begin{tikzpicture}[remember picture,overlay]
\node at (current page.north west)
{\begin{tikzpicture}[remember picture,overlay]
\node[anchor=north west,inner sep=0pt] at (0,5) {\includegraphics[width=\paperwidth]{\thechapterimage}};
\draw[anchor=west] (5cm,-4cm) node [rounded corners=20pt,fill=ocre!10!white,text opacity=1,draw=ocre,draw opacity=1,line width=1.5pt,fill opacity=.6,inner sep=12pt]{\huge\sffamily\bfseries\textcolor{black}{\thechapter. #1\strut\makebox[22cm]{}}};
\end{tikzpicture}};
\end{tikzpicture}}
\par\vspace*{90\p@}
\fi
\fi}

\def\@makeschapterhead#1{
\thispagestyle{empty}
{\centering \normalfont\sffamily
\ifnum \c@secnumdepth >\m@ne
\if@mainmatter
\begin{tikzpicture}[remember picture,overlay]
\node at (current page.north west)
{\begin{tikzpicture}[remember picture,overlay]
\node[anchor=north west,inner sep=0pt] at (0,5) {\includegraphics[width=\paperwidth]{\thechapterimage}};
\draw[anchor=west] (5cm,-4cm) node [rounded corners=20pt,fill=ocre!10!white,fill opacity=.6,inner sep=12pt,text opacity=1,draw=ocre,draw opacity=1,line width=1.5pt]{\huge\sffamily\bfseries\textcolor{black}{#1\strut\makebox[22cm]{}}};
\end{tikzpicture}};
\end{tikzpicture}}
\par\vspace*{90\p@}
\fi
\fi
}
\makeatother

\newcommand{\N}{\mathbb{N}}

\newcommand{\restr}{\upharpoonright}

\newcommand{\opcross}{\mathop{Cross}}

\DeclareMathOperator{\dnrf}{\mathrm{DNR}}
\DeclareMathOperator{\atoms}{\mathrm{atoms}}

\DeclareMathOperator{\true}{\textup{\texttt{T}}}
\DeclareMathOperator{\false}{\textup{\texttt{F}}}

\newcommand{\concat}{\frown}

\newcommand{\dbf}{\mathbf{d}}
\newcommand{\ebf}{\mathbf{e}}
\newcommand{\andd}{\wedge}
\newcommand{\orr}{\vee}
\newcommand{\la}{\langle}
\newcommand{\ra}{\rangle}
\newcommand{\da}{\!\downarrow}
\newcommand{\ua}{\!\uparrow}
\newcommand{\imp}{\rightarrow}

\newcommand{\biimp}{\leftrightarrow}
\newcommand{\Ab}{\mathbb{A}}
\newcommand{\TPb}{\mathbb{U}}
\newcommand{\Tb}{\mathbb{T}}
\newcommand{\Nb}{\mathbb{N}}
\newcommand{\Pb}{\mathbb{P}}
\newcommand{\Qb}{\mathbb{Q}}
\newcommand{\Rb}{\mathbb{R}}

\newcommand{\Psf}{\mathsf{P}}
\newcommand{\Qsf}{\mathsf{Q}}
\newcommand{\Rsf}{\mathsf{R}}
\newcommand{\Tsf}{\mathsf{T}}
\newcommand{\Acal}{\mathcal{A}}

\newcommand{\Ccal}{\mathcal{C}}
\newcommand{\Dcal}{\mathcal{D}}
\newcommand{\Ecal}{\mathcal{E}}
\newcommand{\Fcal}{\mathcal{F}}
\newcommand{\Gcal}{\mathcal{G}}

\newcommand{\Ical}{\mathcal{I}}
\newcommand{\Lcal}{\mathcal{L}}
\newcommand{\Mcal}{\mathcal{M}}
\newcommand{\Ncal}{\mathcal{N}}
\newcommand{\Pcal}{\mathcal{P}}
\newcommand{\Qcal}{\mathcal{Q}}
\newcommand{\Rcal}{\mathcal{R}}
\newcommand{\Scal}{\mathcal{S}}
\newcommand{\Kcal}{\mathcal{K}}
\newcommand{\Tcal}{\mathcal{T}}
\newcommand{\Wcal}{\mathcal{W}}

\newcommand{\cs}{2^\omega}
\newcommand{\str}{2^{<\omega}}
\newcommand{\uh}{{\upharpoonright}}

\renewcommand{\setminus}{\smallsetminus}

\newcommand{\bad}{\mathrm{Bad}}

\newcommand{\range}{\mathrm{ran}}
\newcommand{\seto}{\mathrm{set}}
\newcommand{\parts}{\mathrm{parts}}

\newcommand{\dom}{\mathrm{dom}}

\newcommand{\inter}{\mbox{int}}

\DeclareMathOperator{\red}{\textup{\texttt{red}}}
\DeclareMathOperator{\blue}{\textup{\texttt{blue}}}

\newcommand{\set}[1]{\left\{ #1 \right\}}
\newcommand{\card}[1]{\left| #1 \right|}

\newcommand{\tuple}[1]{\left\langle #1 \right\rangle}

\newcommand{\cond}[1]{\left\{\begin{array}{ll} #1 \end{array}\right.}

\newcommand{\nbd}{\nobreakdash-\hspace{0pt}}

\newcommand{\s}[1]{\ensuremath{\sf{#1}}}
\newcommand{\rans}[1]{#1\mbox{-}\s{RAN}}
\newcommand{\wwkls}[1]{#1\mbox{-}\s{WWKL}}
\newcommand{\gens}[1]{#1\mbox{-}\s{GEN}}
\newcommand{\dnrs}[1]{#1\mbox{-}\s{DNR}}
\newcommand{\sdnrs}[1]{#1\mbox{-}\s{SDNR}}
\newcommand{\rwkls}[1]{#1\mbox{-}\s{RWKL}}

\newcommand{\seqs}[1]{\s{Seq}^{*}(#1)}

\DeclareMathOperator{\er}{\s{ER}}
\DeclareMathOperator{\erp}{\s{\eta \to (\aleph_0,\eta)}^2}
\DeclareMathOperator{\ers}{\s{\eta \to (\eta)^1_{<\infty}}}
\DeclareMathOperator{\erps}{\s{\eta \to (\aleph_0,\eta)^1}}
\DeclareMathOperator{\tto}{\s{TT}}
\DeclareMathOperator{\dsf}{\s{D}}
\DeclareMathOperator{\rca}{\s{RCA}_0}
\DeclareMathOperator{\aca}{\s{ACA}}
\DeclareMathOperator{\hyp}{\s{HYP}}
\DeclareMathOperator{\wkl}{\s{WKL}}
\DeclareMathOperator{\atr}{\s{ATR}}
\DeclareMathOperator{\pioca}{\Pi^1_1\s{-CA}}
\DeclareMathOperator{\wwkl}{\s{WWKL}}
\DeclareMathOperator{\dnr}{\s{DNR}}
\DeclareMathOperator{\dnrzp}{\dnrs{2}}
\DeclareMathOperator{\isig}{\s{I}\Sigma}
\DeclareMathOperator{\ipi}{\s{I}\Pi}
\DeclareMathOperator{\bpi}{\s{B}\Pi}
\DeclareMathOperator{\bsig}{\s{B}\Sigma}
\DeclareMathOperator{\bst}{\bsig^0_2}
\DeclareMathOperator{\ist}{\isig^0_2}

\DeclareMathOperator{\rkl}{\s{RKL}}
\DeclareMathOperator{\rwkl}{\s{RWKL}}
\DeclareMathOperator{\rwwkl}{\s{RWWKL}}

\DeclareMathOperator{\rt}{\s{RT}}
\DeclareMathOperator{\srt}{\s{SRT}}
\DeclareMathOperator{\kl}{\s{KL}}
\DeclareMathOperator{\crt}{\s{CRT}}
\DeclareMathOperator{\rrt}{\s{RRT}}
\DeclareMathOperator{\srrt}{\s{SRRT}}
\DeclareMathOperator{\wsrrt}{\s{wSRRT}}

\DeclareMathOperator{\ipt}{\s{IPT}}
\DeclareMathOperator{\ads}{\s{ADS}}
\DeclareMathOperator{\sads}{\s{SADS}}
\DeclareMathOperator{\cads}{\s{CADS}}
\DeclareMathOperator{\cac}{\s{CAC}}

\DeclareMathOperator{\scac}{\s{SCAC}}

\DeclareMathOperator{\sat}{\s{SAT}}
\DeclareMathOperator{\rsat}{\s{RSAT}}

\DeclareMathOperator{\rcolor}{\s{RCOLOR}}

\DeclareMathOperator{\coh}{\s{COH}}

\DeclareMathOperator{\pizog}{\Pi^0_1\s{G}}
\DeclareMathOperator{\fip}{\s{FIP}}
\DeclareMathOperator{\ts}{\s{TS}}
\DeclareMathOperator{\sts}{\s{STS}}

\DeclareMathOperator{\psrt}{\s{psRT}}

\DeclareMathOperator{\fs}{\s{FS}}
\DeclareMathOperator{\sfs}{\s{SFS}}
\DeclareMathOperator{\emo}{\s{EM}}
\DeclareMathOperator{\semo}{\s{SEM}}

\DeclareMathOperator{\opt}{\s{OPT}}
\DeclareMathOperator{\amt}{\s{AMT}}

\DeclareMathOperator{\ivt}{\s{IVT}}

\newcommand{\Apx}{\operatorname{Apx}}
\newcommand{\Ext}{\operatorname{Ext}}

\begin{document}


\begingroup

\ifdefined\mypublic

\thispagestyle{empty}
\centering
\vspace*{9cm}
\par\normalfont\fontsize{30}{30}\sffamily\selectfont
The reverse mathematics of Ramsey-type theorems\par 

\vspace*{1cm}
{\Huge Ludovic Patey}\par 

\else

\thispagestyle{empty}
\setlength{\TPHorizModule}{1mm}
\setlength{\TPVertModule}{1mm}

\begin{textblock}{20}(185,35)
\includegraphics[scale=0.3]{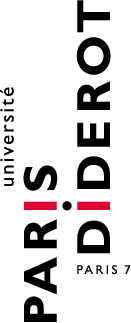}
\end{textblock}
\begin{textblock}{20}(185,15)
\includegraphics[scale=0.25]{logo-irif.png}
\end{textblock}

\begin{center}

{\large \bf\sc Université Paris Diderot -- Paris VII\\
Sorbonne Paris Cité}
\vspace{1cm}

{\large \sc École Doctorale Sciences Mathématiques de Paris Centre}
\vspace*{2cm}

{\LARGE \bf TH\`ESE}\\
\vspace{10pt}
en vue d'obtenir le grade de\\
\vspace{10pt}
{\bf DOCTEUR DE L'UNIVERSITÉ PARIS DIDEROT}\\
\vspace{10pt}
en {\sc Informatique Fondamentale}\\
\vspace{1.5cm}

{\huge\bf Les mathématiques à rebours\\[10pt] de théorèmes de type Ramsey}\\

\vspace{2cm}

{\itshape Présentée et soutenue par}\\[10pt]
{\large\bf Ludovic Patey}\\[10pt]
{\itshape le 26 février 2016}
\vspace{1cm}

{\sc Rapporteurs} :\\[10pt]
\begin{tabular}{r@{\ }ll}
  & M. Alberto {\sc Marcone} & Université d'Udine  \\
  & M. Reed {\sc Solomon} & Université du Connecticut  \\
\end{tabular}
\vspace{1cm}

{\sc Jury} :\\[10pt]
\begin{tabular}{r@{\ }ll}
  & M. Laurent {\sc Bienvenu} & Directeur de thèse\\
  & M. Hugo {\sc Herbelin} & Directeur de thèse \\
  & M. Alberto {\sc Marcone} & Rapporteur  \\
  & M. Serge {\sc Grigorieff} & Examinateur  \\
  & M. Mathieu {\sc Hoyrup} & Examinateur  \\
\end{tabular}

\end{center}

\fi

\endgroup


\chapterimage{chapter_head_1.pdf} 

\pagestyle{empty} 

\setcounter{tocdepth}{1}
\tableofcontents 

\ifdefined\mypublic
\else
\cleardoublepage 
\fi

\pagestyle{fancy} 

\ifdefined\mychapter
\input{\mychapter}
\else

\chapterimage{chapter_head_2.pdf}
\ifdefined\mypublic
\chapter*{Abstract}
\else
\chapter*{Résumé / Abstract}
\fi
\addcontentsline{toc}{chapter}{\textcolor{ocre}{Abstract}}

\ifdefined\mypublic
\bigskip
\else
\begin{center}
\textbf{Les mathématiques à rebours de théorèmes de type Ramsey}
\end{center}

Dans cette thèse, nous investiguons le contenu calculatoire
et la force logique du théorème de Ramsey et de ses conséquences.
Pour cela, nous utilisons les outils des mathématiques à rebours et
de la réduction calculatoire. Nous procédons à une analyse systématique
de divers énoncés de type Ramsey à travers des outils unifiés et minimalistes,
et obtenons une analyse précise de leurs interdépendances.

Nous clarifions notamment le role du nombre de couleurs 
dans le théorème de Ramsey. En particulier, nous montrons que la hiérarchie
du théorème de Ramsey induite par le nombre de couleurs
est strictement croissante au niveau des réductions calculatoires,
et exhibons en mathématiques à rebours une hiérarchie
infinie décroissante de théorèmes de type Ramsey en affaiblissant
les contraintes d'homogénéité.
Ces résultats tendent à montrer que les énoncés de type Ramsey ne sont pas robustes,
c'est-à-dire que de faibles variations des énoncés mènent à des sous-systèmes 
strictement différents.

Enfin, nous poursuivons l'analyse des liens entre les théorèmes de Ramsey
et les arguments de compacité, en étendant le théorème de Liu à de nombreux énoncés
de type Ramsey et en prouvant son optimalité sous différents aspects.

\bigskip
\fi

\begin{center}
\textbf{The reverse mathematics of Ramsey-type theorems}
\end{center}

In this thesis, we investigate the computational content
and the logical strength of Ramsey's theorem and its consequences.
For this, we use the frameworks of reverse mathematics
and of computable reducibility. We proceed to a systematic study
of various Ramsey-type statements under a unified and minimalistic framework
and obtain a precise analysis of their interrelations.

We clarify the role of the number of colors in Ramsey's theorem.
In particular, we show that the hierarchy of Ramsey's theorem induced by the number of colors
is strictly increasing over computable reducibility, and exhibit in reverse mathematics an 
infinite decreasing hiearchy of Ramsey-type theorems by weakening the homogeneity constraints.
These results tend to show that the Ramsey-type statements are not robust, that is,
slight variations of the statements lead to strictly different subsystems. 

Finally, we pursuit the analysis of the links between Ramsey's theorems 
and compacity arguments, by extending Liu's theorem to several Ramsey-type statements
and by proving its optimality under various aspects.

\chapter*{Acknowledgement}
\addcontentsline{toc}{chapter}{\textcolor{ocre}{Acknowledgement}}

I would like to thank my cat.
\footnote{Oh my God! I cannot believe I forgot to thank my Ph.D. advisor.
I would like to thank Laurent Welcome for accepting supervising me along those 3 years.
He has been successively a professor, an advisor and a friend. I learned by him
the basics of computability theory and reverse mathematics, but far beyond that,
he taught me how to become a researcher. My writing of papers, my presentations,
my vision of mathematics is highly influenced by Laurent. My mistakes are my own contribution.}
\footnote{I would like to thank my second Ph.D. advisor, Hugo Herbelin.
Although our mutual research interests early diverged, I really appreciated all the
insightful discussions we had about type theory.}
\footnote{I also would like to thank Alberto Marcone and Reed Solomon
who kindly accepted to serve as referee for my dissertation, and Serge Grigorieff
and Mathieu Hoyrup for being part of my defense commitee.}
\footnote{I realized with horror that Google automatically translated
``Bienvenu'' into ``Welcome''. I woul like to apologize for this awkward mistake.
I would certainly not allow myself to make such a terrible play on words with my advisor's name.
I will ask Google to stop translating ``Bienvenu'' into ``Welcome''.}
\footnote{I would like to give a special thank to Paul Shafer.
Along my Ph.D., he provided me profusion of advice, read carefully lots of my proofs
and answered my questions with a boundless patience. My debt to Paul is beyond calculation.}
\footnote{I also would like to thank the computability team that LIAFA used to have
during the golden age, namely, Serge Grigorieff, Christopher Porter, Benoît Monin, and Antoine Taveneaux.}
\footnote{I had the pleasure to visit the logic groups of UConn, CT during April 2015 and of Notre Dame, IN in November 2015.
In particular, I thank Peter Cholak, Damir Dzhafarov, Gregory Igusa, Manuel Lerman, 
Reed Solomon and Brown Westrick for making them such pleasant stays.}
\footnote{Google denied my request, arguing that ``Bienvenu'' was more often used
as a welcoming statement than to refer to Laurent. I deplore their lack of understanding.}
\footnote{I thank my two recent co-authors Emanuele Frittaion and Keita Yokoyama.
It was a really pleasant collaboration and I learned a lot with them.}
\footnote{More generally, I would like to thank the whole reverse mathematical community
for their warm bienv... welcome. In particular, I would like to thank Wei Wang for our
insightful discussions about forcing arguments, and 
Damir Dzhafarov and Denis Hirschfeldt for our ongoing (and promising) collaboration.}
\footnote{The PPS laboratory is composed of great people with whom I shared very pleasant moments.
In particular, I thank (in alphabetical order) Cyrille Chenavier, Maxime Lucas, Adrien Batmalle, Pierre Vial,
Charles Grellois, Thibaut Girka, Matthieu Boutier, Flavien Breuvart, Lourdes Gonzales, 
Marie Kerjean, Shahin Amini, Jovana Obradovic, Amina Doumane, Guillaume Claret, Ioana Cristescu, and Raphaelle Crubille.}
\footnote{My scholarship at the Ecole Normale Supérieure was probably the best period of my life. 
I met there awesome people, and in particular in my classmates, namely,
Antoine Amarilli, Fabrice Ben Hamouda, Yoann Bourse, Florence Clerc, Wenjie Fang, 
Marc Jeanmougin, Robin Morisset, Pablo Rauzy, and Hang Zhou.}
\footnote{I also met great people during my LMFI master, among whom, Thibaut Gauthier, 
Louis Ioos, Maxime Lucas, Sonia Marin, Pierre Vial, and Marc Verrière.}
\footnote{After some statistics on Google, I must admit that the vast majority of
occurences of ``Bienvenu'' are not citations of Laurent. I will nevertheless keep my high esteem of my Ph.D. advisor.}
\footnote{Some people outside academics played and continue to play an important role in my life. 
I thank Rémy Breuils, Sangmin Shin and Luc Gille.}
\footnote{Last, I would like to give a special thank to my girlfriend Lucile
for her endless support and for making my life brighter.}
\footnote{I now realize that I do not have a cat. Please, forget about him.}

\ifdefined\mypublic
\else
\chapter*{Introduction (Français)}
\addcontentsline{toc}{chapter}{\textcolor{ocre}{Introduction (Français)}}

Le présent document est un rapport de thèse résultant de trois années de recherche sous la supervision
de Laurent Bienvenu et Hugo Herbelin. Nous commençons par une introduction gentille aux domaines
mathématiques dont il est question, puis nous fournissons un résumé détaillé de la thèse ainsi que de ses
principales contributions.

\section*{Calculabilité et mathématiques à rebours}

Cette thèse contribue au domaine des mathématiques à rebours d'un point de vue calculatoire.
Nous introduirons brièvement la calculabilité, qui fournit un environnement de travail solide
pour parler de la notion d'ensemble calculable. Ensuite, nous effectuerons une introduction 
aux mathématiques à rebours avec une orientation calculatoire, en présentant ses principales motivations
et ses défis.

\subsection*{Calculabilité}

Informellement, un ensemble d'entiers naturels est \emph{calculable} si l'on peut décider de manière
effective quels sont les éléments lui appartenant. En ce sens, tout ensemble fini est calculable dans la mesure
où l'on peut apprendre par cœur la liste de ses éléments. En revanche, lorsqu'il s'agit d'ensembles infinis,
il faut trouver une méthode systématique (ou \emph{algorithme}) pour effectuer une telle décision.
La thèse de Church-Turing énonce que la notion d'ensemble calculable est robuste, en ce que tout
modèle de calcul raisonnable définit la même classe d'ensembles calculables. Ainsi, l'on peut
aisément se convaincre qu'un ensemble est calculable en écrivant un algorithme de décision
dans n'importe quel langage de programmation standard.

Bien que définie pour des entiers naturels, la notion d'ensemble calculable s'étend
à n'importe quel type de données par une simple technique de codage.
Toutes les fonctions de codage que nous considérerons sont calculables.
De ce fait, il sera possible de parler indifféremment d'une structure de données
ou de l'ensemble d'entiers naturels correspondant. Par exemple, un couple d'entiers $(a, b)$
peut être codé par l'entier $2^a3^b$. A partir d'une fonction de couplage,
il est possible de définir un codage des chaînes d'entiers, et ainsi de suite.

Combien y a-t-il d'ensembles calculables ? Dès lors que chaque ensemble calculable vient avec 
un algorithme de décision, il ne peut y en avoir qu'un nombre dénombrable. Par un simple argument diagonal,
l'on peut montrer que la collection de tous les ensembles n'est pas dénombrable, et ainsi que la vaste majorité
des ensembles d'entiers naturels est incalculable. L'exemple le plus connu d'un tel ensemble non calculable
est le \emph{problème de l'arrêt}, défini comme l'ensemble des (codes de) programmes qui terminent.

Lorsqu'un ensemble n'est pas calculable, il est naturel de se demander à quel point il ne l'est pas.
Par exemple, pourrions-nous calculer cet ensemble si nous étions capable de résoudre le problème de l'arrêt ?
Il est possible d'étendre la notion d'ensemble calculable en ajoutant la fonction caractéristique $\chi_Y$
d'un ensemble $Y$ parmi les primitives de son langage de programmation préféré. La fonction $\chi_Y$ n'est pas calculable,
et doit être par conséquent considérée comme un \emph{oracle}. Un ensemble~$X$ décidable par un algorithme
dans ce langage étendu est dit \emph{calculable relativement à $Y$}. On l'abrège usuellement par $X \leq_T Y$ (où $T$
tient de Turing). Si $X \leq_T Y$, alors l'ensemble $Y$ est \emph{au moins aussi difficile à calculer que $X$}.
Il est possible de définir une hiérarchie infinie de degrés d'incalculabilité, dans le même esprit que la hiérarchie
d'infinitudes révélée par Cantor.

La notion d'ensemble calculable admet une caractérisation purement logique en terme de formules de l'arithmétique.
Fixons un ensemble infini de \emph{variables de nombres} $x, y, z, \dots$, qui sont supposées prendre pour valeur
des entiers naturels. Nous pouvons former des \emph{expressions numériques} en les utilisant comme des paramètres
dans des additions et multiplications d'entiers. Par exemple, $(3+x) \cdot y$ est une expression numérique,
qui, après avoir remplacé $x$ par 18 et $y$ par 2, s'évalue en l'entier~42. Il est ensuite possible de définir
des \emph{formules} en reliant deux expressions numériques par une égalité $=$ ou une inégalité $\leq$, en utilisant
des quantificateurs sur les entiers $(\exists x)$ (il existe $x$) et $(\forall x)$ (pour tout $x$), des quantificateurs
bornés $(\exists x \leq y)$, $(\forall x \leq y)$, et en les composant avec des connecteurs logiques tels que
$\wedge$, $\vee$, $\neg$, $\imp$, $\equiv$ (et, ou, non, implique, si et seulement si). Par exemple,
la formule $(\forall x)(\exists y \leq x)(x = 2y \vee x = 2y + 1)$ est une formule affirmant que tout nombre
$x$ est soit pair, soit impair.

Une formule $\varphi$ contenant une variable \emph{libre}, c'est-à-dire, telle que $x$ n'apparaît pas sous le champ
d'un quantificateur, peut être associée à un ensemble $S_\varphi$, défini comme l'ensemble de tous les entiers naturels
$n$ tels que la formule $\varphi$ est satisfaite en replaçant~$x$ par~$n$. Par exemple, la formule $(\exists y)(x = 2y)$
définit l'ensemble $S = \{n : (\exists y)(n = 2y)\}$ qui n'est autre que l'ensemble des entiers pairs.
Il est possible de classifier la complexité des formules en fonction de leur alternance de quantificateurs.
Une formule est $\Sigma^0_1$ (resp.\ $\Pi^0_1$) si elle est de la forme $(\exists x)\varphi$ (resp.\ $(\forall x)\varphi$),
où $\varphi$ est une formule ne contenant que des quantificateurs bornés. Un ensemble est dit $\Sigma^0_1$
s'il est défini par une formule $\Sigma^0_1$. Les ensembles $\Pi^0_1$ sont définis similairement.
Intuitivement, les ensembles $\Sigma^0_1$ sont ceux dont les éléments peuvent être \emph{énumérés},
car il suffit de tester chaque valeur de~$x$ jusqu'à en trouver une pour laquelle la formule $\varphi$
correspondante est satisfaite. Un ensemble est $\Delta^0_1$ s'il est à la fois $\Sigma^0_1$ et $\Pi^0_1$.
Un célèbre théorème d'Emil Post énonce que les ensembles $\Delta^0_1$ sont \emph{exactement} les ensembles calculables.
Ce théorème n'est qu'un point de départ de toute une correspondance entre la calculabilité et la définissabilité par des formules.

Nous avons maintenant les bases de calculabilité nécessaires à une introduction aux mathématiques à rebours.

\subsection*{Mathématiques à rebours}

Un \emph{axiome} est un énoncé mathématique pris comme un postulat. Étant donné un ensemble d'axiomes,
souvent appelé \emph{théorie}, et des \emph{règles de déduction}, on peut dériver des conséquences mathématiques
appelées \emph{théorèmes}. D'un point de vue purement mathématique, un axiome et un théorème ne sont rien d'autre
que des énoncés mathématiques. La seule différence est qu'un axiome est auto-justifié, tandis qu'il faut fournir
la preuve d'un théorème.

Il est raisonnable de vouloir qu'une théorie soit \emph{calculatoirement énumerable}, c'est-à-dire,
que l'ensemble de ses axiomes soit $\Sigma^0_1$. En effet, l'on voudrait être capable de savoir
quels axiomes sont disponibles pour prouver un théorème. Heureusement, la vaste majorité des formalisations
de concepts aboutit de fait à des théories calculatoirement énumérables. C'est en particulier le cas de la théorie des
ensembles. Dans son premier théorème d'incomplétude, Gödel affirme qu'une telle théorie, si elle est cohérente
et suffisamment expressive, ne peut être qu'\emph{incomplète}, c'est-à-dire qu'elle contient des énoncés mathématiques
qui ne sont ni prouvables, ni réfutables. Ce théorème exclus tout espoir de trouver un jour un
ensemble calculable d'axiomes qui nous permettrait de décider n'importe quel énoncé mathématique. 
Il en résulte que le programme fondationnel de recherche d'axiomes naturels que nous pourrions ajouter
à la théorie des ensembles est un processus non terminé, et un sujet de recherche toujours actif.

Les mathématiques à rebours sont un vaste programme mathématique cherchant à déterminer
quels axiomes sont nécessaires pour prouver les théorèmes des mathématiques.
L'idée sous-jacente est très simple. Retirez tous les axiomes forts de la théorie des ensembles,
et gardez seulement une théorie de base très faible, dans laquelle presque aucun des théorèmes usuels
n'est prouvable. Cette théorie doit cependant être suffisamment forte pour pouvoir formaliser toutes
les astuces triviales de codage que nous utilisons, et ainsi ne pas dépendre de la présentation
des structures de données utilisées. Ensuite, choisissez un théorème ordinaire~$\Psf$,
utilisé dans la vie de tous les jours, et essayez de le prouver avec un nombre minimal d'axiomes.
Parfois, le théorème $\Psf$ est directement prouvable dans la théorie de base. Dans ce cas, cela signifie
que le théorème ne possède pas une preuve complexe, puisque nous avons choisi une théorie de base très faible.
Lorsque ce n'est pas le cas, nous devons ajouter certains axiomes $A_0, A_1, \dots, A_n$ à la théorie de base
pour prouver $\Psf$. A ce niveau, nous savons que les axiomes $A_0, A_1, \dots, A_n$ sont \emph{suffisants}
pour prouver $\Psf$, mais comment s'assurer qu'ils sont vraiment \emph{nécessaires} ? Il suffit simplement
d'essayer de prouver la réciproque, c'est-à-dire d'ajouter le théorème $\Psf$ à la théorie de base,
et essayer de prouver $A_0 \wedge A_1 \wedge \dots \wedge A_n$. Si nous réussissons, alors nous venons de prouver
que, dans la théorie de base, les axiomes $A_0, A_1, \dots, A_n$ sont les \emph{minimaux} requis pour prouver~$\Psf$.

Le programme des mathématiques à rebours a été initié par Harvey Friedman en 1975. Depuis lors, de nombreux chercheurs
ont contribué à son développement, le contributeur le plus notable étant Stephen Simpson. Des milliers
de théorèmes provenant des domaines principaux des mathématiques ont été étudiés,
notamment en algèbre, en analyse et en topologie. Un phénomène surprenant est apparu aux débuts
des mathématiques à rebours: la plupart des théorèmes étudiés ne nécessite que des axiomes très faibles.
Plus encore, nombre d'entre eux se trouvent être \emph{équivalents} à l'un parmi cinq grands ensembles
d'axiomes, plus connus sous le nom de \emph{Club des Cinq}. Ces cinq ensembles d'axiomes correspondent à
des approches philosophiques bien connues. Une partie de la recherche actuelle consiste à essayer de comprendre
ce phénomène, et en particulier à étudier les théorèmes échappant au Club des Cinq.

Nous allons maintenant entrer dans le détail des mathématiques à rebours et de la théorie de base.
Pour cela, nous devons étendre notre notion de formule à l'\emph{arithmétique du second ordre}. 
Prenons notre précédent langage, augmenté de primitives pour parler d'ensembles d'entiers naturels.
Plus précisément, il nous faut fixer une collection de \emph{variables d'ensembles} $X, Y, Z$,
et étendre notre précédente notion de formule avec la relation de possession $(n \in X)$
et des quantificateurs sur les ensembles $(\forall X)$, $(\exists X)$. Par exemple,
la formule $(\exists X)(\forall x)[x \in X \biimp x < 5]$ énonce que l'ensemble
des entiers naturels inférieurs à 5 existe. Harvey Friedman ayant remarqué que la grande
majorité des théorèmes mathématiques pouvait être énoncé de manière naturelle dans le langage
de l'arithmétique du second ordre, nous pouvons restreindre le programme des mathématiques à rebours
aux énoncés mathématiques de l'arithmétique du second ordre sans trop perdre en généralité.

Dans ce contexte, un axiome précise le comportement des entiers naturels et ensembles d'entiers
naturels, que nous appellerons respectivement \emph{partie du premier ordre} et \emph{partie
du second ordre}. Par exemple, l'énoncé ``$0 \neq 1$'' est un axiome. Plus la théorie de base
aura d'axiomes, moins ses modèles seront chaotiques. La théorie de base $\rca$,
signifiant Axiome de Compréhension Récursif, contient un ensemble d'axiomes du premier ordre
connu sous le nom d'\emph{arithmétique de Robinson}. Ces axiomes décrivent
le comportement des entiers naturels vis-à-vis de l'addition, la multiplication et de l'ordre.
Par le théorème de compacité, il n'est pas possible d'exclure tous les comportements
non-standards des entiers naturels simplement avec des axiomes du premier ordre, mais
l'arithmétique de Robinson s'assure qu'ils vont au moins se comporter comme prévu vis-à-vis
des opérations standards.

$\rca$ contient également des axiomes du second ordre, qui vont décrire quels ensembles
d'entiers naturels sont assurés d'exister. Les axiomes du second ordre de $\rca$ 
peuvent être classifiés selon deux catégories: les axiomes d'\emph{induction} 
et les axiomes de \emph{compréhension}. Étant donnée une formule $\varphi(x)$ 
avec une variable de nombre distinguée $x$, l'axiome d'induction pour $\varphi$ 
énonce que si $\varphi(0)$ est vrai, et $\varphi(n+1)$ est vrai lorsque $\varphi(n)$ 
l'est, alors $\varphi(n)$ est satisfaite pour tout entier naturel $n$. L'axiome de 
compréhension pour $\varphi$ énonce que la collection des $n$'s tels que $\varphi(n)$ 
est satisfaite existe en tant qu'ensemble. On dénote par $\{n : \varphi(n)\}$ cet 
ensemble. Les axiomes d'induction et de compréhension peuvent être vus comme des 
méthodes pour construire de nouveaux ensembles à partir de ceux existants dans les 
modèles. $\rca$ contient le schéma d'induction pour toute formule $\Sigma^0_1$, et 
le schéma de compréhension pour tout prédicat $\Delta^0_1$. Par la correspondance 
entre les prédicats $\Delta^0_1$ et les ensembles calculables, nous n'autorisons 
la construction de nouveaux ensembles que calculatoirement à partir des précédents 
ensembles. En ce sens, $\rca$ capture les \emph{mathématiques calculables}.

De nombreux théorèmes peuvent être vus comme des \emph{problèmes mathématiques}, 
venant avec une classe naturelle d'\emph{instances}, et de \emph{solutions} pour 
chaque instance. Par exemple, le lemme de K\"onig énonce pour tout arbre infini 
à branchement fini l'existence d'un chemin infini à travers l'arbre. Ici, 
l'instance est un arbre infini à branchement fini $T$, et une solution de $T$ 
est un chemin infini à travers $T$. Il existe une notion naturelle de \emph{réduction 
calculatoire} d'un problème $\Psf$ vers un autre problème $\Qsf$ (noté $\Psf \leq_c \Qsf$), 
signifiant que toute instance $X_0$ de $\Psf$ calcule une instance $X_1$ de $\Qsf$ 
telle que toute solution de $X_1$ calcule (avec l'aide de $X_0$) une solution de $X_0$. 
Informellement, $\Psf \leq_c \Qsf$ si nous pouvons résoudre calculatoirement le 
problème mathématique $\Psf$ à partir d'une boite noire permettant de résoudre $\Qsf$. 
Une réduction calculatoire de $\Psf$ vers $\Qsf$ peut être vue comme une preuve 
de $\Qsf \imp \Psf$ dans $\rca$, dans laquelle une seule application de $\Qsf$ 
est autorisée. En ce sens, la réduction calculatoire est plus \emph{précise} 
que les mathématiques à rebours. En particulier, cela permet de révéler de 
subtiles distinctions entre des énoncés mathématiques qui seraient indistinguables 
du point de vue des mathématiques à rebours.

\section*{Résumé de la thèse}

Nous commençons cette section en expliquant plus précisément le contenu de cette thèse. 
Nous détaillerons ensuite ces explications chapitre par chapitre, puis nous terminerons 
avec un résumé des principales contributions de cette thèse.

\subsection*{Le sujet}

Cette thèse porte sur les mathématiques à rebours d'énoncés mathématiques 
venant de la théorie de Ramsey. La théorie de Ramsey est une branche des 
mathématiques étudiant les conditions sous lesquelles une certaine structure 
apparaît dans une collection d'objets suffisamment grande. Par exemple, 
dans tout groupe de six personnes, au moins trois personnes se connaissent 
mutuellement ou se sont étrangères mutuellement. Le plus connu de ces théorèmes 
est peut-être le \emph{théorème de Ramsey}, qui énonce que pour tout $k$-coloriage 
de $[\Nb]^n$ (où $[\Nb]^n$ dénote les $n$-uplets d'entiers naturels), il existe 
un ensemble infini d'entiers naturels~$H$ tel que $[H]^n$ est monochromatique. 
Un tel ensemble est dit \emph{homogène}. En particulier, le théorème de Ramsey 
pour les singletons ($n=1$) n'est autre que la version infinie du \emph{principe des tiroirs}, 
qui énonce que si un nombre infini d'objets est réparti en $k$ boites, 
au moins l'une des boites contient un nombre infini d'objets.

Au cours des deux dernières décennies, la théorie de Ramsey est devenue 
l'un des sujets de recherche les plus importants en mathématiques à rebours. 
Cette théorie fournit un large panel de théorèmes échappant au phénomène 
du Club des Cinq, et dont la force logique est notoirement difficile à évaluer. 
Nous effectuons une étude systématique de nombreuses conséquences du théorème de Ramsey,
parmi lesquelles le théorème d'Erd\H{o}s-Moser, les théorèmes d'ensemble libre, 
d'ensemble mince et le théorème de l'arc-en-ciel de Ramsey, avec les outils 
des mathématiques à rebours et de la réduction calculatoire. En outre, nous 
les comparons à de nombreux théorèmes précédemment étudiés en mathématiques à rebours.

\subsection*{Structure de la thèse}

Cette thèse est divisée en 4 parties: un préambule, la force du théorème de Ramsey, 
des études thématiques, et une conclusion. Nous allons maintenant expliquer 
brièvement le contenu de chaque partie.
\bigskip

Dans le préambule, nous introduisons les concepts de base de la calculabilité que 
nous aurons à manipuler. Puis nous présentons une introduction détaillée des 
mathématiques à rebours et du phénomène du Club des Cinq. Ensuite, nous comparons 
les mathématiques à rebours à différentes notions de réductions introduites récemment, 
et qui mettent l'accent sur certains aspects calculatoires des théorèmes considérés. 
Enfin, nous introduisons la principale technique de séparation utilisée tout au 
long de cette thèse: le forcing effectif. Cette technique est un raffinement de la 
notion de forcing de Cohen dans le but de préserver certaines propriétés calculatoires.

Dans la second partie, nous effectuons une étude systématique du théorème de Ramsey 
et d'un certain nombre de ces conséquences, en mathématiques à rebours et avec la 
réduction calculatoire. Nous considérons principalement des théorèmes bien établis 
et qui ont déjà été introduits en mathématiques à rebours. Parmi eux, nous étudions 
la cohésion, une version stable du théorème de Ramsey pour les paires, le théorème 
d'Erd\H{o}s-Moser, les théorèmes d'ensemble libre, d'ensemble mince et le théorème 
de l'arc-en-ciel de Ramsey. Nous fournissons des preuves simplifiées de séparations 
existantes avec un ensemble d'outils uniformisés et minimalistes, et résolvons de 
nombreuses questions ouvertes.

La troisième partie est plus thématique. Nous étudions différents théorèmes en 
mathématiques à rebours sous certains angles calculatoires. En premier lieu, nous 
étudions l'existence d'instances maximalement difficiles à résoudre (\emph{instances 
universelles}) de théorèmes de type Ramsey, et les degrés de Turing dominant leurs 
solutions. Ensuite, nous montrons l'orthogonalité de différents théorèmes de type 
Ramsey vis-à-vis de la compacité. Nous identifions également une nouvelle classe 
de théorèmes généralisant le théorème de Ramsey pour les paires, et dont la force 
calculatoire est actuellement inconnue. Enfin, nous concevons de nouvelles notions 
de forcing avec de meilleures propriétés de définissabilité, dans l'espoir de prouver 
que les hiérarchies des théorèmes d'ensemble libre, d'ensemble mince et de l'arc-en-ciel 
de Ramsey, sont strictes en mathématiques à rebours.

Dans la conclusion, nous discutons de la naturalité des outils utilisés pour étudier 
les énoncés de type Ramsey en mathématiques à rebours, et essayons de donner de 
nouvelles perspectives sur le rôle des théorèmes de type Ramsey par rapport au phénomène 
du Club des Cinq. Nous fournissons ensuite une liste de questions ouvertes restantes, 
avec leurs motivations et des explications sur la raison de leur difficulté.

\subsection*{Principales contributions}

Au cours de cette thèse, nous résolvons 25 questions ouvertes, venant de 17 papiers 
différents. Nous allons maintenant détailler les principales contributions originales de cette thèse.
\bigskip

\begin{itemize}
	\item La principale contribution est probablement d'ordre méthodologique. En accord avec le fameux proverbe ``\emph{Si tu donnes un poisson à un homme, il mangera un jour. Si tu lui apprends à pêcher, il mangera toujours}'', nous sommes davantage intéressés par le développement d'outils permettant de séparer de manière systématique des théorèmes en mathématiques à rebours, plutôt que par les séparations en tant que telles. Dans~\cite{Patey2015Iterative}, nous avons présenté une simplification d'outils introduits par Lerman, Solomon et Towsner~\cite{Lerman2013Separating}. Nous montrons que ces outils sont suffisamment généralistes pour reprouver différentes séparations avec des arguments plus simples (Corollaire~\ref{cor:emo-wkl-sts2-sads}, Corollaire~\ref{cor:ads-not-implies-scac}) et séparer le théorème de Ramsey pour les pairs du théorème des arbres pour les pairs (Théorème~\ref{thm:tree-theorem-rt22+wkl-not-tt22}), répondant par cela même à une question de Mont\'alban~\cite{Montalban2011Open}.
\bigskip

	\item Parmi les contributions importantes, mentionnons la séparation du théorème de Ramsey pour $k+1$ couleurs du théorème de Ramsey pour $k$ couleurs sous la réduction calculatoire (Théorème~\ref{thm:colors-ramsey-general}). Cette question fut initialement posée sous une forme plus faible par Mileti~\cite{Mileti2004Partition}, puis pour la réduction de Weihrauch par Dorais, Dzhafarov, Hirst, Mileti, et Shafer~\cite{Dorais2016uniform} qui l'ont considérée comme la plus importantes des questions non résolues de leur article. Nous répondons aux deux questions. Hirschfeldt et Jockusch~\cite{Hirschfeldtnotions} on récemment demandé combien d'applications du théorème de Ramsey pour $m$ couleurs sont nécessaires pour résoudre le théorème de Ramsey pour $k$ couleurs en fonction de $m$ et $k$. Nous donnons une réponse précise à travers le Théorème~\ref{thm:colors-ramsey-general}.
	\bigskip
	
	\item Nous révélons que la hiérarchie du théorème d'ensemble fin, basée sur le nombre de couleurs est strictement décroissante en mathématiques à rebours (Théorème~\ref{thm:ts2-omega-models}) et que le théorème d'ensemble libre n'implique pas le théorème de Ramsey pour les paires (Corollaire~\ref{cor:fs-not-imply-rt22}). Cela répond en particulier à de nombreuses questions de Cholak, Guisto, Hirst, Jockusch~\cite{Cholak2001Free} et de Mont\'alban~\cite{Montalban2011Open}.
	\bigskip

	\item La séparation du théorème de Ramsey pour les paires du lemme faible de K\"onig est restée ouverte pendant des décennies, jusqu'à ce que Liu~\cite{Liu2012RT22} la résolve récemment. Nous étendons le théorème de Liu aux théorèmes d'ensemble libre, d'ensemble fin et d'arc-en-ciel de Ramsey, montrant ainsi qu'aucun d'entre eux n'implique le lemme doublement faible de K\"onig sous $\rca$. Ces questions avaient été posées par Hirschfeldt~\cite{Hirschfeldt2015Slicing} et leur réponse a fait intervenir une machinerie élaborée. Nous montrons également l'optimalité du théorème de Liu en prouvant que le théorème de Ramsey pour les paires n'évite pas les 1-énumerations d'ensembles clos (Corollaire~\ref{cor:sts2-not-simultaneous-cb-enum-avoidance}).
	\bigskip

	\item Nous montrons la faiblesse calculatoire du théorème d'Erd\H{o}s-Moser, en prouvant qu'il est borné par des degrés low${}_2$ (Section~\ref{sect:degrees-bounding-em-low2-degree}) et qu'il n'implique pas même le théorème de modèle atomique (Section~\ref{sect:dominating-the-erdos-moser-theorem}) en mathématiques à rebours. Ce dernier résultat renforce les séparations pré-existences de Lerman, Solomon et Towsner~\cite{Lerman2013Separating} et de Wang~\cite{Wang2014Definability}. En particulier, cela répond à une question d'Hirschfeldt, Shore et Slaman~\cite{Hirschfeldt2009atomic} en montrant que la cohesion n'implique pas le théorème de modèle atomique dans $\rca$. 
	\bigskip

	\item Nous développons de nouvelles notions de forcing avec de bonnes propriétés définitionnelles pour différents théorèmes de type Ramsey. L'existence de telles notions de forcing peut être vue comme un pas en avant vers la preuve que les hiérarchies des théorèmes d'ensemble libre, d'ensemble fin et d'arc-en-ciel de Ramsey basées sur la taille de leur $n$-uplets est stricte. La question de savoir si ces hiérarchies sont strictes a été posée par Cholak, Giusto, Hirst et Jockusch~\cite{Cholak2001Free} et reste ouverte à ce jour. En particulier, nous répondons positivement à deux conjectures de Wang~\cite{Wang2014Definability} à travers le Théorème~\ref{thm:coh-preservation-arithmetic-hierarchy} et le Théorème~\ref{thm:em-preserves-arithmetic}.
\end{itemize}

\fi

\ifdefined\mypublic
\chapter*{Introduction}
\addcontentsline{toc}{chapter}{\textcolor{ocre}{Introduction}}
\else
\chapter*{Introduction (English)}
\addcontentsline{toc}{chapter}{\textcolor{ocre}{Introduction (English)}}
\fi

This document is a thesis report resulting from three years of research under the supervision
of Laurent Bienvenu and Hugo Herbelin. We start with a gentle introduction to the mathematical fields
that we deal with here. We then give a detailed summary of the thesis and its main
contributions.

\section*{Computability theory and reverse mathematics}

This thesis contributes to the domain of reverse mathematics
under a computational perspective. We briefly introduce computability theory,
which provides a solid framework to talk about the notion of computable set.
We then provide a computationally-oriented introduction to reverse mathematics
with its main motivations and challenges.

\subsection*{Computability theory}

Informally, a set of natural numbers is \emph{computable} if one can effectively
decide which elements belong to it and which do not. In this sense, every finite
set is computable since one can learn by heart the list of its elements.
However, when dealing with infinite sets, one must find a systematic method (or \emph{algorithm})
to make such a decision. The Church-Turing thesis asserts that the notion of computable set
is robust, in that any reasonable model of computation defines the same class of computable sets.
One can therefore easily get convinced that some set is computable by writing a decision algorithm
in any mainstream programming language.

Although defined over natural numbers, the notion of computable set
extends to any data structure by using a simple coding technique. 
All the coding functions we shall consider are computable and
therefore one can talk transparently about a data structure or its corresponding set
of natural numbers. For example, an ordered pair of integers $(a, b)$ can be coded as the integer~$2^a3^b$.
Once we have a pairing function, we can define a coding of strings, and so on.

How large is the class of computable sets? 
Since each computable set comes with a decision algorithm,
there are only countably many of them. By a simple diagonal argument,
one can show that the collection of all sets is uncountable, and therefore
that the vast majority of sets of natural numbers is incomputable.
The most famous example of such an incomputable set is the \emph{halting problem}, defined as the 
set of all the (codes of) programs which halt.

Whenever a set is not computable, one may naturally wonder how much incomputable it is.
For example, would we be able to compute the set if we were able to solve the halting problem?
One can extend the notion of computable set by adding the characteristic function $\chi_Y$ of a set~$Y$
among the primitives of any programming language. The function~$\chi_Y$ cannot be computed, and must therefore be thought of 
as an \emph{oracle}. A set~$X$ which is decidable by an algorithm in this augmented language
is said to be \emph{computable relative to~$Y$}. This is usually abbreviated by~$X \leq_T Y$
(where T stands for Turing). If $X \leq_T Y$, then the set $Y$ is at least \emph{as hard to compute as} $X$.
One can define an infinite hierarchy of degrees of incomputability, in the same spirit as the hierarchy of infinities
revealed by Cantor.

The notion of computable set admits a purely logical characterization in terms of formulas
in arithmetic. Fix an infinite set of \emph{number variables} $x, y, z, \dots,$
which are intended to range over natural numbers.
We can form \emph{numerical expressions} by using them as placeholders
in some additions and multiplications of integers. For example, 
$(3+x) \cdot y$ is a numerical expression, which, once replaced~$x$ by~18 and~$y$ by~2,
denotes the number 42.
One can then define \emph{formulas} by relating two numerical expressions with the equality $=$
or the inequality $\leq$, using number quantifiers $(\exists x)$ (there exists $x$) and $(\forall x)$ (for all $x$),
bounded quantifiers $(\exists x \leq y)$, $(\forall x \leq y)$, 
and composing them with some logical connectives, such as $\wedge, \vee, \neg, \imp, \biimp$ (and, or, not, implies, if and only if).
For example, the formula $(\forall x)(\exists y \leq x)(x = 2y \vee x = 2y+1)$ is a formula saying that every number~$x$
is an even number or an odd number.

A formula $\varphi$ with one \emph{free} variable~$x$, that is, such that $x$ does not appear in a quantifier,
can be associated a set $S_\varphi$, defined as the set of all the natural number~$n$ such that
the formula $\varphi$ holds whenever~$x$ is replaced by~$n$. For example, the formula
$(\exists y)(x = 2y)$ defines the set~$S = \{ n : (\exists y)(n = 2y) \}$ which is nothing but
the set of even numbers. We can classify the complexity the formulas in function of their alternation of quantifiers.
A formula is $\Sigma^0_1$ (resp.\ $\Pi^0_1$) if it is of the form 
$(\exists x)\varphi$ (resp. $(\forall x)\varphi$), where $\varphi$ is a formula with only bounded quantifiers.
We say that a set is $\Sigma^0_1$ if it is defined by a $\Sigma^0_1$ formula. $\Pi^0_1$ sets are defined accordingly.
Intuitively, $\Sigma^0_1$ sets are those whose elements can be \emph{enumerated},
since one can try each value of~$x$ until we find one for which the corresponding formula $\varphi$ holds.
A set is $\Delta^0_1$ if it is both $\Sigma^0_1$ and $\Pi^0_1$.
A famous theorem from Emil Post asserts that the $\Delta^0_1$ sets are \emph{exactly} the computable sets.
This theorem is only a bootstrap of a whole correspondence between computability and definability by formulas.

We now have the necessary background in computability theory to introduce reverse mathematics.

\subsection*{Reverse mathematics}

An \emph{axiom} is a mathematical statement which is taken as a postulate.
Given a set of axioms, often called a \emph{theory}, and some \emph{deduction rules}, one can derive some mathematical consequences
called \emph{theorems}. From a purely mathematical point of view, an axiom and a theorem are
nothing but mathematical statements. The only difference is that an axiom is self-justified,
whereas a theorem must be given a proof.

It is reasonable to require a theory to be \emph{computably enumerable},
that is, such that the set of axioms is $\Sigma^0_1$. Indeed, one want to be able 
to know which axioms are available to prove a theorem. Thankfully, the vast majority
of formalizations of concepts yield computably enumerable theories.
This is in particular the case for set theory.
In his first incompleteness theorem, G\"odel asserts that every such theory which
is consistent and sufficiently expressive is \emph{incomplete}, that is,
have mathematical statement which are neither provable, nor refutable.
This theorem rules out the hope to find one day a computable set of axioms
which would enable us to decide any mathematical statement.
As a consequence, the foundational search for natural axioms
we could add to set theory is an unfinished process and still an active research subject.

Reverse mathematics is a vast mathematical program that seeks to determine 
which axioms are required to prove theorems of mathematics.
The underlying idea is very simple. Remove all the strong axioms from set theory
and keep only a very weak base theory, in which almost no theorem is provable.
This theory nevertheless has to be strong enough to formalize all the trivial coding
tricks we use, and therefore not to be depending on the presentation of the considered data structures.
Then choose an ordinary theorem~$\Psf$ that is used in everyday life, and try to prove it with the minimum
number of axioms. Sometimes, the theorem~$\Psf$ is directly provable in the base theory.
In this case, this means that the theorem does not have a very complex proof since we chose
a very weak base theory. When it is not the case, we need to add some axioms
$A_0, A_1, \dots, A_n$ to the base theory to prove~$\Psf$. At this stage, we know that
the axioms $A_0, A_1, \dots, A_n$ are \emph{sufficient} to prove~$\Psf$, but how can we ensure
that they are \emph{necessary}? Simply try to prove the reverse implication,
that is, add the theorem~$\Psf$ to the base theory, and try to prove $A_0 \wedge A_1 \wedge \dots \wedge A_n$.
If we succeed, then we just proved that, over the base theory, the axioms $A_0, A_1, \dots, A_n$
are the \emph{minimal} required to prove~$\Psf$.

The reverse mathematical program was founded by Harvey Friedman in 1975.
Since then, many researchers have contributed to this endeavor, 
the most notable contributor being Stephen Simpson.
Thousands of theorems have been studied from core mathematical areas,
such as algebra, analysis, topology, among others.
A surprising phenomenon emerged from the early years of
reverse mathematics: Most theorems studied require very weak axioms. 
Moreover, many of them happen to be \emph{equivalent} to one of five main sets of axioms,
that are referred to as the \emph{Big Five}. These five sets of axioms have been shown to correspond
to well-known philosophical approaches. There is an ongoing research to understand
the reasons of this phenomenon, and in particular which theorems escape the Big Five.

We now go into the details of the framework of reverse mathematics and the base theory.
For this, we need to extend our notion of formula to \emph{second-order arithmetic}. 
Take our previous language, augmented with some primitives to talk about sets of natural numbers.
We therefore fix a set of \emph{set variables} $X, Y, Z$, and extend our previous notion of formula
with ownership ($n \in X$) and quantifiers over sets $(\forall X)$, $(\exists X)$.
For example, the formula $(\exists X)(\forall x)[x \in X \biimp x < 5]$ asserts
that the set of natural numbers smaller than~5 exists.
Harvey Friedman noticed that the large majority of mathematical theorems could be naturally stated
in the language of second-order arithmetic. Therefore, one can restrict the reverse mathematical
program to mathematical statements from second-order arithmetic without loosing too much of generality.

In this setting, an axiom precises the behavior of the natural numbers
and the sets of natural numbers, that we refer to as the \emph{first-order part}
and the \emph{second-order part}, respectively. For example, the statement ``$0 \neq 1$''
is an axiom. The more axioms the base theory has, the less chaotic its models will be.
The base theory $\rca$, standing for Recursive Comprehension Axiom, contains
a set of first-order axioms known as \emph{Robinson arithmetic}. These axioms describe
the behavior of the natural numbers with respect to addition, multiplication
and order. By the compactness theorem, it is not possible to exclude all non-standard behaviors
of the natural numbers with only first-order axioms, but Robinson arithmetic
ensures that they will behave as expected with respect to the standard operations.

$\rca$ also contains second-order axioms which will describe which sets of natural numbers
are ensured to exist. The second-order axioms of~$\rca$ can be classified into two kinds:
the \emph{induction} axioms and the \emph{comprehension} axioms.
Given a formula~$\varphi(x)$ with one distinguished free number variable~$x$,
the induction axiom for~$\varphi$ asserts that if $\varphi(0)$ holds
and $\varphi(n+1)$ holds whenever $\varphi(n)$ holds, then $\varphi(n)$ holds
for every natural number~$n$. The comprehension axiom for~$\varphi$
asserts that the collection of the $n$'s 
such that~$\varphi(n)$ holds exists as a set. We write~$\{ n : \varphi(n) \}$
for this set. The induction and comprehension axioms can be seen as methods 
to build new sets from existing ones in the models.
$\rca$ contains the induction scheme for every $\Sigma^0_1$
formula, and the comprehension scheme for any $\Delta^0_1$ predicate.
By the correspondence between $\Delta^0_1$ predicates and computable sets, 
we allow only to build new sets computably from the previous ones. In this sense,
$\rca$ captures \emph{computable mathematics}.

Many theorems can be seen as \emph{mathematical problems},
coming with a natural class of \emph{instances} and of \emph{solutions} to each instance.
For example, K\"onig's lemma asserts for every infinite, finitely branching tree
the existence of an infinite path through the tree. Here, an instance
is an infinite, finitely branching tree~$T$, and a solution to~$T$
is an infinite path through~$T$. One can naturally define a notion of \emph{computable reduction}
from a problem~$\Psf$ to another problem~$\Qsf$ (written $\Psf \leq_c \Qsf$) by saying that
every instance $X_0$ of~$\Psf$ computes an instance~$X_1$ of~$\Qsf$ such that every solution to~$X_1$
computes (with the help of $X_0$) a solution to~$X_0$. Informally,
$\Psf \leq_c \Qsf$ if we can computably solve the mathematical problem~$\Psf$ from
a $\Qsf$-solver.
A computable reduction from~$\Psf$ to~$\Qsf$ can be seen as a proof that $\Qsf \imp \Psf$ over~$\rca$
where only one use of~$\Qsf$ is allowed. In this sense, computable reduction is more \emph{precise}
than reverse mathematics. In particular, it enables one to reveal some fine distinctions
between mathematical statements which would be indistinguishable from the point of view of reverse mathematics.

\section*{Thesis summary}

We start this section explaining more accurately the content of this thesis. We then
develop those explanations chapter by chapter and we end with a summary of the main
original contribution of the thesis.

\subsection*{The subject}

This thesis is about the reverse mathematics of mathematical statements
coming from Ramsey theory. Ramsey theory is a branch of mathematics 
studying the conditions under which some
structure appears among a sufficiently large collection of objects.
For example, in any group of six people,
there are at least three people who all know each other or three people who
are all strangers to each other. 
Perhaps the most famous such theorem is \emph{Ramsey's theorem},
which asserts that for every $k$-coloring of~$[\Nb]^n$ (where $[\Nb]^n$ denotes the $n$-tuples
of natural numbers), there is an infinite set of natural numbers $H$ such that~$[H]^n$ is given only one color.
Such a set is called \emph{homogeneous}. In particular, Ramsey's theorem for singletons ($n = 1$)
is nothing but the infinite \emph{pigeonhole principle}, which asserts that if an infinite number of items
is spread in $k$ boxes, at least one of the boxes must contain infinitely many items.

In the past two decades, Ramsey theory emerged as one of the most important topics in reverse mathematics.
This theory provides a large class of theorems escaping the Big Five phenomenon, and whose strength is notoriously hard
to gauge. We conduct a systematic study of various consequences of Ramsey's theorem, namely,
the Erd\H{o}s-Moser theorem, the free set, thin set and rainbow Ramsey theorems, among others,
within the frameworks of reverse mathematics and of computable reducibility. 
Moreover, we relate them to several existing theorems studied in reverse mathematics.

\subsection*{Structure of this thesis}

This thesis is divided into 4 parts, namely, a preamble, the strength of Ramsey's theorem,
further topics, and a conclusion. We briefly explain the content of each part.
\bigskip

In the preamble, we introduce the basic concepts about computability theory
that we will deal with. Then, we provide a detailed introduction to reverse mathematics
and the Big Five phenomenon. We then relate reverse mathematics to various recently introduced reducibility notions
which focus on some computational aspects of the theorems we consider.
Last, we introduce the main separation technique that is used all along this thesis, namely,
effective forcing. This technique is a refinement of Cohen's notion of forcing 
in order to preserve some computability-theoretic properties.

In the second part, we conduct a systematic study of Ramsey's theorem
and various of its consequences under reverse mathematics and computable reducibility.
We mainly consider well-established theorems which have already been introduced in reverse mathematics,
namely, cohesiveness, stable Ramsey's theorem for pairs, the Erd\H{o}s-Moser theorem,
the free set, the thin set and the rainbow Ramsey theorem, among others.
We provide simpler proofs of existing separations with a uniform and minimalistic framework,
and answer several open questions.

The third part is more thematic, that is, we study various existing theorems in reverse mathematics
under some specific computational perspective. First, we study the existence of maximally difficult instances
(\emph{universal instances}) of Ramsey-type theorems and the Turing degrees bounding their solutions.
Then, we show the orthogonality of some Ramsey-type theorems with respect to compactness.
We also identify a new class of theorems strengthening Ramsey's theorem for pairs and whose
computational strength is currently unknown. Last, we design new notions of forcing with better definability
properties, with the hope to prove the strictness of the free set, thin set and rainbow Ramsey theorem hierarchies
in reverse mathematics.

In the conclusion, we discuss the naturality of the tools used to study Ramsey-type statements
in reverse mathematics, and try to give some perspective 
on the role of Ramsey-type theorems with respect to the Big Five phenomenon.
We then provide a list of remaining open questions
with their motivation and some insights about the reason of their difficulty.

\subsection*{Main contributions}

Along this thesis, we answer 25 open questions from 17 different papers.
We detail here the main original contributions of this thesis.
\bigskip

\begin{itemize}

	\item The main contribution is probably methodological.
	According to the famous proverb ``\emph{give a man a fish and you feed him for a day; 
	teach a man to fish and you feed him for a lifetime}'', we are more interested in the development
	of tools to separate systematically theorems over reverse mathematics than
	to actual separations. In~\cite{Patey2015Iterative}, we presented a simplification of a framework
	introduced by Lerman, Solomon and Towsner~\cite{Lerman2013Separating}.
	We show that this framework is sufficiently general reprove various separations with a simpler
	argument (Corollary~\ref{cor:emo-wkl-sts2-sads}, Corollary~\ref{cor:ads-not-implies-scac}) 
	and separate Ramsey's theorem for pairs
	from the tree theorem for pairs (Theorem~\ref{thm:tree-theorem-rt22+wkl-not-tt22}), 
	thereby answering a question of Mont\'alban~\cite{Montalban2011Open}.
	\bigskip

	\item Another important contribution is the separation
	of Ramsey's theorem for $k+1$ colors from Ramsey's theorem for $k$ colors
	over computable reducibility (Theorem~\ref{thm:colors-ramsey-general}).
	This question was first asked in a weaker form by Mileti~\cite{Mileti2004Partition},
	and then for Weihrauch reducibility by Dorais, Dzhafarov, Hirst, Mileti and Shafer~\cite{Dorais2016uniform}
	who left it open as a ``chief question''.
	We answer both questions. Hirschfeldt and Jockusch~\cite{Hirschfeldtnotions} recently asked how many
	applications of Ramsey's theorem for~$m$ colors are necessary to solve Ramsey's theorem
	for $k$ colors in function of~$m$ and~$k$. We give a precise answer in Theorem~\ref{thm:colors-ramsey-general}.
	\bigskip

	\item We reveal that the hierarchy of the thin set theorem based on the number of colors is strictly decreasing
	over reverse mathematics (Theorem~\ref{thm:ts2-omega-models}) and that the free set theorem does not imply
	Ramsey's theorem for pairs (Corollary~\ref{cor:fs-not-imply-rt22}). This answers several questions of Cholak, Giusto,
	Hirst and Jockusch~\cite{Cholak2001Free} and of Mont\'alban~\cite{Montalban2011Open}.
	\bigskip

	\item The separation of Ramsey's theorem for pairs from weak K\"onig's lemma was open for decades,
	until Liu~\cite{Liu2012RT22} recently solved it. We extend Liu's theorem to the free set, thin set and the rainbow Ramsey theorem,
	therefore showing that none of them imply weak weak K\"onig's lemma over~$\rca$.
	These questions were asked by Hirschfeldt~\cite{Hirschfeldt2015Slicing} 
	and their answers require an involved machinery.
	We also show the optimality of Liu's theorem by proving that Ramsey's theorem for pairs does not
	avoid 1-enumerations of closed sets (Corollary~\ref{cor:sts2-not-simultaneous-cb-enum-avoidance}).
	\bigskip

	\item We show the computational weakness of the Erd\H{o}s-Moser theorem
	by proving that it admits a low${}_2$ bounding degree (Section~\ref{sect:degrees-bounding-em-low2-degree})
	and that it does not even imply the atomic model theorem (Section~\ref{sect:dominating-the-erdos-moser-theorem}) 
	in reverse mathematics.
	The later result strengthens the previous separations of Lerman, Solomon and Towsner~\cite{Lerman2013Separating}
	and of Wang~\cite{Wang2014Definability}. In particular, it answers a question of 
	Hirschfeldt, Shore and Slaman~\cite{Hirschfeldt2009atomic} by showing that cohesiveness
	does not imply the atomic model theorem over~$\rca$.
	\bigskip

	\item We develop new notions of forcing with good definitional properties
	for various Ramsey-type theorems. The existence of such notions of forcing
	can be seen as a step towards the resolution of the strictness of the free set, thin set,
	and rainbow Ramsey theorem hierarchies based on the number of their tuples.
	The question of the strictness of those hierarchies is asked by Cholak, Giusto,
	Hirst and Jockusch~\cite{Cholak2001Free} and remains open. 
	In particular, we answer positively two conjectures of Wang~\cite{Wang2014Definability} 
	through Theorem~\ref{thm:coh-preservation-arithmetic-hierarchy} and Theorem~\ref{thm:em-preserves-arithmetic}.
	
\end{itemize}

\chapter*{A coffee break}
\addcontentsline{toc}{chapter}{\textcolor{ocre}{A coffee break}}

Everyone knows that a mathematician transforms 
coffee\footnote{The choice of coffee is not crucial in this explanation. You may even replace
it with your favorite beverage: tea, chocolate, orange juice, ...
Note that the notion of beverage has to be defined carefully:
hemlock is a beverage, but it can be drunk only once.}
into theorems. This fact has been empirically verified since ever, and in particular
non vacuously since coffee has been invented.
However, coffee is an expensive ressource, especially in those
troubled times of economic crisis. One may therefore
naturally wonder whether coffee is really necessary for a mathematician
to produce theorems. Water may be sufficient, and much cheaper.

How can we ensure that coffee is optimal for producing theorems?
An attempt would consist of taking a theorem and, using the mathematician,
trying to extract some coffee from it.
Assuming Lavoisier's Law of Conservation of Mass,
coffee cannot be created \emph{ex nihilo}.
In case of success, we can therefore reasonably deduce that the coffee was extracted
from the theorem, and that we \emph{really} need coffee to make a theorem.

We must however be careful with the choice of the mathematician.
The insights we may get from our coffee extraction strongly depends on his skills.
If the mathematician can brew coffee, we cannot establish any link
between coffee and the theorem, in that he may not have used the theorem
in his recipe of coffee\footnote{Beware if you order coffee in touristic places:
rumor has it that some locals replace theorems by cheap lemmas
in the preparation of coffee, in order to increase their margins.}.
On the other hand, a too cautious mathematician will require to drink
only a specific brand of coffee to produce his theorems. 
We do not want to be bothered by those distinctions, and want our insights
to be independent of the precise definition of ``coffee''.
We will therefore choose a mathematician who is able to perform
only basic operations so that we can establish a link between what we provide him
and what he produces, but who is sufficiently skilled not to depend on
the presentation of the input.

Assuming coffee is optimal for the production of a theorem,
how many cups do we need to prove a theorem?
Can the same mathematician produce every theorem if we give him
a sufficient amount of coffee, or do we need the whole
scientific community?
The author tries to answer those questions
for a particular class of theorems coming from Ramsey's theory.
This theory informally states that we can always find some structure
out of a sufficiently large amount of disorder. By a quirk of fate, Ramsey's theory 
introduces disorder in the structured world that used to be reverse mathematics.

In this thesis, you may find action, suspence, hope, disappointment,
hope again, but above all, you will find... theorems.
Some of the presented results are trivial, some are rather technical. 
The author just hopes you will find his theorems tasty enough to prepare a good coffee...

\part{Preamble}

\chapter{Computability theory}

\index{computable!set}
\index{computable!function}
Many physical phenomena can be apparented to computation
and used to solve mathematical problems.
Their analysis often leads to new paradigms of computation,
with the dual hope to increase our understanding of the physical world, 
and to obtain more efficient and maybe \emph{more expressive} computational models.
DNA synthesis and quantum phenomena are case in point examples
of such approaches.

However, among all the computational behaviors found in nature hitherto,
no one has been shown to be more expressive than a computer.
Is this remark an empirical observation, or is there a philosophical justification
to the unability to find more expressive models of computation?
More generally, what kind of mathematical problems can we solve by
physical means? Is there a natural class of functions
that we can reasonably consider as \emph{effectively computable}?

\section{The Church-Turing thesis}

The early 1930's have seen the emergence of three equivalent models of computation discovered 
independently. Herbrand and G\"odel~\cite{Herbrand1931consistency} proposed a system of equations with an evaluation
mechanism in order to remedy the lack of expressive power of primitive recursive functions.
On his side, Church introduced the $\lambda$-calculus, a calculus of functions, and formulated
the thesis that his calculus captured exactly the effectively computable functions.
However, it was hard at the time to be convinced without some manipulation that $\lambda$-calculus
was expressive enough to capture every computational process. It is only
when Turing introduced his now-called Turing machines~\cite{Turing1936computable} in a seminal paper
providing all the evidences of the naturality of his model of computation
that Church's thesis became universally accepted.

\begin{thesis}[Church-Turing thesis]
Every effectively computable function is recursive.
\end{thesis}

Apart from its philosophical meaning, the Church-Turing thesis has practical consequences
in the development of computability theory. It is a common practice, when proving the existence
of a Turing machine computing some property, to simply provide an informal description
of how to compute this property, and use the Church-Turing thesis to deduce that
the property is indeed computable in a formal sense.
For this reason, we shall remain voluntarily vague concerning the choice of the computational model over which
we will build computability theory.

Fix any \emph{reasonable} programming language, and define a \emph{G\"odel numbering}
of all programs, that is, an effective coding of the programs into natural numbers.
By reasonable, we mean any programming language which can be converted to and from 
a Turing machine. In particular, mainstream programming language such as C, Java
and PHP are reasonable.

\index{p.c.|see {partial computable}}
\index{partial computable}

\begin{definition}[Partial computable function]
The \emph{$e$th partial computable (p.c.) function}~$\Phi_e$ is defined by~$P$ if $e$ is the G\"odel number of a program~$P$,
and is the empty program otherwise.
\end{definition}

\index{Turing machine}
\index{$W_e$}
\index{$\Phi_e$}
We call the number~$e$ the \emph{Turing index}, or simply the~\emph{index} of
the p.c.\ function~$\Phi_e$. Since Turing machines provide a suited computational model
to analyse the complexity of algorithms in terms of time and space, they are usually adopted
as the underlying model supporting computability theory. For this reason, $\Phi_e$ is often called a \emph{Turing machine}
and partial computable functions inherit the standard terminology of Turing machines.
In particular, we say that~$\Phi_e$ \emph{diverges on~$x$} and write~$\Phi_e(x) \uparrow$
if~$x$ is not in the domain of~$\Phi_e$. If~$\Phi_e$ does not diverge on~$x$, then it \emph{converges} to the value~$\Phi_e(x)$ 
and we write~$\Phi_e(x) \downarrow$.
We use~$\Phi_e(x) \simeq \Phi_i(y)$ to say that either both functions diverge, or they both converge to the same value.
We write~$W_e$ for the domain of the $e$th partial computable function~$\Phi_e$.
Such a set is called \emph{computably enumerable} as its elements
can be enumerated by executing the program of~$\Phi_e$ and adding an element~$x$
whenever the program halts on the input~$x$.

We also want to define what it means for a set to be computable.
One may think of a set of integers~$X$ as a \emph{decision problem}
whose questions are of the form~``is~$n$ in~$X$?''.

\index{decidable}

\begin{definition}
A set of integers~$X$ is \emph{decidable} if its characteristic function is total computable.
\end{definition}

For any reasonable programming language, the corresponding G\"odel numbering
is \emph{admissible} in the sense of Rogers~\cite{Rogers1967Theory}, and the listing
of partial computable functions satisfies the following standard lemmas and theorems.

\begin{lemma}[Padding Lemma]
For each index~$e$ and threshold~$m$, one may effectively obtain some index~$i > m$
such that~$\Phi_i \simeq \Phi_e$.
\end{lemma}
\begin{proof}[Proof idea]
Add to the program of G\"odel number~$e$ some useless instructions
to artificially increase its G\"odel number.
\end{proof}

Altough used in many arguments, 
the padding lemma should be seen as a curse rather than
a desired property. Indeed, any listing can be transformed into
one satisfying the padding lemma by using the bijection
between $\Nb \times \Nb$ and $\Nb$. What the padding lemma really tells us is that
there is no way to avoid this property, and in particular to
make an effective listing of all partial computable functions without redundancy.

The following theorem is due to Kleene~\cite{Kleene1938notation}.
It is an integrated feature of some alternative models of computation
such as the $\lambda$-calculus~\cite{Barendregt1984Lambda}.

\index{s-n-m theorem}
\begin{theorem}[S-m-n theorem]
For each partial computable function~$f(e, x)$,
one may effectively obtain a computable function~$q$
such that~$\Phi_{q(e)}(x) \simeq f(e,x)$ for each~$x$.
\end{theorem}
\begin{proof}[Proof idea]
Hardcode the parameter~$e$ into the program corresponding to the function~$f$. 
\end{proof}

The S-m-n theorem informally states that given a partial computable function,
one can fix some parameters and transform it into an effective listing 
of partial computable functions of the remaining parameters.

The last theorem about the enumeration of all partial computable functions,
also due to Kleene~\cite{Kleene1938notation}, is certainly the most difficult to interpret.

\index{Kleene's recursion theorem}
\begin{theorem}[Kleene's recursion theorem]
For any partial computable function~$f(e, x)$,
there is an index~$e$ such that~$\Phi_e(x) \simeq f(e,x)$ for each~$x$.
\end{theorem}

Kleene's recursion theorem informally states that 
we can define self-referential algorithms, that is,
algorithms in which we assume we are provided their Turing index.
Altough simple, the proof of Kleene's recursion theorem
is relatively obscure, and one generally prefers
to avoid using this theorem when unecessary.

\section{Definability and the arithmetic hierarchy}

Since there are uncountably many sets of integers, but only countably many programs,
some sets of integers are non-computable. Turing~\cite{Turing1936computable} gave one of the first examples
of unsolvable decision problems.

\index{halting problem}
\begin{theorem}[Halting problem]\label{thm:intro-comp-halting-problem}
The \emph{halting set} $K = \{e : \Phi_e(e) \downarrow \}$ is undecidable.
\end{theorem}
\begin{proof}[Proof idea]
Let~$\Phi_e(n) = 1 - \chi_K(n)$. Is~$e$ in~$K$?
\end{proof}

The existence of non-computable sets is a starting point to define the notion of \emph{relative computability}.
Given a set~$X$, we can add to our programming language the characteristic function of~$X$ as a primitive of the language.
We then denote by~$\Phi^X_0, \Phi^X_1, \dots$ the enumeration of the $X$-partial computable functions, that is,
the partial computable functions which can be computed with the help of the \emph{oracle}~$X$.
The notion of~$X$-c.e.\ set is defined accordingly and we denote by~$W^X_e$ the domain of~$\Phi^X_e$.

\index{$\leq_T$}
\index{$\equiv_T$}
\index{reduction!Turing}
\index{Turing!reduction}
\index{Turing!degree}
\begin{definition}[Turing reducibility]
A set~$X \subseteq \Nb$ \emph{computes} a set~$Y$ (written~$Y \leq_T X$) if $Y = \Phi^X_e$ for some index~$e$.
We write~$X \equiv_T Y$ to say that~$X \leq_T Y$ and~$Y \leq_T X$.
The \emph{Turing degree}~$\deg(X)$ is the collection~$\{ Y \subseteq \Nb : Y \equiv_T X \}$.
\end{definition}

The notion of Turing reducibility extends to Turing degrees. We denote the Turing degrees
by bold lower case letters~$\dbf, \ebf, \dots$ and call~$\mathbf{0}$ the degree of the computable sets.
Letting the \emph{effective join} $X \oplus Y$ be the set~$\{2n : n \in X\} \cup \{2n+1 : n \in Y\}$,
two Turing degrees~$\dbf = \deg(X)$ and~$\ebf=\deg(Y)$ have a least upper-bound~$\dbf \cup \ebf = \deg(X \oplus Y)$.

Theorem~\ref{thm:intro-comp-halting-problem} tells us that
$K \not \leq_T \emptyset$. The halting problem can be relativized to any oracle, and leads to the notion of Turing jump
which plays a central role in computability theory.

\index{Turing!jump}
\index{jump operator}
\begin{definition}[Jump operator]
The \emph{jump} of a set $X$ is the set~$X' = \{e : \Phi^X_e(e) \downarrow\}$.
\end{definition}

The Turing jump can be iterated as follows: Given a set~$X$, $X^{(0)} = X$ and~$X^{(n+1)} = (X^{(n)})'$.
In particular, $\emptyset' = K$. From now on, we drop the notation~$K$ for the halting set, and simply use~$\emptyset'$.

\subsection{The arithmetic hierarchy}

According to Hirschfeldt~\cite{Hirschfeldt2015Slicing}, computability theory
is above all about the relationship between computation and definability.
Post's theorem (see below) is one of the core theorems establishing this correspondance.
It provides bridges between computability and provability.

\index{arithmetic hierarchy}
\begin{definition}[Arithmetic hierarchy]
A set $X$ is~$\Sigma^0_n$ if there is a computable relation~$R(x_1,\dots, x_n, y)$ such that
$y \in X$ iff
\[
\underbrace{\exists x_1 \forall x_2 \exists x_3 \forall x_4 \dots Q_nx_n}_{n \mbox{ alterning quantifiers}}R(x1_1, \dots, x_n, y)
\]
where~$Q_n$ is $\exists$ if~$n$ is odd and~$\forall$ if~$n$ is even.
The definition of~$X$ being $\Pi^0_n$ is the same, except that the leading quantifier is a~$\forall$.
\end{definition}

It is easy to see that a set~$X$ is $\Sigma^0_n$ iff $\overline{X}$ is $\Pi^0_n$. A set~$X$
is~$\Delta^0_n$ if it is both~$\Sigma^0_n$ and~$\Pi^0_n$.

\index{Post theorem}
\begin{theorem}[Post~\cite{Post1948Degrees}]
A set is~$\Sigma^0_{n+1}$ iff it is $\emptyset^{(n)}$-c.e., and is $\Delta^0_{n+1}$
iff it is $\emptyset^{(n)}$-computable.
\end{theorem}

In particular, $\Delta^0_1$ sets are the computable sets, $\Sigma^0_1$ sets are the c.e.\ sets
and~$\Delta^0_2$ sets are the sets computable from the halting problem.

Given an $X$-partial computable function~$\Phi^X_e$, we can define
the~$s$th approximation~$\Phi^X_e[s]$ of~$\Phi^X_e$ as the partial computable function
obtained by running the program of~$\Phi^X_e$ in at most $s$ steps. If within this time,
the program halts on some input~$x$, then~$\Phi^X_e[s](x) \downarrow = \Phi^X_e(x)$.
Otherwise, $\Phi^X_e[s](x) \uparrow$. We denote by~$W^X_e[s]$ the domain of~$\Phi^X_e[s]$.
The set~$W^X_e[s]$ is computable uniformly in~$e$ and~$s$, and $W^X_e[s] \subseteq W^X_e$
for each~$s$.
Shoenfield~\cite{Shoenfield1959degrees} generalized the notion of approximation to any $\Delta^0_2$ set.

\index{Shoenfield's limit lemma}
\begin{lemma}[Shoenfield's limit lemma]
A set~$X$ is $\Delta^0_2$ iff there is a computable~$\{0,1\}$-valued function~$g$
such that~$\lim_s g(n, s)$ exists and~$\lim_s g(n, s) = \chi_X(n)$ for all~$n$.
\end{lemma}

The function~$g$ is called a~\emph{$\Delta^0_2$ approximation} of~$X$.
When~$X$ is~$\Delta^0_2$ we usually fix a $\Delta^0_2$ approximation~$g$
and write~$X[s] = \{n \leq s : g(n, s) = 1 \}$.

\subsection{Low and high degrees}

If~$X \leq_T Y$, then~$X' \leq_T Y'$. However, it is not true that
if~$X <_T Y$, then~$X' <_T Y'$.
In particular, there are some non-computable sets~$X$ such that~$X' =_T \emptyset'$.

\index{low set}
\index{high set}
\index{lown set@low${}_n$ set}
\begin{definition}
A set~$X$ is \emph{low over~$Y$} if~$(X \oplus Y)' \leq_T Y'$.
A set~$X$ is \emph{high over~$Y$} if~$Y'' \leq_T = (X \oplus Y)'$.
A set~$X$ is \emph{low} (\emph{high}) if it is low over~$\emptyset$ (high over~$\emptyset$).
\end{definition}

Intuitively, a set is low if it is indistinguishable from a computable
set from the point of view of the halting set. A set is high
if its jump is maximally difficult to decide. In particular,
every low set is~$\Delta^0_2$. By the previous assertion, there exist some non-computable low sets.
The low sets are downward-closed under the Turing reducibility.
Note that lowness is a degree-theoretic property.
An important feature of lowness is that it is preserved under relativization.

\begin{lemma}
If~$X$ is low over a low set~$Y$, then~$X$ is low.
\end{lemma}
\begin{proof}
$X' \leq_T (X \oplus Y)' \leq_T Y' \leq \emptyset'$ 
\end{proof}

The notions of lowness and highness can be generalized to arbitrary jumps.
Thus, a set~$X$ is \emph{low${}_n$ over~$Y$} if~$(X \oplus Y)^{(n)} \leq_T Y^{(n)}$.
The notion of high${}_n$ness is defined accordingly.
The high degrees admit a simple characterization that we shall present in the next section.

\section{Domination and hyperimmunity}

There are many ways to encode some non-computable information
into a set~$X$. One way consists in using the sparsity of~$X$
to compute some \emph{fast-growing function}.
The~\emph{principal function $p_X$} of a set~$X = \{n_0 < n_1 < \dots \}$
is the function defined by~$p_X(i) = n_i$. The sparser 
the set~$X$ is, the faster its principal function will grow.

\index{modulus function}
\begin{example}
The \emph{modulus} function~$m_e$ of a c.e.\ set~$W_e$
is defined for each~$n$ as the least stage~$s$ such that~$(W_e \cap [0, n)) \setminus W_e[s] = \emptyset$.
Any function growing faster than the modulus function of~$W_e$ computes~$W_e$.
In particular, any sufficiently sparse set computes the halting set.
\end{example}

The notions of domination and hyperimmunity are suitable tools to study the sparsity of a set.
They will naturally appear in many constructions over this thesis.

\index{domination}
\index{dominating function}
\index{hyperimmune!set}
\index{hyperimmune!function}
\begin{definition}[Domination, hyperimmunity]
A function~$f$ \emph{dominates} a function~$g$
if $f(x) \geq g(x)$ for \emph{all but finitely many~$x$}.
A function~$f$ is \emph{hyperimmune} if it is not dominated by any computable function.
\end{definition}

A set~$X$ is \emph{hyperimmune} if its principal function~$p_A$ is hyperimmune.
One must think of hyperimmune sets as having large holes from time to time.
This intuition can be made formal by the following characterization of hyperimmune sets.

\begin{lemma}
Fix a computable list of all finite sets~$D_0, D_1, \dots$
A set~$X$ is hyperimmune iff for every computable function~$f$ such that the sets~$D_{f(0)}, D_{f(1)}, \dots$ are mutually
disjoint, there is some~$i$ such that~$X \cap D_{f(i)} = \emptyset$.
\end{lemma}

A Turing degree~$\dbf$ is hyperimmune if it contains a hyperimmune function,
otherwise $\dbf$ is \emph{hyperimmune-free}.
Every computable function~$f$ is dominated by the computable function~$g(x) = f(x)+1$.
Therefore, there is no computable hyperimmune function and~$\mathbf{0}$ is a hyperimmune-free degree.
Miller and Martin~\cite{Miller1968degrees} showed the existence of a non-zero hyperimmune-free degree. However,
there is no such $\Delta^0_2$ degree by the following theorem.

\begin{theorem}[Miller and Martin~\cite{Miller1968degrees}]
Every non-zero~$\Delta^0_2$ degree is hyperimmune.
\end{theorem}

We shall see in section~\ref{sect:intro-forcing-effective-forcing} that hyperimmunity is a weak form
of genericity. The following theorem shows that highness is a dual notion of hyperimmunity.

\begin{theorem}[Martin~\cite{Martin1966Classes}]
A set~$X$ is high iff it computes a function dominating every computable functions.
\end{theorem}

\section{PA degrees}

In this section, we are interested in the following problem:
Given an infinite, finitely branching tree, how complicated is it 
to find an infinite path through the tree?
In order to make precise the definition of a tree and of a path,
we need to introduce some notation.

\index{string}
\index{sequence}
\index{real}
\index{binary string}
Fix a set~$X$.
A \emph{string} over~$X$ is an ordered tuple of integers $a_0, \dots, a_{n-1} \in X$.
The empty string is written~$\varepsilon$. A \emph{sequence} over~$X$
is an infinite listing of integers $a_0, a_1, \dots \in X$.
The sets $X^s$, $X^{<s}$ and~$X^{<\Nb}$ are the sets of strings of length~$s$, strictly smaller than~$s$
and finite length, respectively. When $X$ is not specified, we shall take~$X = \Nb$.
A \emph{binary string} is a string in~$2^{<\Nb}$ and a \emph{real} is a sequence in~$2^\Nb$.
We often identify a real~$X$ with the set of integers~$\{n : X(n) = 1\}$.

\index{tree}
\index{tree!binary}
\index{path}
\begin{definition}[Tree, path]
A tree $T \subseteq \Nb^{<\Nb}$ is a set downward-closed under the prefix relation.
The tree~$T$ is \emph{finitely branching} if every node~$\sigma \in T$
has finitely many immediate successors.
A \emph{binary} tree is a tree~$T \subseteq 2^{<\Nb}$.
A sequence $P$ is a \emph{path} through~$T$ if for every $\sigma \prec P$, that is, every initial segment~$\sigma$ of~$P$,
$\sigma \in T$.
\end{definition}

\index{Konig's lemma@K\"onig's lemma}
\index{weak Konig's lemma@weak K\"onig's lemma}
Given a tree~$T$, we denote by~$[T]$ the collection of its infinite paths.
K\"onig's lemma asserts that every infinite, finitely branching tree
has an infinite path. We call weak K\"onig's lemma the restriction of K\"onig's lemma
to infinite binary trees.

First, notice that the complexity of a path through a tree depends on the effectiveness of the tree~$T$ itself,
since every real~$X$ computes the infinite binary tree~$T_X = \{ \sigma \in 2^{<\Nb} : \sigma \prec X \}$
whose unique path is~$X$. We shall therefore restrict ourselves to computable trees.

The standard proof of K\"onig's lemma is the following: Given an infinite, finitely
branching tree, since the root has finitely many immediate successors, at least one of the successors
induces an infinite subtree. Choose one of them, and apply the same procedure on the induced subtree.
The standard proof of K\"onig's lemma is not effective, in that it requires to decide which node induce an infinite subtree.
Kreisel showed that there is no effective proof of weak K\"onig's lemma in the following sense.

\begin{theorem}[Kreisel~\cite{Kreisel1953variant}]
There is an infinite binary tree with no computable path.
\end{theorem}

For now, we shall restrict ourselves to binary trees.

\subsection{$\Pi^0_1$ classes}

\index{pi01 class@$\Pi^0_1$ class}
\index{class!$\Pi^0_1$}
A $\Pi^0_1$ class is a collection of sets of the form~$\{X \in 2^\Nb : (\forall n)\varphi(X \uh n)\}$
for some $\Pi^0_1$ formula~$\varphi$. It is easy to see that for every $\Pi^0_1$ class $\Ccal$,
there is a binary tree~$T$ such that~$\Ccal = [T]$. A $\Pi^0_1$ class is \emph{non-empty}
if it is the collection of paths through an infinite binary tree.

Weak K\"onig's lemma can be seen as a collection of problems
parameterized by their tree. One may question the relevance of studying a class of problems
as a single statement, since the complexity of finding a path through a computable tree
depends on the tree. However, there is a computable binary tree of \emph{maximal complexity},
in the following sense.

\index{class!universal}
\begin{definition}[Universality]
A non-empty $\Pi^0_1$ class~$\Ccal$ is \emph{universal} if for every non-empty
$\Pi^0_1$ class~$\Dcal$, every member of~$\Ccal$ computes some member of~$\Dcal$.
\end{definition}

Intuitively, finding a member of a universal $\Pi^0_1$ class is as hard as 
finding a member of any non-empty $\Pi^0_1$ class.
There are many natural examples of universal $\Pi^0_1$ classes.
By G\"odel's second incompleteness theorem~\cite{Goedel1931Uber}, 
there is no computable completion of Peano arithmetic.
Computability theory gives us a more precise understanding of
how hard it is to find such a completion.

\begin{theorem}[Scott~\cite{Scott1962Algebras}]
The $\Pi^0_1$ class of completions of Peano arithmetic is universal.
\end{theorem}

We call~\emph{PA} the degrees of completion of Peano arithmetic.
By Solovay [unpublished] (see~\cite{Downey2010Algorithmic}), the PA degrees are closed upward under the Turing reduction.
The notion of PA degree can be relativized to any set~$X$.
We write~$P \gg X$ to say that~$P$ is of PA degree relative to~$X$,
that is, if every non-empty $\Pi^{0,X}_1$ class has a $P$-computable member.
In particular, $P$ computes~$X$ as witnessed by the $X$-computable tree~$T_X$.

\subsection{Basis theorems}

The common structure of $\Pi^0_1$ classes
enables one to deduce general upper bounds on the complexity of 
finding a member in them. A \emph{basis} for $\Pi^0_1$ classes
states that each non-empty $\Pi^0_1$ class has a member satisfying some property.
In particular, there is a PA degree satisfying this property.
The three basis theorems presented in this section
are proven by Jockusch and Soare~\cite{Jockusch197201}.

The first basis theorem surprisingly states
that if some set~$X$ is non-computable,
being able to compute a completion of Peano arithmetic
is of no help to compute~$X$.

\index{basis theorem!cone avoidance}
\index{cone avoidance|see {basis theorem}}
\begin{theorem}[Cone avoidance basis theorem]\label{thm:intro-comp-cone-avoidance-basis}
If~$X$ is non-computable, every non-empty~$\Pi^0_1$ class
has a member which does not compute~$X$.
\end{theorem}

The cone avoidance basis theorem should not be interpreted
as ``the completions of Peano arithmetic carry no incomputable information''.
In fact, this extra information is of behavioral nature, as shows
the following lemma generalizing~\cite[Lemma 4.2]{Cholak2001strength}.

\begin{lemma}[Folklore]\label{lem:intro-comp-pa-choose}
Fix two sets~$P$ and~$X$. Then $P \gg X$ iff 
for every uniform sequence of pairs of $\Pi^{0,X}_1$ statements~$(\gamma_{e,0}, \gamma_{e,1})_{e \in \Nb}$
there is a $P$-computable function~$f$ such that~$\gamma_{e,f(i)}$ is true
whenever $\gamma_{e,0}$ or~$\gamma_{e,1}$ is true.
\end{lemma}

Lemma~\ref{lem:intro-comp-pa-choose} is another core
example of the correspondance between computability and definability.
The next basis theorem is, according to Cenzer~\cite{Cenzer199901},
\emph{perhaps the most cited result in the theory of $\Pi^0_1$ classes}.
It shows that we can always find a member which is indistiguishable
from the computable sets from the point of view of the halting set.

\index{basis theorem!low}
\index{low basis theorem|see {basis theorem}}
\begin{theorem}[Low basis theorem]
Every non-empty~$\Pi^0_1$ class has a member of low degree.
\end{theorem}

Last, every non-empty $\Pi^0_1$ class has a member~$P$ such that
the $P$-computable functions are all dominated by some computable functions.

\index{basis theorem!hyperimmune-free}
\index{hyperimmune-free basis theorem|see {basis theorem}}
\begin{theorem}[Hyperimmune-free basis theorem]
Every non-empty~$\Pi^0_1$ class
has a member of hyperimmune-free degree.
\end{theorem}

As noted by Hirschfeldt~\cite{Hirschfeldt2015Slicing},
the cone avoidance basis theorem is a consequence of the low basis and the hyperimmune-free
basis theorems. Indeed, given a non-computable set~$X$ and a non-empty $\Pi^0_1$ class~$\Ccal$,
either $X$ is $\Delta^0_2$, in which case no hyperimmune-free degree bounds it,
or it is not $\Delta^0_2$, and by the low basis theorem, there is a member of~$\Ccal$
of low degree, which therefore does not compute~$X$.

We shall prove in further parts several 
other basis theorems for more marginal
computability-theoretic properties.

\subsection{K\"onig's lemma}

We now clarify the links between K\"onig's lemma
and weak K\"onig's lemma. Given a collection of sets~$S$,
we denote by~$\deg(S)$ the set of Turing degrees of elements of~$S$.

The basis theorems stated above are all based on the fact that
we can computably bound the successors of any node.
More precisely, a tree~$T \subseteq \Nb^{<\Nb}$ is \emph{computably bounded}
if there is a computable function~$f$ such that if~$\sigma \in T$,
then~$\sigma(n) \leq f(n)$ for all~$n < |\sigma|$.
Jockusch and Soare~\cite{Jockusch197201} proved that we can code a computably bounded
tree into a binary tree.

\begin{theorem}[Jockusch and Soare]
For every computable, computably bounded tree~$T$, there is a binary
tree~$\hat{T}$ such that~$\deg([\hat{T}]) = \deg([T])$.
\end{theorem}

However, when considering computable, finitely branching trees in
their full generality, the problem of finding a path becomes more difficult.
In particular, Jockusch, Lewis and Remmel~\cite{Jockusch199101} constructed a computable,
finitely branching tree whose unique path computes the halting set.
Since every infinite, finitely branching tree is $\emptyset'$-computably bounded,
every PA degree relative to~$\emptyset'$ computes an infinite path through it.
In the other direction, Jockusch, Lewis and Remmel proved the existence
of such a tree whose paths are all of PA degree relative to~$\emptyset'$
with the following one-to-one correspondence.

\begin{theorem}[Jockusch, Lewis and Remmel]
For every $\emptyset'$-computable infinite binary tree~$T$,
there is a computable, finitely branching tree~$\hat{T}$
such that~$\deg([\hat{T}]) = \deg([T])$.
\end{theorem}

Therefore, according to Hirschfeldt~\cite{Hirschfeldt2015Slicing}, K\"onig's lemma can be considered
has behaving like weak K\"onig's lemma ``one jump up''.

\subsection{Algorithmic randomness}\label{subsect:intro-comp-algorithmic-randomness}

The word ``randomness'' is part of the common language
and can be employed with very different meanings, according to which aspect
of randomness we consider. Algorithmic randomness uses the framework of computability theory
to study the various notions of randomness.
Among them, Martin-L\"of randomness stands out from the crowd
by its numerous characterizations with very different paradigms, namely,
measure theory, Kolmogorof complexity, martingales.

\index{MLR|see {Martin-L\"of random}}
\index{Martin-L\"of!test}
\index{Martin-L\"of!random}
\begin{definition}[Martin-L\"of randomness]
A~\emph{Martin-Löf test} is a sequence~$U_0, U_1, \dots$ of uniformly~$\Sigma^0_1$
classes such that~$\mu(U_i) \leq 2^{-i}$ for every~$i \in \omega$.
A real~$Z$ is \emph{Martin-L\"of random} if for every Martin-L\"of test $U_0, U_1, \dots$,
$Z \not \in \bigcap_i U_i$.
\end{definition}

Martin-L\"of randomness has strong connections with a subclass of weak K\"onig's lemma.

\index{tree!measure}
\begin{definition}
A binary tree~$T \subseteq 2^{<\Nb}$ has \emph{positive measure}
if $\lim_s \frac{|\{\sigma \in T : |\sigma| = s \}|}{2^s} > 0$.
\end{definition}

The restriction of weak K\"onig's lemma to trees of positive measure
is poetically called \emph{weak weak K\"onig's lemma}.
Informally, a tree of positive measure is very bushy,
and each branch is very likely to be infinite.
Therefore, one may expect to obtain an infinite path by a random walk trough 
a tree of positive measure.
The following theorem makes precise the correspondence between weak weak K\"onig's lemma
and Martin-L\"of randomness.

\begin{theorem}[Ku\v{c}era~\cite{Kucera1985Measure}]\ 
\begin{itemize}
	\item[(i)] A Martin-L\"of random real is a path (up to prefix) through a tree~$T \subseteq 2^{<\Nb}$ iff $T$
	has positive measure.
	\item[(ii)] There is a computable binary tree of positive measure
	whose paths are all Martin-L\"of random.
\end{itemize}
\end{theorem}

\chapter{Reverse mathematics}

What are the set existence axioms needed to prove ordinary theorems
in mathematics? How constructive are our theorems?
How to relate two theorems and formalize the intuition
that some theorems are consequences of others?
When will the author stop asking questions?

\index{reverse mathematics}
Reverse mathematics is a vast foundational program
whose goal is to study the logical strength of everyday life theorems.
Before getting into the details, let us introduce some of its main motivations.

\subsection{In search of optimal axioms}

There may be several ways to measure the naturality of an axiom.
A simple way consists in providing a philosophical justification of its truth.
One can also consider the collection of its consequences,
and in particular what desired theorems are unprovable without it.
This method requires to be able to classify theorems
according to the axioms they need. The finer this classification is,
the more insights about the naturality of the axioms we get.
In the best case, we should be able to find axioms
which are \emph{equivalent} to the considered theorem.
This is exactly the purpose of reverse mathematics, which
tries to (and succeeds in) finding the optimal axioms
to prove natural theorems.

\subsection{The constructive perspective}

The 20th century have seen the emergence of constructivism,
with Brower's intuitionism in 1908, Hilbert and Bernay's finitism in 1920,
and Bishop's constructive analysis in 1967 among others.
Mathematicians started to care about the constructive content of their theorems,
and to rebuild the core of mathematics under this perspective.
Type theory and topos theory are examples of alternative
foundational theories focusing on the constructive content of the theorems.
Reverse mathematics follows on from this movement by providing
a framework to analyse the logical strength and the computational
content of ordinary theorems.

\subsection{The interrelationships between theorems}

It is a common mathematical practice to state corollaries
after a theorem. Informally, a corollary is an immediate consequence
of the previous result. Being a corollary~$\Qsf$ of a theorem~$\Psf$ is not an intrinsic property of the statement~$\Qsf$,
but a relationship between~$\Psf$ and~$\Qsf$. Sometimes, the theorem~$\Psf$
can also be deduced from its corollary~$\Qsf$ by elementary means, in which case we may
think of~$\Psf$ and~$\Qsf$ as being equivalent. Reverse mathematics enables one to make formal
the notion of ``immediate consequence'', by prodiving a framework in which
an implication~$\Psf \imp \Qsf$ is interpreted as ``the statement~$\Qsf$ can be logically deduced from~$\Psf$
in a very weak theory''.

\subsection{The search of simpler proofs}

Although a theorem logically reduces to its truth value, 
in the evolving process of mathematical discovery, 
a theorem is a living object.
The proofs of important theorems are regularly refined and
simplified, with the ultimate goal to find a proof that
Paul Erd\H{o}s would call ``from The Book''~\cite{Aigner2014Proofs}.
Such a proof would be natural enough to reveal the essence of the theorem, and weak enough so
that no stronger statement can be proved without modifying the core of the argument.
In his attempt to find the optimal axioms
needed by a theorem, the reverse mathematician
will sometimes need to find new and constructively simpler proofs,
and in this way contributes to give new insights into the theorems.

\section{The strength of theorems}

\index{instance}
\index{solution}
Many theorems of ``ordinary'' mathematics can be formulated as \emph{mathematical problems},
coming with a natural class of~\emph{instances}. Each instance has a collection
of \emph{solutions} associated to it.
In the case of K\"onig's lemma which asserts that every infinite, finitely branching tree has an infinite path,
an instance is an \emph{infinite, finitely branching tree} and a solution is an \emph{infinite path} through the tree.
More precisely, we are interested in the theorems which can be formalized in
the language of second-order arithmetic. This language happens to be expressive
enough to state in a natural way many theorems coming from~\emph{countable mathematics}.
See~Simpson~\cite{Simpson2009Subsystems} for a discussion about the notion of countable mathematics.
Formally, we mainly consider~$\Pi^1_2$ sentences~$\Psf$ of the form
\[
\forall X(\varphi(X) \imp \exists Y \psi(X,Y))
\]
where~$\varphi$ and~$\psi$ are arithmetic formulas, that is, formulas which do not contain any set quantifier.
A set~$X$ such that~$\varphi(X)$ holds is a $\Psf$-instance, and a set~$Y$
such that~$\psi(X,Y)$ holds is a \emph{solution to~$X$}.

Some theorems are more \emph{effective} than others.
Consider the proof of the intermediate value theorem.
Given a continuous real-valued function which is negative 
at $0$ and positive at $1$, one can (non-uniformly) compute a real $x \in (0,1)$
such that $f(x)=0$ essentially by using the usual interval-halving procedure.
On the other hand, we have seen that K\"onig's lemma is non-effective,
in the sense that given a computable infinite, finitely branching tree,
one cannot in general compute an infinite path.
In the same spirit as for the Turing degrees, one may wonder
how non-effective theorems are whenever they do not admit
a computable solution.

There are many ways to understand the notion of strength of a theorem.
In this chapter, we embrace the provability approach, and try to capture
the logical strength of the theorems in terms of set existence axioms.
The main idea is to place oneself in a theory~$\Tsf$ weak enough so that
non-effective theorems are not provable. In this theory, a proof of~$\Psf \imp \Qsf$
can therefore be interpreted as ``the statement~$\Qsf$ can be obtained
from~$\Psf$ by effective (and therefore elementary) means''.

\section{The base theory}

\index{rca@$\rca$}
We now formally define the weak theory which will serve as a base
to study our theorems. Such a theory has to be chosen carefully
so that an implication~$\Psf \imp \Qsf$ provides some insights into the relation
between the statements~$\Psf$ and~$\Qsf$. The base theory is called~$\rca$,
standing for~\emph{recursive comprehension axiom}.

We first describe the first-order part of~$\rca$.
We would like the integers to behave sufficiently nicely, not 
to be bothered by the coding details of the objects we manipulate.
The base theory therefore contains the basic first-order Peano axioms.

\[
\begin{array}{ll}
\forall m (m+1 \neq 0) & \forall m \forall n (m \times (n+1) = (m \times n) + m)\\
\forall m \forall n(m+1 = n+1 \imp m=n) & \forall m \forall n (m < n+1 \biimp (m < n \orr m = n))\\
\forall m (m+0 = m) & \forall m \neg(m < 0)\\
\forall m \forall n(m+(n+1) = (m+n)+1) & \forall m (m \times 0 = 0) \\
\end{array}
\]

\index{scheme!comprehension}
In the second-order part, we want to allow only sets which can be constructed from other sets
in an effective way.
By Post's theorem, a set~$X$ is computable from~$Y$ if it is~$\Delta^{0,Y}_1$-definable.
We therefore add a set-building scheme, the \emph{comprehension scheme} restricted to~$\Delta^0_1$ formulas
with parameters. The $\Delta^0_1$ \emph{comprehension scheme} consists of the universal closures of all formulas of the form
\begin{align*}
\forall n (\varphi(n) \biimp \psi(n)) \imp \exists X \forall n(n \in X \biimp \varphi(n)),
\end{align*}
where $\varphi$ is $\Sigma^0_1$, $\psi$ is $\Pi^0_1$, and $X$ is not free in $\varphi$.

\index{scheme!induction}
\index{isig@$\isig^0_n$|see {scheme induction}}
Finally, we add some induction scheme.
We should also be careful about the amount of induction we want,
since it adds new finite sets in the theory.
Since the logical strength of a theorem is defined according to its set existence axioms,
it is natural to control the finite sets as well. Such a control
enables one to prove conservativity results.
For each $n \in \omega$, the $\Sigma^0_n$ (resp.\ $\Pi^0_n$) \emph{induction scheme}, 
denoted $\isig^0_n$ (resp.\ $\ipi^0_n$), consists of the universal closures of all formulas of the form
\begin{align*}
[\varphi(0) \andd \forall n(\varphi(n) \imp \varphi(n+1))] \imp \forall n \varphi(n),
\end{align*}
where $\varphi$ is $\Sigma^0_n$ (resp.\ $\Pi^0_n$). 
The $\Sigma^0_1$ induction scheme is exactly the amount of induction we want to add to~$\rca$.
Indeed, the $\Sigma^0_1$ induction scheme is equivalent to the \emph{bounded $\Delta^0_1$ comprehension scheme},
that is, the $\Delta^0_1$ comprehension scheme stating the existence of sets
of the form~$\{n \in X \biimp n < b \wedge \varphi(n)\}$.

\index{scheme!bounding}
\index{bsig@$\bsig^0_n$|see {scheme bounding}}
One can prove over $\rca$ the properties of many basic primitive recursive codings, and in particular
the properties of the Cantor pairing function~$\tuple{\cdot, \cdot} : \Nb \times \Nb \to \Nb$.
After some manipulation, the mathematician gets an intuition about which argument can,
and which cannot be formalized over~$\rca$. One of the main pitfalls
of which a beginner (and even a more experimented mathematician) should take care,
is that a finite union of finite sets is not necessarily finite. Those statements
are the \emph{bounding schemes}.
For each $n \in \omega$, the $\Sigma^0_n$ (resp.\ $\Pi^0_n$) \emph{bounding scheme}, denoted $\bsig^0_n$ (resp.\ $\bpi^0_n$), consists of the universal closures of all formulas of the form
\begin{align*}
\forall a[(\forall n < a)(\exists m)\varphi(n,m) \imp \exists b(\forall n < a)(\exists m < b)\varphi(n,m)],
\end{align*}
where $\varphi$ is $\Sigma^0_n$ (resp.\ $\Pi^0_n$).
The $\Sigma^0_n$ bounding scheme is a strict consequence of the~$\Sigma^0_n$ induction scheme.
In particular, $\rca \vdash \bsig^0_1$, but not~$\bsig^0_2$. Such a concern happens
for example when we use the fact that given the $\Delta^0_2$ approximation~$g$ of a $\Delta^0_2$
set, and some integer~$b$, there is some stage~$s$ after which $g(x, \cdot)$ is constant for every $x < b$.

\section{The Big Five phenomenon}

\index{Big Five}
Over the past decades, theorems from all over mathematics have been studied
within the framework of reverse mathematics. A surprising phenomenon
emerged from the early years of reverse mathematics:
Most theorems studied require \emph{very weak set existence axioms}.
Moreover, many of them happen to be \emph{equivalent} to one of five main subsystems over~$\rca$,
that Mont\'alban~\cite{Montalban2011Open} refers to as the Big Five.
These are, in the strictly increasing strength order, $\rca$,
weak K\"onig's lemma ($\wkl$), the arithmetic comprehension axiom scheme ($\aca$),
arithmetic transfinite recursion~($\atr$) and the $\Pi^1_1$ comprehension axiom scheme ($\pioca$).
These subsystems are characterized by their \emph{robustness}, that is,
the invariance of their strength under small changes in their definition.
It is currently unclear whether the existence of the Big Five reflects
a bias in the statements studied in mathematics, or whether
there is some philosophical justification to this structure.
We briefly detail the five subsystems.

\subsection{$\rca$}

We already detailed the definition of~$\rca$. By the equivalence
for a set to be $\Delta^0_1$-definable and to be computable,
$\rca$ can be thought of as capturing \emph{constructive mathematics}.
If we consider the structures in which the first-order part is composed
of the standard integers~$\omega$ with the standard operations, 
$\rca$ admits a minimal model~$\Mcal_0$ (in the inclusion sense),
in which the second-order part consists of all the computable sets.
Therefore, if some theorem~$\Psf$ admits a computable instance with no computable solution,
$\Mcal_0$ is not a model of~$\Psf$, hence~$\rca$ does not prove~$\Psf$.

\subsection{$\wkl$}
\index{wkl@$\wkl$}
Weak K\"onig's lemma is the statement ``Every infinite binary tree has a path''.
$\wkl$ captures \emph{compactness arguments} and is equivalent to several famous
theorems, such as G\"odel's completeness theorem or the Heine/Borel covering lemma.

\subsection{$\aca$}
\index{aca@$\aca$}
The arithmetic comprehension axiom scheme consists of the universal closures of all formulas of the form
\begin{align*}
\exists X \forall n(n \in X \biimp \varphi(n)),
\end{align*}
where $\varphi$ is a formula with no set quantifier and such that~$X$ is not free in $\varphi$.
$\aca$ is sufficiently expressive to prove the vast majority of ordinary theorems over~$\rca$.
It is in particular equivalent to full K\"onig's lemma and to the statement 
``Every set has a Turing jump'' over~$\rca$.
Thus, every model of~$\rca + \aca$ contains the halting set.

\subsection{$\atr$}
\index{atr@$\atr$}
The arithmetic transfinite recursion axiom informally asserts that 
the Turing jump can be iterated along any countable well-ordering,
starting at any set. More precisely, let~$\varphi(n, X)$
be an arithmetic formula with parameters, and let~$\Phi(X) = \{n : \varphi(n, X) \}$.
$\atr$ is the statement ``For every well-order~$(D, \preceq)$,
the sequence of sets~$\{K_n : n \in D\}$ such that
$K_n = \Phi(\{\tuple{m,x} : m \prec n \wedge x \in K_m\})$ exists.''
$\atr$ is sufficiently strong to prove many basic theorems of classical descriptive set theory,
and is equivalent to Ulm's theorem for countable reduced Abelian groups.

\subsection{$\pioca$}
\index{pioca@$\pioca$}
The~$\Pi^1_1$ comprehension scheme is the comprehension scheme restricted to $\Pi^1_1$ formulas,
that is, formulas of the form~$\forall X \theta$ where~$\theta$ is an arithmetic formula.
$\pioca$ is strong enough to develop a good theory of ordinals
and is equivalent to the Cantor-Bendixson theorem. 
\bigskip

In this thesis, we shall focus exclusively on statements below~$\aca$
over~$\rca$. A large part of mathematics can already be developped within this subsystem.

\section{Ramsey's theorem}

In order to put to the test the robustness of the Big Five phenomenon, the reverse mathematics
community started to investigate candidate theorems leading to their own subsystems.
Ramsey's theory provides, among others, a large collection of statements escaping
the Big Five phenomenon.
This theory informally asserts that given some \emph{size}~$s$,
every \emph{sufficiently large} collection of objects has a
sub-collection of size~$s$, whose objects satisfy some \emph{structural properties}.
Perhaps the most famous example of such a statement is \emph{Ramsey's theorem}.
Given a set~$X$, we denote by~$[X]^k$ the set of unordered $k$-tuples over~$X$.

\index{Ramsey's theorem}
\index{homogeneous!for a coloring}
\index{rt@$\rt^n_k$|see {Ramsey's theorem}}
\begin{definition}[Ramsey's theorem]
A subset~$H$ of~$\Nb$ is~\emph{homogeneous} for a coloring~$f : [\Nb]^n \to k$ (or \emph{$f$-homogeneous}) 
if each $n$-tuples over~$H$ are given the same color by~$f$.
$\rt^n_k$ is the statement ``Every coloring $f : [\Nb]^n \to k$ has an infinite $f$-homogeneous set''.
\end{definition}

Jockusch~\cite{Jockusch1972Ramseys} conducted a computable analysis of Ramsey's theorem
and gave precise definitional bounds on the complexity of homogeneous sets.
Simpson~\cite{Simpson2009Subsystems} formalized Ramsey's theorem in the setting of reverse mathematics,
and proved that $\rt^n_k$ is equivalent to~$\aca$ over~$\rca$ for each~$n \geq 3$.

The strength of Ramsey's theorem for pairs has been a long-standing open problem,
until Seetapun~\cite{Seetapun1995strength} proved that~$\rt^2_2$ is strictly weaker than~$\aca$ over~$\rca$.
Cholak, Jockusch and Slaman~\cite{Cholak2001strength} studied extensively Ramsey's theorem for pairs.
In particular, they showed that~$\rt^2_2$ is not a consequence of weak K\"onig's lemma over~$\rca$.
Later, Liu~\cite{Liu2012RT22} clarified the relation between $\rt^2_2$ and the Big Five by proving that~$\rt^2_2$
does not imply~$\wkl$ over~$\rca$.

Every computable instance of Ramsey's theorem for singletons has a computable solution.
However, the statement~$(\forall k)\rt^1_k$ (often written~$\rt^1_{<\infty}$) is not provable over~$\rca$.
In fact, it is equivalent to~$\bst$, the bounding scheme for~$\Sigma^0_2$ formulas~\cite{Cholak2001strength}.
The full statement $(\forall n)(\forall k)\rt^n_k$, written~$\rt$, is strictly stronger than~$\aca$,
again for induction reasons.

The full diagram of the relation between variants of Ramsey's theorem and the Big Five
is summarised in Figure~\ref{fig:intro-rm-ramsey-zoo}. An arrow between
$\Psf$ and~$\Qsf$ means that~$\rca \vdash \Psf \imp \Qsf$. The missing arrows are non-implications.

\usetikzlibrary{arrows}
\usetikzlibrary{decorations.markings}

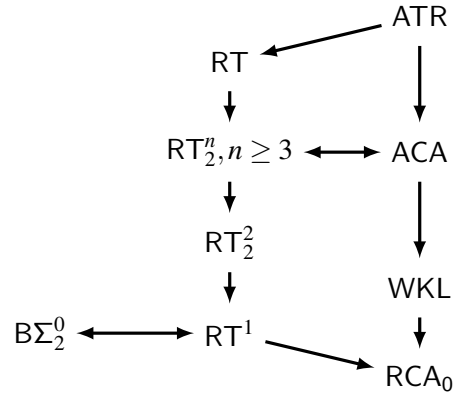
\begin{figure}[htbp]
\begin{center}
\begin{tikzpicture}[x=2.5cm, y=1.2cm, 
	node/.style={minimum size=2em},
	impl/.style={draw,very thick,-latex},
	equiv/.style={draw, very thick, latex-latex}]

	\node[node] (atr) at (2, 4.5) {$\atr$};
	\node[node] (rt22) at (1,2) {$\rt^2_2$};
	\node[node] (rt32) at (1, 3)  {$\rt^n_2, n \geq 3$};
	\node[node] (rt) at (1, 4) {$\rt$};
	\node[node] (rt1) at (1, 1) {$\rt^1$};
	\node[node] (bst) at (0, 1) {$\bst$};
	\node[node] (wkl) at (2, 1.5) {$\wkl$};
	\node[node] (aca) at (2, 3) {$\aca$};
	\node[node] (rca) at (2, 0.5) {$\rca$};
	
	\draw[impl] (atr) -- (rt);
	\draw[impl] (atr) -- (aca);
	\draw[impl] (rt) -- (rt32);
	\draw[equiv] (aca) -- (rt32);
	\draw[impl] (aca) -- (wkl);
	\draw[impl] (wkl) -- (rca);
	\draw[impl] (rt32) -- (rt22);
	\draw[impl] (rt22) -- (rt1);
	\draw[equiv] (rt1) -- (bst);
	\draw[impl] (rt1) -- (rca);

\end{tikzpicture}
\end{center}
\caption{Ramsey's theorem zoo}
\label{fig:intro-rm-ramsey-zoo}
\end{figure}

Due to the complexity of their separations, Ramsey-type statements
received a lot of attention from the reverse mathematics community.
Many consequences of Ramsey's theorem have been investigated in reverse mathematics,
coming from graph theory~\cite{Bovykin2005strength,Lerman2013Separating,PateyDominating}, 
model theory~\cite{Hirschfeldt2009atomic,Conidis2008Classifying}, set theory~\cite{Downey2012finite,Cholak2015Any}, combinatorics~\cite{Csima2009strength,FriedmanFom53free,Cholak2001Free,Wang2014Some}
and order theory~\cite{Hirschfeldt2007Combinatorial} among others.

\chapter{Reducibilities}\label{chap:introduction-reducibilities}

There are many ways to understand the strength of a theorem.
We have already seen the provability approach, which consists
in comparing the logical consequences of the theorems
in terms of set existence axioms.
We will now present the computational approach
which comes up with various reducibility notions,
depending on which specific constructive properties we focus on.
More precisely, we will present three reducibility notions
and relate them to reverse mathematics.

\section{Computable reducibility}

One of the simplest way to define the relative computable strength between two theorems
consists of using the natural reducibility notion induced by the Turing reduction.
For this, take again the interpretation of a theorem~$\Psf$ as a mathematical problem with instances and solutions.
Informally, a problem~$\Qsf$ is \emph{simpler} than another problem~$\Psf$ if the ability
to solve any~$\Psf$-instance gives the ability to solve any~$\Qsf$-instance. 

In other words, the theorem~$\Psf$ can be seen as a \emph{blackbox},
where the input is a $\Psf$-instance, and the output is \emph{any} valid solution. 
Therefore, $\Qsf$ is simpler than~$\Psf$ if we can \emph{simulate} a blackbox for~$\Qsf$ from any blackbox for~$\Psf$.
We want some effectiveness restriction on the simulation to give some computational meaning to this relation.
In particular, the $\Psf$-instance has to be constructed computably from the~$\Qsf$-instance,
and a solution to the~$\Qsf$-instance has to be transformed into a solution to the~$\Psf$-instance by constructive means.

\index{reduction!computable}
\index{computable reducibility}
\index{computable reducibility!strong}
\index{$\leq_c$}
\index{$\leq_{sc}$}
\begin{definition}[Computable reducibility] Fix two~$\Pi^1_2$ statements~$\Psf$ and $\Qsf$.
$\Qsf$ is \emph{computably reducible} to~$\Psf$ (written $\Qsf \leq_c \Psf$)
if every~$\Qsf$-instance~$I$ computes a~$\Psf$-instance~$J$ such that for every solution~$X$ to~$J$,
$X \oplus I$ computes a solution to~$I$.
\end{definition}

Computable reducibility is illustrated in Figure~\ref{fig:intro-reduc-comp-reduc}.
Note that the solution to the~$\Qsf$-instance~$I$ can be computed with the help of the instance~$I$.
One can imagine a reducibility notion where the solution to~$I$ has to be computable
from the solution to the~$\Psf$-instance. This leads to the notion of~\emph{strong computable reducibility}.
Various reducibility notions have a strong version.

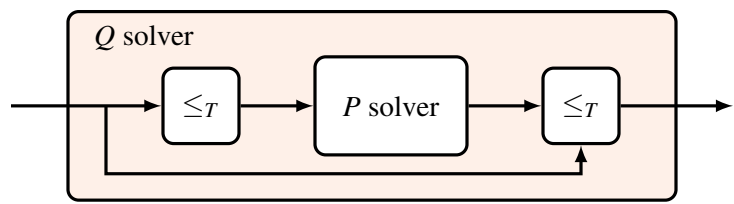
\begin{figure}[htbp]
\begin{center}
\begin{tikzpicture}[x=2.5cm, y=0.9cm, 
	node/.style={minimum size=2em},
	box/.style={draw,shape=rectangle, very thick,rounded corners},
	arrow/.style={draw, very thick, -latex}]
 
	\node[box,fill=ocre!10,minimum width=8cm,minimum height=2.5cm] at (0.9,0) {};
	\node[box,fill=white,minimum width=2cm,minimum height=1.3cm] (P) at (1,0) {$P$ solver};
	\node[box,fill=white,minimum width=1cm,minimum height=1cm] (firstR) at (0,0) {$\leq_T$};
	\node[box,fill=white,minimum width=1cm,minimum height=1cm] (secondR) at (2,0) {$\leq_T$};
	\node[node] at (-0.3,1) {$Q$ solver};

	\draw[arrow] (-1,0) -- (firstR);
	\draw[arrow] (firstR) -- (P);
	\draw[arrow] (P) -- (secondR);
	\draw[arrow] (secondR) -- (2.8, 0);
	\draw[arrow] (-0.5,0) -- (-0.5,-1) -- (2,-1) -- (secondR);
\end{tikzpicture}
\end{center}
\caption{Computable reducibility}
\label{fig:intro-reduc-comp-reduc}
\end{figure}

Computable reducibility provides in some sense a more \emph{fine-grained}
analysis of the relations between theorems than reverse mathematics.
Indeed, in reverse mathematics, a principle can be applied an arbitrary number of times.
For example, the statement~``Every set has a Turing jump'' is equivalent to~``Every set has a Turing double jump''
over~$\rca$, since for any model~$\Mcal$ of the former statement and any set~$X \in \Mcal$, $X' \in \Mcal$,
and hence~$X'' \in \Mcal$.
Therefore, an implication from~$\Psf$ to~$\Qsf$ over~$\rca$
may involve many~$\Psf$-instances to obtain a solution to one $\Qsf$-instance.
On the other hand, computable reducibiliy allows the use of only one $\Psf$-instance in a reduction from~$\Qsf$ to~$\Psf$.

This difference is useful to reveal some subtle distinctions between statements
that would be collapsed in reverse mathematics. In particular, Ramsey's theorem
does not depend on the number of its colors in reverse mathematics as witnessed by a simple
color blindness argument, whereas we shall see that $\rt^n_{k+1} \not \leq_c \rt^n_k$ for any~$n, k \geq 2$.

Many implications proofs over~$\rca$ happen to be computable reductions.
However, a proof of~$\Qsf \leq_c \Psf$ is not sufficient to deduce that~$\rca \vdash \Psf \imp \Qsf$.
Indeed, the argument has to be formalizable over~$\rca$, and in particular one must care
about the amount of induction used to prove its validity.
On the other hand, proving that~$\Qsf$ does not computably reduce to~$\Psf$ is simpler than
separating~$\Qsf$ from~$\Psf$ over standard models, since in the $\Qsf \not \leq_c \Psf$ case, one has to build a $\Qsf$-instance~$I$
which diagonalizes against only the~$I$-computable $\Psf$-instances. Computable non-reducibility
is therefore often used as a preliminary separation, which is then generalized
to a separation over~$\rca$. This is in particular the approach of Lerman, Solomon and Towsner~\cite{Lerman2013Separating}.

\section{Weihrauch reducibility}

Some proofs are more uniform than others. Let us take back the example
of the intermediate value theorem.
Given a continuous real-valued function~$f$ which is negative 
at $0$ and positive at $1$, the computation of a real~$x \in (0,1)$ such that~$f(x) = 0$
is split into two cases.
In the first case, there is a rational~$x \in (0, 1)$ such that~$f(x) = 0$.
Every rational is computable, so $f$ admits a computable solution, even though the real~$x$ is not effectively provided.
In the second case, we compute a solution to~$f$ essentially by using an interval-halving procedure.
The standard proof of the intermediate value theorem involves a case analysis,
and in this sense, is not uniform. 
Moreover, Brattka and Gherardi~\cite{Brattka2011Effective}
proved that the intermediate value theorem does not admit a uniform proof.

\emph{Weihrauch reducibility} is a refinement of computable reducibility,
where the transformation of a $\Qsf$-instance into a $\Psf$-instance has to be done
uniformly in the instance, and so has to be the computation of a $\Qsf$-solution.

\index{reduction!Weihrauch}
\index{Weihrauch reducibility}
\index{Weihrauch reducibility!strong}
\index{uniform reducibility|see {Weihrauch reducibility}}
\index{$\leq_W$}
\index{$\leq_{sW}$}
\begin{definition}[Weihrauch reducibility]
Fix two~$\Pi^1_2$ statements~$\Psf$ and $\Qsf$.
$\Qsf$ is \emph{Weihrauch reducible} to~$\Psf$ (written $\Qsf \leq_W \Psf$)
if there are two Turing functionals~$\Phi$ and~$\Gamma$ such that 
for every $\Qsf$-instance~$I$, $\Phi^I$ is a $\Psf$-instance~$J$ such that for every solution~$X$,
$\Gamma^{I \oplus X}$ is a solution to~$I$.
\end{definition}

Weihrauch reducibility has been used by Dorais, Dzhafarov, Hirst, Mileti and Shafer~\cite{Dorais2016uniform} 
to study extensively the uniform content of combinatorial theorems. It is currently an active
research subject~\cite{DzhafarovStrong,Hirschfeldtnotions,Rakotoniaina2015Computational}.

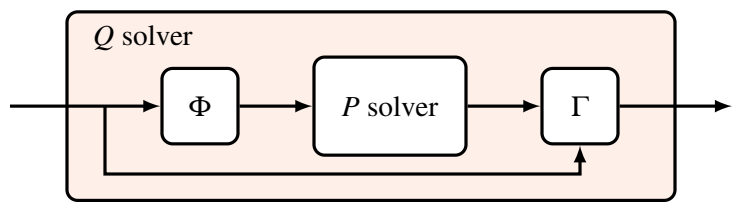
\begin{figure}[htbp]
\begin{center}
\begin{tikzpicture}[x=2.5cm, y=0.9cm, 
	node/.style={minimum size=2em},
	box/.style={draw,shape=rectangle, very thick,rounded corners},
	arrow/.style={draw, very thick, -latex}]
 
	\node[box,fill=ocre!10,minimum width=8cm,minimum height=2.5cm] at (0.9,0) {};
	\node[box,fill=white,minimum width=2cm,minimum height=1.3cm] (P) at (1,0) {$P$ solver};
	\node[box,fill=white,minimum width=1cm,minimum height=1cm] (firstR) at (0,0) {$\Phi$};
	\node[box,fill=white,minimum width=1cm,minimum height=1cm] (secondR) at (2,0) {$\Gamma$};
	\node[node] at (-0.3,1) {$Q$ solver};

	\draw[arrow] (-1,0) -- (firstR);
	\draw[arrow] (firstR) -- (P);
	\draw[arrow] (P) -- (secondR);
	\draw[arrow] (secondR) -- (2.8, 0);
	\draw[arrow] (-0.5,0) -- (-0.5,-1) -- (2,-1) -- (secondR);
\end{tikzpicture}
\end{center}
\caption{Weihrauch reducibility}
\label{fig:intro-reduc-uniform-reduc}
\end{figure}

Although the flavor of Weihrauch reducibility is very similar to the one of computable reducibility,
the constructive properties they focus on is very different. In particular, Weihrauch reducibility
enables one to establish distinctions between statements which are provable over~$\rca$,
whereas computable reducibility cannot since $\rca$ has a model whose second-order part is exactly
the collection of computable sets. For example, the intermediate value theorem is provable over~$\rca$,
but the proof is not uniform (see Brattka and Gherardi~\cite{Brattka2011Effective}).

\section{Computable entailment}

Although the historical motivation of reverse mathematics is proof-theoretic,
many researchers in the field adopt the computational point of view.
In particular, they do not feel concerned about the amount of induction needed
to formalize their arguments. Starting from this observation, Shore~\cite{Shore2010Reverse}
proposed a computability-theoretic reduction with a model-theoretic flavor closely related to reverse mathematics.

The models of reverse mathematics are of the form~$\Mcal = (M, \Ical, +, \times, <, 0, 1, \in)$
where~$M$ and~$\Ical$ are the first-order and second-order parts, respectively,
$+$ and~$\times$ are binary functions on~$M$, $<$ is a binary relation on~$M$
and~$0$ and~$1$ are members of~$M$. An~\emph{$\omega$-structure} 
is a structure~$\Mcal = (\omega, \Ical, +_\omega, \times_\omega, <_\omega, 0, 1, \in)$
where~$\omega$ are the standard integers, coming with the standard operations. 
An $\omega$-structure is therefore fully specified by its second-order part~$\Ical$.
The second-order parts of $\omega$-models of~$\rca$ admits a simple purely computability-theoretic characterization.

\index{Turing!ideal}
\begin{definition}[Turing ideal]
A \emph{Turing ideal} $\Ical$ is a collection of subsets of~$\omega$ which is closed under
\begin{itemize}
	\item[(i)] the Turing reduction: $(\forall X \in \Ical)(\forall Y \leq_T X)[Y \in \Ical]$
	\item[(ii)] the effective join: $(\forall X, Y \in \Ical)[X \oplus Y \in \Ical]$.
\end{itemize}
\end{definition}

Friedman~\cite{Friedman1974Some} proved that the second-order parts of~$\omega$-models of~$\rca$
are exactly the Turing ideals. Therefore, we say that a second-order formula $\varphi$
\emph{holds} in a Turing ideal~$\Ical$ (written~$\Ical \models \varphi$)
if it holds in the $\omega$-structure whose second-order part is~$\Ical$.
In particular, a $\Pi^1_2$ statement~$\Psf$ holds in~$\Ical$ if every~$\Psf$-instance in~$\Ical$
has a solution in~$\Ical$.

\index{computable entailment}
\index{$\leq_\omega$}
\begin{definition}[Computable entailment]
Fix two second-order formulas~$\varphi$ and~$\psi$. We say that 
$\varphi$ \emph{computably entails} $\psi$ if every Turing ideal satisfying~$\varphi$
also satisfies~$\psi$.
\end{definition}

In other words, $\varphi$ computably entails~$\psi$ if whenever~$\varphi$ holds
in an $\omega$-model of~$\rca$, then so does~$\psi$.
Note that the notion of computable entailment is defined over every second-order formulas.
It is therefore a more general framework than computable and Weihrauch reducibilities.
Hirschfeldt and Jockusch~\cite{Hirschfeldtnotions} used the notation~$\Qsf \leq_\omega \Psf$ to express that the~$\Pi^1_2$
statement~$\Psf$ computably entails the $\Pi^1_2$ statement~$\Qsf$.
In this thesis, we shall consider only computable entailments for~$\Pi^1_2$ statements.

Although recently introduced as a formal reduction, this reducibility notion was already used \emph{de facto}
for a long time, in particular to separate statements over~$\rca$. Indeed,
any proof of $\Qsf \not \leq_\omega \Psf$ is a proof of~$\rca + \Psf \not \vdash \Qsf$.
Computable entailment and computable reducibility coincide for some statements like~$\wkl$.

\begin{theorem}[\cite{Hirschfeldtnotions}]\label{thm:intro-reduc-equiv-entail}
For any $\Pi^1_2$ statement~$\Psf$, $\Psf \leq_c \wkl$ iff $\Psf \leq_\omega \wkl$.
\end{theorem}
\begin{proof}
Suppose that~$\Psf \leq_\omega \wkl$. Let~$I$ be a $\Psf$-instance,
and let~$T$ be an $I$-computable tree such that every path is of PA degree relative to~$I$.
By Scott~\cite{Scott1962Algebras}, for every path~$S$ through~$T$, $S$ bounds a Turing ideal~$\Ical$ 
containing~$I$ in which~$\wkl$ holds. Since~$\Psf \leq_\omega \wkl$, $\Psf$ holds in~$\Ical$,
so $\Ical$ contains a solution to~$I$ and therefore every path through~$T$ computes a solution to~$I$.
\end{proof}

Theorem~\ref{thm:intro-reduc-equiv-entail} is essentially due to the fact 
that every set of PA degree bounds an $\omega$-model of~$\wkl$.
Such a property is also true for any notion for which
a Van Lambalgen-like theorem holds. This is in particular the case of Martin-L\"of randomness~\cite{VanLambalgen1990axiomatization}
and 1-genericity~\cite{Yu2006Lowness} (defined in section~\ref{sect:intro-forcing-effective-forcing}).

Hirschfeldt and Jockusch~\cite{Hirschfeldtnotions} introduced a \emph{generalized Weihrauch reducibility}
which is the dual notion of computable entailment for Weihrauch reductions. This reduction enables
one to express a relation between statements which needs more than one application, but for which
the multiple instances are uniformly built. This is in particular the case for
the reduction from~$\rt^n_{k+1}$ to~$\rt^n_k$.

\section{Separating principles}

Consider two $\Pi^1_2$ statements~$\Psf$ and~$\Qsf$.
Unfolding the definition, in order to prove that $\Psf$ does not imply~$\Qsf$ over $\rca$,
one needs to create a model~$\Mcal$ satisfying~$\rca + \Psf$
but not~$\Qsf$. In its full generality, the model~$\Mcal$ can be non-standard,
but one usually wants to construct an $\omega$-structure,
so that the separation is not due to some induction technicalities,
but also holds with our standard understanding of the integers.
The usual construction of an $\omega$-model of $\rca + \Psf$
is done by building a Turing ideal in which $\Psf$ holds, as follows.

\begin{quote}
Start with a \emph{topped} Turing ideal $\Ical_0$, that is,
a Turing ideal of the form~$\{Z : Z \leq_T C \}$ for some set~$C$.
The Turing ideal of all computable sets is commonly chosen.
Assuming that as stage~$s$, we have a topped Turing ideal~$\Ical_s = \{ Z : Z \leq_T C_s \}$
for some set~$C_s$,
\begin{itemize}
	\item[(i)] Pick the \emph{next} $\Psf$-instance~$X \in \Ical_s$ with no solution in~$\Ical_s$, for a reasonable
	ordering so that each $\Psf$-instance will receive attention at some stage.
	\item[(ii)] Construct a solution~$Y$ to~$X$. Such a solution exists by the classical proof of~$\Psf$.
	\item[(iii)] Set~$\Ical_{s+1} = \{ Z : Z \leq_T Y \oplus C_s \}$. 
\end{itemize}
Then, go to the next stage. Finally, take~$\Ical = \bigcup_s \Ical_s$.
\end{quote}

Furthermore, if one want the resulting Turing ideal~$\Ical$ not to satisfy~$\Qsf$,
one need to be a bit more careful during the construction. First, there must be
a $\Qsf$-instance~$I$ with no~$I$-computable solution, otherwise every statement
computably entails~$\Qsf$. Start with~$\Ical_0 = \{Z : Z \leq_T I\}$.
By the choice of the~$\Qsf$-instance, $\Qsf$ does not hold in~$\Ical_0$.
Ideally (no pun intended), one would like to preserve the property that~$\Qsf$ does not hold in~$\Ical_s$
at any stage~$s$. However, sometimes, this invariant is not strong enough, and one adds solutions to
the~$\Qsf$-instance~$I$ during the step~(iii). One therefore needs to maintain a stronger invariant,
which leads to the notion of \emph{preservation of a weakness property}.

\index{preservation}
\index{preservation!of weakness}
\begin{definition}[Weakness preservation]
A \emph{weakness property} is a collection $\Pcal$ of subsets of $\omega$ which
is downward-closed under the Turing reduction.
A $\Pi^1_2$ statement~$\Psf$ admits \emph{$\Pcal$ preservation} for some weakness property~$\Pcal$,
if for every set~$C \in \Pcal$ and every $C$-computable $\Psf$-instance~$X$,
there is a solution~$Y$ to~$X$ such that~$Y \oplus C \in \Pcal$.
\end{definition}

We have already seen two important examples of weakness properties, namely, lowness and hyperimmune-freeness.
By the low and hyperimmune-free basis theorems,
$\wkl$ admits both lowness and hyperimmune-freeness preservation.
$\Pcal$ preservation is designed so that it is possible to tweak
the previous Turing ideal construction to maintain the invariant
that~$\Ical_s \subseteq \Pcal$.

\begin{lemma}\label{lem:intro-reduc-preservation-included}
If~$\Psf$ admits $\Pcal$ preservation, then for every set~$C \in \Pcal$,
$\Psf$ holds in a Turing ideal~$\Ical \subseteq \Pcal$ containing~$C$.
\end{lemma}
\begin{proof}
Follow the previous construction, starting with $\Ical_0 = \{ Z : Z \leq_T C \}$.
At any stage~$s$, maintain the invariant that~$\Ical_s = \{ Z : Z \leq_T C_s \}$
for some set~$C_s \in \Pcal$. Therefore, $\bigcup_s \Ical_s \subseteq \Pcal$.
This is possible since by $\Pcal$ preservation for~$\Psf$,
one can choose at step (ii) a solution~$Y$ to~$X$ such that~$Y \oplus C_s \in \Pcal$.
\end{proof}

In particular, $\Pcal$ preservation is strong enough to obtain
separations over $\omega$-models, and \emph{a fortiori}
separations over~$\rca$.

\begin{lemma}\label{lem:intro-reduc-preservation-separation}
If~$\Psf$ admits $\Pcal$ preservation but not~$\Qsf$,
then $\Qsf \not \leq_\omega \Psf$.
\end{lemma}
\begin{proof}
Let~$C \in \Pcal$ and~$I \leq_T C$ be a $\Qsf$-instance
such that for every solution~$Y$ to~$I$, $Y \oplus C \not \in \Pcal$.
By Lemma~\ref{lem:intro-reduc-preservation-included},
$\Psf$ holds in a Turing ideal~$\Ical \subseteq \Pcal$
such that~$C \in \Ical$. In particular, $I \in \Ical$.
Suppose for the sake of contradiction that $\Ical$ contains a solution~$Y$ to~$I$.
Then~$Y \oplus C \in \Ical$ by the closure under the effective join,
contradicting our choices of~$C$ and~$I$.
\end{proof}

Weakness preservation admits a dual presentation in terms of \emph{avoidance}.

\index{avoidance}
\begin{definition}[Avoidance]
Let~$\Qcal$ be a collection of sets upward-closed under the Turing reducibility.
A $\Pi^1_2$ statement~$\Psf$ admits \emph{$\Qcal$ avoidance} if
for every set~$C \not \in \Qcal$ and every $C$-computable $\Psf$-instance~$X$,
there is a solution~$Y$ to~$X$ such that~$Y \oplus C \not \in \Qcal$.
\end{definition}

The two notions are trivially equivalent, since $\Qcal$ avoidance is~$\{Z : Z \not \in \Qcal\}$ preservation
and~$\Pcal$ preservation is~$\{Z : Z \not \in \Pcal\}$ avoidance.
There is however an ontological difference between those two notions.
Weakness properties such as lowness talk intuitively about how easy to compute a set is,
whereas the sets satisfying an avoidance property can be arbitrarily complicated.
In general, one will prefer to use a presentation such that the considered property is countable.
For example, there are countably many low sets, so it is more natural to define lowness preservation rather than non-lowness avoidance.
On the other hand, the upper cone of a non-computable set is countable, so one will rather use cone avoidance to build
a Turing ideal avoiding the halting set.

We shall often consider families of preservations or avoidances.

\index{avoidance!of cone}
\begin{example}[Cone avoidance]
A $\Pi^1_2$ statement~$\Psf$ admits \emph{cone avoidance} if for every pair of sets~$A$ and~$C$ such that $A \not \leq_T C$,
every~$C$-computable~$\Psf$-instance~$X$ has a solution~$Y$ such that~$A \not \leq_T Y \oplus C$.
In other words, $\Psf$ admits cone avoidance if for every set~$A$,
$\Psf$ admits~$\Qcal_A$ avoidance, where~$\Qcal_A = \{Z : A \leq_T Z\}$.
\end{example}

\chapter{Effective forcing}\label{chap:introduction-effective-forcing}

We now introduce the notion of effective forcing, which is one of the main tools to construct solutions
to theorems in reverse mathematics.
The forcing framework is a technique introduced by Cohen in set theory
to prove consistency and independence results~\cite{Cohen1963independence}. His technique was a major breakthrough
and  has been successfully re-applied in set theory and computability theory over the following years.
In this chapter, we present an effectivization of the forcing framework to fit the purposes
of computability theory.

Consider the following problem: Let~$\Mcal$ be a \emph{universe} satisfying some set of properties~$T$.
By ``universe'', we mean a collection of objects. Typically, $\Mcal$ will be a Turing ideal, or a structure,
and~$T$ will be a collection of sentences such that~$\Mcal \models T$.
We would like to add some elements to~$\Mcal$, so that the new universe~$\Ncal$
satisfies some extra properties~$S$, while preserving~$T$.
This problem basically occurs when we try to separate a statement~$\Psf$ from a statement~$\Qsf$,
by constructing a Turing ideal in which $\Psf$ holds, but not~$\Qsf$.
At some stage~$s$, we want to add a solution to some chosen~$\Psf$-instance, 
while avoiding adding solutions to a $\Qsf$-instance.
The forcing framework enables one to create new objects while controlling 
the propagation of the properties from the \emph{ground} universe to the new one.

\section{Basic notions}

\index{forcing}
\index{condition}
\index{forcing!extension}
A \emph{forcing notion} is a partial order~$(P, \preceq)$, whose elements
are called \emph{conditions}. A condition~$c$ informally represents
a partial approximation of the object we are constructing. It often comes
with a natural notion of \emph{satisfiability}. An object~$G$ satisfies a condition~$c$
if it belongs to the candidate objects approximated by~$c$.
A condition~$d$ \emph{extends} a condition~$c$ (written~$d \preceq c$)
if the approximation of~$d$ is more precise than the one of~$c$.
In other words, $d \preceq c$ if the collection of objects satisfying the
condition~$d$ is a subset of the ones satisfying~$c$, hence the reverse notation for 
the extension relation.

Given a formula~$\varphi(G)$ with only one distinguished set parameter~$G$, 
we can represent the property~$\varphi(G)$ by the set~$D_\varphi \subseteq P$ of all conditions~$c$
such that~$\varphi(H)$ holds for every set~$H$ satisfying~$c$.
A set~$D \subseteq P$ is \emph{dense} if every condition in~$P$ has an extension in~$D$.
Therefore, if the set~$D_\varphi$ for some property~$\varphi(G)$ is dense, 
then whatever the partial approximation~$c$ we have chosen so far,
it is always possible to find a refinement of the approximation which furthermore satisfies the property~$\varphi(G)$.

\begin{example}[Tree forcing]\label{ex:intro-forcing-tree-forcing}
Let~$T \subseteq 2^{<\omega}$ be a computable, infinite binary tree for which we want to construct an infinite path.
Consider the forcing notion~$(P, \preceq)$ where~$P$ is the collection of all computable infinite subtrees of~$T$.
A set~$G$ satisfies a condition~$U \subseteq T$ if it is a path through~$U$. A condition~$V$
extends~$U$ if~$V \subseteq U$. Note that~$[V] \subseteq [U]$. 
By taking a subtree~$V$, one removes some paths and therefore refines the approximation of~$U$.
For each index~$e$, consider the following sets
\[
D_{e,\downarrow} = \{ U \in P : (\exists s)(\forall \sigma \in U^{[s]})\Phi^\sigma_e(e) \downarrow\}
\hspace{20pt} D_{e, \uparrow} = \{ U \in P : (\forall \sigma \in U)\Phi^\sigma_e(e) \uparrow \}
\]
where~$U^{[s]} = \{\sigma \in U : |\sigma| = s \}$. The set~$D_e = D_{e,\downarrow} \cup D_{e,\uparrow}$
is dense. In particular, if~$U \in D_e$, then we have \emph{decided} the termination of~$\Phi^G_e(e)$
for every set~$G$ satisfying~$U$.
\end{example}

\index{filter}
\index{generic!filter}
\index{generic!set}
A \emph{filter} is a set~$F \subseteq P$ which is closed upward under the extension relation,
i.e., if $d \in \Fcal$ and~$d \preceq c$, then~$c \in F$, and such that any two elements
are compatible, that is, if~$c, d \in F$, then there is some~$e \in F$ such that $e \preceq c$ and~$e \preceq d$.
A filter is~\emph{$\Dcal$-generic} for a countable collection of dense sets~$\Dcal$
if it intersects every element of~$\Dcal$. Such a collection~$\Dcal$ represents the collection
of properties we would like to satisfy by the constructed object.

\begin{theorem}\label{thm:intro-forcing-filter-existence}
Let $(P, \preceq)$ be a notion of forcing, let~$\Dcal$ be a countable collection
of dense subsets of~$P$, and let~$c \in P$. There is a $\Dcal$-generic filter containing~$c$.
\end{theorem}

The collection~$\Dcal$ is often not explicited and we argue that
every \emph{sufficiently generic} filter satisfies some property if, whenever we put sufficiently many
dense sets in~$\Dcal$, every $\Dcal$-generic filter satisfies the property.
When the filter is sufficiently generic, there is often a unique object satisfying the intersection
of all its members, in which case we skip the filter notation and simply say that a set~$G$
is sufficiently generic for a forcing notion.

\begin{example}[Tree forcing]
Take the forcing notion of Example~\ref{ex:intro-forcing-tree-forcing}.
Consider the following set for each string~$\sigma \in 2^{<\omega}$:
\[
D_\sigma = \{U \in P : (\forall \tau \in U)[\tau \preceq \sigma \vee \sigma \preceq \tau] \}
\]
For each~$n$, the set~$D_n = \bigcup_{\sigma \in 2^n} D_\sigma$ is dense. In particular, every sufficiently generic filter~$F$
will intersect $D_\sigma$ for each~$n$ and exactly one~$\sigma \in 2^n$.
Therefore, the filter~$F$ will induce a unique real~$G_F = \bigcup \{ \sigma \in 2^{<\omega} : D_\sigma \cap F \neq \emptyset \}$.
\end{example}

\section{The forcing relation}

\index{forcing!relation}
One of the main features of the forcing framework is the ability to control
the properties of the resulting object during the construction.
We have interpreted a property $\varphi(G)$ as the set~$D_\varphi \subseteq P$
of all conditions~$c$ such that the property holds for every set satisfying~$c$.
This interpretation is however too restrictive, and in particular, the set~$D_\varphi \cup D_{\neg \varphi}$ is not dense in
general. For example, with the forcing notion defined in Example~\ref{ex:intro-forcing-tree-forcing},
the set
\[
\{U \in P : (\forall H \in [U])[\Phi^H_e \emph{ is total}] \vee (\forall H \in [U])[\Phi^H_e \emph{ is partial}] \}
\]
is not dense, although for every set~$G$ sufficiently generic, $\Phi^G_e$ will be either total or partial.

We therefore widen our satisfaction of a property, and say that a condition~$c$ \emph{forces} 
a formula~$\varphi(G)$ (written~$c \Vdash \varphi$) if~$\varphi(G)$ holds
for every sufficiently generic filter containing~$c$. In particular, the forcing relation satisfies
the following enjoyable property.

\begin{theorem}\label{thm:intro-forcing-relation-sound-complete}
Let~$F$ be a sufficiently generic filter, $G$ be the corresponding generic set
and~$\varphi$ be any formula where~$G$ is the only free variable.
Then $\varphi(G)$ holds iff there is a condition~$c \in F$ such that~$c \Vdash \varphi$.
\end{theorem}

This forcing relation is defined in a semantic way, but it also admits a purely syntactic definition.
We denote by~$c \Vdash^{*} \varphi(G)$ the syntactic forcing relation, which is defined inductively as follows.
The base case, that is, the forcing relation over~$\Sigma^0_0$ formulas, depends on the forcing notion we consider and in particular
the way we obtain a set~$G$ from a generic filter. It will always be clear in the context.
In Example~\ref{ex:intro-forcing-tree-forcing},~$U \Vdash^{*} \varphi$ holds for some $\Sigma^0_0$
formula~$\varphi(G)$ if there is some threshold~$t$ such that~$\varphi(\sigma)$ holds for every~$\sigma \in U$
of length~$t$. The forcing relation is extended to arbitrary formulas by the following rules.
\begin{itemize}
	\item[1.] $c \Vdash^{*} \neg \varphi$ iff $d \not \Vdash^{*} \varphi$ for all~$d \preceq c$
	\item[2.] $c \Vdash^{*} \varphi \wedge \psi$ iff $c \Vdash^{*} \varphi$ and~$c \Vdash^{*} \psi$
	\item[3.] $c \Vdash^{*} (\exists n)\varphi(n)$ iff for every~$d \preceq c$, there is some~$n \in \omega$
	and some~$e \preceq d$ such that~$e \Vdash^{*} \varphi(n)$.
\end{itemize}
\smallskip

By the following theorem, the syntactic and the semantic forcing relations coincide.

\begin{theorem}
Let $c$ be a condition and $\varphi$ be any formula where~$G$ is the only free variable.
Then $c \Vdash \varphi$ iff $c \Vdash^{*} \varphi$.
\end{theorem}

\section{Effective forcing}\label{sect:intro-forcing-effective-forcing}

\index{effective forcing}
\index{ngeneric@$n$-generic}
\index{weakly n-generic@weakly $n$-generic}
Theorem~\ref{thm:intro-forcing-filter-existence} admits an effective analogue if we impose effectiveness restrictions
on the collection~$\Dcal$. Given a set~$D \subseteq P$, let~$D^{\bot} = \{c \in P : (\forall d \preceq c) d \not \in D\}$.
We cannot effectively enumerate all dense~$\Sigma^0_n$ sets. However, for every~$\Sigma^0_n$ set~$D$,
the set~$D \cup D^{\bot}$ is dense. We say that a filter~$F \subseteq P$ is \emph{$n$-generic}
if it intersects~$D \cup D^{\bot}$ for every $\Sigma^0_n$ set~$D \subseteq P$. 
The filter~$F$ \emph{meets} $D$ if it intersects~$D$, and~\emph{avoids} $D$
if it intersects~$D^{\bot}$.
A filter~$F \subseteq P$ is \emph{weakly $n$-generic} if it intersects every dense $\Sigma^0_n$ set~$D \subseteq P$.

\index{Cohen forcing}
Cohen notion of forcing~$(2^{<\omega}, \succeq)$ is of particular interest. 
The conditions are strings~$\sigma \in 2^{<\omega}$, and a string~$\tau$ extends~$\sigma$
if $\sigma$ is a prefix of~$\tau$. Cohen genericity is a very useful notion
as it represents, together with Martin-L\"of randomness, the informal idea of~\emph{typical set}. 
Genericity is typical in a Baire category sense, while Martin-L\"of randomness is typical in a measure-theoretic sense.
Weakly 1-genericity for Cohen forcing coincides with the notion of hyperimmunity that we already met.

\begin{theorem}[Kurtz~\cite{Kurtz1982Randomness,Kurtz1983Notions}]
A degree~$\dbf$ bounds a Cohen weakly 1-generic real iff $\dbf$ is hyperimmune.
\end{theorem}
\begin{proof}[Proof idea]
Given a string~$\sigma$, let~$p_\sigma$ denote its principal function.
For every computable function, the dense set~$D_f = \{\sigma : (\exists n)p_\sigma(n) > f(n)\}$
witnesses the hyperimmunity of weakly 1-generic reals. 
In the other direction, take a hyperimmune function~$f$ and construct a weakly 1-generic real by a greedy algorithm, 
looking at stage~$s$ at the~$f(s)$th approximation of the c.e.\ sets of strings.
\end{proof}

Cohen forcing is a very simple notion of forcing, both conceptually
and effectively. In general, one will want to use the most effective notion of forcing possible,
that is, a forcing notion whose conditions and extension relation are effectively describable, 
in order to make the generic set preserve more properties.

\begin{example}[Tree forcing]\label{ex:intro-forcing-tree-forcing-cone}
Suppose we want to prove the cone avoidance basis theorem (Theorem~\ref{thm:intro-comp-cone-avoidance-basis}).
Let~$T$ be a computable, infinite binary tree and~$A$ be a non-computable set.
Consider the forcing notion of Example~\ref{ex:intro-forcing-tree-forcing}.
It suffices to prove that for every index~$e$, the following set is dense.
\[
D_e = \{ U \in P : (\exists n)(\forall \sigma \in U)[\Phi_e^{\sigma}(n) \downarrow \imp \Phi_e^{\sigma}(n) \neq A(n)] \}
\]
Fix some infinite binary tree~$U \in P$ and let~$\Gamma_n = \{ \sigma \in 2^{<\omega} : \Phi_e^\sigma(n) \uparrow \}$.
By convention on the notation~$\Phi^\sigma_e(\cdot)$, the set~$U \cap \Gamma_n$ is a computable binary tree. We claim that one of the two following must hold:
\begin{itemize}
	\item[1.] $\Phi^\sigma_e(n) \downarrow \neq A(n)$ for some~$n$ and some extendible $\sigma \in U$
	\item[2.] $U \cap \Gamma_n$ is infinite for some~$n$
\end{itemize}
In the first case, the tree $U^{[\sigma]} = \{\tau \in U : \tau \succeq \sigma\}$ is in~$D_e$,
and in the second case, $U \cap \Gamma_n \in D_e$. If none of 1 and 2 hold, then
we can compute~$A$, by a procedure which on input~$n$, searches for some threshold~$t$ and some value~$x$
such that $\Phi^\sigma_e(n) \downarrow = x$ for every string~$\sigma \in U$ of length~$t$.
This contradicts the assumption that~$A$ is not computable.
\end{example}

The core of the argument is the following. We want to diagonalize against a complex object (here, the set~$A$)
while working within a weak forcing notion (computable trees). If there is no way to diagonalize against the object,
then the forcing notion gives us some grasp on it, and therefore the object admits a simple description.
In Example~\ref{ex:intro-forcing-tree-forcing-cone}, if we used arbitrary infinite binary trees 
instead of computable ones, there would be no contradiction
in obtaining an $U$-computation of~$A$.

The importance of an effective description of the partial order is even more visible when
we consider the forcing relation. In order to construct a set~$G$ of low${}_n$ degree,
that is, such that~$G^{(n)} \leq_T \emptyset^{(n)}$, one need to decide $\Sigma^0_n$ formulas
in a $\emptyset^{(n)}$-effective construction. In particular,
one want to $\emptyset^{(n)}$-effectively decide whether, given a condition~$c \in P$
and a $\Sigma^0_n$ formula~$\varphi(G)$, there is an extension~$d \Vdash^{*} \varphi$
or $d \Vdash^{*} \neg \varphi$. The relation~$d \Vdash^{*} \varphi$ for $\Sigma^0_0$ formulas
is usually decidable. However, the third inductive rule makes the relation~$d \Vdash^{*} \varphi$ 
too complex when dealing with higher rank formulas.
We need to define a stronger relation of forcing which has better definitional properties.
Define~$c \Vdash^{\circ} \varphi$ in the same way as~$c \Vdash^{*} \varphi$ for $\Sigma^0_0$ formulas,
and add the following inductive rules.
\begin{itemize}
	\item[(a)] $c \Vdash^{\circ} (\forall n)\varphi(n)$ iff $d \Vdash^{\circ} \varphi(n)$ for all~$n \in \omega$ and all~$d \preceq c$
	\item[(b)] $c \Vdash^{\circ} (\exists n)\varphi(n)$ iff $c \Vdash^{\circ} \varphi(n)$ for some~$n \in \omega$.
\end{itemize}
\smallskip

The relation~$\Vdash^{\circ}$ is known as the \emph{strong forcing relation} and was introduced by Cohen~\cite{Cohen1963independence}. 
Note that, unlike the previous forcing relation, the strong forcing relation is not semantic.
In particular, there are some logically equivalent formulas~$\varphi$ and~$\psi$
such that~$c \Vdash^{\circ} \varphi$ and~$c \not \Vdash^{\circ} \psi$. See Mummert~\cite{MummertIs}
for a discussion. Thankfully, the strong forcing relation also satisfies Theorem~\ref{thm:intro-forcing-relation-sound-complete}.
It relates to the weak forcing relation as follows.

\begin{theorem}
Let $c$ be a condition and $\varphi$ be any formula where~$G$ is the only free variable.
Then $c \Vdash \varphi$ iff $c \Vdash^{\circ} \neg \neg \varphi$.
\end{theorem}

Looking at Cohen forcing~$(2^{<\omega}, \succeq)$, is it clear that the set~$\{\sigma : \sigma \Vdash^{\circ} \varphi\}$
for some~$\Sigma^0_n$ ($\Pi^0_n$) formula~$\varphi$ is $\Sigma^0_n$ ($\Pi^0_n$).
However, when considering more complicated forcing notions, the complexity of the forcing relation
can increase due to the complexity of the quantification over all extensions of a given condition.
This is in particular the case for \emph{Mathias forcing}, which we introduce below.

From now on, we will consider only the strong forcing relation and write it~$\Vdash$.

\section{Mathias forcing}\label{sect:introduction-forcing-mathias-forcing}

Mathias forcing is a central notion of forcing in the computable analysis of Ramsey's theorem
and its consequences. Variants of Mathias forcing have been successfully used to separate
various Ramsey-type theorems.

\index{Mathias forcing}
\index{forcing!Mathias}
\index{condition!Mathias}
\begin{definition}[Mathias forcing]
A \emph{Mathias condition} is an ordered pair~$(F, X)$, where~$F$ is a finite set of integers
and $X$ is an infinite set such that~$max(F) < min(X)$.
A condition~$d = (E, Y)$ \emph{extends} $c = (F, X)$ (written~$d \leq c$)
if $F \subseteq E$, $Y \subseteq X$ and $E \setminus F \subset X$.
A set $G$ \emph{satisfies} a Mathias condition $(F, X)$
if $F \subset G$ and $G \setminus F \subseteq X$.
\end{definition}

One may think of a Mathias condition~$(F, X)$ 
as a finite approximation~$F$ of the generic set~$G$,
together with an infinite \emph{reservoir}~$X$,
that is, a set of candidate elements we may add to $F$ later on.
The restriction~$max(F) < min(X)$ ensures that~$F$ is an initial segment
of the set~$G$ (when we view sets as binary strings/sequences).

The shape of a notion of forcing strongly depends on the nature of the object
one wants to construct. It is therefore common to use slight variations of Mathias forcing
to force solutions to Ramsey-type statements. A reservoir is not only a restriction
of the domain, but also a warrant of the extensibility of a condition.

\begin{example}[Mathias forcing for~$\rt^2_2$]\label{ex:intro-forcing-mathias-rt22}
Let~$f : [\omega]^2 \to 2$ be a computable coloring for which we want to construct
an infinite $f$-homogeneous set. First, note that a finite $f$-homogeneous set~$F$ may not be extendible
into an infinite one. For example, a set~$F$ which is $f$-homogeneous for color~$0$,
and such that~$f(x, y) = 1$ for some~$x \in F$ and cofinitely many~$y$ is not extendible.

The suitable forcing notion to obtain extendible conditions is the following variant of Mathias forcing:
The conditions are tuples~$(F_0, F_1, X)$ such that $(F_0, X)$ and~$(F_1, X)$ form Mathias conditions.
The conditions furtheremore satisfy the following property:
\[
(\forall i < 2)(\forall y \in X)[F_i \cup \{y\} \mbox{ is $f$-homogeneous for color } i]
\]
This property over the reservoir ensures that at least one of the finite $f$-homogeneous sets~$F_0$ or~$F_1$
is extendible into an infinite one. To see this, apply~$\rt^2_2$ over the coloring~$f$ restricted to the reservoir~$X$,
to obtain an infinite $f$-homogeneous set~$H \subseteq X$, say for color~$i$. Then the set~$F_i \cup H$
is also $f$-homogeneous.
\end{example}

As explained in the previous section, one usually wants to impose some effectiveness
property over the notion of forcing. In the case of Mathias forcing, the reservoir~$X$
often satisfies some \emph{weakness property}.
Suppose we want to prove that some statement~$\Psf$ admits $\Pcal$ preservation
for some weakness property~$\Pcal$. Let~$C \in \Pcal$ and~$I$ be a $C$-computable $\Psf$-instance.
One usually works with a Mathias forcing $(F, X)$ such that~$X \oplus C \in \Pcal$.
In particular, the absence of a solution~$Y$ to~$I$ such that~$Y \oplus C \in \Pcal$
ensures the extendibility of all finite approximations.

\begin{example}[Mathias forcing for~$\rt^2_2$]
Suppose we want to prove that~$\rt^2_2$ admits cone avoidance.
Let~$A$ be a non-computable set, and let~$f : [\omega]^2 \to 2$ be a computable coloring.
If there is an infinite $f$-homogeneous set~$H$ such that~$A \not \leq_T H$,
then we are done. So suppose it is not the case.
Consider the variant of Mathias forcing introduced in Example~\ref{ex:intro-forcing-mathias-rt22},
where we furthermore impose that~$A \not \leq_T X$.

For every condition~$c = (F_0, F_1, X)$ and every~$i < 2$, there is an extension~$d = (E_0, E_1, Y)$
such that~$|E_i| > |F_i|$. By symmetry, say~$i = 0$. Ask whether there is some~$x \in Y$ such that
the set~$Y = \{ y \in X : f(x, y) = 0 \}$ is infinite. If so, the condition~$(F_0 \cup \{x\}, F_1, Y \setminus [0, x])$
is the desired extension. If there is no such~$x \in Y$, then for almost all~$x$ and~$y \in Y$,
$f(x,y) = 1$. By thinning out the set~$X$, one may obtain an $X$-computable $f$-homogeneous set~$H$.
In particular, $A \not \leq_T H$, contradicting our assumption.

Of course, we need to prove that the sets forcing the generic~$G$ not to compute~$A$
are dense, but this goes beyond an introduction to effective forcing.
\end{example}

Mathias forcing is a good notion of forcing to control properties
below the first jump. In particular, the forcing relation over~$\Sigma^0_1$
formulas is $\Sigma^0_1$. However, the situation gets more complicated when controlling
iterated jumps. Cholak, Dzhafarov and Hirst~\cite{Cholak2014Generics} studied the complexity
of the forcing relation for computable Mathias forcing, that is, a Mathias forcing
whose reservoirs are computable.

\begin{theorem}[\cite{Cholak2014Generics}]
Let~$(F, X)$ be a computable Mathias condition and let~$\varphi(G)$ be a formula with exactly one free set variable.
\begin{itemize}
	\item[(1)] If~$\varphi$ is~$\Sigma^0_0$ then the relation~$(F, X) \Vdash \varphi$ is computable.
	\item[(2)] If~$\varphi$ is~$\Pi^0_1$, $\Sigma^0_1$ or~$\Sigma^0_2$ then so is the relation~$(F, X) \Vdash \varphi$.
	\item[(3)] For~$n \geq 2$, if $\varphi$ is~$\Pi^0_n$ then the relation~$(F, X) \Vdash \varphi$ is~$\Pi^0_{n+1}$.
	\item[(4)] For~$n \geq 3$, if~$\varphi$ is~$\Sigma^0_n$ then the relation~$(F, X) \Vdash \varphi$ is~$\Sigma^0_{n+1}$.
\end{itemize}
\end{theorem}

The extra complexity of the forcing relation is essentially due to the fact that the sentence~``$Y$
is an infinite subreservoir of~$X$'' is hard to describe.
The separation between Ramsey-type statements can be usually done by
computability-theoretic properties which admit a diagonalization
below the first jump. In this case, Mathias forcing is a suited notion of forcing.
However, in some cases, one needs to control iterated jumps, in particular
when dealing with hiearchies such as the free set and the thin set hierarchies.
We shall see in Chapter~\ref{chap-controlling-iterated-jumps} alternative notions of forcing
for various Ramsey-type theorems, whose complexity of the forcing relation
admits the same bounds as for Cohen forcing.

Finally, Cholak, Jockusch and Slaman~\cite{Cholak2001strength} 
noted that every sufficiently computable Mathias generic produces a set of high degree.

\begin{lemma}[\cite{Cholak2001strength}]
Every set sufficiently generic for computable Mathias forcing
is of high degree.
\end{lemma}
\begin{proof}
For every computable function~$f$, the set~$D_f = \{ (F, X) : p_{F \cup X} \mbox{ dominates } f \}$
is dense. In particular, for every infinite set~$G$ satisfying some condition~$(F, X) \in D_f$,
$p_G$ dominates~$f$. 
\end{proof}

One can however restrict the amount of genericity of the filter to obtain sets
of lower degree, as did Cholak, Jockusch and Slaman for Ramsey's theorem for pairs.

\part{The strength of Ramsey's theorem}

\chapter{A proof of Ramsey's theorem}

Ramsey's theory claims that true chaos does not exist\footnote{Ramsey never went into my room.},
in the sense that one can always find some substructure in a sufficiently large collection of arbitrary data.
Ramsey's theorem is perhaps the most famous example of such a statement and asserts
that every coloring of the~$n$-tuples of integers into $k$ colors has an infinite monochromatic set.

The strength of Ramsey's theorem and its consequences is notoriously hard to tackle.
Its investigation led to a number of new separation techniques and 
contributed to increase the global understanding of reverse mathematics.
The effective analysis of Ramsey’s theorem was started by Jockusch~\cite{Jockusch1972Ramseys}
who gave precise bounds on the definitional complexity of solutions to Ramsey's theorem.

\begin{theorem}[Jockusch]\label{thm:intro-rev-ramsey-bounds} Fix some~$n, k \geq 2$.
\begin{itemize}
	\item Every computable coloring $f : [\Nb]^n \to k$ has a $\Pi^0_n$ homogeneous set.
	\item There is a computable coloring $f : [\Nb]^n \to k$ with no $\Sigma^0_n$ homogeneous set.
\end{itemize}
\end{theorem}

In particular, Theorem~\ref{thm:intro-rev-ramsey-bounds} shows that~$\rt^n_k$ is a consequence of~$\aca$ in standard models since
$\aca$ can be understood as the statement ``Every set has a Turing $n$th jump''.
There are many ways to prove Ramsey's theorem. The simplest one is an inductive argument using \emph{prehomogeneous sets}.

\index{prehomogeneous set}
\begin{definition}
Let~$f : [\Nb]^{n+1} \to k$ be a coloring. A set~$P$ is \emph{$f$\nbd prehomogeneous}
if for every~$\sigma \in [P]^n$ and every~$x,y \geq max(\sigma)$, $f(\sigma, x) = f(\sigma, y)$.
\end{definition}

In other words, a set~$P$ is $f$-prehomogenous if the color of a tuple
over~$P$ do not depend on its last component. In particular, every $f$-homogeneous
set is $f$-prehomogeneous. From a computational point of view, Jockusch~\cite{Jockusch1972Ramseys}
showed that K\"onig's lemma ($\kl$) proves the existence of a prehomogeneous set.
Simpson~\cite{Simpson2009Subsystems} formalized Jockusch's argument in~$\rca$.

\begin{lemma}[($\rca + \kl$, $\leq_{sc} \kl$)]\label{lem:strength-kl-prehom}
Every coloring~$f : [\Nb]^{n+1} \to k$ has an infinite $f$\nbd prehomogenous set.
\end{lemma}
\begin{proof}
Let~$T \subseteq \Nb^{<\Nb}$ be the set of all strings~$a_0 < a_1 < \dots < a_{t-1}$
such that $\{a_0, \dots, a_{t-1}\}$ is $f$-prehomogeneous and
for each~$i < t$, $a_i$ is the least $a > a_{i-1}$ such that
$f(\sigma, a_i) = f(\tau, a_i)$ for each~$\sigma, \tau \in [\{a_0, \dots, a_{i-1}\}]^n$ (where~$a_{-1} = -1$).
The set~$T$ is an infinite, $\Delta^{0,T}_1$ finitely branching tree
whose paths are $f$-prehomogeneous.
\end{proof}

Since K\"onig's lemma holds in~$\rca + \aca$, we can use the notion of prehomogeneous
set to reduce~$\rt^{n+1}_k$ to~$\rt^n_k$.

\begin{theorem}\label{thm:strength-aca-rtn}
For every standard~$n$, $\rca \vdash \aca \imp \rt^n_{<\infty}$.
\end{theorem}
\begin{proof}
Assume by induction over~$n$ that $\rca \vdash \aca \imp \rt^n_{<\infty}$.
Let~$f : [\Nb]^{n+1} \to k$ be a coloring. By Lemma~\ref{lem:strength-kl-prehom},
there is an infinite $f$-prehomogeneous set~$P$. Let~$\tilde{f} : [P]^n \to k$
be such that for each~$\sigma \in [P]^n$, $\tilde{f}(\sigma) = f(\sigma, x)$ for
the least~$x \in P$ such that~$x > max(\sigma)$. By~$\rt^n_{<\infty}$,
there is an infinite $\tilde{f}$-homogeneous set~$H$. The set~$H$ is $f$-homogeneous.
\end{proof}

Note that the induction in Lemma~\ref{thm:strength-aca-rtn} lives in the meta theory.
In fact, Simpson~\cite{Simpson2009Subsystems} proved that the statement
$\rt$ does not hold in~$\rca + \aca$. Of course, $\rt$ holds in any $\omega$-model of~$\aca$.

Is the proof of~$\rt^n_{<\infty}$ optimal? The use of prehomogeneous sets
as intermediary objects is done without loss of optimality since every homogeneous
set is prehomogeneous. However, the resort to K\"onig's lemma to prove the existence
of a prehomogeneous set may be unnecessary. In fact, Hirschfeldt and Jockusch~\cite{Hirschfeldtnotions}
showed that whenever~$n \geq 3$, the use of K\"onig's lemma is optimal, in the following sense.

\begin{theorem}[$\rca$]
For every set~$X$, there is a $\Delta^{0,X}_1$ coloring~$f : [\Nb]^3 \to 2$ such that every $f$-prehomogeneous set
is of PA degree relative to the jump of~$X$.
\end{theorem}

In particular, $\rca \vdash \rt^3_2 \imp \kl$ and $\kl \leq_{sc} \rt^3_2$.
It follows that, $\rt^n_k$ and $\rt^n_{<\infty}$ are equivalent to~$\aca$ over~$\rca$ for every~$n \geq 3$
and~$k \geq 2$.

\section{Cohesiveness and stability}

Can we do better for Ramsey's theorem for pairs? Let us have a closer look
at the standard proof of~$\rt^2_2$ using prehomogeneous sets.
Let~$f : [\Nb]^2 \to k$ be a coloring and let~$P$ be an infinite $f$-prehomogeneous set.
The set~$P$ induces a $\Delta^{0,f \oplus P}_1$ partition~$A_0, \dots, A_{k-1}$ of~$P$ defined by
\[
A_i = \{ x \in P : (\forall y \in P)[y > x \imp f(x, y) = i] \}
\] 
Note that any infinite subset of one of the~$A$'s is a solution to~$f$.
Apply~$\rt^1_k$ to obtain an infinite $f$-homogeneous set.
In their seminal paper, Cholak, Jockusch and Slaman~\cite{Cholak2001strength} had the idea
to use a stronger variant of the pigeonhole principle in order to weaken the properties
expected from the set~$P$.

\index{d@$\mathsf{D}^n_k$}
\begin{definition}
$\mathsf{D}^n_k$ is the statement ``For every~$\Delta^0_n$ $k$-partition~$A_0, \dots, A_{k-1}$ of~$\Nb$,
there is an infinite subset of one of the~$A$'s''.
\end{definition}

By Shoenfield's limit lemma, for every $\Delta^0_2$ such $k$-partition~$A_0, \dots, A_{k-1}$ of~$\Nb$,
there is a computable function~$f : [\Nb]^2 \to k$ such that for every~$x$, $\lim_y f(x, y)$ exists
and~$x \in A_i$ iff $\lim_y f(x,y) = i$. The function~$f$ can be seen as a particular
$\rt^2_k$-instance such that $\lim_y f(x,y)$ exists for every~$x$. Such a function is called~\emph{stable}
and we denote by~$\srt^2_k$ the statement~$\rt^2_k$ restricted to stable functions.
In particular, every infinite $f$-homogeneous set for color~$i$ is an infinite subset of~$A_i$.
On the other hand, we can $f$-computably thin out any infinite subset of~$A_i$ into an infinite $f$-homogeneous set
for color~$i$. The statements~$\mathsf{D}^2_k$ and~$\srt^2_k$ are therefore computably equivalent.
Cholak, Jockusch and Slaman~\cite{Cholak2001strength} formalized the equivalence over~$\rca$. 
However, they used the~$\Sigma^0_2$ bounding scheme ($\bst$)
to compute an infinite $f$-homogeneous set from a subset of one of the~$A$'s.
Later, Chong, Lempp and Yang~\cite{Chong2010role} showed that~$\mathsf{D}^2_2$ implies~$\bst$ over~$\rca$
using an involved argument.

\begin{theorem}\label{thm:strength-d22-srt22}
For every~$k$,
$\rca \vdash \mathsf{D}^2_k \biimp \srt^2_k$ and~$\mathsf{D}^2_k =_c \srt^2_k$.
\end{theorem}

Given a coloring~$f : [\Nb]^2 \to k$, the restriction~$f : [P]^2 \to k$ is stable
for every $f$-prehomogeneous set~$P$. In fact, every set which is \emph{eventually} prehomogeneous
satisfies this property. Let~$\crt^2_2$ be the statement asserting for every coloring~$f : [\Nb]^2 \to 2$
the existence of an infinite set~$C$ such that~$f : [C]^2 \to 2$ is stable.
One immediately sees that $\rca \vdash \rt^2_2 \biimp [\crt^2_2 \wedge \srt^2_2]$.
However, we shall rather use the following seemingly stronger statement.
Given two sets~$A$ and~$B$, the notation $A \subseteq^{*} B$ means that
the set~$A$ is included in~$B$, \emph{up to finite changes}.

\index{cohesiveness}
\index{coh@$\coh$|see {cohesiveness}}
\index{p-cohesive}
\begin{definition}[Cohesiveness]
An infinite set $C$ is $\vec{R}$-cohesive for a sequence of sets $R_0, R_1, \dots$
if for each $i \in \omega$, $C \subseteq^{*} R_i$ or $C \subseteq^{*} \overline{R_i}$.
A set $C$ is \emph{p-cohesive} if it is $\vec{R}$-cohesive where
$\vec{R}$ is an enumeration of all primitive recursive sets.
$\coh$ is the statement ``Every uniform sequence of sets $\vec{R}$
has an $\vec{R}$-cohesive set.''
\end{definition}

Cohesiveness implies~$\crt^2_2$ over~$\rca$ as shows the following simple reduction.

\begin{lemma}\label{lem:strength-coh-to-crt22}
$\rca \vdash \coh \imp \crt^2_2$ and $\crt^2_2 \leq_{sW} \coh$
\end{lemma}
\begin{proof}
Let~$f : [\Nb]^2 \to 2$ be a coloring.
Define the~$\Delta^{0,f}_1$ sequence of sets~$R_0, R_1, \dots$ for each~$x$ as
$R_x = \{ y \in \Nb : f(x, y) = 1 \}$
and let~$C$ be an $\vec{R}$-cohesive set. The function~$f : [C]^2 \to 2$ is stable.
\end{proof}

At first sight, cohesiveness seems to be more than what we require
in the proof of Lemma~\ref{lem:strength-coh-to-crt22}.
Indeed, any infinite set~$C$ such that for every~$x \in C$, $C \subseteq^{*} R_x$ or~$C \subseteq^{*} \overline{R}_x$
would be sufficient. Mileti~\cite{Mileti2004Partition} proved that cohesiveness can chosen without loss of optimality. 

\begin{lemma}\label{lem:strength-rt22-to-coh}
$\rca \vdash \rt^2_2 \imp \coh$ and~$\coh \leq_{sW} \rt^2_2$.
\end{lemma}
\begin{proof}[Proof idea]
Let~$R_0, R_1, \dots$ be a sequence of sets. By adding some dummy sets,
we can assume that for each~$x \neq y$, there is some~$i$ such that~$R_i(x) \neq R_i(y)$.
Let~$i(x,y)$ be the least such~$i$. Define~$f : [\Nb]^2 \to 2$ for each~$x < y$
by~$f(x, y) = 1$ iff~$x \in R_{i(x,y)}$. Any infinite $f$-homogeneous set is~$\vec{R}$-cohesive.
However, one must be careful with the induction to carry the proof over~$\rca$.
\end{proof}

Later, Hirschfeldt and Shore~\cite{Hirschfeldt2007Combinatorial} clarified the links between $\crt^2_2$ and~$\coh$ 
by proving that they coincide over~$\rca + \bst$.
Putting Theorem~\ref{thm:strength-d22-srt22},
Lemma~\ref{lem:strength-coh-to-crt22} and Lemma~\ref{lem:strength-rt22-to-coh} 
together, we obtain the following equivalence.

\begin{theorem}\label{thm:strength-rt22-coh-d22}
$\rca \vdash \rt^2_2 \biimp [\coh \wedge \dsf^2_2]$. 
\end{theorem}

\section{Effective constructions}

The effective analysis of Ramsey's theorem for pairs now reduces to 
the analysis of cohesiveness and the $\dsf^2_2$ statement.
The former one will be extensively studied in Chapter~\ref{chap:degrees-unsolvability-cohesiveness}.
We now proceed to an effective analysis of $\dsf^2_2$.
Most of the constructions we will encounter in this thesis will be $\Pcal$ preservations
for some weakness property~$\Pcal$ which is a genericity notion,
that is, for which we can prove that every sufficiently Mathias generic
preserves $\Pcal$. For example, hyperimmunity and cone avoidance are genericity notions
whereas lowness is not. For any such weakness property $\Pcal$, we do not care about the effectiveness
of the construction. In particular, a proof of $\Pcal$ preservation for~$\dsf^2_2$
often happens to be a proof of strong $\Pcal$ preservation for~$\rt^1_2$.

\index{preservation!of weakness (strong)}
\begin{definition}[Strong weakness preservation]
A $\Pi^1_2$ statement~$\Psf$ admits \emph{strong $\Pcal$ preservation} for some weakness property~$\Pcal$,
if for every set~$C \in \Pcal$ and every (non-necessarily $C$-computable) $\Psf$-instance~$X$,
there is a solution~$Y$ to~$X$ such that~$Y \oplus C \in \Pcal$.
\end{definition}

Strong $\Pcal$ preservation represents a \emph{structural} weakness of the statement~$\Psf$
in that it is not possible to encode in a $\Psf$-instance the amount of information to escape~$\Pcal$.
On the other hand, the standard $\Pcal$ preservation represents an \emph{effective} weakness,
in that the weakness of the solutions may come from the effectiveness restrictions on the~$\Psf$-instance.
Some statements like~$\wkl$ are effectively weak while being structurally strong.
Indeed, by the cone avoidance basis theorem, $\wkl$ admits cone avoidance, whereas 
there is an infinite binary tree whose unique path computes the halting set.

Since any $\dsf^2_2$-instance can be seen as a non-effective~$\rt^1_2$-instance,
if~$\rt^1_2$ admits strong~$\Pcal$ preservation then $\dsf^2_2$ admits $\Pcal$ preservation.
In particular, the proofs of $\Pcal$ preservation for~$\rt^2_2$ are often based
on the following lemma which holds by Theorem~\ref{thm:strength-rt22-coh-d22}.

\begin{lemma}
If~$\rt^1_2$ admits strong~$\Pcal$ preservation and~$\coh$ admits $\Pcal$ preservation,
then~$\rt^2_2$ admits $\Pcal$ preservation.
\end{lemma}

We now present a proof sketch of strong $\Pcal$ avoidance for~$\rt^1_2$ 
for an arbitrary collection $\Pcal$ of subsets of~$\omega$.
Of course, such a proof will contain some holes which can be completed only for some particular~$\Pcal$'s.
Let~$A_0 \cup A_1 = \omega$ be an arbitrary $\rt^1_2$-instance.
Fix some set~$C$ belonging to some weakness property~$\Pcal$
and suppose that there is no infinite subset~$H$ of one of the~$A$'s
such that~$H \oplus C \in \Pcal$, as otherwise we are done.
We want to build a set~$G$ such that~$G \cap A_0$ and~$G \cap A_1$
are both infinite, and either~$(G \cap A_0) \oplus C \in \Pcal$
or~$(G \cap A_1) \oplus C \in \Pcal$.
For this, we use a variant of Mathias forcing~$(F, X)$ where $X \oplus C \in \Pcal$.
First, notice that Mathias generics are cohesive.

\begin{lemma}\label{lem:strength-mathias-generic-cohesive}
If~$G$ is sufficiently Mathias generic, then $G$ is cohesive.
\end{lemma}
\begin{proof}
Let~$R$ be a computable set and~$c = (F, X)$ be a Mathias condition.
The sets~$X \cap R$ and~$X \cap \overline{R}$ are both $X \oplus C$-computable,
and one of them is infinite. Therefore, either~$(F, X \cap R)$
or~$(F, X \cap \overline{R})$ is a valid extension forcing~$G \subseteq^{*} R$
in the former case, and~$G \subseteq^{*} \overline{R}$ in the latter case.
\end{proof}

Lemma~\ref{lem:strength-mathias-generic-cohesive} is one of the main reasons
why it is so difficult to separate Ramsey-type statements from~$\coh$.
Indeed, many Ramsey-type theorems are constructed using variants of Mathias forcing
whose generics are cohesive. It is still possible to prove that some statement~$\Psf$
does not imply~$\coh$ using Mathias forcing, but then some effectiveness
restriction has to be put on the overall construction.

Assuming that there is no infinite subset~$H$ of one of the~$A$'s
such that~$H \oplus C \in \Pcal$, we can prove that $G \cap A_0$
and~$G \cap A_1$ are both infinite if~$G$ is sufficiently generic.

\begin{lemma}\label{lem:strength-weakness-assumption-infinite} 
For every condition~$c = (F, X)$ and every~$i < 2$,
there is an extension~$d = (E, Y)$ to~$c$ such that~$|E \cap A_i| > |F \cap A_i|$.
\end{lemma}
\begin{proof}
Fix $c$ and~$i$. First suppose that the set~$X \cap A_i$ is empty. In this case~$X \subseteq A_{1-i}$,
is infinite, and is such that~$X \oplus C \in \Pcal$, contradicting our assumption.
So suppose~$X \cap A_i \neq \emptyset$, and let~$n \in X \cap A_i$.
The condition~$d = (F \cup \{n\}, X \setminus [0, n])$ is the desired extension.
\end{proof}

At this stage, it is important to clarify the difference between 
a generic for this notion of forcing and a solution to the $\rt^1_2$-instance~$A_0, A_1$.
Every sufficiently generic filter for this variant of Mathias forcing induces a single set~$G$
that we call the~\emph{generic}.
However, we are in fact constructing two sets~$G \cap A_0$ and~$G \cap A_1$,
which aim to be \emph{solutions} to the~$\rt^1_2$-instance for color~$0$ and~$1$, respectively.
One may wonder why we build two solutions in parallel, since Lemma~\ref{lem:strength-weakness-assumption-infinite}
shows that both can be forced to be infinite. The need for two solutions becomes clear
when trying to force the solutions to preserve~$\Pcal$.
Indeed, given a set~$H$, one can usually ensure that~$H \in \Pcal$ 
by satisfying a countable collection of requirements~$\Rcal^H_0, \Rcal^H_1, \dots$
In the case of strong~$\Pcal$ preservation for~$\rt^1_2$, we satisfy a countable collection
of disjunctive requirements $\Qcal_{0,0}, \Qcal_{0,1}, \Qcal_{1,0}, \dots$ defined for each~$e_0, e_1$ by
\[
 \Qcal_{e_0, e_1} = \Rcal^{(G \cap A_0) \oplus C}_{e_0} \vee \Rcal^{(G \cap A_1) \oplus C}_{e_1}
\]

By a simple pairing argument, if the $\Qcal$-requirements are all satisfied,
then all the~$\Rcal^{G \cap A_i}$-requirements are satisfied for some~$i < 2$,
in which case~$G \cap A_i$ is a solution to the~$A$'s such that~$(G \cap A_i) \oplus C \in \Pcal$.

\begin{example}[Cone avoidance]
Suppose we want to prove that~$\rt^1_2$ admits strong cone avoidance.
One can split the property~$A \not \leq_T H$ into the countable collection of requirements
$\Rcal^H_0, \Rcal^H_1, \dots$ where $\Rcal^H_e$ is ``$\Phi^H_e \neq A$''. 
Therefore, to prove that~$\rt^1_2$ admits strong cone avoidance, one will try to satisfy the following
requirements for each~$e_0, e_1 \in \omega$.
\[
 \Qcal_{e_0, e_1} : \Phi^{(G \cap A_0) \oplus C}_{e_0} \neq A \vee \Phi^{(G \cap A_1) \oplus C}_{e_1} \neq A
\]
If all the~$\Qcal$-requirements are satisfied, then either~$A \not \leq_T (G \cap A_0) \oplus C$
or~$A \not \leq_T (G \cap A_1) \oplus C$.
\end{example}

We now describe how to decide $\Sigma^{0,C}_1$ formulas over the solutions.
Let~$\varphi_0(H)$ and~$\varphi_1(H)$ be two $\Sigma^{0,C}_1$ formulas with one distinguished set parameter.
Given a condition~$c = (F, X)$, we would like to decide whether there is an extension~$d = (E, Y)$
forcing~$\varphi_0(G \cap A_0) \vee \varphi_1(G \cap A_1)$ or forcing~$\neg \varphi_1(G \cap A_0) \vee \neg \varphi_1(G \cap A_1)$.
By ``forcing'', we mean that the property holds for every set~$G$ satisfying the condition~$d$.

There exist two main arguments to decide $\Sigma^{0,C}_1$ formulas. Both seem
to be equi-expressive whenever the overall construction has no effectiveness restriction.
Indeed, all the existing weakness preservation proofs for generic weakness notions can
use the former and the latter argument interchangeably.

In both cases, we assume that~$\wkl$ admits $\Pcal$ preservation. 
In fact, this assumption is too strong and we will present in Chapter~\ref{chap:ramsey-type-konig-lemma}
the Ramsey-type weak K\"onig's lemma which is a weakening of $\wkl$ sufficient to carry both arguments.

\section{Seetapun-style forcing}

This argument has been introduced by Seetapun~\cite{Seetapun1995strength} to
prove that Ramsey's theorem for pairs admits cone avoidance. 
It has been recently generalized by Dzhafarov~\cite{DzhafarovStrong} who extracted the core of its combinatorics.

Fix a condition~$c = (F, X)$ and pick an $X \oplus C$-computable sequence of finite sets of maximal length
$E_0 < E_1 < \dots$ over~$X$, such that~$\varphi_0((F \cap A_0) \cup E_s)$ holds for each~$s$.
Suppose the sequence~$E_0 < E_1 < \dots$ finite. Then there is some stage~$s$ and a threshold~$t$ such that
for every finite set $E \subseteq X \cap (t, +\infty)$, $\varphi_0((F \cap A_0) \cup E)$ does not hold.
In this case, the condition~$d = (F, X \cap (t, +\infty))$ forces~$\varphi_0(G \cap A_0)$ not to hold
and therefore forces~$\neg \varphi_0(G \cap A_0) \vee \neg \varphi_1(G \cap A_1)$.
So suppose now that the sequence is infinite. 

If there is some~$s \in \omega$ such that $E_s \subseteq A_0$, 
then the condition~$d = (F \cup E_s, X \setminus [0, max(E_s)])$ forces~$\varphi_0(G \cap A_0)$ to hold
and therefore forces $\varphi_0(G \cap A_0) \vee \varphi_1(G \cap A_1)$.
From now on, we suppose that for every~$s \in \omega$, $E_s \cap A_1 \neq \emptyset$.
In particular, the set~$A_1$ is not hyperimmune relative to~$X \oplus C$.

Let~$T \subseteq X^{<\omega}$ be the $X \oplus C$-computable, finitely branching tree
such that~$\sigma \in T$ iff $\sigma(i) \in E_s$ for each~$i < |\sigma|$
and for every~$E \subseteq \{\sigma(i) : i < |\sigma|-1\}$, $\varphi_1((F \cap A_1) \cup E)$ does not hold.

If the tree $T$ is finite, then there is some leaf~$\sigma \in T$ such that~$\range(\sigma) \subseteq A_1$.
Let~$E \subseteq \range(\sigma)$ be such that~$\varphi_1((F \cap A_1) \cup E)$ holds.
The condition~$d = (F \cup E, X \setminus [0, max(E)])$ forces~$\varphi_0(G \cap A_1)$ to hold
and therefore forces $\varphi_0(G \cap A_0) \vee \varphi_1(G \cap A_1)$.
Last, suppose that the tree~$T$ is infinite.

By $\Pcal$ preservation for~$\wkl$, there is an infinite path~$P$ through~$T$
such that~$P \oplus C \not \in \Pcal$. The set~$Y = \range(P)$
is an infinite subset of~$X$ such that~$\varphi_1((F \cap A_1) \cup E)$ does not hold
for every finite set~$E \subseteq Y$.
The condition~$d = (F, Y)$ forces~$\varphi_1(G \cap A_1)$ not to hold, and therefore
forces~$\neg \varphi_0(G \cap A_0) \vee \neg \varphi_1(G \cap A_1)$.
\bigskip

Seetapun's argument has a trial and error flavor, in that we first try to 
find an extension making~$\varphi_0(G \cap A_0)$ hold. If we fail trying such an extension,
then we exploit this failure to build an extension making~$\varphi_1(G \cap A_1)$ hold.

\section{CJS-style forcing}

This argument has been introduced by Cholak, Jockusch and Slaman~\cite{Cholak2001strength}
to control the first jump of solutions to Ramsey's theorem for pairs. It has been reused
by Dzhafarov and Jockusch~\cite{Dzhafarov2009Ramseys} to prove strong cone avoidance for~$\rt^1_2$.
Fix a condition~$c = (F, X)$ and ask the following question:

\begin{quote}
Is it the case that for every 2-partition~$Z_0 \cup Z_1 = X$, there is 
some side~$i < 2$ and a finite set~$E \subseteq Z_i$ such that~$\varphi_i((F \cap A_i) \cup E)$ holds?
\end{quote}

Suppose the answer is yes. In particular, taking~$Z_0 = X \cap A_0$ and~$Z_1 = X \cap A_1$,
there is some~$i < 2$ and a finite set~$E \subseteq X \cap A_i$ such that~$\varphi_i((F \cap A_i) \cup E)$ holds.
The condition~$d = (F \cup E, X \setminus [0, max(E)])$ is an extension forcing~$\varphi_i(G \cap A_i)$ to hold,
and therefore forcing~$\varphi_0(G \cap A_0) \vee \varphi_1(G \cap A_1)$.

Suppose now that the answer is no. Then, the $\Pi^{0,X \oplus C}_1$ class~$\Ccal$
of all sets~$Z_0 \oplus Z_1$ such that~$Z_0 \cup Z_1 = X$ is a 2-partition
and for every side~$i < 2$ and every finite set~$E \subseteq Z_i$,
$\varphi_i((F \cap A_i) \cup E)$ does not hold is not empty.
Since~$\wkl$ admits $\Pcal$ preservation, there is a set~$Z_0 \oplus Z_1 \in \Ccal$
such that~$Z_0 \oplus Z_1 \oplus X \oplus C \in \Pcal$.
As~$Z_0 \cup Z_1 = X$, one of the~$Z$'s is infinite, say~$Z_i$.
The condition~$d = (F, Z_i)$ is an extension of~$c$
forcing~$\varphi_i(G \cap A_i)$ not to hold, and therefore forcing
$\neg \varphi_0(G \cap A_0) \vee \neg \varphi_1(G \cap A_1)$.
\bigskip

The intuition behind this argument is the following: Given a condition~$c = (F, X)$,
we would like to know whether there is some side~$i < 2$ and a finite set~$E \subseteq X \cap A_i$
such that $\varphi_i((F \cap A_i) \cup E)$ holds.
The~$\rt^1_2$-instance~$A_0 \cup A_1 = \omega$ is non-effective, and therefore non-accessible.
All we know is $F \cap A_0$ and~$F \cap A_1$ since it requires only a finite amount of information from the~$A$'s.

As we cannot access the~$\rt^1_2$-instance $A_0, A_1$, we over-approximate the question
and ask whether for \emph{every} $\rt^1_2$-instance~$Z_0, Z_1$,
there is some finite set~$E$ which is compatible with the~$Z$'s and such that~$\varphi_i((F \cap A_i) \cup E)$ holds.
If the answer is yes, then in particular, the answer is yes for the good $\rt^1_2$-instance~$A_0, A1$.
However, if the answer is no, then the $\rt^1_2$-instance witnessing the failure may not be $A_0, A_1$
and therefore at first sight, we cannot deduce that $\varphi_i(G \cap A_i)$ will not hold for both~$i < 2$.

This is where we use an essential feature of Ramsey-type statements: the ability
to make a set a solution to multiple instances at the same time. Indeed, if we let~$B_0, B_1$
be the $\rt^1_2$-instance witnessing the negative answer, making the solutions~$G \cap A_0$
and~$G \cap A_1$ be also solutions to~$B_0$, $B_1$ from now on ensures
that~$\varphi_i(G \cap A_i)$ will not hold for some~$i < 2$.

One may wonder why we are able to access the witness of failure~$B_0, B_1$
whereas is can be \emph{any} $\rt^1_2$-instance. In fact, we would be happy
with any $\rt^1_2$-instance witnessing this failure, and this class admits 
a simple description since the formulas are~$\Sigma^{0,C}_1$.
Among this class, at least one witness instance has a simple description,
and can be accessed with the help of any PA degree.

Last, notice that by a compactness argument, the main question is~$\Sigma^{0,X \oplus C}_1$
and therefore can be decided in the jump of~$X \oplus C$.
In the case of a positive answer, the extension can be found $A_0 \oplus A_1 \oplus X \oplus C$-effectively,
and therefore in the jump of~$X \oplus C$ if~$A_0 \oplus A_1$ is $\Delta^{0,C}_2$.
In the case of a negative answer, we need to choose which one over~$Z_0$ and~$Z_1$ is infinite,
which requires the computational power of a PA degree relative to the jump of~$X \oplus C$.

\chapter{The colors in Ramsey's theorem}

The strength of Ramsey's theorem is known to remain the same when
changing the number of colors in the setting of reverse mathematics.
Indeed, given some coloring~$f : [\omega]^n \to k^2$, we can define another coloring
$g : [\omega]^n \to k$ by merging colors together by blocks of size~$k$. After one application of~$\rt^n_k$
to the coloring~$g$, we obtain an infinite set~$H$ over which~$f$ uses at most~$k$ different colors.
Another application of~$\rt^n_k$ gives an infinite $f$-homogeneous set.
This standard proof of~$\rca \vdash \rt^n_k \imp \rt^n_{k^2}$ involves two applications
of~$\rt^n_k$. In this chapter, we use computable reducibility to show that multiple applications are really necessary
to reduce~$\rt^n_k$ to~$\rt^n_\ell$ whenever~$k > \ell$ and~$n \geq 2$.

More generally, Hirschfeldt and Jockusch~\cite{Hirschfeldtnotions} asked how many applications of~$\rt^n_\ell$ are needed
in a proof that~$\rca \vdash \rt^n_\ell \imp \rt^n_k$ whenever~$k > \ell$ and~$n \geq 2$.
They introduced a refinement~$\Qsf \leq_\omega^n \Psf$ of computable entailment 
in which at most~$n$ applications of~$\Psf$ are allowed. This notion can be easily formulated within
the game-theoretic framework~\cite{Hirschfeldtnotions} (see section~\ref{sect:colors-ramsey-reduction-games}). 

Given two integers~$u, \ell \geq 1$, we let~$\pi(u, \ell)$ denote the unique~$a \geq 1$ such that
$u = a \cdot \ell - b$ for some~$b \in [0, \ell)$. Informally, $\pi(u, \ell)$ is the minimal number of pigeons
we can ensure in at least one pigeonhole, given~$u$ pigeons and~$\ell$ pigeonholes.
Given~$k \geq 1$ and~$\ell \geq 2$, we define~$m_{k, \ell}$ inductively as follows.
First, $m_{1, \ell} = 0$. Assuming that~$m_{s, \ell}$ is defined for every~$s < k$,
let~$m_{k, \ell} = 1+m_{\pi(k, \ell), \ell} \mbox{ whenever } k \geq 2$. Note that~$m_{\pi(k, \ell)}$
is already defined since $\pi(k, \ell) < k$ whenever~$k, \ell \geq 2$.
In this chapter, we answer the question of Hirschfeldt and Jockusch by proving the following theorem.

\begin{theorem}\label{thm:colors-ramsey-general} For every~$k > \ell \geq 2$
\begin{multicols}{2}
\begin{itemize}
	\item[(i)] $\rt^n_k \leq_\omega^2 \rt^n_\ell$ for every~$n \geq 3$
	\item[(ii)] $\srt^n_k \not \leq_c \rt^n_\ell$ for every~$n \geq 2$
	\item[(iii)] $\rt^2_k \leq_\omega^{m_{k, \ell}} \rt^2_\ell$
	\item[(iv)] $\srt^2_k \not \leq_\omega^{m_{k, \ell}-1} \rt^2_\ell$
\end{itemize}
\end{multicols}
\end{theorem}

We first prove items (i) and~(iii). The negative results will be proven
at the end of this chapter.


\begin{proof}[Proof of Theorem~\ref{thm:colors-ramsey-general} item (i)]
Let~$f : [\omega]^n \to k$ be an instance of~$\rt^n_k$.
By Jockusch~\cite{Jockusch1972Ramseys}, there is an $f$-computable $\rt^n_\ell$-instance~$g : [\omega]^n \to \ell$
such that every infinite $g$-homogeneous set $f$-computes the Turing jump of~$f$.
Let~$H_g$ be an infinite $g$-homogeneous set.
By Hirschfeldt and Jockusch~\cite{Hirschfeldtnotions}, there is an $f \oplus H_g$-computable 
$\rt^n_\ell$-instance $h : [\Nb]^n \to \ell$ such that every infinite $h$-homogeneous set~$H_h$
$f \oplus H_g$-computes a set~$P$ of PA degree relative to~$(f \oplus H_g)^{(n-2)}$ and
therefore relative to~$f^{(n-1)}$. By Lemma~A.2 in~\cite{Hirschfeldtnotions},
$P$ computes an infinite $f$-homogeneous set, and therefore
so does~$f \oplus H_g \oplus H_h$.
\end{proof}

\begin{proof}[Proof of Theorem~\ref{thm:colors-ramsey-general} item (iii)]
We prove that~$\rt^2_k \leq_\omega^{m_{k,\ell}} \rt^2_\ell$
by induction over~$k \geq 1$.
In the case~$k = 1$, $\rt^2_k$ is computably true and therefore we are done.
Fix~$k \geq 2$ and suppose it holds for all~$k' < k$.
Let~$f : [\omega]^2 \to k$ be a coloring.
Define the coloring~$g : [\omega]^2 \to \ell$ by~$g(x, y) = f(x, y) \mod \ell$.
For any infinite $g$-homogeneous set~$H_g$, the coloring~$f$ uses
at most~$\pi(k,\ell)$ colors over~$H_g$. Therefore, there is a $f \oplus H_g$-computable coloring~$h : [\omega]^2 \to \pi(k, \ell)$
such that every infinite $h$-homogeneous set $f \oplus H_g$-computes an infinite $f$-homogeneous set.
By induction hypothesis, since $\pi(k, \ell)$ is smaller than~$k$ whenever~$k, \ell \geq 2$,
$\rt^2_{\pi(k, \ell)} \leq_\omega^{m_{\pi(k, \ell), \ell}} \rt^2_\ell$.
By definition of~$m_{k, \ell}$, $m_{k, \ell} = 1 + m_{\pi(k, \ell), \ell}$ so
$\rt^2_{\pi(k, \ell)} \leq_\omega^{m_{k, \ell}-1} \rt^2_\ell$.
Therefore $\rt^2_k \leq_\omega^{m_{k,\ell}} \rt^2_\ell$.
\end{proof}

The remainder of this chapter is devoted to the proof of items (ii) and~(iv)
of Theorem~\ref{thm:colors-ramsey-general}.

\section{Partition of the integers and strong computable reduciblity}

We start our analysis with partitions of integers.
Of course, every computable partition has an infinite computable homogeneous set,
so we need to consider non-effective partitions and strong computable reducibility.
The study of~$\rt^1_k$ over strong reducibility has close connections with
cohesiveness.
Dzhafarov~\cite{Dzhafarov2014Cohesive} proved that~$\coh \not \leq_{sc} \mathsf{D}^2_{<\infty}$ by iterating the following theorem.

\begin{theorem}[Dzhafarov~\cite{Dzhafarov2014Cohesive}]\label{thm:colors-ramsey-dzhafarov-partitions}
For every~$k \geq 2$ and $\ell < 2^k$, there is a finite sequence~$R_0, \dots, R_{k-1}$
such that for all partitions~$A_0 \cup \dots \cup A_{\ell-1} = \omega$ hyperarithmetical in~$\vec{R}$,
there is an infinite subset of some~$A_j$ that computes no~$\vec{R}$-cohesive set.
\end{theorem}

Hirschfeldt and Jockusch noticed in~\cite{Hirschfeldtnotions} that the proof of Theorem~\ref{thm:colors-ramsey-dzhafarov-partitions}
can be slightly modified to obtain a proof that~$\rt^1_k \not \leq_{sc} \rt^1_\ell$ whenever~$k > \ell \geq 2$.
Mont\'alban asked whether the hyperarithmetic effectiveness restriction
can be removed from Dzhafarov's theorem. 
We give a positive answer, which has been proved independently by Hirschfeldt \& Jockusch~\cite{Hirschfeldtnotions}.
Moreover, we show that~$\vec{R}$ can be chosen to be low. 
More precisely, we will prove in this section the following theorem, from which
we deduce several corollaries about cohesiveness and~$\rt^1_k$.

\begin{theorem}\label{thm:colors-ramsey-rt1-hyperimmunity}
Fix some~$k \geq 1$ and $\ell \geq 2$, some set~$I$ and $k$ $I$-hyperimmune sets~$B_0, \dots, B_{k-1}$.
For every~$\ell$-partition~$A_0 \cup \dots \cup A_{\ell-1} = \omega$,
there exists an infinite subset~$H$ of some~$A_i$ such
that $\pi(k, \ell)$ sets among the~$B$'s are $I \oplus H$-hyperimmune.
\end{theorem}

We will postpone the proof of Theorem~\ref{thm:colors-ramsey-rt1-hyperimmunity}
until after Corollary~\ref{cor:colors-ramsey-srt2-not-reduc-colors}.
Using the existence of a low
$k$-partition~$B_0 \cup \dots \cup B_{k-1} = \omega$ such that
$\overline{B_j}$ is hyperimmune for every~$j < k$,
we deduce the following corollary.

\begin{corollary}\label{cor:colors-ramsey-d2-general-colors}
For every~$k > \ell \geq 2$, there is a low $k$-partition~
$B_0 \cup \dots \cup B_{k-1} = \omega$
such that for all $\ell$-partitions~$A_0  \cup \dots \cup A_{\ell-1} = \omega$,
there is an infinite subset~$H$ of some~$A_i$ and a pair~$j_0 < j_1 < k$
such that every infinite~$H$-computable set intersects both~$B_{j_0}$ and~$B_{j_1}$.
\end{corollary}
\begin{proof}
Fix some~$k > \ell \geq 2$ and a low $k$-partition $B_0 \cup \dots \cup B_{k-1} = \omega$
such that $\overline{B_j}$ is hyperimmune for every~$j < k$.
Since~$k > \ell \geq 2$, $\pi(k, \ell) \geq 2$.
Therefore, by Theorem~\ref{thm:colors-ramsey-rt1-hyperimmunity}, for every $\ell$-partition
$A_0  \cup \dots \cup A_{\ell-1} = \omega$, there is an infinite subset~$H$ of some~$A_i$ and a pair~$j_0 < j_1 < k$
such that~$\overline{B_{j_0}}$ and~$\overline{B_{j_1}}$ are~$H$-hyperimmune.
In particular, every infinite $H$-computable set intersects both~$B_{j_0}$ and~$B_{j_1}$.
\end{proof}

The positive answer to Mont\'alban's question is an immediate consequence
of the previous corollary.

\begin{corollary}\label{cor:colors-ramsey-dzhafarov-comb}
For every~$k \geq 2$ and $\ell < 2^k$, there is a finite sequence
of low sets~$R_0, \dots, R_{k-1}$
such that for all partitions~$A_0 \cup \dots \cup A_{\ell-1} = \omega$,
there is an infinite subset of some~$A_i$ that computes no~$\vec{R}$-cohesive set.
\end{corollary}
\begin{proof}
Given~$k \geq 2$ and~$\ell < 2^k$, fix the low $2^k$-partition $(B_\sigma : \sigma \in 2^k)$
whose existence is stated by Corollary~\ref{cor:colors-ramsey-d2-general-colors}.
For each~$i < k$, define~$R_i = \bigcup_{\sigma(i) = 1} B_\sigma$.
Note that by disjointness of the $B$'s, $\overline{R_i} = \bigcup_{\sigma(i) = 0} B_\sigma$.
By choice of the $B$'s, for all $\ell$-partitions~$A_0 \cup \dots \cup A_{\ell-1} = \omega$, there is
an infinite subset~$H$ of some~$A_j$ and a pair~$\sigma <_{lex} \tau \in 2^k$
such that every infinite~$H$-computable set intersects both~$B_\sigma$ and~$B_\tau$.
Let~$i < k$ be the least bit such that~$\sigma(i) \neq \tau(i)$.
As~$\sigma <_{lex} \tau$, $\sigma(i) = 0$ and~$\tau(i) = 1$. By definition of~$R_i$,
$B_{\tau} \subseteq R_i$ and~$B_{\sigma} \subseteq \overline{R_i}$. 
Therefore no infinite~$H$-computable set is homogeneous for~$R_i$. In particular 
no infinite~$H$-computable set is~$\vec{R}$-cohesive.
\end{proof}

The construction of the~$B$'s is done uniformly in~$k$.
We can therefore deduce the following corollary.

\begin{corollary}
There exists a sequence of low sets~$R_0, R_1, \dots$
such that every finite partition of~$\omega$ has an infinite subset in one of its parts
which does not compute an~$\vec{R}$-cohesive set.
\end{corollary}

The effectiveness of~$B$ in the statement of Corollary~\ref{cor:colors-ramsey-d2-general-colors} enables us
to deduce computable non-reducibility results about stable Ramsey's theorem for pairs,
thanks to the computable equivalence between~$\srt^2_\ell$ and the statement~$\mathsf{D}^2_\ell$.

\begin{corollary}\label{cor:colors-ramsey-srt2-not-reduc-colors}
For every~$k > \ell \geq 2$, $\srt^2_k \not \leq_c \srt^2_\ell$.
\end{corollary}
\begin{proof}
Fix $k > \ell \geq 2$. By Corollary~\ref{cor:colors-ramsey-d2-general-colors},
there is a~$\Delta^0_2$ $k$-partition~
$B_0 \cup \dots \cup B_{k-1} = \omega$
such that for all $\ell$-partitions~$A_0, \dots, A_{\ell-1}$ of~$\omega$,
there is an infinite subset~$H$ of some~$A_i$ which does not compute an infinite
subset of any~$B_j$.
By Cholak et al.~\cite{Cholak2001strength}, for every stable computable function~$f : [\omega]^2 \to \ell$,
there exists a~$\Delta^0_2$~$\ell$-partition~$A_0 \cup \dots \cup A_{\ell-1} = \omega$ such that every
infinite subset of a part computes an infinite~$f$-homogeneous set. Therefore, for every such function~$f$,
there exists an infinite $f$-homogeneous set which does not compute an infinite subset of any~$B_j$.
By Shoenfield's limit lemma~\cite{Shoenfield1959degrees}, the~$\Delta^0_2$ approximation~$g : [\omega]^2 \to k$ of the~$k$-partition~
$B_0 \cup \dots \cup B_{k-1} = \omega$ is a stable computable function and every infinite~$g$-homogeneous set with color~$j$ is an infinite subset of~$B_j$.
\end{proof}

We now turn to the proof of Theorem~\ref{thm:colors-ramsey-rt1-hyperimmunity}.
We shall prove it by induction over~$\ell$, using a forcing construction whose
forcing conditions are Mathias conditions~$(F, X)$ where~$X$ is an infinite set 
such that the~$B$'s are $X \oplus I$-hyperimmune. The case where~$\ell = 1$ trivially holds
since $\pi(k, 1) = k$.

\subsection{Forcing limitlessness}

For every~$\ell$-partition~$A_0 \cup \dots \cup A_{\ell-1} = \omega$, 
we want to satisfy the following scheme of requirements to ensure that~$G \cap A_i$
is infinite for each~$i < \ell$.
\[
\Qcal_p : (\exists n_0, \dots, n_{\ell-1} > p)[n_0 \in G \cap A_0 \wedge \dots \wedge n_{\ell-1} \in G \cap A_{\ell-1}]
\]
Of course, all requirements may not be satisfiable
if some part~$A_i$ is finite. Usually, a forcing argument starts with the assumption that the instance 
is non-trivial, that is, does not admit a solution with the desired properties (cone avoiding, low, ...).
In order to force the solution to be infinite,
it suffices to ensure that the reservoirs satisfy the desired properties, and therefore 
cannot be a solution to a non-trivial instance.

In our case, we say that an~$\ell$-partition~$A_0 \cup \dots \cup A_{\ell-1}$ is~\emph{non-trivial} if
there is no infinite set~$H$ included in the complement of one of the~$A$'s and such that
the $B$'s are $H \oplus I$-hyperimmune. 
The following lemma states that we can focus on non-trivial partitions without loss of generality.

\begin{lemma}\label{lem:colors-ramsey-partition-reduction}
For every trivial $\ell$-partition $A_0 \cup \dots \cup A_{\ell-1}$, 
there is an infinite set~$H \subseteq A_i$ for some~$i < \ell$ such that $\pi(k, \ell)$ sets among the $B$'s are $H \oplus I$-hyperimmune.
\end{lemma}
\begin{proof}
Let~$G = \{n_0 < n_1 < \dots \}$ be an infinite subset of $\overline{A}_i$ for some~$i < \ell$
such that the~$B$'s are $G \oplus I$-hyperimmune.
Define the~$(\ell-1)$-partition~$(C_j : j \neq i)$ by setting~$C_j = \{ s \in \omega : n_s \in A_j \}$ for each~$j \neq i$.
By induction hypothesis, there exists an infinite set~$H_0 \subseteq C_j$ for some~$j \neq i$
such that $\pi(k, \ell - 1)$ sets among the $B$'s are $H_0 \oplus G \oplus I$-hyperimmune.
Note that $\pi(k, \ell-1) \leq \pi(k, \ell)$.
The set $H = \{ n_s : s \in H_0 \}$ is an $H_0 \oplus G$-computable subset of $A_j$
and $\pi(k, \ell)$ sets among the $B$'s are $H \oplus J$-hyperimmune.
\end{proof}

Notice that the proof of Lemma~\ref{lem:colors-ramsey-partition-reduction}
uses the induction hypothesis with a different context, namely, $G \oplus I$ instead of $I$.
This is where we needed to use the relativized version of the theorem in the proof.
A condition~$c = (F, X)$ \emph{forces $\Qcal_p$}
if there exists some~$n_0, \dots, n_{m-1} > p$ such that~$n_i \in F \cap A_j$
for each~$i < \ell$. Therefore, if~$G$ satisfies~$c$ and~$c$ forces~$\Qcal_p$,
then~$G$ satisfies the requirement~$\Qcal_p$.
We now prove that the set of conditions forcing~$\Qcal_p$ is dense for each~$p \in \omega$.
Thus, every sufficiently generic filter will induce an infinite solution.

\begin{lemma}\label{lem:colors-ramsey-coh-strong-reduc-force-Q}
For every condition~$c$ and every~$p \in \omega$, there is an extension forcing~$\Qcal_p$.
\end{lemma}
\begin{proof}
Fix some~$p \in \omega$. It is sufficient to show that given a condition~
$c = (F, X)$ and some~$i < \ell$,
there exists an extension~$d_0 = (E, Y)$ and some integer~$n_i > p$
such that~$n_i \in E \cap A_i$.
By iterating the process for each~$i < \ell$, we obtain the desired extension~$d$.
By definition of non-triviality, $A_i$ is co-immune in~$X$
and therefore $X \cap A_i$ is infinite.
Take any~$n_i \in X \cap A_i \cap (p, +\infty)$.
The condition $d_0 = (F \cup \{n_i\}, X \setminus [0, n_i])$ is the desired extension.
\end{proof}

\subsection{Forcing non-homogeneity}

The second scheme of requirements aims at ensuring that
for some~$i < \ell$,
at least~$\pi(k, \ell)$ sets among the~$B$'s are $(G \cap A_i) \oplus I$-hyperimmune.
The requirements are of the following form for each~$j < k$
and each tuple of indices~$\vec{e} = e_0, \dots, e_{\ell-1}$.
\[
\Rcal_{\vec{e}, j} : \hspace{10pt} \Rcal^{A_0, B_j}_{e_0} \vee \dots \vee \Rcal^{A_{\ell-1}, B_j}_{e_{\ell-1}}
\]
where $\Rcal_e^{A, B}$ is the statement ``$\Phi_e^{(G \cap A) \oplus I}$ does not dominate $p_B$''.

We claim that if all the requirements are satisfied, then~$(G \cap A_i)$ has the desired property
for some~$i < \ell$. Indeed, if for some fixed~$j < k$, all the requirements~$\Rcal_{\vec{e}, j}$
are satisfied, then by the usual pairing argument, there is some~$i < \ell$ such that $B_j$
is $(G \cap A_i) \oplus I$-hyperimmune. So if all the requirements are satisfied, then by the pigeonhole
principle, there is some~$i < \ell$ such that~$\pi(k, \ell)$ sets among the~$B$'s are $(G \cap A_i) \oplus I$-hyperimmune.

A condition~\emph{forces $\Rcal_{\vec{e}, j}$} if every set~$G$ satisfying this condition also satisfies the requirement~$\Rcal_{\vec{e}}$.
The following lemma is the core of the forcing argument.

\begin{lemma}\label{lem:colors-ramsey-coh-strong-reduc-force-R}
For every condition~$c = (F, X)$, every~$j < k$ and every tuple of Turing indices~$\vec{e}$, 
there exists an extension~$d = (E, Y)$ forcing~$\Phi_{e_i}^{(G \cap A_i) \oplus I}$ not to dominate~$p_{B_j}$
for some~$i < \ell$.
\end{lemma}
\begin{proof}
Let~$f$ be the partial~$X \oplus I$-computable function which on input~$x$, 
searches for a finite set of integers~$U$ such that for every $\ell$-partition
$Z_0 \cup \dots \cup Z_{\ell-1} = X$, there is some~$i < \ell$
and some set~$E \subseteq Z_i$ such that $\Phi_{e_i}^{((F \cap A_i) \cup E) \oplus I}(x) \downarrow \in U$.
If such a set~$U$ is found, then~$f(x) = max(U)+1$, otherwise~$f(x) \uparrow$. We have two cases.
\begin{itemize}
	\item Case 1: The function~$f$ is total. By~$X \oplus I$-hyperimmunity of~$B_j$,
	$f(x) \leq p_{B_j}(x)$ for some~$x$. Let~$U$ be the finite set witnessing~$f(x) \downarrow$.
	Letting~$Z_i = X \cap A_i$ for each~$i < \ell$, there is some~$i$ and some finite set $E \subseteq X \cap A_i$
	such that~$\Phi_{e_i}^{((F \cap A_i) \cup E) \oplus I}(x) \downarrow \in U$.
	The condition~$d = (F \cup E, X \setminus [0, max(E)])$ is an extension forcing
	$\Phi_{e_i}^{(G \cap A_i) \oplus I}(x) < f(x) \leq p_{B_j}(x)$.

	\item Case 2: There is some~$x$ such that~$f(x) \uparrow$.
	By compactness, the~$\Pi^{0,X \oplus I}_1$ class $\Ccal$ of sets~$Z_0 \oplus \dots \oplus Z_{\ell-1}$
	such that~$Z_0 \cup \dots \cup Z_{\ell-1} = X$ and for every~$i < \ell$
	and every set~$E \subseteq Z_i$, $\Phi_{e_i}^{((F \cap A_i) \cup E) \oplus I}(x) \uparrow$
	is non-empty. By the hyperimmune-free basis theorem~\cite{Jockusch197201},
	there is some~$\ell$-partition $Z_0 \oplus \dots \oplus Z_{\ell-1} \in \Ccal$
	such that all the~$B$'s are $Z_0 \oplus \dots \oplus Z_{\ell-1} \oplus X \oplus I$-hyperimmune.
	Let~$i < \ell$ be such that~$Z_i$ is infinite. The condition~$d = (F, Z_i)$
	is an extension of~$c$ forcing~$\Phi_{e_i}^{(G \cap A_i) \oplus I}(x) \uparrow$.
\end{itemize}
\end{proof}

\subsection{Construction}
We have all necessary ingredients to build an infinite set~$G$
such that each~$G \cap A_i$ is infinite, and such that $\pi(k, \ell)$ sets among the~$B$'s
are $(G \cap A_i) \oplus I$-hyperimmune for some~$i < \ell$. 
Thanks to Lemma~\ref{lem:colors-ramsey-coh-strong-reduc-force-Q} 
and Lemma~\ref{lem:colors-ramsey-coh-strong-reduc-force-R}, define an infinite descending
sequence of conditions~$(\varepsilon, \omega) \geq c_0 \geq \dots$ such that for each~$s \in \omega$,
\begin{itemize}
	\item[(a)] $c_s$ forces~$\Qcal_s$
	\item[(b)] $c_s$ forces~$\Rcal_{\vec{e},j}$ if~$s = \tuple{\vec{e},j}$
\end{itemize}
where~$c_s = (F_s, X_s)$.
Define the set~$G = \bigcup_s F_s$. By (a), $G \cap A_i$ is infinite 
for every~$i < \ell$, and by (b), each requirement~$\Rcal_{\vec{e},j}$
is satisfied. This finishes the proof of Theorem~\ref{thm:colors-ramsey-rt1-hyperimmunity}.

\section{Reducibility to Ramsey's theorem for pairs}\label{sect:colors-ramsey-rt-reducibility}

Dorais et al.~\cite{Dorais2016uniform} asked whether $\rt^n_k \not \leq_c \rt^n_\ell$ for every~$n \geq 2$ and $k > \ell \geq 2$.
Hirschfeldt \& Jockusch~\cite{Hirschfeldtnotions} 
and Rakotoniaina~\cite{Rakotoniaina2015Computational} proved that~$\srt^n_k$ is not Weihrauch reducible
to~$\rt^n_\ell$ whenever~$k > \ell$. We extend the result to computable reducibility.
In the first place, we shall focus on the case~$n = 2$.
For this, we will take advantage of the proof of~$\rt^2_\ell$ that applies the cohesiveness principle
to obtain a stable coloring~$f : [\omega]^2 \to \ell$. This coloring can itself be considered
as the~$\Delta^0_2$ approximation of a~$\emptyset'$-computable $\ell$-partition of~$\omega$,
and therefore as a non-effective instance of~$\rt^1_\ell$.
Any infinite subset of one of its parts computes an infinite set homogeneous for~$f$.

Lerman, Solomon and Towsner~\cite{Lerman2013Separating} introduced a framework
to separate Ramsey-type statements in reverse mathematics.
It turns out that their framework correspond to variants of preservation of hyperimmunity~\cite{Patey2015Iterative}.
We shall see that this framework provides a sufficiently fine-grained analysis of Ramsey's theorem
to answer Hirschfeldt and Jockusch's question.

\index{preservation!of hyperimmunity}
\begin{definition}[Preservation of hyperimmunity]
A $\Pi^1_2$ statement~$\Psf$ \emph{admits preservation of hyperimmunity} if for each set~$Z$, 
each sequence of $Z$-hyperimmune sets~$A_0, A_1, \dots$,
and each $\Psf$-instance $X \leq_T Z$,
there is a solution $Y$ to~$X$ such that the~$A$'s are $Y \oplus Z$-hyperimmune.
\end{definition}

By Lemma~\ref{lem:intro-reduc-preservation-separation},
if some statement~$\Psf$ admits preservation of hyperimmunity
and some other statement~$\Qsf$ does not, then $\rca \wedge \Psf \nvdash \Qsf$.
By the hyperimmune-free basis theorem~\cite{Jockusch197201},
weak K\"onig's lemma admits preservation of hyperimmunity.
Note that any hyperimmune-free degree preserves \emph{all} the hyperimmune
sets simultaneously, which is much stronger than preserving only
countably many of them. Some statements such as cohesiveness
have computable instances with no hyperimmune-free solution, while 
they admit preservation of hyperimmunity.

\begin{theorem}\label{thm:coh-hyperimmunity-preservation}
$\coh$ admits preservation of hyperimmunity.
\end{theorem}

The proof is done by the usual construction of a cohesive set
with Mathias forcing, combined with the following lemma.

\begin{lemma}\label{lem:coh-preservation-lemma}
For every set~$Z$, every $Z$-computable Mathias condition~$(F, X)$, 
every index~$e$ and every $Z$-hyperimmune set~$B$,
there exists an extension~$(E, Y)$ such that $X =^{*} Y$
and which forces~$\Phi_e^{G \oplus Z}$ not to dominate~$p_B$.
\end{lemma}
\begin{proof}
Let~$f$ be the partial $Z$-computable function which
on input~$x$, searches for some finite set~$E \subseteq X$ and some~$n$
such that~$\Phi_e^{(F \cup E) \oplus Z}(x) \downarrow = n$. If such an~$n$
is found, then $f(x) = n$, otherwise~$f(x) \uparrow$.
If there is some~$x$ such that~$f(x) \uparrow$,
then~$(F, X)$ already forces~$\Phi_e^{G \oplus Z}$ to be partial.
Otherwise, $f$ is a $Z$-computable function, and by $Z$-hyperimmunity of~$B$,
there is some~$x$ such that~$f(x) < p_B(x)$.
Let~$E \subseteq X$ be such that~$\Phi_e^{(F \cup E) \oplus Z}(x) \downarrow < p_B(x)$.
The condition~$d = (F \cup E, X \setminus [0, max(E)])$ extends~$c$ and forces
$\Phi_e^{G \oplus Z}(x) < p_B(x)$.
\end{proof}

Note that this theorem is optimal in the sense that every p-cohesive set is hyperimmune.
$\rt^2_\ell$ does not admit preservation of hyperimmunity for any~$\ell \geq 2$.
However, it is able to preserve the hyperimmunity of some sets among an initial sequence
of hyperimmune sets.

\index{preservation!of k hyperimmunities}
\begin{definition}[Preservation of $k$ hyperimmunities]
A $\Pi^1_2$ statement~$\Psf$ \emph{admits preservation of $p$ among $k$ hyperimmunities} 
if for each set~$Z$, each $Z$-hyperimmune sets~$A_0, \dots, \allowbreak A_{k-1}$,
and each $\Psf$-instance $X \leq_T Z$,
there is a solution $Y$ to~$X$ such that at least~$p$ sets among the~$A$'s are $Y \oplus Z$-hyperimmune.
\end{definition}

We say that $\Psf$ admits preservation of~$k$ hyperimmunities if it admits preservation
of~$k$ among $k$ hyperimmunities. Beware, preservation of~$\omega$ among~$\omega$ hyperimmunities
is strictly weaker than preservation of hyperimmunity.
For example, $\rt^2_2$ admits preservation of $\omega$ among $\omega$ hyperimmunities, 
but does not even admit preservation of~$2$ hyperimmunities.
Using Theorem~\ref{thm:colors-ramsey-rt1-hyperimmunity}, we can deduce the following theorem.

\begin{theorem}\label{thm:colors-ramsey-rt2-hyperimmunity}
For every~$k \geq 1$ and~$\ell \geq 2$, $\rt^2_\ell$ admits preservation of~$\pi(k, \ell)$ among~$k$ hyperimmunities.
\end{theorem}
\begin{proof}
Fix $k$ $Z$-hyperimmune sets $B_0, \dots, B_{k-1}$ for some set~$Z$.
Let~$f : [\omega]^2 \to \ell$ be a $Z$-computable coloring
and consider the sequence of sets~$R_0, R_1, \dots$ defined for each~$x \in \omega$ by
\[
R_x = \{ s : f(x, s) = 1 \}
\]
By Theorem~\ref{thm:coh-hyperimmunity-preservation}, there
is an infinite~$\vec{R}$-cohesive set~$C$ such that
the~$B$'s are hyperimmune relative to~$C \oplus Z$.
Let~$\tilde{f} : \omega \to \ell$ be defined by
$\tilde{f}(x) = \lim_{s \in C} f(x, s)$. By Theorem~\ref{thm:colors-ramsey-rt1-hyperimmunity}
relativized to~$C \oplus Z$, there is an infinite $\tilde{f}$-homogeneous set~$H$
such that $\pi(k, \ell)$ among the~$B$'s are $H \oplus C \oplus Z$-hyperimmune.
In particular, $H \oplus C \oplus Z$ computes an infinite $f$-homogeneous set. 
\end{proof}

Using again the existence of a low $k$-partition $B_0 \cup \dots \cup B_{k-1}$
such that~$\overline{B_j}$ is hyperimmune
for every~$j < k$, we deduce the following corollary.

\begin{corollary}\label{cor:colors-ramsey-rt2-general-colors}
For every~$k > \ell \geq 2$, there is a low $k$-partition~
$B_0 \cup \dots \cup B_{k-1} = \omega$
such that each computable coloring $f: [\omega]^2 \to \ell$
has an infinite $f$-homogeneous set~$H$ and a pair~$j_0 < j_1 < k$
such that every infinite~$H$-computable set intersects both~$B_{j_0}$ and~$B_{j_1}$.
\end{corollary}
\begin{proof}
Fix some~$k > \ell \geq 2$ and a low $k$-partition $B_0 \cup \dots \cup B_{k-1} = \omega$
such that $\overline{B_j}$ is hyperimmune for every~$j < k$.
Since~$k > \ell \geq 2$, $\pi(k, \ell) \geq 2$.
Therefore, by Theorem~\ref{thm:colors-ramsey-rt2-hyperimmunity}, for every $\rt^2_\ell$-instance
$f : [\omega]^2 \to \ell$, there is an infinite $f$-homogeneous set~$H$ and a pair~$j_0 < j_1 < k$
such that~$\overline{B_{j_0}}$ and~$\overline{B_{j_1}}$ are~$H$-hyperimmune.
In particular, every infinite $H$-computable set intersects both~$B_{j_0}$ and~$B_{j_1}$.
\end{proof}

Using Corollary~\ref{cor:colors-ramsey-rt2-general-colors} in a relativized form, we can extend
the result to colorings over arbitrary tuples.

\begin{theorem}
For every~$n \geq 2$, and every $k > \ell \geq 2$, there is a~$\Delta^0_n$ $k$-partition~
$B_0 \cup \dots \cup B_{k-1} = \omega$
such that each computable coloring $f: [\omega]^n \to \ell$
has an infinite $f$-homogeneous set~$H$ and a pair~$j_0 < j_1 < k$
such that every infinite~$H$-computable set intersects both~$B_{j_0}$ and~$B_{j_1}$.
\end{theorem}
\begin{proof}
This is proved in a relativized form by induction over~$n$.
The case~$n = 2$ is proved by relativizing Corollary~\ref{cor:colors-ramsey-rt2-general-colors}.
Now assume it holds for some~$n$ in order to prove it for~$n+1$.
Let~$P \gg \emptyset^{(n-1)}$ be such that~$P' \leq \emptyset^{(n)}$.
Such a set exists by the relativized low basis theorem~\cite{Jockusch197201}.
Applying the induction hypothesis to~$P$,
there is a~$\Delta^{0,P}_2$ (hence~$\Delta^0_{n+1}$) $k$-partition~$B_0 \cup \dots \cup B_{k-1} = \omega$
such that each $P$-computable coloring $f: [\omega]^n \to \ell$
has an infinite $f$-homogeneous set~$H$ and a pair~$j_0 < j_1 < k$
such that every infinite~$H \oplus P$-computable set intersects both~$B_{j_0}$ and~$B_{j_1}$.

Let~$f : [\omega]^{n+1} \to \ell$ be a computable coloring.
By Jockusch~\cite[Lemma 5.4]{Jockusch1972Ramseys}, there exists an infinite set~$C$ pre-homogeneous for~$f$
such that~$C \leq_T P$.
Let~$\tilde{f} : [C]^n \to \ell$ be the~$P$-computable coloring defined for each~$\sigma \in [C]^n$
by~$\tilde{f}(\sigma) = f(\sigma, a)$, where~$a \in A$, $a > max(\sigma)$.
Every $\tilde{f}$-homogeneous set is~$f$-homogeneous.
By definition of~$B_0 \cup \dots \cup B_{k-1} = \omega$, 
there exists an infinite~$\tilde{f}$-homogeneous (hence $f$-homogeneous) set~$H$
and a pair~$j_0 < j_1 < k$
such that every infinite~$H \oplus P$-computable set intersects both~$B_{j_0}$ and~$B_{j_1}$.
\end{proof}

Using the fact that $\mathsf{D}^n_k \leq_c \srt^n_k$ for every $n, k \geq 2$, we obtain
the following corollary strengthening the result 
of Hirschfeldt \& Jockusch~\cite{Hirschfeldtnotions} 
and Rakotoniaina~\cite{Rakotoniaina2015Computational}.
This proves item (ii) of Theorem~\ref{thm:colors-ramsey-general}.

\begin{corollary}\label{cor:ramsey-colors-not-computable-reducible}
For every~$n \geq 2$ and every~$k > \ell \geq 2$, $\srt^n_k \not \leq_c \rt^n_\ell$.
\end{corollary}

This answers in particular Question~7.1 of Dorais et al.~\cite{Dorais2016uniform}.
Note that two applications of~$\rt^n_2$ are sufficient to deduce~$\rt^n_k$ in the case~$n \geq 3$
by item (i) of Theorem~\ref{thm:colors-ramsey-general}.
The following corollary answers positively Question~5.5.3 of Mileti~\cite{Mileti2004Partition}.

\begin{corollary}\label{cor:two-ramsey-instances-cannot-be-solved-by-one}
There exists two stable computable functions~$f_1 : [\omega]^2 \to 2$
and~$f_2 : [\omega]^2 \to 2$ such that there is no computable~$g : [\omega]^2 \to 2$
with the property that every set~$H_g$ homogeneous for~$g$ computes
both a set~$H_{f_1}$ homogeneous for~$f_1$ and a set~$H_{f_2}$ homogeneous for~$f_2$.
\end{corollary}
\begin{proof}
By Corollary~\ref{cor:colors-ramsey-rt2-general-colors} with~$\ell = 2$ and~$k = 3$,
there exists a~$\Delta^0_2$ 3-partition~$B_0 \cup B_1 \cup B_2 = \omega$
such that each computable coloring $f: [\omega]^2 \to 2$
has an infinite $f$-homogeneous set~$H$ and a pair~$j_0 < j_1 < 3$
such that every infinite~$H$-computable set intersects both~$B_{j_0}$ and~$B_{j_1}$.
As in Corollary~\ref{cor:colors-ramsey-dzhafarov-comb}, we assume that the~$B$'s are disjoint.
By Shoenfield's limit lemma~\cite{Shoenfield1959degrees}, there exist
two stable computable colorings~$f_1$ and $f_2$ such that~$\lim_s f_1(\cdot, s) = B_0$ 
and $\lim_s f_2(\cdot, s) = B_1$.
If~$j_0 = 0$ (resp. $j_0 = 1$) then~$H$ does not compute an infinite set homogeneous for~$f_1$ (resp. $f_2$).
This completes the proof.
\end{proof}

\section{Ramsey's theorem and reduction games}\label{sect:colors-ramsey-reduction-games}

Hirschfeldt and Jockusch~\cite{Hirschfeldtnotions} introduced the following game-theoretic presentation
of computable entailment.

\index{reduction game}
\begin{definition}[Reduction game]
Given two $\Pi^1_2$ statements~$\Psf$ and~$\Qsf$, the \emph{reduction game} $G(\Qsf \to \Psf)$
is a two-player game that proceeds as follows.

On the first move, Player 1 plays a $\Psf$-instance $X_0$
and Player 2 either plays an $X_0$-computable solution to~$X_0$
and declares victory, in which case the game ends, or responds
with an $X_0$-computable $\Qsf$-instance~$Y_1$.

For $n > 1$, on the~$n$th move, Player 1 
plays a solution $X_{n-1}$ to the~$\Qsf$-instance $Y_{n-1}$.
Then Player 2 plays either an $(\bigoplus_{i < n} X_i)$-computable solution to~$X_0$
and declares victory, in which case the game ends again, or 
plays an $(\bigoplus_{i < n} X_i)$-computable $\Qsf$-instance $Y_n$.
Player 2 wins the game if it ever declares victory. Otherwise Player 1 wins.
\end{definition}

In particular, they showed that $\Psf \leq_\omega \Qsf$ if and only if
Player 2 has a winning strategy for the reduction game~$G(\Qsf \to \Psf)$.
They defined the reduction~$\Psf \leq_\omega^n \Qsf$ to hold if
Player 2 has a winning strategy that guarantees victory in $n+1$ or fewer moves.

We start with a simple lemma which states some properties of 
a finite sequence of integers $a_0, \dots, a_{m_{k, \ell}}$.
Intuitively, $a_n$ is the number of sets among $k$ hyperimmune sets
which can be preservered at the end of the $(n+1)$th move
during the reduction game $G(\rt^2_\ell \to \Psf)$.

\begin{lemma}\label{lem:colors-ramsey-finite-sequence-properties}
Fix some~$k \geq \ell \geq 2$ and let 
$a_0, \dots, a_{m_{k, \ell}}$ be the finite sequence defined by~$a_0 = k$ and~$a_{i+1} = \pi(a_i, \ell)$.
Then
\begin{itemize}
	\item[(i)] $m_{k, \ell} = i + m_{a_i, \ell}$ for each~$i \leq m_{k, \ell}$
	\item[(ii)] $a_i \geq 1$ for each~$i \leq m_{k, \ell}$
	\item[(iii)] $a_{i+1} < a_i$ for each~$i < m_{k, \ell}$
	\item[(iv)] $a_{m_{k, \ell}} = 1$.
\end{itemize}
\end{lemma}
\begin{proof}
Items (i-ii) are proven by induction over~$i$. The case~$i = 0$ is trivial,
and assuming it holds for~$i$, $m_{k, \ell} = i + m_{a_i, \ell} = i + 1 + m_{\pi(a_i, \ell), \ell} = i + 1 + m_{a_{i+1}, \ell}$.
and $a_{i+1} = \pi(a_i, \ell) \geq 1$ since~$a_i, \ell \geq 1$.
Item (iii) is proven as follows: Since~$i < m_{k, \ell}$, by item (i) $a_i \neq 1$.
$a_{i+1} = \pi(a_i, \ell) < a_i$ since~$a_i, \ell > 2$.
Item (iv) is proven from item (i) by taking~$i = m_{k, \ell}$.
Indeed, we obtain $m_{a_{m_{k, \ell}}, \ell} = 0$ and therefore~$a_{m_{k, \ell}} = 1$.
\end{proof}

\begin{proof}[Proof of Theorem~\ref{thm:colors-ramsey-general} item (i)]
Consider the reduction game $G(\rt^2_\ell \to \rt^2_k)$
and let~$a_0, \dots, a_{m_{k, \ell}}$ be the finite sequence defined in Lemma~\ref{lem:colors-ramsey-finite-sequence-properties}.

Let~$B_0 \cup \dots \cup B_{k-1} = \omega$ be a $\Delta^0_2$ $k$-partition
of~$\omega$ such that~$\overline{B_j}$ is hyperimmune for each~$j < k$.
By Shoenfield's limit lemma~\cite{Shoenfield1959degrees}, there is a stable computable function~$f_0 : [\omega]^2 \to k$
such that~$x \in B_j$ iff $\lim_s f_0(x, s) = j$ for each~$x \in \omega$.

We prove by induction over~$n \leq m_{k, \ell}$ that at the end of the $n$th move, if the game has not yet ended,
$a_{n-1}$ sets among the~$B$'s have a $(\bigoplus_{i < n} f_i)$-hyperimmune complement.
Player 1 first plays $f_0$. Since~$f_0$ is computable, $a_0 (=k)$ sets among the $B$'s 
have an $f_0$-hyperimmune complement.
After $n$ moves, $0 \leq n \leq m_{k, \ell}$, suppose that~$a_{n-1}$ sets among the~$B$'s 
have a $(\bigoplus_{i < n} f_i)$-hyperimmune complement. If at the end of the~$n$th move,
Player 2 does not declare victory, then he plays a $(\bigoplus_{i < n} f_i)$-computable $\rt^2_\ell$-instance
$g : [\omega]^2 \to \ell$. By Theorem~\ref{thm:colors-ramsey-rt2-hyperimmunity},
there is an infinite $H_g$-homogeneous set such that $\pi(a_{n-1}, \ell)$ sets among the $a_{n-1}$
remaining $B$'s with a $(\bigoplus_{i < n} f_i)$-hyperimmune complement 
have a $(\bigoplus_{i < n} f_i) \oplus H_g$-hyperimmune complement. Player 1 plays $H_g$.
This finishes the induction.

Suppose now that Player 2 wins in $m_{k, \ell}$ moves or fewer,
and let~$H$ be the $(\bigoplus_{i < n} f_i)$-computable solution to~$f$.
Then at least~$a_{n-1}$ among the~$B$'s have an $H$-hyperimmune complement. By items (iii) and (iv)
of Lemma~\ref{lem:colors-ramsey-finite-sequence-properties}, $a_{n-1} \geq 2$ and therefore
$H$ intersects at least two $B$'s, contradicting the fact that~$H$ is an infinite $f$-homogeneous set.
\end{proof}

\chapter{Degrees of unsolvability of cohesiveness}\label{chap:degrees-unsolvability-cohesiveness}

Cohesiveness plays a central role in reverse mathematics.
It appears naturally in the standard proof of Ramsey's theorem, as a preliminary
step to reduce an instance of Ramsey's theorem over $(n+1)$-tuples into a non-effective instance over $n$-tuples.
An important part of current research about Ramsey-type principles
in reverse mathematics consists in trying to understand whether
cohesiveness is a consequence of stable Ramsey's theorem for pairs,
or more generally whether it is a combinatorial consequence of the infinite pigeonhole principle~\cite{Cholak2001strength,Dzhafarov2014Cohesive,Dzhafarov2012Cohesive,Wang2013Omitting}.
Chong et al.~\cite{Chong2014metamathematics} recently showed using non-standard models that cohesiveness
is not a proof-theoretic consequence of the pigeonhole principle.
However it remains unknown whether stable Ramsey's theorem for pairs 
computationally implies cohesiveness.

Cohesiveness is a $\Pi^1_2$ statement whose instances are sequences of sets~$\vec{R}$
and whose solutions are~$\vec{R}$-cohesive sets. It is natural to wonder about the degrees
of unsolvability of the~$\vec{R}$-cohesive sets according to the sequence of sets~$\vec{R}$.
Mingzhong Cai asked whether whenever a uniformly computable sequence of sets~$R_0$, $R_1, \dots$
has no computable~$\vec{R}$-cohesive set, there exists a non-computable set which does not compute
one. On the opposite direction, one may wonder whether every unsolvable instance of~$\coh$ is maximally difficult.

In this chapter, we establish a pointwise correspondence
between the sets cohesive for a sequence and the sets whose jump computes a member of a~$\Pi^{0, \emptyset'}_1$ class.
Then, using the known interrelations between typical sets and~$\Pi^0_1$ classes,
we give precise genericity and randomness bounds above which no typical set helps computing a cohesive set.
We identify various layers of unsolvability and spot a class of instances sharing many
properties with the universal instance.
As the author~\cite{Patey2015Combinatorial} did about the pigeonhole principle and
weak K\"onig's lemma ($\wkl$), 
we show that some unsolvable instances of~$\coh$ are combinatorial consequences of the pigeonhole principle.

\section{Cohesiveness and K\"onig's lemma}\label{sect:strength-ramsey-cohesiveness-konig-lemma}

Jockusch and Stephan~\cite{Jockusch1993cohesive} studied the degrees of unsolvability of cohesiveness
and proved that~$\coh$ admits a universal instance whose solutions
are the p-cohesive sets. They characterized their degrees as those whose
jump is PA relative to~$\emptyset'$.
We refine their analysis by establishing 
an instance-wise correspondence between the degrees
of the sets cohesive for a sequence, and the degrees whose jump
computes a member of a non-empty $\Pi^{0, \emptyset'}_1$ class.

\index{$\Ccal(\vec{R})$}
\begin{definition}
Fix a uniformly computable sequence of sets~$R_0, R_1, \dots$
For every~$\sigma \in 2^{<\omega}$, let~$R_\sigma = \bigcap_{\sigma(i) = 1} R_i \cap \bigcap_{\sigma(i) = 0} \overline{R}_i$.
For example, $R_{0110} = \overline{R_0} \cap R_1 \cap R_2 \cap \overline{R_3}$.
By convention, $R_\varepsilon = \omega$.
Let $\Ccal(\vec{R})$ be the $\Pi^{0,\emptyset'}_1$ class of binary sequences~$P$
such that for every~$\sigma \prec P$, the set~$R_\sigma$ is infinite.
\end{definition}

Our first lemma shows that the degrees of~$\vec{R}$-cohesive
sets can be characterized by their jumps.
This lemma reveals in particular that low sets fail to solve unsolvable instances
of cohesiveness.

\begin{lemma}\label{lem:coh-tree}
For every uniformly computable sequence of sets~$R_0, R_1, \dots$,
a set computes an~$\vec{R}$-cohesive set if and only if
its jump computes a member of~$\Ccal(\vec{R})$.
\end{lemma}
\begin{proof}
Fix an~$\vec{R}$-cohesive set~$C$.
Let~$P = \bigcup \{ \sigma \in 2^{<\omega} : C \subseteq^{*} R_\sigma \}$.
The sequence~$P$ is infinite and~$C'$-computable as there exists exactly
one string~$\sigma$ of each length such that~$C \subseteq^{*} R_\sigma$.
In particular, for every~$\sigma \prec P$, $R_\sigma$ is infinite,
so~$P$ is a member of~$\Ccal(\vec{R})$.

Conversely, let~$X$ be a set whose jump computes a member~$P$ of~$\Ccal(\vec{R})$.
By Shoenfield's limit lemma~\cite{Shoenfield1959degrees}, there exists an~$X$-computable function~$f(\cdot, \cdot)$
such that for each~$x$, $\lim_s f(x, s) = P(x)$. Define an~$\vec{R}$-cohesive set~$C = \bigcup_s C_s$
$X$-computably by stages $C_0 = \emptyset \subsetneq C_1 \subsetneq \dots$ as follows. 
At stage~$s$, search for some string~$\sigma$ of length~$s$
and some integer~$n \in R_\sigma$ greater than~$s$ such that~$f(x, n) = \sigma(x)$ for each~$x < |\sigma|$.
We claim that such~$\sigma$ and~$n$ must exist, as there exists a threshold~$n_0$
such that for every~$n > n_0$, $f(x, n) = P(x)$ for each~$x < s$.
Let~$\sigma \prec P$ be of length~$s$. By definition of~$P$, $R_\sigma$ is infinite,
so there must exist some $n \in R_\sigma$ which is greater than~$n_0$ and~$s$.
Set~$C_{s+1} = C_s \cup \{n\}$ and go to the next stage.
We now check that~$C = \bigcup_s C_s$ is~$\vec{R}$-cohesive.
For every~$x \in \omega$, there exists a threshold~$n_1$ such that
for every~$n > n_1$, $f(x, n) = P(x)$. By construction, for every element~$n \in C \setminus C_{n_1}$,
$n \in R_\sigma$ for some string~$\sigma$ such that~$\sigma(x) = P(x)$. Therefore~$C \subseteq^{*} R_x$
or~$C \subseteq^{*} \overline{R_x}$.
\end{proof}

Jockusch and Stephan~\cite{Jockusch1993cohesive}
showed the existence of a uniformly computable sequence of sets~$R_0, R_1, \dots$
having no low~$\vec{R}$-cohesive set. We prove that it suffices
to consider any sequence~$\vec{R}$ with no computable~$\vec{R}$-cohesive set to obtain this property.

\begin{corollary}\label{cor:low-helps-not-coh}
A uniformly computable sequence of sets~$R_0, R_1, \dots$ has
a low~$\vec{R}$-cohesive set if and only if it has a computable~$\vec{R}$-cohesive set.
\end{corollary}
\begin{proof}
Let~$X$ be a low $\vec{R}$-cohesive set. By Lemma~\ref{lem:coh-tree},
the jump of~$X$ (hence~$\emptyset'$) computes a member of~$\Ccal(\vec{R})$
and therefore there exists a computable~$\vec{R}$-cohesive set.
\end{proof}

One may naturally wonder about the shape of the~$\Pi^{0,\emptyset'}_1$
classes~$\Ccal(\vec{R})$ for uniformly computable sequences $R_0, R_1, \dots$
The next lemma shows that~$\Ccal(\vec{R})$ can be any~$\Pi^{0,\emptyset'}_1$
class. Together with Lemma~\ref{lem:coh-tree}, it establishes an instance-wise correspondence
between cohesive sets and~$\Pi^{0,\emptyset'}_1$ classes.

\begin{lemma}\label{lem:tree-coh-inv}
For every non-empty~$\Pi^{0,\emptyset'}_1$ class~$\Dcal \subseteq 2^{\omega}$,
there exists a uniformly computable sequence of sets~$R_0, R_1, \dots$
such that~$\Ccal(\vec{R}) = \Dcal$.
\end{lemma}
\begin{proof}
By Shoenfield's limit lemma~\cite{Shoenfield1959degrees}, 
there exists a computable function~$g : 2^{<\omega} \times \omega \to 2$
whose limit exists and such that~$\Dcal$ is the collection of~$X$ such that 
for every~$\sigma \prec X$, $\lim_s g(\sigma, s) = 1$. We can furthermore assume
that whenever~$g(\sigma, s) = 1$, then for every~$\tau \prec \sigma$, $g(\tau, s) = 1$,
and that for every~$s \in \omega$, the set~$U_s = \{ \sigma \in 2^s : g(\sigma, s) = 1\}$ is non-empty.
We define a uniformly computable sequence of sets~$R_0, R_1, \dots$ such that~$\Ccal(\vec{R}) = \Dcal$
by stages as follows. 

As stage~$0$, $R_i = \emptyset$ for every~$i \in \omega$.
Suppose that we have already decided~$R_i \restr n_s$ for every~$i \in \omega$ and some~$n_s \in \omega$.
At stage~$s+1$, we will add elements to~$R_0, \dots, R_s$ so that for each string~$\sigma$ of length~$s+1$,
$R_\sigma \restr [n_s, n_s+p] \neq \emptyset$ if and only if $\sigma \in U_{s+1}$.
To do so, consider the set~$U_{s+1} = \{\sigma_0, \dots, \sigma_p\}$ defined above
and add $\{n_s + i : \sigma_i(j) = 1, i \leq p \}$ to $R_j$ for each~$j \leq s$.
Set~$n_{s+1} = n_s + p + 1$ and go to the next stage.

We claim that $R_\sigma$ is infinite if and only if~$\sigma \prec X$ for some~$X \in \Dcal$.
Assume that~$R_\sigma$ is infinite. By construction, there are infinitely many~$s$
such that~$R_\sigma \restr [n_s, n_s+p] \neq \emptyset$. So there are infinitely many stages~$s$ such that
$\tau \in U_s$ ($g(\tau, s) = 1$) for some~$\tau \succeq \sigma$. By assumption on~$g$,
there are infinitely many~$\tau \succeq \sigma$ such that $g(\tau, s) = 1$ for infinitely many~$s$.
Therefore, by compactness, there exists some~$X \in \Dcal$ such that~$\sigma \prec X$.
Conversely, if~$\sigma \prec X$ for some~$X \in \Dcal$, then there are infinitely many stages~$s$
such that~$\tau \in U_s$ for some~$\tau \succeq \sigma$. At each of these stages,
$R_\sigma \restr [n_s, n_s + p] \supseteq R_\tau \restr [n_s, n_s + p] \neq \emptyset$.
Therefore~$R_\sigma$ is infinite.
\end{proof}

Jockusch et al.\ proved in~\cite{Jockusch199101} that for every~$\Pi^{0,\emptyset'}_1$ class $\Ccal \subseteq 2^{\omega}$,
there exists a~$\Pi^0_1$ class~$\Dcal \subseteq \omega^\omega$ such that~$\deg(\Ccal) = \deg(\Dcal)$,
where~$\deg(\Ccal)$ is the class of degrees of members of~$\Ccal$.
For the reader who is familiar with Weihrauch degrees, what we actually prove here is that
König's lemma is the jump of the cohesiveness principle under Weihrauch reducibility.
Bienvenu [personal communication] suggested the use of Simpson's Embedding Lemma~\cite[Lemma 3.3]{Simpson2007extension}
to prove the reducibility of some unsolvable instances of cohesiveness to various statements.

\begin{lemma}[Bienvenu]\label{lem:sigma3-coh}
For every~$\Sigma^{0,\emptyset'}_3$ class~$\Ecal \subseteq \omega^\omega$ with no $\emptyset'$-computable member,
there exists a uniformly computable sequence of sets~$R_0, R_1, \dots$
with no computable~$\vec{R}$-cohesive set but such that every member of~$\Ecal$ computes an~$\vec{R}$-cohesive set.
\end{lemma}
\begin{proof}
By a relativization of~Lemma~3.3 in~\cite{Simpson2007extension}, there exists a~$\Pi^{0, \emptyset'}_1$ class~$\Dcal$
whose degrees (relative to~$\emptyset'$) are exactly $\deg(\Ecal) \cup PA[\emptyset']$. Therefore~$\Dcal$ has no~$\emptyset'$-computable member and every member of~$\Ecal$ $\emptyset'$-computes a member of~$\Dcal$.
By Lemma~\ref{lem:tree-coh-inv}, there exists a uniformly computable sequence of sets~$R_0, R_1, \dots$
such that~$\Ccal(\vec{R}) = \Dcal$. By Lemma~\ref{lem:coh-tree}, there exists no computable $\vec{R}$-cohesive set,
but every member of~$\Dcal$ (and in particular every member of~$\Ecal$) computes an~$\vec{R}$-cohesive set.
\end{proof}

\section{How genericity helps solving cohesiveness}

A natural first approach in the analysis of the strength of a principle consists
in looking in which way typical sets can help computing
a solution to an unsolvable instance.
The notion of typical set is usually understood in two different ways: using the genericity approach
and the randomness approach. Wang~\cite{Wang2013Omitting} answered the question
of Mingzhong Cai by investigating the solvability of cohesiveness by typical sets.
We now take advantage of the analysis of the previous section to deduce optimal bounds
on how much genericity is needed to avoid solving an unsolvable instance of~$\coh$.

Recall that a real is~\emph{$n$-generic} if it meets or avoids each~$\Sigma^0_n$ set of strings.
A real is~\emph{weakly $n$-generic} if it meets each~$\Sigma^0_n$ dense set of strings.

By Friedberg's jump inversion theorem~\cite{Friedberg1957criterion}, there exists a 1-generic which is of high degree,
and therefore computes a cohesive set for every uniformly computable sequence of sets.
Wang~\cite{Wang2013Omitting} proved that whenever a uniformly computable sequence of sets~$R_0, R_1, \dots$
has no computable~$\vec{R}$-cohesive sets, no weakly 3-generic computes an~$\vec{R}$-cohesive set.
He asked whether there exists a 2-generic computing an~$\vec{R}$-cohesive set.
We prove the optimality of Wang's bound by showing the existence of an unsolvable instance
of~$\coh$ which is solvable by a 2-generic real.

\begin{lemma}\label{lem:cohesiveness-unsolvable-but-2-generic}
There exists a 2-generic real~$G$ together with a
uniformly computable sequence of sets~$R_0, R_1, \dots$
with no computable~$\vec{R}$-cohesive set such that
$G$ computes an~$\vec{R}$-cohesive set.
\end{lemma}
\begin{proof}
Fix any~$\Delta^0_3$ 2-generic real~$G$ and consider the singleton~$\Ecal = \{G\}$.
As no 2-generic is~$\Delta^0_2$, the class $\Ecal$ has no~$\emptyset'$-computable member.
By Lemma~\ref{lem:sigma3-coh}, there exists a uniformly computable sequence of sets
$R_0, R_1, \dots$ with no computable~$\vec{R}$-cohesive set, 
such that~$G$ computes an~$\vec{R}$-cohesive set.
\end{proof}

However, if we slightly increase the unsolvability of the sequence of sets,
no 2-generic real helps computing a set cohesive for the sequence.
A~\emph{1-enum} of a class~$\Ccal \subseteq 2^{<\omega}$
is a sequence of strings~$\sigma_0, \sigma_1, \dots$
such that~$|\sigma_s| = s$ and~$[\sigma_s] \cap \Ccal \neq \emptyset$ for each~$s \in \omega$.
The notion is extensively studied in chapter~\ref{chap-avoiding-enumerations-closed-sets}.

\begin{theorem}
For any uniformly computable sequence of sets~$R_0, R_1, \dots$
such that~$\Ccal(\vec{R})$ has no~$\emptyset'$-computable 1-enum,
no 2-generic real computes an~$\vec{R}$-cohesive set.
\end{theorem}
\begin{proof}
By Jockusch~\cite{JockuschJr1980Degrees}, every $n$-generic set is~GL${}_n$
and in particular, every 2-generic is~GL${}_1$. Therefore,
by Lemma~\ref{lem:coh-tree},
a 2-generic set~$G$ computes an~$\vec{R}$-cohesive set
if and only if there exists some functional~$\Gamma$
such that~$\Gamma^{G \oplus \emptyset'}$ is a member of~$\Ccal(\vec{R})$.
Fix a functional~$\Gamma$ such that~$\Gamma^{G \oplus \emptyset'}$ is total for some
2-generic set~$G$, and define the following $\Sigma^{0, \emptyset'}_1$ set:
\[
W_{bad} = \{ \sigma \in 2^{<\omega} : [\Gamma^{\sigma \oplus \emptyset'}] \cap \Ccal(\vec{R}) = \emptyset \}
\]

We claim that $G$ meets~$W_{bad}$.
Suppose for the sake of contradiction that $G$ avoids $W_{bad}$. By 2-genericity of~$G$, there exists a string~$\sigma \prec G$
with no extension in~$W_{bad}$. We show that there exists a~$\emptyset'$-effective
procedure which computes a 1-enum of~$\Ccal(\vec{R})$, contradicting our hypothesis.

On input~$n$, $\emptyset'$-effectively search for a~$\tau_n \succeq \sigma$
such that~$\Gamma^{\tau_n \oplus \emptyset'} \restr n$ is defined.
Such~$\tau_n$ exists as~$\sigma \prec G$ and~$\Gamma^{G \oplus \emptyset'}$ is total.
As~$\tau_n \not \in W_{bad}$, $[\Gamma^{\tau_n \oplus \emptyset'}] \cap \Ccal(\vec{R}) \neq \emptyset$
and therefore~$(\tau_n : n \in \omega)$ is a $\emptyset'$-computable 1-enum of~$\Ccal(\vec{R})$.
\end{proof}

Note that if we assume that~$G$ is weakly 3-generic and therefore avoids the set
$W_{bad} \cup W_{partial}$ where
\[
W_{partial} = \{ \sigma \in 2^{<\omega} : (\forall \tau \succeq \sigma) |\Gamma^{\tau \oplus \emptyset'}| < |\sigma| \}
\]
then we can furthermore impose that~$\tau_{n+1} \succeq \tau_n$ and~$\emptyset'$-compute a member of~$\Ccal(\vec{R})$.
This suffices to reprove that no weakly 3-generic helps solving an unsolvable intance of~$\coh$.

We now prove a theorem inspired by the proof
of domination closure of p-cohesive degrees by Jockusch and Stephan~\cite{Jockusch1993cohesive}.

\begin{theorem}
For any uniformly computable sequence of sets~$R_0, R_1, \dots$
such that~$\Ccal(\vec{R})$ has no~$\emptyset'$-computable 1-enum,
every~$\vec{R}$-cohesive set is of hyperimmune degree.
\end{theorem}
\begin{proof}
Suppose for the sake of contradiction that there exists some~$\vec{R}$-cohesive set~$C = \{a_0 <  a_1 < \dots \}$
and a computable set~$B = \{b_0 < b_1 < \dots\}$ such that~$(\forall i)(a_i < b_i)$.
For each~$n \in \omega$, let $B_n = \{n, n+1, \dots, b_n\}$. Note that~$a_n \in B_n$ for every~$n$,
and therefore for every length~$s$, there exists a string~$\sigma_s$ of length~$s$
such that~$(\forall^{\infty} n) R_{\sigma_s} \cap B_n \neq \emptyset$.
Let~$\sigma_0, \sigma_1, \dots$ be the~$\emptyset'$-computable sequence of such strings.
We claim that this sequence is a 1-enum of~$\Ccal(\vec{R})$, therefore contradicting our hypothesis.
Indeed, as $(\forall^{\infty} n) R_{\sigma_s} \cap B_n \neq \emptyset$, the set~$R_{\sigma_s}$
is infinite and therefore~$\Ccal(\vec{R}) \cap [\sigma_s] \neq \emptyset$.
\end{proof}

Of course, there exists some uniformly computable sequence of sets~$R_0, R_1, \dots$
with no computable $\vec{R}$-cohesive set but with an~$\vec{R}$-cohesive
set of hyperimmune-free degree. Simply apply Lemma~\ref{lem:sigma3-coh} with~$\Ecal = \{X\}$
where~$X$ is a~$\Delta^0_3$ set of hyperimmune-free degree. Such a set is known
to exist by Miller and Martin~\cite{Miller1968degrees}. The class $\Ecal$ has no $\emptyset'$-computable
member as every $\Delta^0_2$ set is hyperimmune.

\section{How randomness helps solving cohesiveness}

We now explore the interrelations between cohesiveness
and the measure-theoretic paradigm of typicality, namely, algorithmic randomness.
We already defined in section~\ref{subsect:intro-comp-algorithmic-randomness} the notion of Martin-L\"of random real.
The notion naturally relativizes by considering~$\Sigma^0_n$ tests.

\index{n-random@$n$-random}
\index{Martin-L\"of!$\Sigma^0_n$ test}
\index{weakly n-random@weakly $n$-random}
\begin{definition}[Randomness]
A~\emph{$\Sigma^0_n$ (Martin-Löf) test} is a sequence~$U_0, U_1, \dots$ of uniformly~$\Sigma^0_n$
classes such that~$\mu(U_i) \leq 2^{-i}$ for every~$i \in \omega$.
A real~$Z$ is \emph{$n$-random} if for every~$\Sigma^0_n$ test $U_0, U_1, \dots$,
$Z \not \in \bigcap_i U_i$.
A real~$Z$ is \emph{weakly $n$-random} if it is in every~$\Sigma^0_n$ class of measure~1.
\end{definition}

Wang~\cite{Wang2013Omitting} proved that whenever a uniformly computable sequence of sets~$R_0, R_1, \dots$
has no computable~$\vec{R}$-cohesive sets, there exists a Martin-Löf random real computing no~$\vec{R}$-cohesive set.
Thanks to Corollary~\ref{cor:low-helps-not-coh}, 
we know that it suffices to take any low Martin-Löf random real to obtain this property.
Wang asked whether we can always ensure the existence of a 3-random real computing
an~$\vec{R}$-cohesive set whenever the instance is unsolvable. 
The next two lemmas answer this question by proving that it depends on the considered
sequence of sets~$\vec{R}$.

\begin{lemma}
There exists a uniformly computable
sequence of sets~$R_0, R_1, \dots$ with no computable~$\vec{R}$-cohesive set,
but such that every 2-random real computes an~$\vec{R}$-cohesive set.
\end{lemma}
\begin{proof}
Let~$\Dcal$ be a $\Pi^{0,\emptyset'}_1$ class of positive measure
with no~$\emptyset'$-computable member.
By Lemma~\ref{lem:tree-coh-inv},
there exists a uniformly computable sequence of sets~$R_0, R_1, \dots$
such that~$\Ccal(\vec{R}) = \Dcal$.
By Kautz~\cite{Kautz1991Degrees,Kautz1998improved}, every 2-random real is,
up to prefix, a member of~$\Ccal(\vec{R})$.
Therefore, by Lemma~\ref{lem:coh-tree}, every 2-random real
computes an~$\vec{R}$-cohesive set.
\end{proof}

\begin{lemma}\label{lem:randomness-p-cohesive}
No weakly 3-random real computes a p-cohesive set.
\end{lemma}
\begin{proof}
Jockusch and Stephan~\cite{Jockusch1993cohesive} proved that degrees of p-cohesive sets
are those whose jump is PA relative to~$\emptyset'$.
By a relativization of~Stephan~\cite{Stephan2006Martin}, every 2-random real
whose jump is of PA degree relative to~$\emptyset'$ is high.
By Kautz~\cite{Kautz1991Degrees}, no weakly 3-random real is high.
\end{proof}

Avigad et al.~\cite{Avigad2012Algorithmic}
introduced the principle $\wwkls{n}$ stating that every~$\Delta^0_n$ tree of positive measure has a path.
In particular, $\wwkls{1}$ is~$\wwkl$.  
Thanks to Lemma~\ref{lem:randomness-p-cohesive},
for every~$n \in \omega$, one can apply the usual constructions to build an~$\omega$-model
of~$\wwkls{n}$ which does not contain any p-cohesive set and therefore is not a model of~$\coh$.
Pick any~$n$-random $Z$ which does not compute any p-cohesive set
and consider it as an infinite join~$Z_0 \oplus Z_1 \oplus \dots$.
By Van Lambalgen's theorem~\cite{VanLambalgen1990axiomatization}, the~$\omega$-structure whose second-order part
is the Turing ideal~$\{ X : (\exists i) X \leq_T Z_0 \oplus \dots \oplus Z_i \}$
is a model of~$\wwkls{n}$. Moreover it does not contain a p-cohesive set.

\section{How Ramsey-type theorems help solving cohesiveness}

In his paper separating Ramsey's theorem for pairs from
weak K\"onig's lemma, Liu~\cite{Liu2012RT22} proved that every (non-necessarily effective)
set~$A$ has an infinite subset of either it or its complement which is not of PA degree.
The absence of effectiveness conditions on~$A$ shows the combinatorial nature
of the weakness of the infinite pigeonhole principle.
On the other hand, the author~\cite{Patey2015Combinatorial} showed that this weakness depends on the choice
of the instance of~$\wkl$, by constructing a computable tree with no computable path together
with a~$\Delta^0_2$ set~$A$ such that every infinite subset of either~$A$ or~$\overline{A}$
computes a path trough the tree. We answer a similar question for cohesiveness
and study the weakness of the pigeonhole principle for typical partitions.

\begin{lemma}
There exists a~$\Delta^0_3$ (in fact low over~$\emptyset'$) set~$A$ 
and a uniformly computable sequence of sets~$R_0, R_1, \dots$
with no computable~$\vec{R}$-cohesive set, such that
every infinite subset of either~$A$ or~$\overline{A}$
computes an~$\vec{R}$-cohesive set.
\end{lemma}
\begin{proof}
Fix a set~$A$ which is low over~$\emptyset'$ and bi-immune relative to~$\emptyset'$.
The set of the infinite, increasing sequences which form an subset of either~$A$ or~$\overline{A}$ is $\Pi^{0,A}_1$,
hence~$\Pi^{0,\emptyset'}_2$ in the Baire space:
$$
\Ecal = \{ X \in \omega^\omega : (\forall s)[X(s) < X(s+1)] 
	\wedge [(\forall s)(X(s) \in A) \vee (\forall s)(X(s) \in \overline{A})] \}
$$
Moreover, $\Ecal$ has no~$\emptyset'$-computable member by bi-immunity relative to~$\emptyset'$ of~$A$.
Apply Lemma~\ref{lem:sigma3-coh} to complete the proof.
\end{proof}

We showed in a previous section the existence of a uniformly computable
sequence of sets~$R_0, R_1, \dots$ with no computable $\vec{R}$-cohesive set
such that every 2-random real computes an~$\vec{R}$-cohesive set.
The following lemma strengthens this result by constructing an unsolvable instance of~$\coh$
solvable by every infinite subset of any 2-random real.

\index{d.n.c.|see {diagonally non-computable}}
\index{diagonally non-computable}
\begin{definition}[Diagonal non-computability]
A function~$f : \omega \to \omega$ is \emph{diagonally non-computable} relative to~$X$
if for every~$e \in \omega$, $f(e) \neq \Phi^X_e(e)$.
\end{definition}

By Kjos-Hanssen~\cite{Kjos-Hanssen2009Infinite} and Greenberg \& Miller~\cite{Greenberg2009Lowness}, 
a set computes a function d.n.c.\ relative
to~$\emptyset^{(n-1)}$ if and only if it computes an infinite subset of an~$n$-random.

\begin{lemma}
There exists a uniformly computable sequence of sets~$R_0, R_1, \dots$
with no computable~$\vec{R}$-cohesive set, such that
every function d.n.c.\ relative to~$\emptyset'$ computes an~$\vec{R}$-cohesive set.
\end{lemma}
\begin{proof}
The class of functions which are d.n.c.\ relative to~$\emptyset'$ is $\Pi^{0, \emptyset'}_1$ in the Baire space:
$$
\Ecal = \left\{ f \in \omega^\omega : (\forall e)[\Phi^{\emptyset'}_e(e) \uparrow \vee f(e) \neq \Phi^{\emptyset'}_e(e)]\right\}
$$
Moreover, $\Ecal$ has no~$\emptyset'$-computable member.
Apply Lemma~\ref{lem:sigma3-coh} to complete the proof.
\end{proof}

In contrast with this lemma, 
if we require a bit more uncomputability in the~$\vec{R}$-cohesive sets of the sequence~$R_0, R_1, \dots$,
we can ensure the existence of a function d.n.c.\ relative to~$\emptyset'$
which does not compute an~$\vec{R}$-cohesive set.

\begin{theorem}\label{thm:coh-dnc-avoid}
Fix a uniformly computable sequence of sets~$R_0, R_1, \dots$
such that~$\Ccal(\vec{R})$ has no $\emptyset'$-computable 1-enum.
For every set~$X$, there exists a function~$f$ d.n.c.\ relative to~$X$
whose jump does not compute a 1-enum of~$\Ccal(\vec{R})$.
In particular, $f$ does not compute an~$\vec{R}$-cohesive set.
\end{theorem}

\index{forcing!bushy tree}
The proof of Theorem~\ref{thm:coh-dnc-avoid} is done by a bushy tree forcing argument.
See the survey from Khan and Miller~\cite{Khan2015Forcing} for terminology and definitions.
Fix a set~$X$. We will construct a GL${}_1$ function which is d.n.c.\ relative to~$X$.
Our forcing conditions are tuples~$(\sigma, B)$ where~$\sigma \in \omega^{<\omega}$
and $B \subseteq \omega^{<\omega}$ is an upward-closed set $k$-small above~$\sigma$ for some~$k \in \omega$.
A sequence~$f$ \emph{satisfies} a condition~$(\sigma, B)$ if $\sigma \prec f$ and $B$ is small above every initial segment of~$f$.
Our initial condition is~$(\varepsilon, B_{DNC}^X)$ where~
$$
B_{DNC}^X = \{ \sigma \in \omega^{<\omega} : (\exists e) \sigma(e) = \Phi^X_e(e) \}
$$
Therefore every infinite sequence~$f$ satisfying~$(\varepsilon, B_{DNC}^X)$ is d.n.c.\ relative to~$X$.
Thanks to the following lemma, we can prevent~$f \oplus \emptyset'$ from computing
a 1-enum of~$\Ccal(\vec{R})$. As the constructed function~$f$ is GL${}_1$, $f' \leq_T f \oplus \emptyset'$
does not compute a 1-enum of~$\Ccal(\vec{R})$.

\begin{lemma}
For every condition~$c = (\sigma, B)$ and every Turing functional~$\Gamma$,
there exists an extension~$d = (\tau, C)$ forcing~$\Gamma^{f \oplus \emptyset'}$ to be partial
or such that~$\Gamma^{\tau \oplus \emptyset'}$ is not a 1-enum of~$\Ccal(\vec{R})$.
\end{lemma}
\begin{proof}
Suppose that~$B$ is $k$-small above~$\sigma$.
For every~$n \in \omega$,
define the~$\Sigma^{0, \emptyset'}_1$ set~$D_n = \{ \tau \in \omega^{<\omega} : \Gamma^{\tau \oplus \emptyset'}(n) \downarrow \in 2^n\}$.
Make a~$\emptyset'$-effective search for an~$n \in \omega$ such that one of the following holds:
\begin{itemize}
	\item[(a)] $D_n$ is $k2^n$-small above~$\sigma$ for some~$n \in \omega$
	\item[(b)] $D_{n,\rho} = \{ \tau \in \omega^{<\omega} : \Gamma^{\tau \oplus \emptyset'}(n) \downarrow = \rho \}$
	is $k$-big above~$\sigma$ for some string~$\rho \in 2^n$ such that $[\rho] \cap \Ccal(\vec{R}) = \emptyset$.
\end{itemize}
Such an~$n$ exists, as otherwise, for every~$n \in \omega$, $D_n$ is $k2^n$-big above~$\sigma$.
By the smallness additivity property, $D_{n,\rho}$ is $k$-big above~$\sigma$ for some~$\rho \in 2^n$.
For every such string~$\rho$, $[\rho] \cap \Ccal(\vec{R}) \neq \emptyset$. Therefore we can $\emptyset'$-compute
a 1-enum of~$\Ccal(\vec{R})$ by searching on each input~$n$ for some~$\rho$ of length~$n$
such that~$D_{n,\rho}$ is $k$-big above~$\sigma$.

If we are in case~(a), take~$d = (\tau, C \cup D_n)$ as the desired extension.
The condition~$d$ forces~$\Gamma^{f \oplus \emptyset'}$ to be partial.
If we are in case~(b), by the concatenation property, there exists an extension~$\tau \in D_{n,\rho}$
such that~$B$ is still~$k$-small above~$\tau$. The condition $d = (\tau, B)$ is an extension
forcing~$\Gamma^{f \oplus \emptyset'}$ not to be a 1-enum of~$\Ccal(\vec{R})$
as~$\Gamma^{f \oplus \emptyset'}(n) = \Gamma^{\tau \oplus \emptyset'}(n) = \rho$ and~$[\rho] \cap \Ccal(\vec{R}) = \emptyset$.
\end{proof}

Looking at the proof of the previous lemma, we can~$\emptyset'$-decide in which case we are, and then use the knowledge of~$f$
to see which path has been chosen in the bushy tree. The construction therefore yields a~GL${}_1$ sequence.

\chapter{Thin set and free set theorems}\label{chap:thin-set-free-set-theorems}

Simpson~\cite[Theorem III.7.6]{Simpson2009Subsystems} proved
that Ramsey's theorem for~$n$-tuples is equivalent to the arithmetic comprehension
axiom over~$\rca$ whenever~$n \geq 3$. Therefore, the Ramsey hierarchy collapses
at level three. One may wonder about some natural weakenings of Ramsey's theorem over arbitrary tuples
which remain strictly weaker than~$\aca$. 

Ramsey's theorem asserts that every $k$-coloring over~$[\Nb]^n$ has an infinite
monochromatic set~$H$. A natural weakening consists of allowing more colors in the solution~$H$.
Let~$\rt^n_{k, \ell}$ be the statement ``Every $k$-coloring over~$[\Nb]^n$
has an infinite set~$H$ using at most~$\ell$ colors''.
Since~$\rt^n_{k, \ell}$ and~$\rt^n_{m, \ell}$ are equivalent over~$\rca$ whenever~$k, m \geq \ell$,
we shall fix~$\ell$ to be~$k-1$. The statement~$\rt^n_{k, k-1}$ is known as a \emph{thin set theorem}.

\index{thin set theorem}
\index{ts@$\ts^n_k$|see {thin set theorem}}
\begin{definition}[Thin set theorem]
Given a coloring~$f : [\Nb]^n \to k$ (resp.\ $f : [\Nb]^n \to \Nb$), an infinite set~$H$
is \emph{thin} for~$f$ if~$|f([H]^n)| \leq k-1$ (resp.\ $f([H]^n) \neq \Nb$).
For every~$n \geq 1$ and~$k \geq 2$, $\ts^n_k$ is the statement
``Every coloring $f : [\Nb]^n \to k$ has a thin set''
and $\ts^n$ is the statement ``Every coloring $f : [\Nb]^n \to \Nb$ has a thin set''.
\end{definition}

\index{sts@$\sts^n_k$|see {thin set theorem}}
As for~$\rt^2_2$, we shall denote by~$\sts^2_k$ the restriction of~$\ts^2_k$ to stable colorings.
$\ts$ is the statement~$(\forall n)\ts^n$.
The reverse mathematical analysis of the thin set theorem started with Friedman~\cite{FriedmanFom53free,Friedman2013Boolean}.
It has been studied by Cholak et al.~\cite{Cholak2001Free},
Wang~\cite{Wang2014Some} and the author~\cite{Patey2015Combinatorial,Patey2015Degrees} among others.

\section{The strength of the thin set theorem for pairs}

Various proofs involving Ramsey's theorem for pairs and its consequences
can be adapted to the thin set theorem. For instance, Downey, Hirschfeldt,
Lempp and Solomon~\cite{Downey200102} constructed a computable instance of~$\srt^2_2$
with no computable solution. A similar theorem can be obtained
for the thin set theorem simply by slightly changing their proof.

\begin{theorem}[$\rca + \ist$]\label{thm:sts-no-low}
There exists a computable stable coloring $f : [\Nb]^2 \to \Nb$ with no low infinite set thin for $f$.
\end{theorem}
\begin{proof}
This is a straightforward adaptation of the proof of~\cite{Downey200102}.
We assume that definitions and the procedure $P(m)$ has been defined like in the original proof.
Given a stable coloring $f : [\Nb]^2 \to \Nb$, define $A_i = \set{x \in \Nb : (\forall^{\infty} s)f(x, s) \neq i}$.

We need to satisfy the following requirements for all $\Sigma^0_2$ sets $U$, all partial computable
binary functions $\Psi$ and all $i \in \Nb$:
$$
\Rcal_{U, \Psi, i} : U \subseteq A_i \wedge U \in \Delta^0_2 \wedge U \mbox{ infinite } \wedge
\forall n(\lim_s \Psi(n, s) \mbox{ exists}) \imp U' \neq \lim_s \Psi(-,s)
$$

The strategy for satisfying a single requirement $\Rcal_{U, \Psi, i}$ is almost the same.
It begins by choosing an $e \in \Nb$. Whenever a number $x$ enters $U$, 
it enumerates the axiom $\tuple{e, \set{x}}$ for $U'$. Whenever it sees that $\Psi(e, s) \downarrow \neq 1$
for some new number $s$, it \emph{commits} every $x$ for which it has enumerated an axiom $\tuple{e, \set{x}}$
to be assigned color $i$, i.e. starts settings $f(x, t) = i$ for every $t \geq s$.

If $U$ is $\Delta^0_2$ and infinite, and $\lim_s \Psi(e, s)$ exists and is not equal to 1,
then eventually an axiom $\tuple{e, \set{x}}$ for some $x \in U$ is enumerated, in which case
$U'(e) = 1 \neq \lim_s \Psi(e,s)$. On the other hand, if $U \subseteq A_i$ and $\lim_s \Psi(e,s) = 1$
then for all axioms $\tuple{e, \set{x}}$ that are enumerated by our strategy, $x$ is eventually commited to 
be assigned color $i$, which implies that $x \not \in U$. Thus in this case, $U'(e) = 0 \neq \lim_s \Psi(e,s)$.

The global construction is exactly the same as in the original proof.
\end{proof}

We now prove that the thin set theorem implies the \emph{atomic model theorem}.

\index{atomic model theorem}
\index{amt@$\amt$|see {atomic model theorem}}
\index{atom}
\begin{definition}[Atomic model theorem]
A formula $\varphi(x_1, \dots, x_n)$ is an \emph{atom} of a theory $T$ if for each formula $\psi(x_1, \dots, x_n)$
we have $T \vdash \varphi \imp \psi$ or $T \vdash \varphi \imp \neg \psi$ but not both.
A theory $T$ is \emph{atomic} if, for every formula $\psi(x_1, \dots, x_n)$ consistent with $T$,
there is an atom $\varphi(x_1, \dots, x_n)$ extending it, i.e. one such that $T \vdash \varphi \imp \psi$.
A model $\Acal$ of $T$ is \emph{atomic} if every $n$-tuple from $\Acal$ satisfies an atom of $T$. 
$\amt$ is the statement ``Every complete atomic theory has an atomic model''.
\end{definition}

$\amt$ has been introduced as a principle by Hirschfeldt, Shore and Slaman in~\cite{Hirschfeldt2009atomic}.
It is in particular computably equivalent to the statement ``For every $\Delta^0_2$ function~$f$,
there is a function~$g$ not dominated by~$f$''~\cite{Hirschfeldt2009atomic,Conidis2008Classifying}. 
Therefore, the atomic model theorem can be seen as a statement between hyperimmunity 
and hyperimmunity relative to~$\emptyset'$.
Hirschfeldt et al.\ in~\cite{Hirschfeldt2009atomic} proved that~$\wkl$ and~$\amt$ are incomparable over~$\rca$,
and in particular that $\amt$ is a consequence of the stable ascending descending sequence principle 
(defined in chapter~\ref{chap:erdos-moser-theorem}).
We show that the proof can be adapted to the stable thin set theorem for pairs.

\begin{theorem}\label{thm:sts2-amt}
$\rca \vdash \sts^2 \imp \amt$ and $\amt \leq_W \sts^2$.
\end{theorem}
\begin{proof}
We prove that for any atomic theory $T$, there exists a $\Delta^{0,T}_1$ stable
coloring $f : [\Nb]^2 \to \Nb$ such that for any infinite $f$-thin set $H$, 
there is a~$\Delta^{0, H \oplus T}_1$ atomic model of $T$.
The proof is very similar to~\cite[Theorem~4.1]{Hirschfeldt2009atomic}.
We begin again with an atomic theory $T$ and consider the tree $\Scal$ of standard
Henkin constructions of models of $T$. We want to define a stable coloring $f : [\Nb]^2 \to \Nb$
such that any infinite $f$-thin set computes an infinite path $P$ through $\Scal$
that corresponds to an atomic model $\Acal$ of $T$.
Define as in~\cite[Theorem~4.1]{Hirschfeldt2009atomic} 
a monotonic computable procedure $\Phi$ which on tuple $\tuple{x_1, \dots, x_n}$
will return a tuple $\tuple{\sigma_1, \dots, \sigma_n}$ such that $\sigma_{i+1}$ is the
least node of $\Scal$ extending $\sigma_i$ such that we have found no witness that $\sigma_{i+1}$
is not an atom about $c_0, \dots, c_{x_i}$ after a standardized search of $x_{i+1}$ many steps.
$\sigma_1$ is the least node on $\Scal$ mentioning $c_0$ and such that we have not found a witness
that $\sigma_1$ is not an atom about $c_0$ after $x_1$ many steps.

The construction of the coloring $f$ will involve a movable marker procedure.
At each stage $s$, we will ensure to have defined $f$ on $\set{x : x \leq s}$.
For each color $i$, we can associate the set $C_i = \{x : (\forall^{\infty} s)f(x,s) = i\}$.
At stage $s$, we maintain a set $C_{i,s}$ with the intuition that $C_i = \lim_s C_{i,s}$.

For each~$e, i \in \Nb$, the requirement
$\Rcal_{e,i}$ ensures that for any sequence
$x_1, \dots, x_n, d_{e,i,t}$ in $\overline{C_i}$ that is increasing
in natural order, $\sigma_{n+1}$ includes an atom about $c_0, \dots, c_{x_n}$
where $d_{e,i,t}$ is the value of the marker~$d_{e,i}$ associated to~$\Rcal_{e,i}$
at stage~$t$, and $\Phi(x_1, \dots, x_n, d_{e,i,t}) = \tuple{\sigma_1, \dots, \sigma_{n+1}}$.

We say that the requirement $\Rcal_{e,i}$ \emph{needs attention at stage $s$}
if there exists a sequence $x_1, \dots, x_n, d_{e,i,s}$ of elements of $\overline{C_{i,s}}$ increasing in natural order,
such that $\Phi(x_1, \dots, x_n, d_{e,i,s}) = \tuple{\sigma_1, \dots, \sigma_{n+1}}$ and by stage $s$ we have seen
a witness that $\sigma_{n+1}$ does not supply an atom about $c_0, \dots, c_{x_n}$.

At stage $s$, suppose the highest priority requirement needing attention is $\Rcal_{e,i}$.
The strategy \emph{commits to $C_i$} each $x < s$ that are in greater or equal to $d_{e,i,s}$.
We let $d_{e, i,s+1} = s$. All $d_{u, j, s+1}$ become undefined for $\tuple{u,j} > \tuple{e,i}$.
If no requirement needs attention, we let $d_{u,j,s+1} = s$
for the least $\tuple{u,j}$ such that $d_{u,j,s}$ is undefined.
For each $x < s$, set $f(x, s) = i$ if $x$ is committed to be in $C_i$.
Otherwise set $f(x, s) = 0$. We then go to the next stage.

\begin{claim}
The resulting coloring is stable.
\end{claim}
\begin{proof*}
Take any $x \in \Nb$. If no requirement ever commits $x$ to be in some $D_i$
then $x$ is committed at stage $x+1$ to be in $C_0$ and this commitment is never injured,
so $(\forall^\infty s)f(x,s) = 0$. Otherwise by $\isig^0_1$ there is a requirement
of highest priority that commits $x$ to be in some $C_i$.
Say it is $\Rcal_{e,i}$ and it acts to commit $x$ at stage $s$.
This means that $d_{e,i,s} \leq x < s$. Then we set $d_{e,i,s+1} = s$ and never decrease this marker.
No requirement of higher priority will act after stage $s$ on $x$ by our choice of $\Rcal_{e,i}$
and the markers for strategies of lower priority will be initialized after stage $s$ to 
a value greater than $s$. So $x$ will stay for ever in $C_i$. Thus $(\forall^\infty s)f(x,s) = i$.
\end{proof*}

\begin{claim}
Each requirement $\Rcal_{e,i}$ acts finitely often and
$d_{e,i,s}$ will eventually remain fixed. Moreover,
if $d_{e,i,s}$ never changes after stage $t$, then, for any sequence
$x_1, \dots, x_n, d_{e,i,t}$ in $\overline{C_i}$ that is increasing
in natural order, $\sigma_{n+1}$ includes an atom about $c_0, \dots, c_{x_n}$
where $\Phi(x_1, \dots, x_n, d_{e,i,t}) = \tuple{\sigma_1, \dots, \sigma_{n+1}}$.
\end{claim}
\begin{proof*}
We prove it by $\Sigma^0_1$ induction.
Assume that $\Rcal_{e,i}$ acts at stage $s$ and no requirement of higher
priority ever acts again. We then set $d_{e,i,s+1} = s$ and 
now act again for $\Rcal_{e,i}$ only if we discover a new witness as described
in the definition of needing attention. As we never act for any requirement
of higher priority, at any stage $t > s$ the numbers between $d_{e,i,s}$ and $d_{e,i,t}$
will all be committed to~$C_i$. Then the sequences $x_1, \dots, x_n \leq d_{e,i,t}$ in $\overline{C_i}$,
increasing in natural order are sequences $x_1, \dots, x_n \leq d_{e,i,s}$ in $\overline{C_i}$.
Hence their set is bounded. By the same trick as in~\cite[Theorem~4.1]{Hirschfeldt2009atomic},
we can avoid the use of $\bst$ by constructing a single atom extending each $\sigma(x_1, \dots, x_n)$
where $\sigma(x_1, \dots, x_n)$ is the next to last value under $\Phi$.
By $\isig^0_1$, there is a first such atom and a bound on the witnesses needed to show that all smaller candidates
are not such an atom. Once we passed such a stage, no change occurs in $d_{e,i,t}$ and
its value must also be above the stage where all witnesses are found. After such
a stage, $\Rcal_{e,i}$ will never need attention again.
\end{proof*}

The construction of an atomic model of $T$ from any infinite $f$-thin set for color~$i$
is exactly the same as in~\cite[Theorem~4.1]{Hirschfeldt2009atomic}.
\end{proof}

We have seen that the thin set theorem is strong enough to carry out some constructions
involving Ramsey's theorem. We will now prove that it is a strict weakening of Ramsey's theorem.
Furthermore, we shall see that the thin set theorems with bounded range form
a strictly decreasing hierarchy in reverse mathematics. Since the thin set theorem
with unbounded range is a consequence of each of them, it implies none of them
over~$\rca$.

\section{The thin set theorem and strong reducibility}\label{sect:ts-strong-reduc}

We start our analysis of the weakness of the thin set theorem
with partitions of integers like we did with Ramsey's theorem. 
Every computable partition has an infinite 
computable set avoiding one of its parts. The natural
reducibility to consider is therefore strong computable reducibility.
In this section, we show that~$\ts^1_k \not \leq_{sc} \ts^1_{k+1}$.

\subsection{Negative preservation results}

We now prove that the thin set theorems with bounded range
cannot preserve too many hyperimmunities simultaneously.
More precisely, we shall see that~$\ts^2_k$
does not admit preservation of~$k$ hyperimmunities for every~$k \geq 2$.
The following theorem builds a function with stronger properties than necessary.
The property (ii) will be useful to show that a stable version of 
the ascending descending sequence does not admit preservation of 2 hyperimmunities.

\begin{theorem}\label{thm:strength-ramsey-delta2-hyperimmune}
There is a stable computable coloring $f : [\omega]^2 \to \omega$ such that for each~$i$,
\begin{itemize}
	\item[(i)] $\{ x : \lim_s f(x,s) \neq i \}$ is hyperimmune;
	\item[(ii)] $f(x,y) \neq i \wedge f(y,z) \neq i \imp f(x,z) \neq i$ 
		for each each~$x < y < z \in \omega$.
\end{itemize}
\end{theorem}
\begin{proof}
The construction of the function~$f$ is done by a finite injury priority argument
with a movable marker procedure. We want to satisfy the following scheme of requirements for each~$e \in \omega$ and~$i \in \omega$,
where~$B_i = \{ x : \lim_s f(x, s) = i \}$:
\[
\Rcal_{e,i} : \Phi_e \mbox{ does not dominate } p_{\overline{B_i}}
\]

The requirements are given the usual priority ordering.
We proceed by stages, maintaining a sequence of sets~$B_0, B_1, \dots$
which represent the limit of the function~$f$.
At stage 0, $B_{i,0} = \emptyset$ for each~$i$
and $f$ is nowhere defined. Moreover, each requirement~$\Rcal_{e,i}$
is given a movable marker~$m_{e,i}$ initialized to~0.

A strategy for~$\Rcal_{e,i}$ \emph{requires attention at stage~$s+1$}
if $\Phi_{e,s}(m_{e,i}) \downarrow = n$ for some~$n < s$ and~$\overline{B_{i,s}} \cap [m_{e,i}, n] \neq \emptyset$.
The strategy sets~$B_{i,s+1} = B_{i,s} \cup [m_{e,i},n]$, 
and~$B_{j,s+1} = B_{j,s} \setminus [m_{e,i},n]$ for every~$j \neq i$.
Then it is declared~\emph{satisfied} until some strategy of higher priority changes its marker. 
Each marker~$m_{e',i'}$ of strategies of lower priorities is assigned the value~$s+1$.

At stage~$s+1$, assume that~$B_{0,s} \cup \dots \cup B_{k-1,s} = [0,s)$
and that~$f$ is defined for each pair over~$[0,s)$.
For each~$x \in [0,s)$, set~$f(x,s) = i$ for the unique~$i$ such that~$x \in B_{i,s}$.
If some strategy requires attention at stage~$s+1$, take the least one
and satisfy it.
If no such requirement is found, set~$B_{0,s+1} = B_{0,s} \cup \{s\}$
and~$B_{i,s+1} = B_{i,s}$ for $i > 0$.
Then go to the next stage. This ends the construction.

Each time a strategy acts, it changes the markers of strategies of lower priority, and is declared satisfied.
Once a strategy is satisfied, only a strategy of higher priority can injury it.
Therefore, each strategy acts finitely often and the markers stabilize.
It follows that the $B$'s also stabilize and that~$f$ is a stable function.

\begin{claim}
For every~$i \in \omega$ and every~$x < y < z$, $f(x,y) \neq i \wedge f(y,z) \neq i \imp f(x,z) \neq i$.
\end{claim}
\begin{proof*}
Suppose that $f(x,y) \neq i$ but $f(x,z) = i$
for some~$i < k$. Let~$s \leq z$ be the least stage such that~$f(x, t) = i$
for every~$t \in [s+1, z]$. At stage~$s+1$, some strategy~$\Rcal_{e,i}$
moved to~$B_i$ the whole interval~$[m_{e,i},s]$. Since $m_{e',i'} \leq m_{e,i}$
for every strategy~$\Rcal_{e',i'}$ of higher priority, none of the elements in~$[m_{e,i}, s]$ leave~$B_i$ before stage~$z+1$.
As $f(x,y) \neq i$, $y \not \in [s+1,z]$ so $y \in [m_{e,i},s]$.
Therefore $y \in B_{i,z}$ and thus $f(y, z) = i$.
\end{proof*}

\begin{claim}
For every~$e, i \in \omega$, $\Rcal_{e,i}$ is satisfied.
\end{claim}
\begin{proof*}
By induction over the priority order. Let~$s_0$ be a stage after which
no strategy of higher priority will ever act. By construction, $m_{e,i}$ will not change after stage~$s_0$.
If $\Phi_e$ is total and dominates $p_{\overline{B_i}}$, 
$\Phi_e(m_{e,i})$ will eventually halt and output some~$n$
such that $\overline{B_i} \cap [m_{e,i}, n] \neq \emptyset$
and therefore $\Rcal_{e,i}$ will require attention at some stage~$s \geq u$.
As no strategy of higher priority ever acts after stage~$s_0$, $\Rcal_{e,i}$
will receive attention, be satisfied and never be injured.
\end{proof*}

Satisfying $\Rcal_{e,i}$ for every~$e \in \omega$ and $i < k$ guarantees that $f$
has no computable thin set.
This last claim finishes the proof of Theorem~\ref{thm:strength-ramsey-delta2-hyperimmune}.
\end{proof}

Note that we could have interleaved some lowness requirements
to ensure that the function~$f$ is low.

\begin{corollary}\label{cor:ts2-non-preservation}\ 
\begin{itemize}
	\item[(i)] $\sts^2$ does not admit preservation of hyperimmunity.
	\item[(ii)] For every~$k \geq 2$, $\sts^2_k$ does not admit preservation of $k$ hyperimmunities.
\end{itemize}
\end{corollary}
\begin{proof}
We first prove (i).
Let~$f : [\omega]^2 \to \omega$ be as in Theorem~\ref{thm:strength-ramsey-delta2-hyperimmune}
and let~$B_i = \{ x : \lim_s g(x, s) = i \}$.
Any infinite $f$-thin set $H$ for color~$i$ is an infinite subset of~$\overline{B_i}$.
In particular, $\overline{B_i}$ is not $H$-hyperimmune since~$p_H$ dominates $p_{\overline{B_i}}$.

We now prove~(ii).
Let~$f$ be as in case (i), and let~$g : [\omega]^2 \to k$
be the stable computable function defined by~$g(x, s) = f(x, s)$ for~$x < k-1$
and~$g(x, s) = f(k-1, s)$ otherwise. Let~$C_i = \{x : \lim_s g(x, s) = i \}$.
Notice that~$\overline{C_i} \subseteq \overline{B_i}$ for each~$i < k$
and therefore that~$\overline{C_i}$ is hyperimmune for each~$i < k$.
The argument is similar to case (i).
\end{proof}

\subsection{Strong preservation of hyperimmunity}

Because every computable instance of~$\ts^1_k$ having a computable solution,
$\ts^1_k$ admits preservation of $k$ hyperimmunities for every~$k$.
On the other hand, we would like to say that~$\ts^1_k$ does not \emph{combinatorially}
preserve $k$ hyperimmunities since Theorem~\ref{thm:strength-ramsey-delta2-hyperimmune} shows
the existence of a non-effective instance of~$\ts^1_k$ whose solutions do not preserve $k$ hyperimmunities.
This combinatorial notion of preservation is called \emph{strong preservation}.

\index{preservation!of k hyperimmunities (strong)}
\begin{definition}[Strong preservation of $k$ hyperimmunities]\ 
A~$\Pi^1_2$ statement~$\Psf$ \emph{admits strong preservation of $k$ hyperimmunities} 
if for each set $Z$, each list of $Z$-hyperimmune sets $B_0, \dots, B_{k-1}$ and each (arbitrary) $\Psf$-instance $X$, 
there exists a solution $Y$ to~$X$ such that the~$B$'s are $Y \oplus Z$-hyperimmune.
\end{definition}

We have seen through Theorem~\ref{thm:strength-ramsey-delta2-hyperimmune} that for every~$k \geq 2$,
$\ts^1_k$ does not admit strong preservation of~$k$ hyperimmunities.
The following theorem shows the optimality of Theorem~\ref{thm:strength-ramsey-delta2-hyperimmune}.

\begin{theorem}\label{thm:ts1-strong-preservation}
For every~$k \geq 2$, $\ts^1_{k+1}$ admits strong preservation of~$k$ hyperimmunities.
\end{theorem}

Theorem~\ref{thm:ts1-strong-preservation} is proven at the end of this section. We first state a few immediate corollaries.
Putting Theorem~\ref{thm:strength-ramsey-delta2-hyperimmune} and Theorem~\ref{thm:ts1-strong-preservation}
together, we obtain the desired separation over strong computable reducibility.

\begin{corollary}
For every~$\ell > k \geq 2$, $\ts^1_k \not \leq_{sc} \ts^1_\ell$
\end{corollary}

Using the computable equivalence between the problem of finding a infinite set thin
for an~$\Delta^0_2$ $\ell$-partition and~$\sts^2_\ell$, we deduce the following corollary.

\begin{corollary}
For every~$\ell > k \geq 2$, $\sts^2_k \not \leq_{c} \sts^2_\ell$
\end{corollary}
\begin{proof}
Fix~$\ell > k \geq 2$ and consider the~$\Delta^0_2$ $k$-partition~$B_0 \cup \dots \cup B_{k-1} = \omega$
of Theorem~\ref{thm:strength-ramsey-delta2-hyperimmune}.
By Shoenfield's limit lemma~\cite{Shoenfield1959degrees}, there exists a stable computable function~$g : [\omega]^2 \to k$
such that~$B_j = \{ x : \lim_s g(x,s) = j \}$ for each~$j < k$. Every infinite $g$-thin set
is thin for the~$B$'s. Fix any stable computable function~$f : [\omega]^2 \to \ell$
and let~$A_i = \{x : \lim_s f(x,s) = i\}$ for each~$i < m$. By Theorem~\ref{thm:ts1-strong-preservation},
there exists an infinite set~$H$ thin for the~$A$'s which does not compute an infinite set thin for the~$B$'s
(hence for~$g$). As~$H \oplus f$ computes an infinite $f$-thin set~$G$, $f$ has an infinite $f$-thin set
which does not compute an infinite $g$-thin set.
\end{proof}

The remainder of this section is devoted to the proof of Theorem~\ref{thm:ts1-strong-preservation}.
Fix some set~$C$ preserving the hyperimmunity of some sets~$B_0, \dots, B_{k-1}$
and fix some $(k+1)$-partition~$A_0 \cup \dots \cup A_k = \omega$.
We will construct a set~$G$ such that $G \cap \overline{A_i}$
is infinite for each~$i \leq k$ and all the $B$'s are $(G \cap \overline{A_i}) \oplus C$-hyperimmune
for some~$i \leq k$.
Our forcing conditions are Mathias conditions~$(F, X)$ where
$X$ is an infinite set of integers such that all the $B$'s are $X \oplus C$-hyperimmune.

\subsection{Forcing limitlessness}

We want to satisfy the following scheme of requirements to ensure that~$G \cap \overline{A_i}$
is infinite for each~$i \leq k$:
\[
\Qcal_p : (\exists m_0, \dots, m_k > p)[m_0 \in G \cap \overline{A_0} \wedge \dots \wedge m_k \in G \cap \overline{A_k}]
\]
We say that an~$(k+1)$-partition~$A_0 \cup \dots \cup A_k = \omega$ is \emph{non-trivial} 
if there exists no infinite set~$H$ homogeneous for the~$A$'s
such that all the $B$'s are $H \oplus C$-hyperimmune. Of course, every infinite set homogeneous for the~$A$'s
is thin for the~$A$'s, so if the partition~$A_0 \cup \dots \cup A_k = \omega$ is trivial, we succeed.
Therefore we will assume from now on that the partition is non-trivial.
A condition~$c = (F, X)$ \emph{forces $\Qcal_p$}
if there exist some~$m_0, \dots, m_k > p$ such that~$m_i \in F \cap \overline{A_i}$
for each~$i \leq k$. Therefore, if~$G$ satisfies~$c$ and~$c$ forces~$\Qcal_p$,
then~$G$ satisfies the requirement~$\Qcal_p$.
We now prove that the set of conditions forcing~$\Qcal_p$ is dense for each~$p \in \omega$.
Thus, every sufficiently generic filter will induce a set~$G$
such that $G \cap \overline{A_i}$ is infinite for each~$i \leq k$.

\begin{lemma}\label{lem:ts2-reduc-force-Q}
For every condition~$c$ and every~$p$, there is an extension forcing~$\Qcal_p$.
\end{lemma}
\begin{proof}
Fix some~$p \in \omega$. It is sufficient to show that given a condition~
$c = (F, X)$ and some~$i \leq k$,
there exists an extension~$d_0 = (E, Y)$ and some integer~$m_i > p$
that~$m_i \in E \cap \overline{A_i}$.
By iterating the process for each~$i \leq k$, we obtain an extension forcing~$\Qcal_p$.
Suppose for the sake of contradiction that~$X \cap \overline{A_i}$ is finite.
One can then $X$-compute an infinite set~$H \subseteq A_i$, 
contradicting non-triviality of the~$A$'s.
Therefore, there exists an~$m_i \in X \cap \overline{A_i}$ such that $m_i > p$.
The condition $d_0 = (F \cup \{m_i\}, X \setminus [0, m_i])$
is the desired extension.
\end{proof}

\subsection{Forcing preservation}

The second scheme of requirements consists of ensuring that
the sets~$B_0, \dots, B_{k-1}$ are all $(G \cap \overline{A_i}) \oplus C$-hyperimmune for some~$i \leq k$.
The requirements are of the following form for each tuple of indices~$\vec{e} = e_0, \dots, e_k$:
\[
\Rcal_{\vec{e}} : \bigwedge_{j < k} \Rcal_{e_0}^{A_0, B_j}
	\vee \dots \vee \bigwedge_{j < k} \Rcal_{e_k}^{A_k, B_j}
\]
where~$\Rcal_e^{A_i, B_j}$ is the requirement ``$\Phi^{(G \cap \overline{A_i}) \oplus C}_{e_0} \mbox{ does not dominate } p_{B_j}$''.
A condition~\emph{forces $\Rcal_{\vec{e}}$} if every set~$G$ satisfying this condition also satisfies requirement~$\Rcal_{\vec{e}}$.
The following lemma is the core of the forcing argument.

\begin{lemma}\label{lem:ts2-reduc-step}
For every condition~$c = (F, X)$, every~$i_0 < i_1 \leq k$,
every $j < k$ and every vector of indices~$\vec{e}$, there is
an extension~$d$ forcing either~$\Phi^{(G \cap \overline{A_{i_0}}) \oplus C}_{e_{i_0}}$ or
$\Phi^{(G \cap \overline{A_{i_1}}) \oplus C}_{e_{i_1}}$ not to dominate $p_{B_j}$.
\end{lemma}
\begin{proof}
Let~$f$ be the partial $X \oplus C$-computable function which on input~$x$,
searches for a finite set of integers~$U$ such that 
for every $2$-cover~$Z_{i_0} \cup Z_{i_1} = X$, 
there is some $i \in \{i_0,i_1\}$ and some set~$E \subseteq Z_i$ 
such that~$\Phi_{e_i}^{((F \cap \overline{A_i}) \cup E) \oplus C}(x) \downarrow \in U$.
If such a set~$U$ is found, then~$f(x) = max(U)+1$, otherwise~$f(x) \uparrow$.
We have two cases:
\begin{itemize}
	\item Case 1: The function~$f$ is total. By $X \oplus C$-hyperimmunity of~$B_j$,
	$f(x) \leq p_{B_j}(x)$ for some~$x$. Let~$U$ be the finite set witnessing~$f(x) \downarrow$.
	Letting $Z_{i_0} = X \cap \overline{A_{i_0}}$
	and~$Z_{i_1} = X \cap \overline{A_{i_1}}$, there is some~$i \in \{i_0,i_1\}$
	and some finite set~$E \subseteq X \cap \overline{A_i}$ 
	such that~$\Phi_{e_i}^{((F \cap \overline{A_i}) \cup E) \oplus C}(x) \downarrow \in U$.
	The condition~$d = (F \cup E, X \setminus [0, max(E)])$ is an extension
	forcing $\Phi_{e_i}^{(G \cap \overline{A_i}) \oplus C}(x) < f(x) \leq p_{B_j}(x)$.

	\item Case 2: There is some~$x$ such that~$f(x) \uparrow$.
	By compactness, the~$\Pi^{0,X \oplus C}_1$ class $\Ccal$ of sets~$Z_{i_0} \oplus Z_{i_1}$ such that
	$Z_{i_0} \cup Z_{i_1} = X$ and for every~$i \in \{i_0,i_1\}$
	and every set~$E \subseteq Z_i$ $\Phi_{e_i}^{((F \cap \overline{A_i}) \cup E) \oplus C}(x) \uparrow$
	is non-empty. As $\wkl$ admits preservation of $k$ hyperimmunities, there exists some
	2-cover~$Z_{i_0} \oplus Z_{i_1} \in \Ccal$
	such that all the~$B$'s are~$Z_{i_0} \oplus Z_{i_1} \oplus X \oplus C$-hyperimmune.
	Let~$i \in \{i_0,i_1\}$ be such that $Z_i$ is infinite. The condition~$d = (F, Z_i)$
	is an extension of~$c$ forcing~$\Phi_{e_i}^{(G \cap \overline{A_i}) \oplus C}(x) \uparrow$.
\end{itemize}
\end{proof}

As usual, the following lemma iterates Lemma~\ref{lem:ts2-reduc-step}
and uses the fact that~$k+1 > k$ to satisfy the requirement~$\Rcal_{\vec{e}}$.

\begin{lemma}\label{lem:ts2-reduc-force-R}
For every condition~$c$,
and every vector of indices~$\vec{e}$, there exists
an extension~$d$ forcing~$\Rcal_{\vec{e}}$.
\end{lemma}
\begin{proof}
Fix a condition~$c$, and iterate applications of Lemma~\ref{lem:ts2-reduc-step}
to obtain an extension~$d$ such that for each~$j < k$, 
$d$ forces~$\Phi_{e_i}^{(G \cap \overline{A_i}) \oplus C}$ not to dominate $p_{B_j}$ for $k$ different~$i$'s.
By the pigeonhole principle, there exists some~$i \leq k$
such that $d$ forces~$\Phi_{e_i}^{(G \cap \overline{A_i}) \oplus C}$ not to dominate $p_{B_j}$ for each~$j < k$.
Therefore, $d$ forces~$\Rcal_{\vec{e}}$.
\end{proof}

\subsection{Construction}

Thanks to Lemma~\ref{lem:ts2-reduc-force-Q} and Lemma~\ref{lem:ts2-reduc-force-R}, 
define an infinite descending sequence of conditions
$(\emptyset, \omega) \geq c_0 \geq \dots$ such that for each~$s \in \omega$,
\begin{itemize}
	\item[(a)] $c_s$ forces~$\Qcal_s$
	\item[(b)] $c_s$ forces~$\Rcal_{\vec{e}}$ if~$s = \tuple{\vec{e}}$
\end{itemize}
where~$c_s = (F_s, X_s)$. Let~$G = \bigcup_s F_s$.
By (a), $G \cap \overline{A_i}$ is infinite for every~$i \leq k$ and by (b),
$G$ satisfies each requirement~$\Rcal_{\vec{e}}$.
This finishes the proof of Theorem~\ref{thm:ts1-strong-preservation}.

\section{The weakness of the thin set theorem for pairs}\label{sect:preservation-ts2-omega-models}

There is a fundamental difference in the way we proved that~$\rt^1_k \not \leq_{sc} \rt^1_\ell$
and that~$\ts^1_\ell \not \leq_{sc} \ts^1_k$ whenever~$k > \ell$. In the former case,
we have built an instance~$I$ of~$\rt^1_k$ satisfying some hyperimmunity properties,
and used those properties to construct a solution~$X$ to each instance of~$\rt^1_\ell$
which does not compute a solution to~$I$. We did not ensure that those hyperimmunity properties
are preserved relative to the solution~$X$, which prevents us from iterating the construction.
As it happens, those properties are not preserved as multiple applications
of~$\rt^1_\ell$ are sufficient to compute a solution to~$I$.
In the latter case, we proved that~$\ts^1_\ell$ has an instance whose solutions do not preserve some definitional property,
whereas each instance of~$\ts^1_k$ has a solution preserving it. This preservation enables us to
iterate the applications of~$\ts^1_k$ and build $\omega$-structures whose
second-order part is made of sets preserving this property.
We will take advantage of those observations to obtain new separations in reverse mathematics.

In this section, we prove that~$\ts^2_{k+1}$ does not imply~$\ts^2_k$ over~$\rca$
for every~$k \geq 2$. In particular, we answer several questions
asked by Cholak, Giusto, Hirst and Jockusch~\cite{Cholak2001Free} and by Mont\'alban~\cite{Montalban2011Open}
about the relation between~$\rt^2_2$ and~$\ts^2$.
Dorais et al.~\cite{Dorais2016uniform} proved that~$\rca \vdash \ts^n_k \imp \aca$ for~$n \geq 3$ 
whenever~$k$ is not large enough. Therefore we cannot hope to obtain the same separation
result over~$\rca$ for arbitrary tuples. However, we shall see that
$\ts^n_k$ is not computably reducible to~$\ts^n_{k+1}$ for~$n, k \geq 2$.

\begin{theorem}\label{thm:ts2-omega-models}
For every~$k \geq 2$, $\rca \wedge \coh \wedge \wkl \wedge \ts^2_{k+1} \nvdash \sts^2_k$.
\end{theorem}

Cholak et al.~\cite{Cholak2001Free} and Mont\'alban~\cite{Montalban2011Open} asked whether~$\ts^2$ implies~$\rt^2_2$
over~$\rca$. Thanks to Theorem~\ref{thm:ts2-omega-models}, we answer negatively,
noticing that~$\ts^2_2$ is the statement~$\rt^2_2$
and~$\rca \vdash \ts^2_k \imp \ts^2$ for every~$k \geq 2$ (see Dorais et al.~\cite{Dorais2016uniform}).

\begin{corollary}
$\ts^2$ does not imply~$\rt^2_2$ over~$\rca$.
\end{corollary}

Using the standard trick of prehomogeneous sets, we can generalize
from computable non-reducibility over pairs to arbitrary tuples.

\begin{corollary}\label{cor:ts-non-computable-reduction}
For every~$k, n \geq 2$ there exists a~$\Delta^0_n$ $k$-partition
$B_0 \cup \dots \cup B_{k-1} = \omega$ such that every computable
coloring $f : [\omega]^n \to k+1$ has an infinite $f$-thin set
computing no set thin for the~$B$'s. 
\end{corollary}
\begin{proof}
This is proved in a relativized form by induction over~$n \geq 2$.
The case~$n = 2$ is obtained by relativizing the proof of Theorem~\ref{thm:ts2-omega-models},
which indeed shows the existence of a $\Delta^0_2$ $k$-partition $B_0 \cup \dots \cup B_{k-1} = \omega$
such that every computable coloring $f : [\omega]^2 \to k+1$ has an infinite $f$-thin set computing no set thin for the~$B$'s.
Now assume it holds for some~$n$ in order to prove it for~$n+1$.
By the relativized low basis theorem~\cite{Jockusch197201},
let~$P \gg \emptyset^{(n-1)}$ be such that~$P' \leq \emptyset^{(n)}$.
Applying the induction hypothesis to~$P$,
there is a~$\Delta^{0,P}_2$ (hence~$\Delta^0_{n+1}$) $k$-partition~$B_0 \cup \dots \cup B_{k-1} = \omega$
such that each $P$-computable coloring $f: [\omega]^n \to k+1$
has an infinite $f$-homogeneous set~$H$ such that~$H \oplus P$
does not compute an infinite set thin for the~$B$'s.

Let~$f : [\omega]^{n+1} \to k+1$ be a computable coloring.
By Jockusch~\cite[Lemma 5.4]{Jockusch1972Ramseys}, there exists an infinite set~$C$ pre-homogeneous for~$f$
such that~$C \leq_T P$.
Let~$\tilde{f} : [C]^n \to k+1$ be the~$P$-computable coloring defined for each~$\sigma \in [C]^n$
by~$\tilde{f}(\sigma) = f(\sigma, a)$, where~$a \in A$, $a > max(\sigma)$.
Every $\tilde{f}$-thin set is~$f$-thin.
By definition of~$B_0 \cup \dots \cup B_{k-1} = \omega$, 
there exists an infinite~$\tilde{f}$-thin (hence $f$-thin) set~$H$
such that $H \oplus P$ does not compute an infinite set thin for the~$B$'s.
\end{proof}

\begin{corollary}
For every~$k, n \geq 2$, $\sts^n_k \not \leq_c \ts^n_{k+1}$
\end{corollary}

By Shoenfield's limit lemma~\cite{Shoenfield1959degrees}, a stable computable coloring over $(n+1)$-tuples
can be considered as a non-effective coloring over $n$-tuples. This consideration establishes
a bridge between preservation properties for colorings over~$(n+1)$-tuples and strong preservation 
properties for colorings over $n$-tuples. In particular, it enables us
to prove preservation results by induction over~$n$.
The following lemma has been proven by the author in its full generality in~\cite{Patey2015Combinatorial}.
Nevertheless we reprove it in the context of preservation of hyperimmunity.

\begin{lemma}\label{lem:ts-strong-to-weak}
For every~$k,n \geq 1$ and~$\ell \geq 2$, if~$\ts^n_\ell$ admits strong preservation of~$k$ hyperimmunities,
then $\ts^{n+1}_\ell$ admits preservation of~$k$ hyperimmunities.
\end{lemma}
\begin{proof}
Fix any set~$C$, $k$ $C$-hyperimmune sets~$A_0, \dots, A_{k-1}$ and any $C$-computable
coloring $f : [\omega]^{n+1} \to \ell$.
Consider the uniformly~$C$-computable sequence of sets~$\vec{R}$ defined for each~$\sigma \in [\omega]^n$ and~$i < \ell$ by
\[
R_{\sigma,i} = \{s \in \omega : f(\sigma,s) = i\}
\]
As~$\coh$ admits preservation of~$k$ hyperimmunities, there exists
some~$\vec{R}$-cohesive set~$G$ such that $G \oplus C$ preserves hyperimmunity
of the~$A$'s. The cohesive set induces a $(G \oplus C)'$-computable coloring~$\tilde{f} : [\omega]^n \to \ell$ defined by:
\[
(\forall \sigma \in [\omega]^n) \tilde{f}(\sigma) = \lim_{s \in G} f(\sigma,s)
\]
As ~$\ts^n_\ell$ admits strong preservation of $k$ hyperimmunities,
there exists an infinite $\tilde{f}$-thin set~$H$ such that
$H \oplus G \oplus C$ preserves hyperimmunity
of the~$A$'s. $H \oplus G \oplus C$ computes an infinite $f$-thin set.
\end{proof}

Using Theorem~\ref{thm:ts1-strong-preservation} together with Lemma~\ref{lem:ts-strong-to-weak}, 
we deduce the following corollary.

\begin{corollary}\label{cor:ts2-preservation}
For every~$k \geq 2$, $\ts^2_{k+1}$ admits preservation of~$k$ hyperimmunities.
\end{corollary}

We are now ready to prove Theorem~\ref{thm:ts2-omega-models}.

\begin{proof}[Proof of Theorem~\ref{thm:ts2-omega-models}]
Fix some~$k \geq 2$.
By Theorem~\ref{thm:coh-hyperimmunity-preservation} and by the hyperimmune-free basis theorem~\cite{Jockusch197201},
$\coh$ and $\wkl$ admit preservation of $k$ hyperimmunities.
By Corollary~\ref{cor:ts2-preservation}, so does $\ts^2_{k+1}$.
By Corollary~\ref{cor:ts2-non-preservation}, $\sts^2_k$ does not admit preservation of~$k$ hyperimmunities.
An application of Lemma~\ref{lem:intro-reduc-preservation-separation} completes the proof.
\end{proof}

\section{The weakness of the thin set theorem for tuples}\label{sect:preservation-ts-omega-models}

In this section, we extend the preservation of hyperimmunity of the thin set theorem for pairs
to arbitrary tuples, using the same construction pattern as Wang~\cite{Wang2014Some}.
We deduce that $\ts^n_\ell$ does not imply~$\ts^n_k$ over~$\rca$ whenever~$\ell$ is large enough,
which is informally the strongest result we can obtain since Proposition~5.3 in Dorais et al.~\cite{Dorais2016uniform}
states that $\rca \vdash \aca \biimp \ts^n_k$ for~$n \geq 3$ whenever~$k$ is not large enough.

\begin{theorem}\label{thm:ts-arbitrary-preservation}
For every~$k,n \geq 1$, $\ts^n_\ell$ admits strong preservation of $k$ hyperimmunities 
for sufficiently large~$\ell$. 
\end{theorem}

Using the fact that~$\rca \vdash \ts^n_\ell \imp \ts^n$ for every~$n, \ell \geq 2$,
we obtain the following preservation result for~$\ts$.

\begin{corollary}
For every~$k \geq 1$, $\ts$ admits strong preservation of $k$ hyperimmunities.
\end{corollary}

Thanks to the existing preservations of hyperimmunity
and Proposition~2.4 from Wang~\cite{Wang2014Definability}, we deduce the following separations
over~$\omega$-models.

\begin{corollary}
For every~$k \geq 2$, $\rca \wedge \coh \wedge \wkl \wedge \ts^2_{k+1} \wedge \ts \nvdash \sts^2_k$.
\end{corollary}

The remainder of this section is devoted to the proof of Theorem~\ref{thm:ts-arbitrary-preservation}.

\subsection{Proof structure}

We shall follow the proof structure
of strong cone avoidance by Wang~\cite{Wang2014Some}. Fix some~$k \geq 1$.
The induction works as follows:

\begin{itemize}
	\item[(A1)] In section~\ref{sect:preservation-ts2-omega-models}
	we proved that~$\ts^1_{k+1}$ admits strong preservation of $k$ hyperimmunities.
	This is the base case of our induction.
	\item[(A2)] Assuming that for each~$t \in (0,n)$, $\ts^t_{d_t+1}$ admits
	strong preservation of $k$ hyperimmunities, we prove
	that~$\ts^n_{d_{n-1}+1}$ admits preservation of $k$ hyperimmunities.
	This is done by Lemma~\ref{lem:ts-strong-to-weak}.
	\item[(A3)] Then we prove that~$\ts^n_{d_n+1}$ admits strong preservation of
	$k$ hyperimmunities where
	\[
	d_n = d_1 d_{n-1} + \sum_{0 < t < n} d_t d_{n-t}
	\]
\end{itemize}

Properties (A1) and (A2) are already proven. We now prove property (A3). It is again done in several steps.
Fix a coloring $f : [\omega]^n \to d_n+1$ and a set $C$ preserving hyperimmunity of $k$ sets~$A_0, \dots, A_{k-1}$.

\begin{itemize}
	\item[(S1)] First, we construct an infinite set $D \subseteq \omega$ such that $D \oplus C$ preserves hyperimmunity of the~$A$'s
	and a sequence $(I_\sigma : 0 < |\sigma| < n)$ such that for each $t \in (0,n)$ and each $\sigma \in [\omega]^t$
	\begin{itemize}
		\item[(a)] $I_\sigma$ is a subset of $\{0, \dots, d_n\}$ with at most $d_{n-t}$ many elements  
		\item[(b)] $(\exists b)(\forall \tau \in [D \cap (b,+\infty)]^{n-t}) f(\sigma,\tau) \in I_\sigma$
	\end{itemize}

	\item[(S2)]
	Second, we construct an infinite set~$E \subseteq D$ such that $E \oplus C$ preserves hyperimmunity of the~$A$'s
	and a sequence $(I_t : 0 < t < n)$ such that for each $t \in (0,n)$
	\begin{itemize}
		\item[(a)] $I_t$ is a subset of~$\{0,\dots, d_n\}$ of size at most $d_t d_{n-t}$
		\item[(b)] $(\forall \sigma \in [E]^t)(\exists b)(\forall \tau \in [E \cap (b,+\infty)]^{n-t}) f(\sigma,\tau) \in I_t$
	\end{itemize}

	\item[(S3)] Third, we construct a sequence~$(\xi_i \in [E]^{<\omega} : i < \omega)$ such that
	\begin{itemize}
		\item[(a)] The set $G = \bigcup_i \xi_i$ is infinite and $G \oplus C$ preserves hyperimmunity of the~$A$'s
		\item[(b)] $|f([\xi_i]^n)| \leq d_{n-1}$ and~$max(\xi_i) < min(\xi_{i+1})$ for each~$i < \omega$
		\item[(c)] For each~$t \in (0,n)$ and~$\sigma \in [\bigcup_{j < i} \xi_j]^t$,
		$f(\sigma,\tau) \in I_t$ for all~$\tau \in [\bigcup_{j \geq i} \xi_j]^{n-t}$
	\end{itemize}

	\item[(S4)] Finally, we build an infinite set~$H \subseteq G$
	such that~$H \oplus C$ preserves hyperimmunity of the~$A$'s
	and~$|f([H]^n)| \leq d_n$.
\end{itemize}

\subsection{Generalized cohesiveness}

Before proving that~$\ts^n_{d_n+1}$ admits strong preservation of~$k$ hyperimmunities,
we need to prove strong preservation for a generalized notion of cohesiveness
already used by the author in~\cite{Patey2015Combinatorial}.
Cohesiveness can be seen as the problem which takes as an input a coloring of pairs $f : [\omega]^2 \to \ell$
and fixes the first parameter to obtain an infinite sequence of colorings of integers $f_x : \omega \to \ell$ for each $x \in \omega$.
A solution to this problem is an infinite set~$G$ which is eventually homogeneous for each coloring $f_x$.

Going further in this approach, we can consider that cohesiveness is a degenerate
case of the problem which takes as an input a coloring of pairs $f : [\omega]^2 \to \omega$ using infinitely many colors,
and fixes again the first parameter to obtain an infinite sequence of colorings of integers $f_x : \omega \to \omega$.
A solution to this problem is an infinite set~$G$ such that for each color $i$, either eventually the color will be avoided
by $f_x$ over $G$, or $G$ will be eventually homogeneous for~$f_x$ with color $i$.

We can generalize the notion to colorings over tuples $f : [\omega]^n \to \omega$,
seeing $f$ as an infinite sequence of colorings over $t$-uples $f_{\sigma} : [\omega]^t \to \omega$ for each $\sigma \in [\omega]^{n - t}$. We will create a set~$G$ such that at most $d_t$ colors will appear for arbitrarily large pairs over~$G$ for each function $f_{\sigma}$.
This set will be constructed by applying $\ts^t_{d_t+1}$ to $f_{\sigma}$ for each~$\sigma$.

We do not need Theorem~\ref{thm:generalized-cohesivity-strong-avoidance} in its 
full generality to complete our step (S1). However, it
will be useful in a later section for proving that the free set theorem admits
preservation of $k$ hyperimmunities.

\begin{theorem}\label{thm:generalized-cohesivity-strong-avoidance}
Fix a coloring $f : [\omega]^n \to \omega$, some $t \leq n$ and
suppose that
$\ts^s_{d_s+1}$ admits strong preservation of~$k$ hyperimmunities for each $s \in (0, t]$. 
For every set $C$ preserving hyperimmunity of some sets~$A_0, \dots, A_{k-1}$,
there exists an infinite set $G$ such that $G \oplus C$ preserves hyperimmunity of the~$A$'s
and for every $\sigma \in [\omega]^{<\omega}$ such that $n-t \leq \card{\sigma} < n$, 
$$
\card{\set{x : (\forall b)(\exists \tau \in [G \cap (b, +\infty)]^{n-|\sigma|}) f(\sigma, \tau) = x}} \leq d_{n-|\sigma|}
$$
\end{theorem}
\begin{proof}
Our forcing conditions are Mathias conditions $(F, X)$ where $X \oplus C$ preserves hyperimmunity
of the~$A$'s.
Lemma~\ref{lem:coh-preservation-lemma} states that for every set~$G$ which is sufficiently
generic for~$(F,X)$, $G \oplus C$ preserves $k$ hyperimmunities. It suffices therefore to prove the following lemma.

\begin{lemma}\label{lem:gen-coh-function}
For every condition $(F, X)$ and $\sigma \in [\omega]^{<\omega}$ such that $n-t \leq \card{\sigma} < n$,
for every finite set $I$ such that $\card{I} = d_{n-|\sigma|}$, 
there exists an extension $(F, \tilde{X})$ such that 
$$
\{f(\sigma, \tau) : \tau \in [\tilde{X}]^{n-|\sigma|}\} \subseteq I
\hspace{10pt} \mbox{ or } \hspace{10pt} 
I \not \subseteq \{f(\sigma, \tau) : \tau \in [\tilde{X}]^{n-|\sigma|}\}
$$
\end{lemma}
\begin{proof*}
Define the function $g : [X]^{n - |\sigma|} \to I \cup \{\bot\}$
by $g(\tau) = f(\sigma, \tau)$ if $f(\sigma, \tau) \in I$
and $g(\tau) = \bot$ otherwise. By strong preservation of~$k$ hyperimmunities of $\ts^{n - |\sigma|}_{d_{n-|\sigma|}+1}$,
there exists an infinite subset $\tilde{X} \subseteq X$ such that $\tilde{X} \oplus C$
preserves hyperimmunity of the~$A$'s and $\card{\{g(\tau) : \tau \in [\tilde{X}]^{n-|\sigma|}\}} \leq d_{n-|\sigma|}$.
The condition $(F, \tilde{X})$ is the desired extension.
\end{proof*}

Using Lemma~\ref{lem:coh-preservation-lemma}
and Lemma~\ref{lem:gen-coh-function}, one can define an infinite
descending sequence of conditions $(\emptyset, \omega) \geq c_0 \geq c_1 \geq \dots$
such that for each $s \in \omega$
\begin{itemize}
  \item[1.] $c_s = (F_s, X_s)$ with $\card{F_s} \geq s$
  \item[2.] $c_s$ forces $W_e^{G \oplus C} \neq A_i$ if $s = \tuple{e,i}$
  \item[3.] $\{f(\sigma, \tau) : \tau \in [X_s]^{n-|\sigma|}\} \subseteq I$
  or $I \not \subseteq \{f(\sigma, \tau) : \tau \in [X_s]^{n-|\sigma|}\}$ if $s = \tuple{\sigma, I}$
  and $\card{I} = d_{n-|\sigma|}$.
\end{itemize}
The set $G = \bigcup_s F_s$ is an infinite set 
such that $G \oplus C$ preserves hyperimmunity of the~$A$'s. We claim that $G$ satisfies the desired properties.
Fix a $\sigma \in [\omega]^{<\omega}$ such that $n-t \leq \card{\sigma} < n$. Suppose that there exists $d_{n - |\sigma|} + 1$ elements 
$x_0, \dots, x_{d_{n-|\sigma|}}$ such that $(\forall b)(\exists \tau \in [G \cap (b, +\infty)]^{n-|\sigma|}) f(\sigma, \tau) = x_i$
for each $i \leq d_{n-|\sigma|}$.
Let $I = \{x_0, \dots, x_{d_{n-|\sigma|}-1}\}$.
By step $s = \tuple{\sigma, I}$, $G$ satisfies $(F_s, X_s)$ such that $\{f(\sigma, \tau) : \tau \in [X_s]^{n-|\sigma|}\} \subseteq I$
or $I \not \subseteq \{f(\sigma, \tau) : \tau \in [X_s]^{n-|\sigma|}\}$. In the first case it contradicts
the choice of $x_{d_{n-|\sigma|}}$ and in the second case it contradicts the choice of an element of $I$.
This finishes the proof of Theorem~\ref{thm:generalized-cohesivity-strong-avoidance}.
\end{proof}

\subsection{Step (S1) : Construction of the set~$D$}

We start with the construction of an infinite set~$D \subseteq \omega$ such that $D \oplus C$ preserves hyperimmunity of the~$A$'s
and a sequence $(I_\sigma : 0 < |\sigma| < n)$ such that for each $t \in (0,n)$ and each $\sigma \in [\omega]^t$
\begin{itemize}
	\item[(a)] $I_\sigma$ is a subset of $\{0, \dots, d_n\}$ with at most $d_{n-t}$ many elements  
	\item[(b)] $(\exists b)(\forall \tau \in [G \cap (b,+\infty)]^{n-t}) f(\sigma,\tau) \in I_\sigma$
\end{itemize}

Let $D$ be the set constructed in Theorem~\ref{thm:generalized-cohesivity-strong-avoidance}
for $t = n-1$. For each $\sigma \in [\omega]^{<\omega}$ such that $0 < \card{\sigma} < n$,
let
$$
I_\sigma = \{x \leq d_n :  (\forall b)(\exists \tau \in [G \cap (b, +\infty)]^{n-|\sigma|}) f(\sigma, \tau) = x\}
$$
By choice of $D$, the set $I_\sigma$ has at most $d_{n-|\sigma|}$ many elements.
Moreover, for each $y \leq d_n$ such that $y \not \in I_\sigma$, there exists a bound $b_y$
such that $(\forall \tau \in [D \cap (b_y, +\infty)]^{n-|\sigma|}) f(\sigma, \tau) \neq x$.
So taking $b = max(b_y : y \leq d_n \wedge y \not \in I_\sigma)$, we obtain
$$
(\forall \tau \in [D \cap (b,+\infty)]^{n-|\sigma|}) f(\sigma,\tau) \in I_\sigma
$$

\subsection{Step (S2) : Construction of the set~$E$}

We now construct an infinite set~$E \subseteq D$ such that $E \oplus C$ preserves hyperimmunity of the~$A$'s
and a sequence $(I_t : 0 < t < n)$ such that for each $t \in (0,n)$
\begin{itemize}
	\item[(a)] $I_t$ is a subset of~$\{0,\dots, d_n\}$ of size at most $d_t d_{n-t}$
	\item[(b)] $(\forall \sigma \in [E]^t)(\exists b)(\forall \tau \in [E \cap (b,+\infty)]^{n-t}) f(\sigma,\tau) \in I_t$
\end{itemize}

For each $t \in (0,n)$ and~$\sigma \in [\omega]^t$, let $F_t(\sigma) = I_\sigma$.
Using strong preservation of $k$ hyperimmunities
of $\ts^t_{d_t+1}$,
we build a finite sequence $D \supseteq E_1 \supseteq \dots \supseteq E_{n-1}$
such that for each $t \in (0,n)$
\begin{itemize}
	\item[1.] $E_t \oplus C$ preserves hyperimmunity of the~$A$'s
	\item[2.] $|F_t([E_t]^t)| \leq d_t$
\end{itemize}
Let $E = E_{n-1}$ and $I_t = \bigcup F_t([E]^t)$ for each $t \in (0,n)$.
As for each~$\sigma \in [E]^t$, $|F_t(\sigma)| = |I_\sigma| \leq d_{n-t}$, $|I_t| \leq d_t d_{n-t}$, so property (a) holds.
We now check that property (b) is satisfied.
Fix a $\sigma \in [E]^t$. By definition of~$I_t$, $F_t(\sigma) = I_\sigma \subseteq I_t$.
As $E \subseteq D$,
\[
(\exists b)(\forall \tau \in [E \cap (b,+\infty)]^{n-t}) f(\sigma,\tau) \in I_\sigma \subseteq I_t
\]

\subsection{Step (S3) : Construction of the set~$G$}

Given the set~$E$ and the sequence of sets of colors~$(I_t : 0 < t < n)$,
we will construct a sequence~$(\xi_i \in [E]^{<\omega} : i < \omega)$ such that
\begin{itemize}
	\item[(a)] The set $G = \bigcup_i \xi_i$ is infinite and $G \oplus C$ preserves hyperimmunity of the~$A$'s
	\item[(b)] $|f([\xi_i]^n)| \leq d_{n-1}$ and~$max(\xi_i) < min(\xi_{i+1})$ for each~$i < \omega$
	\item[(c)] For each~$t \in (0,n)$ and~$\sigma \in [\bigcup_{j < i} \xi_j]^t$,
	$f(\sigma,\tau) \in I_t$ for all~$\tau \in [\bigcup_{j \geq i} \xi_j]^{n-t}$
\end{itemize}

We construct our set~$G$ by Mathias forcing~$(\sigma, X)$ where
$X$ is an infinite subset of~$E$ such that $X \oplus C$ preserves hyperimmunity
of the~$A$'s. Using property~(b) of~$E$, we can easily
construct an infinite sequence $(\xi_i \in [E]^{<\omega} : i < \omega)$
satisfying properties~(b) and~(c) of step (S3). The following lemma
shows how to satisfy property~(a).

\begin{lemma}\label{lem:ts-lemma-set-g}
Fix a condition~$(\sigma, X)$, some~$e \in \omega$ and some~$j < k$.
There exists an extension~$(\sigma\xi, Y)$
with~$|f([\xi]^n)| \leq d_{n-1}$, forcing~$\Phi_e^{G \oplus C}$ not to dominate $p_{A_j}$.
\end{lemma}
\begin{proof}
Let $h$ be the partial~$X \oplus C$-computable function which on input~$x$,
searches for a finite set of integers~$U$ such that 
for every coloring $g : [X]^n \to d_n+1$, there is a set~$\xi \in [X]^{<\omega}$
such that~$|g([\xi]^n)| \leq d_{n-1}$ and~$\Phi_e^{\sigma\xi \oplus C}(x) \downarrow \in U$.
If such a set~$U$ is found, then $f(x) = max(U) + 1$, otherwise $f(x) \uparrow$.
We have two cases:
\begin{itemize}
	\item Case 1: the function is total. By $X \oplus C$-hyperimmunity of $A_j$,
	$f(x) \leq p_{A_j}(x)$ for some~$x$. Let~$U$ be the finite set witnessing~$f(x) \downarrow$.
	Letting~$g = f$, there exists a set~$\xi \in [X]^{<\omega}$
	such that~$|f([\xi]^n)| \leq d_{n-1}$ and~$\Phi_e^{\sigma\xi \oplus C}(x) \downarrow \in U$.
	The condition~$d = (\sigma\xi, X)$ is an extension forcing $\Phi_e^{G \oplus C}(x) < f(x) \leq p_{A_j}(x)$.

	\item Case 2: There is some~$x$ such that~$f(x) \uparrow$.
	By compactness, the $\Pi^{0, X \oplus C}_1$ class~$\Ccal$ of colorings~$g : [X]^n \to d_n+1$
	such that for every set~$\xi \in [X]^{<\omega}$
	satisfying~$|g([\xi]^n)| \leq d_{n-1}$, $\Phi_e^{\sigma\xi \oplus C}(x) \uparrow$
	is non-empty.
	As $\wkl$ admits preservation of $k$ hyperimmunities,
	there is some coloring~$g \in \Ccal$ such that~$g \oplus X \oplus C$
	preserves hyperimmunity of the~$A$'s.
	By preservation of $k$ hyperimmunities of~$\ts^n_{d_{n-1}+1}$,
	there exists an infinite subset~$Y \subseteq X$ such that
	$Y \oplus C$ preserves hyperimmunity of the~$A$'s
	and~$|g([Y]^n)| \leq d_{n-1}$.
	The condition~$(\sigma, Y)$ forces~$\Phi_e^{G \oplus C}(x) \uparrow$.
\end{itemize}
\end{proof}

Using Lemma~\ref{lem:ts-lemma-set-g} and property (b) of the set~$E$,
we can construct an infinite descending sequence of conditions~$(\epsilon,E) \geq c_0 \geq \dots$
such that for each~$s \in \omega$
\begin{itemize}
	\item[(i)] $\sigma_{s+1} = \sigma_s\xi_s$ with $|\sigma_s| \geq s$ and ~$f([\xi_s]^n) \leq d_{n-1}$
	\item[(ii)] $f(\sigma,\tau) \in I_t$ for each~$t \in (0,n)$, $\sigma \in [\sigma_s]^t$ and~$\tau \in [X]^{n-t}$.
	\item[(iii)] $c_s$ forces~$\Phi_e^{G \oplus C}$ not to dominate $p_{A_j}$ if~$s = \tuple{e,j}$
\end{itemize}
where~$c_s = (\sigma_s, X_s)$. The set~$G = \bigcup_s \sigma_s$
satisfies the desired properties.

\subsection{Step (S4) : Construction of the set~$H$}

Finally, we build an infinite set~$H \subseteq G$
such that~$H \oplus C$ preserves hyperimmunity of the~$A$'s
and~$|f([H]^n)| \leq d_n$.

For each~$i < \omega$, let~$J_i = f([\xi_i]^n)$.
By property (b) of step (S3), $J_i$ is a subset of~$\{0, \dots, d_n\}$
such that $|J_i| \leq d_{n-1}$. For each subset~$J \subseteq \{0, \dots, d_n\}$
of size~$d_{n-1}$, define the set
\[
Z_J = \{ x \in G : (\exists i) x \in \xi_i \wedge f([\xi_i]^n) \subseteq J \}
\]
There exists finitely many such~$J$'s, and the~$Z$'s form a partition of~$G$.
Apply strong preservation of $k$ hyperimmunities of~$\ts^1_{d_1+1}$
to obtain a finite set~$S$ of $J$'s of such that~$|S| \leq d_1$
and an infinite set~$H \subseteq \bigcup_{J \in S} Z_J \subseteq G$
such that~$H \oplus G \oplus C$ preserves hyperimmunity of the~$A$'s.

\begin{lemma}
$|f([H]^n)| \leq d_n$
\end{lemma}
\begin{proof}
As~$H \subseteq G$, any~$\sigma \in [H]^n$ can be decomposed into $\rho^\concat \tau$ for some~$\rho \in [\xi_i]^{<\omega}$
and some~$\tau \in [\bigcup_{j \geq i} \xi_j]^{<\omega}$ with~$|\rho| > 0$.
If~$|\tau| = 0$ then $f(\sigma) \in \bigcup_{J \in S} J$  by definition of~$H$.
If~$|\tau| > 0$, then~$f(\sigma) \in I_{|\rho|}$ by property (c) of step (S3).
In any case
\[
f(\sigma) \in (\bigcup_{J \in S} J) \cup (\bigcup_{t \in (0,n)} I_t)
\]
Recall that $|S| \leq d_1$, $|J| = d_{n-1}$ for each~$J \in S$, and~$|I_t| \leq d_t d_{n-t}$ for each~$t \in (0,n)$.
We obtain therefore the desired inequality.
\end{proof}

This completes property (A3) and the proof of Theorem~\ref{thm:ts-arbitrary-preservation}.

\section{The free set theorem}\label{sect:fs-omega-models}

The free set theorem is a strengthening of the thin set theorem
in which every member of a free set is a witness of thinness of the same set.
Indeed, if~$H$ is an infinite~$f$-free set for some function~$f$,
for every~$a \in H$, $H \setminus \{a\}$ is $f$-thin with witness color~$a$.
See Theorem~3.2 in~\cite{Cholak2001Free} for a formal version of this claim.

\index{free set theorem}
\index{fs@$\fs^n$|see {free set theorem}}
\begin{definition}[Free set theorem]
Given a coloring $f : [\Nb]^n \to \Nb$, an infinite set~$H$
is \emph{$f$-free} if for every~$\sigma \in [H]^n$,
$f(\sigma) \in H \imp f(\sigma) \in \sigma$.
For every~$n \geq 1$, $\fs^n$ is the statement
``Every coloring $f : [\Nb]^n \to \Nb$ has a free set''.
\end{definition}

\index{sfs@$\sfs^n$|see {free set theorem}}
Again, $\sfs^2$ is the restriction of~$\fs^2$ to stable colorings
and $\fs$ is the statement~$(\forall n)\fs^n$.
Cholak et al.~\cite{Cholak2001Free} studied the thin set theorem with infinitely many colors as a weakening
of the free set theorem. Cholak et al.~\cite{Cholak2001Free}
proved that~$\rca \vdash \rt^n_2 \imp \fs^n \imp \ts^n$ for every~$n \geq 2$.
Wang~\cite{Wang2014Some} proved that~$\fs$ (hence~$\ts$) does not imply~$\aca$
over~$\omega$-models. The author~\cite{Patey2015Combinatorial} proved
that~$\fs$ does not imply~$\wkl$ (and in fact weak weak König's lemma) over~$\rca$.

The forcing notions used by Wang in~\cite{Wang2014Definability}
and by the author in~\cite{Patey2015Combinatorial} for constructing solutions to free set instances
both involve the thin set theorem for a finite, but arbitrary number of colors.
These constructions may suggest some relation between~$\fs^n$ and~$\ts^n_k$
for arbitrarily large~$k$, but the exact nature of this relation is currently unclear.

In this section, we use the preservation of $k$ hyperimmunities of the thin set theorem
to deduce similar preservation results for the free set theorem,
and thereby separate~$\fs$ from~$\rt^2_2$ over~$\rca$.
More precisely, we prove the following preservation theorem.

\begin{theorem}\label{thm:fs2-preservation-definitions}
For every~$k \in \omega$,
$\fs$ admits strong preservation of $k$ hyperimmunities.
\end{theorem}

Cholak et al.~\cite{Cholak2001Free}
asked whether any of~$\fs^2$, $\fs^2+\coh$ and
$\fs^2+\wkl$ imply~$\rt^2_2$ and
Hirschfeldt~\cite{Hirschfeldt2015Slicing} asked whether~$\fs^2+\wkl$ implies~$\srt^2_2$.
We answer all these questions negatively with the following corollary.

\begin{corollary}\label{cor:fs2-not-imply-sts2n}
For every~$k \geq 2$,
$\rca \wedge \wkl \wedge \coh \wedge \fs \wedge \ts^2_{k+1} \not \vdash \sts^2_k$.
\end{corollary}

In particular, since the statement~$\ts^2_2$ is nothing but Ramsey's theorem for pairs, we deduce
the following corollary.

\begin{corollary}\label{cor:fs-not-imply-rt22}
$\rca \wedge \fs \nvdash \rt^2_2$.
\end{corollary}

The remainder of this section is devoted to the proof of Theorem~\ref{thm:fs2-preservation-definitions}.
The proof is done by induction over the size of the tuples.
The base case of our induction states that $\fs^0$ admits strong preservation of $k$ hyperimmunities. 
Consider $\fs^0$ as a degenerate case of the free
set theorem, where an instance is a constant~$c$ and a solution to~$c$ is an infinite set~$H$
which does not contain~$c$. Indeed, a function~$f : [\omega]^0 \to \omega$ can be considered as a constant~$c$,
and a set~$H$ is $f$-free if for every~$\varepsilon \in [H]^0$, $f(\varepsilon) \in H \imp f(\varepsilon) \in \varepsilon$.
As~$f(\varepsilon) \not \in \varepsilon$, $f(\varepsilon) = c \not \in H$.
From now on, we will assume that~$\fs^t$ admits strong 
preservation of~$k$ hyperimmunities for every~$t \in [0,n)$.

We start with a lemma similar to Lemma~\ref{lem:ts-strong-to-weak}.

\begin{lemma}\label{lem:fs-strong-to-weak}
For every~$n \geq 1$ and~$k \geq 2$, if~$\fs^n$ admits strong preservation of~$k$ hyperimmunities,
then $\fs^{n+1}$ admits preservation of~$k$ hyperimmunities.
\end{lemma}
\begin{proof}
Fix any set~$C$, $k$ $C$-hyperimmune sets~$A_0, \dots, A_{k-1}$ and any $C$-computable
coloring $f : [\omega]^{n+1} \to \omega$.
Consider the uniformly~$C$-computable sequence of sets~$\vec{R}$ defined for each~$\sigma \in [\omega]^n$ and $y \in \omega$ by
\[
R_{\sigma,y} = \{s \in \omega : f(\sigma,s) = y\}
\]
As~$\coh$ admits preservation of~$k$ hyperimmunities, there exists
some~$\vec{R}$-cohesive set~$G$ such that $G \oplus C$ preserves hyperimmunity
of the~$A$'s. The cohesive set induces a coloring~$\tilde{f} : [\omega]^n \to \omega$ 
defined for each~$\sigma \in [\omega]^n$ by
\[
\tilde{f}(\sigma) = \cond{
	\lim_{s \in G} f(\sigma,s) & \mbox{ if it exists}\\
	0 & \mbox{ otherwise}\\
}
\]
As~$\fs^n$ admits strong preservation of $k$ hyperimmunities,
there exists an infinite $\tilde{f}$-free set~$H$ such that
$H \oplus G \oplus C$ preserves hyperimmunity
of the~$A$'s. In particular,
\[
(\forall \sigma \in [H]^n)(\forall y \in H \setminus \sigma)(\forall^\infty s)f(\sigma, s) \neq y
\]
Thus $H \oplus G \oplus C$ computes an infinite $f$-free set.
\end{proof}

\subsection{Trapped functions}

Although the notion of free set can be defined for every coloring
over tuples of integers, we shall restrict ourselves to a particular kind
of colorings: left trapped functions. The notion of trapped function
has been introduced by Wang in~\cite{Wang2014Some} to prove that~$\fs$ does not imply~$\aca$
over~$\omega$-models. It has been later reused by the author in~\cite{Patey2015Combinatorial} to separate~$\fs$
from~$\wwkl$ over~$\omega$-models. Given a string~$\sigma \in [\omega]^{<\omega}$ and some~$n < |\sigma|$,
we write~$\sigma(n)$ for the value of~$\sigma$ at the~$(n+1)$th position.

\index{left trapped function}
\index{right trapped function}
\begin{definition}
A function $f : [\omega]^n \to \omega$ is \emph{left (resp. right) trapped}
if for every $\sigma \in [\omega]^n$, $f(\sigma) \leq \sigma(n-1)$ (resp. $f(\sigma) > \sigma(n-1)$).
\end{definition}

The following lemma is again a particular case of a more general statement
proven by the author in~\cite{Patey2015Combinatorial}.
It follows from the facts that $\fs^n$ for right trapped functions
is strongly Weihrauch reducible to the diagonally non-computable principle ($\dnr$),
which itself is computably reducible to~$\fs^n$ for left trapped functions.

\begin{lemma}\label{lem:hyperimmunity-preservation-fs-trapped-to-untrapped}
For each $k, n \geq 1$, 
if $\fs^n$ for left trapped functions admits (strong) 
preservation of $k$ hyperimmunities then so does $\fs^k$.
\end{lemma}

It therefore suffices to prove strong preservation
of~$k$ hyperimmunities for left trapped functions.

\subsection{Case of left trapped functions}

In this part, we will prove the following theorem which,
together with Lemma~\ref{lem:hyperimmunity-preservation-fs-trapped-to-untrapped}
is sufficient to prove Theorem~\ref{thm:fs2-preservation-definitions} by induction over~$n$.

\begin{theorem}\label{thm:strong-preservation-fs-left-trapped}
For each~$k, n \geq 1$, if~$\fs^t$ admits strong preservation of~$k$ hyperimmunities
for each~$t \in [0,n)$, then so does~$\fs^n$ for left trapped functions.
\end{theorem}

The two following lemmas will ensure that the reservoirs of our forcing conditions
will have good properties, so that the conditions will be extensible.

\begin{lemma}\label{lem:fs-left-trapped-preserves-smaller}
Suppose that $\fs^t$ admits strong preservation of $k$ hyperimmunities for each $t \in (0,n)$ for some~$k \in \omega$.
Fix a set $C$, some $C$-hyperimmune sets~$A_0, \dots, A_{k-1}$, a finite set $F$ and an infinite set $X$ computable in $C$.
For every function $f : [X]^n \to \omega$ there exists an infinite
set $Y \subseteq X$ such that $Y \oplus C$ preserves hyperimmunity of the~$A$'s and 
$(\forall \sigma \in [F]^t)(\forall \tau \in [Y]^{n-t})f(\sigma, \tau) \not \in Y \setminus \tau$
for each $t \in (0,n)$.
\end{lemma}
\begin{proof}
Fix the finite enumeration $\sigma_0, \dots, \sigma_{m-1}$ of all elements of $[F]^t$ for all $t \in (0,n)$.
We define a finite decreasing sequence of sets~$X = Y_0 \supseteq Y_1 \supseteq \dots \supseteq Y_m$
such that for each~$s < m$
\begin{itemize}
	\item[(a)] all the $A$'s are $Y_{s+1} \oplus C$-hyperimmune
	\item[(b)] $\forall \tau \in [Y_{s+1}]^{n - |\sigma_s|}f(\sigma_s, \tau) \not \in Y_{s+1} \setminus \tau$
\end{itemize}
Given some stage~$s < m$ and some set~$Y_s$,
define the function $f_{\sigma_s} : [Y_s]^{n - |\sigma_s|} \to \omega$ 
by $f_{\sigma_s}(\tau) = f(\sigma_s, \tau)$.
By strong preservation of $k$ hyperimmunities of $\fs^{n - |\sigma_s|}$,
there exists an infinite set $Y_{s+1} \subseteq Y_s$ satisfying (a) and~(b).
We claim that~$Y_m$ satisfies the properties of the lemma.
Fix some~$\sigma \in [F]^t$ and some~$\tau \in [Y_m]^{n - t}$ for some~$t \in (0,n)$.
There is a stage~$s < m$ such that $\sigma = \sigma_s$. Moreover, $\tau \in [Y_{s+1}]^{n - |\sigma_s|}$, 
so by (b), $f(\sigma_s, \tau) \not \in Y_{s+1} \setminus \tau$, therefore~$f(\sigma, \tau) \not \in Y_m \setminus \tau$, 
completing the proof.
\end{proof}

\begin{lemma}\label{lem:fs-cohesivity-strong-preservation}
Suppose that $\ts^t_{d_t+1}$ admits strong preservation of~$k$ hyperimmunities for each $t \in (0,n]$
and $\fs^t$ admits strong preservation of~$k$ hyperimmunities for each $t \in [0,n)$.
For every function $f : [\omega]^n \to \omega$
and every set $C$ preserving hyperimmunity of some sets~$A_0, \dots, A_{k-1}$,
there exists an infinite set $X$ such that $X \oplus C$ preserves hyperimmunity of the~$A$'s 
and for every $\sigma \in [X]^{<\omega}$ such that $0 \leq \card{\sigma} < n$,
$$
(\forall x \in X \setminus \sigma)(\exists b)(\forall \tau \in [X \cap (b, +\infty)]^{n-|\sigma|}) 
  f(\sigma, \tau) \neq x
$$
\end{lemma}
\begin{proof}
Let $X$ be the infinite set constructed in Theorem~\ref{thm:generalized-cohesivity-strong-avoidance}
with $t = n$. For each $s < n$ and $i < d_{n-s}$, let $f_{s,i} : [X]^s \to \omega$ be the function
such that $f_{s,i}(\sigma)$ is the $i$th element of
$\set{x : (\forall b)(\exists \tau \in [X \cap (b, +\infty)]^{n-s}) f(\sigma, \tau) = x}$
if it exists, and 0 otherwise.
Define a finite sequence $X \supseteq X_0 \supseteq \dots \supseteq X_{n-1}$ such that for each $s  < n$
\begin{itemize}
  \item[1.] $X_s \oplus C$ preserves hyperimmunity of the~$A$'s 
  \item[2.] $X_s$ is $f_{s,i}$-free for each $i < d_{n-s}$
\end{itemize}
We claim that $X_{n-1}$ is the desired set. Fix $s < n$ and take any $\sigma \in [X_{n-1}]^s$
and any $x \in X_{n-1} \setminus \sigma$. 
If $(\forall b)(\exists \tau \in [X \cap (b, +\infty)]^{n - s})f(\sigma, \tau) = x$,
then by choice of $X$, there exists an $i < d_{n - s}$ such that
$f_{s, i}(\sigma) = x$, contradicting $f_{s,i}$-freeness of $X_{n-1}$.
So $(\exists b)(\forall \tau \in [X \cap (b, +\infty)]^{n - s})f(\sigma, \tau) \neq x$.
\end{proof}

\begin{proof}[Proof of Theorem~\ref{thm:strong-preservation-fs-left-trapped}]
Fix~$k \geq 2$, some set~$C$, some $C$-hyperimmune sets~$A_0, \dots, A_{k-1}$
and a left trapped coloring~$f : [\omega]^n \to \omega$.
We will construct an infinite $f$-free set~$H$ such that
all the $A$'s is $H \oplus C$-hyperimmune.
Our forcing conditions are Mathias conditions $(F, X)$ such that
\begin{itemize}
  \item[(a)] $X \oplus C$ preserves hyperimmunity of the~$A$'s
  \item[(b)] $(\forall \sigma \in [F \cup X]^n) f(\sigma) \not \in F \setminus \sigma$
  \item[(c)] $(\forall \sigma \in [F \cup X]^t)(\forall x \in (F \cup X) \setminus \sigma)
  (\exists b)(\forall \tau \in [(F \cup X) \cap (b, +\infty)]^{n-t})\\
  f(\sigma, \tau) \neq x$ for each $t \in [0, n)$.
  \item[(d)] $(\forall \sigma \in [F]^t)(\forall \tau \in [X]^{n-t})
  f(\sigma, \tau) \not \in X \setminus \tau$ for each $t \in (0,n)$
\end{itemize}
Properties (c) and (d) will be obtained by 
Lemma~\ref{lem:fs-cohesivity-strong-preservation} 
and Lemma~\ref{lem:fs-left-trapped-preserves-smaller} and are present to maintain the property (b)
over extensions. A set $G$ \emph{satisfies} a condition $(F, X)$ if it is $f$-free and
satisfies the Mathias condition $(F, X)$.
Our initial condition is $(\emptyset, Y)$ where $Y$ is obtained by Lemma~\ref{lem:fs-cohesivity-strong-preservation}.

\begin{lemma}\label{lem:fs2-left-trapped-preserves-1}
For every condition $(F, X)$ there exists an extension
$(H, Y)$ such that $|H| > |F|$.
\end{lemma}
\begin{proof*}
Choose an $x \in X$ such that 
$(\forall \sigma \in [F]^n)f(\sigma) \neq x$ and set $H = F \cup \{x\}$.
By property (c) of $(F, X)$, there exists a $b$
such that
\[
(\forall \sigma \in [F]^{t})
(\forall \tau \in [X \cap (b, +\infty)]^{n-t})f(\sigma, \tau) \neq \{x\} \setminus \sigma
\]
for each $t \in [0, n]$.
By Lemma~\ref{lem:fs-left-trapped-preserves-smaller},
there exists an infinite set $Y \subseteq X \setminus [0, b]$ such that $Y \oplus C$ preserves
hyperimmunity of the~$A$'s and property (d) is satisfied for $(H, Y)$.
We claim that $(H, Y)$ is a valid condition.
Properties (a), (c) and (d) trivially hold. We now check property (b).
By property (b) of $(F, X)$, we only need to check that
$(\forall \sigma \in [F \cup Y]^k)f(\sigma) \neq x$.
This follows from our choice of~$b$.
\end{proof*}

\begin{lemma}\label{lem:fs2-left-trapped-preserves-2}
For every condition $(F, X)$, every~$e \in \omega$ and~$j < k$,
there exists an extension $(H, Y)$
forcing $\Phi_e^{G \oplus C}$ not to dominate $p_{A_j}$. 
\end{lemma}
\begin{proof*}
By removing finitely many elements of $X$, we can assume that
$(\forall \sigma \in [F]^n)f(\sigma) \not \in X$.
Let~$f$ be the partial $X \oplus C$-computable function which on input~$x$
searches for some finite set of integers~$U$ such that
for every left trapped function~$g : [X]^n \to \omega$,
there is a $g$-free set~$E \subseteq X$ such that~$\Phi_e^{(F \cup E) \oplus C}(x) \downarrow \in U$.
If such a set~$U$ is found, $f(x) = 1+max(U)$, otherwise~$f(x) \uparrow$.
We have two cases:
\begin{itemize}
	\item Case 1: the function~$f$ is total. By~$X \oplus C$-hyperimmunity of~$A_j$,
	$f(x) \leq p_{A_j}(x)$ for some~$x$. Let~$U$ be the finite set witnessing~$f(x) \downarrow$.
	By definition of~$f$, taking~$g = f$, there exists a finite $f$-free set~$E$ such that~$a \in W^{(F \cup E) \oplus C}_e$.
	Set~$H = F \cup E$. By property (c) of $(F, X)$, there exists a $b$ such that
	\[
	(\forall \sigma \in [H]^t)(\forall x \in H)
	(\forall \tau \in [X \cap (b, +\infty)]^{n-t}) f(\sigma, \tau) \neq \{x\} \setminus \sigma
	\]
	for each~$t \in [0, n)$.
	By Lemma~\ref{lem:fs-left-trapped-preserves-smaller},
	there exists an infinite set $Y \subseteq X \cap (b, +\infty)$ such that 
	$Y \oplus C$ preserves hyperimmunity of the~$A$'s
	and property (d) is satisfied for $(H, Y)$.
	We claim that $(H, Y)$ is a valid condition.

	Properties (a), (c) and (d) trivially hold. We now check property (b).
	By our choice of $b$, we only need to check that $(\forall \sigma \in [H]^n)
	f(\sigma) \not \in H \setminus \sigma$.
	By property (b) of $(F, X)$, it suffices to check that
	$(\forall \sigma \in [H]^n)f(\sigma) \not \in E \setminus \sigma$.
	By property (d) of $(F, X)$, and our initial assumption on~$X$, we only need to check that
	$(\forall \sigma \in [E]^n)f(\sigma) \not \in E \setminus \sigma$,
	which is exactly $f$-freeness of $E$.

	\item Case 2: there is some~$x$ such that~$f(x) \uparrow$. 
	By compactness, the~$\Pi^{0,X\oplus C}_1$ class~$\Ccal$
	of all left trapped functions~$g : [X]^n \to \omega$
	such that for every~$g$-free set~$E \subseteq X$, $\Phi_e^{(F \cup E) \oplus C}(x) \uparrow$
	is non-empty.
	As~$\wkl$ admits preservation of~$k$ hyperimmunities,
	there exists a left trapped functions~$g \in \Ccal$
	such that $g \oplus X \oplus C$ preserves hyperimmunity of the~$A$'s.
	As~$\fs^n$ admits preservation of~$k$ hyperimmunities,
	there exists some infinite $g$-free set~$Y \subseteq X$ such that~$Y \oplus C$
	preserves hyperimmunity of the~$A$'s. The condition $(F, Y)$ forces $\Phi_e^{G \oplus C}(x) \uparrow$.
\end{itemize}
\end{proof*}

Let~$Y$ be the set constructed in Lemma~\ref{lem:fs-cohesivity-strong-preservation}.
Using Lemma~\ref{lem:fs2-left-trapped-preserves-1}
and Lemma~\ref{lem:fs2-left-trapped-preserves-2}, 
we can define an infinite decreasing sequence of conditions~
$(\emptyset, Y) \geq c_0 \geq \dots$ such that for every~$s \in \omega$
\begin{itemize}
	\item[(i)] $|F_s| \geq s$
	\item[(ii)] $c_s$ forces~$W_e^{G \oplus C} \neq A_j$ if~$s = \tuple{e,j}$
\end{itemize}
where~$c_s = (F_s, X_s)$. Let~$G = \bigcup_s F_s$. By (i), $G$ is infinite
and by~(ii), all the~$A$'s are~$G \oplus C$-hyperimmune.
This completes the proof of Theorem~\ref{thm:strong-preservation-fs-left-trapped}.
\end{proof}

\chapter{The rainbow Ramsey theorem}

Among the consequences of Ramsey's theorem, the rainbow Ramsey theorem
intuitively states the existence of an infinite injective restriction
of any function which is already close to being injective.
We now provide its formal definition.

\index{rainbow Ramsey theorem}
\index{rrt@$\rrt^n_k$|see {rainbow Ramsey theorem}}
\index{bounded function}
\begin{definition}[Rainbow Ramsey theorem]
Fix $n, k \in \N$. A coloring function $f: [\N]^n \to \N$ is \emph{$k$-bounded}
if for every $y \in \N$, $\card{f^{-1}(y)} \leq k$. A set $R$ is a \emph{rainbow} for $f$ (or an \emph{$f$-rainbow})
if $f$ is injective over~$[R]^n$.
$\rrt^n_k$ is the statement ``Every $k$-bounded function $f : [\N]^n \to \N$
has an infinite $f$-rainbow''.
\end{definition}

As usual, $\rrt$ is the statement $(\forall n)(\forall k)\rrt^n_k$.
Galvin noticed that it follows easily from Ramsey's theorem.
Csima and Mileti~\cite{Csima2009strength},
Conidis and Slaman~\cite{Conidis2013Random},
Wang~\cite{WangSome,Wang2013Rainbow} and the author~\cite{Patey2015Somewhere} among others
studied the reverse mathematics of the rainbow Ramsey theorem.

Wang~\cite{WangSome} proved that the rainbow Ramsey theorem is an immediate
consequence of the free set theorem at each level.
More precisely, he showed that $\rca \vdash (\forall n)[\fs^n \imp \rrt^n_2]$.
When not considering the levels of the hierarchies individually, we can obtain a partial reversal.

\begin{theorem}\label{thm:rrt2n-fsn}
For every~$n \geq 1$, $\rca \vdash \rrt^{2n+1}_2 \imp \fs^n$.
\end{theorem}
\begin{proof}
Let~$\tuple{\cdot,\cdot}$ be a bijective coding from~$\{(x, y) \in \Nb^2 : x < y\}$ to~$\Nb$,
such that~$\tuple{x,y} < \tuple{u,v}$ whenever~$x < u$ and~$y < v$. We shall refer to this property as (P1).
We say that a function~$f : [\Nb]^{n} \to \Nb$ is $t$-trapped for some~$t \leq n$
if for every $\vec z \in [\Nb]^n$, $z_{t-1} \leq f(\vec z) < z_t$, 
where~$z_{-1} = -\infty$ and~$z_n = +\infty$.
Wang proved in~\cite[Lemma 4.3]{Wang2014Some} that we can restrict ourselves 
without loss of generality to trapped functions
when $n$ is a standard integer.

Let $f : [\Nb]^{n} \to \Nb$ be a $t$-trapped coloring for some~$t \leq n$. 
We build a $\Delta^{0,f}_1$ 2-bounded coloring $g : [\Nb]^{2n+1} \to \Nb$
such that every infinite $g$-rainbow computes an infinite $f$-thin set.
Given some~$\vec z \in [\Nb]^n$, we write~$\vec z \bowtie_t u$ to denote
the $(2n+1)$-uple
\[
x_0, y_0, \dots, x_{t-1}, y_{t-1}, u, x_t, y_t, \dots, x_{n-1}, y_{n-1}
\]
where~$z_i = \tuple{x_i, y_1}$ for each~$i < n$.
We say that~$\vec z \bowtie_t u$ is \emph{well-formed} if the sequence above is a strictly increasing.

For every $y \in \Nb$ and $\vec{z} \in [\Nb]^n$ such that~$\vec z \bowtie_t y$ is well-formed, 
if $f(\vec z) = \tuple{x,y}$ for some~$x$ such that~$\vec z \bowtie_t x$ is well-formed, 
then set $g(\vec z \bowtie_t y) = g(\vec z \bowtie_t x)$.
Otherwise assign $g(\vec z \bowtie_t y)$ a fresh color.
The function~$g$ is total and 2-bounded.

Let $H = \{ x_0 < y_0 < x_1 < y_1 < \dots \}$ be an infinite $g$-rainbow and let
$H_1 = \{ \tuple{x_i, y_i} : i \in \Nb \}$. We claim that~$H_1$ is $f$-free.
Let~$\vec z \in [H_1]^n$ be such that~$f(\vec z) \in H_1$.
In particular, $f(\vec z) = \tuple{x_i, y_i}$ for some~$i \in \Nb$.
Since~$f$ is $t$-trapped and by (P1), if~$f(\vec z) \neq z_{t-1}$
then~$\vec z \bowtie_t x_i$ and~$\vec z \bowtie_t y_i$ are both well-formed.
Hence~$g(\vec z \bowtie_t x_i) = g(\vec z \bowtie_s y_i)$. Because~$H$ is a $g$-rainbow, either~$x_i$
or~$y_i$ is not in~$H$, contradicting~$\tuple{x_i,y_i} \in H_1$. Therefore~$f(\vec z) = z_{t-1}$.
\end{proof}

\begin{corollary}
$\rrt$ and~$\fs$ coincide over~$\omega$-models.
\end{corollary}

When considering the thin set theorem, we can obtain better bounds.
In particular, the proof is uniform in~$n$, and we can deduce that
the full rainbow Ramsey theorem implies the thin set theorem over~$\rca$.

\begin{theorem}\label{thm:rrt2n-tsn}
$\rca \vdash (\forall n)[\rrt^{n+1}_2 \imp \ts^n]$ and~$\ts^n \leq_{sW} \rrt^{n+1}_2$.
\end{theorem}
\begin{proof}
Fix some~$n \in \Nb$ and let $f : [\Nb]^{n} \to \Nb$ be a coloring. We build a $\Delta^{0,f}_1$
2-bounded coloring $g : [\Nb]^{n+1} \to \Nb$
such that every infinite $g$-rainbow is, up to finite changes, $f$-thin.
For every $y \in \Nb$ and $\vec{z} \in [\Nb]^n$, if $f(\vec z) = \tuple{x,y}$ with $x < y < min(\vec z)$,
then set $g(y, \vec z) = g(x, \vec z)$. Otherwise assign $g(y, \vec z)$ a fresh color.
The function~$g$ is clearly 2-bounded.
Let $H$ be an infinite $g$-rainbow and let~$x, y \in H$ be such that~$x < y$.
Set~$H_1 = H \setminus [0, y]$. We claim that~$H_1$ is $f$-thin with color~$\tuple{x,y}$.
Indeed, for every~$\vec z \in [H_1]^n$, if~$f(\vec z) = \tuple{x,y}$
then $x < y < min(\vec z)$, so~$g(x, \vec z) = g(y, \vec z)$.
This contradicts the fact that~$H$ is a~$g$-rainbow.
\end{proof}

\section{The rainbow Ramsey theorem and randomness}

The thin set, free set and rainbow Ramsey theorem admit the same
definability bounds as the ones proven by Jockusch for Ramsey's theorem~\cite{Jockusch1972Ramseys}. 
In particular, $\Psf^n$ admits a computable instance with no~$\Sigma^0_n$
solution, where~$\Psf^n$ is any of $\rt^n_2$, $\ts^n_2$, $\fs^n$ or~$\rrt^n_2$.
One may wonder whether those problems admit \emph{probabilistic solutions}.

There are many ways to understand the notion of probabilistic solution.
In our case, we shall take the algorithmic randomness approach,
and consider the following computability-theoretic definition.

\index{probabilistic solution}
\begin{definition}
A $\Pi^1_2$ statement~$\Psf$ \emph{admits probabilistic solutions}
if for every $\Psf$-instance~$X$, the measure of oracles~$Z$ such that~$Z \oplus X$
computes a solution to~$X$ is positive.
\end{definition}

\index{n-RAN@$\rans{n}$}
Let~$\rans{n}$ be the statement ``For every set~$X$, there is an $n$-random real relative to~$X$''.
Avigad, Dean and Rute~\cite{Avigad2012Algorithmic} studied the reverse mathematics
of the statement~$\rans{n}$ and proved that~$\rans{1}$ is equivalent to~$\wwkl$ over~$\rca$.
Not to admit probabilistic solutions is sufficient to obtain a separation over standard models.

\begin{lemma}
If~$\Psf$ does not admit probabilistic solutions, then~$\Psf \not \leq_\omega \rans{n}$
for every~$n$.
\end{lemma}
\begin{proof}
Let~$X$ be a $\Psf$-instance such that the measure of oracles $Z$ which $X$-compute a solution to $X$ is null.
In particular, there is an~$n$-random real $R$ relative to~$X$ such that~$R \oplus X$ does not compute a solution to~$X$.
By Van Lambalgen theorem~\cite{VanLambalgen1990axiomatization}, $\rans{n}$ holds in the Turing ideal
$\Ical = \{ Z : (\exists i)[Z \leq_T X \oplus R_0 \oplus \dots \oplus R_i] \}$, where~$R = \bigoplus_i R_i$.
Since every~$C \in \Ical$ is~$R \oplus X$-computable, $\Ical \not \models \Psf$.
\end{proof}

It turns out that very few theorems studied in reverse mathematics admit
probabilistic solutions. One can easily prove that weak statements
such as the atomic model theorem do not admit ones by simply noticing that
the halting set is uniformly almost-everywhere dominating~\cite{Dobrinen2004Almost}.

\begin{theorem}\label{thm:strength-ramsey-amt-nra}
$\amt$ does not admit probabilistic solutions.
\end{theorem}
\begin{proof}
By~\cite{Hirschfeldt2009atomic,Conidis2008Classifying}, $\amt$
is computably equivalent to the statement ``For every $\Delta^0_2$ function~$f$,
there is a function~$g$ not dominated by~$f$''.
Kurtz~\cite{Kurtz1982Randomness} proved the existence of a $\Delta^0_2$ function~$f$ such that for almost
every set~$A$, $f$ dominates every~$A$-computable function.
Let~$T$ be a computable, complete atomic theory such that every atomic model computes 
a function not dominated by~$f$. By choice of~$f$, the measure of oracles which compute
an atomic model of~$T$ is null.
\end{proof}

We shall see in Corollary~\ref{cor:strength-ramsey-srrt32-implies-amt} that the rainbow Ramsey theorem for triples
implies $\amt$ over~$\rca$. Therefore, $\rrt^n_2$ does not admit probabilistic solutions
for~$n \geq 3$. Surprisingly, the case of the rainbow Ramsey theorem for pairs is pretty different.

Csima and Mileti~\cite{Csima2009strength} gave a probabilistic algorithm to compute
solutions to the rainbow Ramsey theorem for pairs.
Conidis and Slaman~\cite{Conidis2013Random} formalized their argument in reverse mathematics
and proved that~$\rrt^2_2$ is a consequence of~$\rans{2}$ over~$\rca$.
We therefore study how randomness and hyperimmunity relate to each other, in order to deduce
preservation of hyperimmunity results about the rainbow Ramsey theorem for pairs
and compare it to the other statements in reverse mathematics.
The following theorem is well-known and has been proven by Agafonov and Levin~\cite{Zvonkin1970complexity}.

\begin{theorem}\label{thm:typicality-randomness-hyperimmunity}
Fix some set~$Z$ and a countable collection of $Z$-hyperimmune sets~$B_0, B_1, \dots$
If~$R$ is sufficiently random relative to~$Z$, the~$B$'s are $R \oplus Z$-hyperimmune.
\end{theorem}
\begin{proof}
It suffices to prove that for every~$e, i \in \omega$, the following class is Lebesgue null.
\[
\Scal_{e,i} = \{X : \Phi_e^{X \oplus Z} \mbox{ dominates } p_{B_i} \}
\]
Suppose it is not the case. There exists~$\sigma \in 2^{<\omega}$ such that
\[
\mu\set{X \in \Scal_{e,i} : \sigma \prec X} > 2^{-|\sigma|-1}
\]
Let~$f$ be the partial $Z$-computable function which
on input~$n$, looks for a clopen set~$[U] \subseteq [\sigma]$
such that $\mu([U]) > 2^{-|\sigma|-1}$ and $\Phi_e^{\tau \oplus Z}(n) \downarrow$ for each~$\tau \in U$,
and outputs~$1+max(\Phi_e^{\tau \oplus Z}(n) : \tau \in U)$. In particular, $[U] \cap \Scal_{e,i} \neq \emptyset$,
so~$f(n) \geq p_{B_i}$.
The function~$f$ is total by choice of~$Z$ and dominates~$B_i$. Contradiction.
\end{proof}

Note that this does not mean that the set~$R$ is hyperimmune-free relative to~$Z$.
In fact, Kautz~\cite{Kautz1991Degrees} proved that the converse holds for~$Z = \emptyset$, that is,
if~$R$ sufficiently random, then it is hyperimmune.

\begin{corollary}
For every~$n$, $\rans{n}$ and~$\rrt^2_2$ admit preservation of hyperimmunity.
\end{corollary}
\begin{proof}
Immediate by Lemma~\ref{thm:typicality-randomness-hyperimmunity} for~$\rans{n}$ and by Csima and Mileti~\cite{Csima2009strength}
for~$\rrt^2_2$.
\end{proof}

In fact, the rainbow Ramsey theorem for pairs admits an exact characterization
in terms of algorithmic randomness. Miller noticed that in the proof of Csima and Mileti~\cite{Csima2009strength}
any function d.n.c.\ relative to~$\emptyset'$ is sufficient to carry the whole construction.
He furthermore obtained a reversal, thereby proving that the rainbow Ramsey theorem for pairs is equivalent
to the diagonal non-computable principle relative to~$\emptyset'$.

\begin{theorem}[Miller~\cite{MillerAssorted}]\label{thm:rrt22-dnrzp-rel}
Fix a set $X$. 
\begin{itemize}
  \item[(i)] There is an $X$-computable 2-bounded coloring $f : [\Nb]^2 \to \Nb$
such that every infinite $f$-rainbow computes (not relative to $X$) a function d.n.c. relative to $X^{'}$.
  \item[(ii)] For every $X$-computable 2-bounded coloring $f : [\Nb]^2 \to \Nb$
  and every function $g$ d.n.c. relative to $X'$, there exists a $g \oplus X$-computable infinite $f$-thin set.
\end{itemize}
\end{theorem}

\index{n-DNR@$\dnrs{n}$}
Following the notation of~\cite{Avigad2012Algorithmic}, 
we denote by~$\dnrs{n}$ the statement ``For every set~$X$, there is a function d.n.c.\ relative to~$X^{(n-1)}$''.
The argument of Miller is easily formalizable in $\rca$, which gives the following corollary.

\begin{corollary}[Miller~\cite{MillerAssorted}]\label{cor:rrt22-dnrzp}
$\rca \vdash \rrt^2_2 \biimp \dnrs{2}$ and $\rrt^2_2 =_u \dnrs{2}$.
\end{corollary}

Kang~\cite{Kang2014Combinatorial} asked whether the rainbow Ramsey theorem for pairs
is a consequence of the thin set theorem for pairs.
We now answer positively, thanks to Miller's characterization of~$\rrt^2_2$.
Note that the implication is strict, since $\amt$ and therefore~$\sts^2$
does not admit probabilistic solutions.

\begin{theorem}\label{thm:ts2-dnrzp}
$\rca \vdash \ts^2 \imp \dnrs{2}$ and~$\dnrs{2} \leq_{sW} \ts^2$.
\end{theorem}
\begin{proof}
We prove that for every set $X$, there is an $X$-computable coloring function $f : [\Nb]^2 \to \Nb$
such that every infinite $f$-thin set computes (not relative to $X$) a function d.n.c. relative to $X^{'}$. 
The structure of the proof is very similar to Theorem~\ref{thm:em-dnrzp}, but instead of diagonalizing against
computing an infinite transitive tournament, we will diagonalize against computing an infinite set avoiding color $i$. 
Applying diagonalization for each color $i$, we will obtain the desired result. 

Let $X$ be a set and $g(.,.)$ be a total $\Delta^{0,X}_1$ function such that $\Phi_e^{X'}(e)=\lim_s g(e,s)$ if the limit exists, and $\Phi^{X'}_e(e) \uparrow$ if the limit does not exist. For each~$e, i, s \in \Nb$,
interpret $g(e,s)$ as the code of a finite set $D_{e,i,s}$ of size $3^{e\cdot i}$. 
We define the coloring~$f$ by $\Sigma_1$-induction as follows. Set $f_0=\emptyset$. At stage~$s+1$, do the following. Start with $f_{s+1}=f_s$. Then, for each $\alpha(e,i) <s$ -- where $\alpha(.,.)$ is the Cantor pairing function, i.e., 
$\alpha(e,i) = \frac{(e+i)(e+i+1)}{2}+e$ -- take the first element $x \in D_{e,i,s} \setminus \bigcup_{(e',i') < (e,i)} D_{e',i',s}$ (notice that these exist by cardinality assumptions on the $D_{e,i,s}$), and if $f_{s+1}(s,x)$ is not already assigned, assign it to color $i$. Finally, for any $z<s$ such that $f_{s+1}(s,z)$ remains undefined, assign any color to it in a predefined way (e.g., for any such pair $\{x,y\}$, set $f_{s+1}(x,y)$ to be $0$). This finishes the construction of $f_{s+1}$. Set $f=\bigcup_s f_s$, which must exist as a set by $\Sigma_1$-induction. 

First of all, notice that $f$ is a coloring function of domain $[\N]^2$, as at the end of stage~$s+1$ of the construction $f(x,y)$ is assigned a value for (at least) all pairs $\{x,y\}$ with $x<s$ and $y<s$. By $\ts^2$, let~$A$ be an infinite $f$-thin set.
Let $i \in \Nb \setminus f([A]^2)$. Let $h(e)$ be the code of the finite set $A_e$ consisting of the first $3^{e\cdot i+1}$ elements of~$A$. We claim that $h(e) \not=\Phi_e^{X'}(e)$ for all~$e$, which would prove $\dnrzp$. Suppose otherwise, i.e., suppose that $\Phi_e^{X'}(e)=h(e)$ for some~$e$. Then there is a stage $s_0$ such that $h(e)=g(e,s)$ for all~$s \geq s_0$ or equivalently $D_{e,i,s}=A_e$ for all~$s \geq s_0$. Let $N_e=\max(A_e)$. The same argument as in the proof of Theorem~\ref{thm:em-dnrzp} shows that for any $s$ be bigger than both $\max(\bigcup_{e,i,s < N_e} D_{e,i,s})$ and $s_0$, the restriction of $f$ to $A_e \cup \{s\}$ does not avoid color $i$, which contradicts the fact that the infinite set $A$ containing $A_e$ avoids color $i$ in~$f$.
\end{proof}

\begin{corollary}
$\rca \vdash \ts^2 \imp \rrt^2_2$ and $\rrt^2_2 \leq_W \ts^2$.
\end{corollary}
\begin{proof}
Immediate by Theorem~\ref{thm:ts2-dnrzp} and Corollary~\ref{cor:rrt22-dnrzp}.
\end{proof}

\section{A stable rainbow Ramsey theorem for pairs}

In this section, we study a stable version of the rainbow Ramsey theorem.
There exist different notions of stability for $k$-bounded functions.

Consider a 2-bounded coloring $f$ of pairs as the history of interactions
between people in an infinite population. $f(x, s) = f(y, s)$ means that $x$ and $y$
interact at time $s$. In this world, $x$ and $y$ get \emph{married} if $f(x, s) = f(y, s)$
for cofinitely many $s$, whereas a person $x$ becomes a \emph{monk} if $f(x, s)$ is a fresh color
for cofinitely many $s$. Finally, a person $x$ is \emph{wise} if for each~$y$, 
either $x$ and~$y$ get married or $x$ and~$y$ eventually
break up forever, i.e., $(\forall y)[(\forall^\infty s) f(x,s) = f(y,s) \vee (\forall^\infty s)f(x,s) \neq f(y,s)]$. 
In particular married people and monks are wise.
Note that 2-boundedness implies that a person~$x$ can get married to at most one~$y$.

$\rrt^2_2$ states that given an world, we can find infinitely many instants where
people behave like monks. However we can weaken our requirement, leading to new principles.

\index{rainbow-stable}
\index{srrt@$\srrt^n_k$}
\index{stable rainbow Ramsey theorem}
\begin{definition}[Stable rainbow Ramsey theorem]\label{def:stable-rainbow-ramsey-theorem}
A coloring $f : [\Nb]^2 \to \Nb$ is 
\emph{rainbow-stable} if for every $x$, one of the following holds:
\begin{itemize}
	\item[(a)] There is a $y \neq x$ such that 
	$(\forall^\infty s)f(x, s) = f(y,s)$
	\item[(b)] $(\forall^{\infty} s) \card{\set{y \neq x : f(x, s) = f(y, s)}} = 0$
\end{itemize}
$\srrt^2_2$ is the statement ``every rainbow-stable 2-bounded coloring $f:[\Nb]^2 \to \Nb$
has a rainbow.''
\end{definition}

Hence in the restricted world of $\srrt^2_2$, everybody either gets married or becomes a monk.
$\srrt^2_2$ is a particular case of $\rrt^2_2$.

The naturality of this version of stability is justified by the existence of various
simple characterizations of the stable rainbow Ramsey theorem for pairs.
We shall later study another version which seems more natural in the sense 
that a stable instance can be obtained from a non-stable one by an application
of the cohesiveness principle. However the latter version does not admit
immediate simple characterizations.

\index{strongly rainbow-stable}
\index{prerainbow}
\begin{definition}
A 2-bounded coloring $f : [\Nb]^2 \to \Nb$ is 
\emph{strongly rainbow-stable} if  
$(\forall x)(\exists y \neq x)(\forall^\infty s)f(x, s) = f(y,s)$
A set $X \subseteq \Nb$ is a \emph{prerainbow} for a 2-bounded coloring $f : [\Nb]^2 \to \Nb$ if 
$(\forall x \in X)(\forall y \in X)(\forall^\infty s \in X)[f(x, s) \neq f(y, s)]$.
\end{definition}

\begin{lemma}[Wang in \cite{WangSome}, $\rca + \bst$]\label{lem:prerainbow-equiv}
Let $f : [\Nb]^2 \to \Nb$ be a $2$-bounded coloring and $X$ be an infinite $f$-prerainbow. 
Then $X \oplus f$ computes an infinite $f$-rainbow $Y \subseteq X$.
\end{lemma}

\begin{theorem}\label{thm:rainbow-stable-vs-strongly}
The following are equivalent over $\rca + \bst$:
\begin{itemize}
  \item[(i)] $\srrt^2_2$
  \item[(ii)] Every strongly rainbow-stable 2-bounded coloring $f : [\Nb]^2 \to \Nb$ has a rainbow.
\end{itemize}
\end{theorem}
\begin{proof}
$(i) \imp (ii)$ is straightforward as any strongly rainbow-stable coloring is rainbow-stable.
$(ii) \imp (i)$: Let $f : [\Nb]^2 \to \Nb$ be a 2-bounded rainbow-stable coloring.
Consider the following collection:
$$
S = \set{ x \in \Nb : (\forall^{\infty} s)(\forall y \neq x)[f(y, s) \neq f(x,s)]}
$$

If $S$ is finite, then take $n \geq max(S)$. The restriction of $f$ to $[n, +\infty)$
is a strongly rainbow-stable 2-bounded coloring and we are done. So suppose $S$ is infinite.
We build a 2-bounded strongly rainbow-stable  coloring $g \leq_T f$ by stages.

At stage $t$, assume $g(x, i)$ is defined for every $x, i < t$.
For every pair $x, y \leq t$ such that $f(x, t) = f(y, t)$,
define $g(x, t) = g(y, t)$. Let $S_t$ be the set of $x \leq t$
such that $g(x, t)$ has not been defined yet.
Writing $S_t = \set{x_1 < x_2 < \dots}$, we set $g(x_{2i}) = g(x_{2i+1})$ for each $i$.
If $S_t$ has an odd number of elements, there remains an undefined value.
Set it to a fresh color.
This finishes the construction. 
It is clear by construction that $g$ is 2-bounded.

\begin{claim}
$g$ is strongly rainbow-stable.
\end{claim}
\begin{proof*}
Fix any $x \in \Nb$. Because $f$ is rainbow-stable, we have two cases:
\begin{itemize}
  \item Case 1: there is a $y \neq x$ such that $(\forall^{\infty} s) f(x, s) = f(y, s)$.
  Let $s_0$ be the threshold such that $(\forall s \geq s_0) f(x, s) = f(y, s)$.
  Then by construction, at any stage $s \geq s_0$, $g(x, s) = g(y, s)$ and we are done.
  
  \item Case 2: $x \in S$. Because $x$ is infinite, it has a successor $y_0 \in S$. 
  By $\bst$, let $s_0$ be the threshold such that for every $y \leq y_0$
  either there is a $z \leq y_0$, $z \neq y$ such that $(\forall s \geq s_0) f(y, s) = f(z, s)$ 
  or $(\forall s \geq s_0)$ $f(y, s)$ is a fresh color.
  Then by construction of $g$, for every $t \geq s_0$, $S_t \uh y = S \uh y$.
  Either $x = x_{2i}$ for some $i$ and then $(\forall t \geq s_0) g(x, t) = g(x_{2i+1}, t)$
  or $x = x_{2i+1}$ for some $i$ and then $(\forall t \geq s_0) g(x, t) = g(x_{2i}, t)$.
\end{itemize}
\end{proof*}

\begin{claim}
Every infinite $g$-prerainbow is an $f$-prerainbow.
\end{claim}
\begin{proof*}
Let $X$ be an infinite $g$-prerainbow and assume for the sake of contradiction
that it is not an $f$-prerainbow. Then there exists two elements $x, y \in X$
such that $(\forall s)(\exists t \geq s)[f(x, t) = f(y, t)]$.
But then because $f$ is rainbow-stable, 
there is a threshold $s_0$ such that $(\forall s \geq s_0)[f(x, s) = f(y, s)]$.

Then by construction of $g$, for every $s \geq s_0$, $g(x, s) = g(y, s)$.
For every $u \in X$ there is an $s \in X$ with $s \geq u, s_0$ such that
$g(x, s) = g(y, s)$ contradicting the fact that $X$ is a $g$-prerainbow.
\end{proof*}

Using Lemma~\ref{lem:prerainbow-equiv}, for any infinite $H$ $g$-prerainbow,
$f \oplus H$ computes an infinite $f$-rainbow. This finishes the proof.
\end{proof}

\bigskip

\subsection{Relation with diagonal non-recursiveness}\label{sect:srrt22-dnr}

It is well-known that being able to compute a d.n.c.\ function is equivalent
to being able to uniformly find a member outside a finite $\Sigma^0_1$ set if we know an upper bound on its size,
and also equivalent to diagonalize against a $\Sigma^0_1$ function. The proof relativizes well 
and is elementary enough to be formalized in $\rca$ (see Theorem~\ref{thm:dnrzn-char-dnr}).

\index{diagonalizing function}
\index{escaping function}
\begin{definition}
Let $(X_e)_{e \in \Nb}$ be a uniform family of finite sets.
An \emph{$(X_e)_{e \in \Nb}$-escaping function} is a function $f : \Nb^2 \to \Nb$
such that $(\forall e)(\forall n)[\card{X_e} \leq n \imp f(e,n) \not \in X_e]$.
Let $h : \Nb \to \Nb$ be a function. An \emph{$h$-diagonalizing function} $f$ is a function $\Nb \to \Nb$
such that $(\forall x)[f(x) \neq h(x)]$.
When $(X_e)_{e \in \Nb}$ and $h$ are clear from context, they may be omitted.
\end{definition}

\begin{theorem}[Folklore]\label{thm:dnrzn-char-dnr}
For every $n \geq 1$, the following are equivalent over $\rca$:
\begin{itemize}
  \item[(i)] $\dnrs{n}$
  \item[(ii)] Any uniform family $(X_e)_{e \in \Nb}$ of $\Sigma^0_n$ finite sets
  has an escaping function.
  \item[(iii)] Any partial $\Delta^0_n$ function has a diagonalizing function.
\end{itemize}
\end{theorem}
\begin{proof}
Fix a set $A$ and some~$n \geq 1$.
\begin{itemize}
  \item $(i) \imp (ii)$:
  Let $(X_e)_{e \in \Nb}$ be a uniform family of finite $\Sigma^{0,A}_{n+1}$ finite sets
  and $f$ be a function d.n.c.\ relative to $A^{(n-1)}$. Define a function $h : \Nb^2 \to \Nb$
  by $h(e, s) = \tuple{f(i_1), \dots, f(i_s)}$ where $i_j$ is the index of the
  partial $\Delta^{0, A}_n$ function which on every input, looks at the $j$th element $k$ of $X_e$
  if it exists, interprets $k$ as an $s$-tuple $\tuple{k_1, \dots, k_s}$ and returns $k_j$.
  The function diverges if no such $k$ exists.
  One easily checks that $h$ is an $(X_e)_{e \in \Nb}$-escaping function.
  
  \item $(ii) \imp (iii)$:
  Let $f : \Nb \to \Nb$ be a partial $\Delta^{0,A}_n$ function.
  Consider the enumeration defined by $X_e = \set{f(e)}$ if it $f(e) \downarrow$ and $X_e = \emptyset$ otherwise.
  This is a uniform family of $\Sigma^{0,A}_n$ finite sets, each of size at most 1. Let $g: \Nb^2 \to \Nb$ be an
  $(X_e)_{e \in \Nb}$-escaping function. Then $h : \Nb \to \Nb$ defined by $h(e) = g(e, 1)$ is
  an $f$-diagonalizing function.  
  
  \item $(iii) \imp (i)$:
  Consider the partial $\Delta^{0,A}_n$ function $f(e) = \Phi^{A^{(n-1)}}_e(e)$. 
  Any $f$-diagonalizing function is d.n.c relative to $A^{(n-1)}$.
\end{itemize}
\end{proof}

In particular, using Miller's characterization of $\rrt^2_2$ by $\dnrzp$, we have the following theorem taking $n=2$:

\begin{theorem}[Folklore]\label{thm:rrt22-char-dnr}
The following are equivalent over $\rca$:
\begin{itemize}
  \item[(i)] $\rrt^2_2$
  \item[(ii)] Any uniform family $(X_e)_{e \in \Nb}$ of $\Sigma^0_2$ finite sets
  has an escaping function.
  \item[(iii)] Any partial $\Delta^0_2$ function has a diagonalizing function.
\end{itemize}
\end{theorem}

In the rest of this section, we will give an equivalent of Theorem~\ref{thm:rrt22-char-dnr}
for $\srrt^2_2$.

\begin{lemma}[$\rca + \bst$]\label{lem:srrt22-char1}
For every $\Delta^0_2$ function $h : \Nb \to \Nb$,
there exists a computable rainbow-stable 2-bounded coloring $c : [\N]^2 \to \Nb$
such that every infinite rainbow $R$ for $c$ computes an $h$-diagonalizing function.
\end{lemma}
\begin{proof}
Fix a $\Delta^0_2$ function $h$ and a uniform family $(D_e)_{e \in \Nb}$ of all finite sets.
We will construct a rainbow-stable 2-bounded coloring $c : [\Nb]^2 \to \Nb$ by
a finite injury priority argument. By Shoenfield's limit lemma, there exists
a total computable function $g(\cdot, \cdot)$ such that $\lim_s g(x, s) = h(x)$ for every $x$.

Our requirements are the following:

\bigskip
$\Rcal_x$: If $\card{D_{\lim_s g(x, s)}} \geq 3x+2$ then $\exists u, v \in D_{\lim_s g(x, s)}$
such that $(\forall^{\infty} s) c(v, s) = c(v, s)$.
\smallskip

We first check that if every requirement is satisfied then we can compute
a function $f : \Nb \to \Nb$ such that $(\forall x)[f(x) \neq h(x)]$ from
any infinite $c$-rainbow. Fix any infinite $c$-rainbow $R$.
Let $f$ be the function which given $x$ returns the index of the set
of the first $3x+2$ elements of $R$. Because of the requirement $\Rcal_x$,
$D_{f(x)} \neq D_{\lim_s g(x,s)}$. Otherwise $\card{D_{f(x)}} = 3x+2$ and
there would be two elements $u, v \in D_{f(x)} \subset R$
such that $(\forall^{\infty} s) c(x, s) = c(y, s)$. So take an element $s \in R$ large enough to witness this fact.
$c(x, s) = c(y, s)$ for $x, y, s \in R$ contradicting the fact that $R$ is a rainbow.
So $D_{f(x)} \neq D_{\lim_s g(x,s)}$ from which we deduce $f(x) \neq \lim_s g(x, s) = h(x)$. 

Our strategy for satisfying a local requirement $\Rcal_x$ is as follows. 
If $\Rcal_x$ receives attention at stage $t$, it checks whether $\card{D_{g(x, t)}} \geq 3x+2$.
If this is not the case, then it is declared satisfied. If $\card{D_{g(x, t)}} \geq 3x+2$,
then it chooses the least two elements $u, v \geq x$, such that $u, v \in D_{g(x, s)}$ 
and $u$ and $v$ are \emph{not} restrained
by a strategy of higher priority and \emph{commits} to assigning a common color.
For any such pair $u, v$, this commitment will remain active as long as the strategy has a restraint
on that element. Having done all this, the local strategy is declared to be satisfied and will not act again
unless either a higher priority puts restraint on $u$ or $v$ or at a further stage $t' > t$, $g(x, t') \neq g(x, t)$.
In both cases, the strategy gets \emph{injured} and has to reset, releasing all its restraints.

To finish stage $t$, the global strategy assigns $c(u, t)$ for all $u \leq t$ as follows:
if $u$ is committed to some assignment of $c(u, t)$ due to a local strategy, define $c(u, t)$ to be this value.
If not, let $c(u, t)$ be a fresh color. This finishes the construction and we now turn to the verification.
It is easy to check that each requirement restrains at most two elements at a given stage.

\begin{claim}
Every given strategy acts finitely often.
\end{claim}
\begin{proof*}
Fix some~$x \in \Nb$.
By~$\bst$ and because $g$ is limit-computable, there exists a stage $s_0$
such that $g(y, s) = g(y, s_0)$ for every~$y \leq x$ and~$s \geq s_0$.
If $|D_{g(x,s_0)}| < 3x+2$, then the requirement is satisfied and does not act any more.
If $|D_{g(x,s_0)}| \geq 3x+2$, then by a cardinality argument, there exists two elements~$u$ and~$v$ in $D_{g(x,s_0)}$
which are not restrained by a strategy of higher priority. Because $D_{g(y, s)} = D_{g(y, s_0)}$ for each~$y \leq x$
and~$s \geq s_0$, no strategy of higher priority will change its restrains and will therefore injure~$\Rcal_x$ after stage~$s_0$.
So $(\forall^{\infty} s) c(u, s) = c(v,s)$ for some $u, v \in D_{\lim_s g(x, s)}$ and requirement $\Rcal_x$ is satisfied.
\end{proof*}

\begin{claim}
The resulting coloring $c$ is rainbow-stable.
\end{claim}
\begin{proof*}
Consider a given element $u \in \Nb$. We distinguish three cases:
\begin{itemize}
  \item Case 1: the element becomes, during the construction, free from any restraint after some stage
$t \geq t_0$ . In this case, by construction, $c(u, t)$ is assigned a fresh color for all $t \geq t_0$.
Then $(\forall^{\infty} s)(\forall v \neq u)[c(u, s) \neq c(v, s)]$.

  \item Case 2: there is a stage $t_0$ at which some restraint is put on $u$ by some local strategy, and this
restraint is never released. In this case, the restraint comes together with a commitment that
all values of $c(u, s)$ and $c(v, s)$ be the same beyond some stage $t_0$ for some fixed $v \neq x$. 
Therefore for all but finitely many stages $s$, $c(u, s) = c(v, s)$.

  \item Case 3: during the construction, infinitely many restraints are put on $u$ and are later released. This
  is actually an impossible case, since by construction only strategies for requirements $\Rcal_y$
  with $y \leq u$ can ever put a restraint on $u$. By~$\bst$,
	there exists some stage after which no strategy ~$\Rcal_y$ acts for every~$y \leq u$
	and therefore the restraints on~$u$ never change again.
\end{itemize}
\end{proof*}
This last claim finishes the proof.
\end{proof}

\begin{lemma}[$\rca + \ist$] \label{lem:char1-srrt22}
For every computable strongly rainbow-stable 2-bounded coloring $f : [\Nb]^2 \to \Nb$
there exists a uniform family $(X_e)_{e \in \Nb}$ of $\Delta^0_2$ finite sets
whose sizes are uniformly $\Delta^0_2$ computable
such that every $(X_e)_{e \in \Nb}$-escaping function
computes a $c$-rainbow.
\end{lemma}
\begin{proof}
Fix any uniform family $(D_e)_{e \in \Nb}$ of finite sets.
Let $f : [\Nb]^2 \to \Nb$ be a 2-bounded rainbow-stable computable coloring.
For an element $x$, define 
$$
\bad(x) = \set{ y \in \Nb : (\forall^\infty s)f(x,s) = f(y,s)}
$$
Notice that $x \in \bad(x)$.
Because $f$ is strongly rainbow-stable, $\bad$ is a $\Delta^0_2$ function.
For a set $S$, $\bad(S) = \bigcup_{x \in S} \bad(x)$.
Define $X_e = \bad(D_e)$. Hence $X_e$ is a $\Delta^0_2$ set, and this uniformly in $e$.
Moreover, $\card{X_e} \leq 2\card{D_e}$ and for every $x$, $\card{\bad(x)} = 2$
so we can $\emptyset'$-compute the size of $X_e$ with the following equality
$$
  \card{X_e} = 2|D_e| - 2 \card{\set{\set{x, y} \subset D_e : \bad(x) = \bad(y)}}
$$

Let $h : \Nb \to \Nb$ be a function satisfying $(\forall e)(\forall n)[\card{X_e} \leq n \imp h(e, n) \not \in X_e]$.
We can define $g : \Nb \to \Nb$ by $g(e) = h(e, 2\card{D_e})$.
Hence $(\forall e)g(e) \not \in X_e$.

We construct a prerainbow $R$ by stages $R_0 (=\emptyset) \subsetneq R_1 \subsetneq R_2, \dots$
Assume that at stage $s$, 
$(\forall \{x,y\} \subseteq R_s)(\forall^{\infty} s)[f(x, s)\neq f(y,s)]$.
Because $R_s$ is finite, we can computably find some index $e$ such that $R_s = D_e$.
Set $R_{s+1} = R_s \cup \set{g(e)}$. By definition, $g(e) \not \in X_e$.
Let $x \in R_s$. Because $g(e) \not \in X_e$, $(\forall^{\infty} s) f(x, s) \neq f(g(e), s)$.
By~$\ist$, the set~$R$ is an $f$-prerainbow.
By Lemma~\ref{lem:prerainbow-equiv} we can compute an infinite $f$-rainbow from $R \oplus f$.
\end{proof}

\begin{theorem}\label{thm:srrt22-characterizations}
The following are equivalent over $\rca + \ist$:
\begin{itemize}
  \item[(i)] $\srrt^2_2$
  \item[(ii)] Any uniform family $(X_e)_{e \in \Nb}$ of $\Sigma^0_2$ finite sets
  whose sizes are uniformly $\Delta^0_2$ has an escaping function.
  \item[(iii)] Any $\Delta^0_2$ function $h : \Nb \to \Nb$ has a diagonalizing function.
\end{itemize}
\end{theorem}
\begin{proof}
$(i) \imp (iii)$ is Lemma~\ref{lem:srrt22-char1} and $(ii) \imp (i)$ follows from Lemma~\ref{lem:char1-srrt22}.
This is where we use~$\ist$.
We now prove $(iii) \imp (ii)$.
Let $(X_e)_{e \in \Nb}$ be a uniform family of $\Sigma^0_2$ finite sets
such that~$|X_e|$ is~$\Delta^0_2$ uniformly in~$e$. For each~$n, i \in \Nb$,
define~$(n)_i$ to be the $i$th component of the tuple whose code is~$n$, if it exists. Define
\[
h(\tuple{e,i}) = \cond{
	(n)_i & \mbox{ where } n \mbox{ is the } i\mbox{th element of } X_e \mbox{ if } i < |X_e|\\
	0 & \mbox{ otherwise}
}
\]
By (iii), let $g : \Nb \to \Nb$ be a total function such that $(\forall e)[g(e) \neq h(e)]$.
Hence for every pair $\tuple{e,i}$ such that~$i \leq |X_e|$, $g(\tuple{e,i}) \neq (n)_i$
where $n$ is the $i$th element of $X_e$.
Define $f : \Nb^2 \to \Nb$ to return on inputs $e$ and $s$ the tuple $\tuple{g(\tuple{e,0}), \dots, g(\tuple{e,s})}$.
Hence if $s \geq \card{X_e}$ then $f(e, s) \neq m$
where $m$ is the $i$th element of $X_e$ for each~$i < |X_e|$. So $f(e,n) \not \in X_e$.
\end{proof}

\begin{corollary}\label{cor:srrt22-omega-model-dnr}
Every $\omega$-model of $\srrt^2_2$ is a model of $\dnr$.
\end{corollary}
\begin{proof}
Let $h : \Nb \to \Nb$ be the $\Delta^0_2$ function which on input $e$ returns $\Phi_e(e)$ if $\Phi_e(e) \downarrow$
and returns 0 otherwise. By (iii) of Theorem~\ref{thm:srrt22-characterizations} there
exists a total function $f : \Nb \to \Nb$ such that $(\forall e)[f(e) \neq h(e)]$.
Hence $(\forall e)[f(e)\neq \Phi_e(e)]$ so $f$ is a d.n.c.\ function. 
\end{proof}

We strengthen Corollary~\ref{cor:srrt22-omega-model-dnr} and prove that in fact $\srrt^2_2$ implies $\dnr$ over~$\rca$.

\begin{theorem}
$\rca \vdash \srrt^2_2 \imp \dnr$
\end{theorem}
\begin{proof}
If $\Phi_{e}(e) \downarrow$ then interpret $\Phi_e(e)$ as the code of a finite set $D_e$ of size $3^{e+1}$
with $min(D_e) > e$.
Let $D_{e,s}$ be the approximation of $D_e$ at stage $s$, i.e. $D_{e, s}$ is the set $\{e+1, \dots, e+3^{e+1}\}$
if $\Phi_{e,s}(e) \uparrow$ and $D_{e,s} = D_e$ if $\Phi_{e,s}(e) \downarrow$.
We will construct a rainbow-stable coloring $f : [\Nb]^2 \to \Nb$ meeting the following requirements for each $e \in \Nb$.
$$
\Rcal_e : \Phi_e(e) \downarrow \imp (\exists a, b \in D_e)(\forall^{\infty} s)f(a, s) = f(b, s)
$$

Before giving the construction, let us explain how to compute a d.n.c.\ function from any infinite $f$-rainbow
if each requirement is satisfied. Let $H$ be an infinite $f$-rainbow. Define the function
$g : \Nb \to \Nb$ which given $e$ returns the code of the $3^{e+1}$ first elements of $H$.
We claim that $g$ is a d.n.c.\ function. Otherwise suppose $g(e) = \Phi_e(e)$ for some $e$.
Then $D_e \subseteq H$, but by $\Rcal_e$, $(\exists a, b \in D_e)(\forall^{\infty} s)f(a, s) = f(b, s)$.
As $H$ is infinite, there exists an $s \in H$ such that $f(a, s) = f(b, s)$, contradicting the 
fact that $H$ is an $f$-rainbow.

We now describe the construction. The coloring $f$ is defined by stages.
Suppose that at stage $s$, $f(u,v)$ is defined for each $u, v < s$.
For each $e < s$ take the first pair $\set{a, b} \in D_{e,s} \setminus \bigcup_{k < e} D_{k,s}$.
Such a pair must exist by cardinality assumption on the $D_{e,s}$. Set $f(a, s) = f(b, s) = i$
for some fresh color $i$. Having done that, for any $u$ not yet assigned, assign $f(u, s)$ a fresh color
and go to stage $s+1$.

\begin{claim}
Each requirement $\Rcal_e$ is satisfied.
\end{claim}
\begin{proof*}
Fix an $e \in \Nb$.
By $\bsig^0_1$ there exists a stage $s$ such that $\Phi_{k,s}(k) = \Phi_k(k)$ for each $k \leq e$.
Then at each further stage $t$, the same par $\set{a, b}$ will be chosen in $D_{e,s}$
to set $f(a, t) = f(b, t)$. Hence if $\Phi_e(e) \downarrow$, there are $a, b \in D_e$ such that
$(\forall^{\infty} s)f(a, s) = f(b, s)$.
\end{proof*}

\begin{claim}
The coloring $f$ is rainbow-stable.
\end{claim}
\begin{proof*}
Fix an element $u \in \Nb$. By $\bsig^0_1$ there is a stage $s$ such that 
$\Phi_{k,s}(k) = \Phi_k(k)$ for each $k < u$. If $u \in \set{a, b}$
for some pair $\set{a,b}$ chosen by a requirement of priority $k < u$ then
at any further stage $t$, $f(u, t) = f(a, t) = f(b,t)$. If $u$ is not chosen by any requirement
of priority $k < u$, then $u$ will not be chosen by any further requirement as $min(D_e) > e$ for each $e \in \Nb$.
So by construction, $f(u,t)$ will be given a fresh color for each $t > s$.
\end{proof*}
\end{proof}

\subsection{K\"onig's lemma and relativized Schnorr tests}\label{srrt22-konig}

D.n.c.\ degrees admit other characterizations in terms of Martin-L\"of tests
and Ramsey-Type K\"onig's lemmas. For the former, it is well-known
that d.n.c.\ degrees are the degrees of infinite subsets of Martin-L\"of randoms \cite{Kjos-Hanssen2009Infinite,Greenberg2009Lowness}.
The latter has been introduced by Flood in~\cite{Flood2012Reverse} under the name $\rkl$ and and renamed into $\rwkl$
in~\cite{Bienvenu2015logical}. It informally states the existence 
of an infinite subset of $P$ or $\overline{P}$ where $P$ is a path through a tree.
We shall 

\index{homogeneous!for a tree}
\begin{definition}
Fix a binary tree $T \subseteq 2^{<\Nb}$ and a $c \in \set{0,1}$.
A string $\sigma \in 2^{<\Nb}$ is \emph{homogeneous for a tree $T$ with color $c$} if there exists a $\tau \in T$
such that $\forall i < \card{\sigma}$, $\sigma(i) = 1 \imp \tau(i) = c$. A set $H$ is \emph{homogeneous for
$T$} if there is a $c \in \set{0,1}$ such that for every initial segment $\sigma$ of $H$, $\sigma$
is homogeneous for $T$ with color $c$.
\end{definition}

Flood proved in~\cite{Flood2012Reverse} the existence of a computable tree of positive
measure such that every infinite  homogeneous set computes a d.n.c.\ function.
We shall study this the Ramsey-type K\"onig's lemma in chapter~\ref{chap:ramsey-type-konig-lemma}.

\index{homogeneous!for a Martin-L\"of test}
\begin{definition}
A \emph{Martin-L\"of test} relative to $X$ is a sequence $(U_i)_{i \in \Nb}$ 
of uniformly $\Sigma^{0,X}_1$ classes such that $\mu(U_n) \leq 2^{-n}$ for all $n$.
A set $H$ is \emph{homogeneous} for a Martin-L\"of test $(U_i)_{i \in \Nb}$
if there exists an $i$ such that $H$ is homogeneous for the tree
corresponding to the closed set~$\overline{U_i}$.
\end{definition}

\begin{theorem}[Flood \cite{Flood2012Reverse}, Bienvenu et al.~\cite{Bienvenu2015logical}]\label{thm:dnrzp-rwwkl}
For every $n \in \Nb$, the following are equivalent over $\rca+\isig^0_{n+1}$:
\begin{itemize}
  \item[(i)] $\dnr[0^{(n)}]$
  \item[(ii)] Every Martin-L\"of test $(U_i)_{i \in \Nb}$ relative to $\emptyset^{(n)}$
  has an infinite homogeneous set.
  \item[(iii)] Every $\Delta^0_{n+1}$ tree of positive measure has an infinite homogeneous set.
\end{itemize}
\end{theorem}

In the rest of this section, we will prove an equivalent theorem
for $\srrt^2_2$.

\index{Schnorr test}
\begin{definition}[Schnorr~\cite{Schnorr1971Zufalligkeit}]
A Martin-L\"of test $(U_n)_{n \in \Nb}$ relative to $X$ 
is a \emph{Schnorr test} relative to $X$ if the measures $\mu(U_n)$
are uniformly $X$-computable.
\end{definition}

\begin{lemma}[$\rca+\bst$]
For every set $A$, every $n \in \Nb$ and every function $f \leq_T A'$
there exists a tree $T \leq_T A'$ such that $\mu(T)$ is an $A'$-computable positive real, $\mu(T) \geq 1 - \frac{1}{2^n}$
and every infinite set homogeneous for $T$ computes a function $g$
such that $g(e) \neq f(e)$ for every~$e$.

Moreover the index for $T$ and for its measure can be found effectively from $n$
and $f$.
\end{lemma}
\begin{proof}
Fix $n \in \Nb$. Let $(D_{e,i})_{e, i \in \Nb}$ be an enumeration of finite sets such that
\begin{itemize}
  \item[(i)] $min(D_{e,i}) \geq i$
  \item[(ii)] $\card{D_{e,i}} = i+2+n$
  \item[(iii)] given an $i$ and finite set $U$ satisfying (i) and (ii),
  one can effectively find an $e$ such that $D_{e, i} = U$.
\end{itemize}

For any canonical index $e$ of a finite set, define $T_e$ to be the downward closure
of the $f$-computable set $\set{\sigma \in \str :  \exists a, b \in D_{f(e), e} : \sigma(a) = 0 \wedge \sigma(b) = 1}$.
The set~$T_e$ exists by~$\bsig^{0,f}_1$, hence~$\bst$.
Define also $T_{\leq e} = \bigcap_{i=0}^e T_e$. It is easy to see that
$$
  \mu(T_e) = 1 - \frac{1}{2^{\card{D_{f(e), e}}-1}}
$$

Fix a $\emptyset'$-computable function $f$. Consider the following tree $T = \bigcap_{i = 0}^{\infty} T_i$.
Because of condition (ii),
$$
\mu(T) \geq 1 - \sum_{i = 0}^{\infty} [1 - \mu(T_i)] 
  = 1 - \sum_{i = 0}^{\infty} \frac{1}{2^{i+1+n}} = 1 - \frac{1}{2^{n}}
$$

\begin{claim}
$T$ is an $f$-computable tree.
\end{claim}
\begin{proof*}
Fix a string $\sigma \in \str$.
$\sigma \in T$ iff $\sigma \in \bigcap_{i=0}^{\infty} T_i$
By definition, $\sigma \in T_i$ iff $\sigma \preceq \tau$
for some $\tau \in \str$ such that there are some elements $a, b \in D_{f(i), i}$
verifying $\tau(a) = 0$ and $\tau(b) = 1$.
When $i \geq \card{\sigma}$, because of conditions (i) and (ii)
there exists $a, b \geq i$ with $a, b \in D_{f(i), i}$
and $\tau \succeq \sigma$ such that $\tau(a) = 0$ and $\tau(b) = 1$.
Hence
$\sigma \in T$ iff $\sigma \in T_{\leq \card{\sigma}}$,
which is an $f$-computable predicate uniformly in $\sigma$.
\end{proof*}

\begin{claim}
$\mu(T)$ is an $f$-computable real.
\end{claim}
\begin{proof*}
Fix any $c \in \Nb$. For any $d \in \Nb$, by condition (ii)
$$
\mu(T_{\leq d}) \geq \mu(T) \geq \mu(T_{\leq d}) - \sum_{i=d}^{\infty} \frac{1}{2^{i+1+n}}
$$
In particular, for $d$ such that $2^{-n} - \sum_{i=0}^{d} \frac{1}{2^{i+1+n}} \leq 2^{-c}$ we have
$$
\card{\mu(T_{\leq d}) - \mu(T)} \leq \sum_{i=d}^{\infty} \frac{1}{2^{i+1+n}} \leq \frac{1}{2^c}
$$
It suffices to notice that $\mu(T_{\leq d})$ is easily $f$-computable
as for $u = max(\bigcup_{i=0}^d D_{f(i), i})$
$$
\mu(T_{\leq d}) = \frac{\card{\set{\sigma \in 2^u : \sigma \in T_{\leq d}}}}{2^u}
$$
\end{proof*}

Let $H$ be an infinite set homogeneous for $T$.

\begin{claim}
$H$ computes a function $g$ such that $g(i) \neq f(i)$ for every $i$.
\end{claim}
\begin{proof*}
Let $g$ be the $H$-computable function which on input $i$
returns an $e \in \Nb$ such that $D_{e, i}$ is the set
of the first $i+2+n$ elements of $H$. Such an element can be effectively found
by condition (iii).

Assume for the sake of contradiction that $g(i) = f(i)$ for some $i$.
Then by definition of being homogeneous for $T$, there exists a $j \in \{0,1\}$ and a $\sigma \in T$
such that $\sigma(u) = j$ whenever $u \in H$. In particular, $\sigma \in T_i$.
So there exists $a, b \in D_{f(i), i} = D_{g(i), i} \subset H$ such that
$\sigma(a) = 0$ and $\sigma(b) = 1$. Hence there exists an $a \in H$ such that $\sigma(a) \neq j$. 
Contradiction.
\end{proof*}

This last claim finishes the proof.
\end{proof}

\begin{corollary}
For every 2-bounded, computable coloring $f : [\Nb]^2 \to \Nb$
there exists a $\emptyset'$-computable tree $T$ of positive $\emptyset'$-computable measure
such that every infinite set homogeneous for $T$
computes an infinite $f$-rainbow.
\end{corollary}

\begin{corollary}\label{cor:schnorr-srrt22}
For every 2-bounded, computable coloring $f : [\Nb]^2 \to \Nb$
there exists a Schnorr test $(U_i)_{i \in \Nb}$ relative to $\emptyset'$ such that
every infinite set homogeneous for $(U_i)_{i \in \Nb}$ computes
an infinite $f$-rainbow.
\end{corollary}

\begin{theorem}[$\rca+\ist$]\label{thm:konig-srrt}
Fix a set $X$. For every $X'$-computable tree $T$ of positive $X'$-computable measure $\mu(T)$
there exists a uniform family $(X_e)_{e \in \Nb}$ of $\Delta^{0, X}_2$ finite sets
whose sizes are uniformly $X'$-computable
and such that every $(X_e)_{e \in \Nb}$-escaping function
computes an infinite set homogeneous for $T$.
\end{theorem}
\begin{proof}
Consider $X$ to be computable for the sake of simplicity.
Relativization is straightforward.
We denote by $(D_e)_{e \in \Nb}$ the canonical enumeration of all finite sets.
Let $T$ be a $\emptyset'$-computable tree of positive $\emptyset'$-computable measure $\mu(T)$.
For each~$s \in \Nb$, let~$T_s$ be the set of strings~$\sigma \in 2^{<\Nb}$ of length~$s$
and let~$\mu_s(T)$ be the first $s$ bits approximation of~$\mu(T)$.
Consider the following set for each finite set $H \subseteq \Nb$ and $k \in \Nb$.
$$
\bad( H, k) = \left\{ n \in \Nb : \mu_{4k}(T \cap \Gamma^0_H \cap \Gamma^0_n) < 2^{-2k}\right\}
$$

First notice that the measure of $T \cap \Gamma^0_H$ (resp. $T \cap \Gamma^0_H \cap \Gamma^0_n$)
is $\emptyset'$-computable  uniformly in $H$ (resp. in $H$ and $n$), 
so one $\bad(H, k)$ is uniformly $\Delta^0_2$. We now prove that $\bad(H,k)$
has a uniform $\Delta^0_2$ upper bound, which is sufficient to deduce that~$|\bad(H, k)|$
is uniformly~$\Delta^0_2$.

Given an $H$ and a $k$, let $\epsilon = 2^{-k-1} - 2^{-2k} - 2^{-4k}$.
We can $\emptyset'$-computably find a length $s = s(H, k)$ such that
$$
\frac{|T_s \cap \Gamma^0_H|}{2^s} - \mu(T \cap \Gamma^0_H)  < \epsilon
$$

\begin{claim}
If $2^{-k} \leq \mu(T \cap \Gamma^0_H)$, then $max(\bad(H,k)) \leq s$
\end{claim}
\begin{proof*}
Fix any $n > s$. By choice of $s$,
$$
\mu(T \cap \Gamma^0_H \cap \Gamma^1_n) \leq \frac{|T_s \cap \Gamma^0_H|}{2^{s+1}}
	\leq \frac{\mu(T \cap \Gamma^0_H)}{2} + \epsilon
$$
Furthermore,
$$
\mu(T \cap \Gamma^0_H \cap \Gamma^0_n) = \mu(T \cap \Gamma^0_H) - \mu(T \cap \Gamma^0_H \cap \Gamma^1_n)
$$
Putting the two together, we obtain
\begin{eqnarray*}
\mu(T \cap \Gamma^0_H \cap \Gamma^0_n) &\geq& \mu(T \cap \Gamma^0_H) - \frac{\mu(T \cap \Gamma^0_H)}{2} - \epsilon\\
	&\geq& \frac{\mu(T \cap \Gamma^0_H)}{2} - \epsilon \geq 2^{-k-1} - \epsilon \geq 2^{-2k} + 2^{-4k}
\end{eqnarray*}
In particular
$$
\mu_{4k}(T \cap \Gamma^0_H \cap \Gamma^0_n) \geq \mu(T \cap \Gamma^0_H \cap \Gamma^0_n) - 2^{-4k} \geq 2^{-2k}
$$
Therefore $n \not \in \bad(H,k)$.
\end{proof*}

For each~$H$ and~$k$, let~$X_{H,k} = \bad(H, k) \cap [0, s(H, k)]$.
The set $X_{H, k}$ is ~$\Delta^0_2$ uniformly in~$H$ and~$k$, and its size is uniformly $\Delta^0_2$.
In addition, by previous claim, if $2^{-k} \leq \mu(T \cap \Gamma^0_H)$ then $\bad(H,k) \subseteq X_{H,k}$.

Let $g : \Pcal_{fin}(\Nb) \times \Nb \times \Nb \to \Nb$ be a total function such that for every finite set $H$
and $k \in \Nb$, $g(H, k, n) \not \in X_{H,k}$ whenever $n \geq \card{X_{H, k}}$. Fix any $k \in \Nb$ 
such that $2^{-k} \leq \mu(T)$.
We construct by $\isig^{0,g}_1$ a set $H$ and a sequence of integers $k_0, k_1, \dots$ by finite approximation as follows.
First let $H_0 = \emptyset$ and~$k_0 = k$. We will ensure during the construction that for all $s$:
\begin{itemize}
  \item[(a)] $\card{H_s} = s$
  \item[(b)] $T \cap \Gamma^0_{H_s}$ has measure at least $2^{-k_s}$
  \item[(c)] $H_s \subseteq H_{s+1}$ and every $n \in H_{s+1} \setminus H_s$ is greater than all elements in $H_s$.
\end{itemize}
Suppose $H_s$ has been defined already. The tree $T \cap \Gamma^0_{H_s}$ has measure at least $2^{-k_s}$
and $\card{\bad(H_s, k_s)}$ has at most $2k_s$ elements. Thus $g(H_s, k_s) \not \in X_{H_s,k_s} \supseteq \bad(H_s, k_s)$.
We set $H_{s+1} = H_s \cup \{g(e, k_s)\}$ and~$k_{s+1}$ be the least integer such that~$2^{-k_{s+1}} \leq 2^{-2k_{s}} - 2^{-4k_s}$.
By definition of~$\bad(H_s, k_s)$, $T \cap \Gamma^0_{H_{s+1}}$ has measure at least
$2^{-2k_s}$ with an approximation of~$2^{-4k_s}$, so has measure at least $2^{-k_{s+1}}$.

Let now $H = \bigcup_s H_s$.

\begin{claim}
$H$ is homogeneous for $T$.
\end{claim}
\begin{proof*}
Suppose for the sake of contradiction that $H$ is not homogeneous for $T$.
This means that there are only finitely many $\sigma \in T$ such that $H$ is homogeneous for $\sigma$.
Therefore for some level $l$,
$\set{\sigma \in T_l \; \mid \; \forall i \in H \ \sigma(i) = 0}=\emptyset$.
Since $H \cap \{0,..,l\} = H_l \cap \{0,..,l\}$, we in fact have 
$\set{\sigma \in T_l \; \mid \; \forall i \in H_l \ \sigma(i) = 0}=\emptyset$.
 
In other words, $T \cap \Gamma^0_{H_l} = \emptyset$ which contradicts property (b) in the definition of $H_l$
ensuring that $T \cap \Gamma^0_{H_l}$ has measure at least $2^{-k_l}$.
Thus $H$ is homogeneous for $T$.
\end{proof*}
\end{proof}

\begin{theorem}
The following are equivalent over $\rca + \ist$:
\begin{itemize}
  \item[(i)] $\srrt^2_2$
  \item[(ii)] Every Schnorr test $(U_i)_{i \in \Nb}$ relative to $\emptyset'$
  has an infinite homogeneous set.
  \item[(iii)] Every $\Delta^0_2$ tree of $\emptyset'$-computable positive measure 
  has an infinite homogeneous set.
\end{itemize}
\end{theorem}
\begin{proof}
$(i) \imp (iii)$ is Theorem~\ref{thm:konig-srrt}
together with Theorem~\ref{thm:srrt22-characterizations}.
$(iii) \imp (ii)$ is obvious
and $(ii) \imp (i)$ is Corollary~\ref{cor:schnorr-srrt22}. 
\end{proof}

Hirschfeldt et al.\ proved in \cite[Theorem 3.1]{Hirschfeldt2008Limit} that for
every $X'$-computable martingale $M$, there is a set low over $X$ on which $M$ does not succeed.
Schnorr proved in \cite{Schnorr1971Zufalligkeit} that for every Schnorr test $C$ relative to $X'$
there exists an $X'$-computable martingale $M$ such that a set does not succeeds on $M$ iff it passes the test $C$.
By Corollary~\ref{cor:schnorr-srrt22}, there exists an $\omega$-model of $\srrt^2_2$ containing only low sets.
However we will prove it more directly under the form of a low basis theorem for $\emptyset'$-computable
trees of $\emptyset'$-computable positive measure. This is an adaptation of \cite[Proposition 2.1]{Barmpalias2012Randomness}.

\begin{theorem}[Low basis theorem for $\Delta^0_2$ trees]\label{thm:low-tree-exact-measure}
Fix a set $X$. Every $X'$-computable tree of $X'$-computable positive measure
has an infinite path $P$ low over $X$ (i.e., such that $(X \oplus P)' \leq_T X'$).
\end{theorem}
\begin{proof}
Fix $T$, an $X'$-computable tree of 
$X'$-computable positive measure $\mu(T)$. We will define
an $X'$-computable subtree $U$ of measure $\frac{\mu(T)}{2}$
such that any infinite path through $T$ is GL${}_1$ over~$X$. It then suffices
to take any $\Delta^{0,X}_2$ path through $U$ to obtain the desired path low over~$X$.

Let $f$ be an $X'$-computable function that on input $e$
returns a stage $s$ after which $e$ goes into $A'$ for at most
measure $2^{-e-2}\mu(T)$ of oracles~$A$. Given $e$ and $s = f(e)$, the oracles $A$
such that $e$ goes into $A'$ after stage $s$ form a $\Sigma^{0,X}_1$ class $V_e$
of measure $\mu(V_e) \leq 2^{-e-2}\mu(T)$.
Thus~$\mu(\bigcap_e \overline{V_e}) \geq 1 - \sum_e 2^{-e-2}\mu(T) \geq 1 - \frac{\mu(T)}{2}$.
Therefore~$\mu(T \cap \bigcap_e \overline{V_e}) \geq \frac{\mu(T)}{2}$.
One can easily restrict $T$ to a subtree $U$ such that $[U] \subseteq \bigcap_e \overline{V_e}$
and $\mu(U) = \frac{\mu(T)}{2}$.
For any path $P \in [U]$ and any $e \in \Nb$, $e \in P' \biimp e \in P'_{f(e)}$.
Hence $P$ is GL${}_1$ over~$X$.
\end{proof}

\begin{corollary}
There exists an $\omega$-model of $\srrt^2_2$ containing only low sets.
\end{corollary}

\begin{corollary}
There exists an $\omega$-model of $\srrt^2_2$ which is neither a model of $\semo$ nor of $\sts^2$.
\end{corollary}
\begin{proof}
If every computable stable tournament had a low infinite subtournament then we could build an $\omega$-model $M$
of $\semo + \sads$ having only low sets, but then $M \models \srt^2_2$ contradicting \cite{Downey200102}.
Moreover, by Theorem~\ref{thm:sts-no-low} any $\omega$-model of $\sts^2$ contains a non-low set.
\end{proof}

In fact we will see later that even $\rrt^2_2$ implies neither $\semo$ nor $\sts^2$ on $\omega$-models.

\subsection{Relations to other principles}

We now relate the stable rainbow Ramsey theorem for pairs
to other existing principles studied in reverse mathematics.
This provides in particular a factorization of existing
implications proofs. For example, both the rainbow Ramsey theorem for pairs
and the stable Erd\H{o}s-Moser theorem are known to imply the omitting partial types principle ($\opt$) over~$\rca$.
In this section, we show that both principles imply~$\srrt^2_2$,
which itself implies~$\opt$ over~$\rca$.
Hirschfeldt \& Shore in~\cite{Hirschfeldt2007Combinatorial}
introduced~$\opt$ and proved its equivalence with~$\hyp$ over~$\rca$.

\begin{theorem}\label{thm:srrt22-imp-opt}
$\rca \vdash \srrt^2_2 \imp \hyp$
\end{theorem}
\begin{proof}[Proof using Cisma \& Mileti construction, $\rca$]
We prove that the construction from Csima \& Mileti in~\cite{Csima2009strength} that $\rca \vdash \rrt^2_2 \imp \hyp$
produces a rainbow-stable coloring. We take the notations and definitions
of the proof of Theorem~4.1 in~\cite{Csima2009strength}. It is therefore essential
to have read it to understand what follows.
Fix an $x \in \Nb$. By $\bsig^0_1$ there exists an $e \in \Nb$ and a stage $t$ after which
$n^k_j$ and $m^k$ will remains stable for any $k \leq e$ and any $j \in \Nb$
and such that $n^e_i \leq x < n^e_{i+1}$ for some $i$.
\begin{itemize}
  \item If $i > 0$ then 
$x$ will be part of no pair $(m,l)$ for any requirement and $f(x, s) = \tuple{x,s}$
will be fresh for cofinitely many $s$.

  \item If $i = 0$ and $n^e_j$ is defined for each $j$
such that $j+1 \leq \frac{(n^e_0 - m^e)^2 - (n^e_0 - m^e)}{2}$
then as there are finitely many such $j$, after some finite stage
$x$ will not be paired any more and $f(x, s) = \tuple{x,s}$
will be fresh for cofinitely many $s$.

  \item If $i = 0$ and $n^e_j$ is undefined for some $j$
  such that $\tuple{m, x} = j+1$ or $\tuple{x, m} = j+1$ for some $m$, then
  $x$ will be part of a pair $(m, l)$ for cofinitely many $s$
  and so there exists an $m$ such that $f(x,s) = f(m, s)$ for cofinitely many $s$.
  
  \item If $i = 0$ and $n^e_j$ is undefined for some $j$
  such that $\tuple{m, x} \neq j+1$ or $\tuple{x, m} \neq j+1$ for any $m$
  then $x$ will not be paired after some stage
  and $f(x, s) = \tuple{x,s}$ will be fresh for cofinitely many $s$.
\end{itemize}
In any case, either $f(x, s)$ is fresh for cofinitely many $s$,
or there is a $y$ such that $f(x, s) = f(y, s)$ for cofinitely many $s$.
So the coloring is rainbow-stable.
\end{proof}

\smallskip

We can also adapt the proof using $\Pi^0_1$-genericity to $\srrt^2_2$.

\begin{proof}[Proof using $\Pi^0_1$-genericity, $\rca + \ist$]
Take any incomplete $\Delta^0_2$ set $P$ of PA degree.
The author proved in~\cite{Patey2015Degrees} the existence of
a $\Delta^0_2$ function $f$ such that $P$ does not compute any $f$-diagonalizing function.

Fix any functional $\Psi$.
Consider the $\Sigma^0_2$ class
$$
U = \set{X \in \cs : (\exists e) \Psi^X(e) \uparrow \vee \Psi^X(e) = f(e)}
$$

Consider any $\Pi^0_1$-generic $X$ such that $\Psi^X$ is total.
Either there exists a $X \in U$ in which case $\Psi^X(e) = f(e)$ hence $\Psi^X$ is not an $f$-diagonalizing function.
Or there exists a $\Pi^0_1$ class $F$ disjoint from $U$ and containing $X$. Any member of $F$ computes an $f$-diagonalizing
function. In particular $P$ computes an $f$-diagonalizing function. Contradiction.
\end{proof}

\begin{corollary}
$\rca \vdash \srrt^2_2 \imp \opt$
\end{corollary}

The following theorem is not surprising as by a relativization of Theorem~\ref{thm:srrt22-imp-opt}
to $\emptyset'$, there exists an $\emptyset'$-computable rainbow-stable coloring of pairs
such that any infinite rainbow computes a function hyperimmune relative to $\emptyset'$.
Csima et al.~\cite{Csima2004Bounding} and Conidis~\cite{Conidis2008Classifying} proved
that $\amt$ is equivalent over $\omega$-models to the statement ``For any $\Delta^0_2$ function $f$, there exists a function $g$
not dominated by $f$''. Hence any $\omega$-model of $\srrt^2_2[\emptyset']$ is an $\omega$-model of $\amt$.
We will prove that the implication holds over $\rca$.

\begin{theorem}\label{thm:srrt2n-stsn}
$\rca \vdash (\forall n)[\srrt^{n+1}_2 \imp \sts^n]$
\end{theorem}
\begin{proof}
Fix some~$n \in \Nb$ and let $f : [\Nb]^{n} \to \Nb$ be a stable coloring. 
If~$n = 1$, then~$f$ has a $\Delta^{0,f}_1$ infinite thin set, so suppose~$n > 1$.
We build a $\Delta^{0,f}_1$
rainbow-stable 2-bounded coloring $g : [\Nb]^{n+1} \to \Nb$
such that every infinite $g$-rainbow is, up to finite changes, $f$-thin.
Construct~$g$ as in the proof of Theorem~\ref{thm:rrt2n-tsn}.
It suffices to check that~$g$ is rainbow-stable whenever $f$ is stable.

Fix some~$x \in \Nb$ and~$\vec z \in [\Nb]^{n-1}$ such that~$x < min(\vec z)$. 
As~$f$ is stable, there exists a stage~$s_0 > max(\vec z)$ after which $f(\vec z, s) = f(\vec z, s_0)$.
Interpret~$f(\vec z, s_0)$ as a tuple~$\tuple{u,v}$.
If~$u \geq v$ or~$v \geq min(\vec z)$ or~$x \not \in \{u,v\}$, then $g(x, \vec z, s)$ will be given
a fresh color for every~$s \geq s_0$.
If~$u < v < min(\vec z)$ and~$x \in \{u,v\}$ (say $x = u$),
then~$g(x, \vec z, s) = g(v, \vec z, s)$ for every~$s \geq v$.
Therefore $g$ is rainbow-stable.
\end{proof}

\begin{corollary}\label{cor:strength-ramsey-srrt32-implies-amt}
$\rca \vdash \srrt^3_2 \imp \amt$
\end{corollary}

\begin{theorem}[$\rca + \bst$]
For every $\Delta^0_2$ function $f$, 
there exists a computable stable coloring $c:[\Nb]^2 \to \Nb$ such that every infinite
$c$-thin set computes an $f$-diagonalizing function.
\end{theorem}
\begin{proof}
Fix a $\Delta^0_2$ function $f$ as stated above.
For any $n \in \Nb$, fix a canonical enumeration $(D_{n, e})_{e \in \Nb}$ of all
finite sets of $n+1$ integers greater than $n$.
We will build a computable stable coloring $c : [\Nb]^2 \to \Nb$ fulfilling 
the following requirements for each $e, i \in \Nb$:

\smallskip
$\Rcal_{e, i}: $ $\exists u \in D_{\tuple{e,i}, f(e)}$ such that $(\forall^{\infty} s)c(u, s) = i$.
\smallskip

We first check that if every requirement is satisfied, then any
infinite $c$-thin set computes an $f$-diagonalizing function.
Let $H$ be an infinite $c$-thin set for color $i$.
Define $h: \Nb \to \Nb$ to be the $H$-computable function which on $e$
returns the value $v$ such that $D_{\tuple{e,i}, v}$
is the set of the $\tuple{e,i}+1$ first elements of $H$ greater than
$\tuple{e,i}$.

\begin{claim}
$h$ is an $f$-diagonalizing function.
\end{claim}
\begin{proof*}
Suppose for the sake of contradiction that $h(e) = f(e)$ for some $e$.
Then $D_{\tuple{e,i}, h(e)} = D_{\tuple{e,i}, f(e)}$.
But by $\Rcal_{e,i}$, $\exists u \in D_{\tuple{e,i}, f(e)}$ such that
$(\forall^{\infty} s)c(u,s) = i$. Then there is an $s \in H$ such that
$c(u, s) = i$, and as $D_{\tuple{e,i}, f(e)} \cup \set{s} \subset H$,
$H$ is not $c$-thin for color~$i$. Contradiction.
\end{proof*}

By Shoenfield's limit lemma, let $g(\cdot, \cdot)$ be the partial approximations of~$f$.
The strategy for satisfying a local requirement $\Rcal_{e,i}$ is as follows.
At stage $s$, it takes the least element $u$ of $D_{\tuple{e,i}, g(x,s)}$ not restrained by a strategy of
higher priority if it exists. Then it puts a restraint on $u$ and \emph{commits}
$u$ to assigning color $i$. For any such $u$, this commitment 
will remain active as long as the strategy has a restraint on that element.
Having done all this, the local strategy is declared to be satisfied and will not act again, unless either a higher priority
puts a restraint on $u$, or releases a $v \in D_{\tuple{e,i}, g(e,s)}$ with $v < u$ or 
at a further stage $t > s$, $g(e,t) \neq g(e, s)$.
In each case, the strategy gets \emph{injured} and has to reset, releasing its restraint.

To finish stage $s$, the global strategy assigns $c(u, s)$ for all $u \leq s$ as follows:
if $u$ is commited to some assignment of $c(u, s)$ due to a local strategy, define $c(u, s)$ to be this value.
If not, let $c(u, t) = 0$. This finishes the construction and we now turn to the verification.
It is easy to check that each requirement restrains at most one element at a given stage.

\begin{claim}
Each strategy $\Rcal_{e,i}$ acts finitely often.
\end{claim}
\begin{proof*}
Fix some strategy~$\Rcal_{e,i}$.
By $\bst$, there is a stage~$s_0$ after which $g(x,s) = f(x)$ for every~$x \leq \tuple{e,i}$.
Each strategy restrains at most one element, 
and the strategies of higher priority will always choose the same elements after stage~$s_0$.
As $\card{D_{\tuple{e,i}, f(e)}} = \tuple{e,i}+1$,
the set of $u \in D_{\tuple{e,i}, f(e)}$ such that no strategy of higher priority puts a restraint on $u$ is non empty
and does not change.
Let $u_{min}$ be its minimal element. By construction, $\Rcal_{e,i}$ will choose $u_{min}$ before stage $s_0$
and will not be injured again.
\end{proof*}

\begin{claim}
The resulting coloring $c$ is stable.
\end{claim}
\begin{proof*}
Fix a $u \in \Nb$.
If $\tuple{e, i} > u$ then $\Rcal_{e,i}$ does not put a restraint on $u$ at any stage.
As each strategy acts finitely often, by $\bst$ there exists a stage $s_0$ after which
no strategy $\Rcal_{e, i}$ with $\tuple{e,i} \leq u$ will act on $u$.
There are two cases: In the first case, at stage $s_0$ the element $u$ is restrained
by some strategy $\Rcal_{e,i}$ with $\tuple{e,i} \leq u$ in which case 
$c(u, s)$ will be assigned a unique color specified by strategy $\Rcal_{e,i}$
for cofinitely many $s$. In the other case, after stage $s_0$, the element $u$ is free from
any restraint, and $c(u, s) = 0$ for cofinitely many $s$.
\end{proof*}
\end{proof}

\begin{corollary}
$\rca + \ist \vdash \sts^2 \imp \srrt^2_2$
\end{corollary}

\begin{theorem}[$\rca$]
For every rainbow-stable 2-bounded coloring $f : [\Nb]^2 \to \Nb$,
there exists an $f$-computable stable tournament $T$ such that every 
infinite transitive subtournament of $T$ computes an $f$-rainbow.
\end{theorem}
\begin{proof}
Use exactly the same construction as in Theorem~3.1 in~\cite{Kang2014Combinatorial}.
We will prove that in case of rainbow-stable colorings, the constructed tournament $T$ is stable.
Fix an $x \in \Nb$. By rainbow-stability, either $f(x, s)$ is a fresh color for cofinitely many $s$,
in which case $T(x, s)$ holds for cofinitely many $s$, or there exists a $y$ such that
$f(y, s) = f(x, s)$ for cofinitely many $s$. If $T(x, y)$ holds then
$T(x, s)$ does not hold and $T(y, s)$ holds for cofinitely many $s$.
Otherwise $T(x, s)$ holds and $T(y, s)$ does not hold for cofinitely many $s$.
Hence $T$ is stable.\end{proof}

\begin{corollary}
$\rca \vdash \semo \imp \srrt^2_2$
\end{corollary}

\begin{question}
Does $\srrt^2_2 + \coh$ imply $\rrt^2_2$ over $\rca$ ?
\end{question}

\section{A weakly-stable rainbow Ramsey theorem for pairs}

Despite the robustness of the stable rainbow Ramsey theorem for pairs
which has been shown to admit several simple characterizations,
rainbow-stability does not seem to be the natural stability notion corresponding
to~$\rrt^2_2$. In particular, it is unknown
whether~$\rca \vdash \coh + \srrt^2_2 \imp \rrt^2_2$.

Wang used in \cite{Wang2014Cohesive} another version of stability for rainbow Ramsey theorems
to prove various results, like the existence of non-PA solution to any instance of $\rrt^3_2$.
This notion leads to a principle between $\rrt^2_2$ and $\srrt^2_2$.

\index{weakly rainbow-stable}
\index{wsrrt@$\wsrrt^n_k$}
\begin{definition}[Weakly stable rainbow Ramsey theorem]
A coloring $f : [\Nb]^2 \to \Nb$ is \emph{weakly rainbow-stable} if
$$
(\forall x)(\forall y)[(\forall^{\infty} s)f(x, s) = f(y,s) \vee (\forall^{\infty} s)f(x, s) \neq f(y,s)]
$$
$\wsrrt^2_2$ is the statement ``every weakly rainbow-stable 2-bounded coloring $f:[\Nb]^2 \to \Nb$
has an infinite rainbow.''
\end{definition}

Weak rainbow-stability can be considered as the ``right'' notion of stability for 2-bounded colorings
as one can extract an infinite weakly rainbow-stable restriction of any 2-bounded coloring using
cohesiveness.
However the exact strength of $\wsrrt^2_2$ is harder to tackle. A characterization
candidate would be computing an infinite subset of a path in a $\emptyset'$-computably
graded $\Delta^0_2$ tree where the notion of computable gradation is taken from the restriction
of Martin-L\"of tests to capture \emph{computably random} reals. 

In this section, we study the weakly-rainbow stable rainbow Ramsey theorem for pairs.
We prove that it is enough be able to escape finite
$\Delta^0_2$ sets to prove $\wsrrt^2_2$. We also separate $\wsrrt^2_2$
from~$\rrt^2_2$ by proving that $\wsrrt^2_2$ contains an $\omega$-model with only low sets.
The question of exact characterizations of $\wsrrt^2_2$ remains open.

It is easy to see that every rainbow-stable coloring is weakly rainbow-stable,
hence $\rca \vdash \wsrrt^2_2 \imp \srrt^2_2$.
Wang proved in \cite[Lemma 4.11]{Wang2014Cohesive} that $\rca \vdash \coh + \wsrrt^n_2 \imp \rrt^n_2$
and that $\wsrrt^2_2[\emptyset']$ has an $\omega$-model
with only low${}_2$ sets. We show through the following theorem that~$\wsrrt^{n+1}_2$ corresponds to the exact strength
of~$\rrt^n_2[\emptyset']$ for every~$n$.

\begin{theorem}\label{thm:wsrrt-jump-of-rrt}
For every standard $n \geq 1$, $\rca + \bst \vdash \wsrrt^{n+1}_2 \biimp \rrt^n_2[\emptyset']$\\
and $\wsrrt^{n+1}_2 =_c \rrt^n_2[\emptyset']$.
\end{theorem}
\begin{proof}
We first prove that $\rca \vdash \wsrrt^{n+1}_2 \imp \rrt^n_2[\emptyset']$.
Fix a $X'$-computable 2-bounded coloring $f : [\Nb]^n \to \Nb$.
Using Shoenfield's limit lemma, there exists an $X$-computable approximation function $h: [\Nb]^{n+1} \to \Nb$
such that $\lim_s h(\vec{x}, s) = f(\vec{x})$ for every $\vec{x} \in [\Nb]^n$.
Let~$\tuple{ \dots }$ be a standard coding of the lists of integers into~$\Nb$
and~$\prec_\Nb$ be a computable total order over~$\Nb^{<\Nb}$.
We define an $X$-computable 2-bounded coloring $g:[\Nb]^{n+1}$ as follows.
$$
g(\vec{x}, s) = \cond{
   \tuple{h(\vec{x}, s), s, 0} & 
   \mbox{ if there is at most one } \vec{y} \prec_\Nb \vec{x} \mbox{ s.t. } h(\vec{y},s) = h(\vec{x}, s)\\
   \tuple{rank_{\prec_\Nb}(\vec{x}), s, 1} & \mbox{ otherwise}
}
$$
(where $rank_{\prec_\Nb}(\vec{x})$ is the position of $\vec{x}$
for any well-order $\prec_\Nb$ over tuples). 
By construction $g$ is 2-bounded and $X$-computable. Also notice that~$g$ is rainbow-stable.

\begin{claim}
Every infinite $g$-rainbow is an $f$-rainbow.
\end{claim}
\begin{proof*}
Let $A$ be an infinite $g$-rainbow. Assume for the sake of contradiction that $\vec{x}, \vec{y} \in [A]^n$
are such that $\vec{y} \prec_\Nb \vec{x}$ and $f(\vec{y}) = f(\vec{x})$.
Fix $t \in \Nb$ such that $h(\vec{z}, s) = f(\vec{z})$ whenever $\vec{z} \preceq_\Nb \vec{x}$
and $s \geq t$. Fix $s$ such that $s \in A$, $s \geq t$ and $s > max(\vec{x})$. Notice that 
since $f$ is 2-bounded and $h(\vec{z}, s) = f(\vec{z})$ for every $\vec{z} \preceq_\Nb \vec{x}$,
we have $g(\vec{z}, s) = \tuple{h(\vec{z}, s), s, 0} = \tuple{f(\vec{z}), s, 0}$ for every $\vec{z} \preceq_\Nb \vec{x}$.
Hence
$$
g(\vec{x}, s) = \tuple{f(\vec{x}), s, 0} = \tuple{f(\vec{y}), s, 0} = g(\vec{y}, s)
$$
contradicting the fact that $A$ is a $g$-rainbow.
\end{proof*}

We now prove that $\rca +\bst \vdash \rrt^n_2[\emptyset'] \imp \wsrrt^{n+1}_2$.
Let $f : [\Nb]^{n+1} \to \Nb$ be a 2-bounded weakly rainbow-stable coloring.
Let $g : [\Nb]^n \to \Nb$ be the 2-bounded coloring which on $\vec{x} \in [\Nb]^n$
will fetch the least $\vec{y} \preceq_n \vec{x}$ such that $(\forall^{\infty} s)f(\vec{x}, s) = f(\vec{y}, s)$
and return color $\tuple{\vec{y}}$. One easily sees that $g$ is $f'$-computable and 2-bounded.
By $\rrt^n_2[\emptyset']$, let $H$ be an infinite $g$-rainbow. 

We claim that $H$ is an $f$-prerainbow.
Suppose for the sake of contradiction that there exists $\vec{x} \preceq_n \vec{y} \in H$
such that $(\forall^{\infty} s)f(\vec{x}, s) = f(\vec{y}, s)$. Then by definition
$g(\vec{x}) = g(\vec{y}) = \tuple{\vec{x}}$ and $H$ is not a $g$-rainbow.
By Lemma~\ref{lem:prerainbow-equiv} and~$\bst$, $f \oplus H$ computes an infinite $f$-rainbow.
\end{proof}

\begin{lemma}[$\rca+\ist$]\label{lem:char1-wsrrt22}
For every computable weakly rainbow-stable 2-bounded coloring $f : [\Nb]^2 \to \Nb$
there exists a uniformly $\Delta^0_2$ family $(X_e)_{e \in \Nb}$ of finite sets
such that every $(X_e)_{e \in \Nb}$-escaping function computes an infinite $f$-rainbow.
\end{lemma}
\begin{proof}
Fix any uniformly $\Delta^0_2$ family $(D_e)_{e \in \Nb}$ of finite sets.
Let $f : [\Nb]^2 \to \Nb$ be a 2-bounded weakly rainbow-stable computable coloring.
For an element $x$, define 
$$
\bad(x) = \set{ y \in \Nb : (\forall^\infty s)c(x,s) = c(y,s)}
$$
Notice that $x \in \bad(x)$.
Because $f$ is weakly rainbow-stable, $\bad$ is a $\Delta^0_2$ function.
For a set $S$, $\bad(S) = \bigcup_{x \in S} \bad(x)$.
Define $X_e = \bad(D_e)$. Hence $X_e$ is a $\Delta^0_2$ set, and this uniformly in $e$.
Moreover, $\card{X_e} \leq 2\card{D_e}$.

Let $h : \Nb \to \Nb$ be a function satisfying $(\forall e)(\forall n)[\card{X_e} \leq n \imp h(e, n) \not \in X_e]$.
We can define $g : \Nb \to \Nb$ by $g(e) = h(e, 2{\card{D_e} \choose 2})$.
Hence $(\forall e)g(e) \not \in X_e$.

We construct a prerainbow $R$ by stages $R_0 (=\emptyset) \subsetneq R_1 \subsetneq R_2, \dots$
as in Lemma~\ref{lem:char1-srrt22}.
Assume that at stage $s$, 
$(\forall \{x,y\} \subseteq R_s)(\forall^{\infty} s)[f(x, s)\neq f(y,s)]$.
Because $R_s$ is finite, we can computably find some index $e$ such that $R_s = D_e$.
Set $R_{s+1} = R_s \cup \set{g(e)}$. By definition, $g(e) \not \in X_e$.
Let $x \in R_s$. Because $g(e) \not \in X_e$, $(\forall^{\infty} s) f(x, s) \neq f(g(e), s)$.
By~$\ist$, the set~$R$ is an $f$-prerainbow.
By Lemma~\ref{lem:prerainbow-equiv} we can compute an infinite $f$-rainbow from $R \oplus f$.
\end{proof}

\subsection{Lowness and bushy tree forcing}

In this section, we prove that the rainbow Ramsey theorem for pairs
restricted to weakly rainbow-stable colorings is strictly weaker
than the full rainbow Ramsey theorem for pairs, by constructing
an~$\omega$-model of~$\wsrrt^2_2$ having only low sets.
As~$\rrt^2_2$ does not admit such a model,
$\wsrrt^2_2$ does not imply $\rrt^2_2$ over~$\rca$.

\begin{theorem}\label{thm:wsrrt22-low}
For every set~$X$ and every weakly rainbow-stable $X$-computable 2-bounded function~$f : [\Nb]^2 \to \Nb$,
there exists an infinite $f$-rainbow low over~$X$.
\end{theorem}

We will use \emph{bushy tree forcing} for building a low solution to a computable
instance of~$\wsrrt^2_2$. This forcing notion has been successfully used for proving
many properties of d.n.c.\ degrees~\cite{Ambos-Spies2004Comparing,Bienvenu2016Diagonally,Khan2015Forcing,Patey2015Ramsey}. 
Indeed, the power of a d.n.c.\ function is known
to be equivalent to finding a function escaping a uniform family of c.e.\ sets~\cite{Kjos-Hanssen2011Kolmogorov},
which is exactly what happens with bushy tree forcing: we build an infinite set by finite approximations,
avoiding a set of bad extensions whose size is computably bounded.
We start by stating the definitions of bushy tree forcing and the basic properties
without proving them. See the survey of Kahn and Miller~\cite{Khan2015Forcing} for a good introduction.

\index{bushy tree}
\begin{definition}[Bushy tree]
Fix a function $h$ and a string $\sigma \in \Nb^{<\Nb}$.
A tree $T$ is  $h$-bushy above $\sigma$ if every $\tau \in T$ is increasing
and comparable with $\sigma$ and whenever $\tau \succeq \sigma$ is not a leaf of~$T$,
it has at least $h(\card{\tau})$ immediate children. 
We call $\sigma$ the \emph{stem} of $T$.
\end{definition}

\index{big set}
\index{small set}
\begin{definition}[Big set, small set]
Fix a function $h$ and some string $\sigma \in \Nb^{<\Nb}$.
A set $B \subseteq \Nb^{<\Nb}$ is \emph{$h$-big above $\Nb$} 
if there exists a finite tree $T$ $h$-bushy above $\sigma$ 
such that all leafs of $T$ are in $B$. 
If no such tree exists, $B$ is said to be \emph{$h$-small above $\sigma$}.
\end{definition}

Consider for example a weakly rainbow-stable 2-bounded function~$f : [\Nb]^2 \to \Nb$.
We want to construct an infinite $f$-prerainbow. We claim that
the following set is~$id$-small above~$\epsilon$, where~$id$ is the identity function:
\[
B_f = \{ \sigma \in \Nb^{<\Nb} : (\exists x, y \in \sigma)(\forall^{\infty} s)f(x,s) = f(y, s)\}
\]
Indeed, given some string~$\sigma \not \in B_f$,
there exists at most~$|\sigma|$ integers~$x$ such that~$\sigma x \in B_f$.
Therefore, given any infinite tree which is $h$-bushy above~$\emptyset$,
at least one of the paths will be an $f$-prerainbow. Also note that because
$f$ is weakly rainbow-stable, the set~$B_f$ is~$\Delta^{0,f}_2$.
We now state some basic properties about bushy tree forcing.


\begin{lemma}[Smallness additivity] \label{lem:smallness-add}
Suppose that $B_1, B_2, \ldots, B_n$ are subsets of $\Nb^{<\Nb}$, 
$g_1$, $g_2$, ..., $g_n$ are functions, and $\sigma \in \Nb^{<\Nb}$.
If $B_i$ is $g_i$-small above~$\sigma$ for all~$i$, then $\bigcup_i B_i$ is $(\sum_i g_i)$-small above $\sigma$. 
\end{lemma}

\begin{lemma}[Small set closure]\label{lem:small-set-closure}
We say that $B  \subseteq \Nb^{<\Nb}$ is \emph{$g$-closed} if whenever $B$ is $g$-big above a string $\rho$ then $\rho \in B$. Accordingly, the \emph{$g$-closure} of any set~$B \subseteq \Nb^{<\Nb}$ is the set $C = \set{\tau \in \Nb^{<\Nb} : B \mbox{ is $g$-big above } \tau}$. If $B$ is $g$-small above a string~$\sigma$, then its closure is also $g$-small above $\sigma$.
\end{lemma}

Note that if~$B$ is a $\Delta^{0,X}_2$ $g$-small set for some computable function~$g$,
so is the $g$-closure of~$B$. Moreover, one can effectively find a $\Delta^{0,X}_2$ index of the $g$-closure of~$B$
given a $\Delta^{0,X}_2$ index of~$B$.
Fix some set~$X$.
Our forcing conditions are tuples~$(\sigma, g, B)$ where
$\sigma$ is an increasing string, $g$ is a computable function and~$B \subseteq \Nb^{<\Nb}$ is a $\Delta^{0,X}_2$
$g$-closed set $g$-small above~$\sigma$.
A condition~$(\tau, h, C)$ \emph{extends} $(\sigma, g, B)$ if~$\sigma \preceq \tau$ and~$B \subseteq C$.
Any infinite decreasing sequence of conditions starting with~$(\epsilon, id, B_f)$
will produce an $f$-prerainbow.

The following lemma is sufficient to deduce the existence of a $\Delta^{0,X}_2$ infinite $f$-prerainbow.

\begin{lemma}\label{lem:wsrrt22-ext}
Given a condition~$(\sigma, g, B)$, one can $X'$-effectively
find some~$x \in \Nb$ such that the condition~$(\sigma x, g, B)$ is a valid extension.
\end{lemma}
\begin{proof}
Pick the first~$x \in \Nb$ greater than $\sigma(|\sigma|)$ such that~$\sigma x \not \in B$.
Such~$x$ exists as there are at most~$g(|\sigma|)-1$ many bad~$x$ by~$g$-smallness of~$B$.
Moreover~$x$ can be found $X'$-effectively as~$B$ is~$\Delta^{0,X}_2$.
By $g$-closure of~$B$, $B$ is $g$-small above~$\sigma x$. Therefore~$(\sigma x, g, B)$
is a valid extension.
\end{proof}

A sequence~$G$ \emph{satisfies} the condition~$(\sigma, g, B)$
if it is increasing, $\sigma \prec G$ and~$B$ is~$g$-small above~$\tau$
for every~$\tau \prec G$.
We say that~$(\sigma, g, B) \Vdash \Phi^{G \oplus X}_e(e) \downarrow$
if~$\Phi^{\sigma \oplus X}_e(e) \downarrow$, and~$(\sigma, g, B) \Vdash \Phi^{G \oplus X}_e(e) \uparrow$
if~$\Phi^{G \oplus X}_e(e) \uparrow$ for every sequence~$G$ satisfying the condition~$(\sigma, g, B)$.
The following lemma decides the jump of the infinite set constructed.

\begin{lemma}\label{lem:wsrrt22-force}
Given a condition~$(\sigma, g, B)$ and an index~$e \in \Nb$,
one can $X'$-effectively find some extension~$d = (\tau, h, C)$ 
such that~$d \Vdash \Phi^{G \oplus X}_e(e) \downarrow$ or
$d \Vdash \Phi^{G \oplus X}_e(e) \uparrow$. 
Moreover, one can~$X'$-decide which of the two holds.
\end{lemma}
\begin{proof}
Consider the following~$\Sigma^{0,X}_1$ set:
\[
D = \{ \tau \in \Nb^{<\Nb} : \Phi^{\tau \oplus X}_e(e) \downarrow \}
\]
The question whether~$D$ is~$g$-big above~$\sigma$ is~$\Sigma^{0,X}_1$ and therefore can be~$X'$-decided.
\begin{itemize}
	\item If the answer is yes, we can~$X$-effectively find a finite tree~$T$ $g$-bushy above~$\sigma$
	witnessing this. As~$B$ is~$\Delta^{0,X}_2$, we can take~$X'$-effectively some leaf~$\tau \in T$.
	By definition of~$T$, $\sigma \prec \tau$. As~$B$ is $g$-closed, $B$ is $g$-small above~$\tau$,
	and therefore~$(\tau, g, B)$ is a valid extension. Moreover~$\Phi^{\tau \oplus X}_e(e) \downarrow$.
	\item If the answer is no, the set~$D$ is $g$-small above~$\sigma$.
	By the smallness additivity property (Lemma~\ref{lem:smallness-add}),
	$B \cup D$ is $2g$-small above~$\sigma$. We can $X$-effectively find a $\Delta^0_2$
	index for its $2g$-closure~$C$. The condition~$(\sigma, 2g, C)$ is a valid extension
	forcing~$\Phi^{G \oplus X}_e(e) \uparrow$.
\end{itemize}
\end{proof}

We are now ready to prove Theorem~\ref{thm:wsrrt22-low}.

\begin{proof}[Proof of Theorem~\ref{thm:wsrrt22-low}]
Fix some set~$X$ and some weakly rainbow-stable $X$-computable 2-bounded function~$f : [\Nb]^2 \to \Nb$.
Thanks to Lemma~\ref{lem:wsrrt22-ext} and Lemma~\ref{lem:wsrrt22-force}, define an infinite decreasing
$X'$-computable sequence of conditions~$c_0 \geq c_1 \geq \dots$ starting with~$c_0 = (\epsilon, id, B_f)$ and
such that for each~$s \in \Nb$,
\begin{itemize}
	\item[(i)] $|\sigma_s| \geq s$
	\item[(ii)] $c_{s+1} \Vdash \Phi^{G \oplus X}_s(s) \downarrow$ or~$c_{s+1} \Vdash \Phi^{G \oplus X}_s(s) \uparrow$
\end{itemize}
where~$c_s = (\sigma_s, g_s, B_s)$. The set~$G = \bigcup_s \sigma_s$
is an $f$-prerainbow. By (i), $G$ is infinite and by~(ii), $G$ is low over~$X$. 
By Lemma~\ref{lem:prerainbow-equiv}, $G \oplus X$ computes an infinite $f$-rainbow.
\end{proof}

\begin{corollary}
There exists an $\omega$-model of $\wsrrt^2_2$ having only low sets.
\end{corollary}

\begin{corollary}
$\wsrrt^2_2$ does not imply~$\rrt^2_2$ over~$\rca$.
\end{corollary}
\begin{proof}
By Corollary~\ref{cor:rrt22-dnrzp}, every model of~$\rrt^2_2$ is a model of~$\dnrs{2}$,
and no function d.n.c.\ relative to~$\emptyset'$ is low.
\end{proof}

\subsection{Relations to other principles}

In this last section, we prove that the rainbow Ramsey theorem
for pairs for weakly rainbow-stable colorings
is a consequence of the stable free set theorem for pairs.
We first need to introduce some useful terminology.

\index{tail rainbow}
\begin{definition}[Wang in~\cite{Wang2014Cohesive}]
Fix a 2-bounded coloring $f : [\Nb]^n \to \Nb$ and~$k \leq n$.
A set~$H$ is a \emph{$k$-tail $f$-rainbow} if $f(\vec{u}, \vec{v}) \neq f(\vec{w}, \vec{x})$
for all~$\vec{u}, \vec{w} \in [H]^{n-k}$ and~ \emph{distinct} $\vec{v}, \vec{x} \in [H]^k$.
\end{definition}

Every 2-bounded coloring~$f : [\Nb]^2 \to \Nb$ admits an infinite $f$-computable
1-tail $f$-rainbow.
Wang proved in~\cite{Wang2014Cohesive} that for every 2-bounded coloring $f : [\Nb]^n \to \Nb$,
every $f$-random computes an infinite 1-tail $f$-rainbow $H$.
We refine this result by the following lemma.

\begin{lemma}[$\rca$]\label{lem:dnc-normal-rrt}
Let $f : [\Nb]^{n+1} \to \Nb$ be a 2-bounded coloring.
Every function d.n.c.\ relative to $f$ computes an infinite 1-tail $f$-rainbow $H$.
\end{lemma}
\begin{proof}
By~\cite{Kjos-Hanssen2011Kolmogorov}, every function d.n.c.\ relative to $f$ computes a function $g$
such that if $|W^f_e| \leq n$ then $g(e, n) \not \in W^f_e$.
Given a finite 1-tail $f$-rainbow $F$,
there exists at most ${\card{F} \choose n}$ elements $x$ such that $F \cup \{x\}$
is not a 1-tail $f$-rainbow. We can define an infinite 1-tail $f$-rainbow $H$ by stages,
starting with $H_0 = \emptyset$. Given a finite 1-tail $f$-rainbow $H_s$ of cardinal $s$,
set $H_{s+1} = H_s \cup \{g(e, {s \choose n})\}$ where $e$ is a Turing index
such that $W^f_e = \{x : H_s \cup \{x\} \mbox{ is not a 1-tail } f\mbox{-rainbow}\}$.
\end{proof}

\begin{theorem}
$\rca + \bst \vdash \sfs^2 \imp \wsrrt^2_2$
\end{theorem}
\begin{proof}
Fix a weakly rainbow-stable 2-bounded coloring $f : [\Nb]^2 \to \Nb$.
As $\rca \vdash \sfs^2 \imp \dnr$, there exists by Lemma~\ref{lem:dnc-normal-rrt}
an infinite 1-tail $f$-rainbow $X$.
We will construct an infinite $X \oplus f$-computable stable coloring 
$g : [X]^2 \to \set{0,1}$ such that every infinite $g$-free set is an $f$-rainbow.
We define the coloring $g : [\Nb]^2 \to \Nb$ by stages as follows.

At stage $s$, assume $g(x, y)$ is defined for every $x, y < s$.
For every pair $x < y < s$ such that $g(x, s) = g(y, s)$, set $g(y, s) = x$.
For the remaining $x < s$, set $g(x, s) = 0$.
This finishes the construction. We now turn to the verification.

\begin{claim}
Every infinite $g$-free set $H$ is an $f$-rainbow.
\end{claim}
\begin{proof*}
Assume for the sake of contradiction that $H$ is not an $f$-rainbow.
Because $X$ is a 1-tail $f$-rainbow and $H \subseteq X$, 
there exists $x, y, s \in H$ such that $c(x, s) = c(y, s)$
with $x < y < s$. As $f$ is 2-bounded, neither $x$ nor $y$ can be part of another pair
$u, v$ such that $f(u, s) = f(v, s)$. So neither $x$ nor $y$ is restrained by another pair
already satisfied, and during the construction we set $g(y, s) = x$.
So $g(y, s) = x$ with $\set{x, y, s} \subset H$, contradicting the freeness of $H$ for $g$.
\end{proof*}

\begin{claim}
The coloring $g$ is stable.
\end{claim}
\begin{proof*}
Fix a $y \in \Nb$. As $f$ is weakly rainbow-stable, we have two cases.
Either there exists an $x < y$ such that $f(y, s) =  f(x, s)$ for cofinitely
many $s$, in which case $g(y, s) = x$ for cofinitely many $s$ and we are done.
Or $f(y, s) \neq f(x, s)$ for each $x < y$ and cofinitely many $s$.
Then by $\bst$, for cofinitely many $s$, $f(y, s) = 0$.
\end{proof*}
\end{proof}

\begin{question}
Does $\sts^2$ imply $\wsrrt^2_2$ over $\rca$ ?
\end{question}

Last, we give an upper bound on the strength of the rainbow Ramsey theorem
which enables us to separate it from the thin set theorem over computable reducibility.

\begin{theorem}\label{thm:generalized-rrt-randomness}
Fix some~$n \geq 2$, and a computable~$\rrt^n_2$-instance~$f : [\omega]^n \to \omega$.
For every~$P \gg \emptyset^{(n-2)}$ and every function~$g$ d.n.c.\ relative to~$P'$,
$P \oplus g$ computes an infinite~$f$-rainbow.
\end{theorem}
\begin{proof}
We prove a relativized version of this theorem by induction over~$n$. 
The case~$n = 2$ is proven by Miller~\cite{MillerAssorted}. He in fact proved the stronger
statement that for every set~$X$, every $X$-computable $\rrt^2_2$-instance~$f$
and every function~$g$ d.n.c.\ relative to~$X'$, $X \oplus g$ computes an infinite $f$-rainbow.

We need to treat the case $n = 3$ independently as well. 
Fix a set~$X$, an $X$-computable $\rrt^3_2$-instance~$f$
and a set $P \gg X'$.  We can assume that~$\omega$ is a 1-tail $f$-rainbow since every $X$-computable $k$-bounded
coloring admits an $X$-computable 1-tail rainbow.
For each~$\sigma, \tau \in [\omega]^n$, let~$R_{\sigma, \tau} = \{ s : f(\sigma, s) = f(\tau, s) \}$.
By Jockusch and Stephan~\cite{Jockusch1993cohesive}, there
is an infinite $\vec{R}$-cohesive set~$C$ such that~$(X \oplus C)' \leq_T P$.
By Theorem~\ref{thm:wsrrt-jump-of-rrt}, there is an $(X \oplus C)'$-computable
$\rrt^2_2$-instance~$h$ such that every infinite $h$-rainbow~$H$ $X \oplus C$-computes
an infinite $f$-rainbow. Using the stronger statement proven by Miller relativized to~$P$,
for every function~$g$ d.n.c.\ relative to~$P'$, $P \oplus g$ computes 
an infinite $h$-rainbow, hence $P \oplus g$-computes an infinite $f$-rainbow.

We now prove the case~$n \geq 4$ by induction.
Fix some set~$X$ and some~$X$-computable $\rrt^n_2$-instance~$f$ such that~$\omega$ is a 1-tail $f$-rainbow.
Fix also a set~$P \gg X^{(n-2)}$ and a function~$g$ d.n.c.\ relative to~$X^{(n-1)}$.
For each~$\sigma, \tau \in [\omega]^n$, let again~$R_{\sigma, \tau} = \{ s : f(\sigma, s) = f(\tau, s) \}$.
By Jockusch and Stephan~\cite{Jockusch1993cohesive}, there
is an infinite $\vec{R}$-cohesive set~$C$ such that~$(X \oplus C)^{(2)} \leq_T Y^{(2)}$.
In particular, $P \gg X^{(n-2)} \geq_T (X \oplus C)^{(n-2)}$ since~$n \geq 4$.
Still by Theorem~\ref{thm:wsrrt-jump-of-rrt}, there is an $(X \oplus C)'$-computable
$\rrt^{n-1}_2$-instance~$h$ such that every infinite $h$-rainbow~$H$ $X \oplus C$-computes
an infinite $f$-rainbow. By induction hypothesis relativized to~$(X \oplus C)'$,
for every function~$g$ d.n.c.\ relative to~$P'$, $P \oplus g$-computes an infinite $h$-rainbow,
hence $P \oplus g$-computes an infinite $f$-rainbow.
\end{proof}

\begin{corollary}
For every~$n \geq 2$, $\sts^n \not \leq_c \rrt^n_2$.
\end{corollary}
\begin{proof}
By Theorem~\ref{thm:sts2-amt}, $\amt \leq_c \sts^2 \leq_c \ts^1[\emptyset']$.
By Csima et al.~\cite{Csima2004Bounding} and Conidis~\cite{Conidis2008Classifying},
$\amt$ is computably equivalent to the statement ``For every $\Delta^0_2$ function~$f$,
there is a function~$g$ not dominated by~$f$''.
By Kurtz~\cite{Kurtz1982Randomness}, for every set~$X$,
there is an $X'$-computable function such that the measure of oracles which $X$-compute a function
not dominated by~$f$ is null.
Therefore, for every set~$X$, there is an $X'$-computable function~$f : \omega \to \omega$
such that the measure of oracles~$Z$ such that~$X \oplus Z$-computes an infinite $f$-thin set is null.
Such a function can also be easily constructed by a direct argument.

Fix a computable $\rrt^n_2$-instance~$g$ for some~$n \geq 2$,
and let~$P \gg \emptyset^{(n-2)}$ be such that~$P' \leq_T \emptyset^{(n-1)}$.
By the previous argument relativized to~$P$, there is a $\emptyset^{(n-1)}$-computable function~$f : \omega \to \omega$
such that the measure of oracles~$Z$ such that~$P \oplus Z$ computes an infinite $\tilde{f}$-thin set is null.
Let~$f : [\omega]^n \to \omega$ be the stable computable function obtained by taking the $\Delta^0_n$ approximation
of the function~$\tilde{f}$. Every infinite $f$-thin set is $\tilde{f}$-thin.
By Ku\v{c}era~\cite{Kucera1985Measure}, the measure of oracles computing a function d.n.c.\ relative to~$P'$ is positive,
so there is some function~$h$ d.n.c.\ relative to~$P'$ such that $P \oplus h$ does not compute an infinite $f$-thin set.
However, by Theorem~\ref{thm:generalized-rrt-randomness}, $P \oplus h$ computes an infinite $g$-rainbow.
\end{proof}

\chapter{The Erd\H{o}s-Moser theorem}\label{chap:erdos-moser-theorem}

The Erd\H{o}s-Moser theorem is a statement from graph theory 
which received a particular interest from reverse mathematical community
as it provides, together with the ascending descending sequence principle,
an alternative proof of Ramsey's theorem for pairs.

\index{Erd\H{o}s-Moser theorem}
\index{em@$\emo$|see {Erd\H{o}s-Moser theorem}}
\index{transitive tournament}
\index{tournament}
\begin{definition}[Erd\H{o}s-Moser theorem] 
A tournament $T$ is an irreflexive binary relation such that for all $x,y \in \omega$ with $x \not= y$, exactly one of $T(x,y)$ or $T(y,x)$ holds. A tournament $T$ is \emph{transitive} if the corresponding relation~$T$ is transitive in the usual sense. 
$\emo$ is the statement ``Every infinite tournament $T$ has an infinite transitive subtournament.''
\end{definition}

Bovykin and Weiermann~\cite{Bovykin2005strength} proved Ramsey's theorem for pairs as follows:
Given a coloring~$f : [\Nb]^2 \to 2$, we can see~$f$ as a tournament~$T$
such that whenever~$x <_\Nb y$, $T(x,y)$ holds if and only if~$f(x,y) = 1$.
Any transitive subtournament~$H$ can be seen as a linear order $(H, \prec)$
such that whenever~$x <_\Nb y$, ~$x \prec y$ if and only if~$f(x,y) = 1$.
Any infinite ascending or descending sequence is $f$-homogeneous.
It is therefore natural to study the ascending descending sequence principle
together with the Erd\H{o}s-Moser theorem.

\index{ascending descending sequence}
\index{ads@$\ads$|see {ascending descending sequence}}
\begin{definition}[Ascending descending sequence]
Given a linear order~$(L, <_L)$, an \emph{ascending} (\emph{descending}) sequence
is a set~$S$ such that for every~$x <_\Nb y \in S$, $x <_L y$ ($x >_L y$).
$\ads$ is the statement ``Every infinite linear order admits an infinite ascending or descending sequence''.
\end{definition}

\index{sads@$\sads$|see {ascending descending sequence}}
$\sads$ is the restriction of~$\ads$ to linear orders of type $\omega+\omega^{*}$.
The above-mentioned argument is easily formalizable over~$\rca$.
Furthermore, both~$\emo$ and~$\ads$ are immediate consequences of Ramsey's theorem
for pairs over~$\rca$. We therefore obtain the following equivalence.

\begin{theorem}[\cite{Bovykin2005strength}]\label{thm:rt22-em-ads}
$\rca \vdash \rt^2_2 \biimp [\emo \wedge \ads]$
\end{theorem}

One can also consider the stable versions of the Erd\H{o}s-Moser theorem
and the ascending descending sequence principle. For the latter statement,
a linear order is stable iff it is of order type~$\omega+\omega^{*}$.
The decomposition of~$\rt^2_2$ into~$\emo$ and~$\ads$ also holds for the stable versions.
The particular shape of the decomposition enables one to make~$\emo$ and~$\semo$ inherit 
several properties of~$\rt^2_2$ and $\srt^2_2$, respectively.

\begin{theorem}\label{thm:strength-ramsey-em-preservation}
If~$\ads$ admits~$\Pcal$ preservation but not~$\rt^2_2$,
then neither does~$\emo$. Moreover, the witnesses of failure for~$\rt^2_2$
and~$\emo$ have the same Turing degree.
\end{theorem}
\begin{proof}
Let~$C \in \Pcal$ and~$f : [\omega]^2 \to 2$ be a $C$-computable coloring
such that for every~$f$-homogeneous set~$H$, $H \oplus C \not \in \Pcal$.
The coloring $f$ can be seen as a tournament $T$
where for each~$x < y$, $T(x, y)$ holds iff $f(x, y) = 1$.
If there is an infinite transitive sub-tournament~$U$ such that~$U \oplus C \in \Pcal$,
then considering~$(U, <_T)$ as a linear order, since~$\ads$ admits~$\Pcal$ preservation,
there is an infinite ascending or descending sequence~$S$ such that~$S \oplus U \oplus C \in \Pcal$.
The set~$S$ is an infinite $f$-homogeneous set such that~$S \oplus C \in \Pcal$,
contradicting our choices of~$C$ and~$f$. Therefore, $C$ and~$T$
witness the failure of~$\Pcal$ preservation for~$\emo$.
\end{proof}

The same theorem holds for the stable version of the statements.
In particular, we can reprove some properties of the Erd\H{o}s-Moser theorem for free.

\begin{corollary}[Kreuzer \cite{Kreuzer2012Primitive}]\label{cor:emo-computable-instance-no-low-solution}
There exists a transitive computable tournament having no low infinite subtournament.
\end{corollary}
\begin{proof}
Downey et al.\ proved in~\cite{Downey200102} that there is a computable instance of~$\srt^2_2$ with no low solution.
On the other hand, Hirschfeldt et al.\ \cite{Hirschfeldt2007Combinatorial} proved $\sads$ admits lowness preservation.
By Theorem~\ref{thm:strength-ramsey-em-preservation}, there is a computable instance of~$\semo$
with no low solution.
\end{proof}

It sometimes happens that the combinatorics of two statements~$\Psf$ and~$\Qsf$
are so incompatible that whenever a $C$-computable $\Qsf$-instance~$X$ as no computable solution,
then every $C$-computable $\Psf$-instance has a solution which does not $C$-compute a solution to~$X$.
In this case, we say that~\emph{$\Psf$ admits $\Qsf$ avoidance}.

\begin{theorem}\label{thm:strength-ramsey-em-principle-preservation}
If $\Qsf \leq_c \rt^2_2$ and $\ads$ admits $\Qsf$ avoidance, then $\Qsf \leq_c \emo$.
\end{theorem}
\begin{proof}
Let $I$ be any instance of $\Qsf$. As $\Qsf \leq_c \rt^2_2$, there exists
an $I$-computable coloring $f : [\Nb]^2 \to 2$ such that for any infinite $f$-homogeneous set $H$,
$I \oplus H$ computes a solution to $I$. The coloring $f$ can be seen as a tournament $T$
where for each~$x < y$, $T(x, y)$ holds iff $f(x, y) = 1$. If $T$ has an infinite sub-tournament
$U$ such that $H \oplus I$ does not compute a solution to $I$, consider $H$ as an $I \oplus H$-computable
stable linear order. Then since~$\ads$ admits $\Qsf$-avoidance, there exists a solution $S$ to $H$
such that $S \oplus H \oplus I$ does not compute a solution to $I$. But $S$ is an infinite $f$-homogeneous set,
contradicting our choice of $f$.
\end{proof}

Using the stable version of Theorem~\ref{thm:strength-ramsey-em-principle-preservation},
we deduce again properties of the Erd\H{o}s-Moser theorem for free.
Rice~\cite{RiceThin} proved that~$\dnr$ is a consequence of~$\emo$ over~$\rca$.
We shall see in the next section that~$\dnr$ even follows from the stable Erd\H{o}s-Moser theorem
over~$\rca$. Meanwhile, we can deduce the following weaker result.

\begin{corollary}\label{cor:omega-sem-dnr}
$\dnr \leq_c \semo$.
\end{corollary}
\begin{proof}
Hirschfeldt et al.\ proved in~\cite{Hirschfeldt2008strength} that $\dnr \leq_c \srt^2_2$
and in \cite{Hirschfeldt2007Combinatorial} that $\ads$ admits $\dnr$ avoidance.
By the stable version of Theorem~\ref{thm:strength-ramsey-em-principle-preservation}, $\dnr \leq_c \semo$.
\end{proof}

\begin{corollary}
$\coh \leq_c \srt^2_2$ if and only if $\coh \leq_c \semo$
\end{corollary}
\begin{proof}
If $\coh \leq_c \semo$, then $\coh \leq_c \srt^2_2$ since $\semo \leq_c \srt^2_2$.
We now prove the reverse implication. 
We showed in section~\ref{sect:strength-ramsey-cohesiveness-konig-lemma}
that we can associate a $\Pi^{0,\emptyset'}_1$ class $\Ccal(\vec{R})$
to any sequence of sets~$R_0$, $R_1$, ..., so that
a degree bounds an~$\vec{R}$-cohesive set if and only if
its jump bounds a member of~$\Ccal(\vec{R})$.
Hirschfeldt et al.~\cite{Hirschfeldt2007Combinatorial} proved that
every $X$-computable instance $I$ of~$\sads$ has a solution $Y$ low over~$X$.
Therefore, if $X$ does not compute an~$\vec{R}$-cohesive set,
then~$X'$ does not compute a member of~$\Ccal(\vec{R})$.
As~$(Y \oplus X)' \leq X'$, $(Y \oplus X)'$ does not compute a member of~$\Ccal(\vec{R})$,
$Y \oplus X$ does not compute an~$\vec{R}$-cohesive set.
In other words, $\sads$ admits $\coh$ avoidance. Conclude by 
the stable version of Theorem~\ref{thm:strength-ramsey-em-principle-preservation}.
\end{proof}

\section{The strength of the Erd\H{o}s-Moser theorem}

We shall see (Corollary~\ref{cor:emo-wkl-sts2-sads}) that~$\emo$
does not imply the stable thin set theorem for pairs over~$\rca$.
The Erd\H{o}s-Moser theorem is not known to imply the cohesiveness principle as well.
However, a simple combinatorial argument shows that it implies the disjunction of~$\sts^2$
and~$\coh$.

\begin{theorem}\label{thm:emo-implies-sts2-or-coh}
$\rca \vdash \emo \rightarrow [\sts^2 \vee \coh]$
\end{theorem}
\begin{proof}
Let $f : [\Nb]^2 \to \Nb$ be a stable coloring and $R_0, R_1, \dots$ be a uniform sequence of sets.
We denote by $\tilde{f}$ the function defined by $\tilde{f}(x) = \lim_s f(x,s)$.
We build a $\Delta^{0, f \oplus \vec{R}}_1$ tournament $T$ such that every infinite transitive subtournament
is either thin for $\tilde{f}$ or is an $\vec{R}$-cohesive set. As every set $H$ thin for $\tilde{f}$
$H \oplus f$-computes a set thin for $f$, we are done.
For each $x, s \in \Nb$, set $T(x,s)$ to hold if one of the following holds:
\begin{itemize}
	\item[(i)] $f(x,s) = 2i$ and $x \in R_i$  
	\item[(ii)] $f(x,s) = 2i+1$ and $x \not \in R_i$
\end{itemize}
Otherwise set $T(s,x)$ to hold. Let $H$ be an infinite transitive subtournament of $T$
which is not $\tilde{f}$-thin. Suppose for the sake of contradiction that $H$ is not $\vec{R}$-cohesive.
Then there exists an $i \in \Nb$ such that $H$ intersects $R_i$ and $\overline{R_i}$ infinitely many times.
As $H$ is not $\tilde{f}$-thin, there exists $x, y \in H$ such that $\tilde{f}(x) = \lim_s f(x,s) = 2i$
and $\tilde{f}(y) = \lim_s f(y,s) = 2i+1$. As $H$ intersects $R_i$ and $\overline{R_i}$ infinitely many times,
there exists $s_0 \in R_i \cap H$ and $s_1 \in \overline{R_i} \cap H$ 
such that $f(x,s_0) = f(x,s_1) = 2i$ and $f(y,s_0) = f(y,s_1) = 2i+1$.
But then $T(x,s_0)$, $T(s_0, y)$, $T(y, s_1)$ and $T(s_1, x)$ hold, forming a 4-cycle
and therefore contradicting the transitivity of~$H$.
\end{proof}

\begin{definition}
Let $T$ be a tournament on a domain $D \subseteq \N$. A $n$-cycle is
a tuple $(x_1, \dots, x_n) \in D^n$ such that 
for every $0 < i < n$, $T(x_i, x_{i+1})$ holds
and $T(x_n, x_1)$ holds.
\end{definition}

Kang~\cite{Kang2014Combinatorial} attributed to Wang a direct proof of $\rca \vdash \emo \imp \rrt^2_2$.
We provide an alternative proof using the characterization of $\rrt^2_2$
by $\dnrzp$ from~Miller.

\begin{theorem}\label{thm:em-dnrzp}
$\rca \vdash \emo \rightarrow \dnrzp$
\end{theorem}
\begin{proof}

Let $X$ be a set. Let $g(.,.)$ be a total $X$-computable function such that $\Phi_e^{X'}(e)=\lim_s g(e,s)$ if the limit exists, and $\Phi^{X'}_e(e) \uparrow$ if the limit does not exist. Interpret $g(e,s)$ as the code of a finite set $D_{e,s}$ of size $3^{e+1}$. We define the tournament~$T$ by $\Sigma_1$-induction as follows. Set $T_0=\emptyset$. At stage~$s+1$, do the following. Start with $T_{s+1}=T_s$. Then, for each $e<s$, take the first pair $\{x,y\} \in D_{e,s} \setminus \bigcup_{k<e} D_{k,s}$ (notice that such a pair exists the by cardinality assumptions on the $D_{e,s}$), and if $T_{s+1}(s,x)$ and $T_{s+1}(s,y)$ are not already assigned, assign them in a way that $(x,y,s)$ forms a $3$-cycle in $T_{s+1}$. Finally, for any $z<s$ such that $T_{s+1}(s,z)$ remains undefined, assign any truth value to it in a predefined way (e.g., for any such pair $\{x,y\}$, set $T_{s+1}(x,y)$ to be true if $x<y$, and false otherwise). This finishes the construction of $T_{s+1}$. Set $T=\bigcup_s T_s$, which must exist as a set by $\Sigma_1$-induction. 

First of all, notice that $T$ is a tournament of domain $[\N]^2$, as at the end of stage~$s+1$ of the construction $T(x,y)$ is assigned a truth value for (at least) all pairs $\{x,y\}$ with $x<s$ and $y<s$. By $\emo$, let $H$ be an infinite, transitive subtournament of~$T$. Let $f(e)$ be the code of the finite set $A_e$ consisting of the first $3^{e+1}$ elements of~$H$. We claim that $f(e) \not=\Phi_e^{X'}(e)$ for all~$e$, which would prove $\dnrzp$. Suppose otherwise, i.e., suppose that $\Phi_e^{X'}(e)=f(e)$ for some~$e$. Then there is a stage $s_0$ such that $f(e)=g(e,s)$ for all~$s \geq s_0$ or equivalently $D_{e,s}=A_e$ for all~$s \geq s_0$. Let $N_e=\max(A_e)$. We claim that for any $s$ be bigger than both $\max(\bigcup_{e,s < N_e} D_{e,s})$ and $s_0$, the restriction of $T$ to $A_e \cup \{s\}$ is not a transitive subtournament, which contradicts the fact that the restriction $H$ of $T$ to the infinite set $A$ containing $A_e$ is transitive. 

To see this, let $s$ be such a stage. At that stage~$s$ of the construction of $T$, a pair $\{x,y\} \in D_{e,s} \setminus \bigcup_{k<e} D_{k,s}$ is selected, and since $D_{e,s}=A_e$, this pair is contained in $A_e$. Furthermore, we claim that $T(s,x)$ and $T(s,y)$ become assigned at that precise stage, i.e., were not assigned before. This is because, by construction of $T$, when the value of some $T(a,b)$ is assigned at a stage~$t$, either $a \leq t$ or $b \leq t$. Thus, if $T(s,x)$ was already assigned at the beginning of stage~$s$, it would have in fact been assigned during or before stage~$x$. However, $x \in A_e$, so $x < N_e$, and at stage~$N_e$ the number $s$, by definition of $N_e$, has not appeared in the construction yet. In particular $T(s,x)$ is not assigned at the end of stage~$x$. This proves our claim, therefore $T(s,x)$ and $T(s,y)$ do become assigned exactly at stage~$s$, in a way -- still by construction -- that $\{x,y,s\}$ form a $3$-cycle for~$T$. Therefore the restriction of $T$ to $A_e \cup \{s\}$ is not a transitive subtournament, which is what we needed to prove. 
\end{proof}

\begin{corollary}[Wang in~\cite{Kang2014Combinatorial}]\label{cor:emo-rrt22}
$\rca \vdash \emo \to \rrt^2_2$
\end{corollary}
\begin{proof}
Immediate by Theorem~\ref{thm:em-dnrzp} and Corollary~\ref{cor:rrt22-dnrzp}.
\end{proof}

We have seen (see Corollary~\ref{cor:omega-sem-dnr})
that every $\omega$-model of~$\semo$ is a model of~$\dnr$.
We now give a direct proof of it and show that it holds over~$\rca$.

\begin{theorem}
$\rca \vdash \semo \rightarrow \dnr$
\end{theorem}
\begin{proof}
This is obtained by small variation of the proof of Theorem~\ref{thm:em-dnrzp}.
Fix a set $X$. Let $g(.,.)$ be a total $X$-computable function such that 
$\Phi^X_e(e)=\lim_s g(e,s)$ if~$\Phi^X_e(e) \downarrow$
and~$\lim_s g(e,s) = 0$ otherwise. Interpret $g(e,s)$ as a code
of a finite set $D_{e,s}$ of size $3^{e+1}$ such that~$min(D_{e,s}) \geq e$ and construct the infinite tournament~$T$
accordingly. The argument for constructing a function d.n.c.\ relative to~$X$
given an infinite transitive subtournament is similar. We will only prove that
the tournament~$T$ is stable.

Fix some~$u \in \Nb$. 
By $\bst$, which is provable from $\semo$ over $\rca$ (see \cite{Kreuzer2012Primitive}),
there exists some stage~$s_0$ after which~$D_{e,s}$ remains constant
for every~$e \leq u$. If~$u$ is part of a pair~$\{x,y\} \subset D_{e,s}$
for some~$s \geq s_0$ and~$e$, then~$e \leq u$ because $min(D_{e,s}) \geq e$.
As the~$D_{e,s}$'s remain constant for each~$e \leq u$, 
the pair~$\{x,y\}$ will be chosen at every stage~$s \geq s_0$
and therefore~$T(u, s)$ will be assigned the same value for every~$s \geq s_0$.
If~$u$ is not part of a pair~$\{x,y\}$, it will always be assigned
the default value at every stage~$s \geq s_0$.
In both cases, $T(u,s)$ stabilizes at stage~$s_0$.
\end{proof}

\section{The combinatorics of the Erd\H{o}s-Moser theorem}\label{sect:combinatorics-emo}

The standard way of building an infinite object by forcing consists of defining an increasing
sequence of finite approximations, and taking the union of them. Unlike
$\coh$ where every finite set can be extended to an infinite cohesive set,
some finite transitive subtournaments may not be extensible to an infinite one.
We therefore need to maintain some extra properties which will guarantee that
the finite approximations are extendible.
The nature of these properties constitute the core of the combinatorics of $\emo$.

Lerman et al.~\cite{Lerman2013Separating} proceeded to an analysis
of the Erd\H{o}s-Moser theorem.
They showed in particular that it suffices to ensure that the finite transitive subtournament~$F$
has infinitely many \emph{one-point extensions}, that is, infinitely many elements~$x$ such that
$F \cup \{x\}$ is transitive, to extend~$F$ to an infinite transitive subtournament (see~\cite[Lemma 3.4]{Lerman2013Separating}).
This property is sufficient to add elements one by one to the finite approximation.
However, when adding elements by block, we shall maintain a stronger invariant. We will require that
the reservoir is included in a minimal interval of the finite approximation~$F$.
In this section, we reintroduce the terminology of Lerman et al.~\cite{Lerman2013Separating}
and give a presentation of the combinatorics of the Erd\H{o}s-Moser theorem
motivated by its computational analysis.

\index{minimal interval}
\begin{definition}[Minimal interval]
Let $R$ be an infinite tournament and $a, b \in R$
be such that $R(a,b)$ holds. The \emph{interval} $(a,b)$ is the
set of all $x \in R$ such that $R(a,x)$ and $R(x,b)$ hold.
Let $F \subseteq R$ be a finite transitive subtournament of $R$.
For $a, b \in F$ such that $R(a,b)$ holds, we say that $(a,b)$
is a \emph{minimal interval of $F$} if there is no $c \in F \cap (a,b)$,
i.e., no $c \in F$ such that $R(a,c)$ and $R(c,b)$ both hold.
\end{definition}

Fix a computable tournament~$R$, and consider a pair $(F, X)$ where
\begin{itemize}
	\item[(i)] $F$ is a finite $R$-transitive set representing the \emph{finite approximation}
	of the infinite $R$-transitive subtournament we want to construct
	\item[(ii)] $X$ is an infinite set disjoint from $F$, included in a minimal interval of~$F$
	and such that $F \cup \{x\}$ is $R$-transitive for every $x \in X$. In other words,
	$X$ is an infinite set of one-point extensions.
	Such a set $X$ represents the \emph{reservoir}, that is, a set of candidate
	elements we may add to~$F$ later on.
\end{itemize}

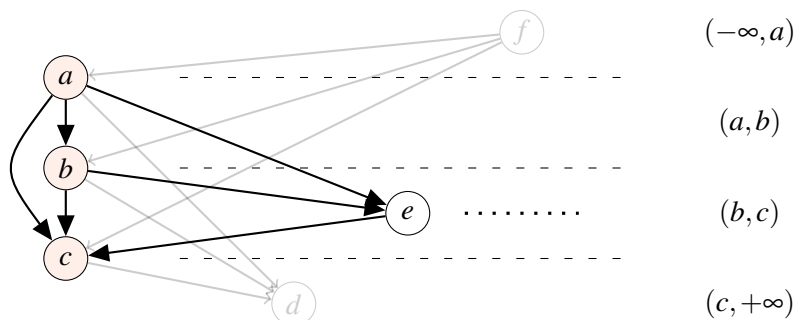
\begin{figure}[htbp]
\begin{center}
\begin{tikzpicture}[x=1.5cm, y=1.2cm, 
		node/.style={circle, draw, fill=ocre!10, inner sep=0pt, minimum size=1.5em}, 
		arrow/.style={black,thick,->,-triangle 45},
		good/.style={fill=white},
		hiddenarrow/.style={black,thick,->, opacity=0.2},
		hidden/.style={opacity=0.2}
]

  \node[node] (a) at (1, 2) {$a$};
	\node[node] (b) at (1, 1) {$b$};
	\node[node] (c) at (1, 0) {$c$};
	\node[node,good,hidden] (d) at (3, -0.5) {$d$};
	\node[node,good] (e) at (4, 0.5) {$e$};
	\node[node,good,hidden] (f) at (5, 2.5) {$f$};

	\node at (7, 2.5) {$(-\infty, a)$};
	\node at (7, 1.5) {$(a, b)$};
	\node at (7, 0.5) {$(b, c)$};
	\node at (7, -0.5) {$(c, +\infty)$};

	\draw[arrow] (a) -- (b);
	\draw[arrow] (b) -- (c);
	\draw[arrow] (a)  .. controls (0.4,1) .. (c);

	\draw[hiddenarrow] (a) -- (d);
	\draw[hiddenarrow] (c) -- (d);
	\draw[hiddenarrow] (b) -- (d);

	\draw[arrow] (a) -- (e);
	\draw[arrow] (e) -- (c);
	\draw[arrow] (b) -- (e);

	\draw[hiddenarrow] (f) -- (a);
	\draw[hiddenarrow] (f) -- (c);
	\draw[hiddenarrow] (f) -- (b);

	\draw[very thick, loosely dotted] (4.5,0.5) -- (5.5,0.5);

	\draw[loosely dashed] (2, 2) -- (6, 2);
	\draw[loosely dashed] (2, 1) -- (6, 1);
	\draw[loosely dashed] (2, 0) -- (6, 0);
\end{tikzpicture}
\end{center}
\caption{In this figure, $F = \{a, b, c\}$ is a transitive set,
$X = \{d, e, f, \dots \}$ a set of one-point extensions, $(b, c) = \{e, \dots \}$ a minimal interval of~$F$ 
and $(F, X \cap (b, c))$ an EM condition. The elements~$d$ and~$f$ are not part of the minimal interval~$(b, c)$.} 
\end{figure}

\index{$\to_R$}
The infinite set~$X$ ensures extensibility of the finite set~$F$ into an infinite $R$-transitive
subtournament. Indeed, by applying the Erd\H{o}s-Moser theorem to $R$ over the domain $X$, there exists an infinite $R$-transitive
subtournament $H \subseteq X$. One easily checks that $F \cup H$ is $R$-transitive.
The pair $(F, X)$ is called an Erd\H{o}s-Moser condition in~\cite{Patey2015Degrees}.
A set~$G$ \emph{satisfies} an EM condition~$(F, X)$ if it is $R$-transitive and satisfies the Mathias condition~$(F, X)$.
In order to simplify notation, given a tournament $R$ and two sets~$E$ and~$F$,
we denote by $E \to_R F$ the formula $(\forall x \in E)(\forall y \in F) R(x,y)$.

Suppose now that we want to add a finite number of elements of~$X$ into $F$ to obtain
a finite $T$-transitive set $\tilde{F} \supseteq F$,
and find an infinite subset $\tilde{X} \subseteq X$ such that $(\tilde{F}, \tilde{X})$
has the above mentioned properties. We can do this in a few steps:

\begin{itemize}
	\item[1.] Choose a finite (non-necessarily $R$-transitive) set $E \subset X$.
	\item[2.] Any element $x \in X \setminus E$ induces a 2-partition $\tuple{E_0, E_1}$ of $E$
	by setting $E_0 = \{y \in E : R(y, x) \}$ and $E_1 = \{y \in E : R(x, y)\}$.
	Consider the coloring $f$ which associates to any element of $X \setminus E$ the corresponding 2-partition $\tuple{E_0, E_1}$ of $E$.
	\item[3.]
	As~$E$ is finite, there exists finitely many 2-partitions of~$E$, so $f$ colors $X \setminus E$ with
	finitely many 2-partitions. By Ramsey's theorem for singletons applied to~$f$, there exists a 2-partition $\tuple{E_0, E_1}$ of $E$
	together with an infinite subset $\tilde{X} \subseteq X \setminus E$ such that for every $x \in \tilde{X}$, $f(x) = \tuple{E_0, E_1}$.
	By definition of~$f$ and~$E_i$, $E_0 \to_R \tilde{X} \to_R E_1$.
	
	\item[4.] Take any $R$-transitive subset $F_1 \subseteq E_i$ for some~$i < 2$ and set $\tilde{F} = F \cup F_1$.
	The pair $(\tilde{F}, \tilde{Y})$ satisfies the required properties (see Lemma~\ref{lem:emo-cond-valid} for a proof).
\end{itemize}

In a computational point of view, if we start with a computable condition~$(F, X)$, that is, where~$X$ is a computable set,
we end up with a computable extension~$(\tilde{F}, \tilde{Y})$.
Remember that our goal is to define a $\Delta^0_2$ function~$f$ which will dominate every $G$-computable function
for some solution~$G$ to~$R$. For this, we need to be able to~$\emptyset'$-decide
whether~$\Phi^G_e(n) \downarrow$ or $\Phi^G_e(n) \uparrow$ for every solution~$G$ to~$R$ satisfying 
some condition~$(F, X)$. 
More generally, given some~$\Sigma^0_1$ formula $\varphi$, we focus on the computational power required to decide a question of the form

\smallskip
{\itshape
Q1: Is there an $R$-transitive extension $\tilde{F}$ of $F$ in $X$ such that $\varphi(\tilde{F})$ holds?
}
\smallskip

Trying to apply naively the algorithm above requires a lot of computational power.
In particular, step 3 requires to choose a true formula among finitely many $\Pi^{0, X}_2$ formulas.
Such a step needs the power of PA degree relative to the jump of~$X$.
We shall apply the same trick as for cohesiveness, consisting in not trying to choose a true $\Pi^{0,X}_2$ formula,
but instead parallelizing the construction. Given a finite set $E \subset X$, 
instead of finding an infinite subset $\tilde{Y} \subset X \setminus E$
whose members induce a 2-partition of $E$, we will construct
as many extensions of $(F, X)$ as there are 2-partitions of~$E$. The question now becomes

\smallskip
{\itshape
Q2: Is there a finite set $E \subseteq X$ such that for every 2-partition $\tuple{E_0, E_1}$ of~$E$,
there exists an $R$-transitive subset $F_1 \subseteq E_i$ for some $i < 2$ such that $\varphi(F \cup F_1)$ holds?
}
\smallskip

This question is $\Sigma^{0,X}_1$, which is good enough for our purposes.
If the answer is positive, we will try the witness $F_1$ associated to each 2-partition of $E$ in parallel.
Note that there may be some 2-partition $\tuple{E_0, E_1}$ of~$E$
such that the set $Y = \{ x \in X \setminus E : E_0 \to_R \{x\} \to_R E_1 \}$ is finite,
but this is not a problem since there is \emph{at least}
one good 2-partition such that the corresponding set is infinite. 
The whole construction yields again a tree of pairs~$(F, X)$.

If the answer is negative, we want to ensure that
$\varphi(\tilde{F})$ will not hold at any further stage of the construction.
For each~$n \in \omega$, let $H_n$ be the set of the $n$ first elements of~$X$.
Because the answer is negative, for each~$n \in \omega$, there exists a 2-partition $\tuple{E_0, E_1}$
of~$H_n$ such that for every $R$-transitive subset $F_1 \subseteq E_i$ for any $i < 2$, $\varphi(F \cup F_1)$ does not hold.
Call such a 2-partition an \emph{avoiding} partition of $H_n$. 
Note that if $\tuple{E_0, E_1}$ is an avoiding partition of $H_{n+1}$, then $\tuple{E_0 \uh n, E_1 \uh n}$
is an avoiding partition of $H_n$. So the set of avoiding 2-partitions of some $H_n$
forms an infinite tree~$T$. Moreover, the predicate ``$\tuple{E_0, E_1}$ is an avoiding partition of $H_n$''
is $\Delta^{0, H_n}_1$ so the tree $T$ is $\Delta^{0,X}_1$. The collection of the infinite paths through $T$
forms a non-empty $\Pi^{0,X}_1$ class $\Ccal$ defined as the collection of 2-partitions $Z_0 \cup Z_1 = X$
such that for every $i < 2$ and every $R$-transitive subset $F_1 \subseteq Z_i$, $\varphi(F \cup F_1)$
does not hold. Apply weak K\"onig's lemma to obtain a 2-partition~$Z_0 \cup Z_1 = X$
such that for every finite $R$-transitive subset $F_1$ of any of its parts, $\varphi(F \cup F_1)$ does not hold.

\section{The weakness of the Erd\H{o}s-Moser theorem}\label{sect:weakness-erdos-moser-theorem}

We now study the computational weakness of the Erd\H{o}s-Moser theorem
using the framework of preservation of hyperimmunity. The forcing notion
to construct solutions to the Erd\H{o}s-Moser theorem shares many properties
with the one for Ramsey's theorem for pairs. The main difference is that
only one object is built, namely, a transitive subtournament, whereas
in the case of Ramsey's theorem for pairs, a set homogeneous for color~0
and another one for color~1 are built simultaneously. 
It follows that there is no disjunction in the requirements of the Erd\H{o}s-Moser theorem
and that countably many hyperimmunities can be preserved simultaneously.

\begin{theorem}\label{thm:em-hyperimmunity}
$\emo$ admits preservation of hyperimmunity.
\end{theorem}

Before proving Theorem~\ref{thm:em-hyperimmunity}, we introduce
the notion of \emph{Erd\H{o}s Moser condition} and prove some basic combinatorial lemmas.

\index{Erd\H{o}s-Moser condition}
\index{EM condition|see {Erd\H{o}s-Moser condition}}
\begin{definition}
An \emph{Erd\H{o}s Moser condition} (EM condition) for an infinite tournament~$T$
is a Mathias condition $(F, X)$ where
\begin{itemize}
	\item[(a)] $F \cup \{x\}$ is $T$-transitive for each $x \in X$
	\item[(b)] $X$ is included in a minimal $T$-interval of $F$.
\end{itemize}
\end{definition}

The extension relation is the usual Mathias extension.
EM conditions have good properties for tournaments as stated by the following lemmas.
Given a tournament $T$ and two sets $E$ and $F$,
we denote by $E \to_T F$ the formula $(\forall x \in E)(\forall y \in F) T(x,y) \mbox{ holds}$.

\begin{lemma}\label{lem:emo-cond-beats}
Fix an EM condition $(F, X)$ for a tournament $T$.
For every $x \in F$, $\{x\} \to_T X$ or $X \to_T \{x\}$.
\end{lemma}
\begin{proof}
Fix an $x \in F$. Let $(u,v)$ be the minimal $T$-interval containing $X$,
where $u, v$ may be respectively $-\infty$ and $+\infty$.
By definition of interval, $\{u\} \to_T X \to_T \{v\}$.
By definition of minimal interval, $T(x,u)$ or $T(v,x)$ holds.
Suppose the former holds.
By transitivity of $F \cup \{y\}$ for every $y \in X$, $T(x,y)$ holds, therefore $\{x\} \to_T X$.
In the latter case, by symmetry, $X \to_T \{x\}$.
\end{proof}

\begin{lemma}\label{lem:emo-cond-valid}
Fix an EM condition $c = (F, X)$ for a tournament $T$, 
an infinite subset $Y \subseteq X$ and a finite $T$-transitive set $F_1 \subset X$ such that
$F_1 < Y$ and $[F_1 \to_T Y \vee Y \to_T F_1]$.
Then $d = (F \cup F_1, Y)$ is a valid extension of $c$.
\end{lemma}
\begin{proof}
Properties of a Mathias condition for $d$ are immediate.
We prove property~(a). Fix an $x \in Y$.
To prove that $F \cup F_1 \cup \{x\}$ is $T$-transitive,
it suffices to check that there exists no 3-cycle in $F \cup F_1 \cup \{x\}$.
Fix three elements $u < v < w \in F \cup F_1 \cup \{x\}$.
\begin{itemize}
	\item Case 1: $\{u, v, w\} \cap F \neq \emptyset$. Then $u \in F$ as $F < F_1 < \{x\}$
	and $u < v < w$. If $v \in F$ then using the fact that $F_1 \cup \{x\} \subset X$
	and property (a) of condition $c$, $\{u, v, w\}$ is $T$-transitive.
	If $v \not \in F$, then by Lemma~\ref{lem:emo-cond-beats},
	$\{u\} \to_T X (\supseteq F \cup \{x\})$ or $X \to_T \{u\}$
	hence $\{u\} \to_T \{v,w\}$ or $\{v,w\} \to_T \{u\}$ so $\{u, v, w\}$ is $T$-transitive.
	
	\item Case 2: $\{u, v, w\} \cap F = \emptyset$. Then at least $u, v \in F_1$ because $F_1 < \{x\}$.
	If $w \in F_1$, then $\{u,v,w\}$ is $T$-transitive by $T$-transitivity of $F_1$.
	Otherwise, as $F_1 \to_T Y$ or $Y \to_T F_1$, $\{u,v\} \to_T \{w\}$
	or $\{w\} \to_T \{u,v\}$ and $\{u,v,w\}$ is $T$-transitive.
\end{itemize}
We now prove property (b).
Let $(u,v)$ be the minimal $T$-interval of $F$ in which $X$ (hence $Y$) is included by property (b) of
condition $c$. $u$ and $v$ may be respectively $-\infty$ and $+\infty$.
By assumption, either $F_1 \to_T Y$ or $Y \to_T F_1$. As $F_1$ is a finite $T$-transitive set,
it has a minimal and a maximal element, say $x$ and $y$. 
If $F_1 \to_T Y$ then $Y$ is included in the $T$-interval $(y, v)$.
Symmetrically, if $Y \to_T F_1$ then $Y$ is included in the $T$-interval $(u, x)$.
To prove minimality for the first case, assume that some $w$ is in the interval $(y, v)$.
Then $w \not \in F$ by minimality of the interval $(u,v)$ w.r.t. $F$,
and $w \not \in F_1$ by maximality of $y$.
Minimality for the second case holds by symmetry. 
\end{proof}

We are now ready to prove Theorem~\ref{thm:em-hyperimmunity}.

\begin{proof}[Proof of Theorem~\ref{thm:em-hyperimmunity}]
Let~$B_0, B_1, \dots$ be a countable sequence of~$C$-hyperimmune sets
for some set~$C$, and let~$T$ be a $C$-computable tournament.
We will build an infinite transitive subtournament~$G$
such that the~$B$'s are $G \oplus C$-hyperimmune.
For this, we use Mathias forcing~$(F, X)$ where~$X$
is an infinite set such that the~$B$'s are $X \oplus C$-hyperimmune.
The following lemma is the core of the argument.

\begin{lemma}\label{lem:emo-preserve-hyperimmunity}
Fix a condition~$c = (F,X)$ and two integers $e,i \in \omega$.
There is an extension~$d$ of~$c$ forcing~$\Phi_e^{G \oplus C}$ not to dominate $p_{B_i}$.
\end{lemma}
\begin{proof*}
Let~$f$ be the partial $X \oplus C$-computable function which on input~$x$
searches for a finite set of integers~$U$ such that for every 2-partition $Z_0 \cup Z_1 = X$,
there is some~$j < 2$ and a finite $T$-transitive set~$E \subseteq Z_j$ such that
$\Phi^{(F \cup E) \oplus C}_e(x) \downarrow \in U$. If such a set~$U$ is found,
$f(x) = 1+max(U)$, otherwise~$f(x) \uparrow$. We have two cases.
\begin{itemize}
	\item Case 1: $f$ is total. By $X \oplus C$-hyperimmunity of~$B_i$, there
	is some~$x$ such that~$f(x) \leq p_{B_i}(x)$. Let~$U$ be the finite set witnessing~$f(x) \downarrow$.
	By compactness, there exists a finite set~$H \subset X$ such that for every partition~$H_0 \cup H_1 = H$,
	there is some~$j < 2$ and a $T$-transitive set~$E \subseteq H_j$
	such that~$\Phi^{(F \cup E) \oplus C}_e(x) \downarrow \in U$.
	Each element~$y \in X$ induces a partition~$H_0 \cup H_1 = H$ such that~$H_0 \to_T \{y\} \to_T H_1$.
	There exists finitely many such partitions, so by the infinite pigeonhole principle,
	there exists an~$X$-computable infinite set~$Y \subset X$ and a partition~$H_0 \cup H_1 = H$
	such that~$H_0 \to_T Y \to_T H_1$. Let~$j < 2$ and~$E \subseteq H_j$ be the $T$-transitive set
	such that $\Phi^{(F \cup E) \oplus C}_e(x) \downarrow \in U$.
	By Lemma~\ref{lem:emo-cond-valid}, $(F \cup E, Y)$ is a valid extension
	forcing~$\Phi^{G \oplus C}_e(x) \downarrow \leq p_{B_i}(x)$.

	\item Case 2: there is some~$x$ such that~$f(x) \uparrow$. By compactness,
	the~$\Pi^{0,X \oplus Z}_1$ class~$\Ccal$ of sets $Z_0 \oplus Z_1$ such that~$Z_0 \cup Z_1 = X$ and
	for every~$j < 2$ and every $T$-transitive set~$E \subseteq Z_j$,
	$\Phi^{(F \cup E) \oplus C}_e(x) \uparrow$. 
	By preservation of hyperimmunity of~$\wkl$, there exists some partition~$Z_0 \oplus Z_1 \in \Ccal$
	such that the~$B$'s are~$Z_0 \oplus Z_1 \oplus Z$-hyperimmune. The set~$Z_j$
	is infinite for some~$j < 2$ and the condition~$(F, Z_i)$ is an EM extension
	forcing~$\Phi^{G \oplus C}_e(x) \uparrow$.
\end{itemize}
\end{proof*}

Thanks to Lemma~\ref{lem:emo-cond-valid} and Lemma~\ref{lem:emo-preserve-hyperimmunity}, 
define an infinite descending sequence of conditions
$(\emptyset, \omega) \geq c_0 \geq \dots$ such that for each~$s \in \omega$,
\begin{itemize}
	\item[(a)] $|F_s| \geq s$
	\item[(b)] $c_s$ forces~$\Phi^{G \oplus C}_e$ not to dominate $p_{B_i}$ if~$s = \tuple{e,i}$
\end{itemize}
where~$c_s = (F_s, X_s)$. Let~$G = \bigcup_s F_s$.
By (a), $G$ is infinite and by (b), the~$B$'s are $G \oplus C$-hyperimmune.
This finishes the proof of Theorem~\ref{thm:em-hyperimmunity}.
\end{proof}

\begin{corollary}\label{cor:sads-not-preservation-2-hyperimmunities}
Neither $\ads$ nor $\sads$ admit preservation of 2 hyperimmunities.
\end{corollary}
\begin{proof}
By the stable version of Theorem~\ref{thm:rt22-em-ads}, $\rca \vdash \srt^2_2 \biimp [\semo \wedge \sads]$.
If~$\sads$ admitted preservation of 2 hyperimmunities, then so would~$\rt^2_2$,
contradicting Corollary~\ref{cor:ts2-non-preservation}.

In a more direct way, let~$f : [\omega]^2 \to 2$
be the stable computable function of Theorem~\ref{thm:strength-ramsey-delta2-hyperimmune}.
Let~$L = (\omega, <_L)$ be defined for each~$x <_\omega y$
by~$x <_L y$ iff $f(x, y) = 0$, and~$x >_L y$ otherwise.
By item (ii) of Theorem~\ref{thm:strength-ramsey-delta2-hyperimmune}, $L$ is a linear order
of order type~$\omega+\omega^{*}$ and by item (i) of Theorem~\ref{thm:strength-ramsey-delta2-hyperimmune},
the $\omega$ part and the~$\omega^{*}$ part are both hyperimmune.
The $\omega$ (resp.\ $\omega^{*}$) part is hyperimmune-free relative to any infinite ascending (resp.\ descending) sequence.
\end{proof}

The following corollary has been proven by Wang~\cite{Wang2014Definability} using the notion
of preservation of non-c.e.\ definitions.

\begin{corollary}[Wang~\cite{Wang2014Definability}]\label{cor:emo-wkl-sts2-sads}
$\rca \wedge \emo \wedge \coh \wedge \wkl \nvdash \sts^2 \vee \sads$.
\end{corollary}
\begin{proof}
By Theorem~\ref{thm:em-hyperimmunity}, Theorem~\ref{thm:coh-hyperimmunity-preservation} 
and the hyperimmune-free basis theorem~\cite{Jockusch197201},
$\emo$, $\coh$ and~$\wkl$ admit preservation of hyperimmunity. 
Let~$f : [\omega]^2 \to \omega$ be the stable computable function of Theorem~\ref{thm:strength-ramsey-delta2-hyperimmune},
let~$B_i = \{x : \lim_s f(x, s) = i\}$ and let~$L = (\omega, <_L)$ be defined as in Corollary~\ref{cor:sads-not-preservation-2-hyperimmunities}. The complement of the~$B$'s are all hyperimmune.
The set~$\overline{B_i}$ is hyperimmune-free relative to any infinite $f$-thin set for color~$i$,
and $\overline{B_1}$ (resp.\ $\overline{B_0}$) is hyperimmune-free relative to any infinite ascending
(resp.\ descending) sequence. Lemma~\ref{lem:intro-reduc-preservation-separation} enables us to conclude.
\end{proof}

Note that by Theorem~\ref{thm:emo-implies-sts2-or-coh}, every model of~$\rca \wedge \emo$ which is not a model of~$\sts^2$
is a model of~$\coh$. Therefore, Theorem~\ref{thm:em-hyperimmunity} implies Theorem~\ref{thm:coh-hyperimmunity-preservation}.

There has been recently a lot of literature around the weakness of the Erd\H{o}s-Moser theorem.
Lerman, Solomon and Towsner~\cite{Lerman2013Separating} proved that~$\rca \wedge \emo \nvdash \sads$.
The author~\cite{Patey2013note} refined their proof to obtain~$\rca \wedge \emo \nvdash \sts^2$.
Wang~\cite{Wang2014Definability} enhanced these results by proving that
$\rca \wedge \emo \wedge \coh \wedge \wkl \nvdash \sts^2 \vee \sads$.
Recently, the author~\cite{PateyDominating} strengthened all the above-mentioned results
by proving that~$\emo \wedge \coh \wedge \wkl$ does not even imply the atomic model theorem ($\amt$).
This is a strictly stronger result since $\amt$ is a consequence of both~$\sads$ and~$\sts^2$ over~$\rca$.

Even though the deep combinatorics remain the same, the argument to separate the Erd\H{o}s-Moser
theorem from~$\amt$ is significantly more complicated. Indeed, the atomic model theorem
is a genericity notion. In particular, for every computable, complete atomic theory,
every set of hyperimmune degree relative to~$\emptyset'$ computes an atomic model.
Therefore, every solution sufficiently generic for the notion of forcing of the Erd\H{o}s-Moser
theorem will compute solutions to the atomic model theorem. One needs to restrict the amount of genericity
of the solutions by making the overall construction effective.

Over the next sections, we will prove computable non-reducibility results, which are therefore strictly weaker
than a separation over~$\rca$, but already contain the main ideas of the general construction. See~\cite{PateyDominating}
for the whole proof.

\section{Dominating cohesive sets}\label{sect:dominating-erdos-moser-dominating-cohesive}

Before proving that the atomic model theorem does not computably
reduce to the Erd\H{o}s-Moser theorem theorem,
we illustrate the key features of our construction by 
showing that $\amt$ does not computable reduce to~$\coh$.
The remainder of this section is devoted to the proof of the following theorem.

\begin{theorem}\label{thm:amt-comp-reduc-coh}
$\amt \not \leq_c \coh$
\end{theorem}

In order to prove Theorem~\ref{thm:amt-comp-reduc-coh},
we need to construct a $\Delta^0_2$ function $f$ such that 
for every uniformly computable sequence of sets~$R_0, R_1, \dots$,
there is an $\vec{R}$-cohesive set~$G$ such that every~$G$-computable
function is dominated by~$f$. Thankfully, Jockusch \& Stephan~\cite{Jockusch1993cohesive}
proved that for every such sequence of sets~$\vec{R}$, 
every p-cohesive set computes an infinite $\vec{R}$-cohesive set.
The sequence of all primitive recursive sets is therefore called a \emph{universal instance}.
Hence we only need to build a $\Delta^0_2$ function~$f$ and a p-cohesive set~$G$
such that every $G$-computable function is dominated by~$f$ to obtain Theorem~\ref{thm:amt-comp-reduc-coh}.

Given some uniformly computable sequence of sets~$R_0, R_1, \dots$,
the usual construction of an $\vec{R}$-cohesive set~$G$ is done by a computable Mathias forcing.
The forcing conditions are pairs $(F,X)$, where~$F$ is a finite set representing the finite
approximation of~$G$ and~$X$ is an infinite, computable reservoir such that~$max(F) < min(X)$.
The construction of the~$\vec{R}$-cohesive set is obtained by building
an infinite, decreasing sequence of Mathias conditions, starting with~$(\emptyset, \omega)$
and interleaving two kinds of steps.
Given some condition~$(F,X)$,
\begin{itemize}
	\item[(S1)] the \emph{extension} step consists in taking an element $x$ from $X$ and adding it to~$F$,
	therefore forming the extension $(F \cup \{x\}, X \setminus [0,x])$;
	\item[(S2)] the \emph{cohesiveness} step consists of deciding which one of $X \cap R_i$
	and $X \cap \overline{R}_i$ is infinite, and taking the chosen one as the new reservoir.
\end{itemize}
The first step ensures that the constructed set~$G$ will be infinite, whereas
the second step makes the set $G$ $\vec{R}$-cohesive.
Looking at the effectiveness of the construction, the step (S1) is computable,
assuming we are given some Turing index of the set~$X$.
The step (S2), on the other hand, requires to decide which one of two computable sets
is infinite, knowing that at least one of them is. This decision
requires the computational power of a PA degree relative to~$\emptyset'$ (see \cite[Lemma 4.2]{Cholak2001strength}).
Since we want to build a $\Delta^0_2$ function~$f$ dominating every $G$-computable function,
we would like to make the construction of~$G$ $\Delta^0_2$. Therefore the step (S2) has to be revised.

\subsection{Effectively constructing a cohesive set}

The above construction leads to two observations.
First, at any stage of the construction, the reservoir~$X$ of the Mathias condition~$(F, X)$
has a particular shape. Indeed, after the first application of stage~(S2), 
the set $X$ is, up to finite changes, of the form $\omega \cap R_0$
or $\omega \cap \overline{R_0}$. After the second application of (S2), it is in one of the following forms: $\omega \cap R_0 \cap R_1$,
$\omega \cap R_0 \cap \overline{R}_1$, $\omega \cap \overline{R}_0 \cap R_1$,
$\omega \cap \overline{R}_0 \cap \overline{R}_1$, and so on. More generally, given some string~$\sigma \in 2^{<\omega}$,
we can define~$R_\sigma$ inductively as follows:
First, $R_\varepsilon = \omega$, and then, if $R_\sigma$ has already been defined for some string $\sigma$ of length~$i$,
$R_{\sigma 0} = R_\sigma \cap \overline{R}_i$ and~$R_{\sigma 1} = R_\sigma \cap R_i$.
By the first observation, we can replace Mathias conditions by pairs ~$(F, \sigma)$, where $F$ is a finite set
and $\sigma \in 2^{<\omega}$. The pair~$(F, \sigma)$ denotes the Mathias condition~$(F, R_\sigma \cap (max(F), +\infty))$.
A pair $(F, \sigma)$ is \emph{valid} if $R_\sigma$ is infinite.
The step (S2) can be reformulated as choosing, given some valid condition $(F, \sigma)$, which one of $(F, \sigma 0)$
and $(F, \sigma 1)$ is valid.

\begin{center}
\begin{tikzpicture}[x=2.6cm, y=2cm, 
		node/.style={}, 
		arrow/.style={black, thick,->},
		continue/.style={loosely dotted, thick},
		hidden/.style={opacity=0.2}
]

  \node[node] (a0) at (1, 1) {$(\emptyset, \varepsilon)$};
	\node[node] (a00) at (2, 1) {$(\{x\}, \varepsilon)$};
	\node[node] (a000) at (3, 1.5) {$(\{x\}, 0)$};
	\node[node] (a001) at (3, 0.5) {$(\{x\}, 1)$};
	\node[node] (a0000) at (4, 1.5) {$(\{x, y\}, 0)$};
	\node[node] (a0010) at (4, 0.5) {$(\{x, z\}, 1)$};
	\node[node] (a00000) at (5, 2) {$(\{x, y\}, 00)$};
	\node[node] (a00001) at (5, 1.5) {$(\{x, y\}, 01)$};
	\node[node] (a00100) at (5, 0.5) {$(\{x, z\}, 10)$};
	\node[node] (a00101) at (5, 0) {$(\{x, z\}, 11)$};

	\draw[arrow] (a0) to  node [above, black] {\tiny (S1)} (a00);
	\draw[arrow] (a00) -- (a000);
	\draw[arrow] (a00) -- (a001);
	\draw[arrow] (a001) to  node [below] {\tiny (S1)} (a0010);
	\draw[arrow] (a000) to  node [above] {\tiny (S1)} (a0000);
	\draw[arrow] (a0010) -- (a00100);
	\draw[arrow] (a0010) to  node [below] {\tiny (S2)} (a00101);
	\draw[arrow] (a0000) to  node [above] {\tiny (S2)} (a00000);
	\draw[arrow] (a0000) -- (a00001);

	\node[node] (s20) at (2.5,1) {\tiny (S2)};

	\draw[continue] (5.5, 2) -- (6, 2);
	\draw[continue] (5.5, 1.5) -- (6, 1.5);
	\draw[continue] (5.5, 0.5) -- (6, 0.5);
	\draw[continue] (5.5, 0) -- (6, 0);
\end{tikzpicture}
\end{center}

Second, we do not actually need to decide which one of~$R_{\sigma 0}$ and~$R_{\sigma 1}$ is infinite
assuming that~$R_{\sigma}$ is infinite. Our goal is to dominate every~$G$-computable function with a $\Delta^0_2$ function $f$.
Therefore, given some $G$-computable function~$g$, it is sufficient to find a finite set $S$ of candidate values for~$g(x)$ 
and make~$f(x)$ be greater than the maximum of $S$. Instead of choosing which one of~$R_{\sigma 0}$ and~$R_{\sigma 1}$ is infinite,
we will explore both cases in parallel. The step (S2) will split some condition $(F, \sigma)$
into two conditions~$(F, \sigma 0)$ and $(F, \sigma 1)$. Our new forcing conditions are therefore tuples $(F_\sigma : \sigma \in 2^n)$
which have to be thought of as $2^n$ parallel Mathias conditions $(F_\sigma, \sigma)$ for each~$\sigma \in 2^n$.
Note that $(F_\sigma, \sigma)$ may not denote a valid Mathias condition in general since $R_\sigma$ may be finite.
Therefore, the step (S1) becomes~$\Delta^0_2$, since we first have to check whether $R_\sigma$ is non-empty
before picking an element in~$R_\sigma$. The whole construction is $\Delta^0_2$ and yields a $\Delta^0_2$ infinite
binary tree $T$. In particular, any degree PA relative to $\emptyset'$ bounds an infinite path though $T$ and therefore bounds a $G$-cohesive set. However, the degree of the set $G$ is not sensitive in our argument. We only care about the effectiveness of
the tree~$T$.

\subsection{Dominating the functions computed by a cohesive set}

We have seen in the previous section how to make the construction of a cohesive set more effective
by postponing the choices between forcing~$G \subseteq^{*} R_i$ and~$G \subseteq^{*} \overline{R}_i$
to the end of the construction. We now show how to dominate every $G$-computable function
for every infinite path~$G$ through the $\Delta^0_2$ tree constructed in the previous section.
To do this, we will interleave a third step deciding whether $\Phi^G_e(n)$ halts, and if so, collecting the 
candidate values of~$\Phi^G_e(n)$.
Given some Mathias precondition~$(F, X)$ (a precondition is a condition 
where we do not assume that the reservoir is infinite) and some $e,x \in \omega$, we can $\Delta^0_2$-decide 
whether there is some set $E \subseteq X$ such that~$\Phi^{F \cup E}_e(x) \downarrow$.
If this is the case, then we can effectively find this a finite set~$E \subseteq X$ and
compute the value~$\Phi^{F \cup E}_e(x)$. If this is not the case, then for every infinite set~$G$ satisfying
the condition~$(F, X)$, the function~$\Phi^G_e$ will not be defined on input~$x$. In this case,
our goal is vacuously satisfied since $\Phi^G_e$ will not be a function and therefore
we do not need do dominate~$\Phi^G_e$. 
Let us go back to the previous construction.
After some stage, we have constructed a condition~$(F_\sigma : \sigma \in 2^n)$ inducing a finite tree of depth~$n$. 
The step (S3) acts as follows for some~$x \in \omega$:
\begin{itemize}
	\item[(S3)] Let~$S = \{0\}$. For each~$\sigma \in 2^n$ and each~$e \leq x$, decide whether 
	there is some finite set~$E \subseteq R_\sigma \cap (max(F_\sigma), +\infty)$ such that~$\Phi^{F_\sigma \cup E}_e(x) \downarrow$.
	If this is the case, add the value of $\Phi^{F_\sigma \cup E}_e(x)$ to $S$ and set~$\tilde{F}_\sigma = F_\sigma \cup E$, otherwise set $\tilde{F}_\sigma = F_\sigma$.
	Finally, set~$f(x) = max(S)+1$ and take~$(\tilde{F}_\sigma : \sigma \in 2^n)$ as the next condition.
\end{itemize}
Note that the step (S3) is $\Delta^0_2$-computable uniformly in the condition~$(F_\sigma : \sigma \in 2^n)$.
The whole construction therefore remains~$\Delta^0_2$ and so does the function~$f$.
Moreover, given some $G$-computable function~$g$, there is some Turing index~$e$ such that~$\Phi^G_e = g$.
For each $x \geq e$, the step (S3) is applied at a finite stage and decides whether~$\Phi^G_e(x)$ halts or not
for every set satisfying one of the leaves of the finite tree. In particular, this is the case for the set $G$ and 
therefore $\Phi^G_e(x) \in S$. By definition of $f$, $f(x) \geq max(S) \geq \Phi^G_e(x)$. Therefore $f$ dominates the function~$g$.

\subsection{The formal construction}

Let~$R_0, R_1, \dots$ be the sequence of all primitive recursive sets.
We define a $\Delta^0_2$ decreasing sequence of conditions~$(\emptyset, \varepsilon) \geq c_0 \geq c_1 \dots$
such that for each~$s \in \omega$
\begin{itemize}
	\item[(i)] $c_s = (F^s_\sigma : \sigma \in 2^s)$ and~$|F^s_\sigma| \geq s$ if~$R_\sigma \cap (max(F^s_\sigma), +\infty) \neq \emptyset$.
	\item[(ii)] For every~$e \leq s$ and every~$\sigma \in 2^s$, either $\Phi^{F^s_\sigma}_e(s) \downarrow$
	or $\Phi^G_e(s) \uparrow$ for every set~$G$ satisfying $(F^s_\sigma, R_\sigma)$.
\end{itemize}
Let~$P$ be a path through the tree $T = \{ \sigma \in 2^{<\omega} : R_\sigma \mbox{ is infinite} \}$
and let~$G = \bigcup_s F^s_{P \restr s}$. By (i), for each~$s \in \omega$, $|F^s_{P \restr s}| \geq s$
since~$R_{P \restr s}$ is infinite. Therefore the set~$G$ is infinite.
Moreover, for each~$s \in \omega$, the set~$G$ satisfies the condition~$(F^{s+1}_{P \restr s+1}, R_{P \restr s+1})$,
so~$G \subseteq^{*} R_{P \restr s+1} \subseteq R_s$ if~$P(s) = 1$ and
$G \subseteq^{*} R_{P \restr s+1} \subseteq \overline{R}_s$ if~$P(s) = 0$. 
Therefore~$G$ is $\vec{R}$-cohesive.

For each~$s \in \omega$, let~$f(s) = 1 + max(\Phi^{F^s_\sigma}_e(s) : \sigma \in 2^s, e \leq s)$.
The function~$f$ is $\Delta^0_2$. We claim that it dominates every $G$-computable function.
Fix some~$e$ such that~$\Phi^G_e$ is total. For every~$s \geq e$, let~$\sigma = P \restr s$. By (ii), either
$\Phi^{F^s_\sigma}_e(s) \downarrow$ or $\Phi^G_e(s) \uparrow$ for every
set~$G$ satisfying $(F^s_\sigma, R_\sigma)$. Since $\Phi^G_e(s) \downarrow$,
the first case holds. By definition of~$f$, $f(s) \geq \Phi^{F^s_\sigma}_e(s) = \Phi^G_e(s)$.
Therefore~$f$ dominates the function $\Phi^G_e$. This completes the proof of Theorem~\ref{thm:amt-comp-reduc-coh}.

\section{Dominating the Erd\H{o}s-Moser theorem}\label{sect:dominating-the-erdos-moser-theorem}

We now strengthen the analysis of the previous section by
proving that the atomic model theorem is not computably reducible to the Erd\H{o}s-Moser theorem.
Theorem~\ref{thm:amt-comp-reduc-coh} is an immediate consequence of this result since
$[\amt \vee \coh] \leq_c \emo$ (see~\cite{Patey2015Somewhere}).

\begin{theorem}\label{thm:amt-comp-reduc-em-coh}
$\amt \not \leq_c \emo$
\end{theorem}

Just as we did for cohesiveness, we will show how to build solutions to $\emo$ through $\Delta^0_2$ constructions 
postponing the $\Pi^0_2$ choices to the end.

\subsection{Enumerating the computable infinite tournaments}

Proving that some principle~$\Psf$ does not computably reduce to~$\Qsf$
requires to create a $\Psf$-instance~$X$ such that \emph{every} $X$-computable $\Qsf$-instance
has a solution~$Y$ such that~ $Y \oplus X$ does not compute a solution to~$X$.
In the case of $\amt \not \leq_c \coh$, we have been able to restrict ourselves to only one instance of $\coh$,
since Jockusch \& Stephan~\cite{Jockusch1993cohesive} showed it admits a universal instance.
It is currently unknown whether the Erd\H{o}s-Moser theorem admits a universal instance, that is, a computable infinite tournament
such that for every infinite transitive subtournament $H$ and for every computable infinite tournament $T$,
$H$ computes an infinite transitive $T$-subtournament. See~\cite{Patey2015Degrees} for an extensive study of the existence
of universal instances for principles in reverse mathematics.

Since we do not know whether $\emo$ admits a universal instance, we will need to diagonalize against
the solutions to every computable $\emo$-instance. In fact, we will prove a stronger result. We will construct
a $\Delta^0_2$ function~$f$ and an infinite set~$G$ which is eventually transitive simultaneously for every computable infinite tournament,
and such that $f$ dominates every $G$-computable function. There exists no computable sequence of sets
containing all computable sets. Therefore it is not possible to computably enumerate every infinite computable tournament.
However, one can define an infinite, computable, binary tree such that every infinite path
computes such a sequence. 
See the notion of sub-uniformity defined by Mileti in~\cite{Mileti2004Partition} for details.
By the low basis theorem, there exists a low set bounding a sequence containing, 
among others, every infinite computable tournament.
As we shall prove below,  for every set~$C$ and every uniformly $C$-computable sequence of infinite tournaments~$\vec{R}$,
there exists a set~$G$ together with a $\Delta^{0, C}_2$ function $f$ such that
\begin{itemize}
	\item[(i)] $G$ is eventually $R$-transitive for every $R \in \vec{R}$
	\item[(ii)] If $\Phi^{G \oplus C}_e$ is total, then it is dominated by $f$ for every $e \in \omega$.
\end{itemize}
Thus it suffices to choose~$C$ to be our low set and $\vec{R}$ to be a uniformly $C$-computable sequence
of infinite tournaments containing every computable tournament to deduce the existence of a set~$G$ together
with a $\Delta^0_2$ function $f$ such that 
\begin{itemize}
	\item[(i)] $G$ is eventually $R$-transitive for every infinite, computable tournament $R$
	\item[(ii)] If $\Phi^{G \oplus C}_e$ is total, then it is dominated by $f$ for every $e \in \omega$
\end{itemize}

By the computable equivalence between $\amt$ and the escape property,
there exists a computable atomic theory $T$ such that every atomic model computes
a function~$g$ not dominated by~$f$. If $\amt \leq_c \emo$, then there exists
an infinite, computable tournament~$R$ such that every infinite $R$-transitive subtournament 
computes a model of~$T$, hence computes a function~$g$ not dominated by~$f$.
As the set~$G$ is, up to finite changes, an infinite $R$-transitive subtournament,
$G$ computes such a function~$g$, contradicting our hypothesis. Therefore $\amt \not \leq_c \emo$.

\subsection{Cover classes}

\index{k-cover@$k$-cover}
In this part, we introduce some terminology about classes of $k$-covers.
Recall that a $k$-cover is a $k$-partition whose parts are not required to be pairwise disjoint.

\smallskip
\index{cover class}
\emph{Cover class}.
We code a $k$-cover $Z_0 \cup \dots \cup Z_{k-1}$ of some set $X$ as the set $Z = \bigoplus_{i < k} Z_i$.
We will identify a $k$-cover with its code.
A \emph{$k$-cover class} of some set~$X$ is a tuple $\tuple{k, X, \Ccal}$
where $\Ccal$ is a collection of codes of $k$-covers of~$X$. 
We will be interested in $\Pi^0_1$ $k$-cover classes.
For the simplicity of notation, we may use the same letter~$\Ccal$ to denote both a $k$-cover class~$(k, X, \Ccal)$
and the actual collection of $k$-covers~$\Ccal$. We then write $dom(\Ccal)$ for $X$
and $parts(\Ccal)$ for~$k$.

\smallskip
\emph{Restriction of a cover}. Given some $k$-cover $Z = Z_0 \oplus \dots \oplus Z_{k-1}$ of some set~$X$ and given some set~$Y \subseteq X$, we write $Z \restr Y$ for the $k$-cover $(Z_0 \cap Y) \oplus \dots \oplus (Z_{k-1} \cap Y)$ of~$Y$.
Similarly, given some cover class~$(k, X, \Ccal)$ and some set~$Y \subseteq X$, we denote by $\Ccal \restr Y$
the cover class~$(k, Y, \Dcal)$ where $\Dcal = \{ Z \restr Y : Z \in \Ccal \}$.
Given some part~$\nu$ of $\Ccal$ and some set~$E$, we write~$\Ccal^{[\nu, E]}$
for the cover class~$(k, X, \Dcal)$ where 
$\Dcal = \{ Z_0 \oplus \dots \oplus Z_{k-1} \in \Ccal : E \subseteq Z_\nu \}$.

\smallskip
\emph{Refinement}. The collection of cover classes can be given a natural partial order as follows.
Let~$m \geq k$ and $f : m \to k$. An $m$-cover $V_0 \oplus \dots \oplus V_{m-1}$ of $Y$ \emph{$f$-refines}
a $k$-cover $Z_0 \oplus \dots \oplus Z_{k-1}$ of $X$ if $Y \subseteq X$ and $V_\nu \subseteq Z_{f(\nu)}$ for each~$\nu < m$.
Given two cover classes $(k, X, \Ccal)$ and~$(m, Y, \Dcal)$
and some function $f : m \to k$, we say that $\Dcal$ \emph{$f$-refines} $\Ccal$
if for every $V \in \Dcal$, there is some~$Z \in \Ccal$ such that $V$ $f$-refines $Z$.
In this case, we say that \emph{part $\nu$ of $\Dcal$ refines part~$f(\nu)$ of~$\Ccal$}.

\smallskip
\emph{Acceptable part}. 
We say that part $\nu$ of $\Ccal$ is \emph{acceptable} if there exists some $Z_0 \oplus \dots \oplus Z_{k-1} \in \Ccal$
such that $Z_\nu$ is infinite. Part $\nu$ of $\Ccal$ is \emph{empty} if 
for every $Z_0 \oplus \dots \oplus Z_{k-1} \in \Ccal$, $Z_\nu = \emptyset$.
Note that if $\Ccal$ is non-empty and $dom(\Ccal)$ is infinite, then $\Ccal$ has at least one acceptable part.
Moreover, if~$\Dcal \leq_f \Ccal$ and part~$\nu$ of~$\Dcal$ is acceptable, then so is part~$f(\nu)$ of $\Ccal$.
The converse does not hold in general.

\subsection{The forcing notion}

We now get into the core of our forcing argument by defining
the forcing notion which will be used to build an infinite set eventually
transitive for every infinite computable tournament.
Fix a set $C$ and a uniformly $C$-computable sequence of infinite tournaments $R_0, R_1, \dots$
We construct our set~$G$ by a forcing whose conditions are tuples $(\alpha, \vec{F}, \Ccal)$ where
\begin{itemize}
	\item[(a)] $\Ccal$ is a non-empty $\Pi^{0,C}_1$ $k$-cover class of $[t, +\infty)$ 
	for some $k, t \in \omega$ ; $\alpha \in t^{<\omega}$
	\item[(b)] $(F_\nu \setminus [0, \alpha(i))) \cup \{x\}$ is $R_i$-transitive for every $Z_0 \oplus \dots \oplus Z_{k-1} \in \Ccal$,
	every $x \in Z_\nu$, every $i < |\alpha|$ and each $\nu < k$
	\item[(c)] $Z_\nu$ is included in a minimal $R_i$-interval of $F_\nu \setminus [0, \alpha(i))$
	for every $Z_0 \oplus \dots \oplus Z_{k-1} \in \Ccal$, every $i < |\alpha|$ and each $\nu < k$.
\end{itemize}
A condition $d = (\beta, \vec{E}, \Dcal)$ \emph{extends} 
$c = (\alpha, \vec{F}, \Ccal)$ (written $d \leq c$) if $\beta \succeq \alpha$
and there exists a function $f : parts(\Dcal) \to parts(\Ccal)$ such that the following holds:
\begin{itemize}
	\item[(i)] $(E_\nu, dom(\Dcal))$ Mathias extends $(F_{f(\nu)}, dom(\Ccal))$ for each $\nu < parts(\Dcal)$ 
	\item[(ii)] $\Dcal$ $f$-refines $\Ccal^{[f(\nu), E_\nu \setminus F_{f(\nu)}]}$ for each $\nu < parts(\Dcal)$
\end{itemize}

One may think of a condition $(\alpha, \vec{F}, \Ccal)$ where $k = parts(\Ccal)$
as $k$ parallel Mathias conditions which are,
up to finite changes, Erd\H{o}s-Moser conditions simultaneously for the tournaments $R_0, \dots, R_{|\alpha|-1}$.
Given some $i < |\alpha|$, the value $\alpha(i)$ indicates at which point the sets $\vec{F}$
start being $R_i$-transitive.
More precisely, for every part $\nu < k$ and every $k$-cover $Z_0 \oplus \dots \oplus Z_{k-1} \in \Ccal$,
$(F_\nu \setminus [0, \alpha(i)), Z_\nu)$ is an Erd\H{o}s-Moser condition for $R_i$ for each $i < |\alpha|$.
Indeed, because of clause~(i), the elements $E_\nu \setminus F_{f(\nu)}$ added to $E_\nu$ come from $dom(\Ccal)$
and because of clause~(ii), these elements must come from the part $f(\nu)$ of the class~$\Ccal$,
otherwise $\Ccal^{[f(\nu), E_\nu \setminus F_{f(\nu)}]}$ would be empty and so would be $\Dcal$.

Of course, there may be some parts~$\nu$ of $\Ccal$ which are non-acceptable, that is, such that $Z_\nu$ is finite
for every $k$-cover $Z_0 \oplus \dots \oplus Z_{k-1} \in \Ccal$. However, by the infinite pigeonhole principle, 
$Z_\nu$ must be infinite for at least one $\nu < k$.
Choosing $\alpha$ to be in $t^{<\omega}$ instead of $\omega^{<\omega}$
ensures that all elements added to~$\vec{F}$ will have to be $R_i$-transitive
simultaneously for each~$i < |\alpha|$, as the elements are taken from $dom(\Ccal)$
and therefore are greater than the threshold $\alpha(i)$ for each $i < |\alpha|$.
By abuse of terminology, we may talk about a \emph{part of a condition}
when talking about a part of its corresponding cover class.

We start with a few basic lemmas reflecting the combinatorics described in 
the section~\ref{sect:combinatorics-emo}.
They are directly adapted from the basic properties of an Erd\H{o}s-Moser condition
proven in~\cite{Patey2015Degrees}.
The first lemma states that each element of the finite transitive tournaments $\vec{F}$ behaves
uniformly with respect to the elements of the reservoir, that is, is beaten by every element
of the reservoir or beats all of them.

\begin{lemma}\label{lem:em-comp-reduc-uniform-behaviour}
For every condition~$c = (\alpha, \vec{F}, \Ccal)$,
every $Z_0 \oplus \dots \oplus Z_{k-1} \in \Ccal$, 
every part $\nu$ of $\Ccal$, every $i < |\alpha|$ and every $x \in F_\nu \setminus [0, \alpha(i))$, 
either $\{x\} \to_{R_i} Z_\nu$ or $Z_\nu \to_{R_i} \{x\}$.
\end{lemma}
\begin{proof}
By property (c) of the condition~$c$, there exists a minimal $R_i$-interval
$(u, v)$ of $F_\nu \setminus [0, \alpha(i))$ containing $Z_\nu$.
Here, $u$ and $v$ may be respectively $-\infty$ and $+\infty$.
By definition of an interval, $\{u\} \to_{R_i} Z_\nu \to_{R_i} \{v\}$.
By definition of a minimal interval, $R_i(x, u)$ or $R_i(v, x)$ holds.
Suppose the former holds. By transitivity of $F_\nu \setminus [0, \alpha(i))$,
for every $y \in Z_\nu$, $R_i(x, y)$ holds, since both $R_i(x, u)$ and~$R_i(u, y)$ hold. 
Therefore $\{x\} \to_{R_i} Z_\nu$. In the latter case, by symmetry, $Z_\nu \to_{R_i} \{x\}$.
\end{proof}

The second lemma is the core of the combinatorics of the Erd\H{o}s-Moser theorem. It provides
sufficient properties to obtain a valid extension of a condition. Properties (i) and (ii)
are simply the definition of an extension. Properties (iii) and (iv) help propagating
properties (b) and (c) from a condition to its extension. We shall see empirically that 
properties (iii) and (iv) are simpler to check than (b) and (c), 
as the former properties match exactly the way we add elements to our finite tournaments $\vec{F}$. 
Therefore, ensuring that these properties
are satisfied usually consists of checking that we followed the standard process of adding elements
to~$\vec{F}$. 

\begin{lemma}\label{lem:em-comp-reduc-sufficient-cond-ext}
Fix a condition~$c = (\alpha, \vec{F}, \Ccal)$ where $\Ccal$ is a $k$-cover class of~$[t, +\infty)$. 
Let $E_0, \dots, E_{m-1}$ be finite sets, $\Dcal$ be a non-empty $\Pi^{0,C}_1$ $m$-cover class of $[t', +\infty)$
for some~$t' \geq t$ and $f : m \to k$ be a function such that for each~$i < |\alpha|$ and $\nu < m$,
\begin{itemize}
	\item[(iii)] $E_\nu$ is $R_i$-transitive
	\item[(iv)] $V_\nu \to_{R_i} E_\nu$ or $E_\nu \to_{R_i} V_\nu$ for each $V_0 \oplus \dots \oplus V_{m-1} \in \Dcal$
\end{itemize}
Set $H_\nu = F_{f(\nu)} \cup E_\nu$ for each $\nu < m$.
If properties (i) and (ii) of an extension are satisfied for~$d = (\alpha, \vec{H}, \Dcal)$ with witness $f$,
then~$d$ is a valid condition extending~$c$.
\end{lemma}
\begin{proof}
All we need is to check properties (b) and (c) for~$d$ in the definition of a condition.
We prove property (b). Fix an $i < |\alpha|$, some part $\nu$ of $\Dcal$, and an $x \in V_\nu$
for some $V_0 \oplus \dots \oplus V_{m-1} \in \Dcal$. In order to prove that 
$(F_{f(\nu)} \cup E_\nu) \setminus [0, \alpha(i)) \cup \{x\}$
is $R_i$-transitive, it is sufficient to check that the set contains no 3-cycle.
Fix three elements $u < v < w \in (F_{f(\nu)} \cup E_\nu) \setminus [0, \alpha(i)) \cup \{x\}$.
\begin{itemize}
	\item Case 1: $\{u, v, w\} \cap F_{f(\nu)} \setminus [0, \alpha(i)) \neq \emptyset$. 
	Then $u \in F_{f(\nu)} \setminus [0, \alpha(i))$ as $F_{f(\nu)} < E_\nu < \{x\}$ and $u < v < w$.
	By property (ii), there is some $Z_0 \oplus \dots \oplus Z_{k-1} \in \Ccal$ such that $E_\nu \cup \{x\} \subseteq Z_{f(\nu)}$.
	If $v \in F_{f(\nu)}$, then by property (b) of the condition~$c$ on~$Z_{f(\nu)}$, $\{u, v, w\}$ is $R_i$-transitive.
	If $v \not \in F$, then by Lemma~\ref{lem:em-comp-reduc-uniform-behaviour}, $\{u\} \to_{R_i} Z_{f(\nu)}$
	or $Z_{f(\nu)} \to_{R_i} \{u\}$, so $\{u, v, w\}$ is $R_i$-transitive since~$v, w \in Z_{f(\nu)}$.

	\item Case 2: $\{u, v, w\} \cap  F_{f(\nu)} \setminus [0, \alpha(i)) = \emptyset$. 
	Then at least $u, v \in E_\nu$ because $E_\nu < \{x\}$.
	If $w \in E_\nu$ then $\{u, v, w\}$ is $R_i$-transitive by $R_i$-transitivity of $E_\nu$.
	In the other case, $w = x \in V_\nu$. As $E_\nu \to_{R_i} V_\nu$ or $V_\nu \to_{R_i} E_\nu$,
	$\{u, v\} \to_{R_i} \{w\}$ or $\{w\} \to_{R_i} \{u, v\}$ and $\{u, v, w\}$ is $R_i$-transitive.
\end{itemize}

We now prove property (c) for $d$. Fix some $V_0 \oplus \dots \oplus V_{m-1} \in \Dcal$, 
some part $\nu$ of~$\Dcal$ and some $i < |\alpha|$.
By property (ii), there is some~$Z_0 \oplus \dots \oplus Z_{k-1} \in \Ccal$ such that $E_\nu \cup V_\nu \subseteq Z_{f(\nu)}$.
By property (c) of the condition~$c$, $Z_{f(\nu)}$ (and so $V_\nu$) is included in a minimal $R_i$-interval $(u, v)$ of 
$F_{f(\nu)} \setminus [0, \alpha(i))$.
Here again, $u$ and $v$ may be respectively $-\infty$ and $+\infty$. 
By assumption, either $E_\nu \to_{R_i} V_\nu$ or $V_\nu \to_{R_i} E_\nu$. As $E_\nu$ is a finite $R_i$-transitive set,
it has a minimal and a maximal element, say~$x$ and~$y$. If $E_\nu \to_{R_i} V_\nu$
then $V_\nu$ is included in the $R_i$-interval $(y, v)$.
Symmetrically, if $V_\nu \to_{R_i} E_\nu$ then 
$V_\nu$ is included in the $R_i$-interval $(u, x)$.
To prove minimality for the first case, assume that some $w$ is in the interval $(y, v)$.
Then $w \not \in F_{f(\nu)} \setminus [0, \alpha(i))$ by minimality of the interval $(u, v)$ with respect to 
$F_{f(\nu)} \setminus [0, \alpha(i))$, and $w \not \in E_\nu$ by maximality of~$y$.
Minimality for the second case holds by symmetry.
\end{proof}

Now we have settled the necessary technical lemmas, we start proving
lemmas which will be directly involved in the construction of the transitive subtournament.
The following simple progress lemma states that we can always find an extension of a condition
in which we increase both the finite approximations corresponding to the acceptable parts 
and the number of tournaments for which we are transitive simultaneously.
Moreover, this extension can be found uniformly.

\begin{lemma}[Progress]\label{lem:em-comp-reduc-ext}
For every condition~$c = (\alpha, \vec{F}, \Ccal)$ and every $s \in \omega$,
there exists an extension $d = (\beta, \vec{E}, \Dcal)$ such that $|\beta| \geq s$ and
$|E_\nu| \geq s$ for every acceptable part $\nu$ of~$\Dcal$.
Furthermore, such an extension can be found $C'$-effectively, uniformly in~$c$ and~$s$.
\end{lemma}
\begin{proof}
Fix a condition $c = (\alpha, \vec{F}, \Ccal)$.
First note that for every $\beta \succeq \alpha$ such that $\beta(i) > max(F_\nu : \nu < parts(\Ccal))$
whenever $|\alpha| \leq i < |\beta|$, $(\beta, \vec{F}, \Ccal)$ is a condition extending~$c$.
Therefore it suffices to prove that for every such condition~$c$ and every part $\nu$ of $\Ccal$,
we can $C'$-effectively find a condition~$d = (\alpha, \vec{H}, \Dcal)$ refining~$c$
with witness~$f : parts(\Dcal) \to parts(\Ccal)$ such that $f$ forks only parts refining part $\nu$ of $\Ccal$,
and either every such part $\mu$ of $\Dcal$ is empty or $|H_\mu| > |F_\nu|$.
Iterating the process finitely many times enables us to conclude.

Fix some part $\nu$ of $\Ccal$ and let~$\Dcal$ be the collection of $Z_0 \oplus \dots \oplus Z_{k-1} \in \Ccal$
such that $Z_\nu = \emptyset$. We can $C'$-decide whether or not $\Dcal$ is empty.
If $\Dcal$ is non-empty, then $(\alpha, \vec{F}, \Dcal)$ is a valid extension of~$c$
with the identity function as witness and such that part $\nu$ of $\Dcal$ is empty.
If $\Dcal$ is empty, we can $C'$-computably find some $Z_0 \oplus \dots \oplus Z_{k-1} \in \Ccal$
and pick some~$x \in Z_\nu$.
Consider the $C$-computable $2^{|\alpha|}$-partition $(X_\rho : \rho \in 2^{|\alpha|})$ of $\omega$ defined by
$$
X_\rho = \{ y \in \omega : (\forall i < |\alpha|)[R_i(y, x) \leftrightarrow \rho(i) = 1] \} 
$$
Let $\tilde{\Dcal}$ be the cover class refining $\Ccal^{[\nu, x]}$ such that
part $\nu$ of $\tilde{\Dcal}$ has $2^{|\alpha|}$ forks induced by the
$2^{|\alpha|}$-partition~$\vec{X}$. Define $\vec{H}$ by
$H_\mu = F_\mu$ if $\mu$ refines a part different from $\nu$,
and $H_\mu = F_\nu \cup \{x\}$ if $\mu$ refines part $\nu$ of~$\Ccal$.
The forking according to~$\vec{X}$ ensures that property (iv) of Lemma~\ref{lem:em-comp-reduc-sufficient-cond-ext} holds.
By Lemma~\ref{lem:em-comp-reduc-sufficient-cond-ext}, $d = (\alpha, \vec{H}, \tilde{\Dcal})$ is a valid extension of~$c$.
\end{proof}

\subsection{The strategy}

Thanks to Lemma~\ref{lem:em-comp-reduc-ext}, we can define an infinite, $C'$-computable
decreasing sequence of conditions $(\varepsilon, \emptyset, \{\omega\}) \geq c_0 \geq c_1 \geq \dots$
such that for each~$s \in \omega$, 
\begin{itemize}
	\item[1.] $|\alpha_s| \geq s$.
	\item[2.] $|F_{s, \nu}| \geq s$ for each acceptable part~$\nu$ of~$\Ccal_s$
\end{itemize}
where $c_s = (\alpha_s, \vec{F}_s, \Ccal_s)$.
As already noticed, if some acceptable part $\mu$ of $\Ccal_{s+1}$ refines some part $\nu$ of~$\Ccal_s$,
part $\nu$ of~$\Ccal_s$ is also acceptable.
Therefore, the set of acceptable parts forms an infinite, finitely branching $C'$-computable tree~$\Tcal$.
Let $P$ be any infinite path through~$\Tcal$. 
The set $H(P) = (\bigcup_s F_{s, P(s)})$ is infinite,
and $H(P) \setminus [0, \alpha_{i+1}(i))$ is $R_i$-transitive for each $i \in \omega$.

Our goal is to build a $C'$-computable function dominating every function computed
by $H(P)$ for at least one path $P$ trough~$\Tcal$. However, it requires too much
computational power to distinguish acceptable parts from non-acceptable ones,
and even some acceptable part may have only finitely many extensions. Therefore,
we will dominate the functions computed by~$H(P)$ for \emph{every} path $P$ trough~$\Tcal$.
 
At a finite stage, a condition contains finitely many parts, each one representing
the construction of a transitive subtournament.
As in the construction of a cohesive set, it suffices to check one by one whether 
there exists an extension of our subtournaments which will
cause a given functional to terminate at a given input.
In the next subsection, we develop the framework necessary to decide such a termination
at a finite stage.

\subsection{Forcing relation}

As a condition $c = (\alpha, \vec{F}, \Ccal)$ corresponds to the construction of multiple
subtournaments $F_0, F_1, \dots$ at the same time, the forcing relation will depend on which
subtournament we are considering. In other words, the forcing relation depends on the part $\nu$ of~$\Ccal$
we focus on.

\begin{definition}\label{def:em-comp-reduc-forcing-relation}
Fix a condition $c = (\alpha, \vec{F}, \Ccal)$, a part~$\nu$ of~$\Ccal$ and two integers~$e$, $x$.
\begin{itemize}
	\item[1.] $c \Vdash_\nu \Phi_e^{G \oplus C}(x) \uparrow$ if $\Phi_e^{(F_\nu \cup F_1) \oplus C}(x) \uparrow$
	for all $Z_0 \oplus \dots \oplus Z_{k-1} \in \Ccal$ and all subsets $F_1 \subseteq Z_\nu$
	such that $F_1$ is $R_i$-transitive simultaneously for each $i < |\alpha|$.
	\item[2.] $c \Vdash_\nu \Phi_e^{G \oplus C}(x) \downarrow$ if $\Phi_e^{F_\nu \oplus C}(x) \downarrow$.
\end{itemize}
\end{definition}

The forcing relations defined above satisfy the usual forcing properties.
In particular, let $c_0 \geq c_1 \geq \dots$ be an infinite decreasing
sequence of conditions. This sequence induces an infinite, finitely branching tree of acceptable parts~$\Tcal$.
Let~$P$ be an infinite path trough~$\Tcal$. If 
$c_s \Vdash_{P(s)} \Phi_e^{G \oplus C}(x) \uparrow$ (resp. $c_s \Vdash_{P(s)} \Phi_e^{G \oplus C}(x) \downarrow$)
at some stage~$s$, then $\Phi_e^{H(P) \oplus C}(x) \uparrow$ (resp. $\Phi_e^{H(P) \oplus C}(x) \downarrow$).

Another important feature of this forcing relation is that we can decide $C'$-uniformly in its parameters
whether there is an extension forcing~$\Phi^{G \oplus C}_e(x)$ to halt or to diverge. 
Deciding this relation with little computational power is useful because our $C'$-computable dominating function will
need to decide the termination $\Gamma^{G \oplus C}(x)$ to check whether it has to dominate the value 
outputted by $\Gamma^{G \oplus C}(x)$.

\begin{lemma}\label{lem:em-comp-reduc-force-dense}
For every condition~$c = (\alpha, \vec{F}, \Ccal)$ and every pair of integers $e, x \in \omega$,
there exists an extension~$d = (\alpha, \vec{H}, \Dcal)$ such that for each part~$\nu$ of~$\Dcal$
$$
d \Vdash_\nu \Phi_e^{G \oplus C}(x) \uparrow \hspace{10pt} \vee \hspace{10pt} d \Vdash_\nu \Phi_e^{G \oplus C}(x) \downarrow
$$
Furthermore, such an extension can be found $C'$-effectively, uniformly in~$c$, $e$ and~$x$.
\end{lemma}
\begin{proof}
Given a condition~$c$ and two integers $e, x \in \omega$,
let $I_{e,x}(c)$ be the set of parts $\nu$ of~$c$
such that $c \not \Vdash_\nu \Phi_e^{G \oplus C}(x) \downarrow$ and $c \not \Vdash_\nu \Phi_e^{G \oplus C}(x) \uparrow$.
Note that $I_{e,x}(c)$ is $C'$-computable uniformly in~$c$, $e$ and~$x$.
It suffices to prove that given such a condition~$c$ and a part~$\nu \in I_{e,x}(c)$, one can $C'$-effectively
find an extension~$d$ with witness $f$ such that $f(I_{e,x}(d)) \subseteq I_{e,x}(c) \setminus \{\nu\}$.
Applying iteratively the operation enables us to conclude.

Fix a condition~$c = (\alpha, \vec{F}, \Ccal)$ where $\Ccal$ is a $k$-cover class, and fix some part~$\nu \in I_{e,x}(c)$.
The strategy is the following: either we can fork part~$\nu$ of $\Ccal$ into enough parts so that we 
force~$\Phi_e^{G \oplus C}(x)$ to diverge
on each forked part, or we can find an extension forcing $\Phi_e^{G \oplus C}(x)$ to converge on part~$\nu$ without forking.
Hence, we ask the following question.

\smallskip
{\itshape
Q2: Is it true that for every $k$-cover~$Z_0 \oplus \dots \oplus Z_{k-1} \in \Ccal$,
for every~$2^{|\alpha|}$-partition $\bigcup_{\rho \in 2^\alpha} X_\rho = Z_\nu$,
there is some~$\rho \in 2^{|\alpha|}$ and some finite set~$F_1$ which is $R_i$-transitive
for each~$i < |\alpha|$ simultaneously, and such that~$\Phi_e^{(F_\nu \cup F_1) \oplus C}(x) \downarrow$?
}
\smallskip

If the answer is no, then by forking the part~$\nu$ of~$\Ccal$ into $2^{|\alpha|}$ parts,
we will be able to force~$\Phi_e^{G \oplus C}(x)$ to diverge.
Let~$m = k+2^{|\alpha|}-1$ and define the function~$f : m \to k$ by $f(\mu) = \mu$ 
if $\mu < k$ and $f(\mu) = \nu$ otherwise.
Let~$\Dcal$ be the collection of all~$m$-covers $V_0 \oplus \dots \oplus V_{m-1}$ which $f$-refine
some $Z_0 \oplus \dots \oplus Z_{k-1} \in \Ccal$ and such that for every part $\mu$ of $\Dcal$ 
$f$-refining part~$\nu$ of $\Ccal$ and every subset $F_1 \subseteq V_\mu$
which is $R_i$-transitive simultaneously for each~$i < |\alpha|$, $\Phi_e^{F_\nu \cup F_1}(x) \uparrow$.
Note that~$\Dcal$ is a $\Pi^{0,C}_1$ $m$-cover class $f$-refining $\Ccal$. Moreover~$\Dcal$
is non-empty since the answer to~{\itshape Q2} is no.
Let~$\vec{E}$ be defined by $E_\mu = F_\mu$ if $\mu < k$
and $E_\mu = F_\nu$ otherwise. The condition~$d = (\alpha, \vec{E}, \Dcal)$ extends~$c$
with witness~$f$. For every part~$\mu$ of $\Dcal$ $f$-refining part $\nu$ of $\Ccal$, $d \Vdash_\mu \Phi_e^{G \oplus C}(x) \uparrow$,
therefore $f(I_{e,x}(d)) \subseteq I_{e,x}(c) \setminus \{\nu\}$.

Suppose now that the answer is yes. By compactness, we can $C'$-effectively find a finite set~$E \subseteq Z_\nu$
for some~$Z_0 \oplus \dots \oplus Z_{k-1} \in \Ccal$ such that for every $2^{|\alpha|}$-partition $(E_\rho : \rho \in 2^{|\alpha|})$
of $E$, there is some $\rho \in 2^{|\alpha|}$ and some set $F_1 \subseteq E_\rho$ which is $R_i$-transitive
simultaneously for each $i < |\alpha|$ and such that $\Phi_e^{(F_\nu \cup F_1) \oplus C}(x) \downarrow$.
There are finitely many $2^{|\alpha|}$-partitions of $E$. Let~$n$ be the number of such partitions. 
These partitions induce a finite $C$-computable $n$-partition of~$dom(\Ccal)$
defined for each $(E_\rho : \rho \in 2^{|\alpha|})$ by
$$
X_{\tuple{E_\rho : \rho \in 2^{|\alpha|}}} = \left\{ y \in dom(\Ccal) : (\forall i < |\alpha|) 
	\cond{
	\mbox{ if } \rho(i) = 0 \mbox{ then } E_\rho \to_{R_i} \{y\} \\ 
	\mbox{ if } \rho(i) = 1 \mbox{ then } \{y\} \to_{R_i} E_\rho} \right\}
$$

Let~$\tilde{\Dcal}$ be the $\Pi^{0,C}_1$ $(k+n-1)$-cover class refining $\Ccal^{[\nu, E]}$
and such that part~$\nu$ of~$\Ccal^{[\nu, E]}$ is refined accordingly to the above partition of~$dom(\Ccal)$.
Let~$f : k+n-1 \to k$ be the refining function witnessing it. 
Define~$\vec{H}$ as follows. For every part~$\mu$ of $\Dcal$, refining part~$\nu$ of $\Ccal^{[\nu, E]}$,
by definition of~$\tilde{\Dcal}$, there is some~$2^{|\alpha|}$-partition $\tuple{E_\rho : \rho \in 2^{|\alpha|}}$ of~$E$
such that for every $V_0 \oplus \dots V_{k+n-2} \in \tilde{\Dcal}$, $V_\mu \subseteq X_{\tuple{E_\rho : \rho \in 2^{|\alpha|}}}$.
By choice of~$E$, there exists some set $F_1 \subseteq E_\rho$ for some~$\rho \in 2^{|\alpha|}$
which is $R_i$-transitive simultaneously for each $i < |\alpha|$ and such that
$\Phi_e^{(F_\nu \cup F_1) \oplus C}(x) \downarrow$.
This set $F_1$ can be found $C'$-effectively. Set $H_\mu = F_\nu \cup F_1$.
For every part~$\mu$ of $\tilde{\Dcal}$ which refines some part~$\xi$ of $\Ccal^{[\nu, E]}$ different from~$\nu$,
set~$H_\mu = F_\xi$.
By Lemma~\ref{lem:em-comp-reduc-sufficient-cond-ext}, $d = (\alpha, \vec{H}, \tilde{\Dcal})$ is a valid condition
extending~$c$. Moreover, for every part $\mu$ of~$\tilde{\Dcal}$ refining part~$\nu$ of $\Ccal$,
$d \Vdash_\mu \Phi_e^{G \oplus C}(x) \downarrow$. Therefore $f(I_{e,x}(d)) \subseteq I_{e,x}(c) \setminus \{\nu\}$.
\end{proof}

\subsection{Construction}

We are now ready to construct our infinite transitive subtournament~$H(P)$ together
with a $C'$-computable function~$f$ dominating every~$H(P) \oplus C$-computable function.
Thanks to Lemma~\ref{lem:em-comp-reduc-ext} and Lemma~\ref{lem:em-comp-reduc-force-dense}, we can $C'$-compute an infinite
descending sequence of conditions $(\epsilon, \emptyset, 1^{<\omega}) \geq c_0 \geq c_1 \geq \dots$
such that at each stage $s \in \omega$,
\begin{itemize}
	\item[1.] $|\alpha_s| \geq s$
	\item[2.] $|F_{s, \nu}| \geq s$ for each acceptable part~$\nu$ of~$\Ccal_s$
	\item[3.] $c_s \Vdash_\nu \Phi_e^{G \oplus C}(x) \downarrow$ or $c_s \Vdash_\nu \Phi_e^{G \oplus C}(x) \uparrow$
	for each part~$\nu$ of~$\Ccal_s$ if $\tuple{e, x} = s$
\end{itemize}
where $c_s = (\alpha_s, \vec{F}_s, \Ccal_s)$.
Property 1 ensures that the resulting set with be eventually transitive
for every tournament in~$\vec{R}$. Property~2 makes the subtournaments infinite.
Last, property 3 enables us to $C'$-decide at a finite stage whether a functional terminates on a given
input, with the transitive subtournament as an oracle.

Define the $C'$-computable function $f : \omega \to \omega$ as follows:
On input~$x$, the function~$f$ looks at all stages~$s$ such that $s = \tuple{e,x}$ for
some $e \leq x$. For each such stage~$s$, and each part~$\nu$ in~$\Ccal_s$,
the function $C'$-decides whether $c_s \Vdash_\nu \Phi^{G \oplus C}_e(x) \downarrow$
or $c_s \Vdash_\nu \Phi^{G \oplus C}_e(x) \uparrow$. 
In the first case, $f$ computes the value $\Phi^{F_{s, \nu} \oplus C}_e(x)$.
Having done all that, $f$ returns a value greater than the maximum of the computed values.

Fix any infinite path~$P$ trough the infinite tree $\Tcal$ of the acceptable parts induced
by the infinite descending sequence of conditions. 
We claim that $f$ dominates every function computed by~$H(P) \oplus C$.
Fix any Turing index $e \in \omega$ such that $\Phi_e^{H(P) \oplus C}$ is total.
Consider any input~$x \geq e$ and the corresponding stage $s = \tuple{e,x}$. 
As $\Phi_e^{H(P) \oplus C}$ is total, $c_s \not \Vdash_{P(s)} \Phi_e^{G \oplus C}(x) \uparrow$,
hence by property 3, $c_s \Vdash_{P(s)} \Phi_e^{G \oplus C}(x) \downarrow$.
By construction, $f(x)$ computes the value of $\Phi_e^{F_{s,P(s)} \oplus C}(x)$ and returns
a greater value. As $F_{s,P(s)}$ is an initial segment of $H(P)$, 
$\Phi_e^{F_{s,P(s)} \oplus C}(x) = \Phi_e^{H(P) \oplus C}(x)$
and therefore $f(x) > \Phi_e^{H(P) \oplus C}(x)$.
This completes the proof of~$\amt \not \leq_c \emo$.

\chapter{A Ramsey-type König's lemma}\label{chap:ramsey-type-konig-lemma}

This chapter is a joint work with Laurent Bienvenu and Paul Shafer.

Weak K\"onig's lemma informally captures compactness arguments~\cite{Simpson2009Subsystems}.
It is involved in many constructions of solutions to Ramsey-type statements, e.g.,
cone avoidance~\cite{Seetapun1995strength,Wang2014Some} or control of the jump~\cite{Cholak2001strength,Patey2016Controlling}.
It is natural to wonder whether the use of compactness is really necessary
to prove Ramsey's theorem. The question of whether $\rt^2_2$ implies~$\wkl$ over~$\rca$
has been open for two decades, until Liu~\cite{Liu2012RT22} solved it by proving that PA degrees
are not a combinatorial consequence of~$\rt^1_2$.

Recently, Flood~\cite{Flood2012Reverse} clarified the relation between Ramsey-type theorems
and $\wkl$, by introducing a Ramsey-type variant of weak K\"onig's lemma ($\rwkl$).
Informally, seing a set as a 2-coloring of the integers, 
for every $\Pi^0_1$ class of 2-colorings, $\rwkl$ states
the existence of an infinite set homogeneous for one of them.
The exact statement of~$\rwkl$ has to be done with some care, as we do not want 
to state the existence of a member of the $\Pi^0_1$ class.

\index{Ramsey-type weak K\"onig's lemma}
\index{rwkl@$\rwkl$|see {Ramsey-type weak K\"onig's lemma}}
\index{homogeneous!for a tree}
\begin{definition}[Ramsey-type weak K\"onig's lemma]\label{def-RWKL}
A set $H \subseteq \Nb$ is \emph{homogeneous} for a $\sigma \in 2^{<\Nb}$ if $(\exists c < 2)(\forall i \in H)(i < |\sigma| \imp \sigma(i)=c)$, and a set $H \subseteq \Nb$ is \emph{homogeneous} for an infinite tree $T \subseteq 2^{<\Nb}$ if the tree $\{\sigma \in T : \text{$H$ is homogeneous for $\sigma$}\}$ is infinite.  $\rwkl$ is the statement ``for every infinite subtree of $2^{<\Nb}$, there is an infinite homogeneous set.''
\end{definition}

Flood~\cite{Flood2012Reverse} proved that~$\rwkl$ is a strict consequence of both~$\srt^2_2$ and~$\wkl$.
The Ramsey-type weak K\"onig's lemma is in fact a consequence of the stable Erd\H{o}s-Moser theorem
over~$\rca$, as proven independently by Flood and Towsner~\cite{Flood2014Separating}
and Bienvenu, Patey and Shafer~\cite{Bienvenu2015logical}.

\begin{theorem}\label{thm:sem-imp-rwkl}
$\rca \vdash \semo \imp \rwkl$.
\end{theorem}
\begin{proof}
Let $T \subseteq 2^{<\Nb}$ be an infinite tree.  For each $s \in \Nb$, let $\sigma_s$ be the leftmost element of~$T^s$.  We define a tournament $R$ from the tree $T$.  For $x < s$, if $\sigma_s(x)=1$, then $R(x,s)$ holds and $R(s,x)$ fails; otherwise, if $\sigma_s(x) = 0$, then $R(x,s)$ fails and $R(s,x)$ holds. This tournament $R$ is essentially the same as the coloring $f(x,s)=\sigma_s(x)$ defined by Flood in his proof that $\rca \vdash \srt^2_2 \imp \rwkl$ (\cite{Flood2012Reverse}~Theorem~5), where he showed that $f$ is stable.  By the same argument, $R$ is stable. 

Apply $\semo$ to $R$ to get an infinite transitive sub-tournament $U$.  Say that a $\tau \in U^{<\Nb}$ satisfies~$(\star)$ if $\range(\tau)$ is not homogenous for $T$ with color $1$ and $(\forall k < |\tau|)R(\tau(k),\tau(k+1))$.  Consider a hypothetical $\tau \in U^{<\Nb}$ satisfying $(\star)$.  There must be a $k < |\tau|$ such that $R(s,\tau(k))$ for cofinitely many $s$.  This is because otherwise there would be infinitely many $s$ such that $(\forall k < |\tau|)R(\tau(k),s)$ and hence infinitely many $s$ for which $\range(\tau)$ is homogeneous for $\sigma_s$ with color $1$, contradicting that $\range(\tau)$ is not homogeneous for $T$ with color $1$.  From the facts that $R(s,\tau(k))$ for cofinitely many $s$, that $(\forall k < |\tau|)R(\tau(k),\tau(k+1))$, and that $U$ is transitive, we conclude that $R(s,\tau(|\tau|-1))$ for cofinitely many $s$.

The proof now breaks into two cases.  First, suppose that the $\tau(|\tau|-1)$ for the $\tau \in U^{<\Nb}$ satisfying $(\star)$ are unbounded.  Then, because $(\star)$ is a $\Sigma^0_1$ property of $U$, there is an infinite set $X$ consisting of numbers of the form $\tau(|\tau|-1)$ for $\tau \in U^{<\Nb}$ satisfying $(\star)$.  As argued above, every $x \in X$ satisfies $R(s,x)$ for cofinitely many $s$.  Thus we can thin out $X$ to an infinite set $H$ such that $(\forall x,y \in H)(x < y \imp R(y,x))$. Thus~$H$ is homogeneous for $T$ with color $0$ because $H$ is homogeneous for $\sigma_y$ with color $0$ for every $y \in H$.

Second, suppose that the $\tau(|\tau|-1)$ for the $\tau \in U^{<\Nb}$ satisfying $(\star)$ are bounded, say by $m$.  Then $H = U \setminus \{0,1,\dots,m\}$ is homogeneous for $T$ with color $1$.  To see this, suppose not.  Then there is a finite $V \subseteq H$ that is not homogeneous for $T$ with color $1$.  Let $\tau \in V^{<\Nb}$ be the enumeration of $V$ in the order given by $R$:  $(\forall k < |\tau|)R(\tau(k),\tau(k+1))$.  Then $\tau$ satisfies $(\star)$, but $\tau(|\tau|-1) > m$.  This is a contradiction.
\end{proof}

\section{A Ramsey-type weak weak K\"onig's lemma}

Flood~\cite{Flood2012Reverse} proved that~$\rwkl$ implies~$\dnr$ over~$\rca$ and asked whether the implication is strict.
We shall prove in section~\ref{sect:strength-ramsey-strength-rwkl} that it is the case.
In this section, we clarify the relation between the Ramsey-type weak K\"onig's lemma and~$\dnr$.
Just as $\wkl$ can be weakened to $\wwkl$ by restricting to trees of positive measure, so can $\rwkl$ be weakened to $\rwwkl$ by restricting to trees of positive measure.

\index{Ramsey-type weak weak K\"onig's lemma}
\index{rwwkl@$\rwwkl$|see {Ramsey-type weak weak K\"onig's lemma}}
\begin{definition} \label{def:RWWKL}
$\rwwkl$ is the statement ``for every subtree of $2^{<\Nb}$ with positive measure, there is an infinite homogeneous set.''
\end{definition}

Applying $\rwwkl$ to a tree in which every path is Martin-L\"of random yields an infinite subset of a Martin-L\"of random set, and every infinite subset of every Martin-L\"of random set computes a $\dnrf$ function.  In fact, computing an infinite subset of a Martin-L\"of random set is equivalent to computing a $\dnr$ function, as the following theorem states.

\begin{theorem}[Kjos-Hanssen~\cite{Kjos-Hanssen2009Infinite}, Greenberg and Miller~\cite{Greenberg2009Lowness}]\label{thm:miller-dnc}
For every $A \in 2^\omega$, $A$ computes a $\dnrf$ function if and only if $A$ computes an infinite subset of a Martin-L\"of random set.
\end{theorem}

Theorem~\ref{thm:miller-dnc} also relativizes:  a set $A$ computes a $\dnrf(X)$ function if and only if it computes an infinite subset of a set that is Martin-L\"of random relative to $X$.  Thus one reasonably expects that $\dnr$ and $\rwwkl$ are equivalent over $\rca$.  This is indeed the case, as we show.  The proof makes use of the following computability-theoretic lemma, which reflects a classical fact concerning diagonally non-computable functions.

\begin{lemma}\label{lemma:dnr-rwwkl_helper}
The statement ``for every set $X$ there is a function $g \colon \Nb^3 \imp \Nb$ such that $\forall e,k,n(g(e,k,n) > n \andd (|W_e^X| < k \imp g(e,k,n) \notin W_e^X))$'' is provable in $\rca + \dnr$.
\end{lemma}

\begin{proof}
Fix a sequence of functions $(b_k)_{k \in \Nb}$ such that, for each $k \in \Nb$, $b_k$ maps $\Nb$ onto $\Nb^k$ in such a way that $b_k^{-1}(\vec x)$ is infinite for every $\vec x \in \Nb^k$.  Let $c \colon \Nb \imp \Nb$ be a function such that, for all $e,i,k,x \in \Nb$, $\Phi_{c(e,i,k)}^X(x) = b_k(y)(i)$ for the $(i+1)$\textsuperscript{th} number $y$ enumerated in $W_e^X$ if $|W_e^X| \geq i+1$; and $\Phi_{c(e,i,k)}^X(x)\ua$ otherwise.  Let $f$ be diagonally non-computable relative to $X$.  Define $g$ by letting $g(e,k,n)$ be the least $x > n$ such that $b_k(x) = \la f(c(e,0,k)), f(c(e,1,k)), \dots, f(c(e,k-1,k)) \ra$.  Suppose for the sake of contradiction that $|W_e^X| < k$ but that $g(e,k,n) \in W_e^X$.  Then $g(e,k,n)$ is the $(i+1)$\textsuperscript{th} number enumerated into $W_e^X$ for some $i+1 < k$.  Hence $\Phi_{c(e,i,k)}^X(c(e,i,k)) = b_k(g(e,k,n))(i)$.  However, by the definition of $g$, $b_k(g(e,k,n))(i) = f(c(e,i,k))$.  Thus $f(c(e,i,k)) = \Phi_{c(e,i,k)}^X(c(e,i,k))$, contradicting that $f$ is $\dnrf$ relative to $X$.
\end{proof}

Notice that in the statement of the above lemma, $W_e^X$ need not exist as a set.  Thus `$|W_e^X| < k$' should be interpreted as `$\forall s(|W_{e,s}^X| < k)$,' where $(W_{e,s}^X)_{s \in \Nb}$ is the standard enumeration of $W_e^X$.
The following theorem has been obtained independently by Flood and Towsner~\cite{Flood2014Separating}
and Bienvenu, Patey and Shafer~\cite{Bienvenu2015logical}.

\begin{theorem}\label{thm:dnr-rwwkl}
$\rca \vdash \dnr \biimp \rwwkl$.
\end{theorem}

\begin{proof}
The direction $\rwwkl \rightarrow \dnr$ is implicit in Flood's proof that $\rca \vdash \rwkl \imp \dnr$ (\cite{Flood2012Reverse}~Theorem~8).  Indeed, Flood's proof uses the construction of a tree of positive measure due to Jockusch~\cite{Jockusch197401}. The proof that $\dnr \imp \rwwkl$ is similar to the original proof of Theorem~\ref{thm:miller-dnc}.  However, some adjustments are needed as the original argument uses techniques from measure theory and algorithmic randomness which can only be formalized within $\wwkl$.  We instead use explicit combinatorial bounds.  

Assume $\dnr$, and consider a tree $T$ of measure $\geq 2^{-c}$ for some $c$, which we can assume to be $\geq 3$ (the reason for this assumption will become clear). For a given set~$H \subseteq \Nb$ and a value $v \in \{0,1\}$, let $\Gamma^v_H = \{\sigma \in \str : (\forall i \in H)(\sigma(i) = v)\}$, and abbreviate $\Gamma^v_{\{n\}}$ by $\Gamma^v_{n}$.  For a tree $T$ and a constant $c$, let $\bad(n,T,c)$ be the $\Sigma^0_1$ predicate `$\mu(T \cap \Gamma^0_n) <  2^{-2c}$.'  In the following claim, $\{n : \bad(n,T,c)\}$ need not \emph{a priori} exist as a set, so `$|\{n : \bad(n,T,c)\}| < 2c$' should be interpreted in the same manner as `$|W_e^X| < k$' in the statement of Lemma~\ref{lemma:dnr-rwwkl_helper}.

\begin{claim}
If $c \geq 3$ and $\mu(T) \geq 2^{-c}$, then $|\{n : \bad(n,T,c)\}| < 2c$.
\end{claim}

\begin{proof}
Suppose for the sake of contradiction that $|\{n : \bad(n,T,c)\}| \geq 2c$, and let $B$ be the first $2c$ elements enumerated in $\{n : \bad(n,T,c)\}$.  For each $n \in B$, the tree $T \cap \Gamma^0_n$ has measure $<2^{-2c}$, which implies that $(\forall n \in B)(\exists i)(|T^i \cap \Gamma^0_n| < 2^{i-2c})$ (recall that $T^i$ is the set of strings in $T$ of length $i$).  By $\bsig^0_1$, let $N_0$ be such that $(\forall n \in B)(\exists i < N_0)(|T^i \cap \Gamma^0_n| < 2^{i-2c})$, and observe that $(\forall n \in B)(\forall j > N_0)(|T^j \cap \Gamma^0_n| < 2^{j-2c})$.  Let $N = N_0 + \max(B)$.

On the one hand,
\begin{align*}
\left|T^N \cap \bigcup_{n \in B} \Gamma^0_n\right| = |T^N \setminus \Gamma^1_B| \geq  |T^N| - |\Gamma^1_B \cap \{0,1\}^N| \geq  2^{N-c} - 2^{N-2c}.
\end{align*}

On the other hand,
\begin{align*}
\left|T^N \cap \bigcup_{n \in B} \Gamma^0_n\right| = \left| \bigcup_{n \in B} T^N \cap \Gamma^0_n\right| \leq \sum_{n \in B} |T^N \cap \Gamma^0_n| \leq 2c \cdot 2^{N-2c}.
\end{align*}

Putting the two together, we get that $ 2^{N-c} - 2^{N-2c} \leq 2c \cdot 2^{N-2c}$, which is a contradiction for $c \geq 3$. 
\end{proof}

Let $g$ be as in Lemma~\ref{lemma:dnr-rwwkl_helper} for $X=T$.  Given a (canonical index for a) finite set $F$ and a $c$, we can effectively produce an index $e(F,c)$ such that $\forall n(n \in W_{e(F,c)}^T \biimp \bad(n, T \cap \Gamma^0_F,c))$.  Computably construct an increasing sequence $h_0 < h_1 < h_2 < \dots$ of numbers by letting, for each $s \in \Nb$, $H_s = \{h_i : i < s\}$ and $h_s = g(e(H_s,c \cdot 2^s),c \cdot 2^{s+1},\max(H_s \cup \{0\}))$.  Using $\Pi^0_1$-induction ($\ipi^0_1$), we prove that $\forall s(\mu(T \cap \Gamma^0_{H_s}) \geq 2^{-c \cdot 2^s})$.  For $s=0$, this is simply the assumption $\mu(T) \geq 2^{-c}$.  Assuming $\mu(T \cap \Gamma^0_{H_s}) \geq 2^{-c \cdot 2^s}$, the claim implies that $|W_{e(H_s,c \cdot 2^s)}^T| < c \cdot 2^{s+1}$.  Thus $h_s = g(e(H_s,c \cdot 2^s),c \cdot 2^{s+1},\max(H_s \cup \{0\})) \notin W_{e(H_s,c \cdot 2^s)}$, and therefore $\neg \bad(h_s,T \cap \Gamma^0_{H_s},c \cdot 2^s)$.  This means that $\mu(T \cap \Gamma^0_{H_s} \cap \Gamma^0_{h_s}) \geq 2^{-c \cdot 2^{s+1}}$, which is what we wanted because $\Gamma^0_{H_s} \cap \Gamma^0_{h_s} = \Gamma^0_{H_{s+1}}$.

Let $H = \{h_s : s \in \Nb\}$, which exists by $\Delta^0_1$ comprehension because the sequence $h_0 < h_1 < h_2 < \dots$ is increasing.  We show that $H$ is homogeneous for $T$.  Suppose for the sake of contradiction that $H$ is not homogeneous for $T$.  This means that there are only finitely many $\sigma \in T$ such that $H$ is homogeneous for $\sigma$. Therefore at some level $s$, $\{\sigma \in T^s : (\forall i \in H)(\sigma(i) = 0)\}=\emptyset$.  As $H \cap \{0,1,\dots,s\} \subseteq H_s$, we in fact have that $\{\sigma \in T^s : (\forall i \in H_s)(\sigma(i) = 0)\}=\emptyset$.  In other words, $T \cap \Gamma^0_{H_s} = \emptyset$, which contradicts $\mu(T \cap \Gamma^0_{H_s}) \geq 2^{-c \cdot 2^s}$.  Thus $H$ is homogeneous for $T$.
\end{proof}

\section{Ramsey-type statisfiability}

One can conceivably consider a Ramsey-type variant of any $\Pi^1_2$ statement $\forall X \exists Y \varphi(X,Y)$ so long as one can provide a reasonable formulation of what it means for a set $Z$ to be consistent with a $Y$ such that $\varphi(X,Y)$.  For example, in the case of $\rwkl$, we think of a set $H$ as being consistent with a path through an infinite tree $T \subseteq \str$ if $H$ is homogeneous for $T$.  We are interested in analyzing the strengths of Ramsey-type variants of statements that are equivalent to $\wkl$ over $\rca$.  Several such statements have trivial Ramsey-type variants.  For example, $\rca$ proves that for every pair of injections $f,g \colon \Nb \imp \Nb$ with disjoint ranges, there is an infinite set $X$ consistent with being a separating set for the ranges of $f$ and $g$ because $\rca$ proves that there is an infinite subset of the range of $f$.  The obvious Ramsey-type variant of Lindenbaum's lemma (every consistent set of sentences has a consistent completion) is also easily seen to be provable in $\rca$.  For the remainder of this paper, we consider non-trivial Ramsey-type variants of the compactness theorem for propositional logic and of graph coloring theorems.  Many of these variants are equivalent to $\rwkl$, which we take as evidence that $\rwkl$ is robust.

\index{finitely satisfiable}
\begin{definition}
A set $C$ of propositional formulas is \emph{finitely satisfiable} if every finite $C_0 \subseteq C$ is satisfiable (i.e., has a satisfying truth assignment).  We denote by $\sat$ the compactness theorem for propositional logic, which is the statement ``every finitely satisfiable set of propositional formulas is satisfiable.''
\end{definition}

It is well-known that $\sat$ is equivalent to $\wkl$ over $\rca$ (see \cite{Simpson2009Subsystems}~Theorem IV.3.3).

If $C$ is a set of propositional formulas, then let $\atoms(C)$ denote the set of propositional atoms appearing in the formulas in $C$.  Strictly speaking, $\rca$ does not in prove that $\atoms(C)$ exists for every set of propositional formulas $C$.  However, in $\rca$ we can rename the atoms appearing in a set of propositional formulas $C$ in such a way as to produce an equivalent set of propositional formulas $C'$ for which $\atoms(C')$ does exist.  Indeed, we may assume that $\atoms(C) = \Nb$ whenever $\atoms(C)$ is infinite.  Thus for simplicity, we always assume that $\atoms(C)$ exists as a set.

\index{homogeneous!for a set of formulas}
\begin{definition}
Let $C$ be a set of propositional formulas.  A set $H \subseteq \atoms(C)$ is \emph{homogeneous for $C$} if there is a $c \in \{\true, \false\}$ such that every finite $C_0 \subseteq C$ is satisfiable by a truth assignment $\nu$ such that $(\forall a \in H)(\nu(a) = c)$.
\end{definition}

As is typical, we identify $\true$ with $1$ and $\false$ with $0$.

\index{Ramsey-type satisfiability}
\index{rsat@$\rsat$|see {Ramsey-type satisfiability}}
\begin{definition}[Ramsey-type satisfiability]
$\rsat$ is the statement ``for every finitely satisfiable set $C$ of propositional formulas with $\atoms(C)$ infinite, there is an infinite $H \subseteq \atoms(C)$ that is homogeneous for $C$.''
\end{definition}

We first show that $\rca \vdash \rsat \imp \rwkl$.  In fact, we show that the restriction of $\rsat$ to what we call \emph{$2$-branching} clauses implies $\rwkl$ over $\rca$.  This technical restriction is useful for the proof of Theorem~\ref{thm:rcolor3-rsat} in our analysis of Ramsey-type graph coloring principles.

\index{literal}
Recall that a propositional formula $\ell$ is called a \emph{literal} if either $\ell = a$ or $\ell = \neg a$ for some propositional atom $a$ and that a \emph{clause} is a disjunction of literals.

\index{branching clause@2-branching clause}
\index{rsat2branching@$\rsat_{\textup{2-branching}}$}
\begin{definition}
Let $\{a_i : i \in \Nb\}$ be an infinite set of propositional atoms.  A set $C$ of clauses is called \emph{$2$-branching} if, for every clause $\ell_0 \orr \ell_1 \orr \cdots \orr \ell_{n-1} \in C$ and every $i < n$, the literal $\ell_i$ is either $a_i$ or $\neg a_i$.  $\rsat_{\textup{2-branching}}$ is $\rsat$ restricted to $2$-branching clauses.
\end{definition}

\begin{proposition}\label{prop:RWKL2BranchRSAT}
$\rca \vdash \rsat_{\textup{2-branching}} \imp \rwkl$.
\end{proposition}
\begin{proof}
Let $A = \{a_i : i \in \Nb\}$ be a set of propositional atoms, and to each string $\sigma \in 2^{<\Nb}$ associate the clause $\theta_\sigma = \bigvee_{i<|\sigma|}\ell_i$, where $\ell_i = a_i$ if $\sigma(i) = 0$ and $\ell_i = \neg a_i$ if $\sigma(i) = 1$.  Let $T \subseteq 2^{<\Nb}$ be an infinite tree.  Let $C = \{\theta_\sigma : \sigma \notin T\}$, and observe that $C$ is $2$-branching.  We show that $C$ is finitely satisfiable.  Given $C_0 \subseteq C$ finite, choose $n$ large enough so that the atoms appearing in the clauses in $C_0$ are among $\{a_i : i < n\}$.  As $T$ is infinite, choose a $\tau \in T$ of length $n$.  Define a truth assignment $t \colon \{a_i : i < n\} \imp \{\true, \false\}$ by $t(a_i) = \tau(i)$.  Now, if $\theta$ is a clause in $C_0$, then $\theta = \theta_\sigma = \bigvee_{i<|\sigma|}\ell_i$ for some $\sigma \notin T$ with $|\sigma| < n$.  Thus there is an $i < n$ such that $\sigma(i) \neq \tau(i)$ (because $\tau \in T$ and $\sigma \notin T$), from which we see that $t(\ell_i)=\true$ and hence that $t(\theta_\sigma)=\true$.  Thus $t$ satisfies $C_0$.

By $\rsat_{\textup{2-branching}}$, let $H_0 \subseteq A$ and $c \in \{\true,\false\}$ be such that $H_0$ is homogeneous for $C$ with truth value $c$.  Let $H = \{i \in \Nb : a_i \in H_0\}$.  We show that $H$ is homogeneous for a path through $T$ with color $c$.  Given $n \in \Nb$, we want to find a $\tau \in T$ such that $|\tau| = n$ and $(\forall i < |\tau|)(i \in H \imp \tau(i) = c)$.  Thus let $t \colon \{a_i : i < n\} \imp \{\true, \false\}$ be a truth assignment satisfying $C_0 = \{\theta_\sigma : \sigma \notin T \andd |\sigma| = n\}$ such that $(\forall a \in \{a_i : i < n\} \cap H_0)(t(a) = c)$.  Let $\tau \in 2^n$ be defined by $\tau(i) = t(a_i)$ for all $i<n$.  Notice that $(\forall i < |\tau|)(i \in H \imp \tau(i) = c)$ and that $t(\theta_\tau) = \false$.  If $\tau \notin T$, then $\theta_\tau \in C_0$, contradicting that $t$ satisfies $C_0$.  Thus $\tau \in T$ as desired.
\end{proof}

\section{Ramsey-type graph coloring principles}

\index{locally k-colorable graph@locally $k$-colorable graph}
\index{colorable graph@$k$-colorable graph}
Let $k \in \Nb$, and let $G = (V, E)$ be a graph.  A function $f \colon V \imp k$ is a \emph{$k$-coloring} of $G$ if $(\forall x,y \in V)((x,y) \in E \imp f(x) \neq f(y))$.  A graph is \emph{$k$-colorable} if it has a $k$-coloring, and a graph is \emph{locally $k$-colorable} if every finite subgraph is $k$-colorable.  A simple compactness argument proves that every locally $k$-colorable graph is $k$-colorable.  In the context of reverse mathematics, we have the following well-known equivalence.

\begin{theorem}[see~\cite{Hirst1990Marriage}]\label{thm-ColoringIsWKL}
For every $k \in \omega$ with $k \geq 2$, the following statements are equivalent over $\rca$:
\begin{itemize}
\item[(i)] $\wkl$
\item[(ii)] Every locally $k$-colorable graph is $k$-colorable.
\end{itemize}
\end{theorem}

In light of Theorem~\ref{thm-ColoringIsWKL}, we define Ramsey-type analogs of graph coloring principles and compare them to Ramsey-type weak K\"onig's lemma.

\index{homogeneous!for a graph}
\index{Ramsey-type graph coloring}
\index{rcolor@$\rcolor_k$|see {Ramsey-type graph coloring}}
\begin{definition}[Ramsey-type graph coloring]
Let $G = (V,E)$ be a graph.  A set $H \subseteq V$ is {\itshape homogeneous for $G$} if every finite $V_0 \subseteq V$ induces a subgraph that is $k$-colorable by a coloring that colors every vertex in $V_0 \cap H$ with color $0$.
$\rcolor_k$ is the statement  ``for every infinite, locally $k$-colorable graph $G = (V,E)$, 
there is an infinite $H \subseteq V$ that is homogeneous for $G$.''
\end{definition}

The goal of this section is to obtain the analog of Theorem~\ref{thm-ColoringIsWKL} with $\rwkl$ in place of $\wkl$ and with $\rcolor_k$ in place of the statement ``every locally $k$-colorable graph is $k$-colorable.''  We are able to obtain this analog for all standard $k \geq 3$ instead of all standard $k \geq 2$.  The case $k = 2$ remains open.  Showing the forward direction, that $\rca \vdash \rwkl \rightarrow \rcolor_k$ (indeed, that $\rca \vdash \rwkl \rightarrow \rcolor_k$), is straightforward.

\begin{lemma}\label{RWKLprovesLRCOLOR}
For every $k \in \omega$, $\rca \vdash \rwkl \rightarrow \rcolor_k$.
\end{lemma}

\begin{proof}
Let $G = (V,E)$ be an infinite graph such that every finite $V_0 \subseteq V$ induces a $k$-colorable subgraph.  Enumerate $V$ as $(v_i)_{i \in \Nb}$, and let $T \subseteq k^{<\Nb}$ be the tree
\begin{align*}
T = \{\sigma \in k^{<\Nb} : (\forall i,j < |\sigma|)((v_i,v_j) \in E \imp \sigma(i) \neq \sigma(j))\}.
\end{align*}
$T$ exists by $\Delta^0_1$ comprehension and is downward closed.  $T$ is infinite because for any $n \in \Nb$, any $k$-coloring of the subgraph induced by $\{v_i : i < n\}$ corresponds to a string in the tree of length $n$.  Apply $\rwkl_k$ (which follows from $\rca + \rwkl$) to $T$ to get an infinite set $H_0 \subseteq \Nb$ and a color $c < k$ such that $H_0$ is homogeneous for a path through $T$ with color $c$.  Let $H = \{v_i : i \in H_0\}$.  We show that every finite $V_0 \subseteq V$ induces a subgraph that is $k$-colorable by a coloring that colors every $v \in V_0 \cap H$ with color $0$.  Let $V_0 \subseteq V$ be finite, let $n = \max\{i+1 : v_i \in V_0\}$, and let $\sigma \in T$ be such that $|\sigma| = n$ and such that $H_0$ is homogeneous for $\sigma$ with color $c$.  Then the coloring of $V_0$ given by $v_i \mapsto \sigma(i)$ is a $k$-coloring of $V_0$ that colors the elements of $V_0 \cap H$ with color $c$.  Swapping colors $0$ and $c$ thus gives a $k$-coloring of $V_0$ that colors the elements of $V_0 \cap H$ with color $0$.
\end{proof}

We now prove that $\rca \vdash \rcolor_3 \imp \rwkl$ (Theorem~\ref{thm:rcolor3-rsat} below).  Our proof factors through the Ramsey-type satisfiability principles and is a rather elaborate exercise in circuit design.  The plan is to prove that $\rca \vdash \rcolor_3 \imp \rsat_{\textup{2-branching}}$, then appeal to Proposition~\ref{prop:RWKL2BranchRSAT}.  Given a $2$-branching set of clauses $C$, we compute a locally $3$-colorable graph $G$ such that every set homogeneous for $G$ computes a set that is homogeneous for $C$.  $G$ is built by connecting $\emph{widgets}$, which are finite graphs whose colorings have desirous properties.  A widget $W(\vec v)$ has distinguished vertices $\vec v$ through which we connect the widget to the larger graph.  These distinguished vertices can also be regarded, in a sense, as the inputs and outputs of the widget.

In an $\rcolor_3$ instance built out of widgets according to an $\rsat_{\textup{2-branching}}$ instance, some of the vertices code literals so that the colorings of these coding vertices code truth assignments of the corresponding literals in such a way that a homogeneous set for the $\rsat_{\textup{2-branching}}$ instance can be decoded from a homogeneous set for the graph that contains only coding vertices.  However, we have no control over what vertices appear in an arbitrary homogeneous set.  Therefore, we must build our graph so that the color of \emph{every} vertex gives information about the color of \emph{some} coding vertex.

When we introduce a widget, we prove a lemma concerning the three key aspects of the widget's operation:  soundness, completeness, and reversibility.  By soundness, we mean conditions on the $3$-colorings of the widget, which we think of as input-output requirements for the widget.  By completeness, we mean that the widget is indeed $3$-colorable and, moreover, that $3$-colorings of certain sub-widgets extend to $3$-colorings of the whole widget.  By reversibility, we mean that the colors of some vertices may be deduced from the colors of other vertices.

To aid the analysis of our widgets, we introduce a notation for the property that a coloring colors two vertices with the same color.

\index{$=_\nu$}
\begin{notation}
Let $G = (V,E)$ be a graph, let $a, b \in V$, and let $\nu \colon V \imp k$ be a $k$-coloring of $G$.  We write $a =_\nu b$ if $\nu(a) = \nu(b)$.
\end{notation}

The graph $G$ that we build from widgets has three distinguished vertices $0$, $1$, and $2$, connected as a triangle.  The intention of these vertices is to code truth values.  If $v$ is a vertex coding a literal $\ell$, then $(v,2)$ is an edge in $G$, and, for a $3$-coloring $\nu$, we interpret $v =_\nu 0$ as $\ell$ is false and $v =_\nu 1$ as $\ell$ is true.  Our widgets often include vertices $0$, $1$, and $2$.

\begin{widget}
$R_{\substack{x \mapsto y \\ y \mapsto z}}(a,u)$ is the following widget.
\begin{center}
\begin{tikzpicture}[x=1cm, y=1cm, node/.style={circle, draw, minimum size=2em}, widget/.style={rectangle, draw, minimum size=2em}]
	
	\node[node] (x) at (0, 0) {$x$};
	\node[node] (y) at (1, 1) {$y$};
	\node[node] (z) at (0, 2) {$z$};
	
	\node[node] (v) at (3, 1) {$v$};
	\node[node] (a) at (4, 2) {$a$};
	\node[node] (u) at (4, 0) {$u$};
		
	\draw (x) -- (y) -- (z) -- (x);	
	\draw (a) -- (u) -- (v) -- (a);

	\draw (x) -- (u);
	\draw (y) -- (v);
	\draw (z) -- (a);
		
\end{tikzpicture}
\end{center}
\end{widget}

\begin{lemma}\label{lem:RWidget}{\ }
\begin{itemize}
\item[(i)]  Let $\nu$ be a $3$-coloring of $R_{\substack{x \mapsto y \\ y \mapsto z}}(a,u)$.  If $a =_\nu x$ then $u =_\nu y$, and if $a =_\nu y$ then $u =_\nu z$.

\item[(ii)] Every $3$-coloring of the subgraph of $R_{\substack{x \mapsto y \\ y \mapsto z}}(a,u)$ induced by $\{x,y,z,a\}$ can be extended to a $3$-coloring of $R_{\substack{x \mapsto y \\ y \mapsto z}}(a,u)$.

\item[(iii)] In every $3$-coloring of $R_{\substack{x \mapsto y \\ y \mapsto z}}(a,u)$, the color of each vertex in $\{u, v\}$ determines the color of $a$.
\end{itemize}
\end{lemma}

\begin{proof}
The lemma follows from examining the two possible (up to permutations of the colors) $3$-colorings of $R_{\substack{x \mapsto y \\ y \mapsto z}}(a,u)$:
\begin{align}
a &=_\nu x & v &=_\nu z & u &=_\nu y\\
a &=_\nu y & v &=_\nu x & u &=_\nu z.
\end{align}
We see $(i)$ immediately.  For $(ii)$, if $a =_\nu x$, then color the widget according to the first coloring; and if $a =_\nu y$, then color the widget according to the second coloring.  For $(iii)$, if $u =_\nu y$ or $v =_\nu z$, then $a =_\nu x$; and if $u =_\nu z$ or $v =_\nu x$, then $a =_\nu y$.
\end{proof}

The intention is that, in $R_{\substack{x \mapsto y \\ y \mapsto z}}(a,u)$, the vertices $x$, $y$, and $z$ are some permutation of the vertices $0$, $1$, and $2$.  For example, $R_{\substack{0 \mapsto 1 \\ 1 \mapsto 2}}(a,u)$ is the instance of this widget where $x=0$, $y=1$, and $z=2$.  The notation `$R_{\substack{0 \mapsto 1 \\ 1 \mapsto 2}}(a,u)$' is evocative of Lemma~\ref{lem:RWidget}~$(i)$.  Thinking of $a$ as the widget's input and of $u$ as the widget's output, Lemma~\ref{lem:RWidget}~$(i)$ says that the widget maps $0$ to $1$ and maps $1$ to $2$.

\begin{widget}
$U_{x,y,z}(\ell,b,u)$ is the following widget.
\begin{center}
\begin{tikzpicture}[x=1cm, y=1cm, node/.style={circle, draw, minimum size=2em}, widget/.style={rectangle, draw, minimum size=2em}]
	
	\node[node] (x) at (-1, -1) {$x$};
	\node[node] (y) at (1, -1) {$y$};
	\node[node] (z) at (0, -0) {$z$};
	
	\node[node] (l) at (1, 1) {$\ell$};
	\node[node] (lb) at (-1,1) {$\bar\ell$};
	\node[node] (b) at (4, 1) {$b$};
	
	\node[widget] (R) at (1, 2) {$R_{\substack{x \mapsto y \\ y \mapsto z}}(\ell, r)$};	
	\node[node] (r) at (1, 3) {$r$};
	\node[node] (d) at (4, 2) {$d$};
	\node[node] (u) at (0, 4) {$u$};
		
	\draw (x) -- (y) -- (z) -- (x);
	\draw (z) -- (l) -- (lb) -- (z);
	\draw (lb) -- (u);
	\draw (l) -- (R) -- (r) -- (u);
	\draw (l) -- (d) -- (u);
	\draw (b) -- (d);
	
\end{tikzpicture}
\end{center}
\end{widget}

In the diagram above, the box labeled `$R_{\substack{x \mapsto y \\ y \mapsto z}}(\ell, r)$' represents an $R_{\substack{x \mapsto y \\ y \mapsto z}}(\ell, r)$ sub-widget.  The vertices $\ell$ and $r$ are the same as those appearing inside $R_{\substack{x \mapsto y \\ y \mapsto z}}(\ell, r)$.  They have been displayed to show how they connect to the rest of the $U_{x,y,z}(\ell,b,u)$ widget.  The vertices $x$, $y$, and $z$ are also the same as the corresponding vertices appearing inside $R_{\substack{x \mapsto y \\ y \mapsto z}}(\ell, r)$, and some of the edges incident to them (for example, the edge $(x, r)$) have been omitted to improve legibility.

The properties of $U_{x,y,z}(\ell,b,u)$ highlighted by the next lemmas may seem ill-motivated at first.  We explain their significance after the proofs.

\begin{lemma}\label{lem:U3Widget}{\ }
\begin{itemize}
\item[(i)]  Every $3$-coloring $\nu$ of the subgraph of $U_{x,y,z}(\ell,b,u)$ induced by $\{x,y,z,\ell,b\}$ can be extended to a $3$-coloring of $U_{x,y,z}(\ell,b,u)$.

\item[(ii)]  If $\nu$ is a $3$-coloring of $U_{x,y,z}(\ell,b,u)$ in which $\ell =_\nu x$ and $b = _\nu y$, then $u =_\nu x$.

\item[(iii)]  Every $3$-coloring $\nu$ of the subgraph of $U_{x,y,z}(\ell,b,u)$ induced by $\{x,y,z,\ell,b\}$ in which $\ell =_\nu x$ and $b \neq_\nu y$ can be extended to a $3$-coloring of $U_{x,y,z}(\ell,b,u)$ in which $u =_\nu z$.

\item[(iv)]  Every $3$-coloring $\nu$ of the subgraph of $U_{x,y,z}(\ell,b,u)$ induced by $\{x,y,z,\ell,b\}$ in which $\ell =_\nu y$ can be extended to a $3$-coloring of $U_{x,y,z}(\ell,b,u)$ in which $u =_\nu y$.
\end{itemize}
\end{lemma}

\begin{proof}
For $(i)$, let $\nu$ be a $3$-coloring of the subgraph of $U_{x,y,z}(\ell,b,u)$ induced by $\{x,y,z,\ell,b\}$.  
\begin{itemize}
\item If $\ell =_\nu x$ and $b =_\nu x$, then color the widget so that $\bar\ell =_\nu y$, $r =_\nu y$, $d =_\nu y$, and $u =_\nu z$.

\item If $\ell =_\nu x$ and $b =_\nu y$, then color the widget so that $\bar\ell =_\nu y$, $r =_\nu y$, $d =_\nu z$, and $u =_\nu x$.

\item If $\ell =_\nu x$ and $b =_\nu z$, then color the widget so that $\bar\ell =_\nu y$, $r =_\nu y$, $d =_\nu y$, and $u =_\nu z$.

\item If $\ell =_\nu y$ and $b =_\nu x$, then color the widget so that $\bar\ell =_\nu x$, $r =_\nu z$, $d =_\nu z$, and $u =_\nu y$.

\item If $\ell =_\nu y$ and $b =_\nu y$, then color the widget so that $\bar\ell =_\nu x$, $r =_\nu z$, $d =_\nu x$, and $u =_\nu y$.

\item If $\ell =_\nu y$ and $b =_\nu z$, then color the widget so that $\bar\ell =_\nu x$, $r =_\nu z$, $d =_\nu x$, and $u =_\nu y$.
\end{itemize}
In each of the above cases, the sub-widget $R_{\substack{x \mapsto y \\ y \mapsto z}}(\ell, r)$ is colored according to Lemma~\ref{lem:RWidget}.

For $(ii)$, let $\nu$ be a $3$-coloring of $U_{x,y,z}(\ell,b,u)$ in which $\ell =_\nu x$ and $b =_\nu y$.  Then it must be that $\bar\ell =_\nu y$ and $d =_\nu z$, and therefore it must be that $u =_\nu x$.

Item $(iii)$ can be seen to hold by inspecting the first and third colorings in the proof of $(i)$.  

Item $(iv)$ can be seen to hold by inspecting the last three colorings in the proof of $(i)$.
\end{proof}

\begin{lemma}\label{lem:U3WidgetDecode}
Let $\nu$ be a $3$-coloring of $U_{x,y,z}(\ell,b,u)$.  If $w$ is $\bar\ell$, $u$, or any vertex appearing in the $R_{\substack{x \mapsto y \\ y \mapsto z}}(\ell, r)$ sub-widget that is not $x$, $y$, or $z$, then the color of $w$ determines the color of $\ell$.  Moreover,
\begin{itemize}
\item if $d =_\nu x$, then $\ell =_\nu y$;
\item if $d =_\nu y$, then $\ell =_\nu x$;
\item if $d =_\nu z$, then $b \neq_\nu z$.
\end{itemize}
\end{lemma}

\begin{proof}
Let $\nu$ be a $3$-coloring of $U_{x,y,z}(\ell,b,u)$.  It is easy to see that if $\bar\ell =_\nu x$, then $\ell =_\nu y$ and that if $\bar\ell =_\nu y$, then $\ell =_\nu x$.  If $w$ is a vertex in $R_{\substack{x \mapsto y \\ y \mapsto z}}(\ell, r)$ that is not $x$, $y$, or $z$, then the color of $w$ determines the color of $\ell$ by Lemma~\ref{lem:RWidget}~$(iii)$.  For $u$, if $u =_\nu x$ or $u =_\nu z$ it cannot be that $\ell =_\nu y$ because then $\bar\ell =_\nu x$ and, by Lemma~\ref{lem:RWidget}~$(i)$, $r =_\nu z$.  On the other hand, if $u =_\nu y$, it cannot be that $\ell =_\nu x$ because then $\bar\ell =_\nu y$.  Thus if $u =_\nu x$ or $u =_\nu z$, then $\ell =_\nu x$; and if $u =_\nu y$, then $\ell =_\nu y$.  It is easy to see that if $d =_\nu x$ then $\ell =_\nu y$, that if $d =_\nu y$ then $\ell =_\nu x$, and that if $d =_\nu z$ then $b \neq_\nu z$ because $\ell$ and $b$ are neighbors of $d$.
\end{proof}

Consider a clause $\ell_0 \orr \ell_1 \orr \cdots \orr \ell_{n-1}$.  The idea is to code truth assignments that satisfy the clause as $3$-colorings of a graph constructed by chaining together widgets of the form $U_{x,y,z}(\ell_i,b,u)$.  Let $\nu$ be a $3$-coloring of $U_{x,y,z}(\ell_i,b,u)$.  The color of the vertex $\ell_i$ represents the truth value of the literal $\ell_i$:  $\ell_i =_\nu x$ is interpreted as $\ell_i$ is false, and $\ell_i =_\nu y$ is interpreted as $\ell_i$ is true.  The color of the vertex $b$ represents the truth value of $\ell_0 \orr \ell_1 \orr \cdots \orr \ell_{i-1}$ as well as the truth value of the literal $\ell_{i-1}$:  $b =_\nu x$ is interpreted as $\ell_0 \orr \ell_1 \orr \cdots \orr \ell_{i-1}$ is true but $\ell_{i-1}$ is false; $b =_\nu y$ is interpreted as $\ell_0 \orr \ell_1 \orr \cdots \orr \ell_{i-1}$ is false (and hence also as $\ell_{i-1}$ is false); and $b =_\nu z$ is interpreted as $\ell_{i-1}$ is true (and hence also as $\ell_0 \orr \ell_1 \orr \cdots \orr \ell_{i-1}$ is true).  Similarly, the color of the vertex $u$ represents the truth value of $\ell_0 \orr \ell_1 \orr \cdots \orr \ell_i$ as well as the truth value of the literal $\ell_i$.  However, the meanings of the colors are permuted:  $u =_\nu x$ is interpreted as $\ell_0 \orr \ell_1 \orr \cdots \orr \ell_i$ is false (and hence also as $\ell_i$ is false); $u =_\nu y$ is interpreted as $\ell_i$ is true (and hence also as $\ell_0 \orr \ell_1 \orr \cdots \orr \ell_i$ is true); and $u =_\nu z$ is interpreted as $\ell_0 \orr \ell_1 \orr \cdots \orr \ell_i$ is true but $\ell_i$ is false.  Lemma~\ref{lem:U3Widget} tells us that $U_{x,y,z}(\ell_i,b,u)$ properly implements this coding scheme.  Lemma~\ref{lem:U3Widget}~$(ii)$ says that if a $3$-coloring codes that $\ell_i$ is false and that $\ell_0 \orr \ell_1 \orr \cdots \orr \ell_{i-1}$ is false, then it must also code that $\ell_0 \orr \ell_1 \orr \cdots \orr \ell_i$ is false.  Lemma~\ref{lem:U3Widget}~$(iii)$ says that if $\nu$ is a $3$-coloring of the subgraph of $U_{x,y,z}(\ell_i,b,u)$ induced by $\{x,y,z,\ell_i,b\}$ coding that $\ell_i$ is false and that $\ell_0 \orr \ell_1 \orr \cdots \orr \ell_{i-1}$ is true, then $\nu$ can be extended to a $3$-coloring of $U_{x,y,z}(\ell_i,b,u)$ coding that $\ell_0 \orr \ell_1 \orr \cdots \orr \ell_i$ is true.  The reader may worry that here it is also possible to extend $\nu$ to incorrectly code that $\ell_0 \orr \ell_1 \orr \cdots \orr \ell_i$ is false, so we assure the reader that this is irrelevant.  What is important is that it is possible to extend $\nu$ to code the correct information.  Lemma~\ref{lem:U3Widget}~$(iv)$ says that if $\nu$ is a $3$-coloring of the subgraph of $U_{x,y,z}(\ell_i,b,u)$ induced by $\{x,y,z,\ell_i,b\}$ coding that $\ell_i$ is true, then $\nu$ can be extended to a $3$-coloring of $U_{x,y,z}(\ell_i,b,u)$ coding that $\ell_0 \orr \ell_1 \orr \cdots \orr \ell_i$ is true.  Lemma~\ref{lem:U3WidgetDecode} helps us deduce the colors of literal-coding vertices from the colors of auxiliary vertices and hence helps us compute a homogeneous set for a set of clauses from a homogeneous set for a graph.

The next widget combines $U_{x,y,z}(\ell,b,u)$ widgets into widgets coding clauses.

\begin{widget}
$D(\ell_0,\ell_1,\dots,\ell_{n-1})$ is the following widget.
\begin{center}
\begin{tikzpicture}[x=1cm, y=1cm, node/.style={circle, draw, minimum size=2em}, widget/.style={rectangle, draw, minimum size=2em}]
	
	\node[node] (0) at (-1, 1) {$0$};
	\node[node] (1) at (1, 1) {$1$};
	\node[node] (2) at (0, 0) {$2$};
	
	\draw (0) -- (1) -- (2) -- (0);

	\node[node] (l0) at (0,-1) {$\ell_0$};
	
	\node[widget] (U1) at (0, -2) {$U^1(\ell_1', \ell_0, u_1)$};
	\node[node] (l1p) at (-2, -2) {$\ell_1'$};
	\node[widget] (R1) at (-4, -2) {$R^1(\ell_1, \ell_1')$};
	\node[node] (l1) at (-6, -2) {$\ell_1$};
	\node[node] (u1) at (0, -3) {$u_1$};
				
	\draw (2) -- (l0) -- (U1) -- (u1);
	\draw (l1) -- (R1) -- (l1p) -- (U1);
	
	\node[widget] (U2) at (0, -4) {$U^2(\ell_2', u_1, u_2)$};
	\node[node] (l2p) at (-2, -4) {$\ell_2'$};
	\node[widget] (R2) at (-4, -4) {$R^2(\ell_2, \ell_2')$};
	\node[node] (l2) at (-6, -4) {$\ell_2$};
	\node[node] (u2) at (0, -5) {$u_2$};
				
	\draw (u1) -- (U2) -- (u2);
	\draw (l2) -- (R2) -- (l2p) -- (U2);
	
	\node[widget] (U3) at (0, -6) {$U^3(\ell_3, u_2, u_3)$};
	\node[node] (l3) at (-2, -6) {$\ell_3$};
	\node[node] (u3) at (0, -7) {$u_3$};
				
	\draw (u2) -- (U3) -- (u3);
	\draw (l3) -- (U3);
	
	\node[node] (un2) at (0, -9) {$u_{n-2}$};
	\node[widget] (Un1) at (0, -10) {$U^{n-1}(\ell_{n-1}', u_{n-2}, u_{n-1})$};
	\node[node] (ln1p) at (-3, -10) {$\ell_{n-1}'$};
	\node[widget] (Rn1) at (-6, -10) {$R^{n-1}(\ell_{n-1}, \ell_{n-1}')$};
	\node[node] (ln1) at (-9, -10) {$\ell_{n-1}$};
	\node[node] (un1) at (0, -11) {$u_{n-1}$};
	\node[node] (x) at (0, -12) {$x$};

	\draw[dashed] (u3) -- (un2);				
	\draw (un2) -- (Un1) -- (un1) -- (x);
	\draw (ln1) -- (Rn1) -- (ln1p) -- (Un1);
	
\end{tikzpicture}
\end{center}
The widget also contains the edge $(2,\ell_i)$ for each $i < n$, which we omitted from the diagram to keep it legible.  For $0 < i < n$, the sub-widget $U^i(\ell_i', u_{i-1}, u_i)$ is $U_{0,1,2}(\ell_i, u_{i-1}, u_i)$ if $i \equiv 0 \mod 3$, is $U_{2,0,1}(\ell_i', u_{i-1}, u_i)$ if $i \equiv 1 \mod 3$ (with $\ell_0$ in place of $u_0$ when $i=1$), and is $U_{1,2,0}(\ell_i', u_{i-1}, u_i)$ if $i \equiv 2 \mod 3$.  For $0 < i < n$, the sub-widget $R^i(\ell_i,\ell_i')$ is $R_{\substack{1 \mapsto 0 \\ 0 \mapsto 2}}(\ell_i,\ell_i')$ if $i \equiv 1 \mod 3$ and is $R_{\substack{0 \mapsto 1 \\ 1 \mapsto 2}}(\ell_i,\ell_i')$ if $i \equiv 2 \mod 3$.  If $i \equiv 0 \mod 3$, then there is just the vertex $\ell_i$ instead of the subgraph

\begin{center}
\begin{tikzpicture}[x=1cm, y=1cm, node/.style={circle, draw, minimum size=2em}, widget/.style={rectangle, draw, minimum size=2em}]
	
	\node[node] (l) at (0, 0) {$\ell_i$};
	\node[widget] (R) at (2, 0) {$R^i(\ell_i, \ell_i')$};
	\node[node] (lp) at (4, 0) {$\ell_i'$};
	
	\draw (l) -- (R) -- (lp);
		
\end{tikzpicture}.
\end{center}
The vertex $x$ is $0$ if $n-1 \equiv 0 \mod 3$, is $2$ if $n-1 \equiv 1 \mod 3$, and is $1$ if $x \equiv 2 \mod 3$.  Note that the vertex $x$ is thus drawn twice because it is identical to one of $0$, $1$, $2$.  For clarity, we also point out that in the case of $D(\ell_0)$, the widget is simply
\begin{center}
\begin{tikzpicture}[x=1cm, y=1cm, node/.style={circle, draw, minimum size=2em}, widget/.style={rectangle, draw, minimum size=2em}]
	
	\node[node] (0) at (0, 0) {$0$};
	\node[node] (1) at (1, 1) {$1$};
	\node[node] (2) at (2, 0) {$2$};
	
	\draw (0) -- (1) -- (2) -- (0);

	\node[node] (l0) at (1,-1) {$\ell_0$};
	
	\draw (0) -- (l0) -- (2);
	
\end{tikzpicture}.
\end{center}
\end{widget}

\begin{lemma}\label{lem:DWidget}{\ }
\begin{itemize}
\item[(i)] Every $3$-coloring $\nu$ of the subgraph of $D(\ell_0,\ell_1,\dots,\ell_{n-1})$ induced by $\{0,1,2,\ell_0,\ell_1,\dots,\ell_{n-1}\}$ in which $\ell_i =_\nu 1$ for some $i < n$ can be extended to a $3$-coloring of $D(\ell_0,\ell_1,\dots,\ell_{n-1})$.

\item[(ii)] There is no $3$-coloring $\nu$ of $D(\ell_0,\ell_1,\dots,\ell_{n-1})$ in which $\ell_0 =_\nu \ell_1 =_\nu \cdots =_\nu \ell_{n-1} =_\nu 0$.
\end{itemize}
\end{lemma}

\begin{proof}
For~$(i)$, let $\nu$ be a $3$-coloring of the subgraph induced by $\{0,1,2,\ell_0,\ell_1,\dots,\ell_{n-1}\}$ in which $\ell_i =_\nu 1$ for some $i < n$.  For each $i < n$, let $D_i(\ell_0,\ell_1,\dots,\ell_{n-1})$ denote the subgraph of $D(\ell_0,\ell_1,\dots,\ell_{n-1})$ induced by $0$, $1$, $2$, and the vertices appearing in $R^j(\ell_j,\ell_j')$ and $U^j(\ell_j',u_{j-1},u_j)$ for all $j \leq i$.  That is, if $i < n-1$, then $D_i(\ell_0,\ell_1,\dots,\ell_{n-1})$ is $D(\ell_0,\ell_i,\dots,\ell_i)$ without the edge between $u_i$ and $x$; and if $i=n-1$, then $D_i(\ell_0,\ell_1,\dots,\ell_{n-1})$ is $D(\ell_0,\ell_i,\dots,\ell_{n-1})$.  Item~$(i)$ is then the instance $i = n-1$ of the following claim.

\begin{claim}
For all $i < n$, $\nu$ can be extended to a $3$-coloring of $D_i(\ell_0,\ell_1,\dots,\ell_{n-1})$.  Moreover, if $\ell_j =_\nu 1$ for some $j \leq i$, then $\nu$ can be extended to a $3$-coloring of $D_i(\ell_0,\ell_1,\dots,\ell_{n-1})$ in which $\nu(u_i)$ codes this fact.  That is, if $i \equiv 0 \mod 3$, then $u_i \neq_\nu 0$; if $i \equiv 1 \mod 3$, then $u_i \neq_\nu 2$; and if $i \equiv 2 \mod 3$, then $u_i \neq_\nu 1$ (for $i=0$, interpret $u_0$ as $\ell_0$).
\end{claim}

\begin{proof}
By induction on $i < n$.  For $i=0$, $D_0(\ell_0,\ell_1,\dots,\ell_{n-1})$ is the subgraph of induced by $\{0,1,2,\ell_0\}$, which is $3$-colored by $\nu$ by assumption.  Clearly if $\ell_0 =_\nu 1$, then $\ell_0 \neq_\nu 0$.  Now suppose that $\nu$ has been extended to a $3$-coloring of $D_{i-1}(\ell_0,\ell_1,\dots,\ell_{n-1})$.  For the sake of argument, suppose that $i \equiv 1 \mod 3$ (the $i \equiv 0 \mod 3$ and $i \equiv 2 \mod 3$ cases are symmetric), and suppose that if $\ell_j =_\nu 1$ for some $j \leq i-1$, then $u_{i-1} \neq_\nu 0$.  First suppose that $\ell_i =_\nu 0$.  As $R^i(\ell_i,\ell_i') = R_{\substack{1 \mapsto 0 \\ 0 \mapsto 2}}(\ell_i,\ell_i')$, apply Lemma~\ref{lem:RWidget}~$(i)$ to extend $\nu$ to $R^i(\ell_i,\ell_i')$ so that $\ell_i' =_\nu 2$.  By Lemma~\ref{lem:U3Widget}~$(i)$, it is possible to extend $\nu$ to $U^i(\ell_i',u_{i-1},u_i)$.  Furthermore, if $\ell_j =_\nu 1$ for some $j \leq i-1$, then $u_{i-1} \neq_\nu 0$.  In this situation, by Lemma~\ref{lem:U3Widget}~$(iii)$, it is possible to extend $\nu$ to $U^i(\ell_i',u_{i-1},u_i) = U_{2,0,1}(\ell_i',u_{i-1},u_i)$ so that $u_i =_\nu 1$ (and hence $u_i \neq_\nu 2$).  Now suppose that $\ell_i =_\nu 1$.  As $R^i(\ell_i,\ell_i') = R_{\substack{1 \mapsto 0 \\ 0 \mapsto 2}}(\ell_i,\ell_i')$, apply Lemma~\ref{lem:RWidget}~$(i)$ to extend $\nu$ to $R^i(\ell_i,\ell_i')$ so that $\ell_i' =_\nu 0$.  By Lemma~\ref{lem:U3Widget}~$(iv)$, it is possible to extend $\nu$ to $U^i(\ell_i',u_{i-1},u_i) = U_{2,0,1}(\ell_i',u_{i-1},u_i)$ so that $u_i =_\nu 0$ (and hence $u_i \neq_\nu 2$).
\end{proof}

For~$(ii)$, suppose for the sake of contradiction that $\nu$ is a $3$-coloring of $D(\ell_0,\ell_1,\dots,\ell_{n-1})$ in which $\ell_0 =_\nu \ell_1 =_\nu \cdots =_\nu \ell_{n-1} =_\nu 0$.  We prove by induction on $i < n$ that $u_i =_\nu 0$ if $i \equiv 0 \mod 3$, $u_i = 2$ if $i \equiv 1 \mod 3$, and $u_i =_\nu 1$ if $i \equiv 2 \mod 3$ (again $u_0$ is interpreted as $\ell_0$).  Item~$(ii)$ follows from the case $i=n-1$ because this gives the contradiction $u_{n-1} =_\nu x$.  For $i = 0$, $\ell_0 =_\nu 0$ by assumption.  Now consider $0 < i < n$, assume for the sake of argument that $i \equiv 1 \mod 3$ (the $i \equiv 0 \mod 3$ and $i \equiv 2 \mod 3$ cases are symmetric), and assume that $u_{i-1} =_\nu 0$.  By Lemma~\ref{lem:RWidget}~$(i)$ for the widget $R^i(\ell_i,\ell_i') = R_{\substack{1 \mapsto 0 \\ 0 \mapsto 2}}(\ell_i,\ell_i')$, we have that $\ell_i' =_\nu 2$.  Thus $U^i(\ell_i',u_{i-1},u_i) = U_{2,0,1}(\ell_i',u_{i-1},u_i)$, $\ell_i' =_\nu 2$, and $u_{i-1} =_\nu 0$, so it must be that $u_i =_\nu 2$ by Lemma~\ref{lem:U3Widget}~$(ii)$.
\end{proof}

\begin{lemma}\label{lem:DWidgetDecode}
Let $\nu$ be a $3$-coloring of $D(\ell_0,\ell_1,\dots,\ell_{n-1})$.  If $0<i<n$ and $w$ is a vertex appearing in an $R^i(\ell_i,\ell_i')$ sub-widget or a $U^i(\ell_i',u_{i-1},u_i)$ sub-widget that is not $0$, $1$, or $2$, then the color of $w$ determines either the color of $\ell_i$ or the color of $\ell_{i-1}$.
\end{lemma}

\begin{proof}
Consider a $3$-coloring $\nu$ of $D(\ell_0,\ell_1,\dots,\ell_{n-1})$, an $i$ with $0<i<n$, and a vertex $w$ in an $R^i(\ell_i,\ell_i')$ sub-widget or a $U^i(\ell_i',u_{i-1},u_i)$ sub-widget that is not $0$, $1$, or $2$.  If $w$ appears in $R^i(\ell_i,\ell_i')$, then the color of $w$ determines the color of $\ell_i$ by Lemma~\ref{lem:RWidget}~$(iii)$.  If $w$ appears in $U^i(\ell_i',u_{i-1},u_i)$, then there are a few cases.  If $w$ is not $u_{i-1}$ or $d$, then the color of $w$ determines the color $\ell_i'$ by Lemma~\ref{lem:U3WidgetDecode}, which we have just seen determines the color of $\ell_i$ (or $\ell_i'$ is $\ell_i$ in the case $i \equiv 0 \mod 3$).  Consider $w = u_{i-1}$.  If $i=1$, then $u_{i-1}$ is really $\ell_0$, and of course the color of $\ell_0$ determines the color of $\ell_0$.  Otherwise, $i>1$, $u_{i-1}$ appears in the sub-widget $U^{i-1}(\ell_{i-1}',u_{i-2},u_{i-1})$, and hence the color of $u_{i-1}$ determines the color of $\ell_{i-1}$.

Lastly, consider $w=d$.  $U^i(\ell_i',u_{i-1},u_i)$ is $U_{x,y,z}(\ell_i',u_{i-1},u_i)$, where $x$, $y$, and $z$ are some permutation of $0$, $1$, and $2$.  If $d =_\nu x$ or $d =_\nu y$, then this determines the color of $\ell_i'$ by Lemma~\ref{lem:U3WidgetDecode}, which in turn determines the color of $\ell_i$.  Otherwise $d =_\nu z$, meaning that $u_{i-1} \neq_\nu z$ by Lemma~\ref{lem:U3WidgetDecode}.  If $i=1$, then $z=1$, $u_0$ is really $\ell_0$, and we conclude that $\ell_0 =_\nu 0$.  If $i > 1$, then $U^{i-1}(\ell_{i-1}', u_{i-2}, u_{i-1})$ is $U_{y,z,x}(\ell_{i-1}', u_{i-2}, u_{i-1})$ and, by examining the proof of Lemma~\ref{lem:U3WidgetDecode}, $u_{i-1} \neq_\nu z$ implies that $\ell_{i-1}' =_\nu y$, which in turn determines the color of $\ell_{i-1}$.
\end{proof}

To code a conjunction of two clauses $\ell_0 \orr \ell_1 \orr \cdots \orr \ell_{n-1}$ and $s_0 \orr s_1 \orr \cdots \orr s_{m-1}$, we overlap the widgets $D(\ell_0,\ell_1,\dots,\ell_{n-1})$ and $D(s_0,s_1,\dots,s_{m-1})$ by sharing the vertices pertaining to the longest common prefix of $\ell_0,\ell_1,\dots,\ell_{n-1}$ and $s_0,s_1,\dots,s_{m-1}$.  For example, consider the clauses $\ell_0 \orr \ell_1 \orr \ell_2 \orr \ell_3 \orr \ell_4$ and $\ell_0 \orr \ell_1 \orr s_2 \orr s_3$, where $\ell_2 \neq s_2$.  We overlap $D(\ell_0,\ell_1,\ell_2,\ell_3,\ell_4)$ and $D(\ell_0,\ell_1,s_2,s_3)$ as follows:

\begin{center}
\begin{tikzpicture}[x=1cm, y=1cm, node/.style={circle, draw, minimum size=2em}, widget/.style={rectangle, draw, minimum size=2em}]
	
	\node[node] (0) at (-1, 1) {$0$};
	\node[node] (1) at (1, 1) {$1$};
	\node[node] (2) at (0, 0) {$2$};
	
	\node[node] (l0) at (0,-1) {$\ell_0$};
			
	\draw (0) -- (1) -- (2) -- (0);
	
	\draw (2) -- (l0);
	
	\node[widget] (U1) at (-0, -2) {$U^1(\ell_1', \ell_0, u_1)$};
	\node[node] (l1p) at (-2, -2) {$\ell_1'$};
	\node[widget] (R1) at (-4, -2) {$R^1(\ell_1, \ell_1')$};
	\node[node] (l1) at (-6, -2) {$\ell_1$};
	\node[node] (u1) at (0, -3) {$u_1$};
	
	\draw (l0) -- (U1) -- (u1);
	\draw (l1) -- (R1) -- (l1p) -- (U1);
	
	\node[widget] (U2l) at (-2, -4) {$U^2(\ell_2', u_1, u_2)$};
	\node[node] (l2p) at (-4, -4) {$\ell_2'$};
	\node[widget] (R2l) at (-6, -4) {$R^2(\ell_2, \ell_2')$};
	\node[node] (l2) at (-6, -5) {$\ell_2$};
	\node[node] (u2) at (-2, -5) {$u_2$};
	
	\node[widget] (U2s) at (2, -4) {$U^2(s_2', u_1, v_2)$};
	\node[node] (s2p) at (4, -4) {$s_2'$};
	\node[widget] (R2s) at (6, -4) {$R^2(s_2, s_2')$};
	\node[node] (s2) at (6, -5) {$s_2$};
	\node[node] (v2) at (2, -5) {$v_2$};

	\draw (u1) -- (U2l) -- (u2);
	\draw (l2) -- (R2l) -- (l2p) -- (U2l);
	
	\draw (u1) -- (U2s) -- (v2);
	\draw (s2) -- (R2s) -- (s2p) -- (U2s);
	
	\node[widget] (U3l) at (-2, -6) {$U^3(\ell_3', u_2, u_3)$};
	\node[node] (l3) at (-4, -6) {$\ell_3$};
	\node[node] (u3) at (-2, -7) {$u_3$};
	
	\node[widget] (U3s) at (2, -6) {$U^3(s_3', v_2, v_3)$};
	\node[node] (s3) at (4, -6) {$s_3$};
	\node[node] (v3) at (2, -7) {$v_3$};

	\draw (u2) -- (U3l) -- (u3);
	\draw (l3) -- (U3l);
	
	\draw (v2) -- (U3s) -- (v3);
	\draw (s3) -- (U3s);
	
	\node[widget] (U4l) at (-2, -8) {$U^4(\ell_4', u_3, u_4)$};
	\node[node] (l4p) at (-4, -8) {$\ell_4'$};
	\node[widget] (R4l) at (-6, -8) {$R^4(\ell_4, \ell_4')$};
	\node[node] (l4) at (-6, -9) {$\ell_4$};
	\node[node] (u4) at (-2, -9) {$u_4$};

	\draw (u3) -- (U4l) -- (u4);
	\draw (l4) -- (R4l) -- (l4p) -- (U4l);
	
	\node[node] (2p) at (-2, -10) {$2$};
	\node[node] (0p) at (2, -8) {$0$};
	
	\draw (u4) -- (2p);
	\draw (v3) -- (0p);
	
\end{tikzpicture}
\end{center}

\begin{theorem}\label{thm:rcolor3-rsat}
$\rca \vdash \rcolor_3 \imp \rwkl$.
\end{theorem}

\begin{proof}
We prove $\rca \vdash \rcolor_3 \imp \rsat_{\textup{2-branching}}$.  The theorem follows by Proposition~\ref{prop:RWKL2BranchRSAT}.

Let $C$ be a $2$-branching and finitely satisfiable set of clauses over an infinite set of atoms $A = \{a_i : i \in \Nb\}$.  We assume that no clause in $C$ is a proper prefix of any other clause in $C$ by removing from $C$ every clause that has a proper prefix also in $C$.  We build a locally $3$-colorable graph $G$ such that every infinite homogeneous set for $G$ computes an infinite homogeneous set for $C$.  To start, $G$ contains the vertices $0$, $1$, and $2$, as well as the literal-coding vertices $a_i$ and $\neg a_i$ for each atom $a_i \in A$.  These vertices are connected according to the diagram below.

\begin{center}
\begin{tikzpicture}[x=1cm, y=1cm, node/.style={circle, draw, minimum size=2em}]

	\node[node] (0) at (-1,1) {$0$};
	\node[node] (1) at (-1, -1)  {$1$};
	\node[node] (2) at (0, 0) {$2$};
	
	\node[node] (a0) at (2,1) {$a_0$};
	\node[node] (na0) at (4,1) {$\neg a_0$};
	\node[node] (a1) at (2,-1) {$a_1$};
	\node[node] (na1) at (4,-1) {$\neg a_1$};
	
	\draw (0) -- (1) -- (2) -- (0);
	\draw (2) -- (a0) -- (na0) -- (2);
	\draw (2) -- (a1) -- (na1) -- (2);
	\draw[dashed] (3,-1.5) -- (3,-3);
	
\end{tikzpicture}
\end{center}

Now build $G$ in stages by considering the clauses in $C$ one-at-a-time.  For clause $\ell_0 \orr \ell_1 \orr \cdots \orr \ell_{n-1}$, find the previously appearing clause $s_0 \orr s_1 \orr \cdots \orr s_{m-1}$ having the longest common prefix with $\ell_0 \orr \ell_1 \orr \cdots \orr \ell_{n-1}$.  Then add the widget $D(\ell_0,\ell_1,\dots,\ell_{n-1})$ by overlapping it with $D(s_0,s_1,\dots,s_{m-1})$ as described above.  In $D(\ell_0,\ell_1,\dots,\ell_{n-1})$, for each $i < n$, the vertex $\ell_i$ is the vertex $a_i$ if the literal $\ell_i$ is the literal $a_i$, and the vertex $\ell_i$ is the vertex $\neg a_i$ if the literal $\ell_i$ is the literal $\neg a_i$.  The vertices appearing in the sub-widgets $R^i(\ell_i,\ell_i')$ and $U^i(\ell_i',u_{i-1},u_i)$ for $i$ beyond the index at which $\ell_0 \orr \ell_1 \orr \cdots \orr \ell_{n-1}$ differs from $s_0 \orr s_1 \orr \cdots \orr s_{m-1}$ are chosen fresh, except for $0$, $1$, $2$, and the literal-coding vertices $\ell_i$.  This completes the construction of $G$.

\begin{claim}
$G$ is locally $3$-colorable.
\end{claim}

\begin{proof}
Let $G_0$ be a finite subgraph of $G$.  Let $s$ be the latest stage at which a vertex in $G_0$ appears, and let $C_0 \subseteq C$ be the set of clauses considered up to stage $s$.  By extending $G_0$, we may assume that it is the graph constructed up to stage $s$.

By the finite satisfiability of $C$, let $t \colon \atoms(C_0) \imp \{\true, \false\}$ be a truth assignment satisfying $C_0$.  The truth assignment $t$ induces a $3$-coloring $\nu$ on the literal-coding vertices in $G_0$.  First define $\nu$ on the truth value-coding vertices by $\nu(0) = 0$, $\nu(1)=1$, and $\nu(2)=2$.  If $t(\ell)$ is defined for the literal $\ell$, then set $\nu(\ell) = t(\ell)$ (identifying $0$ with $\false$ and $1$ with $\true$).  If $\ell$ is a literal-coding vertex in $G_0$ on which $t$ is not defined, then set $\nu(\ell)=1$ if $\ell$ is a positive literal and set $\nu(\ell)=0$ if $\ell$ is a negative literal.  For each clause $\ell_0 \orr \ell_1 \orr  \cdots \orr \ell_{n-1}$ in $C_0$, extend $\nu$ to a $3$-coloring of $G_0$ by coloring each widget $D(\ell_0,\ell_1,\dots,\ell_{n-1})$ according to the algorithm implicit in the proof of Lemma~\ref{lem:DWidget}~$(i)$.  The hypothesis of Lemma~\ref{lem:DWidget}~$(i)$ is satisfied because $t$ satisfies $C_0$, so for each clause $\ell_0 \orr \ell_1 \orr \cdots \orr \ell_{n-1}$ in $C_0$, there is an $i < n$ such that $\ell_i =_\nu 1$.  Overlapping widgets $D(\ell_0,\ell_1,\dots,\ell_{n-1})$ and $D(s_0,s_1,\dots,s_{m-1})$ are colored consistently because the colors of the shared vertices depend only on the colors of the literal-coding vertices corresponding to the longest common prefix of the two clauses.
\end{proof}

Apply $\rcolor_3$ to $G$ to get an infinite homogeneous set $H$.  We may assume that $H$ contains exactly one of the truth value-coding vertices $0$, $1$, or $2$.  Call this vertex $c$.

Consider a vertex $w \in H$ that is not $c$.  The vertex $w$ appears in some widget $D(\ell_0, \ell_1, \dots, \ell_{n-1})$, and, by Lemma~\ref{lem:DWidgetDecode}, from $w$ we can compute an $i < n$ and a $c_i \in \{0,1\}$ such that $\ell_i =_\nu c_i$ whenever $\nu$ is a $3$-coloring of $D(\ell_0, \ell_1, \dots, \ell_{n-1})$ in which $w =_\nu c$.  Moreover, for each literal $\ell$, we can compute a bound on the number of vertices $w$ in the graph whose color determines the color of $\ell$.  Still by Lemma~\ref{lem:DWidgetDecode}, if $w$ appears in an $R^i(\ell_i,\ell_i')$ sub-widget or a $U^i(\ell_i',u_{i-1},u_i)$ sub-widget, then the color of $w$ determines either the color of $\ell_i$ or the color of $\ell_{i-1}$.  Thus the vertices whose colors determine the color of $\ell_i$ only appear in $R^i(\ell_i,\ell_i')$, $U^i(\ell_i',u_{i-1},u_i)$, $R^{i+1}(\ell_{i+1},\ell_{i+1}')$, and $U^{i+1}(\ell_{i+1}',u_i,u_{i+1})$ sub-widgets.  The fact that $C$ is a $2$-branching set of clauses and our protocol for overlapping the $D(\ell_0, \ell_1, \dots, \ell_{n-1})$ widgets together imply that, for every $j > 0$, there are at most $2^j$ sub-widgets of the form $R^j(\ell_j,\ell_j')$ and at most $2^j$ sub-widgets of the form $U^j(\ell_j',u_{j-1},u_j)$.  This induces the desired bound on the number of vertices whose colors determine the color of $\ell_i$.

Thus from $H$ we can compute an infinite set $H'$ of pairs $\la \ell, c_\ell \ra$, where each $\ell$ is a literal-coding vertex and each $c_\ell$ is either $0$ or $1$, such that every finite subgraph of $G$ is $3$-colorable by a coloring $\nu$ such that $(\forall \la \ell, c_\ell \ra \in H')(\ell =_\nu c_\ell)$.  Modify $H'$ to contain only pairs $\la a, c_a \ra$ for positive literal-coding vertices $a$ by replacing each pair of the form $\la \neg a, c_{\neg a} \ra$ with $\la a, 1-c_{\neg a} \ra$.  Now apply the infinite pigeonhole principle to $H'$ to get an infinite set $H''$ of positive literal-coding vertices $a$ and a new $c \in \{0,1\}$ such that the corresponding $c_a$ is always $c$.  We identify a positive literal-coding vertex $a$ with the corresponding atom and show that $H''$ is homogeneous for $C$.

Let $C_0 \subseteq C$ be finite.  Let $G_0$ be the finite subgraph of $G$ containing $\{0,1,2\}$, the literal-coding vertices whose atoms appear in the clauses in $C_0$, and the $D(\ell_0,\ell_1,\dots,\ell_{n-1})$ widgets for the clauses $\ell_0 \orr \ell_1 \orr \cdots \orr \ell_{n-1}$ in $C_0$.  By the homogeneity of $H''$ for $G$, there is a $3$-coloring $\nu$ of $G_0$ such that $a =_\nu c$ for every $a \in H''$.  From $\nu$, define a truth assignment $t$ on $\atoms(C_0)$ by $t(a) = \true$ if $a =_\nu 1$ and $t(a) = \false$ if $a =_\nu 0$.  This truth assignment satisfies every clause $\ell_0 \orr \ell_1 \orr \cdots \orr \ell_{n-1}$ in $C_0$.  The $3$-coloring $\nu$ must color the widget $D(\ell_0,\ell_1,\dots,\ell_{n-1})$, so by Lemma~\ref{lem:DWidget}~$(ii)$, it must be that $\ell_i =_\nu 1$ for some $i < n$.  Then $t(\ell_i) = \true$ for this same $i$, so $t$ satisfies $\ell_0 \orr \ell_1 \orr \cdots \orr \ell_{n-1}$.  Moreover, $t(a)$ is the truth value coded by $c$ for every $a \in H''$, so $H''$ is indeed an infinite homogeneous set for $C$.
\end{proof}

It follows that $\rwkl$ and $\rcolor_k$ are equivalent for every fixed $k \geq 3$.

\begin{corollary}\label{cor:RCOLORkisRWKL}
For every $k \geq 3$, $\rca \vdash \rwkl \biimp \rcolor_k$.
\end{corollary}

\begin{proof}
Fix $k \in \omega$ with $k \geq 3$.  $\rca \vdash \rwkl \imp \rcolor_k$ by Lemma~\ref{RWKLprovesLRCOLOR}.  
It is easy to see that $\rca \vdash \rcolor_k \imp \rcolor_3$.  Given a locally $3$-colorable graph $G$, augment $G$ by a clique $C$ containing $k-3$ fresh vertices, and put and edge between every vertex in $C$ and every vertex in $G$.  The resulting graph $G'$ is locally $k$-colorable, and every infinite set that is $k$-homogeneous for $G'$ is also $3$-homogeneous for $G$.  Finally, $\rca \vdash \rcolor_3 \imp \rwkl$ by Theorem~\ref{thm:rcolor3-rsat}.
\end{proof}

The question of the exact strength of $\rcolor_2$ remains open.  
We were unable to determine if $\rcolor_2$ implies $\rwkl$ or even if $\rcolor_2$ implies $\dnr$.

\section{The strength of $\rwkl$}\label{sect:strength-ramsey-strength-rwkl}

In this section, we prove various non-implications concerning $\rwkl$ and $\rcolor_2$.  The main result is that $\rca + \wwkl \nvdash \rcolor_2$ (Theorem~\ref{thm-WWKLDoesNotProveRCOLOR2}).  From this it follows that $\rca + \dnr \nvdash \rwkl$, which answers Flood's question of whether or not $\rca \vdash \dnr \imp \rwkl$ from~\cite{Flood2012Reverse}.  We also show that $\rca + \cac \nvdash \rcolor_2$ (Theorem~\ref{cac-not-imply-rcolor2}), where~$\cac$ stands for the chain antichain principle, defined below.  Note that it is immediate that $\rca + \cac \nvdash \rwkl$ because $\rca + \rwkl \vdash \dnr$ (by~\cite{Flood2012Reverse}) but $\rca + \cac \nvdash \dnr$ (by~\cite{Hirschfeldt2007Combinatorial}).  We do not know if $\rca \vdash \rcolor_2 \imp \dnr$, so we must give a direct proof that $\rca + \cac \nvdash \rcolor_2$.

In summary, the situation is thus.  $\wkl$ and $\rt^2_2$ each imply $\rwkl$ and therefore each imply $\rcolor_2$.  However, if $\wkl$ is weakened to $\wwkl$, then it no longer implies $\rcolor_2$.  Similarly, if $\rt^2_2$ is weakened to $\cac$, then it no longer implies $\rcolor_2$.

We begin our analysis of $\rcolor_2$ by constructing an infinite, computable, bipartite graph with no infinite, computable, homogeneous set.  It follows that $\rca \nvdash \rcolor_2$.  The graph we construct avoids potential infinite, r.e., homogeneous sets in a strong way that aids our proof that $\rca + \cac \nvdash \rcolor_2$.

\index{column-wise homogeneous set}
\begin{definition}
Let $G = (V,E)$ be an infinite graph.  A set $W \subseteq V^2$ is \emph{column-wise homogeneous for $G$} if $W^{[x]}$ is infinite for infinitely many $x$ (where $W^{[x]} = \set{y : \tuple{x,y} \in W}$ is the $x$\textsuperscript{th} column of $W$), and $\forall x \forall y(y \in W^{[x]} \imp \text{$\set{x,y}$ is homogeneous for $G$})$.
\end{definition}

\begin{lemma}\label{lem-RCOLOR2notCE}
There is an infinite, computable, bipartite graph $G = (\omega, E)$ such that no r.e.\ set is column-wise homogeneous for $G$.
\end{lemma}

\begin{proof}
The construction proceeds in stages, starting at stage $0$ with $E = \emptyset$.  We say that \emph{$W_e$ requires attention at stage $s$} if $e < s$ and there is a least pair $\tuple{x,y}$ such that
\begin{itemize}
\item $e < x < y < s$,
\item $y \in W_{e,s}^{[x]}$,
\item $x$ and $y$ are not connected to each other, and
\item neither $x$ nor $y$ is connected to a vertex $\leq e$.
\end{itemize}
At stage $s$, let $e$ be least such that $W_e$ requires attention at stage $s$ and has not previously received attention.  $W_e$ then receives attention by letting $\tuple{x,y}$ witness that $W_e$ requires attention at stage $s$, letting $u$ and $v$ be the least isolated vertices $> s$, and adding edges $(x,u)$, $(u,v)$, and $(v,y)$ to $E$.  This completes the construction.

We verify the construction.  We first show that $G$ is acyclic by showing that it is acyclic at every stage.  It follows that $G$ is bipartite because a graph is bipartite if and only if it has no odd cycles.  All vertices are isolated at the beginning of stage $0$, hence $G$ is acyclic at the beginning of stage~$0$.  By induction, suppose $G$ is acyclic at the beginning of stage~$s$.  If no $W_e$ requires attention at stage~$s$, then no edges are added at stage $s$, hence $G$ is acyclic at the beginning of stage~$s+1$.  If some least $W_e$ requires attention at stage $s$, then during stage~$s$ we add a length-$3$ path connecting the connected components of the $x$ and $y$ such that $\tuple{x,y}$ witnesses that $W_e$ requires attention at stage~$s$.  This action does not add a cycle because by the definition of requiring attention,~$x$ and~$y$ are not connected at the beginning of stage~$s$.  Hence $G$ is acyclic at the beginning of stage~$s+1$.

We now show that, for every $e$, if there are infinitely many $x$ such that $W_e^{[x]}$ is infinite, then there are an~$x$ and a~$y$ with $y \in W_e^{[x]}$ and $\set{x,y}$ not homogeneous for $G$.  If $W_e$ receives attention, then there is a length-$3$ path between an~$x$ and a~$y$ with $y \in W_e^{[x]}$, in which case $\set{x,y}$ is not homogeneous for $G$.  Thus it suffices to show that if $W_e^{[x]}$ is infinite for infinitely many~$x$, then $W_e$ requires attention at some stage.

Suppose that $W_e^{[x]}$ is infinite for infinitely many $x$, and suppose for the sake of contradiction that $W_e^{[x]}$ never requires attention.  Let $s_0$ be a stage by which every $W_i$ for $i < e$ that ever requires attention has received attention.  The graph contains only finitely many edges at each stage, so let $x_0$ be an upper bound for the vertices that are connected to the vertices $\leq e$ at stage $s_0$.  Notice that when some $W_i$ receives attention, the vertices connected at that stage are not connected to vertices $\leq i$.  Therefore once all the $W_i$ for $i < e$ that ever require attention have received attention, no vertex that is not connected to a vertex $\leq e$ is ever connected to a vertex $\leq e$.  In particular, no vertex $\geq x_0$ is ever connected to a vertex $\leq e$.  Now let $x > x_0$ be such that $W_e^{[x]}$ is infinite, and let $s_1 > s_0$ be a stage by which every $W_i$ for $i < x$ that ever requires attention has received attention.  Let $y_0$ be an upper bound for the vertices that are connected to $x$ and the vertices $\leq e$ at stage $s_1$, and again note that no vertex $\geq y_0$ is ever connected to $x$ or a vertex $\leq e$.  As $W_e^{[x]}$ is infinite, let $s > s_1$ be a stage at which there is a $y > y_0$ with $x < y < s$ and $y \in W_{e,s}^{[x]}$.  This $y$ is not connected to $x$, and neither $x$ nor $y$ is connected to a vertex $\leq e$, so $W_e$ requires attention at stage $s$, a contradiction.
\end{proof}

\begin{proposition}\label{prop-RCADoesNotProveRCOLOR2}
$\rca \nvdash \rcolor_2$.
\end{proposition}
\begin{proof}
Consider the $\omega$-model of $\rca$ whose second-order part consists of exactly the computable sets.  The graph $G$ from Lemma~\ref{lem-RCOLOR2notCE} is in the model because $G$ is computable.  However, the model contains no homogeneous set for $G$ because if $H$ were an infinite, computable, homogeneous set, then $\{\la x, y \ra : x, y \in H\}$ would be a computable, column-wise homogeneous set, thus contradicting Lemma~\ref{lem-RCOLOR2notCE}.
\end{proof}

The notion of \emph{restricted $\Pi^1_2$ conservativity} helps separate Ramsey-type weak K\"{o}nig's lemma and the Ramsey-type coloring principles from a variety of weak principles.

\index{restricted $\Pi^1_2$ sentence}
\index{restricted $\Pi^1_2$ conservative}
\begin{definition}[\cite{Hirschfeldt2007Combinatorial,Hirschfeldt2009atomic}]
\item A sentence is \emph{restricted $\Pi^1_2$} if it is of the form 
$\forall A(\Theta(A) \imp \exists B(\Phi(A,B)))$, where $\Theta$ is arithmetic and $\Phi$ is $\Sigma^0_3$.
A theory $T$ is \emph{restricted $\Pi^1_2$ conservative} over a theory $S$ 
if $S \vdash \varphi$ whenever $T \vdash \varphi$ and $\varphi$ is restricted $\Pi^1_2$.
\end{definition}

\begin{theorem}[\cite{Hirschfeldt2007Combinatorial,Hirschfeldt2009atomic}]\label{thm-rPi12conservative}
$\coh$ and~$\amt$ are restricted $\Pi^1_2$ conservative over $\rca$.
\end{theorem}

$\rcolor_2$ is a restricted $\Pi^1_2$ sentence, so we immediately have that neither $\coh$ nor $\amt$ implies $\rcolor_2$ over $\rca$.  Consequently, over $\rca$, $\coh$, $\amt$ and $\opt$ are all incomparable with $\rwkl$ and with $\rcolor_2$.

\begin{theorem}
$\rwkl$ and~ $\rcolor_2$ are incomparable with each of $\coh$, $\amt$, and $\opt$ over $\rca$.
\end{theorem}
\begin{proof}
By~\cite{Hirschfeldt2009atomic}, $\amt$ implies~$\opt$ over $\rca$.
Thus we need only show that neither $\rca + \coh$ nor $\rca + \amt$ prove $\rcolor_2$ 
and that $\rca + \rwkl$ proves neither $\coh$ nor $\opt$.  Observe that $\rcolor_2$ is restricted $\Pi^1_2$ sentence, so we have that neither $\rca + \coh$ nor $\rca + \amt$ proves $\rcolor_2$ by Proposition~\ref{prop-RCADoesNotProveRCOLOR2} and Theorem~\ref{thm-rPi12conservative}.  $\rca + \rwkl$ proves neither $\coh$ nor $\opt$ because $\rca + \wkl$ proves $\rca + \rwkl$ and $\rca + \wkl$ proves neither $\coh$~\cite{Hirschfeldt2007Combinatorial} nor $\opt$~\cite{Hirschfeldt2009atomic}.
\end{proof}

\index{chain antichain}
\index{cac@$\cac$|see {chain antichain}}
\index{partial order}
\index{scac@$\scac$|see {chain antichain}}
\begin{definition}[Chain-antichain]
A \emph{partial order} $P = (P, \leq_P)$ consists of a set $P$ together with a reflexive, antisymmetric, transitive, binary relation $\leq_P$ on $P$.  A \emph{chain} in $P$ is a set $S \subseteq P$ such that $(\forall x, y \in S)(x \leq_P y \orr y \leq_P x)$.  An \emph{antichain} in $P$ is a set $S \subseteq P$ such that $(\forall x, y \in S)(x \neq y \imp x |_P y)$ (where $x |_P y$ means that $x \nleq_P y \andd y \nleq_P x$).  A partial order $(P, \leq_P)$ is \emph{stable} if either $(\forall i \in P)(\exists s)[(\forall j > s)(j \in P \imp i \leq_P j) \orr (\forall j > s)(j \in P \imp i \mid_P j)]$ or $(\forall i \in P)(\exists s)[(\forall j > s)(j \in P \imp i \geq_P j) \orr (\forall j > s)(j \in P \imp i \mid_P j)]$.  $\cac$ is the statement ``every infinite partial order has an infinite chain or an infinite antichain.''  $\scac$ is the restriction of $\cac$ to stable partial orders.
\end{definition}

We now adapt the proof that $\rca + \cac \nvdash \dnr$ in~\cite{Hirschfeldt2007Combinatorial} to prove that $\rca + \cac \nvdash \rcolor_2$.  We build an $\omega$-model of $\rca + \scac + \coh$ that is not a model of $\rcolor_2$ by alternating between adding chains or antichains to stable partial orders and adding cohesive sets without ever adding an infinite set homogeneous for the graph from Lemma~\ref{lem-RCOLOR2notCE}.

\begin{lemma}\label{lem-SCACForcing}
Let $X$ be a set, let $G = (V, E)$ be a graph computable in $X$ such that no column-wise homogeneous set for $G$ is r.e.\ in $X$, and let $P = (P, \leq_P)$ be an infinite, stable partial order computable in $X$.  Then there is an infinite $C \subseteq P$ that is either a chain or an antichain such that no column-wise homogeneous set for $G$ is r.e.\ in $X \oplus C$.
\end{lemma}

\begin{proof}
For simplicity, assume that $X$ is computable.  The proof relativizes to non-computable $X$.  As $P$ is stable, assume for the sake of argument that $P$ satisfies $(\forall i \in P)(\exists s)[(\forall j > s)(j \in P \imp i \leq_P j) \orr (\forall j > s)(j \in P \imp i \mid_P j)]$.  The case with $\geq_P$ in place of $\leq_P$ is symmetric.  Also assume that there is no computable, infinite antichain $C \subseteq P$, for otherwise we are done.

Let $U = \set{i \in P : (\exists s)(\forall j > s)(j \in P \imp i \leq_P j)}$.  The fact that there is no computable, infinite antichain in $P$ implies that $U$ is infinite.  Let $F = (F, \sqsubseteq)$ be the partial order consisting of all $\sigma \in U^{<\omega}$ that are increasing in both $<$ and $\leq_P$, where $\tau \sqsubseteq \sigma$ if $\tau \succeq \sigma$.  Let $H$ be sufficiently generic for $F$, and notice that $H$ (or rather, the range of $H$, which is computable from $H$ as $H$ is increasing in $<$) is an infinite chain in $P$.  Suppose for the sake of contradiction that $W_e^H$ is column-wise homogeneous for $G$.  Fix a $\sigma \preceq H$ such that
\begin{align*}
\sigma \Vdash \forall x \forall y (y \in (W_e^H)^{[x]} \imp \text{$\set{x,y}$ is homogeneous for $G$}).
\end{align*}
Define a partial computable function $\tau \colon \omega^2 \imp P^{<\omega}$ by letting $\tau(x,i) \in P^{<\omega}$ be the string with the least code such that $\tau(x,i) \supseteq \sigma$, that $\tau(x,i)$ is increasing in both $<$ and $\leq_P$, and that $|(W_e^{\tau(x,i)})^{[x]}| > i$.  From here there are two cases:

Case~1: There are infinitely many pairs $\tuple{x,i}$ such that both $\tau(x,i)$ is defined and there is a $y \in (W_e^{\tau(x,i)})^{[x]}$ with $\set{x,y}$ not homogeneous for $G$.  The last element of such a $\tau(x,i)$ is in $P \setminus U$ because otherwise $\tau(x,i) \in F$ and $\tau(x,i) \preceq \sigma$, contradicting that $\sigma \Vdash \forall x \forall y (y \in (W_e^H)^{[x]} \imp \text{$\set{x,y}$ is homogeneous for $G$})$.  Thus the set $C$ consisting of the last elements of such strings $\tau(x,i)$ is an infinite r.e.\ subset of $P \setminus U$.  As elements $i$ of $P \setminus U$ have the property $(\exists s)(\forall j > s)(j \in P \imp i \mid_P j)$, we can thin $C$ to an infinite r.e.\ antichain in $P$ and hence to an infinite computable antichain in $P$, a contradiction.

Case~2: There are finitely many pairs $\tuple{x,i}$ such that $\tau(x,i)$ is defined and there is a $y \in (W_e^{\tau(x,i)})^{[x]}$ with $\set{x,y}$ not homogeneous for $G$.  In this case, let $x_0$ be such that if $x > x_0$ and $\tau(x,i)$ is defined, then $(\forall y \in (W_e^{\tau(x,i)})^{[x]})(\text{$\set{x,y}$ is homogeneous for $G$})$.  Notice that if $|(W_e^H)^{[x]}| > i$, then there is a $\tau$ with $\sigma \preceq \tau \preceq H$ such that $|(W_e^\tau)^{[x]}| > i$.  Hence if $(W_e^H)^{[x]}$ is infinite, then $\tau(x,i)$ is defined for all $i$.  Thus let 
\begin{align*}
W = \set{\tuple{x, \max(W_e^{\tau(x,i)})^{[x]}} : x > x_0 \andd i \in \omega \andd \text{$\tau(x,i)$ is defined}}.
\end{align*}
Then $W$ is an r.e.\ set that is column-wise homogeneous for $G$, a contradiction.

Thus there is no column-wise homogeneous set for $G$ that is r.e.\ in $H$.  Therefore (the range of) $H$ is our desired chain $C$.
\end{proof}

\begin{lemma}\label{lem-COHForcing}
Let $X$ be a set, let $G = (V, E)$ be a graph computable in $X$ such that no column-wise homogeneous set for $G$ is r.e.\ in $X$, and let $\vec R = (R_i)_{i \in \omega}$ be a sequence of sets uniformly computable in $X$.  Then there is an infinite set $C$ that is cohesive for $\overline R$ such that no column-wise homogeneous set for $G$ is r.e.\ in $X \oplus C$.
\end{lemma}

\begin{proof}
For simplicity, assume that $X$ is computable.  The proof relativizes to non-computable $X$.

We force with computable Mathias conditions $(D,L)$, where $D \subseteq \omega$ is finite, $L \subseteq \omega$ is infinite and computable, and every element of $D$ is less than every element of $L$.  The order is $(D_1, L_1) \sqsubseteq (D_0, L_0)$ if $D_0 \subseteq D_1$, $L_1 \subseteq L_0$, and $D_1 \setminus D_0 \subseteq L_0$.  Let $H$ be sufficiently generic.  Then $H$ is an infinite cohesive set for $\vec R$ (as in, for example, Section~4 of~\cite{Cholak2001strength}).

Suppose for the sake of contradiction that $W_e^H$ is column-wise homogeneous for $G$.  Let $(D,L)$ be a condition such that $D \subseteq H \subseteq L$ and
\begin{align*}
(D,L) \Vdash \forall x \forall y (y \in (W_e^H)^{[x]} \imp \text{$\set{x,y}$ is homogeneous for $G$}).
\end{align*}
Let
\begin{align*}
W = \set{\tuple{x,y} : \exists E(\text{$E$ is finite} \andd D \subseteq E \subseteq L \andd \tuple{x,y} \in W_e^E)}.
\end{align*}
$W$ is an r.e.\ set, and $\forall x \forall y (y \in W^{[x]} \imp \text{$\set{x,y}$ is homogeneous for $G$}$).  To see the second statement, suppose there is a $\tuple{x,y} \in W$ such that $\set{x,y}$ is not homogeneous for $G$, and let $E$ witness $\tuple{x,y} \in W$.  Then $(E, L \setminus E) \preceq (D,L)$, but $(E, L \setminus E) \Vdash (y \in (W_e^H)^{[x]} \andd \text{$\set{x,y}$ is not homogeneous for $G$})$, a contradiction.  Finally, $W \supseteq W_e^H$ because if $\tuple{x,y} \in W_e^H$, then there is a finite $E$ with $D \subseteq E \subseteq L$ such that $\tuple{x,y} \in W_e^E$, in which case $\tuple{x,y} \in W$.  Thus $W$ is an r.e.\ set that is column-wise homogeneous for $G$.  This contradicts the lemma's hypothesis.  Therefore no column-wise homogeneous set for $G$ is r.e.\ in $H$, so $H$ is the desired cohesive set.
\end{proof}

\begin{theorem}\label{cac-not-imply-rcolor2}
$\rca + \cac \nvdash \rcolor_2$
\end{theorem}
\begin{proof}
Iterate and dovetail applications of Lemma~\ref{lem-SCACForcing} and Lemma~\ref{lem-COHForcing} to build a collection of sets $\Scal$ such that $(\omega, \Scal) \vDash \rca + \coh + \scac$, the graph $G$ from Lemma~\ref{lem-RCOLOR2notCE} is in $\Scal$, and no set that is r.e.\ in any set in $\Scal$ is column-wise homogeneous for $G$.  Then $(\omega, \Scal) \vDash \cac$ by \cite{Hirschfeldt2007Combinatorial}, and $(\omega, \Scal) \nvDash \rcolor_2$ by the same argument as in Proposition~\ref{prop-RCADoesNotProveRCOLOR2}.
\end{proof}

We conclude by proving that $\rca + \dnr \nvdash \rwkl$, thereby answering \cite{Flood2012Reverse}~Question~9.  In fact, we prove the stronger result $\rca + \wwkl \nvdash \rcolor_2$.  This is accomplished by building a computable bipartite graph $G$ such that the measure of the set of oracles that compute homogeneous sets for $G$ is $0$.  It follows that there is a Martin-L\"of random $X$ that does not compute a homogenous set for $G$, and a model of $\rca + \wwkl + \neg \rcolor_2$ is then easily built from the columns of $X$.

Recall that, in the context of a bipartite graph $G = (V,E)$, a set $H \subseteq V$ is $2$-homogeneous for $G$ if no two vertices in $H$ are connected by an odd-length path in $G$.  Here we simply say that such an $H$ is \emph{$G$-homogeneous} (or just \emph{homogeneous}).  Likewise, if $H \subseteq V$ contains two vertices that are connected by an odd-length path in $G$, then $H$ is \emph{$G$-inhomogeneous} (or just \emph{inhomogeneous}).

\begin{theorem}\label{RCOLOR2hasNRA}
There is a computable bipartite graph $G = (\omega, E)$ such that the measure of the set of oracles that enumerate homogeneous sets for $G$ is $0$. 
\end{theorem}

\begin{proof}
By Lebesgue density considerations (see, for example, \cite{Nies2009Computability}~Theorem~1.9.4), if a positive measure of oracles enumerate infinite homogeneous sets for a graph $G$, then
\begin{align*}
(\forall \epsilon > 0)(\exists e)[\mu\{X : \text{$W_e^X$ is infinite and $G$-homogeneous}\} > 1-\epsilon].
\end{align*}
Thus it suffices to build $G$ to satisfy the following requirement $R_e$ for each $e \in \omega$:
\begin{align*}
R_e:\; \mu\{X : \text{$W_e^X$ is infinite and $G$-homogeneous}\} \leq 0.9.
\end{align*}

Let us first give a rough outline of the construction.  Observe our construction must necessarily produce a graph $G$ that does not contain an infinite connected component.  If $G$ has an infinite connected component, then that component contains a vertex $v$ such that infinitely many vertices are connected to $v$ by an even-length path.  These vertices that are at an even distance from $v$ can be effectively enumerated, and they form a homogeneous set.  Thus our graph $G$ must be a union of countably many finite connected components.  Each stage of the construction adds at most finitely many edges, and thus at each stage of the construction all but finitely many vertices are isolated. For each~$e$, our plan is the following.  We monitor the action of $W_e^X$ for all oracles~$X$ until we see a sufficient measure of $X$'s produce enough vertices (in a sense to made precise). Then, the idea is to satisfy the requirement by adding edges to these vertices in a way that defeats about half (in the measure-theoretic sense) of the oracles~$X$. This is done by a two-step process. Requirement $R_e$ acts by either type~I or type~II actions, the second type following the first type.  In a type~I action, $R_e$ locks some finite number of vertices, thereby preventing lower priority requirements from adding edges to these locked vertices.  In a type~II action, $R_e$ merges finitely many of $G$'s connected components into one connected component by adding some new edges while maintaining that $G$ is a bipartite graph. This merging is made in a way which ensures that for a sufficient measure of oracles~$X$, $W_e^X$ is inhomogeneous for the resulting graph.

We now present the construction in full detail.  At stage $s$, we say that
\begin{itemize}
\item \emph{$R_e$ requires type I attention} if $R_e$ has no vertices locked and there are strings of length $s$ witnessing that
\begin{align*}
\mu\{X : (\exists x \in W_{e,s}^X)(\text{$x$ is not connected to any $v$ locked by $R_k$ for any $k < e$})\} > 0.9;
\end{align*}

\item \emph{$R_e$ requires type II attention} if it currently has locked vertices due to a type~I action, has never acted according to type II, and there are strings of length $s$ witnessing that
\begin{align*}
\mu\{X : (\exists y \in W_{e,s}^X)(\text{$y$ is not connected to any $v$ locked by $R_k$ for any $k \leq e$})\} > 0.9;
\end{align*}

\item \emph{$R_e$ requires attention} if $R_e$ requires type I attention or requires type II attention.
\end{itemize}

At stage $0$, $E = \emptyset$, and no requirement has locked any vertices.

At stage $s+1$, let $e < s$ be least such that $R_e$ requires attention (if there is no such $e$, then go on to the next stage).  If $R_e$ requires type I attention, let $x_0, x_1, \dots, x_{n-1}$ be vertices that are not connected to any $v$ locked by $R_k$ for any $k < e$ and such that the strings of length $s$ witness that $\mu\{X : (\exists i < n)(x_i \in W_{e,s}^X)\} > 0.9$.  $R_e$ locks the vertices $x_0, x_1, \dots, x_{n-1}$.  All requirements $R_k$ for $k > e$ unlock all of their vertices.

If $R_e$ requires type II attention, let $y_0, y_1, \dots, y_{m-1}$ be vertices that are not connected to any $v$ locked by $R_k$ for any $k \leq e$ and such that the strings of length $s$ witness that $\mu\{X : (\exists j < m)(y_j \in W_{e,s}^X)\} > 0.9$.  Let $x_0, x_1, \dots, x_{n-1}$ be the vertices that are locked by $R_e$.  First we merge the connected components of the $x_i$'s into a single connected component and the connected components of the $y_j$'s into a single connected component.  To do this, let $a, b, c, d > s$ be fresh vertices, and add the edges $(a,b)$ and $(c,d)$.  The graph is currently bipartite, so for each $i < n$ add either the edge $(x_i, a)$ or $(x_i, b)$ so as to maintain a bipartite graph.  Similarly, merge the connected components of the $y_j$'s by adding either the edge $(y_j, c)$ or $(y_j, d)$ for each $j < m$.  The component of the $x_i$'s is disjoint from the component of the $y_j$'s because the $y_j$'s were chosen not to be connected to the $x_i$'s.  Thus both the graph $G_1$ obtained by adding the edge $(a,c)$ and the graph $G_2$ obtained by adding the edge $(a,d)$ are bipartite.  Each pair $\{x_i, y_j\}$ is homogeneous for exactly one of $G_1$ and $G_2$, and the strings of length $s$ witness that
\begin{align*}
\mu\{X : (\exists i < n)(\exists j < m)(x_i \in W_{e,s}^X \andd y_j \in W_{e,s}^X)\} > 0.8
\end{align*}
and therefore that
\begin{align*}
\mu\{X : \text{$W_{e,s}^X$ is either $G_1$-inhomogeneous or $G_2$ inhomogeneous}\} > 0.8.
\end{align*}
Thus the strings of length $s$ either witness that
\begin{align*}
\mu\{X : \text{$W_{e,s}^X$ is $G_1$-inhomogeneous}\} > 0.4,
\end{align*}
in which case we extend to $G_1$ by adding the edge $(a,c)$,
or that
\begin{align*}
\mu\{X : \text{$W_{e,s}^X$ is $G_2$-inhomogeneous}\} > 0.4,
\end{align*}
in which case we extend to $G_2$ by adding the edge $(a,d)$.  This completes the construction.

To verify the construction, we first notice that $G$ is bipartite because it is bipartite at every stage.  Furthermore, $G$ is computable because if an edge $(u,v)$ is added at stage $s$, either $u > s$ or $v > s$.  Thus to check whether an edge $(u,v)$ is in $G$, it suffices to check whether the edge has been added by stage $\max(u,v)$. 

We now verify that every requirement is satisfied.  Suppose that $R_e$ acts according to type II at some stage $s+1$.  Then $R_e$ is satisfied because we have ensured that 
\begin{align*}
 \mu \{X : W_e^X~  \text{is $G$-inhomogeneous}\} > 0.4 
\end{align*}
and thus that
\begin{align*}
  \mu \{X : W_e^X~  \text{is $G$-homogeneous}\} \leq 0.6.
\end{align*}

We prove by induction that, for every $e \in \omega$, $R_e$ is satisfied and there is a stage past which $R_e$ never requires attention.  Consider $R_e$.  If $\mu\{X : \text{$W_e^X$ is infinite}\} \leq 0.9$, then $R_e$ is satisfied and $R_e$ never requires attention.  So assume that $\mu\{X : \text{$W_e^X$ is infinite}\} > 0.9$.  By induction, let $s_0$ be a stage such that no $R_k$ for $k < e$ ever requires attention at a stage past $s_0$.  If $R_e$ has locked vertices at stage $s_0$, then these vertices remain locked at all later stages because no higher priority $R_k$ ever unlocks them.  If $R_e$ does not have locked vertices at stage $s_0$, then let $s_1 \geq s_0$ be least such that the strings of length $s_1$ witness that $R_e$ requires type~I attention.  Such an $s_1$ exists because $\mu\{X : \text{$W_e^X$ is infinite}\} > 0.9$ and because the finite set of vertices that are connected to vertices locked by the $R_k$ for $k < e$ have stabilized by stage $s_0$.  $R_e$ then requires and receives type I attention at stage $s_1$, and the vertices that $R_e$ locks at stage $s_1$ are never later unlocked.  So there is a stage $s_1 \geq s_0$ by which $R_e$ has locked a set of vertices that are never unlocked.  If $R_e$ has acted according to type II by stage $s_1$, then $R_e$ is satisfied and never requires attention past stage $s_1$.  If $R_e$ has not acted according to type II by stage $s_1$, let $s_2 \geq s_1$ be least such that the strings of length $s_2$ witness that $R_e$ requires type~II attention.  Such an $s_2$ exists because $\mu\{X : \text{$W_e^X$ is infinite}\} > 0.9$ and because, past stage $s_1$, no requirement except $R_e$ can act to connect a vertex to a vertex locked by an $R_k$ for a $k \leq e$.  $R_e$ then requires and receives type~II attention at stage $s_2$.  Hence $R_e$ is satisfied, and $R_e$ never requires attention at a later stage.  This completes the proof.
\end{proof}

\begin{theorem}\label{thm-WWKLDoesNotProveRCOLOR2}
$\rca + \wwkl \nvdash \rcolor_2$.
\end{theorem}

\begin{proof}
Let $G$ be the computable graph from Theorem~\ref{RCOLOR2hasNRA}.  There are measure $1$ many Martin-L\"of random sets, but only measure $0$ many sets compute homogeneous sets for $G$.  Thus let $X$ be a Martin-L\"of random set that does not compute a homogeneous set for $G$, and let $\mathfrak{M}$ be the structure whose first-order part is $\omega$ and whose second-order part is $\{Y : \exists k (Y \leq_T \bigoplus_{i<k}X^{[i]})\}$.  It is well-known that $\mathfrak{M} \vDash \rca + \wwkl$, which one may see by appealing to van Lambalgen's theorem (see \cite{Downey2010Algorithmic}~Section 6.9) and the equivalence between $\wwkl$ and $\rans{1}$.  Moreover, $\mathfrak{M} \nvDash \rcolor_2$ because $\mathfrak{M}$ contains the bipartite graph $G$, but it does not contain any homogeneous set for $G$.
\end{proof}

It now follows that $\rca + \dnr \nvdash \rwkl$.  This has been proved independently by Flood and Towsner~\cite{Flood2014Separating} using
the techniques introduced by Lerman, Solomon, and Towsner~\cite{Lerman2013Separating}.  Recently, the author~\cite{Patey2015Ramsey} enhanced the separation of $\dnr$ and $\rwkl$ by proving that for every computable order $h$, there is an $\omega$-model of the statement ``for every $X$ there is a function that is $\dnrf$ relative to $X$ and bounded by $h$'' that is not a model of $\rcolor_2$.  This answers a question in~\cite{Flood2014Separating}.

\begin{corollary}\label{cor:DNRNoImpRWKL}
$\rca + \dnr \nvdash \rwkl$.
\end{corollary}

\begin{proof}
This follows from Theorem~\ref{thm-WWKLDoesNotProveRCOLOR2} because $\rca \vdash \wwkl \imp \dnr$ and $\rca \vdash \rwkl \imp \rcolor_2$.
\end{proof}

\part{Further topics}

\chapter{Degrees bounding principles}\label{chap-degrees-bounding-principles}

Many theorems studied in reverse mathematics are~$\Pi^1_2$ statements
and can be seen as collections of problems parameterized by their instances.
The study of theorems which are collections of problems rather than the study of the instances individually
is justified by the existence of a single argument, namely, the proof of the theorem, to construct
a solution to each instance. An effective analysis of the argument enables one
to give a uniform bound on the solutions to the instances, even tough
some instances may be harder to solve than others.

Some theorems such as cohesiveness or weak K\"onig's lemma happen to have maximally difficult instances,
that is, instances such that any solution to them can produce a solution to any other instance. Such an
instance is then called~\emph{universal}. 

\index{universal instance}
\index{P-bounding@$\Psf$-bounding}
\begin{definition}[Universal instance]
Given a statement $\Psf$, a degree $\textbf{d}$ \emph{bounds} $\Psf$ (written $\textbf{d} \gg_{\Psf} \emptyset$)
if every computable instance $X$ of $\Psf$ has a $\textbf{d}$-computable solution.
A statement $\Psf$ admits a \emph{universal instance} if there is a computable
$\Psf$-instance $X$ such that every solution to $X$ bounds $\Psf$.
\end{definition}

The notation $\textbf{d} \gg \emptyset$ historically means that the degree $\textbf{d}$
is PA and therefore is equivalent to $\textbf{d} \gg_{\wkl} \emptyset$.
It is well-known that $\wkl$ admits a universal instance 
-- e.g. take the $\Pi^0_1$ class of completions of Peano arithmetics --.
A few principles have been proven to admit universal instances, namely, weak K\"onig's lemma~\cite{Odifreddi1992Classical}, 
König's lemma~\cite{Hirschfeldt2015Slicing}, cohesiveness~\cite{Jockusch1993cohesive},
the Ramsey-type weak weak König's lemma~\cite{Bienvenu2015logical}, 
the finite intersection property ($\fip$) \cite{Downey2012finite}, 
the omitting partial type theorem ($\opt$) \cite{Hirschfeldt2009atomic},
or even the rainbow Ramsey theorem for pairs~\cite{MillerAssorted}, but most of the principles
studied in reverse mathematics do not admit one.

The standard way to prove that some statement~$\Psf$ does not admit a universal
instance consists of showing that every computable~$\Psf$-instance has a solution
satisfying some weakness property~$\Pcal$, but that no degree in~$\Pcal$ bounds~$\Psf$.
More generally, the following lemma holds.

\begin{lemma}
Fix some weakness property~$\Pcal$ and let~$\Psf$ and~$\Qsf$ be two~$\Pi^1_2$ statements such that the following holds.
\begin{itemize}
	\item[(i)] Every computable~$\Psf$-instance has a solution in~$\Pcal$.
	\item[(ii)] For every set~$X \in \Pcal$, there is a computable~$\Qsf$-instance with no $X$-computable solution
\end{itemize}
Then no $\Pi^1_2$ statement~$\Rsf$ such that~$\Qsf \leq_c \Rsf \leq_c \Psf$
admits a universal instance.
\end{lemma}
\begin{proof}
Let~$\Rsf$ be a $\Pi^1_2$ statement such that~$\Qsf \leq_c \Rsf \leq_c \Psf$.
Let~$I$ be a computable~$\Rsf$-instance. We will prove that it is not universal.
Since~$\Rsf \leq_c \Psf$, $I$ has a solution~$Y \in \Pcal$,
and since~$\Qsf \leq_c \Rsf$, there is a computable~$\Rsf$-instance~$J$ with no~$Y$-computable solution.
Therefore, the solution~$Y$ to~$I$ does not bound~$\Psf$ and~$I$ is not a universal instance.
\end{proof}

Mileti~\cite{Mileti2004Partition} studied the degrees bounding Ramsey's theorem for pairs
and proved in particular that stable Ramsey theorem for pairs admits no bound of low${}_2$ degree.
Therefore, no statement computably in between Ramsey's theorem for pairs and
its stable version admits a universal instance.
Independently, Hirschfeldt \& Shore proved in~\cite{Hirschfeldt2007Combinatorial} 
that the stable ascending descending sequence principle admits no bound of low degree.
Hence neither $\sads$ nor the stable chain antichain principle admit a universal instance.
In this chapter, we generalize these results and study the degree boundings Ramsey-type theorems.

\section{Degrees bounding the atomic model theorem}

The atomic model theorem is a statement of model theory admitting
a simple, purely computability-theoretic characterization over $\omega$-models.
This statement happens to have a weak computational content and is therefore a consequence
of many other principles in reverse mathematics. For those reasons,
the atomic model theorem is a good candidate for factorizing proofs 
of properties which are closed upward by the consequence relation.
Recall that the principle $\amt$ has been proven in~\cite{Hirschfeldt2009atomic,Conidis2008Classifying}
to be computably equivalent to the following principle:

\index{escape property}
\begin{definition}[Escape property]
For every $\Delta^0_2$ function $f$, there exists a function~$g$ such that $f(x) \leq g(x)$ for infinitely many $x$.
\end{definition}

This equivalence does not hold over $\rca$ as, unlike $\amt$, 
the escape property implies $\ist$ over $\bst$~\cite{Hirschfeldt2009atomic}.
Hirschfeldt and Shore proved in~\cite{Hirschfeldt2007Combinatorial} 
that the stable ascending descending sequence principle admits no bound of low degree.
Using the characterization of $\amt$ in terms of the escape property, 
we can easily deduce the following theorem.

\begin{theorem}\label{thm:amt-bounding}
There is no low${}_2$ $\Delta^0_2$ degree bounding $\amt$.
\end{theorem}

Note that Theorem~\ref{thm:amt-bounding} is a strengthening of the result of 
Hirschfeldt and Shore~\cite{Hirschfeldt2007Combinatorial}
since Hirschfeldt et al.~\cite{Hirschfeldt2009atomic} proved that~$\amt$ is a consequence of $\sads$ over computable reducibility.
Theorem~\ref{thm:amt-bounding} can be easily
proven using the following characterization of $\Delta^0_2$ low${}_2$ sets in terms of domination:

\begin{lemma}[Martin, \cite{Martin1966Classes}]\label{lem:shypzp-lowt}
A set $A \leq_T \emptyset'$ is low${}_2$ iff there exists an $f \leq_T \emptyset'$ dominating every $A$-computable function.
\end{lemma}
\begin{proof}
A set $A$ is low${}_2$ iff $\emptyset'$ is high relative to $A$.
We conclude the lemma from the observation that a set $X$ is high relative to a set $A \leq_T \emptyset'$
iff it computes a function dominating every $A$-computable function.
\end{proof}

As explained Conidis in~\cite{Conidis2008Classifying},
Theorem~\ref{thm:amt-bounding} cannot be extended to every low${}_2$ sets.
Indeed, Soare~\cite{Conidis2008Classifying} constructed a low${}_2$ set bounding
the escape property using a forcing argument. So there exists
a low${}_2$ degree bounding $\amt$.

\begin{corollary}\label{cor:amt-universal}
No principle $\Psf$ having an $\omega$-model with only low sets
and such that $\amt \leq_c \Psf$ admits a universal instance.
\end{corollary}
\begin{proof}
Suppose for the sake of contradiction that $\Psf$ has a universal instance $U$
and an $\omega$-model $\Mcal$ with only low sets. As $U$ is computable, $U \in \Mcal$.
Let $X \in \Mcal$ be a (low) solution to $U$. In particular, $X$ is low${}_2$ and $\Delta^0_2$,
so by Lemma~\ref{lem:shypzp-lowt} and the computable equivalence of $\amt$ and the escape property, 
there exists a computable instance~$Y$ of~$\amt$ 
such that $X$ does not compute a solution to~$Y$. 
As $\amt \leq_c \Psf$, there exists a $Y$-computable (hence computable) instance~$Z$
of~$\Psf$ such that every solution to~$Z$ computes a solution to $Y$. Thus $X$ does not compute a solution to~$Z$,
contradicting the universality of~$U$.
\end{proof}

Hirschfeldt et al.\ proved in~\cite{Hirschfeldt2007Combinatorial}
the existence of an $\omega$-model of $\sads$ and $\scac$ with only low sets. 
Therefore we obtain another proof that neither $\sads$ nor $\scac$ admits a universal instance.
The result was first proven in~\cite{Hirschfeldt2007Combinatorial} using an ad-hoc notion of reducibility.

\begin{corollary}
None of $\amt$, $\sads$ and $\scac$ admit a universal instance.
\end{corollary}

The previous argument cannot directly be applied to $\srt^2_2$, $\semo$ or $\sts^2$
as none of those principles admit an $\omega$-model with only 
low sets~\cite{Downey200102,Kreuzer2012Primitive,Patey2015Somewhere}.
However Corollary~\ref{cor:amt-universal} can be extended to principles such that every computable
instance has a $\Delta^0_2$ low${}_2$ solution. It is currently unknown whether
every $\Delta^0_2$ set admits a $\Delta^0_2$ low${}_2$ infinite subset in either it or its complement.
A positive answer would lead to a proof that $\srt^2_2$, $\semo$ and $\sts^2$ have no universal instance,
and more importantly, would provide an $\omega$-model of $\srt^2_2$ that is not a model of $\dnrs{2}$
as explained in~\cite{Hirschfeldt2008strength}. We shall see in a later section that none of $\srt^2_2$,
$\semo$ and $\sts^2$ admits a universal instance.

\section{No low${}_2$ degree bounds $\sts^2$ or~$\sads$}

Mileti originally proved in~\cite{Mileti2004Partition} that no principle $\Psf$ having an $\omega$-model
with only low${}_2$ sets and satisfying $\srt^2_2 \leq_c \Psf$ admits a universal instance, and deduced
that none of $\srt^2_2$ and $\rt^2_2$ admit one.
In this section, we reapply his argument to much weaker statements and derive non-universality
results to a large range of principles in reverse mathematics.
Thin set theorem and ascending descending sequence are example of statements
weak enough to be a consequence of many others, and surprisingly strong enough
to diagonalize against low${}_2$ sets.
As does Mileti in~\cite{Mileti2004Partition}, we deduce the following theorem
by the application of his proof technique.

\begin{theorem}\label{thm:sts-semo-sads-bounding}
There exists no low${}_2$ degree bounding any of $\sts^2$ or $\sads$.
\end{theorem}

\begin{corollary}\label{cor:sts-semo-sads-universal}
No principle $\Psf$ having an $\omega$-model with only low${}_2$ sets
and such that any of $\sts^2$, $\sads$ is computably reducible to $\Psf$
admits a universal instance.
\end{corollary}

The proof of Theorem~\ref{thm:sts-semo-sads-bounding} is split into three lemmas.
Lemma~\ref{lem:universal-instance} provides a general way of obtaining
bounding and universality results, assuming the ability of a principle to diagonalize
against a particular set. 
Lemma~\ref{lem:sts2-universal} and Lemma~\ref{lem:sads-universal}
state the desired diagonalization for respectively $\sts^2$ and $\sads$.

\begin{corollary}\label{thm:rt22-sts2-universal}
None of the following principles admits a universal instance:
$\rt^2_2$, $\rwkls{2}$, $\fs^2$, $\ts^2$, $\cac$, $\ads$ and their
stable versions.
\end{corollary}
\begin{proof}
Each of the above mentioned principles is a consequence of $\rt^2_2$ over $\rca$
and computably implies either $\sads$ or $\sts^2$.
See~\cite{Flood2012Reverse} for $\rwkls{2}$, \cite{Cholak2001Free} for~$\fs^2$ and $\ts^2$,
and \cite{Hirschfeldt2007Combinatorial} for $\cac$ and $\ads$.
By Theorem~3.1 of \cite{Cholak2001strength}, there exists an $\omega$-model of $\rt^2_2$
having only low${}_2$ sets. The result now follows from Corollary~\ref{cor:sts-semo-sads-universal}.
\end{proof}

In order to prove Theorem~\ref{thm:sts-semo-sads-bounding},
we need the following theorem proven by Mileti. It simply consists of applying a relativized 
version of the low basis theorem to a $\Pi^0_1$ class of completions of the enumeration
of all partial computable sets.

\begin{theorem}[Mileti, Corollary 5.4.5 of \cite{Mileti2004Partition}]\label{thm:mileti-f-compute}
For every set $X$, there exists $f : \omega^2 \to \set{0,1}$ low over $X$ such that for every
$X$-computable set $Z$, there exists an $e \in \omega$ with $Z = \set{a \in \omega : f(e,a) = 1}$.
\end{theorem}

Before going into the core lemmas, we show how to obtain Theorem~\ref{thm:sts-semo-sads-bounding}
from Lemma~\ref{lem:sts2-universal} and Lemma~\ref{lem:sads-universal}.

\begin{lemma}\label{lem:universal-instance}
Fix an $n \in \omega$ and two principles $\Psf$ and $\Qsf$ such that $\Psf \leq_c \Qsf$.
Suppose that for any $f : \omega^2 \to \set{0,1}$ satisfying $f'' \leq_T \emptyset^{(n+2)}$,
there exists a computable instance $I$ of $\Psf$ such that for each $e \in \omega$, 
if $\{a \in \omega : f(e,a) = 1\}$ is infinite then it is not a solution to~$I$. Then the following holds:
\begin{itemize}
  \item[(i)] For any degree $\mathbf{d}$ low${}_2$ over $\emptyset^{(n)}$ there is a computable instance
  $U$ of $\Psf$ such that $\mathbf{d}$ does not bound a solution to $U$.
  \item[(ii)] There is no degree low${}_2$ over $\emptyset^{(n)}$ bounding $\Psf$.
  \item[(iii)] If every computable instance $I$ of $\Qsf$ has a solution low${}_2$ over $\emptyset^{(n)}$,
  then $\Qsf$ has no universal instance.
\end{itemize}
\end{lemma}
\begin{proof}\ 
\begin{itemize}
  \item[(i)]
Consider any set $X$ of degree low${}_2$ over $\emptyset^{(n)}$.
By Theorem~\ref{thm:mileti-f-compute}, there exists a function $f : \omega^2 \to \set{0,1}$
low over $X$, hence low${}_2$ over $\emptyset^{(n)}$, such that any $X$-computable set $Z$
is of the form $\{a \in \omega : f(e,a) = 1\}$ for some $e \in \omega$. 
Take a computable instance $I$ of $\Psf$ having no solution of the form $\{a \in \omega : f(e, a) = 1\}$
for any $e \in \omega$. Then $X$ does not compute a solution to $I$.
  \item[(ii)] Immediate from (i).
  \item[(iii)] Take any computable instance $U$ of $\Qsf$. By assumption, $U$ has a solution $X$ low${}_2$ over $\emptyset^{(n)}$.
  By (i), there exists an instance $I$ of $\Psf$ such that $X$ does not compute a solution to~$I$.
	As $\Psf \leq_c \Qsf$,
	there exists an $I$-computable (hence computable) instance $J$ of $\Qsf$ such that any solution to~$J$
	computes a solution to~$I$. Then $X$ does not compute a solution to~$J$,
	hence $U$ is not a universal instance.
\end{itemize}
\end{proof}

We will prove the following lemmas which, together with Lemma~\ref{lem:universal-instance},
are sufficient to deduce Theorem~\ref{thm:sts-semo-sads-bounding}.

\begin{lemma}\label{lem:sts2-universal}
Fix a set $X$. Suppose $f : \omega^2 \to \set{0,1}$ satisfies $f'' \leq_T X^{''}$.
There exists an $X$-computable stable coloring $g : [\omega]^2 \to \omega$ such that
for all $e \in \omega$, if $\{a \in \omega : f(e, a) = 1\}$ is infinite then it
is not thin for $g$.
\end{lemma}

\begin{lemma}\label{lem:sads-universal}
Fix a set $X$. Suppose $f : \omega^2 \to \set{0,1}$ satisfies $f'' \leq_T X^{''}$.
There exists a stable $X$-computable linear order $L$ such that for all $e \in \omega$,
if $\{a \in \omega : f(e, a) = 1\}$ is infinite then it is neither an ascending
nor a descending sequence in~$L$.
\end{lemma}

Before proving the two remaining lemmas, we relativize the results to colorings
over arbitrary tuples.

\begin{theorem}\label{thm:stsn-bounding}
For any $n$, there exists no degree low${}_2$ over $\emptyset^{(n)}$
bounding $\sts^{n+2}$.
\end{theorem}
\begin{proof}
Apply Lemma~\ref{lem:sts2-universal} relativized to $X = \emptyset^{(n)}$
together with Lemma~\ref{lem:universal-instance}.
Simply notice that if $f : [\omega]^n \to \omega$ is a $\emptyset'$-computable
coloring, the computable coloring $g : [\omega]^{n+1} \to \omega$ obtained by an application
of Shoenfield's limit lemma is such that every infinite set thin for $g$ is thin for~$f$.
\end{proof}

\begin{corollary}\label{cor:stsn-universal}
For any $n$, no principle $\Psf$ having an $\omega$-model 
with only low${}_2$ over $\emptyset^{(n)}$ sets
and such that $\sts^{n+2} \leq_c \Psf$ admits a universal instance.
\end{corollary}
\begin{proof}
Same reasoning as Corollary~\ref{cor:sts-semo-sads-universal}
using the remark in the proof of Theorem~\ref{thm:stsn-bounding}.
\end{proof}

\begin{theorem}
For any $n$, none of
$\rt^{n+2}_2$, $\rwkls{(n+2)}$, $\fs^{n+2}$, $\ts^{n+2}$ and their stable versions
admits a universal instance.
\end{theorem}
\begin{proof}
Fix an $n \in \omega$. Each of the above cited principles $\Psf$
satisfies $\sts^{n+2} \leq_c \Psf$ and is a consequence of $\rt^{n+2}_2$ over $\omega$-models.
Cholak et al.~\cite{Cholak2001strength} proved the existence of an $\omega$-model
of $\rt^{n+2}_2$ having only low${}_2$ over $\emptyset^{(n)}$ sets. 
Apply Corollary~\ref{cor:stsn-universal}.
\end{proof}

We now turn to the proofs of Lemma~\ref{lem:sts2-universal},
and Lemma~\ref{lem:sads-universal}.

\begin{proof}[Proof of Lemma~\ref{lem:sts2-universal}]
We prove it in the case when~$X = \emptyset$. The general case follows by a straightforward relativization.
For each $e \in \omega$, let $Z_e = \set{a \in \omega : f(e,a) = 1}$.
The proof is very similar to \cite[Theorem~5.4.2.]{Mileti2004Partition}.
We build a $\emptyset'$-computable function $c : \omega \to \omega$
such that for all $e \in \omega$, if $Z_e$ is infinite then it
is not thin for~$c$. Given such a function $c$, we can then apply Shoenfield's limit lemma
to obtain a stable computable function $h : [\omega]^2 \to \omega$ such that
for each $x \in \omega$, $\lim_s h(x,s) = c(x)$.
Every set thin for $h$ is thin for $c$, and therefore for all $e \in \omega$, 
if $Z_e$ is infinite then it
is not thin for~$h$.

Suppose by Kleene's fixpoint theorem that
we are given a Turing index $d$ of the function $c$ as computed relative to~$\emptyset'$.
The construction is done by a finite injury priority argument
satisfying the following requirements for each $e, i \in \omega$:

\bigskip
$\Rcal_{e, i}$ : $Z_e$ is finite or $(\exists a)[f(e,a) = 1 \mbox{ and } \Phi^{\emptyset'}_d(a) = i]$
\bigskip

The requirements are ordered in a standard way, that is, ordering pairs in lexicographic order.
Notice that each of these requirement is $\Sigma^f_2$, and furthermore we can effectively find an index
for each as such. Therefore, for each $e$ and $i \in \omega$, we can effectively find an integer
$m_{e,i}$ such that $R_{e,i}$ is satisfied if and only if $m_{e,i} \in f''$. 
By Shoenfield's limit Lemma relativized to $\emptyset'$ and low${}_2$-ness of $f$, 
there exists a $\emptyset'$-computable function $g : \omega^2 \to 2$ such that for all $m$,
we have $m \in f'' \biimp \lim_s g(m, s) = 1$ and $m \not \in f'' \biimp \lim_s g(m, s) = 0$.
Notice that for all $e$ and $i \in \omega$, $R_{e,i}$ is satisfied if and only if $\lim_s g(m_{e,i}, s) = 1$.

At stage $s$, assume we have defined $c(u)$ for every $u < s$.
If there exists a least strategy $\Rcal_{e,i}$ (in priority order) with $\tuple{e,i} < s$
such that $g(m_{e,i}, s) = 0$, set $c(s) = i$. Otherwise set $c(s) = 0$.
This ends the construction. We now turn to the verification.

\begin{claim}
Every requirement $\Rcal_{e,i}$ is satisfied.
\end{claim}
\begin{proof*}
By induction over ordered pairs $\tuple{e,i}$
in lexicographic order. Suppose that $R_{e',i'}$ is satisfied for all $\tuple{e',i'} < \tuple{e,i}$, but
$\Rcal_{e,i}$ is not satisfied. Then there exists a threshold $t \geq \tuple{e,i}$ such that $g(m_{e',i'}, s) = 1$
for all $\tuple{e',i'} < \tuple{e,i}$ and $g(m_{e,i}, s) = 0$ whenever $s \geq t$.
By construction, $c(s) = i$ for every $s \geq t$.
As $Z_e$ is infinite, there exists an element $s \in Z_e$ such that $c(s) = i$,
so $Z_e$ is not thin for $c$ with witness $i$ and therefore $\Rcal_{e,i}$ is satisfied. Contradiction.
\end{proof*}
\end{proof}

\begin{proof}[Proof of Lemma~\ref{lem:sads-universal}]
Again, we prove the lemma in the case when~$X = \emptyset$.
For each $e \in \omega$, let $Z_e = \set{a \in \omega : f(e,a) = 1}$.
The proof is very similar to \cite[Theorem~5.4.2.]{Mileti2004Partition}.
We build a $\Delta^0_2$ set $U$ together with a stable computable linear order $L$
such that $U$ is the $\omega$ part of $L$, that is, 
$U$ is the collection of elements $L$-below cofinitely many other elements.
We furthermore ensure that for each $e \in \omega$, if $Z_e$ is infinite, then 
it intersects both $U$ and $\overline{U}$. Therefore, if $Z_e$ is infinite,
it is neither an ascending, nor a descending sequence in $L$ as otherwise
it would be included in either $U$ or $\overline{U}$.

Assume by Kleene's fixpoint theorem that we are given the Turing index $d$ of~$U$ as computed relative to~$\emptyset'$.
The set $U$ is built by a finite injury priority construction
with the following requirements for each $e \in \omega$:

\begin{itemize}
  \item $\Rcal_{2e}$ : $Z_e$ is finite or $(\exists a)[f(e,a) = 1 \mbox{ and } \Phi^{\emptyset'}_d(a) = 1]$
  \item $\Rcal_{2e+1}$ : $Z_e$ is finite or $(\exists a)[f(e,a) = 1 \mbox{ and } \Phi^{\emptyset'}_d(a) = 0]$
\end{itemize}

Notice again that each of these requirement is $\Sigma^f_2$, and furthermore we can effectively find an index
for each as such. Therefore, for each $i \in \omega$, we can effectively find an $m_i$ such that $R_i$
is satisfied if and only if $m_i \in f''$. By two applications of Shoenfield's limit Lemma and
low${}_2$-ness of $f$, there exists a computable function $g : \omega^3 \to 2$ such that for all $m \in \omega$,
we have $m \in f'' \biimp \lim_t \lim_s g(m, s, t) = 1$ and $m \not \in f'' \biimp \lim_t \lim_s g(m, s, t) = 0$.
Notice that for all $i \in \omega$, 
$$
R_i \mbox{ is satisfied } \biimp \lim_t \lim_s g(m_i, s, t) = 1
$$

At stage $0$, $U_0 = \emptyset$ and every integer is a \emph{leader} and \emph{follows} itself.
We say that $\Rcal_i$ \emph{requires attention for $u$ at stage $s$}
if $i \leq u \leq s$, $u$ is \emph{leader} and $g(m_i, s, u) = 0$.
At stage $s+1$, assume we have decided $u <_L v$ or $u >_L v$ for every $u, v < s$.
Set $u <_L s$ if $u \in U_s$ and $u >_L s$ if $u \not \in U_s$.
Initially set $U_{s+1} = U_s$. For each leader $u \leq s$ which has not been claimed at stage $s+1$
and for which some requirement $\Rcal_i$, $i < u$ requires attention, say that the least such $\Rcal_i$ \emph{claims} $u$
and act as follows.
\begin{itemize}
	\item[(a)] If $i = 2e$ and $u \not \in U_s$, then add $[u,s]$ to $U_{s+1}$,
	where the interval~$[u,s]$ is taken in the usual order on~$\omega$ and not in~$<_L$.
	Elements of $[u+1, s]$ \emph{follow~$u$} and are no more considered as leaders from now on and at any further stage.
	\item[(b)] If $i = 2e+1$ and $u \in U_s$, then remove $[u,s]$ from $U_{s+1}$.
	Similarly, elements of $[u+1, s]$ are no more leaders and \emph{follow~$u$}.
\end{itemize}
Then go to the next leader $u \leq s$. This ends the construction.
An immediate verification shows that at every stage, 
\begin{itemize}
	\item if $u$ stops being a leader it never becomes again a leader
	\item if $u$ follows $v$ then $v \leq u$, $v$ is a leader, every $w$ between $v$ and $u$ follows $v$
	and thus $u$ will never follow any $w > v$.
\end{itemize}
So the leader that $u$ follows eventually stabilizes. Moreover, because $g$ is limit-computable,
each leader eventually stops increasing its number of followers and therefore
there are infinitely many leaders.

\begin{claim}
$L$ is a linear order.
\end{claim}
\begin{proof*}
As $L$ is a tournament, it suffices to check there is no 3-cycle.
By symmetry, we check only the case where $u <_L s <_L v <_L u$ forms a 3-cycle
with $s$ the maximal element in $<_\omega$ order. By construction, this means that
$u \in U_s$, $v \not \in U_s$.
If $u <_\omega v$, then $u \not \in U_v$ and so there exists a leader $w \leq_\omega u$ and an even number $i \leq w$
such that $\Rcal_i$ requires attention for $w$ at a stage $t \geq v$. Case (a) of the construction
applies and the interval $[w+1, t]$ is included $U$ at least until stage $s$. 
As $v \in [w+1, t]$, $v \in U_s$ contradicting our hypothesis.
Case $u >_\omega v$ is symmetric.
\end{proof*}

\begin{claim}
$U$ is $\Delta^0_2$.
\end{claim}
\begin{proof*}
Suppose for the sake of contradiction that there exists a least element $u$
entering $U$ and leaving it infinitely many times.
Such a~$u$ must be a leader, otherwise it would not be the least one.
Let $\Rcal_i$ be the least requirement claiming $u$ infinitely many times.
As $\lim_s g(m_i, s, u)$ exists, it will claim $u$ cofinitely many times and therefore
$u$ will be in $U$ or in $\overline{U}$ cofinitely many times. Contradiction.
\end{proof*}

It immediately follows that $L$ is stable.

\begin{claim}
Every requirement $\Rcal_i$ is satisfied.
\end{claim}
\begin{proof*}
By induction over $R_i$ in priority order. 
Suppose that $R_j$ is satisfied for all $j < i$, but $\Rcal_i$ is not satisfied.
Then there exists a threshold $t_0 \geq i$ such that $\lim_s g(m_j, s, t) = 1$
for all $j < i$ and $\lim_s g(m_i, s, t) = 0$ whenever $t \geq t_0$.

Then for every leader $u \geq t_0$, $\Rcal_i$ will claim $u$ cofinitely many times,
and therefore $u$ will be in $U$ if $i$ is even and in $\overline{U}$ if $i$ is odd.
As every element follows the least leader below itself, every $v$ above the least leader
greater than $t_0$ will be in $U$ if $i$ is even and in $\overline{U}$ if $i$ is odd.
So if $Z_e$ is infinite, there will be such a $v \in Z_e$ satisfying $\Rcal_i$.
Contradiction.
\end{proof*}
\end{proof}

\section{A low${}_2$ degree bounding the Erd\H{o}s-Moser theorem}\label{sect:degrees-bounding-em-low2-degree}

The Erd\H{o}s-Moser theorem does not fall into the scope of Theorem~\ref{thm:sts-semo-sads-bounding}
since neither~$\sads$ nor~$\sts^2$ computably reduces to~$\emo$ (see Corollary~\ref{cor:emo-wkl-sts2-sads}).
In fact, the converse holds: there is a low${}_2$ degree bounding the Erd\H{o}s-Moser theorem.
We provide in the next subsections two different proofs
of the existence of a low${}_2$ degree bounding $\emo$.
More precisely, we construct a low${}_2$ set $G$
which is, up to finite changes, transitive for every infinite computable tournament.

By Theorem~\ref{thm:emo-implies-sts2-or-coh}, $[\sts^2 \vee \coh] \leq_c \emo$.
Therefore every low${}_2$ degree bounding $\emo$ also bounds $\coh$.
The proof does not seem adaptable to prove that $\coh$ is a consequence of $\emo$
even in $\omega$-models. However we can prove a weaker statement:

\begin{lemma}
For every set $X$, there exists an infinite $X$-computable tournament $T$
such that for every infinite $T$-transitive subtournament $U$,
$U \subseteq^{*} X$ or $U \subseteq^{*} \overline{X}$.
\end{lemma}
\begin{proof}
Fix a set $X$. We define a tournament $T$ as follows:
For each $a < b$, set $T(a, b)$ to hold iff $a \in X$ and $b \in X$ or $a \not \in X$ and $b \not \in X$.
Suppose for the sake of contradiction that $U$ is an infinite transitive subtournament of $T$
which intersects infinitely often $X$ and $\overline{X}$.
Take any $a, c \in U \cap X$ and $b, d \in U \cap \overline{X}$ such that $a < b < c < d$.
Then $T(a, c)$, $T(c, b)$, $T(b, d)$ and $T(d, a)$ hold contradicting transitivity of $U$.
\end{proof}

Using the previous lemma, the constructed set $G$ must be cohesive
and therefore provides another proof of the existence of a low${}_2$ cohesive set.
Finally, we can deduce a statement slightly weaker than Theorem~\ref{thm:stsn-bounding}
simply by the existence of a low${}_2$ degree bounding $\emo$.

\begin{lemma}
There exists a set $C$ such that
there is no low${}_2$ over $C$ degree $\textbf{d} \gg_{\sads} C$.
\end{lemma}
\begin{proof}
Fix a low${}_2$ set $C \gg_{\emo} \emptyset$
and a set $X$ low${}_2$ over $C$. By low${}_2$-ness of $C$, $X$ is low${}_2$.
Consider the stable coloring $f : [\omega]^2 \to 2$ constructed by Mileti in~\cite[Corollary 5.4.5]{Mileti2004Partition},
such that $X$ computes no infinite $f$-homogeneous set.
We can see $f$ as a stable tournament $T$ such that for each $x < y$,
$T(x,y)$ holds iff $f(x,y) = 1$. As $C \gg_{\emo} \emptyset$, there exists an infinite $C$-computable
transitive subtournament $U$ of $T$. $U$ is a stable linear order
such that every infinite ascending or descending sequence is $f$-homogeneous.
Therefore $X$ computes no infinite ascending or descending sequence in $U$.
\end{proof}

The following question remains open:

\begin{question}
Does $\emo$ admit a universal instance?
\end{question}

\subsection{A low${}_2$ degree bounding $\emo$ using first jump control}

The following theorem uses the proof techniques introduced in~\cite{Cholak2001strength} for producing
low${}_2$ sets by controlling the first jump. It is done in the same spirit
as Theorem~3.6 in~\cite{Cholak2001strength}.

\begin{theorem}\label{thm:low2-degree-bounding-em-first-jump}
For every set $P \gg \emptyset'$, there exists a set $G \gg_{\emo} \emptyset$
such that $G' \leq_T P$.
\end{theorem}
\begin{proof}
We will prove Theorem~\ref{thm:low2-degree-bounding-em-first-jump} using a variation
of Erd\H{o}s-Moser conditions. Recall that an EM condition for an infinite tournament~$T$
is a Mathias condition $(F, X)$ where
\begin{itemize}
	\item[(a)] $F \cup \{x\}$ is $T$-transitive for each $x \in X$
	\item[(b)] $X$ is included in a minimal $T$-interval of $F$.
\end{itemize}
The main properties of an EM condition are proven in section~\ref{sect:weakness-erdos-moser-theorem}
under Lemma~\ref{lem:emo-cond-beats} and Lemma~\ref{lem:emo-cond-valid}.

Let $C$ be a low set such that there exists a uniformly $C$-computable enumeration
$\vec{T}$ of infinite tournaments containing every computable tournament. Note that $P \gg C'$.
Our forcing conditions are tuples $(\sigma, F, X)$ where $\sigma \in \omega^{<\omega}$ and the following holds:
\begin{itemize}
	\item[(a)] $(F,X)$ forms a Mathias condition and $X$ is a set low over $C$.
	\item[(b)] $(F \setminus [0, \sigma(\nu)], X)$ is an EM condition for $T_\nu$ for each $\nu < |\sigma|$.
\end{itemize}
A condition $(\tilde{\sigma}, \tilde{F}, \tilde{X})$ \emph{extends} a condition $(\sigma, F, X)$
if $\sigma \preceq \tilde{\sigma}$ and $(\tilde{F},\tilde{X})$ Mathias extends $(F,X)$.
A set $G$ \emph{satisfies} the condition $(\sigma, F, X)$ if
$G \setminus [0, \sigma(\nu)]$ is $T_\nu$-transitive for each $\nu < |\sigma|$ and $G$ satisfies the Mathias condition $(F,X)$.
An \emph{index} of a condition $(\sigma, F, X)$ is a code of the tuple $\tuple{\sigma, F, e}$ where $e$ is a lowness index of~$X$.

The first lemma simply states that we can ensure 
that $G$ will be infinite and eventually transitive for each tournament in $\vec{T}$.

\begin{lemma}\label{lem:infinite-em-first-jump-control}
For every condition $c = (\sigma, F, X)$ and every $i, j \in \omega$, one can
$P$-compute an extension
$(\tilde{\sigma}, \tilde{F}, \tilde{X})$ such that $|\tilde{\sigma}| \geq i$ and $|\tilde{F}| \geq j$
uniformly from $i$, $j$ and an index of~$c$. 
\end{lemma}
\begin{proof}
Let $x$ be the first element of $X$. As $X$ is low over $C$, $x$ can be found $C'$-computably
from a lowness index of $X$.
The condition $(\tilde{\sigma}, F, X)$ is a valid extension of $c$ 
where $\tilde{\sigma} = \sigma^\frown x \dots x$ so that $|\tilde{\sigma}| \geq i$.
It suffices to prove that we can $C'$-compute an extension $(\tilde{\sigma}, \tilde{F}, \tilde{X})$
with $|\tilde{F}| > |F|$ and iterate the process.
Define the computable coloring $g : X \to 2^{|\tilde{\sigma}|}$
by $g(s) = \rho$ where $\rho \in 2^{|\tilde{\sigma}|}$ such that $\rho(\nu) = 1$ iff $T_\nu(x, s)$ holds.
One can find uniformly in $P$ a $\rho \in 2^{|\tilde{\sigma}|}$ such that the following $C$-computable set is infinite:
$$
Y = \{ s \in X \setminus \{x\} : g(s) = \rho\} 
$$
By Lemma~\ref{lem:emo-cond-valid}, $((F \cup \{x\}) \setminus [0,\tilde{\sigma}(\nu)], Y)$
is a valid EM extension for $T_\nu$. As $Y$ is low over $C$,
$(\tilde{\sigma}, F \cup \{x\}, Y)$ is a valid extension for~$c$.
\end{proof}

It remains to be able to decide $e \in (G \oplus C)'$ uniformly in $e$. 
We first need to define a forcing relation.

\begin{definition}Fix a condition $c = (\sigma, F, X)$ and two integers $e$ and $x$.
\begin{itemize}
	\item[1.] $c \Vdash \Phi^{G \oplus C}_e(x) \uparrow$
	if $\Phi^{(F \cup F_1) \oplus C}_e(x) \uparrow$ for all finite subsets $F_1 \subseteq X$
	such that $F_1$ is $T_\nu$-transitive simultaneously for each $\nu < |\sigma|$.
	\item[2.] $c \Vdash \Phi^{G \oplus C}_e(x) \downarrow$
	if $\Phi^{F \oplus C}_e(x) \downarrow$.
\end{itemize}
\end{definition}

Note that the way we defined our forcing relation $c \Vdash \Phi^{G \oplus C}_e(x) \uparrow$ differs
slightly from the ``true'' forcing notion $\Vdash^{*}$ inherited by the notion of satisfaction of $G$.
The true forcing definition of this statement is the following:

$c \Vdash^{*} \Phi^{G \oplus C}_e(x) \uparrow$ if $\Phi^{(F \cup F_1) \oplus C}_e(x) \uparrow$ 
for all finite \emph{extensible} subsets $F_1 \subseteq X$
such that $F_1$ is $T_\nu$-transitive simultaneously for each $\nu < |\sigma|$, i.e., 
for all finite subsets $F_1 \subseteq X$ such that there exists an extension
$d = (\tilde{\sigma}, F \cup F_1, \tilde{X})$.

However $c \Vdash^{*} \Phi^{G \oplus C}_e(x) \uparrow$ is not a $\Pi^0_1$ statement
whereas $c \Vdash \Phi^{G \oplus C}_e(x) \uparrow$ is.
In particular the fact that $c \not \Vdash \Phi^{G \oplus C}_e(x) \uparrow$
does not mean that $c$ has an extension forcing its negation.
This subtlety is particularly important in Lemma~\ref{lem:force-jump-em-first-jump-control}.
The following lemma gives a sufficient constraint, namely being included in a part of a particular partition,
on finite transitive sets to ensure that they are \emph{extensible}.

\begin{lemma}\label{lem:part1-em-first-jump-control}
Let $c = (\sigma, F, X)$ be a condition and $E \subseteq X$ be a finite set. 
There exists a $2^{|\sigma|}$ partition $(E_\rho : \rho \in 2^{|\sigma|})$ of $E$ and an infinite set $Y \subseteq X$
low over $C$ such that $E < Y$ and for all $\rho \in 2^{|\sigma|}$ and $\nu < |\sigma|$,
if $\rho(\nu) = 0$ then $E_\rho \to_{T_\nu} Y$ and if $\rho(\nu) = 1$ then $Y \to_{T_\nu} E_\rho$.

Moreover this partition and a lowness index of $Y$ can be uniformly $P$-computed
from an index of $c$ and the set $E$.
\end{lemma}
\begin{proof}
Given a set $E$, define $P_E$ to be the finite set of ordered $2^{|\sigma|}$-partitions of $E$,
that is,
$$
P_E = \{ (E_\rho : \rho \in 2^{|\sigma|}) : \bigcup_{\rho \in 2^{|\sigma|}} E_\rho = E 
\mbox{ and } \rho \neq \xi \imp E_\rho \cap E_\xi = \emptyset\}
$$
Define the $C$-computable coloring $g : X \to P_E$ by $g(x) = (E^x_\rho : \rho \in 2^{|\sigma|})$ where
$E^x_\rho = \{ a \in E : (\forall \nu < |\sigma|)[T_\nu(a, x) \mbox{ holds iff } \rho(\nu) = 0]\}$.
On can find uniformly in $P$ a partition $(E_\rho : \rho \in 2^{|\sigma|})$
such that the following $C$-computable set is infinite:
$$
Y = \{ x \in X \setminus E : g(x) = (E_\rho : \rho \in 2^{|\sigma|}) \}
$$
By definition of $g$, for all $\rho \in 2^{|\sigma|}$ and $\nu < |\sigma|$,
if $\rho(\nu) = 0$ then $E_\rho \to_{T_\nu} Y$ and if $\rho(\nu) = 1$ then $Y \to_{T_\nu} E_\rho$.
\end{proof}

We are now ready to prove the key lemma of this forcing, stating
that we can $P$-decide whether or not $e \in G'$ for any $e \in \omega$.

\begin{lemma}\label{lem:force-jump-em-first-jump-control}
For every condition $(\sigma, F, X)$ and every $e \in \omega$, there exists
an extension $d = (\tilde{\sigma}, \tilde{F}, \tilde{X})$ such that one of the following holds:
\begin{itemize}
	\item[1.] $d \Vdash \Phi^{G \oplus C}_e(e) \downarrow$
	\item[2.] $d \Vdash \Phi^{G \oplus C}_e(e) \uparrow$
\end{itemize}
This extension can be $P$-computed uniformly 
from an index of $c$ and $e$.
Moreover there is a $C'$-computable procedure to decide which case holds from an index of $d$.
\end{lemma}
\begin{proof}
Let $k = |\sigma|$.
Using a $C'$-computable procedure, we can decide from an index of $c$ and $e$ whether
there exists a finite set $E \subset X$ such that for every $2^k$-partition
$(E_i : i < 2^k)$ of $E$, there exists an $i < 2^k$ and a subset $F_1 \subseteq E_i$
$T_\nu$-transitive simultaneously for each $\nu < k$ and
satisfying $\Phi_e^{(F \cup F_1) \oplus C}(e) \downarrow$.
\begin{itemize}
	\item[1.] If such a set $E$ exists, it can be $C'$-computably found. 
	By Lemma~\ref{lem:part1-em-first-jump-control}, one can $P$-computably find
	a $2^k$-partition $(E_\rho : \rho \in 2^k)$ of $E$ and a set $Y \subseteq X$ low over $C$
	such that for all $\rho \in 2^k$ and $\nu < k$,
	if $\rho(\nu) = 0$ then $E_\rho \to_{T_\nu} Y$ and if $\rho(\nu) = 1$ then $Y \to_{T_\nu} E_\rho$.
	We can $C'$-computably find a $\rho \in 2^k$ and a set $F_1 \subseteq E_\rho$
	which is $T_\nu$-transitive simultaneously for each $\nu < k$ and
	satisfying $\Phi_e^{(F \cup F_1) \oplus C}(e) \downarrow$. 
	By Lemma~\ref{lem:emo-cond-valid}, $(F \setminus [0,\sigma(\nu)]) \cup F_1, Y)$
	is a valid EM extension of $(F \setminus [0, \sigma(\nu)], X)$ for $T_\nu$ for each $\nu < k$.
	As $Y$ is low over $C$, $(\sigma, F \cup F_1, Y)$ 
	is a valid extension of~$c$ forcing $\Phi_e^{G \oplus C}(e) \downarrow$.

	\item[2.] If no such set exists, then by compactness, the $\Pi^{0,C}_1$ class
	of all $2^k$-partitions $(X_i : i < 2^k)$ of $X$ such that for every $i < 2^k$ 
	and every finite set $F_1 \subseteq X_i$ which is $T_\nu$-transitive simultaneously for each $\nu < k$,
	$\Phi_e^{(F \cup F_1) \oplus C}(e) \uparrow$ is non-empty. In other words,
	the $\Pi^{0,C}_1$ class of all $2^k$-partitions $(X_i : i < 2^k)$ of~$X$
	such that for every $i < 2^k$, $(\sigma, F, X_i) \Vdash \Phi_e^{G \oplus C}(e) \uparrow$
	is non-empty. By the relativized low basis theorem, there exists
	a $2^k$-partition $(X_i : i < 2^k)$ of $X$ low over $C$. Furthermore, a lowness index for this partition
	can be uniformly $C'$-computably found. Using~$P$, one can find
	an $i < 2^k$ such that $X_i$ is infinite.
	$(\sigma, F, X_i)$ is a valid extension of~$c$ forcing $\Phi_e^{G \oplus C}(e) \uparrow$.
\end{itemize}
\end{proof}

Using Lemma~\ref{lem:infinite-em-first-jump-control} 
and Lemma~\ref{lem:force-jump-em-first-jump-control}, one can $P$-compute an infinite
decreasing sequence of conditions $c_0 = (\epsilon, \emptyset, \omega) \geq c_1 \geq \dots$
such that for each $s > 0$
\begin{itemize}
	\item[1.] $|\sigma_s| \geq s$, $|F_s| \geq s$
	\item[2.]	$c_s \Vdash \Phi_s^{G \oplus C}(s) \downarrow$ or $c_s \Vdash \Phi_s^{G \oplus C}(s) \uparrow$
\end{itemize}
where $c_s = (\sigma_s, F_s, X_s)$.
The resulting set $G = \bigcup_s F_s$ is $T_\nu$-transitive up to finite changes for each $\nu \in \omega$
and $G' \leq_T P$.
\end{proof}

\subsection{A low${}_2$ degree bounding $\emo$ using second jump control}
We now use the second proof technique used in~\cite{Cholak2001strength} for producing a low${}_2$ set.
It consists of directly controlling the second jump of the produced set.

\begin{theorem}\label{thm:low2-degree-bounding-em-second-jump}
There exists a low${}_2$ degree bounding $\emo$.
\end{theorem}
\begin{proof}
Similar to Theorem~\ref{thm:low2-degree-bounding-em-first-jump}, 
we fix a low set $C$ such that there exists a uniformly $C$-computable enumeration
$\vec{T}$ of infinite tournaments containing every computable tournament. In particular $P \gg C'$.

Our forcing conditions are the same as in Theorem~\ref{thm:low2-degree-bounding-em-first-jump}.
We can release the constraints of infinity and lowness over $C$ for $X$
in a condition $(\sigma, F, X)$. This gives the notion of \emph{precondition}.
The forcing relations extend naturally to preconditions.

\begin{definition}\label{def:em-smallness}
Fix a finite set of Turing indices $\vec{e}$.
A condition $(\sigma, F,X)$ is \emph{$\vec{e}$-small} if there exists a number~$x$
and a sequence $(\sigma_i, F_i, X_i : i < n)$ such that for each $i < n$
\begin{itemize}
	\item[(i)] $(\sigma_i, F_i, X_i)$ is a precondition extending $c$
	\item[(ii)] $(X_i : i < n)$ is a partition of $X \cap (x,+\infty)$
	\item[(iii)] $max(X_i) < x \mbox{ or } (\sigma_i, F \cup F_i, X_i) \Vdash (\exists e \in \vec{e})(\exists y < x)\Phi_e^{G \oplus C}(y) \uparrow$
\end{itemize}
A condition is \emph{$\vec{e}$-large} if it is not $\vec{e}$-small.
\end{definition}

A condition $(\tilde{\sigma}, \tilde{F}, \tilde{X})$ is a \emph{finite extension}
of $(\sigma, F, X)$ if $\tilde{X} =^{*} X$. Finite extensions
do not play the same fundamental role as in the original forcing in~\cite{Cholak2001strength}
as adding elements to the set $F$ may require to remove infinitely many elements
of the promise set~$X$ to obtain a valid extension.
We nevertheless prove the following classical lemma.

\begin{lemma}\label{lem:em-finite-extension-second-jump}
Fix an $\vec{e}$-large condition $c = (\sigma, F, X)$.
\begin{itemize}
	\item[1.] If $\vec{e'} \subseteq \vec{e}$ then $c$ is $\vec{e'}$-large.
	\item[2.] If $d$ is a finite extension of $c$ then $d$ is $\vec{e}$-large.
\end{itemize}
\end{lemma}
\begin{proof}
Clause 1 is trivial as $\vec{e}$ appears only in a universal quantification in the definition
of $\vec{e}$-largeness. We prove clause 2.
Let $d = (\tilde{\sigma}, \tilde{F}, \tilde{X})$ be an $\vec{e}$-small finite extension of $c$.
We will prove that $c$ is $\vec{e}$-small.
Let $x \in \omega$ and $(\sigma_i, F_i, X_i : i < n)$ witness $\vec{e}$-smallness of $d$.
Let $y = max(x, X \setminus \tilde{X})$. 
For each $i < n$, set $\tilde{X}_i = X_i \cap (y, +\infty)$.
Then $y$ and $(\sigma_i, F_i, \tilde{X}_i : i < n)$ witness $\vec{e}$-smallness of $c$.
\end{proof}

\begin{lemma}\label{lem:decide-smallness-em-second-jump}
There exists a $C''$-effective procedure to decide, given an index of
a condition~$c$ and a finite set of Turing indices $\vec{e}$,
whether $c$ is $\vec{e}$-large.
Furthermore, if $c$ is $\vec{e}$-small, there exists sets $(X_i : i < n)$ low over~$C$
witnessing this, and one may $C'$-compute a value of $n$, $x$, lowness indices
for $(X_i : i < n)$ and the corresponding sequences $(\sigma_i, F_i, X_i : i < n)$
which witness that $c$ is $\vec{e}$-small.
\end{lemma}
\begin{proof}
Fix a condition~$c = (\sigma, F, X)$ 
The predicate ``$(\sigma, F, X)$ is $\vec{e}$-small'' can be expressed
as a $\Sigma^0_2$ statement
$$
(\exists z)(\exists Z)P(z, Z, F, X, \vec{\nu}, \vec{e})
$$
where $P$ is a $\Pi^{0,C}_1$ predicate. Here $z$ codes $n$ and $x$,
and $Z$ codes $(X_i : i < n)$. The predicate $(\exists Z)P(z, Z, F, X, \sigma, \vec{e})$
is $\Pi^{0,C \oplus X}_1$ by compactness. As $X$ is low over $C$ and $F$ and $\sigma$ are finite,
one can compute a $\Delta^{0, C}_2$ index for the same predicate $P$
with parameter $z$, an index of $c$ and $\vec{e}$,
from a lowness index for $X$, $F$ and $\sigma$.
Therefore there exists a $\Sigma^{0,C}_2$ statement with parameters
an index of $c$ and $\vec{e}$ which holds iff $c$ is $\vec{e}$-small.

If $c$ is $\vec{e}$-small, there exists sets $(X_i : i < n)$ low over $X$ (hence low over $C$) witnessing
it by the low basis theorem relativized to $C$. By the uniformity of the proof of the low basis theorem,
one can compute lowness indices of $(X_i : i < n)$ uniformly from a lowness index of~$X$.
\end{proof}

As the extension produced in Lemma~\ref{lem:infinite-em-first-jump-control} is not a finite extension,
we need to refine it to ensure largeness preservation.

\begin{lemma}\label{lem:infinite-em-second-jump-control}
For every $\vec{e}$-large condition $c = (\sigma, F, X)$ and every $i, j \in \omega$, 
one can $P$-compute an $\vec{e}$-large extension
$(\tilde{\sigma}, \tilde{F}, \tilde{X})$ such that $\tilde{\sigma} \geq i$ and $|\tilde{F}| \geq j$
uniformly from an index of $c$, $i$, $j$ and $\vec{e}$. 
\end{lemma}
\begin{proof}
Let $x$ be the first element of $X$. As $X$ is low over $C$, $x$ can be found $C'$-computably
from a lowness index of $X$.
The condition $d = (\tilde{\sigma}, F, X)$ is a valid extension of $c$ 
where $\tilde{\sigma} = \sigma^\frown x \dots x$ so that $|\tilde{\sigma}| \geq i$.
As $d$ is a finite extension of $c$, it is $\vec{e}$-large by Lemma~\ref{lem:em-finite-extension-second-jump}.
It suffices to prove that we can $C'$-compute an $\vec{e}$-large extension $(\tilde{\sigma}, \tilde{F}, \tilde{X})$
with $|\tilde{F}| > |F|$ and iterate the process.
Define the $C$-computable coloring $g : X \to 2^{|\tilde{\sigma}|}$ as in Lemma~\ref{lem:infinite-em-first-jump-control}.
For each $\rho \in 2^{|\tilde{\sigma}|}$, define the following set:
$$
Y_\rho = \{ s \in X \setminus \{x\} : g(s) = \rho\} 
$$
There must be a $\rho \in 2^{|\tilde{\sigma}|}$ such that $Y_\rho$ is infinite 
and $(\tilde{\sigma}, F \cup \{x\}, Y_\rho)$ is $\vec{e}$-large,
otherwise the witnesses of $\vec{e}$-smallness for each $\rho \in 2^{|\tilde{\sigma}|}$ would witness $\vec{e}$-smallness of~$c$.
By Lemma~\ref{lem:decide-smallness-em-second-jump}, one can $C''$-find a $\rho \in 2^{|\tilde{\sigma}|}$ such that
$(\tilde{\sigma}, F \cup \{x\}, Y_\rho)$ is $\vec{e}$-large.
As seen in Lemma~\ref{lem:infinite-em-first-jump-control}, $(\tilde{\sigma}, F, \{x\}, Y_\rho)$ is a valid extension.
\end{proof}

The following lemma is a refinement of Lemma~\ref{lem:part1-em-first-jump-control}
controlling largeness preservation.

\begin{lemma}\label{lem:part1-em-second-jump-control}
Let $c = (\sigma, F, X)$ be an $\vec{e}$-large condition and $E \subseteq X$ be a finite set. 
There is a $2^{|\sigma|}$ partition $(E_\rho : \rho \in 2^{|\sigma|})$ of $E$ and an infinite set $Y \subseteq X$
low over $C$ such that $E < Y$ and
\begin{itemize}
	\item[1.] for all $\rho \in 2^{|\sigma|}$ and $\nu < |\sigma|$, if $\rho(\nu) = 0$ 
	then $E_\rho \to_{T_\nu} Y$ and if $\rho(\nu) = 1$ then $Y \to_{T_\nu} E_\rho$.
	\item[2.] $(\sigma, F \cup F_1, Y)$ is an $\vec{e}$-large condition extending $c$ for every $\rho \in 2^{|\sigma|}$
	and every finite set $F_1 \subseteq E_\rho$ which is $T_\nu$-transitive for each $\nu < |\sigma|$
\end{itemize}
Moreover this partition and a lowness index of $Y$ can be uniformly $C''$-computed
from an index of $c$ and the set $E$.
\end{lemma}
\begin{proof}
Given a set $E$, recall from Lemma~\ref{lem:part1-em-first-jump-control} 
that $P_E$ is the finite set of ordered $2^k$-partitions of $E$.
Define again the computable coloring $g : X \to P_E$ by $g(x) = (E^x_\rho : \rho \in 2^{|\sigma|})$ where
$E^x_\rho = \{ a \in E : (\forall \nu < |\sigma|)[T_\nu(a, x) \mbox{ holds iff } \rho(\nu) = 0]\}$.
If for each partition $(E_\rho : \rho \in 2^{|\sigma|})$, there exists a $\rho \in 2^{|\sigma|}$
and a $F_1 \subseteq E_\rho$ which is $T_\nu$-transitive simultaneously for each $\nu < |\sigma|$
and such that $(\sigma, F \cup F_1, Y)$ is $\vec{e}$-small where
$$
Y = \{ x \in X \setminus E : g(x) = (E_\rho : \rho \in 2^{|\sigma|}) \}
$$
Then we could construct a witness of $\vec{e}$-smallness of $c$ using smallness witnesses
of $(\sigma, F \cup F_1, Y)$ for each partition $(E_\rho : \rho \in 2^{|\sigma|})$.
Therefore there must exist a partition $(E_\rho : \rho \in 2^{|\sigma|})$
such that $Y$ is infinite and $d = (\sigma, F \cup F_1, Y)$ is $\vec{e}$-large for every $\rho \in 2^{|\sigma|}$
and every $F_1 \subseteq E_\rho$ which is $T_\nu$-transitive for each $\nu < |\sigma|$.

By Lemma~\ref{lem:decide-smallness-em-second-jump}, such a partition can be found $C''$-computably.
By definition of $g$, for all $\rho \in 2^{|\sigma|}$ and $\nu < k$,
if $\rho(\nu) = 0$ then $E_\rho \to_{T_\nu} Y$ and if $\rho(\nu) = 1$ then $Y \to_{T_\nu} E_\rho$.
Therefore, by Lemma~\ref{lem:emo-cond-valid}, $((F \setminus [0,\sigma(\nu)]) \cup F_1, Y)$
is a valid EM extension of $(F \setminus [0, \sigma(\nu)], X)$ for $T_\nu$ for each $\nu < |\sigma|$, 
so $d$ is a valid condition.
\end{proof}

\begin{lemma}\label{lem:em-large-to-instance}
Suppose that $c = (\sigma, F, X)$ is $\vec{e}$-large.
For every $y \in \omega$ and $e \in \vec{e}$, 
there exists an $\vec{e}$-large extension $d$
such that $d \Vdash \Phi_e^{G \oplus C}(y) \downarrow$.
Furthermore, an index for $d$ can be $C'$-computed uniformly in an index of $c$, $e$ and $y$.
\end{lemma}
\begin{proof}
Let $k = |\sigma|$.
As $c$ is $\vec{e}$-large, then by a compactness argument,
there exists a finite set $E \subset X$ such that for every $2^k$-partition
$(E_i : i < 2^k)$ of $E$, there exists an $i < k$
and a finite subset $F_1 \subseteq E_i$ which is $T_\nu$-transitive
simultaneously for each $\nu < k$, and $\Phi_e^{(F \cup F_1) \oplus C}(y) \downarrow$.
Moreover this set $E$ can be $C'$-computably found.
By Lemma~\ref{lem:part1-em-second-jump-control},
on can uniformly $C''$-find a partition $(E_\rho : \rho \in 2^k)$
of $E$ and a lowness index for an infinite set $Y \subseteq X$
low over~$C$ such that
\begin{itemize}
	\item[1.] for all $\rho \in 2^k$ and $\nu < k$, if $\rho(\nu) = 0$ 
	then $E_\rho \to_{T_\nu} Y$ and if $\rho(\nu) = 1$ then $Y \to_{T_\nu} E_\rho$.
	\item[2.] $(\sigma, F \cup F_1, Y)$ is an $\vec{e}$-large condition extending~$c$ for every $\rho \in 2^k$
	and every finite set finite set $F_1 \subseteq E_\rho$ which is $T_\nu$-transitive for each $\nu < k$
\end{itemize}
We can then produce by a $C'$-computable search
a $\rho \in 2^k$ and a finite set $F_1 \subseteq E_\rho$
which is $T_\nu$-transitive for each $\nu < k$ and such that
$\Phi_e^{(F \cup F_1) \oplus C}(y) \downarrow$.
By Lemma~\ref{lem:emo-cond-valid}, $((F \setminus [0,\sigma(\nu)]) \cup F_1, Y)$
is a valid EM extension of $(F \setminus [0,\sigma(\nu)], X)$ for $T_\nu$ for each $\nu < k$.
As $Y$ is low over~$C$, $(\sigma, F \cup F_1, Y)$ is a valid $\vec{e}$-large extension.
\end{proof}

\begin{lemma}\label{lem:large-small-force-em-second-jump}
Suppose that $c = (\sigma, F, X)$ is $\vec{e}$-large
and $(\vec{e} \cup \{u\})$-small. There exists a $\vec{e}$-large 
extension $d$ such that $d \Vdash \Phi_u^{G \oplus C}(y) \uparrow$ for some $y \in \omega$.
Furthermore one can find an index for $d$ by applying
a $C''$-computable function to an index of $c$, $\vec{e}$ and $u$.
\end{lemma}
\begin{proof}
By Lemma~\ref{lem:decide-smallness-em-second-jump}, we may choose the sets $(X_i : i < n)$
witnessing that $c$ is $(\vec{e} \cup \{u\})$-small
to be low over $C$. Fix the corresponding $x$ and $(\sigma_i, F_i : i < n)$.
Consider the $i$'s such that $(\sigma_i, F_i, X_i) \Vdash \Phi_u^{G \oplus C}(y) \uparrow$ for some $y < x$.
As $c$ is $\vec{e}$-large, there must be such an~$i < n$ such that
$(\sigma_i, F_i, X_i)$ is an $\vec{e}$-large condition. By Lemma~\ref{lem:decide-smallness-em-second-jump}
we can find $C''$-computably such an~$i < n$.
$(\sigma_i, F_i, X_i)$ is the desired extension. 
\end{proof}

Using the previous lemmas, we can $C''$-compute an infinite descending sequence
of conditions $c_0 = (\epsilon, \emptyset, \omega) \geq c_1 \geq \dots$
together with an infinite increasing sequence of Turing indices $\vec{e}_0 = \emptyset \subseteq \vec{e}_1 \subseteq \dots$ 
such that for each $s > 0$
\begin{itemize}
	\item[1.] $|\sigma_s| \geq s$, $|F_s| \geq s$, $c_s$ is $\vec{e}_s$-large
	\item[2.] Either $s \in \vec{e}_s$ or $c_s \Vdash \Phi_s^{G \oplus C}(y) \uparrow$ for some $y \in \omega$
	\item[3.] $c_s \Vdash \Phi_e^{G \oplus C}(x) \downarrow$ if $s = \tuple{e,x}$ and $e \in \vec{e}_s$\end{itemize}
where $c_s = (\sigma_s, F_s, X_s)$.
The resulting set $G = \bigcup_s F_s$ is $T_\nu$-transitive up to finite changes 
simultaneously for each $\nu \in \omega$
and $G'' \leq_T C'' \leq_T \emptyset''$.
\end{proof}

\section{No incomplete $\Delta^0_2$ degree bounds~$\srrt^2_2$}

Mileti proved in~\cite{Mileti2004Partition} that the only $\Delta^0_2$ degree bounding
$\srt^2_2$ is $\mathbf{0}'$. Using the fact that every $\Delta^0_2$ set has an infinite incomplete
$\Delta^0_2$ subset in either it or its complement \cite{Hirschfeldt2008strength},
we obtain another proof that $\srt^2_2$ admits no universal instance.
 
In this section, we extend Mileti's theorem by proving that the only $\Delta^0_2$ degree bounding
a stable version of the rainbow Ramsey theorem for pairs ($\srrt^2_2$) is $\mathbf{0}'$
and deduce several results about which stable statements admit a universal instance.

Recall that~$\srrt^2_2$ is the restriction of the rainbow Ramsey theorem for pairs
in which everybody either gets married or becomes a monk. 
See Definition~\ref{def:stable-rainbow-ramsey-theorem} for a formal presentation.
The stable rainbow Ramsey theorem for pairs admits several computably equivalent characterizations.
We are in particular interested in a characterization which can be seen as a stable notion of diagonal non-computability
and proven in Theorem~\ref{thm:srrt22-characterizations}.

\begin{definition}
Given a function $f : \omega \to \omega$, a function $g$ is \emph{$f$-diagonalizing}
if $(\forall x)[f(x) \neq g(x)]$.
$\sdnrs{2}$ is the statement ``Every $\Delta^0_2$
function $f : \omega \to \omega$ has an $f$-diagonalizing function''.
\end{definition}

The following theorem extends Mileti's result to $\sdnrs{2}$.
As $\sdnrs{2}$ is computably below many stable principles,
we shall deduce a few more non-universality results.

\begin{theorem}\label{thm:incomplete-srrt22}
For every $\Delta^0_2$ incomplete set $X$,
there exists a $\Delta^0_2$ function $f : \omega \to \omega$
such that $X$ computes no $f$-diagonalizing function.
\end{theorem}

\begin{corollary}
A $\Delta^0_2$ degree $\mathbf{d}$ bounds $\srrt^2_2$ iff $\mathbf{d} = \mathbf{0}'$.
\end{corollary}
\begin{proof}
As $\srrt^2_2 \leq_c \srt^2_2$, any computable instance of $\srrt^2_2$ has a $\Delta^0_2$ solution.
So~$\mathbf{0}'$ bounds $\srrt^2_2$. If $\mathbf{d}$ is incomplete, then by Theorem~\ref{thm:incomplete-srrt22}
and by $\srrt^2_2 =_c \sdnrs{2}$, there is a computable instance of $\srrt^2_2$ such that $\mathbf{d}$ bounds no solution.
\end{proof}

\begin{corollary}
No statement $\Psf$ such that $\srrt^2_2 \leq_c \Psf \leq_c \srt^2_2$ admits a universal instance.
\end{corollary}
\begin{proof}
By \cite[Corollary 4.6]{Hirschfeldt2008strength} every $\Delta^0_2$ set or its complement
has an incomplete $\Delta^0_2$ infinite subset. As $\Psf \leq_c \srt^2_2 \leq_c D^2_2$,
every computable instance $U$ of $\Pcal$ has a $\Delta^0_2$ incomplete solution $X$.
By Theorem~\ref{thm:incomplete-srrt22}, there exists a 
computable coloring $f : [\omega]^2 \to \omega$ such that $X$ computes no infinite $f$-rainbow.
As $\srrt^2_2 \leq_c \Psf$, there exists a computable instance of $\Psf$ such that $X$ does not compute a solution
to it. Hence $U$ is not a universal instance of $\Psf$.
\end{proof}

\begin{corollary}
None of $\srrt^2_2$, $\semo$, $\sts^2$ and $\sfs^2$ admits a universal instance.
\end{corollary}

\begin{proof}[Proof of Theorem~\ref{thm:incomplete-srrt22}]
The proof is an adaptation of \cite[Theorem 5.3.7]{Mileti2004Partition}.
Suppose that $D$ is a $\Delta^0_2$ incomplete set. We will construct a $\Delta^0_2$
function $f : \omega \to \omega$ such that $D$ does not compute any $f$-diagonalizing function.
We want to satisfy the following requirements for each $e \in \omega$:

\bigskip
$\Rcal_e$ : If $\Phi^D_e$ is total, then there is an $a$  such that 
$\Phi^D_e(a) = f(a)$.
\bigskip

For each $e \in \omega$, define the partial function $u_e$ by letting $u_e(a)$ be the use of $\Phi_e^D$ on input
$a$ if $\Phi_e^D(a)\da$ and letting $u_e(a) \ua$ otherwise.
We can assume w.l.o.g.\ that whenever $u_e(a) \da$ then $u_e(a) \geq a$. Also define a computable partial
function $\theta$ by letting $\theta(a) = (\mu t)[a \in \emptyset'_t]$ if $a \in \emptyset'$ and $\theta(a) \ua$ otherwise.

The local strategy for satisfying a single requirement $\Rcal_e$ works as follows. If $\Rcal_e$ receives
attention at stage $s$, this strategy does the following. First it identifies a number $a \geq e$ that is \emph{not}
restrained by strategies of higher priority such that the following conditions holds:
\begin{itemize}
  \item[(i)] $\Phi_{e,s}^{D_s}(a) \downarrow$
  \item[(ii)] $u_{e,s}(a) < max(0, \theta_s(a))$
\end{itemize}
If no such number $a$ exists, the strategy does nothing. Otherwise it puts a restraint on $a$
and \emph{commits} to assigning $f_s(a) = \Phi_{e,s}^{D_s}(a)$.
For any such $a$, this commitment will remain active as long as the strategy has a restraint on this element.
Having done all this, the local strategy is declared to be satisfied and will not act again unless 
either a strategy of higher priority puts restraints on $a$, or the value of $u_{e,s}(a)$ or $\theta_s(a)$ changes.
In both cases the strategy gets \emph{injured} and has to reset, releasing all its restraints.

To finish stage $s$, the global strategy assigns values $f_s(y)$ for all $y \leq s$ as follows:
if $y$ is commited to some value assignment of $f_s(y)$ due to a local strategy, then define $f_s(y)$ to be this value.
If not, let $f_s(y) = 0$. This finishes the construction and we now turn to the verification.

For each $e, a \in \omega$, let $Z_{e, a} = \set{s \in \omega : \Rcal_e \mbox{ restrains } a \mbox{ at stage } s}$.

\begin{claim}\label{claim:incomplete-srrt22-1}
For each $e, a \in \omega$, 
\begin{itemize}
  \item[(a)] if $\Phi^D_e(a) \uparrow$ then $Z_{e, a}$ is finite;
  \item[(b)] if $\Phi^D_e(a) \downarrow$ then $Z_{e, a}$ is either finite or cofinite.
\end{itemize}
\end{claim}
\begin{proof*}
By induction on the priority order. We consider $Z_{e, a}$, assuming that for all $\Rcal_{e'}$ of higher priority,
the set $Z_{e', a}$ is either finite or cofinite. First notice that $Z_{e, a} = \emptyset$
if $a < e$ or $a \not \in \emptyset'$, so we may assume that $a \geq e$ and $a \in \emptyset'$. 
If there exists $e' < e$ such that $Z_{e', a}$
is cofinite, then $Z_{e, a}$ is finite because at most one requirement may claim $a$ at a given stage.
Suppose that $Z_{e', a}$ is finite for all $e' < e$. Fix $t_0$
such that for all $e' < e$ and $s \geq t_0$ $\Rcal_{e'}$ does not restrain $a$ at stage $s$
and $\theta_s(a) = \theta(a)$.

Suppose that $\Phi_e^D(a) \uparrow$. Fix $t_1 \geq t_0$ such that 
$D(b) = D_s(b)$ for all $b \leq \theta(a)$ and all $s \geq t_1$.
Then for all $s \geq t_1$, if $\Phi_{e,s}^{D_s}(a) \downarrow$ 
then we must have $u_{e,s}(a) > \theta(a)$ because otherwise $\Phi^D_e(a) \downarrow$.
It follows that for all $s \geq t_1$, requirement $\Rcal_e$ does not restrain $a$ at stage $s$. Hence $Z_{e, a}$ is finite.

Suppose now that $\Phi^D_e(a) \downarrow$.
Fix $t_1 \geq t_0$ such that for all $s \geq t_1$ we have 
$\Phi^{D_s}_{e,s}(a) \downarrow$ and $D_s(c) = D(c)$ for every $c \leq u_e(a)$.
For every $s \geq t_1$, $u_{e,s}(a) = u_{e,t_1}(a)$ and $\theta_s(a) = \theta_{t_1}(a)$ for each $i \leq a$.
So properties (i) and (ii) will either hold at each stage $s \geq t_1$, or not hold at each stage $s \geq t_1$.
Therefore $Z_{e,a}$ is either finite or cofinite.
\end{proof*}

\begin{claim}
Each requirement $\Rcal_e$ is satisfied. 
\end{claim}
\begin{proof*}
Suppose that $\Phi_e^D$ is total for some $e \in \omega$.
We will prove that $\Phi_e^D$ is not an $f$-diagonalizing function.
Let $A = \set{a \geq e : (\forall e' < e)Z_{e', a} \mbox{ is finite}}$.
Notice that $A$ is cofinite since for each $e' < e$,
there is at most one $a$ such that $Z_{e', a}$ is cofinite.

If for all but finitely many $k \in \omega$, 
we have $k \in \emptyset' \imp k \in \emptyset'_{u_e(k)}$, 
then $\emptyset' \leq_T u_e \leq_T D$, contrary to hypothesis. 
Thus we may let $a$ be the least element of $\{k \in A : k \in \emptyset' \setminus \emptyset'_{u_e(k)}\}$.
We then have
\begin{itemize}
  \item[(1)] $a \geq e$, $\Phi_e^D(a) \downarrow$, $\theta(a) > u_e(a)$
  \item[(2)] For all $e' < e$, there exists $t$ such that $\Rcal_{e'}$
  does not claim $a$ at any stage $s \geq t$.
\end{itemize}
Therefore we may fix $t \geq a$ such that for all $s \geq t$, we have $\Phi^{D_s}_{e,s}(a) \downarrow$,
$\theta_s(a) = \theta(a)$, $u_{e,s}(a) = u_e(a)$, 
and for each $e' < e$, $R_{e'}$ does not claim $a$ at stage $s$. 
Thus, for every $s \geq t$, requirement $\Rcal_e$
claims $a' \leq a$ at stage $s$. Since $Z_{e,i}$ is either finite or cofinite for each $i \leq a$,
it follows that $Z_{e,a}$ is cofinite. By the above argument, we must have
$\Phi^D_e(a) \downarrow$, and by construction, $f(a) = \Phi^D_e(a)$. Therefore $\Rcal_e$ is satisfied.
\end{proof*}

\begin{claim}
The resulting function $f_s$ is $\Delta^0_2$.
\end{claim}
\begin{proof*}
Consider a particular element $a$. By construction, if $e > a$ then $Z_{e,a} = \emptyset$. 
By Claim~\ref{claim:incomplete-srrt22-1}, we have then two cases:
Either $Z_{e, a}$ is finite for all $e \leq a$, in which case for all but finitely many
$s$, $f_s(a) = 0$, or $Z_{e,a}$ is cofinite for some $e$. Then there is a stage $s$ at which
requirement $\Rcal_e$ has committed $f_s(a) = \Phi_e^D(a)$ for assignment
and has never been injured. Thus $f$ is $\Delta^0_2$.
\end{proof*}

This last claim finishes the proof of Theorem~\ref{thm:incomplete-srrt22}.
\end{proof}

\chapter{Avoiding enumerations of closed sets}\label{chap-avoiding-enumerations-closed-sets}

The celebrated first theorem of Liu~\cite{Liu2012RT22} separated Ramsey's theorem for pairs
from weak K\"onig's lemma over computable entailment, therefore showing that the use
of a full compactness argument was not necessary in the proof of~$\rt^2_2$.
Far from closing the topic, Liu's theorem opened two new research directions,
each of them trying to answer a different question:
\smallskip

\begin{itemize}
	\item \emph{What is the precise among of compactness needed in the proof of Ramsey's theorem?}
Flood~\cite{Flood2012Reverse} started the investigation of the question by introducing a Ramsey-type weak K\"onig's lemma,
which happens to be the exact amount of compactness needed in various Ramsey-type statements, such as the Erd\H{o}s-Moser theorem
(see Theorem~\ref{thm:sem-imp-rwkl}).
Further work has been done by Flood and Towsner~\cite{Flood2014Separating} and Bienvenu, Shafer and the author~\cite{Bienvenu2015logical}.
This approach has been developped in chapter~\ref{chap:ramsey-type-konig-lemma}.
\smallskip

	\item \emph{How orthogonal are Ramsey's theorem and compactness arguments?}
Liu's theorem showed that Ramsey's theorem for pairs does not imply the existence of a member
of a particular $\Pi^0_1$ class, namely, the class of the completions of Peano arithmetic.
For which $\Pi^0_1$ classes or more generally for which closed sets of the Baire space is it the case?
Liu~\cite{Liu2015Cone} partially answered this question by proving that Ramsey's theorem for pairs does not imply the existence
of a member of any closed set admitting no computable constant-bound enumeration (defined below).
In this chapter, we push further his investigations by extending this avoidance property to various Ramsey-type
hierarchies such as the free set, the thin set and the rainbow Ramsey theorems.
We also show the optimality of the constant-bound enumeration avoidance by proving
that Ramsey's theorem for pairs does not admit $k$-enumeration avoidance for a fixed~$k$,
and in particular that $\rt^2_2$ implies the existence of members of some special $\Pi^0_1$ classes.
\end{itemize}
\smallskip

All the properties we shall consider in this chapter are special cases of avoidance.
This notion has been defined in chapter~\ref{chap:introduction-reducibilities}, but we now recall the definition.

\begin{definition}[Avoidance]
Let~$\Ccal \subseteq \omega^\omega$ be a set upward-closed under the Turing reducibility.
A $\Pi^1_2$ statement~$\Psf$ admits \emph{$\Ccal$ avoidance} if
for every set~$C \not \in \Ccal$ and every $C$-computable $\Psf$-instance~$X$,
there is a solution~$Y$ to~$X$ such that~$Y \oplus C \not \in \Ccal$.
\end{definition}

The notion of $\Ccal$ avoidance is extended to arbitrary sets by taking their upward-closure
under the Turing reducibility. Moreover, we say that~$\Psf$ admits \emph{strong $\Ccal$ avoidance}
if the~$\Psf$-instance is not required to be $C$-computable. Strong $\Ccal$ avoidance
shows the \emph{structural} weakeness of the statement~$\Psf$, whereas $\Ccal$ avoidance
reflects its \emph{effective} weakness.

Of course, only computably true statements admit $\Ccal$ avoidance for 
an arbitrary set~$\Ccal \subseteq \omega^\omega$. Indeed, let $X$ be a $\Psf$-instance with
no $X$-computable solution and let~$\Ccal$ be the collection of all the solutions to $X$.
If $\Psf$ admits $\Ccal$ avoidance, then by Lemma~\ref{lem:intro-reduc-preservation-included}, 
$\Psf$ holds in a Turing ideal~$\Ical$
containing~$X$ and such that~$\Ical \cap \Ccal = \emptyset$, contradiction.
We must therefore restrict this schema of avoidance to some particular classes of sets.

\section{Avoiding members of a closed set}

In his paper extending the separation of Ramsey's theorem for pairs from weak weak K\"onig's lemma,
Liu~\cite{Liu2015Cone} asked whether whenever an \emph{arbitrary} tree~$T$ has no computable
member, any set $A$ has an infinite subset in either in it or its complement which still does not compute
a path throught~$T$. In this section, we answer negatively and give a general classification
of the theorems in reverse mathematics which admit such a property.

\index{avoidance!of path}
\begin{definition}[Path avoidance]
A $\Pi^1_2$ statement~$\Psf$ admits (strong) \emph{path avoidance}
if it admits (strong) $\Ccal$ avoidance for every closed set~$\Ccal \subseteq \omega^\omega$.
\end{definition}

Unfolding the definition, a $\Pi^1_2$ statements~$\Psf$ admits path avoidance
if for every set~$C$, every closed set~$\Ccal \subseteq \omega^\omega$ with no $C$-computable member,
and every $C$-computable instance~$X$, there is a solution~$Y$ to~$X$ such that $\Ccal$
has no $Y \oplus C$-computable member.
The notion of path avoidance is defined for every closed set of the Baire space.
However, it happens that whenever a principle is shown not to admit path avoidance,
the closed set witnessing the failure belongs to the Cantor space.

\subsection{Cohen genericity}

As usual, a starting point in the analysis of the behavior of the statements of reverse mathematics
with respect to a computability-theoretic property consists of looking at how this property interacts
with typical sets. In this section, we interpret the notion of typicality by Cohen genericity.
Cohen genericity has been introduced in chapter~\ref{chap:introduction-effective-forcing}.

\begin{theorem}\label{thm:genericity-path-avoidance}
Fix a set $C$ computing no member some closed set $\Ccal \subseteq \omega^\omega$.
If $G$ is a real sufficiently Cohen generic, then $G \oplus C$ computes no member of $\Ccal$.
\end{theorem}
\begin{proof}
Given a Turing index $e$, consider the $\Sigma^{0,\Ccal}_2$ sets of strings
$$
D_e = \{\sigma \in 2^{<\omega} : (\exists n)(\forall \tau \succeq \sigma)\Phi_e^{\tau \oplus C}(n) \uparrow\}
$$
$$
H_e = \{\sigma \in 2^{<\omega} : [\Phi_e^{\sigma \oplus C}] \cap \Ccal = \emptyset \}
$$
It suffices to prove that the set $D_e \cup H_e$ is dense.
Let $\sigma \in 2^{<\omega}$. Suppose there exists no finite extension $\tau \in D_e$.
Then for every extension $\tau \succ \sigma$ and every $n \in \omega$, there is an extension~$\rho \succeq \tau$
such that $\Phi^{\rho \oplus C}_e(n) \downarrow$.
Define a $C$-computable sequence of binary strings~$\sigma_0 \prec \sigma_1 \prec \dots$ as follows.
At stage~0, $\sigma_0 = \sigma$. At stage~$s+1$, let~$\sigma_{s+1}$ be the first
string extending~$\sigma_s$ such that~$\Phi^{\sigma_{s+1} \oplus C}_e(s) \downarrow$. 
Such a string exists as~$\sigma_s \succeq \sigma$ and therefore~$\sigma_s \not \in D_e$.
We claim that~$\sigma_s \in H_e$ for some stage~$s \in \omega$. 
If this is not the case, let~$G = \bigcup_s \sigma_s$. 
The real~$G$ is $C$-computable and~$\Phi_e^{G \oplus C}$ is a member of~$\Ccal$, contradiction.
\end{proof}

\begin{corollary}\label{thm:opt-amt-pizog-path-avoidance}
$\opt$, $\amt$ and $\pizog$ admit path avoidance.
\end{corollary}
\begin{proof}
Hirschfeldt et al.~\cite{Hirschfeldt2009atomic} proved that~$\opt$ and $\amt$
are both consequences of $\pizog$, which itself is a restricted notion of Cohen genericity.
\end{proof}

\subsection{The arithmetic hierarchy}

By Simpson's embedding lemma~\cite[Lemma 3.3]{Simpson2007extension} 
(see Corollary~\ref{cor:rt12-no-strong-path-avoidance}), there exists an effectively closed set 
$\Ccal \subseteq 2^\omega$ with no computable member,
and a set $A$ such that every infinite subset in either $A$ or its complement computes a member of $\Ccal$.
Therefore, every degree $\dbf$ such that $A$ is c.e. or co-c.e. relative to $\dbf$ computes a member of $\Ccal$.
However, when considering $\Delta^0_2$ approximations, we never have enough computational power to compute a member
of $\Ccal$, as stated by the following theorem.

\begin{theorem}\label{thm:enum-avoidance-delta2-approx}
Fix a real $C$ computing no member of some closed set $\Ccal \subseteq \omega^\omega$.
For every real $A$, there exists a real $X$ such that $A$ is $\Delta^{0,X}_2$
and $X \oplus C$ computes no member of $\Ccal$.
\end{theorem}
\begin{proof}
For a given real $A$, we build a limit-computable function $f^{\infty} : \omega^2 \to 2$ such that
$f^{\infty} \oplus C$ computes no member of~$\Ccal$ and
$(\forall x)\lim_s f^{\infty}(x,s) = A(x)$.
By Schoenfield's limit lemma, the jump of $f^{\infty}$ computes $A$.
Our forcing conditions are tuples $(g, n)$ such that
$g$ is a finite partial approximation of $f^{\infty}$ and $n$ is an integer.
A condition $(h,m)$ \emph{extends} $(g,n)$ if
\begin{itemize}
	\item[(a)] $dom(g) \subseteq dom(h)$ and $(\forall (x,s) \in dom(g)) g(x,s) = h(x,s)$
	\item[(b)] $m \geq n$ and $(\forall (x,s) \in dom(h) \setminus dom(g))[ s < n \imp h(x,s) = A(s)]$
\end{itemize}
Informally, property (a) states that~$h$ is a function extending $g$
and property (b) forces the $n$ first columns of $g$ to converge to $A \restr n$.
Therefore making $n$ grow arbitrarily large will ensure that the constructed function~$f^{\infty}$ is a $\Delta^0_2$
approximation of $A$. Our initial condition is $(\emptyset, 0)$.
The following lemma states that we can force the constructed function~$f^{\infty}$ to be total.

\begin{lemma}\label{lem:enum-avoidance-delta2-ext}
For every condition $(g,n)$ and every $t \in \omega$, there exists an extension $(h,m)$
such that $m > n$ and $[0, t]^2 \subseteq dom(g)$
\end{lemma}
\begin{proof*}
Let $h$ be the function over domain $[0,t]^2 \cup dom(g)$
defined by $h(x,s) = g(x,s)$ for $(x,s) \in dom(g)$ and 
$h(x,s) = A(x)$ for $(x,s) \not \in dom(g)$. 
$(h, n+1)$ is a  valid extension of $(g,n)$.
\end{proof*}

A function~$f^\infty : \omega^2 \to 2$ \emph{satisfies} a condition~$(g, n)$
if~$(\forall (x,s) \in dom(g)) g(x,s) = f^\infty(x,s)$
and $(\forall (x,s) \in \omega^2 \setminus dom(g))[s < n \imp f^\infty(x,s) = A(s)]$.
In other words, for every finite approximation $h$ of~$f^\infty$ such that~$dom(h) \geq dom(g)$,
$(h, n)$ is a valid extension of~$(g, n)$. Note that~$f^{\infty}$ may not be limit-computable,
and that if~$f^{\infty}$ satisfies~$(g,n)$ and~$m < n$, then~$f^{\infty}$ satisfies~$(g, m)$.
A condition~$c$ \emph{forces~$\Phi^{f^{\infty} \oplus C}_e$ to be partial}
if $\Phi^{f^{\infty} \oplus C}_e$ is partial for every function~$f^{\infty}$
satisfying~$c$.

\begin{lemma}\label{lem:enum-avoidance-delta2-force}
For every condition $(g,n)$ and every $e \in \omega$, there exists an extension $(h,m)$
forcing $\Phi^{f^{\infty} \oplus C}_e$ to be partial, or $[\Phi^{h \oplus C}_e] \cap \Ccal = \emptyset$.
\end{lemma}
\begin{proof*}
If there is an extension~$(h,m)$ forcing $\Phi^{f^{\infty} \oplus C}_e$ to be partial
or such that $\Phi_e^{h \oplus C} \restr n = \sigma$ for some string~$\sigma \in \omega^n$
such that~$\Ccal \cap [\sigma] = \emptyset$, then we are done.
So suppose that it is not the case. We will describe how to $C$-compute a member of~$\Ccal$
and derive a contradiction. Define a $C$-computable sequence
of conditions~$(g, n) = (g_0, n) \geq (g_1, n) \geq \dots$ as follows:
Given some condition~$(g_i, n)$, let~$(g_{i+1}, n)$ be the least extension such that
$\Phi_e^{g_{i+1} \oplus C}(i) \downarrow$. Such an extension exists
as otherwise~$(g_i, n)$ would force~$\Phi^{f^{\infty} \oplus C}_e$ to be partial.
Let~$f^{\infty} = \bigcup_i g_i$. The function~$f^{\infty}$ has been constructed
$C$-computably in such a way that~$\Phi^{f^{\infty} \oplus C}_e$ is total and a member of~$\Ccal$.
This contradicts the assumption that~$C$ does not compute a member of~$\Ccal$.
\end{proof*}

Let~$\Fcal = \{c_0, c_1, \dots\}$ be a sufficiently generic filter containing $(\emptyset, 0)$,
where~$c_s = (g_s, n_s)$. The filter~$\Fcal$ yields a unique partial function~$f^\infty = \bigcup_s g_s$.
By Lemma~\ref{lem:enum-avoidance-delta2-ext}, the function~$f^\infty$ is total,
and by definition of a forcing condition, $f^\infty$ is a $\Delta^0_2$ approximation of the real~$A$.
By Lemma~\ref{lem:enum-avoidance-delta2-force}, $f^{\infty} \oplus C$ computes no member of~$\Ccal$.
\end{proof}

\begin{corollary}\label{cor:coh-path-avoidance}
$\coh$ admits path avoidance.
\end{corollary}
\begin{proof}
Fix a real~$C$ computing no member of some closed set $\Ccal \subseteq \omega^\omega$
and let $R_0, R_1, \dots$ be a uniformly $C$-computable sequence of reals.
By Theorem~\ref{thm:enum-avoidance-delta2-approx}, there exists a real~$X$ such that~$X \oplus C$ computes no member of~$\Ccal$
and the jump of~$X$ computes~$\emptyset''$.
Jockusch and Stephan~\cite{Jockusch1993cohesive} proved that if~$R_0, R_1, \dots$
is a uniform sequence of reals, for any real~$X$ whose jump if of PA degree relative to the jump of~$\vec{R}$,
$X \oplus \vec{R}$ computes an infinite $\vec{R}$-cohesive real.
Therefore~$X \oplus C$ computes an infinite~$\vec{R}$-cohesive real.
\end{proof}

\begin{corollary}
For every real $A$ and every non-computable real $B$,
there exists a real $X$ such that $A \in \Delta^{0,X}_2$
but $X \not \geq_T B$.
\end{corollary}
\begin{proof}
Apply Theorem~\ref{thm:enum-avoidance-delta2-approx} with~$\Ccal = \{B\}$
to obtain a real $X$ such that $A \in \Delta^{0,X}_2$
and $X$ computes no member of $\Ccal$, hence $X \not \geq_T B$.
\end{proof}

We shall see in Corollary~\ref{cor:lots-of-statements-no-strong-path-avoidance} that 
$\coh$ does not admit strong path avoidance since~$\rt^1_2$ does not.

\subsection{The embedding lemma}

The following application of Simpson's embedding lemma is very useful
for proving that some principle does not admit path avoidance.

\begin{lemma}\label{lem:simpson-embedding-path}
If some principle~$\Psf$ has a computable (resp.\ arbitrary) instance with no computable solution
and such that its collection of solutions is a~$\Sigma^0_3$ subset of $\omega^\omega$,
then $\Psf$ does not admit (strong) path avoidance.
\end{lemma}
\begin{proof}
We prove it in the case of path avoidance.
Let~$X$ be a computable~$\Psf$-instance with no computable solution,
and let~$\Ccal \subseteq \omega^\omega$ be its set of solutions.
By Lemma~3.3 in Simpson~\cite{Simpson2007extension}, there exists an effectively closed set~$\Dcal \subseteq 2^\omega$
whose degrees are exactly the PA degrees and the degrees of members of~$\Ccal$.
Since~$X$ has no computable solution, $\Dcal$ has no computable member. Every solution to~$X$
is a member of~$\Ccal$ and thus computes a member of~$\Dcal$. Therefore~$\Psf$ does not admit path avoidance.
\end{proof}

Note that the witness of failure of path avoidance is an effectively closed set.
The three following corollaries are direct applications of Lemma~\ref{lem:simpson-embedding-path}.

\begin{corollary}\label{cor:rt12-no-strong-path-avoidance}
$\rt^1_2$ does not admit strong path avoidance.
\end{corollary}
\begin{proof}
Let~$A$ be a $\Delta^0_2$ bi-immune set.
The collection of its infinite homogeneous sets a~$\Pi^0_2$ subset of~$\omega^\omega$:
\[
\Ccal = \{ X \in \omega^\omega :
	(\forall i)[X(i) <_{\N} X(i+1) \wedge X(i) \in A \biimp X(i+1) \in A] \}
\]
Apply Lemma~\ref{lem:simpson-embedding-path}.
\end{proof}

Of course, if~$\Qsf \leq_c \Psf$ and~$\Qsf$ does not admit path avoidance,
then so does~$\Psf$. We therefore want to prove that very weak principles
do not admit path avoidance to obtain the same conclusion for many statements belonging to the reverse mathematics zoo.

\begin{corollary}\label{cor:dnr-path-avoidance}
$\dnr$ does not admit path avoidance.
\end{corollary}
\begin{proof}
The collection of d.n.c.\ functions is a $\Pi^0_1$ subset of~$\omega^\omega$
with no computable member:
\[
\Ccal = \{ f \in \omega^\omega : (\forall e,s)[\Phi_{e,s}(e) \downarrow \imp \Phi_{e,s}(e) \neq f(e)] \}
\]
Apply Lemma~\ref{lem:simpson-embedding-path}.
\end{proof}

\begin{corollary}\label{cor:sads-path-avoidance}
$\sads$ does not admit path avoidance.
\end{corollary}
\begin{proof}
Tennenbaum~\cite{Rosenstein1982Linear} constructed a computable linear order of order type $\omega + \omega^{*}$
with no computable ascending or descending sequence.
Given a linear order~$\Lcal$, the collection of its infinite ascending or descending sequences is
a~$\Pi^0_1$ subset of~$\omega^\omega$:
\[
\Ccal = \{ X \in \omega^\omega :
	(\forall i)[X(i) <_\Lcal X(i+1)] \vee (\forall i)[X_i >_\Lcal X(i+1)] \}
\]
Apply Lemma~\ref{lem:simpson-embedding-path}.
\end{proof}

The following lemma shows that avoidance is closed
downward under computable reducibilty. As many proofs of reductions
in reverse mathematics are in fact computable reductions, this
lemma has many applications.

\begin{lemma}\label{lem:computably-reducible-avoidance}
\label{lem:strong-computably-reducible-strong-avoidance}
If $\Psf$ is (strongly) computably reducible to~$\Qsf$ and 
$\Qsf$ admits (strong) $\Ccal$ avoidance, then so does $\Psf$.
\end{lemma}
\begin{proof}
We prove it in the case of computable reducibility. The strong case is similar.
Let~$C$ be a real computing no member of~$\Ccal$
and let $I$ be a $C$-computable instance of $P$.
As $P \leq_c Q$, there exists an $I$-computable instance $J$ of $Q$ such that
for every solution $X$ to $J$, $X \oplus I$ computes a solution to $I$.
By $\Ccal$ avoidance of $Q$, there exists a solution $X$ to $J$
such that $X \oplus C$ computes no member of $\Ccal$.
$X \oplus C$ computes a solution $Y$ to $I$, but computes no member of~$\Ccal$.
\end{proof}

\begin{corollary}\label{cor:lots-of-statements-no-path-avoidance}
None of~$\rt^2_2$ $\ads$, $\cac$, $\emo$, $\ts^2$ $\rrt^2_2$ admit path avoidance.
\end{corollary}
\begin{proof}
By Hirschfeldt et al.~\cite{Hirschfeldt2008strength}, $\dnr \leq_c \srt^2_2$.
By Hirschfeldt \& Shore~\cite{Hirschfeldt2007Combinatorial}, $\sads \leq_c \ads \leq_c \cac$.
By Rice~\cite{RiceThin}, $\dnr \leq_c \ts^2$.
By Miller~\cite{MillerAssorted}, $\dnr \leq_c \rrt^2_2$.
By the author~\cite{Patey2015Somewhere}, $\dnr \leq_c \emo$.
Conclude by Lemma~\ref{lem:computably-reducible-avoidance},
Corollary~\ref{cor:dnr-path-avoidance} and Corollary~\ref{cor:sads-path-avoidance}.
\end{proof}

\begin{corollary}\label{cor:lots-of-statements-no-strong-path-avoidance}
$\coh$ does not admit strong path avoidance.
\end{corollary}
\begin{proof}
Immediate by Corollary~\ref{cor:rt12-no-strong-path-avoidance}, 
Lemma~\ref{lem:strong-computably-reducible-strong-avoidance} and the fact that~$\rt^1_2 \leq_{sc} \coh$.
\end{proof}

\subsection{Simultaneous path avoidance}

The notion of path avoidance expresses the ability for a principle to avoid
computing a member of a (boldface) $\Pi^0_1$ set of the Baire space. We now
see that the notion of avoidance for $F_\sigma$ sets coincides with path avoidance. 

\index{avoidance!of path (strong)}
\index{avoidance!of path (simultaneous)}
\begin{definition}[Simultaneous path avoidance]
Fix a countable collection of closed sets $\Ccal_0, \Ccal_1, \dots \subseteq \omega^{\omega}$.
A principle $\Psf$ admits \emph{(strong) path avoidance} for $\vec{\Ccal}$
if it admits (strong) $\bigcup_i \Ccal_i$ avoidance.
A principle~$\Psf$ admits \emph{(strong) simultaneous path avoidance}
if it admits (strong) path avoidance for~$\vec{\Ccal}$ for every
countable collection of closed sets~$\Ccal_0, \Ccal_1, \dots \subseteq \omega^{\omega}$
\end{definition}

The notion of Muchnik reducibility is very useful to substitute closed sets by some others
while preserving the notion of avoidance.

\index{reduction!Muchnik}
\index{Muchnik reducibility}
\begin{definition}[Muchnik reducibility]
Let $\Ccal$ and $\Dcal$ be two classes of reals.
$\Ccal$ is \emph{Muchnik reducible} to $\Dcal$ (denoted by $\Ccal \leq_w \Dcal$)
if for every $X \in \Dcal$, there exists a $Y \in \Ccal$ such that $Y \leq_T X$.
\end{definition}

\begin{lemma}\label{lem:path-avoidance-simultaneous}
Let $\Ccal_0, \Ccal_1, \dots \subseteq \omega^\omega$ be a countable collection 
of closed sets such that~$\Ccal_i$ has no computable member for each $i$.
There exists a closed set $\Dcal \subseteq \omega^\omega$ such that
$\Dcal$ and~$\bigcup_i \Ccal_i$ are Muchnik equivalent. Moreover, if the $\Ccal$'s belong
the Cantor space, then so does~$\Dcal$.
\end{lemma}
\begin{proof}
We may assume that some~$\Ccal_i$ is non-empty as otherwise,
$\Dcal = \emptyset$ is a trivial solution.
Let~$X$ be a member of some~$\Ccal_i$ and define~$\Dcal$ as follows:
$$
\Dcal = \{ \sigma^\concat(i+1 \mod 2)^\concat Z : \sigma^\concat i \prec X \wedge Z \in \Ccal_{|\sigma|}\}
$$
The set~$\Dcal$ is closed and Muchnik equivalent to~$\bigcup_i \Ccal_i$.
\end{proof}

\begin{corollary}
If a principle~$\Psf$ admits (strong) path avoidance,
then it admits (strong) simultaneous path avoidance.
\end{corollary}
\begin{proof}
We prove it in the case of path avoidance. 
Let~$C$ be a set computing no member of~$\bigcup_i \Ccal_i$
for some countable collection of sets~$\Ccal_0, \Ccal_1, \dots \subseteq \omega^\omega$,
and let~$X$ be a $C$-computable instance of~$\Psf$.
By Lemma~\ref{lem:path-avoidance-simultaneous}, there exists a closed set of reals~$\Dcal$
Muchnik equivalent to~$\bigcup_i \Ccal_i$. Therefore $C$ computes no member of~$\Dcal$.
By path avoidance of~$\Psf$, there is a solution $Y$ to $X$ such that~$Y \oplus C$
computes no member of~$\Dcal$ and therefore computes no member of~$\bigcup_i \Ccal_i$.
\end{proof}

\section{Enumeration avoidance}

We have seen in the previous section that Ramsey's theorem for pairs
does not admit path avoidance even for special $\Pi^0_1$ classes
(Corollary~\ref{cor:lots-of-statements-no-path-avoidance}).
In his follow-up paper~\cite{Liu2015Cone}, Liu identified a stronger property enjoyed by the $\Pi^0_1$ class
of completions of Peano arithmetic and for which Ramsey's theorem for pairs admits avoidance.

\index{k-enumeration@$k$-enumeration}
\index{k-enum@$k$-enum|see {$k$-enumeration}}
\index{constant-bound enumeration}
\index{avoidance!of c.b-enum}
\index{avoidance!of c.b-enum (strong)}
\begin{definition}[Constant-bound enumeration]\ 
\begin{itemize}
	\item[(i)]
A \emph{$k$-enumeration} (or $k$-enum) of a set $\Dcal \subseteq \omega^{\omega}$ is a sequence $D_0, D_1, \dots$
of finite sets of strings such that for each $n \in \omega$, $|D_n| \leq k$, $(\forall \sigma \in D_n)|\sigma| = n$
and $\Dcal \cap [D_n] \neq \emptyset$.
A \emph{constant-bound enumeration} (or c.b-enum) of $\Dcal$ is a $k$-enum of $\Dcal$ for some $k \in \omega$.
	\item[(ii)] Fix a collection of sets~$\Ccal_0, \Ccal_1, \dots \subseteq \omega^\omega$.
	A $\Pi^1_2$ statement admits \emph{(strong) c.b-enum avoidance} for $\vec{\Ccal}$
if it admits (strong) $\Dcal$ avoidance, where $\Dcal$ is the set of reals
wich code a c.b-enum of some~$\Ccal_i$.
\end{itemize}
\end{definition}

We can define the notion of (strong) $k$-avoidance accordingly for every~$k$.
Note that the notion of $k$-avoidance for every~$k$ and the notion of c.b-enum avoidance
do not necessarily coincide. In fact, we shall see that c.b-enum avoidance is strictly weaker.
The following trivial lemma shows that the existence of a computable c.b-enum for an effectively closed
set of reals is not invariant in the Muchnik degrees.

\begin{lemma}\label{lem:special-effectively-closed-enum}
For every effectively closed set $\Ccal \subseteq 2^\omega$, 
there exists an effectively closed set $\Dcal \subseteq 2^\omega$
Muchnik equivalent to $\Ccal$ with a computable 1-enum.
\end{lemma}
\begin{proof}
Let $T$ be a computable tree such that $[T] = \Ccal$.
The set $\Dcal = \{\sigma^\concat Z : \sigma \in T \wedge Z \in \Ccal\}$
is effectively closed and Muchnik equivalent to $\Ccal$.
For every $\sigma \in T$, $[\sigma] \cap \Dcal \neq \emptyset$, therefore
we can compute a 1-enum of $\Dcal$ by returning on input $n$
a string of length $n$ in $T$.
\end{proof}

However, when considering the uniform version of Muchnik reducibility, namely,
Medvedev reducibility, then the degrees of the c.b-enums are preserved in the special
case of compact sets.

\index{reduction!Medvedev}
\index{Medvedev reducibility}
\begin{definition}[Medvedev reducibility]
Let $\Ccal \subseteq \omega^\omega$ and $\Dcal \subseteq \omega^\omega$ be two sets of sequences.
We say that $\Ccal$ is \emph{Medvedev reducible} to $\Dcal$ (denoted by $\Ccal \leq_s \Dcal$)
if there exists a Turing functional $\Gamma$ such that $\Gamma^X \in \Ccal$ for every sequence~$X \in \Dcal$.
\end{definition}

\begin{lemma}\label{lem:cbe-trough-medvedev}
Let $\Ccal \subseteq \omega^\omega$ be a set of sequences Medvedev below a compact set of sequences $\Dcal \subseteq \omega^\omega$.
For every $k \in \omega$, every $k$-enum of $\Dcal$ computes a $k$-enum of~$\Ccal$.
\end{lemma}
\begin{proof}
Let $\Gamma$ be the Turing functional witnessing  the Medvedev reduction from $\Ccal$ to $\Dcal$.
We prove it by induction over $k$.
Let $(D_i : i \in \omega)$ be a $k$-enum of $\Dcal$.
Suppose that there exists a $\sigma \in 2^{<\omega}$
such that $\Dcal \cap [\sigma] = \emptyset$
and for infinitely many $i \in \omega$,
$\sigma \preceq \tau$ for some $\tau \in D_i$.
Then $k > 1$ and we can compute a $(k-1)$-enum of $\Dcal$
by computably finding on input $i$ a $j > i$
such that $\sigma \preceq \tau$ for some $\tau \in D_j$
and returning $E_i = \{\sigma \restr i : \sigma \in D_j \setminus \tau\}$.
$\vec{E}$ is a $(k-1)$-enum of $\Dcal$ and by induction hypothesis,
it computes a $(k-1)$-enum of $\Ccal$, so \emph{a fortiori} a $k$-enum of $\Ccal$.

So suppose there exists no such $\sigma$.
This means that for every $i \in \omega$,
there exists a $j > i$ such that $\Dcal \cap [\sigma \restr i] \neq \emptyset$
for each $\sigma \in D_j$.
As $\Ccal \subseteq h^\omega$, by the pigeonhole principle $\Gamma$ 
will produce arbitrarily large $k$-tuples of initial segments
of members of~$\Ccal$. We compute a $k$-enum of $\Ccal$ as follows:
For each $i \in \omega$, let $E_i = \{\Gamma^{\sigma} \restr i : \sigma \in D_j\}$
for some $j$ such that $\Gamma^{\sigma} \restr i$ is defined on each $\sigma \in D_j$.
Such $E_i$ has been shown to exist and can be found computably in $\vec{D}$.
As $[\sigma] \cap \Dcal \neq \emptyset$ for some $\sigma \in D_j$, 
$[\Gamma^{\sigma} \restr i] \cap \Ccal \neq \emptyset$, 
hence $(\exists \tau \in E_i)\Ccal \cap [\tau] \neq \emptyset$
hence and $\vec{E}$ is a valid $k$-enum of~$\Ccal$.
\end{proof}

Like we did for path avoidance, we can define the notion of simultaneous c.b-enum avoidance.

\index{avoidance!of c.b-enum (simultaneous)}
\begin{definition}[Simultaneous c.b-enum avoidance] Fix a $\Pi^1_2$ statement~$\Psf$.
	$\Psf$ admits \emph{(strong) simultaneous c.b-enum avoidance}
	if it admits (strong) c.b-enum avoidance for every countable 
	collection of sets~$\Ccal_0, \Ccal_1 \dots \subseteq 2^\omega$.
\end{definition}

First, notice that unlike path avoidance, we did not require the sets to be closed.
Indeed, a set and its topological closure have the same constant-bound enumerations.
Also notice that we defined the notion of c.b-enum avoidance over the Cantor space.
In fact, this restriction is expressive enough to obtain c.b-enum avoidance over compact sets
of the Baire space, since for every compact set~$\Ccal \subseteq \omega^\omega$,
one can find a closed set~$\Dcal \subseteq 2^\omega$ such that the degrees of the c.b-enums
of~$\Ccal$ and of~$\Dcal$ coincide~\cite{Patey2015Combinatorial}.

We now relate c.b-enum avoidance and simultaneous c.b-enum avoidance
as we did for path avoidance. See~\cite{Patey2015Combinatorial} for an extensive
study of the relations between path avoidance, $k$-avoidance, c.b-enum avoidance
and their simultaneous variants. We start by proving that the notions of c.b-enum avoidance
and simultaneous avoidance do not coincide. Moreover, there is a whole hierarchy of avoidance relations
based on how many closed sets can be avoided simultaneously.

\begin{theorem}\label{thm:not-cbe-avoidance-simultaneous-arbitrary}
There exists a countable collection of closed sets $\Ccal_0, \Ccal_1, \dots \subseteq 2^\omega$ together
with a $\Delta^0_2$ function $f : \omega \to \omega$ and a 1-enum $(\rho_i : i \in \omega)$ such that
\begin{itemize}
	\item[(i)] $\bigcup_{j \neq i} \Ccal_j$ has no computable c.b-enum for each~$i$
	\item[(ii)] $[\rho_i] \cap \Ccal_{f(i)} \neq \emptyset$ for each $i$
\end{itemize}
\end{theorem}
\begin{proof}
Fix a non-computable $\Delta^0_2$ set $X$ and a computable sequence $X_0, X_1, \dots$
of reals pointwise converging to~$X$. We build the closed sets of reals~$\vec{\Ccal}$ by forcing.
Our forcing conditions are tuples $(k, \Ccal_0, \dots, \Ccal_{k-1}, E_0, \dots,  E_{k-1})$ where
\begin{itemize}
	\item[(a)] $\bigcup_{j \neq i} \Ccal_j$ are closed sets containing $X$ and with no computable c.b-enum for each~$i < k$
	\item[(b)] $E_i$ is a finite set of strings for each~$i < k$
	\item[(c)] $(\bigcup_{j \neq i} \Ccal_j) \cap [E_i] = \emptyset$ for each~$i < k$
	\item[(d)] $(\forall s)(\exists i < k)([X_s \restr s] \not \subseteq [\bigcup_{j \neq i} E_j])$

\end{itemize}
A condition $(m, \tilde{\Ccal}_0, \dots, \tilde{\Ccal}_{m-1}, \tilde{E}_0, \dots, \tilde{E}_{m-1})$ \emph{extends}
a condition $(k, \Ccal_0, \dots, \Ccal_{k-1}, E_0, \dots, E_{k-1})$ if $m \geq k$,
$\Ccal_i \subseteq \tilde{\Ccal}_i$ and $E_i \subseteq \tilde{E}_i$ for each~$i < k$.
The set~$E_i$ is a forbidden open set for~$\bigcup_{j \neq i} \Ccal_j$. To force $\bigcup_{j \neq i} \Ccal_j$
not to have computable c.b-enum, we shall put strings in it.
Our initial condition is $(2, \{X\}, \{X\}, \emptyset, \emptyset)$ which is valid
since every c.b-enum of a singleton~$\{X\}$ computes~$X$.
Note that given some condition $c = (k, \Ccal_0, \dots, \Ccal_{k-1}, E_0, \dots, E_{k-1})$,
the condition~$(k+1, \Ccal_0, \dots, \Ccal_{k-1}, \{X\}, E_0, \dots, E_{k-1}, \emptyset)$
is a valid extension of~$c$.

We want our forcing to be $\emptyset'$-effective
to obtain a $\Delta^0_2$ function $f : \omega \to \omega$
such that property (ii) holds. Given some condition~$c = (k, \Ccal_0, \dots, \Ccal_{k-1}, E_0, \dots, E_{k-1})$,
a \emph{code of $c$} is a tuple $\tuple{k, e_0, \dots, e_{k-1}, E_0, \dots, E_{k-1}}$
such that for each~$i < k$, $\Phi^{\emptyset'}_{e_i}$ is the characteristic function
of the set of strings~$\sigma \in 2^{<\omega}$ such that~$[\sigma] \cap \Ccal_i \neq \emptyset$.
Note that a condition may not have a code in general,
but our initial condition $(2, \{X\}, \{X\}, \emptyset, \emptyset)$ has one.
We will show that we can $\emptyset'$-effectively find an infinite decreasing sequence of extensions
having codes and forcing the desired properties.

\begin{lemma}\label{lem:not-cbe-avoidance-simultaneous-arbitrary-enum}
For every condition $c = (k, \Ccal_0, \dots, \Ccal_{k-1}, E_0, \dots, E_{k-1})$ and $s \in \omega$, there exists an extension
$d = (k, \tilde{\Ccal}_0, \dots, \tilde{\Ccal}_{k-1}, E_0, \dots, E_{k-1})$  
and some~$i < k$ such that $[X_s \restr s] \cap \tilde{\Ccal_i} \neq \emptyset$.
Moreover, one can $\emptyset'$-effectively find a code of~$d$ given a code of~$c$.
\end{lemma}
\begin{proof*}
By property (d) of the condition~$c$ there is some~$i < k$ such that 
$[X_s \restr s] \not \subseteq [\bigcup_{j \neq i} E_j]$. Let~$E = \bigcup_{j \neq i} E_j$.
As $E$ is finite, there exists a finite $\tau \succ X_s \restr s$ such that
$[\tau] \cap [E] = \emptyset$. Moreover, those~$i$ and~$\tau$ can be $\emptyset'$-effectively found.
Let $\tilde{\Ccal}_i = \Ccal_i \cup \{\tau^\concat Z : Z \in \Ccal_i\}$
and let $\tilde{\Ccal}_j = \Ccal_j$ for each~$j \neq i$.
The closed set~$\tilde{\Ccal}_i$ is Medvedev above~$\Ccal_i$. Therefore, for each~$j < k$, 
$\bigcup_{r \neq j} \tilde{\Ccal}_r$ is Medvedev above $\bigcup_{r \neq j} \Ccal_r$
and by Lemma~\ref{lem:cbe-trough-medvedev} and property (a) of condition~$c$, it admits no computable c.b-enum.
The condition $d = (k, \tilde{\Ccal}_0, \dots, \tilde{\Ccal}_{k-1}, E_0, \dots, E_{k-1})$
satisfies therefore  properties (a), (b) and (d). We check property~(c).
If $(\bigcup_{r \neq j} \tilde{\Ccal}_r) \cap [E_j] \neq \emptyset$ for some~$j < k$, 
then by property (c) of the condition~$c$, $j \neq i$ and $(\bigcup_{r \neq j} \Ccal_r) \cap [E_j] = \emptyset$.
As $(\bigcup_{r \neq j} \tilde{\Ccal}_r) \subseteq (\bigcup_{r \neq j} \Ccal_r) \cup [\tau]$,
we obtain $[\tau] \cap [E_j] \neq \emptyset$, contradiction.
Hence property (c) holds and $d$ is a valid extension of $c$.
The Turing index of the characteristic function of the strings extensible in~$\tilde{\Ccal}_i$
can be effectively found from the Turing index of the characteristic function of the strings extensible in~$\Ccal_i$.
Therefore the condition~$d$ has a code, which can be~$\emptyset'$-effectively found from a code of~$c$.
\end{proof*}

\begin{lemma}\label{lem:not-cbe-avoidance-simultaneous-arbitrary-force}
For every condition $c = (k, \Ccal_0, \dots, \Ccal_{k-1}, E_0, \dots, E_{k-1})$,
every~$i < k$ and every $e \in \omega$,
there exists an extension $d = (k, \Ccal_0, \dots, \Ccal_{k-1}, \tilde{E}_0, \dots, \tilde{E}_{k-1})$ such that
if  $\Phi_e$ is an $e$-enum then $(\exists n)\Phi_e(n) \subset \tilde{E}_i$.
Moreover, one can $\emptyset'$-effectively find a code of~$d$ given a code of~$c$.
\end{lemma}
\begin{proof*}
Let~$F = \bigcup_{j \neq i} E_j$ and let $u = max(|\sigma| : \sigma \in F)$.
We can $\emptyset'$-effectively find some stage $t > u$ such that $X_t \restr u = X \restr u$.
By Lemma~\ref{lem:not-cbe-avoidance-simultaneous-arbitrary-enum}, we can assume
that for every $s < t$, there is some~$j < k$ such that $[X_s \restr s] \cap \Ccal_j \neq \emptyset$.
As by property (a) of the condition~$c$, $\bigcup_{j \neq i} \Ccal_j$ admits no computable c.b-enum, 
there exists some $n > t + e$ such that either~$\Phi_e(n) \uparrow$, or $[\Phi_e(n)] \cap \bigcup_{j \neq i} \Ccal_j = \emptyset$. 
We can~$\emptyset'$-decide in which case we are. In the first case, we take~$c$ as the desired extension.
Set $\tilde{E}_i = E_i \cup \Phi_e(n)$ and~$\tilde{E}_j = E_j$ for each~$j \neq i$.
Properties (a), (b) and (c) hold for the condition 
$d = (k, \Ccal_0, \dots, \Ccal_{k-1}, \tilde{E}_0, \dots, \tilde{E}_{k-1})$.
We now check property~(d).

Suppose for the sake of contradiction that for some $s$, for every~$j < k$, $[X_s \restr s] \subseteq [\bigcup_{r \neq j} \tilde{E}_r]$.
In particular, $[X_s \restr s] \subseteq [F]$.
In this case $s < t$, otherwise $[X_s \restr s] \subseteq [X \restr u]$.
But then $[X \restr u] \cap [F] \neq \emptyset$
and as $u = max(|\sigma| : \sigma \in F)$, $[X \restr u] \subseteq [E_j]$ for some~$j < k$, contradicting
the fact that~$X \in \bigcup_{r \neq j} \Ccal_r$ and property (c) of the condition~$c$.
By property (d) of the condition~$c$,
there exists some~$j < k$ such that $[X_s \restr s] \not \subseteq [\bigcup_{r \neq j} E_r]$.
Let~$\mu$ be the Lebesgue measure. 
Since~$t > max(|\sigma| : \sigma \in F), t > s$ and $[X_s \restr s] \not \subseteq [\bigcup_{r \neq j} E_r]$,
$\mu([X_s \restr s] \setminus [\bigcup_{r \neq j} E_r]) \geq 2^{-t}$.
Since~$\Phi_e$ is an e-enum and~$n > t + e$, $\mu([\Phi_e(n)]) \leq e \times 2^{-t-e} < 2^{-t}$.
Therefore, $\mu([X_s \restr s] \setminus ([\bigcup_{r \neq j} E_r] \cup [\Phi_e(e)])) > 0$
so $[X_s \restr s] \not \subseteq [\bigcup_{r \neq j} \tilde{E}_r]$, contradiction.
\end{proof*}

Thanks to Lemma~\ref{lem:not-cbe-avoidance-simultaneous-arbitrary-enum} 
and Lemma~\ref{lem:not-cbe-avoidance-simultaneous-arbitrary-force}, 
we build an infinite $\emptyset'$-computable decreasing sequence of conditions
$c_0 = (\{X\}, \{X\}, \emptyset, \emptyset) \geq c_1 \geq c_2 \geq \dots$ together with their codes,
such that for each $s \in \omega$, assuming $c_s = (k_s, \Ccal_{0,s}, \dots, \Ccal_{k_s-1,s}, E_{0,s}, \dots, E_{k_s-1,s})$,
\begin{itemize}
	\item[(i)] $k_s \geq s$
	\item[(ii)] If $\Phi_s$ is a total $s$-enum, then $(\forall i < k_s)(\exists n)\Phi_s(n) \subset E_{i,s}$
	\item[(iii)] $[X_s \restr s] \cap \bigcup_{i < k_s} \Ccal_{i,s} \neq \emptyset$
\end{itemize}
This way, taking $\Ccal_i = \bigcup_{s \geq i} \Ccal_{i,s}$,
we obtain two closed sets admitting no computable c.b-enum by (ii)
and such that $s \mapsto X_s \restr s$ is a computable 1-enum of $\bigcup_i \Ccal_i$ by (ii).
This completes the proof of Theorem~\ref{thm:not-cbe-avoidance-simultaneous-arbitrary}.
\end{proof}

In the following corollary, $n$ c.b-enum avoidance is the restriction
of simultaneous avoidance to sequences of $n$ sets of reals.

\begin{corollary}\label{cor:sts2-not-simultaneous-cb-enum-avoidance}
$\sts^2_n$ does not admit $n$ c.b-enum avoidance for every~$n \geq 2$.
In particular, $\srt^2_2$ does not admit 2 c.b-enum avoidance.
Moreover, $\sts^2$ does not admit simultaneous c.b-enum avoidance.
\end{corollary}

Liu defined in~\cite{Liu2015Cone} c.b-enum avoidance for 
any increasing sequence (in inclusion order) of sets of reals.
We prove that this apparently stronger notion of avoidance is 
in fact reducible to c.b-enum avoidance.

\begin{lemma}\label{lem:cbe-avoidance-increasing-inherit}
Let $\Ccal_0 \subseteq \Ccal_1 \subseteq \dots \subseteq 2^\omega$ be an increasing countable collection of sets of reals
with no computable c.b-enum.
There exists a set $\Dcal \subseteq 2^\omega$ Medvedev below each $\Ccal_i$
such that $\Dcal$ has no computable c.b-enum.
\end{lemma}
\begin{proof}
Fix a set $X \not \leq_T \emptyset'$ and let
$$
\Dcal = \{ \sigma^\concat(1-i)^\concat Z : \sigma^\concat i \prec X \wedge Z \in \Ccal_{|\sigma|}\}
$$
The set $\Dcal$ is Medvedev below each $\Ccal_i$.
We prove that there exists no computable c.b-enum of $\Dcal$.
Fix a computable $k$-enum $(D_i : i \in \omega)$ of $\Dcal$.
By thinning out $\vec{D}$, we can obtain a computable $k$-enum $(E_i : i \in \omega)$ of $\Dcal$
together with a finite set of strings (with possible duplications) $\sigma_0, \dots, \sigma_{r-1}$ for some $r \leq k$
and a computable injective function~$g : \omega \times r \to 2^{<\omega}$ such that
\begin{itemize}
	\item[(i)] $(\forall i < r)\sigma_i \not \prec X \wedge (\sigma_i \uh |\sigma_i|-1) \prec X$
	\item[(ii)] $(\forall i \in \omega)(\forall j < r)[g(i, j) \in E_i \wedge \sigma_j \prec g(i,j)]$
	\item[(iii)] if $\sigma \not \prec X$ then there
	are finitely many $i$ such that $\sigma \prec \tau$ for some $\tau \in E_i \setminus \{g(i,j) : j < r\}$. 
\end{itemize}
If $r = k$ then let $n = max(\{|\sigma_j| : j < r\})$.
For each~$i \in \omega$ and~$j < k$, let~$f(i,j)$ be the unique string~$\rho$ of length~$i$
such that~$\sigma_j^\concat \rho \prec g(n+i,j)$ and let~$F_i = \{g(i,j) : j < k\}$. 
We claim that the sequence~$\vec{F}$ is a $k$-enum of~$\Ccal_n$.
Indeed, since~$\vec{E}$ is a $k$-enum of~$\Dcal$, for each~$i$, there exists a $\tau \in E_{i+n}$
such that~$[\tau] \cap \Dcal \neq \emptyset$. Since~$f$ is injective, there is some~$j < k$
such that~$\tau = g(i+n, j)$.
By construction of $\Dcal$, $[f(i, j)] \cap \Ccal_{|\sigma_j|-1} \neq \emptyset$
so $[f(i, j)] \cap \Ccal_n \neq \emptyset$ since~$\Ccal_n \supseteq \Ccal_{|\sigma_j|-1}$. 

If $r < k$ then consider for each $i$ the non-empty set $F_i = E_i \setminus \{g(i,j) : j < r\}$.
For every $m > max(|\sigma_j| : j < r)$, $(\forall^{\infty} i)(\forall \tau \in F_i)\tau \uh m \prec X$.
Therefore we can $\emptyset'$-compute $X$, contradicting our choice of $X$.
\end{proof}

Although some principles admit c.b-enum but not simultaneous c.b-enum avoidance, 
they can simultaneously avoid computing a c.b-enum of
all effectively closed set with no computable c.b-enum.

\begin{lemma}
Let $\Ccal_0, \Ccal_1, \dots \subseteq 2^\omega$ be a countable collection of effectively closed sets with no computable c.b-enum.
There exists a (non-effectively) closed set $\Dcal \subseteq 2^\omega$ Medvedev below each $\Ccal_i$
such that $\Dcal$ has no computable c.b-enum.
\end{lemma}
\begin{proof}
Fix a set $X \not \leq \emptyset'$ and define $\Dcal$ as in Lemma~\ref{lem:cbe-avoidance-increasing-inherit}.
$\Dcal$ is Medvedev below each $\Ccal_i$.
We prove by induction over $k$ that there exists no computable $k$-enum of $\Dcal$.
Fix a computable $k$-enum $(D_i : i \in \omega)$ of $\Dcal$.
If there exists a $\sigma \not \prec X$ and infinitely many $i$
such that $\sigma \prec \tau$ for some $\tau \in D_i$. As $\sigma \not \prec X$,
there exists $\rho, \nu \in 2^{<\omega}$ and $j \in \{0,1\}$ such that
$\nu^\concat j \prec X$ and $\sigma = \rho^\concat (1-j)^\concat \nu$.
If there exists infinitely many $i$ such that $\sigma \prec \tau$ for some $\tau \in D_i$
and $\Ccal_{|\rho|} \cap [\xi] = \emptyset$ where $\tau = \rho^\concat (1-j)^\concat \xi$, 
then we can computably find infinitely many such $\tau$
and compute a $(k-1)$-enum by enumerating $D_i \setminus \tau$ for each such $i$.
If there are finitely many such $i$, then we can compute a 1-enum of $\Dcal$
by enumerating each such $\tau$.
So suppose that for every $\sigma \not \prec X$, there exist finitely many $i$
such that $\sigma \prec \tau$ for some $\tau \in D_i$. Then the jump of $\vec{D}$
computes $X$, contradicting $X \not \leq_T \emptyset'$.
\end{proof}

Before closing this section about constant-bound enumeration avoidance,
we prove that every principle admitting c.b-enum avoidance admits
simultaneous cone avoidance. Of course, given a countable collection of 
non-computable reals $A_0, A_1, \dots \subseteq \omega$,
the set of reals $\Ccal = \{A_i : i \in \omega\}$ has no computable member but may
have a computable 1-enum. For example, fix any non-computable real $A$
and set $A_i$ to be the real $A$ prefixed by $i$ zeros. 
No $A_i$ is computable, but corresponding set of reals $\Ccal$ will have a trivial 1-enum
consisting of the sequence of all 0-strings.
Thankfully, we can construct
an increasing sequence of sets of reals $\Ccal_0 \subseteq \Ccal_1 \subseteq \dots \subseteq 2^\omega$
such that computing a c.b-enum of $\vec{\Ccal}$ is equivalent to computing one of the reals~$A_i$.

\begin{lemma}\label{lem:strong-enum-avoidance-seetapun}
Fix a set~$C$ and let $A_0, A_1, \dots \subseteq \omega$ be a countable collection non-$C$-computable reals.
There exists a countable collection of closed sets of reals $\Ccal_0 \subseteq \Ccal_1 \subseteq \dots$
with no $C$-computable c.b-enum and such that for every $n$, $A_n$ computes a c.b-enum of~$\Ccal_n$.
\end{lemma}
\begin{proof}
By induction over $n$. Case $n = 0$ is satisfied by defining $\Ccal_0 = \{A_0\}$.
As every c.b-enum of $\{A_0\}$ computes $A_0$, there exists no $C$-computable c.b-enum of $\Ccal_0$.
Suppose we have defined $\Ccal_n$ and consider $A_{n+1}$. We have two cases.
Suppose first that $A_{n+1}$ $C$-computes a c.b-enum of $\Ccal_n$. In this case, set $\Ccal_{n+1} = \Ccal_n$ and we are done.

Suppose now that $A_{n+1}$ $C$-computes no c.b-enum of $\Ccal_n$. Set $\Ccal_{n+1} = \Ccal_n \cup \{A_{n+1}\}$.
If there exists a $C$-computable c.b-enum $(D_i : i < \omega)$ of $\Ccal_{n+1}$, then $A_{n+1} \not \in [D_i]$
for infinitely many $i$, otherwise it would be, up to finite change, a $C$-computable c.b-enum of $\{A_{n+1}\}$
 and would $C$-compute $A_{n+1}$.
So $A_{n+1}$ $C$-computes a c.b-enum of $\Ccal_n$ by looking on input $i$ to the least $j \geq i$
such that $A_{n+1} \not \in D_j$ and returning $D_j \restr i$. This contradicts our hypothesis.
\end{proof}

\begin{corollary}
If a $\Pi^1_2$ statement $\Psf$ admits (strong) c.b-enum avoidance, then it admits (strong) simultaneous cone avoidance.
\end{corollary}
\begin{proof}
We prove it in the case of c.b-enum avoidance.
Let~$C$ be a set, let $A_0, A_1, \dots \subseteq \omega$ be a countable collection non-$C$-computable reals
and let~$X$ be a $C$-computable $\Psf$-instance.
Let~$\Ccal_0 \subseteq \Ccal_1 \subseteq \dots \subseteq 2^\omega$ be the countable collection of closed sets
constructed in Lemma~\ref{lem:strong-enum-avoidance-seetapun}.
By Lemma~\ref{lem:cbe-avoidance-increasing-inherit}, there is a closed set~$\Dcal \subseteq 2^\omega$
Medvedev below each~$\Ccal_i$ and with no $C$-computable c.b-enum.
By c.b-enum avoidance of~$\Psf$, there is an solution~$Y$ to~$X$ such that~$\Dcal$
has no $Y \oplus C$-computable c.b-enum. Suppose for the sake of contradiction 
that $A_i \leq_T Y \oplus C$ for some~$i$. By Lemma~\ref{lem:strong-enum-avoidance-seetapun},
$Y \oplus C$ computes a c.b-enum of~$\Ccal_i$, so $Y \oplus C$ computes a c.b-enum of~$\Dcal$
by Lemma~\ref{lem:cbe-trough-medvedev}, contradiction.
\end{proof}

\section{A framework for enumeration avoidance}

The proofs of c.b-enum avoidance are usually complicated.
Liu~\cite{Liu2015Cone} proved that~$\rt^1_2$ admits strong c.b-enum avoidance 
using an involved machinery. Moreover $\rt^1_2$ is often taken 
as a bootstrap for proving the same property for whole hierarchies
such as the free set and the thin set theorems~\cite{Wang2014Some}.
In this section, we develop some tools to simplify such proofs.
For this, we introduce some reducibility notions
enabling us to \emph{propagate} c.b-enum avoidance from one statement to another.
These reductions happen to be simpler to prove than a direct argument.
Next, we generalize and abstract the lemmas proven by Liu~\cite{Liu2015Cone}
in order to freely reapply them with other statements.

\index{reduction!path}
\index{reduction!enum}
\index{path reducibility}
\index{enum reducibility}
\begin{definition}[Path and enum reducibility]
Fix two $\Pi^1_2$ statements~$\Psf$ and~$\Qsf$.
\begin{itemize}
	\item[(i)] $\Qsf$ is \emph{path reducible} to~$\Psf$ (written $\Qsf \leq_{path} \Psf$)
	if for every closed set of reals $\Ccal \subseteq 2^\omega$, 
	if $\Psf$ admits $\Ccal$ avoidance then so does~$\Qsf$
	\item[(ii)] $\Qsf$ is \emph{enum reducible} to~$\Psf$ (written $\Qsf \leq_{enum} \Psf$)
	if for every countable collection of sets of reals $\Ccal_0, \Ccal_1, \dots \subseteq 2^\omega$, 
	if $\Psf$ admits 1-enum avoidance for $\vec{\Ccal}$ then so does~$\Qsf$
\end{itemize}
\end{definition}

The strong versions of the reducibilities are defined accordingly, and are written
$\Qsf \leq_{spath} \Psf$ and~$\Qsf \leq_{senum} \Psf$, respectively.
The enum reduction, unlike the path reduction, is defined for a countable collection of sets.
The path reduction could have been defined similarly without changing the reduction,
whereas restricting the enum reduction to only one set of reals weakens strictly
the reducibility notion. These notions of reducibility are designed so that they
enjoy the two following lemmas. See~\cite{Patey2015Combinatorial} for a proof.

\begin{lemma}\label{lem:1-enum-path-bridge}
Let~$\Psf$ and~$\Qsf$ be two principles such that $\Qsf \leq_{path} \Psf$.
If $\Psf$ admits c.b-enum (resp. simultaneous c.b-enum, $n$ c.b-enum, 1-enum) avoidance
then so does~$\Qsf$.
The same statement holds if $\Qsf \leq_{spath} \Psf$ and we replace avoidance by strong avoidance.
\end{lemma}

\begin{lemma}\label{lem:cbenum-1-enum-bridge}
Let~$\Psf$ and~$\Qsf$ be two principles such that $\Qsf \leq_{enum} \Psf$.
If $\Psf$ admits c.b-enum (resp. simultaneous c.b-enum, $n$ c.b-enum) avoidance
then so does~$\Qsf$.
The same statement holds if $\Qsf \leq_{senum} \Psf$ and we replace avoidance by strong avoidance.
\end{lemma}

We need to extend the notion of~$k$-partition of the integers
to colorings over arbitrary tuples.
The forcing in Liu's theorem involved $\Pi^0_1$ classes of ordered $k$-partitions of $\omega$.
Those partitions correspond to the sets which are simultaneously homogeneous for a finite number of 2-colorings of the integers.
For example, three functions $g_0, g_1, g_2 : \omega \to 2$ induce the 6-partition
$$
X^{g_0}_0 \cap X^{g_1}_0, X^{g_0}_0 \cap X^{g_2}_0, X^{g_1}_0 \cap X^{g_2}_0,
X^{g_0}_1 \cap X^{g_1}_1, X^{g_0}_1 \cap X^{g_2}_1, X^{g_1}_1 \cap X^{g_2}_1
$$
where $X^g_i$ is the set of the integers $x$ such that $g(x) = i$.
In our case, we will not manipulate colorings of integers but of tuples of integers.
The sets homogeneous for a function $g : [\omega]^n \to k$ do not form a partition of the integers.
This is why have to make explicit the formulas expressing the homogeneity constraints.

\index{coloring atom}
\index{pseudo $k$-partition}
\begin{definition}[Coloring formula] Fix some~$d \geq 1$ and some finite set~$S$.
\begin{itemize}
	\item[1.] A \emph{coloring $d$-atom} over~$S$ is a pair $(g,J)$ (written $g[J]$) 
where $g$ is a function symbol and $J \subset S$ is a set of size $d$. 
A \emph{coloring $d$-formula} over~$S$ is a (possibly empty) conjunction of coloring $d$-atoms over~$S$.
We denote by~$\varepsilon$ the empty conjunction of coloring $d$-atoms.

	\item[2.] A \emph{valuation} of a set of coloring $d$-formulas over~$S$
with function symbols~$g_0, \dots, g_{t-1}$ is a function~$\pi$ with $\dom(\pi) \supseteq \{g_0, \dots, g_{t-1} \}$
and such that for every~$g \in \dom(\pi)$, $\pi(g)$ is a finite set $J \subset S$ of size $d$.

	\item[3.] A valuation~$\pi$ \emph{satisfies} a coloring $d$-formula~$\varphi = g_0[J_0] \wedge \dots \wedge g_{t-1}[J_{t-1}]$
	(written $\pi \models \varphi$) if~$\pi(g_i) = J_j$ for each~$i < t$.

	\item[4.] A \emph{pseudo $k$-partition} of coloring $d$-formulas is an ordered $k$-set of coloring $d$-formulas
$(\varphi_\nu : \nu < k)$ such that for every valuation~$\pi$, $\pi \models \varphi_\nu$ for some~$\nu < k$.
\end{itemize}
\end{definition}

In particular, the singleton~$\{\varepsilon\}$ is trivially a pseudo 1-partition.
Given a coloring formula~$\varphi = g_0[J_0] \wedge \dots \wedge g_k[J_k]$,
we write $\dom(\varphi)$ for the set~$\{g_0, \dots, g_k\}$. The domain
of a pseudo $k$-partition is the union of the domain of its coloring $d$-formulas.

For some fixed~$n$ and~$d$, a natural interpretation of a function symbol $g$ is a function~$f : [\omega]^n \to d+1$. 
In this interpretation, a set~$H$ \emph{satisfies} a coloring atom~$g[J]$
if~$f([H]^n) \subseteq J$. Accordingly, a set~$H$ satisfies a coloring $d$-formula
$\varphi = g_0[J_0] \wedge \dots \wedge g_{t-1}[J_{t-1}]$
within the interpretation $g_i \mapsto f_i$ if~$f_i([H]^n) \subseteq J_i$ for each~$i$.
If~$(\varphi_\nu : \nu < k)$ is a pseudo $k$-partition of coloring $d$-formulas and $\kappa$
is a function from~$\dom(\vec{\varphi})$ to functions of type~$[\omega]^n \to d+1$,
then every set~$H$ which is $f$-thin simultaneous for each~$f \in \mathrm{ran}(\kappa)$
satisfies $\varphi_\nu$ for some~$\nu < k$.
We now prove some closure properties.

\begin{lemma}\label{lem:artn-pa-avoidance-coloring-type-1}
For every pseudo $k$-partition of coloring $d$-formulas $\vec{\varphi} = (\varphi_\nu : \nu < k)$ over a finite set~$S$,
every $\mu < k$ and every function symbol $g$, the set
$\vec{\psi} = (\varphi_\nu : \nu \neq \mu) \cup (\varphi_\mu \wedge g[I] : I \subseteq S \wedge |I| = d)$ is a 
pseudo $(k+{|S| \choose d}-1)$-partition of coloring $d$-formulas.
\end{lemma}
\begin{proof}
Fix some valuation~$\pi$ with $\dom(\pi) \supseteq \dom(\vec{\varphi}) \cup \{g\}$.
As $(\varphi_\nu : \nu < k)$ is a pseudo $k$-partition, 
there exists a $\nu < k$ such that~$\pi \models \varphi_\nu$.
If $\mu \neq \nu$, then we are done since~$\varphi_\nu \in \vec{\psi}$.
If $\mu = \nu$, then $\pi \models \varphi_\mu \wedge g[\pi(g)]$
and we are also done since~$\varphi_\mu \wedge g[\pi(g)] \in \vec{\psi}$.
\end{proof}

We now need to redefine a few notions introduced by Liu in~\cite{Liu2015Cone}.
In the following, a \emph{$k$-cover} of a set $X$ is a $k$-tuple of sets $X_0, \dots, X_{k-1}$
such that $X_0 \cup \dots \cup X_{k-1} = X$. We do not require the sets $X_i$ to be pairwise disjoint.

The main combinatorial property used in the proof that a statement~$\Psf$ admits PA avoidance
is the ability to construct sets which are solutions to multiple $\Psf$-instances simultaneously.
Indeed, fix a notion of forcing to construct a solution~$G$ to a $\Psf$-instance~$X$.
Given a condition~$c$ and a Turing functional~$\Gamma$, when trying to prevent $\Gamma^G$ 
from being a $\{0,1\}$-valued DNC function, 
one can usually define for each~$e$ and~$i < 2$ the $\Pi^0_1$ class $\Ccal_{e,i}$ of $\Psf$-instances~$\tilde{X}$
such that for every extension~$d$ to~$c$ compatible with~$\tilde{X}$, either~$d \Vdash \Gamma^G(e) \uparrow$
or~$d \Vdash \Gamma^G(e) \downarrow \neq i$. If $\Ccal_{e,i} = \emptyset$ for some~$e$ such that~$\Phi_e(e) \downarrow = i$
then we are done. Otherwise, we deduce the existence of one~$e$ such that~$\Ccal_{e,0}$ and~$\Ccal_{e,1}$ are both non-empty.
In this case, we would like to make~$G$ a solution to $X$, $\tilde{X}_0$ and~$\tilde{X}_1$ simultaneously,
where~$\tilde{X}_i \in \Ccal_{e,i}$, so that we force~$\Gamma^G(e) \uparrow$. 
The conjunction of multiple $\Psf$-instances has to be done with some care to guarantee the existence of a solution.
The following notion of supporter informally describes the ``valid'' conjunctions of constraints.

\index{supporter}
\begin{definition}[Supporter]
Fix some integers $k$ and $q$.
A \emph{$k$-supporter} $\vec{\Kcal}$ of $\{1, \dots, q\}$ is $k$-tuple
$(\Kcal_\nu : \nu < k)$ where $\Kcal_\nu = \{K_{\nu,i} : i < q_\nu\}$
such that each $K_{\nu, i}$ is a subset of $\{1, \dots q\}$ and for every ordered $k$-cover
$(P_\nu : \nu < k)$ of $\{1, \dots, q\}$, there exists some $\Kcal_\nu$ and some $K_{\nu,i} \in \Kcal_\nu$
such that $\Kcal_{\nu,i} \subseteq P_\nu$.
\end{definition}

A $k$-supporter $\vec{\Kcal} = (\Kcal_\nu : \nu < k)$ of $\{1, \dots, q\}$
enables us to compose $q$ pseudo $k$-partitions of coloring $d$-formulas
$\vec{\varphi}^1 = (\varphi^1_\nu : \nu < k), \dots, \vec{\varphi}^q = (\varphi^q_\nu : \nu < k)$ 
as follows:
$$
\opcross(\vec{\varphi}^1, \dots, \vec{\varphi}^q, \vec{\Kcal}) = \left\{ 
\bigwedge_{i \in K_{\nu, j}} \varphi^i_\nu : K_{\nu, j} \in \Kcal_\nu, \nu < k 
 \right\}
$$

\begin{lemma}\label{lem:partition-cross-from-supporter}
Let $\vec{\varphi}^1 = (\varphi^1_\nu : \nu < k), \dots, \vec{\varphi}^q = (\varphi^q_\nu : \nu < k)$ 
be $q$ pseudo $k$-partitions of coloring $d$-formulas,
let $\vec{\Kcal} = (\Kcal_\nu : \nu < k)$ be a $k$-supporter of $\{1, \dots, q\}$
and let $K' = \sum_{\nu < k} |\Kcal_\nu|$. Then $\opcross(\vec{\varphi}^1, \dots, \vec{\varphi}^q, \vec{\Kcal})$
is a pseudo $K'$-partition of coloring $d$-formulas.
\end{lemma}
\begin{proof}
Fix a valuation~$\pi$ with~$\dom(\pi) \supseteq \bigcup_i \dom(\vec{\varphi}^i)$.
For every~$i \in (0, q]$, since~$(\varphi^i_\nu : \nu < k)$ is a pseudo $k$-partition of coloring $d$-formulas,
there is some~$\nu_i < k$ such that~$\pi \models \varphi^i_{\nu_i}$.
This induces an ordered $k$-partition $(P_\nu : \nu < k)$ of $\{1, \dots, q\}$
where $P_\nu = \{i \in \{1,\dots, q\} : \nu_i = \nu\}$.
By definition of being a $k$-supporter of $\{1, \dots, q\}$,
there exists some $\Kcal_\nu$ and some $\Kcal_{\nu,j} \in \Kcal_\nu$
such that $\Kcal_{\nu,j} \subseteq P_\nu$. By definition of $P_\nu$ and of the cross operator,
$$
\bigwedge_{i \in \Kcal_{\nu,j}} \varphi^i_{\nu_i} \mbox{ is the same as } \bigwedge_{i \in \Kcal_{\nu,j}} \varphi^i_{\nu}
\mbox{ which is in } \opcross(\vec{\varphi}^1, \dots, \vec{\varphi}^q, \vec{\Kcal})
$$
and $\pi \models \bigwedge_{i \in \Kcal_{\nu,j}} \varphi^i_{\nu_i}$.
Hence $\pi \models \psi$ for some~$\psi \in \opcross(\vec{\varphi}^1, \dots, \vec{\varphi}^q, \vec{\Kcal})$.
\end{proof}

When working with c.b-enums, we will end-up with $\Pi^0_1$ classes of $\Psf$-instances
$\Ccal_{V_0}, \Ccal_{V_1}, \dots$ where~$V_0, V_1, \dots$ are clopen sets,
such that for every extension~$d$ to~$c$
which is compatible with some~$\Psf$-instance~$\tilde{X} \in \Ccal_{V_i}$,
$d$ forces $\Gamma^G$ either to diverge or to intersect~$V_i$.
If we make~$G$ be simultaneously a solution to $\Psf$-instances from
$\Ccal_{V_0}, \dots, \Ccal_{V_q}$ for some clopen sets such that~$\bigcap_{i \leq q} V_i = \emptyset$,
then we will force~$\Gamma^G$ to diverge.
The following notion of~$k$-disperse sequence of clopen sets formalizes this idea.

\index{disperse sequence}
\begin{definition}[Disperse sequence]
Fix some integers $k$ and $q$.
A sequence of $q$ clopen sets $V^{(1)}, \dots, V^{(q)}$
is \emph{$k$-disperse} if for every ordered $k$-cover
$(P_\nu : \nu < k)$ of $\{1, \dots, q\}$, there exists a $\nu < k$
such that $\bigcap_{i \in P_\nu} V^{(i)} = \emptyset$.
\end{definition}

\begin{lemma}\label{lem:cbenum-avoidance-supporter-from-disperse}
Let $(e_\nu : \nu < k)$ be $k$ natural numbers and let $k' = \sum_{\nu < k} e_\nu$.
If $V^{(1)}, \dots, V^{(q)}$ is a $k'$-disperse sequence of clopen sets,
then $\vec{\Kcal} = \{ \Kcal_\nu : \nu < k \}$ where
$$
\Kcal_\nu = \{ K \subseteq \{ 1, \dots, q \} : \{V^{(i)}\}_{i \in K} \mbox{ is an } e_\nu\mbox{-disperse sequence}\}
$$
is a $k$-supporter of $\{1, \dots, q\}$.
\end{lemma}
\begin{proof}
Suppose for the sake of contradiction that there exists a $k$-cover $(P_\nu : \nu < k)$
of $\{1, \dots, q\}$ such that for all $\nu < k$, $P_\nu \not \in K_\nu$, i.e.,
$\{V^{(i)}\}_{i \in P_\nu}$ is not an $e_\nu$-disperse sequence of clopen sets.
Then for each $\nu < k$, there exists an $e_\nu$-cover $(P_{\nu, j} : j < e_\nu)$ of $P_\nu$
such that $(\forall j < e_\nu)(\bigcap_{i \in P_{\nu, j}} V^{(i)} \neq \emptyset)$.
However then $(P_{\nu, j} : j < e_\nu, \nu < k)$ is a $k'$-cover of $\{1, \dots, q\}$
that contradicts the assumption that $V^{(1)}, \dots, V^{(q)}$ is a $k'$-disperse
sequence of clopen sets.
\end{proof}

In particular, we can reprove that the pointwise conjunction
of $k+1$ pseudo $k$-partitions of coloring $d$-formulas
is again a pseudo partition.

\begin{lemma}\label{lem:art2-pa-avoidance-coloring-type-2}
Let $\vec{\psi}^0 = (\varphi^0_\nu : \nu < k), \dots, \vec{\psi}^k = (\varphi^{k}_\nu : \nu < k)$ 
be $k+1$ pseudo $k$-partitions of coloring $d$-formulas.
The set $\vec{\psi} = \{ \varphi^i_\nu \wedge \varphi^j_\nu : i < j \leq k, \nu < k \}$ 
is a pseudo $(k{k+1 \choose 2})$-partition of coloring $d$-formulas.
\end{lemma}
\begin{proof}
First, notice that~$\vec{\psi} = \opcross(\vec{\varphi}^0, \dots, \vec{\varphi}^k, \vec{\Kcal})$,
where~$\vec{\Kcal} = \{\Kcal_\nu : \nu < k\}$ is defined by
$$
\Kcal_\nu = \{ \{i, j\} : i < j \leq k \}
$$
Thanks to Lemma~\ref{lem:partition-cross-from-supporter}, it suffices to prove that~$\vec{\Kcal}$ is a 
$k$-supporter of~$\{0, \dots, k\}$.
Fix some $k$-cover $(P_\nu : \nu < k)$ of~$\{0, \dots, k\}$.
For each~$i \leq k$, let~$\nu_i < k$ be such that~$i \in P_{\nu_i}$.
By the pigeonhole principle, there are some~$i < j \leq k$ such that~$\nu_i = \nu_j$.
Hence~$\{i, j\} \subseteq P_{\nu_i}$. Since~$\{i, j\} \in \Kcal_{\nu_i}$,
we conclude.
\end{proof}

Given a set~$\Ccal \subseteq 2^\omega$ and some $n \in \omega$, define 
$$
C_n = \{ \rho \in 2^n : [\rho] \cap \Ccal \neq \emptyset \}
$$

The following last lemma is the core of the forcing argument.
It asserts that whenever the context is too weak with respect to 
some computability-theoretic notion, then either we will miss some valid elements,
or we will capture some invalid ones.

\begin{lemma}\label{lem:artn-pa-avoidance-either-disagree-correct-or-pairwise-incompatible}.
For every set~$D$ computing no c.b-enum of~$\Ccal$ and every $\Sigma^{0,D}_1$ formula~$\varphi(V)$ where~$V$ is a clopen variable,
one of the following must hold.
\begin{itemize}
	\item[1.] $\varphi(C_n)$ holds for some~$n \in \omega$.
	\item[2.] For every $k \in \omega$, there exists a $k$-disperse sequence of clopen sets $V^{(1)}, \dots, V^{(m)}$
	such that for every $i = 1, \dots, m$, $\varphi(V^{(i)})$ does not hold.
\end{itemize}
\end{lemma}
\begin{proof}
Define the following $D$-c.e. set.
$$
E = \{ W \subseteq 2^{<\omega} : (\forall \rho, \sigma \in W)|\rho| = |\sigma| \wedge  \varphi(W) \}
$$
Suppose case 1 does not hold. In other words $C_n \not \in E$ for every $n \in \omega$.
We prove that for every $k \in \omega$ and almost every $n \in \omega$, the following
is a $k$-disperse:
$$
\Wcal_n = \{ W \subseteq 2^{<\omega} : (\forall \rho \in W)|\rho| = n \wedge W \not \in E\}
$$
Note that $\Wcal_n$ is co-$D$-c.e.\ uniformly in $n$. Fix some~$k \in \omega$.
Let $\Wcal_{n,t}$ denote $\Wcal_n$ at stage $t$.
We have $\Wcal_{n,t+1} \subseteq \Wcal_{n,t}$. Therefore if there exists a $k$-cover $(P_\nu : \nu < k)$
of $\Wcal_n$ such that $(\forall \nu < k)\bigcap_{W \in P_\nu} W \neq \emptyset$, then this cover
can be found in a finite amount of time. Furthermore $C_n \in P_\nu$ for some $\nu < k$,
so 
$$
(\forall \rho \in \bigcap_{W \in P_\nu} W) [\rho] \cap \Ccal \neq \emptyset
$$
It follows that if there exists infinitely many $n$ such that such a $k$-cover exists,
we can $D$-computably find infinitely of them and define the $D$-computable enumeration $h$
which on input $n$ returns $(\rho_\nu : \nu < k)$ such that there exists some $t, m \geq n$ 
and a $k$-cover $(P_\nu : \nu < k)$ of $\Wcal_{m, t}$ such that
$(\forall \nu < k)\bigcap_{W \in P_\nu} W \neq \emptyset$
and $\rho_\nu$ is the leftmost string in $\bigcap_{W \in P_\nu} W$.
This contradicts the fact that $D$ computes no c.b-enum of $\Ccal$.
\end{proof}

\section{The weakness of Ramsey's theorem for pairs}

By Liu's theorem~\cite{Liu2015Cone}, Ramsey's theorem for pairs admits c.b-enum
avoidance. Since the notion of avoidance is downward-closed under
computable entailment, all its consequences admit c.b-enum avoidance.
On the other hand, by a consequence of Simpson's embedding lemma~\cite[Lemma 3.3]{Simpson2007extension},
various weak statements do not admit path avoidance, even for effectively
closed sets (Corollary~\ref{cor:lots-of-statements-no-path-avoidance}).
Furthermore, some statements such as the thin set theorem for pairs do not 
admit simultaneous path avoidance (Corollary~\ref{cor:sts2-not-simultaneous-cb-enum-avoidance}).
In this section, we study consequences of Ramsey's theorem for pairs
which admit stronger notions of avoidance such as 1-enum avoidance
and simultaneous c.b-enum avoidance. This enables us to reprove some existing separations
over computable entailment.

Some lemmas become pretty standard in the analysis of the combinatorics
of Ramsey's theorem and its consequences. They are independently reproven
for each new notion of avoidance. We state them in their most general form,
without providing a proof when they are too trivial.
The first lemma simply reflects the fact that strong avoidance is not sensitive
to the number of applications of a statement.

\begin{lemma}\label{lem:rt1inf-strong-enum-avoidance}
If $\rt^1_2$ admits strong $\Ccal$ avoidance for some set $\Ccal \subseteq \omega^\omega$, then so does 
$\rt^1_{<\infty}$.
\end{lemma}

\begin{corollary}[\cite{Liu2015Cone}]
$\rt^1_{<\infty}$ admits strong c.b-enum avoidance.
\end{corollary}

The next lemma is a part of the inductive proof of Ramsey's theorem
and has been reproven in~\cite[Corollary~1.6]{Liu2012RT22} for PA avoidance,
in \cite[Corollary~5.1]{Liu2015Cone} for c.b-enum avoidance, in~\cite[Theorem 4.11]{Patey2015strength}
for fairness preservation, among others.

\begin{lemma}\label{lem:coh-rtn2-rtnplus12}
If $\coh$ admits $\Ccal$ avoidance and $\rt^n_2$ strong $\Ccal$ avoidance,
then $\rt^{n+1}_2$ admits $\Ccal$ avoidance.
\end{lemma}
\begin{proof}
Let $C$ be a set computing no member of $\Ccal$ and $f : [\omega]^{n+1} \to 2$
be a $C$-computable coloring function. For each $\sigma \in [\omega]^n$, let
$R_\sigma = \{ y : f(\sigma,y) = 1 \}$. By $\Ccal$ avoidance of $\coh$ applied to $\vec{R}$,
there exists an infinite set $X$ such that $X \oplus C$ computes no member of $\Ccal$
and $\lim_{s \in X} f(\sigma, s)$ exists for each $\sigma \in [\omega]^n$.
Let $\tilde{f} : [\omega]^n \to 2$ be the function defined by $f(\sigma) = \lim_{s \in X} f(\sigma, s)$.
By strong $\Ccal$ avoidance of $\rt^n_2$, there exists an infinite set $Y \subseteq X$ and an $i \in \{0,1\}$ such that
$(\forall \sigma \in [Y]^n)\tilde{f}(\sigma) = i = \lim_{s \in X} f(\sigma, s)$ and $Y \oplus X \oplus C$ computes no
member of $\Ccal$. $Y \oplus X \oplus C$ computes an infinite set $H$ such that $f([H]^{n+1}) = i$.
\end{proof}

The following lemma is a typical example of a relative avoidance
which is known not to be enjoyed by the premisse, but which is
useful to propagate weaker notions of avoidance,
such as c.b-enum avoidance.

\begin{corollary}\label{cor:rt22-enum-avoidance}
If $\rt^1_2$ admits strong path avoidance for some set~$\Ccal \subseteq \omega^\omega$,
then $\rt^2_2$ admits path avoidance for~$\Ccal$.
\end{corollary}
\begin{proof}
It follows from Lemma~\ref{lem:coh-rtn2-rtnplus12}
and Theorem~\ref{thm:coh-strong-enum-avoidance}.
\end{proof}

\begin{corollary}[\cite{Liu2015Cone}]\label{cor:rt22-cb-enum-avoidance}
$\rt^2_2$ admits c.b-enum avoidance.
\end{corollary}
\begin{proof}
Apply Lemma~\ref{lem:1-enum-path-bridge} to Corollary~\ref{cor:rt22-enum-avoidance},
using strong c.b-enum avoidance of $\rt^1_2$.
\end{proof}

\subsection{Cohesiveness}

Cohesiveness has been proven to admit path avoidance (Corollary~\ref{cor:coh-path-avoidance}).
On the other hand, it does not admit strong path avoidance since $\rt^1_2 \leq_{sc} \coh$
and by Corollary~\ref{cor:rt12-no-strong-path-avoidance}. There is however a strong relation
between Ramsey's theorem for singletons and cohesiveness. Indeed, the latter
can be constructed by sequentially applying the former. 
This relation is formalized through the following theorem.

\begin{theorem}\label{thm:coh-strong-enum-avoidance}
$\coh \leq_{spath} \rt^1_2$
\end{theorem}
\begin{proof}
Let $\Ccal \subseteq \omega^{\omega}$ be a closed set with no $C$-computable member for some set~$C$,
and let $\vec{R}$ be a countable sequence of sets.
Our forcing conditions are tuples $(F, X)$ forming a Mathias condition,
with the additional requirement that $\Ccal$ has no $X \oplus C$-computable member.
Our initial condition is $(\emptyset, \omega)$. We can easily force our satisfying
sets to be infinite.

\begin{lemma}\label{lem:coh-strong-enum-avoidance-force}
For every condition $c = (F, X)$ and every $e, \in \omega$,
there exists an extension $(\tilde{F}, \tilde{X})$ of~$c$
forcing $\Phi_e^{G \oplus C}$ not to be a member of $\Ccal$.
\end{lemma}
\begin{proof*}
Suppose for the sake of contradiction that there is no extension of~$c$ forcing $\Phi_e^{G \oplus C}$
to be partial or $\Phi_e^{G \oplus C} \uh |\sigma| = \sigma$ for some $\sigma \in 2^{<\omega}$ 
such that $[\sigma] \cap \Ccal = \emptyset$. We show how to $X \oplus C$-compute a member of~$\Ccal$.
Define an $X \oplus C$-computable sequence of sets~$F_0 \subseteq F_1 \subseteq \dots \subseteq X$
such that~$\Phi_e^{(F \cup F_i) \oplus C}(i) \downarrow$ and $\forall x \in F_{i+1} \setminus F_i$, $x \geq max(F_i)$.
Such a sequence exists since there is no extension of~$c$ forcing $\Phi_e^{G \oplus C}$ to be partial.
We claim that the set~$Y$ defined by~$Y(i) = \Phi_e^{(F \cup F_i) \oplus C}(i)$ is a member of~$\Ccal$.
If not, then there is some~$i$ such that $\Ccal \cap [Y \uh i] = \emptyset$.
In this case, $(F \cup  F_i, X \setminus [0, max(F_i)])$ is an extension of~$c$ forcing 
$\Phi_e^{G \oplus C} \uh |\sigma| = \sigma$ for some $\sigma \in 2^{<\omega}$ 
such that $[\sigma] \cap \Ccal = \emptyset$, contradiction.
\end{proof*}

\begin{lemma}\label{lem:coh-strong-enum-avoidance-cohesive}
For every condition $c = (F, X)$ and every $e, i \in \omega$,
there exists an extension $(\tilde{F}, \tilde{X})$ of~$c$
such that $\tilde{X} \subseteq R_i$ or $\tilde{X} \subseteq \overline{R_i}$.
\end{lemma}
\begin{proof*}
Consider the coloring $f : X \to \{0,1\}$ such that
$f(x) = 1$ iff $x \in R_i$. By strong 1-enum avoidance of $\rt^1_2$ for $\vec{\Ccal}$,
there exists an infinite subset $\tilde{X} \subseteq X$
such that $\tilde{X} \oplus C$ does not compute a 1-enum of $\vec{\Ccal}$
and $\tilde{X} \subseteq R_i$ or $\tilde{X} \subseteq \overline{R_i}$.
$(F, \tilde{X})$ is the desired extension.
\end{proof*}

Let~$\Fcal = \{c_0, c_1, \dots\}$ be a sufficiently generic filter containing~$(\emptyset, \omega)$,
where~$c_s = (F_s, X_s)$. The filter~$\Fcal$ yields a unique infinite set~$G = \bigcup_s F_s$.
By Lemma~\ref{lem:coh-strong-enum-avoidance-cohesive}, $G$ is $\vec{R}$-cohesive
and by Lemma~\ref{lem:coh-strong-enum-avoidance-force},
$\Ccal$ has no $G \oplus C$-computable member.
\end{proof}

Wang~\cite{Wang2014Cohesive} proved that $\coh$ admits strong PA avoidance
with a similar argument using strong PA avoidance of~$\rt^1_2$.
Thanks to Theorem~\ref{thm:coh-strong-enum-avoidance}, we can prove the following stronger corollary.

\begin{corollary}
$\coh$ admits strong c.b-enum avoidance.
\end{corollary}
\begin{proof}
By strong c.b-enum avoidance of $\rt^1_2$,
Theorem~\ref{thm:coh-strong-enum-avoidance}
and Lemma~\ref{lem:1-enum-path-bridge}.
\end{proof}

\subsection{The Erd\H{o}s-Moser theorem}

The Erd\H{o}s-Moser theorem shares a lot of the combinatorial
aspects of Ramsey's theorem for pairs. The most notable difference
is that the notion of forcing seen on Section~\ref{sect:weakness-erdos-moser-theorem} yields only one set.
Unlike Ramsey's theorem for pairs whose requirements are disjunctions,
the requirements for the Erd\H{o}s-Moser theorem can be interleaved,
which enables us to prove that the Erd\H{o}s-Moser admits simultaneous c.b-enum avoidance.
As for cohesiveness, $\rt^1_2 \leq_{sc} \emo$ (Lemma~\ref{lem:strong-computably-reducible-strong-avoidance}), 
hence $\emo$ does not admit strong simultaneous c.b-enum avoidance.
We now prove that the Erd\H{o}s-Moser admits 1-enum avoidance relative to the Ramsey-type K\"onig's lemma.
The question whether it admits 1-enum avoidance remains open.

\begin{theorem}\label{thm:em-enum-avoidance}
$\emo \leq_{enum} \rwkl$.
\end{theorem}
\begin{proof}
Since~$\coh$ admits 1-enum avoidance, it suffices to prove the result for stable tournaments.
Fix a countable sequence of sets~$\Ccal_0, \Ccal_1, \dots \subseteq 2^\omega$
for which~$\rwkl$ admits 1-enum avoidance.
Let~$C$ be a set computing no 1-enum of~$\Ccal_i$ for any~$i$,
and let $T$ be a $C$-computable stable tournament.

We will construct an infinite $T$-transitive subtournament by using Erd\H{o}s-Moser conditions, 
defined in section~\ref{sect:weakness-erdos-moser-theorem}.
Recall that an EM condition for~$T$ is a Mathias condition $(F, X)$ where
\begin{itemize}
	\item[(a)] $F \cup \{x\}$ is $T$-transitive for each $x \in X$
	\item[(b)] $X$ is included in a minimal $T$-interval of $F$.
\end{itemize}
We furthermore impose that the~$\Ccal$'s have no $X \oplus C$-computable 1-enum.
A set $G$ \emph{satisfies} a condition $(F, X)$ if it is $T$-transitive and satisfies the
Mathias condition $(F, X)$. Our initial condition is $(\emptyset, \omega)$.
The first lemma shows that we can force the transitive subtournament to be infinite.

\begin{lemma}\label{lem:emo-avoid-enum-2}
For every condition~$c = (F, X)$, there is an extension~$(\tilde{F}, \tilde{X})$
such that~$|\tilde{F}| > |F|$.
\end{lemma}
\begin{proof*}
Let~$x \in X$. Since~$T$ is stable, there is some~$n$ such that~$\{x\} \to_T X \cap [n,+\infty)$
or $X \cap [n,+\infty) \to_T \{x\}$. By Lemma~\ref{lem:emo-cond-valid},
$d = (F \cup \{x\}, X \cap [n,+\infty)$ is a valid extension.
\end{proof*}

\begin{lemma}\label{lem:emo-avoid-enum-3}
For every condition $c = (F, X)$ and every $e,i \in \omega$,
there exists an extension $(\tilde{F}, \tilde{X})$ of~$c$
forcing $\Phi_e^{G \oplus C}$ not to be a 1-enum of $\Ccal_i$ where $G$ is the forcing variable. 
\end{lemma}
\begin{proof*}
Suppose there exists a string $\sigma \in 2^{<\omega}$ such that $[\sigma] \cap \Ccal_i = \emptyset$
and a finite set $E \subset X$ such that for every 2-partition~$E_0 \cup E_1 = E$,
there exists a finite $T$-transitive $F' \subseteq E_j$ for some~$j < 2$ such 
such that $\Phi_e^{(F \cup F') \oplus C}(|\sigma|) \downarrow = \sigma$.
Then consider the 2-partition $E_0 \cup E_1 = E$ defined by~$E_0 = \{ x \in E : (\forall^{\infty} s) T(x,s) \}$
and~$E_1 = \{x \in E : (\forall^\infty s)T(s,x) \}$.  Let~$F' \subseteq E_i$
be such that $\Phi_e^{(F \cup F') \oplus C}(|\sigma|) \downarrow = \sigma$.
In particular, there is some~$n \in \omega$ such that $F' \to_T X \cap [n,+\infty)$
or~$X \cap [n,+\infty) \to_T F'$, so by Lemma~\ref{lem:emo-cond-valid},
the condition $(F \cup F', X \cap [n, +\infty))$ is a valid extension
forcing $\Phi_e^{G \oplus C}$ not to be a 1-enum of $\Ccal_i$.

So suppose there is no such $\sigma \in 2^{<\omega}$.
For each $\sigma \in 2^{<\omega}$, let $\Tcal_{\sigma}$ denote the collection 
of the sets~$Z$ such that for every finite $T$-transitive set $F' \subseteq Z$ or~$F' \subseteq \overline{Z}$, 
$\Phi_e^{(F \cup F') \oplus C}(\card{\sigma}) \uparrow$
or $\Phi_e^{(F \cup F') \oplus C}(\card{\sigma}) \neq \sigma$.
Note that $\Tcal_{\sigma}$ are uniformly $\Pi^{0,X \oplus C}_1$ classes.
Because the previous case does not hold, then by compactness $\Tcal_\sigma \neq \emptyset$
for each $\sigma$ such that $\Ccal_i \cap [\sigma] = \emptyset$.
The set $\{\sigma : \Tcal_\sigma = \emptyset\}$ is $X \oplus C$-c.e. If for each $u \in \omega$,
there exists a $\sigma \in 2^u$ such that $\Tcal_\sigma = \emptyset$ then 
$X \oplus C$ computes a 1-enum of $\Ccal_i$, contradicting our hypothesis.
So there must be a $u$ such that $\Tcal_{\sigma} \neq \emptyset$ for each $\sigma \in 2^u$.

Thanks to 1-enum avoidance of $\rwkl$ for $\vec{\Ccal}$, define a finite decreasing sequence 
$X = X_0 \supseteq \dots \supseteq X_{2^u-1} = \tilde{X}$ such that for each $\sigma \in 2^u$
\begin{itemize}
	\item[1.] $X_\sigma$ is homogeneous for a path in $\Fcal_\sigma$.
	\item[2.] $X_\sigma \oplus C$ computes no 1-enum of any of the $\Ccal$'s.
\end{itemize}
We claim that $(F, \tilde{X})$ is an extension forcing $\Phi_e^{G \oplus C}(u) \uparrow$ or $\Phi_e^{G \oplus C}(u) \not \in 2^u$.
Suppose for the sake of contradiction that there exists a $\sigma \in 2^u$ and a set $G$ satisfying $(F, \tilde{X})$
such that $\Phi_e^{G \oplus C}(u) \downarrow = \sigma$.
By continuity, there exists a finite set $F' \subseteq G$ such that $\Phi_e^{(F \cup F') \oplus C}(u) \downarrow = \sigma$.
The set $F'$ is $T$-transitive by definition of satisfaction of $(F, \tilde{X})$.
It suffices to show that $F' \subseteq Z$ or~$F' \subseteq \overline{Z}$ for some~$Z \in \Tcal_\sigma$
to obtain a contradiction. This is immediate since~$\tilde{X}$ is homogeneous for a path in~$\Tcal_\sigma$.
\end{proof*}

Let~$\Fcal = \{c_0, c_1, \dots\}$ be a sufficiently generic filter containing $(\emptyset, \omega)$,
where~$c_s = (F_s, X_s)$. The filter~$\Fcal$ yields a unique set~$G = \bigcup_s F_s$.
By definition of a condition, the set~$G$ is a transitive subtournament of~$T$.
By Lemma~\ref{lem:emo-avoid-enum-2}, $G$ is infinite and by Lemma~\ref{lem:emo-avoid-enum-3}, 
$G \oplus C$ computes no 1-enum of~$\Ccal_i$ for any~$i \in \omega$.
\end{proof}

\begin{corollary}\label{cor:emo-simultaneous-cbenum-avoidance}
$\emo$ admits simultaneous c.b-enum avoidance.
\end{corollary}
\begin{proof}
Apply Lemma~\ref{lem:cbenum-1-enum-bridge} to Theorem~\ref{thm:em-enum-avoidance}
and simultaneous c.b-enum avoidance of~$\rwkl$ (Theorem~\ref{thm:ts-rwkl-pi01-enum-avoidance}).
\end{proof}

Since~$\rt^1_2 \leq_{sc} \emo$, $\emo$ admits neither strong 1-enum avoidance,
nor strong 2 c.b-enum avoidance. However, we can prove that~$\emo$
admits strong 1-enum avoidance relative to~$\rt^1_2$. This enables us to prove
in particular that $\emo$ admits strong c.b-enum avoidance,
and strong cone avoidance.

\begin{theorem}\label{thm:em-strong-enum-avoidance}
$\emo \leq_{senum} \rt^1_2$.
\end{theorem}

In order to prove Theorem~\ref{thm:em-strong-enum-avoidance},
we need to introduce some notation.

\begin{definition}
A \emph{$\oplus_k$-tournament} is a set $\vec{T} = T_0 \oplus \dots \oplus T_{k-1}$ such that
each $T_i$ is a tournament.
One might think of a $\oplus_k$-tournament as a conjunction of tournaments. 
Thus, notions over tournaments can be naturally
extended to $\oplus_k$-tournaments. For example, a set $U$ is a \emph{subtournament} 
of a $\oplus_k$-tournament $\vec{T}$ if it is a subtournament
of $T_i$ for each $i < k$.
\end{definition}

\begin{proof}[Proof of Theorem~\ref{thm:em-strong-enum-avoidance}]
Let~$\Ccal_1, \Ccal_2, \dots \subseteq 2^\omega$ be a countable collection of sets
for which $\rt^1_2$ admits strong 1-enum avoidance.
Fix a set $C$ computing no 1-enum of $\vec{\Ccal}$ and let $T$ be an infinite tournament.
We will construct an infinite $T$-transitive subtournament by a forcing 
whose conditions are tuples $(k, F, X, \vec{U})$ such that
\begin{itemize}
  \item[(a)] $\vec{U}$ is a $\oplus_k$-tournament
  \item[(b)] $X \oplus C$ does not compute a 1-enum of $\vec{\Ccal}$
  \item[(c)] $(F, X)$ is an EM condition for each $U_i \in \vec{U}$
\end{itemize}
A condition $(m, F', X', \vec{U'})$ \emph{extends} another condition $(k, F, X, \vec{U})$
if $(F', X')$ Mathias extends $(F, X)$, $m \geq k$ and
$\{U_i : i < k\} \subseteq \{U'_i : i < m\}$. 
A set $G$ \emph{satisfies} a condition $(k, F, X, \vec{U})$ if it is $\vec{U}$-transitive and satisfies the
Mathias condition $(F, X)$. Our initial condition is $(1, \emptyset, \omega, T)$.
The first lemma shows that every sufficiently generic filter yields an infinite set.

\begin{lemma}\label{lem:emo-strong-avoid-enum-1}
For every condition $(k, F, X, \vec{U})$, there exists an extension
$(k, \tilde{F}, \tilde{X}, \vec{U})$ such that $|\tilde{F}| > \card{F}$.
\end{lemma}
\begin{proof*}
Take any $x \in X$.
Let $f : X \to 2^k$ be the coloring defined by $f(y) = \sigma_y$
where $\card{\sigma_y} = k$ and for each $i < k$, $\sigma_y(i) = 1$ iff $U_i(x, y)$ holds.
By strong 1-enum avoidance of $\rt^1_{<\infty}$ for $\vec{\Ccal}$, there exists an infinite set $\tilde{X}$
and a $\sigma \in 2^k$ such that
$$
(\forall i < k)(\forall y \in \tilde{X})(U_i(x,y) \mbox{ holds } \biimp \sigma(i) = 1)
$$
and $\tilde{X} \oplus C$ does not compute a 1-enum of $\vec{\Ccal}$.
By Lemma~\ref{lem:emo-cond-valid}, $(F \cup \{x\}, \tilde{X})$ is a valid EM extension
for $U_i$ for each $i < k$ so
$(k, F \cup \{x\}, \tilde{X}, \vec{U})$ is a valid extension.
\end{proof*}

\begin{lemma}\label{lem:part1-strong-avoid-enum}
Fix a set $C$ computing no 1-enum of $\vec{\Ccal}$.
Let $X$ be an infinite $C$-computable set and $\vec{T}$ be a $\oplus_k$-tournament. 
For each finite subset $E \subseteq X$, there is
a $2^k$ partition $E = \bigcup_{\sigma \in 2^k} E_\sigma$ and an infinite set $Y \subseteq X$ such that
$E < Y$, $Y \oplus C$ does not compute a 1-enum of $\vec{\Ccal}$ and for all $\sigma \in 2^k$ and $i < k$,
if $\sigma(i) = 0$ then $E_\sigma \to_{T_i} Y$ and if $\sigma(i) = 1$ then $Y \to_{T_i} E_\sigma$.
\end{lemma}
\begin{proof*}
Given a set $E$, define $P_E$ to be the finite set or ordered $2^k$-partitions of $E$,
i.e.
$$
P_E = \left\{ \tuple{E_\sigma : \sigma \in 2^k} : \bigcup_{\sigma \in 2^k} E_\sigma = E 
\mbox{ and } \sigma \neq \tau \imp E_\sigma \cap E_\tau = \emptyset \right\}
$$
Define the coloring $g : X \to P_E$ by $g(x) = \tuple{E^x_\sigma : \sigma \in 2^k}$ where
$$
E^x_\sigma = \{ a \in E : (\forall i < k) T_i(a, x) \mbox{ holds iff } \sigma(i) = 0\}
$$
By strong 1-enum avoidance of $\rt^1_{<\infty}$ for $\vec{\Ccal}$, there exists an infinite set $Y \subseteq X$
homogeneous for $g$ such that $X \oplus Y$ does not compute a 1-enum of $\vec{\Ccal}$. 
Let $\tuple{E_\sigma : \sigma \in 2^k}$ be the color. 
By removing finitely many elements of $X$, we can ensure that $E < Y$ and 
by definition of $g$, for all $\sigma \in 2^k$ and $i < k$,
if $\sigma(i) = 0$ then $E_\sigma \to_{T_i} Y$ and if $\sigma(i) = 1$ then $Y \to_{T_i} E_\sigma$.
\end{proof*}

\begin{lemma}\label{lem:emo-strong-avoid-enum-2}
For every condition $(k, F, X, \vec{U})$ and every $e, i \in \omega$,
there exists an extension $(m, \tilde{F}, \tilde{X}, \vec{V})$
forcing $\Phi_e^{G \oplus C}$ not to be a 1-enum of $\Ccal_i$ where $G$ is the forcing variable. 
\end{lemma}
\begin{proof*}
Suppose there exists a string $\sigma \in 2^{<\omega}$ such that $[\sigma] \cap \Ccal_i = \emptyset$
and a finite set $E \subset X$
such that for each $2^k$-partition $E = E_0 \cup \dots \cup E_{2^k-1}$, there is a $j < 2^k$
and a $\vec{U}$ transitive set $F' \subseteq E_j$ such that $\Phi_e^{(F \cup F') \oplus C}(|\sigma|) \downarrow = \sigma$.
Take the partition $E = E_0 \cup \dots \cup E_{2^k-1}$ and the infinite set $\tilde{X} \subseteq X$
guaranteed by Lemma~\ref{lem:part1-strong-avoid-enum}. Fix a $j < 2^k$
and an $\vec{U}$-transitive set $F' \subseteq E_j$ such that $\Phi^{(F \cup F') \oplus C}_e(|\sigma|) \downarrow = \sigma$.
By Lemma~\ref{lem:emo-cond-valid}, $(F \cup F', \tilde{X})$ is a valid EM condition for $U_i$
for each $i < k$ so  $(k, F \cup F', \tilde{X}, \vec{U})$ is a valid extension 
and forces $\Phi_e^{G \oplus C}$ not to be a 1-enum of $\Ccal_i$.

So suppose there is no such $\sigma \in 2^{<\omega}$ and finite set $E \subset X$.
For each $\sigma \in 2^{<\omega}$, let $\Tcal_{\sigma}$ denote the collection 
of $\oplus_k$-tournaments $\vec{W}$ satisfying conditions (c) and (d)
such that for each finite set $E \subset X$, there exists a $2^k$-partition $E = E_0 \cup \dots \cup E_{2^k-1}$ such that
for every $j < 2^k$ and $\vec{W}$-transitive set $F' \subseteq E_j$, $\Phi_e^{(F \cup F') \oplus C}(\card{\sigma}) \uparrow$
or $\Phi_e^{(F \cup F') \oplus C}(\card{\sigma}) \neq \sigma$.
Note that $\Tcal_{\sigma}$ are uniformly $\Pi^{0,X \oplus C}_1$ classes.
Because above case does not hold, $\vec{U} \in \Tcal_{\sigma}$ for each $\sigma$ such that $\Ccal_i \cap [\sigma] = \emptyset$.
The set $\{\sigma : \Tcal_\sigma = \emptyset\}$ is $X \oplus C$-c.e. If for each $u \in \omega$,
there exists a $\sigma \in 2^u$ such that $\Tcal_\sigma = \emptyset$ then 
$X \oplus C$ computes a 1-enum of $\Ccal_i$, contradicting our hypothesis.
So there must be a $u$ such that $\Tcal_{\sigma} \neq \emptyset$ for each $\sigma \in 2^u$.
 
Given a $\sigma \in 2^u$, let $\vec{V}_\sigma \in \Tcal_\sigma$.
Define the (non-computable) predicate 
$Q(E, E_0, \dots, E_{2^k-1})$ which holds iff for each $j < 2^k$
and $\vec{V}_\sigma$-transitive set $F' \subseteq E_j$, 
$\Phi_e^{(F \cup F') \oplus C}(u) \uparrow$ or $\Phi_e^{(F \cup F') \oplus C}(u) \neq \sigma$.
For each $m \in \omega$, let $S(m)$ be the set of all $2^k$-partitions
$E_0 \cup \dots \cup E_{2^k-1}$ of the $m$ first elements $E$ of $X$
such that $Q(E, E_0, \dots, E_{2^k-1})$ holds. By definition of $\Tcal_\sigma$,
$S(m)$ is non-empty for each $m \in \omega$. Moreover, if $Q(E, E_0, \dots, E_{2^k-1})$
holds then so does $Q(E \restr s, E_0 \restr s, \dots, E_{2^k-1} \restr s)$.
Therefore $S$ is an infinite finitely branching tree.
Every infinite path in $S$ is a $2^k$-partition $X^\sigma_0 \cup \dots \cup X^\sigma_{2^k-1}$ of $X$ such that for every $j < 2^k$,
and every $\vec{V}_\sigma$-transitive set $F' \subseteq X^\sigma_j$, $\Phi_e^{(F \cup F') \oplus C}(u) \uparrow$
or $\Phi_e^{(F \cup F') \oplus C}(u) \neq \sigma$. By strong 1-enum avoidance of $\rt^1_{<\infty}$ for $\vec{\Ccal}$,
there exists a $j < 2^k$ and an infinite set $X_\sigma \subseteq X_j$ such that
$X_\sigma \oplus C$ computes no 1-enum of $\vec{\Ccal}$.

By repeating the process for each $\sigma \in 2^u$,
we obtain an infinite set $\tilde{X} \subseteq X$ such that
$\tilde{X} \oplus C$ computes no 1-enum of $\vec{\Ccal}$ and
for every $(\bigoplus_{\sigma \in 2^u} \vec{V}_\sigma)$-transitive $F' \subseteq \tilde{X}$,
$\Phi_e^{(F \cup F') \oplus C}(u) \uparrow$ or $\Phi_e^{(F \cup F') \oplus C}(u) \downarrow \not \in 2^u$.
$((2^u+1)k, F, \tilde{X}, \vec{U} \bigoplus_{\sigma \in 2^u} \vec{V}_\sigma)$ is the desired extension.
\end{proof*}

Let~$\Fcal = \{c_0, c_1, \dots\}$ be a sufficiently generic filter containing~$(1, \emptyset, \omega, T)$,
where~$c_s = (k_s, F_s, X_s, \vec{U}_s)$. The filter~$\Fcal$ yields a unique set~$G = \bigcup_s F_s$.
By definition of a forcing condition, the set~$G$
is a transitive subtournament of~$T$. By Lemma~\ref{lem:emo-strong-avoid-enum-1},
$G$ is infinite and by Lemma~\ref{lem:emo-strong-avoid-enum-2}, 
$G \oplus C$ computes no 1-enum of~$\vec{\Ccal}$.
\end{proof}

\begin{corollary}
$\emo$ admits strong c.b-enum avoidance.
\end{corollary}
\begin{proof}
Apply Lemma~\ref{lem:cbenum-1-enum-bridge} to Theorem~\ref{thm:em-strong-enum-avoidance},
knowing that $\rt^1_2$ admits strong c.b-enum avoidance.
\end{proof}

\subsection{The rainbow Ramsey theorem for pairs}

Among the consequences of Ramsey's theorem for pairs,
the rainbow Ramsey theorem for pairs surprisingly
admits a very nice computability-theoretic characterization
in terms of algorithmic randomness.
Miller~\cite{MillerAssorted} proved that it is computably equivalent
to the relativized diagonally non-computable principle,
which itself corresponds to the assertion of the existence of an infinite subset
of a 2-random.

Given a closed set of positive measure, the measure
of oracles which do not compute a c.b-enum of this closed set is null.
In particular, neither weak weak K\"onig's lemma
nor its relativized statements admit c.b-enum avoidance.
We will however see that the rainbow Ramsey theorem for pairs
admits 1-enum avoidance. This results cannot be proven 
for higher exponents since by Theorem~\ref{thm:rrt2n-tsn}, $\sts^2 \leq_{sc} \rrt^3_2$.
Therefore $\rrt^3_2$ does not admit simultaneous c.b-enum avoidance and \emph{a fortiori} does not admit 1-enum avoidance.
We start by proving  that the rainbow Ramsey theorem for singletons
admits strong 1-enum avoidance and then apply the standard inductive
argument.

\begin{theorem}\label{thm:rrt12-strong-enum-avoidance}
$\rrt^1_2$ admits strong 1-enum avoidance.
\end{theorem}

\begin{definition}
A \emph{2-bounded $\oplus_k$-function} is a set $\vec{f} = f_0 \oplus \dots \oplus f_{k-1}$ such that
each $f_i$ is a coding of a 2-bounded coloring over integers.
One might think of an 2-bounded $\oplus_k$-function as a conjunction of 2-bounded functions. 
Thus notions over functions can be naturally
extended to 2-bounded $\oplus_k$-functions: -- e.g. A set $F$ is a \emph{rainbow} 
for a 2-bounded $\oplus_k$-function $\vec{f}$ if it is an $f_i$-rainbow for each $i < k$ --. 
\end{definition}

\begin{proof}[Proof of Theorem~\ref{thm:rrt12-strong-enum-avoidance}]
Let $C$ be a set computing no 1-enum of $\Ccal$ for some set $\Ccal \subseteq 2^\omega$
and $f : \omega \to \omega$ be a 2-bounded coloring.
Our forcing conditions are tuples $(k, F, X, \vec{g})$ such that
\begin{itemize}
  \item[(a)] $\vec{g}$ is a normal 2-bounded $\oplus_k$-function
  \item[(b)] $X$ is an infinite set such that $F < X$ and $X \oplus C$ computes no 1-enum of $\Ccal$
  \item[(c)] $F$ is a finite $\vec{g}$-rainbow.
\end{itemize}
A set $G$ \emph{satisfies} a condition $(k, F, X, \vec{g})$ if it satisfies
the Mathias condition $(F, X)$ and $G$ if $g_i$-free for each $i < k$.
Our initial condition is $(1, \emptyset, \omega, f)$.
A condition $(m, F', X', \vec{g}')$ \emph{extends} another condition $(k, F, X, \vec{g})$
if $(F', X')$ Mathias extends $(F, X)$,  $m \geq k$ and $(\forall i < k)g_i = g'_i$.

\begin{lemma}\label{lem:rrt12-avoid-enum-1}
For every condition $(k, F, X, \vec{g})$ there exists an extension
$(k, \tilde{F}, \tilde{X}, \vec{g})$ such that $|H| > \card{F}$.
\end{lemma}
\begin{proof*}
Take $x \in X \setminus \bigcup_i g_i(F)$.
$F \cup \{x\}$ is a $\vec{g}$-rainbow, hence $(k, F \cup \{x\}, X \setminus [0, x], \vec{g})$ is the desired extension.
\end{proof*}

\begin{lemma}\label{lem:rrt12-avoid-enum-2}
For every condition $(k, F, X, \vec{g})$ and every $e \in \omega$,
there exists an extension $(m, \tilde{F}, \tilde{X}, \vec{h})$
forcing $\Phi_e^{G \oplus C}$ not to be a 1-enum of $\Ccal$, where $G$ is the forcing variable. 
\end{lemma}
\begin{proof*}
Suppose there exists a $\sigma \in 2^{<\omega}$ such that $[\sigma] \cap \Ccal = \emptyset$
and a finite set $F' \subseteq X$ such that $F \cup F'$ is $g_i$-free for each $i < k$ and 
$\Phi_e^{(F_0 \cup F') \oplus C}(|\sigma|) \downarrow = \sigma$.
$(k, F \cup F', X \setminus [0, max(F')], \vec{g})$ is a condition forcing $\Phi_e^{G \oplus C}$ not to be a 1-enum of $\Ccal$.

Suppose there is no such finite set $F' \subset X$.
For each $\sigma \in 2^{<\omega}$, let $\Fcal_{\sigma}$ denote the collection 
of 2-bounded $\oplus_k$-functions $\vec{h}$ such that $F$ is $\vec{h}$-free and
for each finite set $F' \subset X$ such that $F \cup F'$ is $h_j$-free for each $j < k$,
either $\Phi_e^{(F \cup F') \oplus C}(\card{\sigma}) \uparrow$ or $\Phi_e^{(F \cup F') \oplus C}(\card{\sigma}) \neq \sigma$.
Note that $\Fcal_\sigma$ are uniformly $\Pi^{0,X \oplus C}_1$ classes.
Because the above case does not hold, $\vec{g} \in \Fcal_\sigma$ for each $\sigma$ such that $\Ccal \cap [\sigma] = \emptyset$.
The set $\{\sigma : \Fcal_\sigma = \emptyset\}$ is $X \oplus C$-c.e. If for each $u \in \omega$
there exists a $\sigma \in 2^u$ such that $\Fcal_\sigma = \emptyset$ then $X \oplus C$ computes a 1-enum of $\Ccal$,
contradicting our hypothesis.
So there must be an $u \in \omega$ such that $\Fcal_\sigma \neq \emptyset$ for each $\sigma \in 2^u$.

For each $\sigma \in 2^u$, let $\vec{h}_\sigma \in \Fcal_\sigma$.
$((2^u+1)k, F, X, \vec{g} \bigoplus_{\sigma \in 2^u} \vec{h}_\sigma)$
is a condition forcing $\Phi_e^{G \oplus C}(u) \uparrow$ or $\Phi_e^{G \oplus C}(u) \downarrow \not \in 2^u$.
\end{proof*}

Let~$\Fcal = \{c_0, c_1, \dots\}$ be a sufficiently generic filter 
containing $(1, \emptyset, \omega, f)$, where $c_s = (k_s, F_s, X_s, \vec{g}_s)$.
The filter~$\Fcal$ yields a unique set~$G = \bigcup_s F_s$. By Lemma~\ref{lem:rrt12-avoid-enum-1},
the set~$G$ is infinite. By definition of a forcing condition, $G$
is an $f$-rainbow, and by Lemma~\ref{lem:rrt12-avoid-enum-2}, $G \oplus C$ computes no 1-enum of $\Ccal$.
\end{proof}

\begin{lemma}\label{lem:dnr-reducible-to-rrt12zp}
$\dnr \leq_{sc} \rrt^1_2[\emptyset']$
\end{lemma}
\begin{proof}
Fix a set $X$ and a canonical enumeration of all finite sets $(D_i : i \in \omega)$. 
We construct a 2-bounded coloring $f : \omega \to \omega$
such that for every $e \in \omega$,
if $\Phi^X_e(e) \downarrow$ and $D_{\Phi^X_e(e)}$ has at least $2(e+1)$ elements,
then either $D_{\Phi^X_e(e)} \cap [0,e] \neq \emptyset$ or it is not an $f$-rainbow.
We first show how, given an infinite $f$-rainbow $H$, we compute a function $g$ d.n.c. relative to~$X$.
For every $e \in \omega$, $g(e) = i$ where $D_i$ are the first $2(e+1)$ elements of $H$ greater than $e$.
Suppose for the sake of contradiction that $g(e) = \Phi^X_e(e)$ for some $e$.
Then $D_{\Phi^X_e(e)}$ is not an $f$-rainbow and therefore $D_{\Phi^X_e(e)} \neq D_{g(e)}$.
Contradiction.

We now detail the construction of $f$ by stages.
At stage $0$, $dom(f_0) = \emptyset$. Suppose that at stage $s$,
$[0,s) \subseteq dom(f_s)$ and $|dom(f_s)| \leq 3s$.
If $\Phi^X_s(s) \downarrow$ and $|D_{\Phi^X_s(s)}| \geq 2(s+1)$ and has no element before $s$,
then by cardinality, there exist $u, v \in D_{\Phi^X_s(s)} \setminus dom(f_s)$.
Set $f(u) = f(v)$ and give a fresh color to $f(s)$ if $s \not \in dom(f_s)$.
Then go to stage $s+1$.
$f = \bigcup_s f_s$ is the desired coloring. Note that $f$ is $X'$-computable. 
\end{proof}

\begin{corollary}\label{cor:dnr-strong-1-enum-avoidance}
$\dnr$ admits strong 1-enum avoidance.
\end{corollary}
\begin{proof}
By Theorem~\ref{thm:rrt12-strong-enum-avoidance},
Lemma~\ref{lem:dnr-reducible-to-rrt12zp} and Lemma~\ref{lem:strong-computably-reducible-strong-avoidance}.
\end{proof}

\begin{corollary}\label{cor:rrt22-enum-avoidance}
$\rrt^2_2$ admits 1-enum avoidance.
\end{corollary}
\begin{proof}
By Miller~\cite{MillerAssorted}, $\rrt^2_2 =_c \dnr[\emptyset']$.
By Corollary~\ref{cor:dnr-strong-1-enum-avoidance}, $\dnr[\emptyset']$ admits strong 1-enum avoidance,
so \emph{a fortiori} 1-enum avoidance.
Lemma~\ref{lem:computably-reducible-avoidance} enables us to conclude.
\end{proof}

\section{The weakness of the thin set hierarchy}

Wang~\cite{Wang2014Some} first showed, by proving that the thin set and the free set theorems
admit strong cone avoidance,  that increasing the number of colors in the output set
changes fundamentally the combinatorics of Ramsey's theorem
and gives a strictly weaker statement.
We strengthened Wang's result in chapter~\ref{chap:thin-set-free-set-theorems} by showing that the thin set
and the free set theorems admit preservation of arbitrary many hyperimmunities.
In this section, we prove a similar result about simultaneous c.b-enum avoidance.

The argument is however significantly more involved since we cannot make a free use of 
weak K\"onig's lemma which does not admit c.b-enum avoidance.
Flood~\cite{Flood2012Reverse} clarified the situation for Ramsey's theorem for pairs by introducing the 
Ramsey-type weak K\"onig's lemma ($\rwkl$) which happens to be the right amount of compactness
needed for~$\rt^2_2$ and even weaker statements such as the Erd\H{o}s-Moser theorem.
In order to extend c.b-enum avoidance to the thin set and free set theorems,
we need to generalize the Ramsey-type weak K\"onig's lemma to arbitrary colorings.

Fix an enumeration of all $n$-tuples of integers~$t_0, t_1, \dots \subseteq [\omega]^n$.
Every string $\sigma \in k^\omega$ can be interpreted as a function
$f_\sigma : \{t_i : i < |\sigma|\} \to k$ defined by~$f_\sigma(t_i) = \sigma(i)$.
An infinite sequence~$S \in k^\omega$ is then interpreted as the unique function~$f_S : [\omega]^n \to k$
such that~$f_S(t_i) = S(i)$ for each~$i \in \omega$.

\index{pi01(ts)@$\Pi^0_1(\ts^n_k)$}
\begin{definition}
For every~$n, k \geq 2$, $\Pi^0_1(\ts^n_k)$ is the statement ``For every infinite tree~$T \subseteq k^{<\omega}$,
there is an infinite set~$H$ and some~$c < k$ such that for every length~$\ell \in \omega$,
there is some~$\sigma \in T$ of length~$\ell$ such that~$f_\sigma(t_i) \neq c$ for each~$i < |\sigma|$
such that~$t_i \in [H]^n$''.
\end{definition}

Informally, for every non-empty $\Pi^0_1$ class $\Pcal$ of $\ts^n_k$-instances,
$\Pi^0_1(\ts^n_k)$ asserts the existence of an infinite $f$-thin set for some~$f \in \Pcal$.
The statement $\Pi^0_1(\ts^n_k)$ is formulated so that it does not imply the existence
of the function~$f \in \Pcal$. In particular, $\Pi^0_1(\ts^1_2)$ is nothing but the Ramsey-type
weak K\"onig's lemma.
Fix some~$m \geq 1$ and define~$d_n$ inductively as follows.
$$
d_1 = m \hspace{20pt} \tilde{d}_1 = 1 \hspace{10pt} \mbox{ and } \hspace{10pt} 
	\tilde{d}_n = \sum_{0 < s < n} d_s d_{n-s} \hspace{20pt} d_n = (m+1)\tilde{d_n} \hspace{10pt} \mbox{ for } n > 1
$$

\begin{theorem}\label{thm:ts-rwkl-pi01-enum-avoidance}Fix some~$n \geq 1$
\begin{itemize}
	\item[(i)] $\rwkl$ admits simultaneous c.b-enum avoidance.
	\item[(ii)] $\Pi^0_1(\ts^n_{\tilde{d}_n+1})$ admits $m$ c.b-enum avoidance.
	\item[(iii)] $\ts^n_{d_n+1}$ admits strong $m$ c.b-enum avoidance.
\end{itemize}
\end{theorem}

Before proving Theorem~\ref{thm:ts-rwkl-pi01-enum-avoidance}, we state
an immediate corollary. As usual, strong avoidance for $n$-tuples
becomes avoidance for~$(n+1)$-tuples thanks to cohesiveness.

\begin{lemma}
If~$\coh$ admits $\Ccal$ avoidance and~$\ts^n_k$ strong~$\Ccal$ avoidance,
then~$\ts^{n+1}_k$ admits $\Ccal$ avoidance.
\end{lemma}

\begin{corollary}
If~$\ts^n_k$ admits strong path avoidance for some set~$\Ccal \subseteq \omega^\omega$,
then $\ts^{n+1}_k$ admits path avoidance for~$\Ccal$.
\end{corollary}

\begin{corollary}\label{cor:tsm-cb-enum-avoidance}
$\ts^{n+1}_{d_n+1}$ admits $m$ c.b-enum avoidance.
\end{corollary}

Finally, notice that ~$\tilde{d}_n = (m+1)^{2n-2}a_n$, where $a_1 = 1$ and
$$
a_n = \sum_{i=1}^{n-1} a_i a_{n-i}
$$
The reader will have recognized \emph{Catalan numbers}~\cite{Koshy2009Catalan}.
Therefore we can obtain an explicit growth for~$d_n$:
$$
d_n = (m+1)\tilde{d_n} = \frac{(m+1)^{2n-1} {2n \choose n}}{n+1}
$$

We now turn to the proof of Theorem~\ref{thm:ts-rwkl-pi01-enum-avoidance}.
Fix a countable collection of sets~$\Ccal_0, \Ccal_1, \dots \subseteq 2^\omega$.
The proof of Theorem~\ref{thm:ts-rwkl-pi01-enum-avoidance} is done in several steps
by a mutual induction as follows:
\bigskip

\begin{itemize}
	\item[(A1)] Assuming that for each~$t \in (0,n)$, $\ts^t_{d_t+1}$ admits strong c.b-enum avoidance for $\vec{\Ccal}$,
	we prove that $\Pi^0_1(\ts^n_{\tilde{d}_n+1})$ admits c.b-enum avoidance for~$\vec{\Ccal}$.
	\item[(A2)] Assuming that $\rwkl$ admits $m$ c.b-enum avoidance,
	we prove that~$\ts^1_{d_1+1}$ admits strong $m$ c.b-enum avoidance.
	\item[(A3)] Assuming that for each~$t \in (0,n)$, $\ts^t_{d_t+1}$ admits strong 1-enum avoidance for $\vec{\Ccal}$
	and that~$\Pi^0_1(\ts^n_{\tilde{d}_n+1})$ admits 1-enum avoidance for~$\vec{\Ccal}$,
	we prove that~$\ts^n_{d_n+1}$ admits strong 1-enum avoidance for~$\vec{\Ccal}$.
\end{itemize}

\bigskip
By (A1), we deduce that $\rwkl$ admits simultaneous c.b-enum avoidance.
Moreover, 1-enum avoidance can be replaced by c.b-enum avoidance in step (A3) 
by Lemma~\ref{lem:cbenum-1-enum-bridge}.

\subsection{Strong path avoidance of generalized cohesiveness}

As we did for preservation of hyperimmunity, we need to
state a few theorems asserting the existence of sets
satisfying some cohesiveness properties generalized to the thin set theorem.
The proof is very similar to the proof of Theorem~\ref{thm:generalized-cohesivity-strong-avoidance},
and is obtained by replacing Lemma~\ref{lem:coh-preservation-lemma} by Lemma~\ref{lem:coh-strong-enum-avoidance-force}.

\begin{theorem}\label{thm:generalized-cohesivity-strong-enum-avoidance}
Fix a coloring $f : [\omega]^n \to \omega$, some $t \leq n$ and
a closed set~$\Ccal \subseteq \omega^\omega$ for which
$\ts^s_{d_s+1}$ admits strong path avoidance for each $s \in (0, t]$. 
For every set $C$ computing no member of~$\Ccal$,
there exists an infinite set $G$ such that $G \oplus C$ computes no member of~$\Ccal$
and for every $\sigma \in [\omega]^{<\omega}$ such that $n-t \leq \card{\sigma} < n$, 
$$
\card{\set{x : (\forall b)(\exists \tau \in [G \cap (b, +\infty)]^{n-|\sigma|}) f(\sigma, \tau) = x}} \leq d_{n-|\sigma|}
$$
\end{theorem}

When considering a function $f : [\omega]^n \to k$ and taking $t = n-1$, we 
obtain a set similar to the one constructed in section 3.1 in~\cite{Wang2014Some}.

\begin{theorem}\label{thm:enum-avoidance-gen-coh}
Fix a coloring $f : [\omega]^n \to k$ and a closed set~$\Ccal \subseteq \omega^\omega$
for which $\ts^s_{d_s+1}$ admits strong path avoidance
for each $s \in (0, n)$. For every set $C$ which does not compute a member of $\Ccal$,
there exists an infinite set $G$ such that $G \oplus C$ computes no member of $\Ccal$
and a sequence $(I_\sigma : 0 < |\sigma| < n)$ such that for each $\ell \in (0,n)$ and each $\sigma \in [\omega]^\ell$
\begin{itemize}
	\item[(a)] $I_\sigma$ is a subset of $\{0, \dots, k-1\}$ with at most $d_{n-\ell}$ elements  
	\item[(b)] $(\exists b)(\forall \tau \in [G \cap (b,+\infty)]^{n-\ell}) f(\sigma,\tau) \in I_\sigma$
\end{itemize}
\end{theorem}
\begin{proof}
Let $G$ be the set constructed by Theorem~\ref{thm:generalized-cohesivity-strong-enum-avoidance}
for $t = n-1$. For each $\sigma \in [\omega]^{<\omega}$ such that $0 < \card{\sigma} < n$,
let
$$
I_\sigma = \{x < k :  (\forall b)(\exists \tau \in [G \cap (b, +\infty)]^{n-|\sigma|}) f(\sigma, \tau) = x
$$
By choice of $G$, $|I_\sigma|$ has at most $d_{n-|\sigma|}$ many elements.
Moreover, for each $y < k$ such that $y \not \in I_\sigma$, there exists a bound $b_y$
such that $(\forall \tau \in [G \cap (b_y, +\infty)]^{n-|\sigma|}) f(\sigma, \tau) \neq y$.
So taking $b = max\{b_y : y < k \wedge y \not \in I_\sigma\}$, we have
$$
(\forall \tau \in [G \cap (b,+\infty)]^{n-|\sigma|}) f(\sigma,\tau) \in I_\sigma
$$
\end{proof}

\begin{theorem}\label{thm:gen-coh-art}
Fix a coloring $f : [\omega]^n \to k$ and a closed set~$\Ccal \subseteq \omega^\omega$
for which $\ts^s_{d_s+1}$ admits strong path avoidance
for each $s \in (0, n)$. For every set $C$ which does not compute a member of $\Ccal$,
there exists an infinite set $G$ such that $G \oplus C$ computes no member of $\Ccal$
and a finite set  $(\Ical_s : 0 < s < n)$ such that for each $s \in (0,n)$
\begin{itemize}
	\item[(a)] $\Ical_s$ is a finite set of at most $d_s$ sets of colors,
	and $|I| \leq d_{n-s}$ for each $I \in \Ical_s$. 
	\item[(b)] $(\forall \sigma \in [G]^s)(\exists b)(\exists I \in \Ical_s)
	(\forall \tau \in [G \cap (b,+\infty)]^{n-s}) f(\sigma,\tau) \in I$
\end{itemize}
\end{theorem}
\begin{proof}
Let $X$ be the infinite set and $(I_\sigma : 0 < |\sigma| < n)$
be the infinite sequence constructed in~Theorem~\ref{thm:enum-avoidance-gen-coh}.
For each $s \in (0,n)$ and~$\sigma \in [G]^s$, let $F_s(\sigma) = I_\sigma$.
Using strong path avoidance of $\ts^s_{d_s+1}$ for $\Ccal$,
we build a finite sequence $X \supseteq X_1 \supseteq \dots \supseteq X_{n-1}$
such that for each $s \in (0,n)$
\begin{itemize}
	\item[1.] $X_s \oplus C$ computes no member of~$\Ccal$
	\item[2.] $|F_s([X_s]^s)| \leq d_s$
\end{itemize}
Let $G = X_{n-1}$ and $\Ical_s = F_s([G]^s)$ for each $s \in (0,n)$.
We now check that property (b) is satisfied.
Fix a $\sigma \in [G]^s$. Because $G \subseteq X$,
$(\exists b)(\forall \tau \in [G \cap (b,+\infty)]^{n-s}) f(\sigma,\tau) \in I_\sigma$.
So $F_s(\sigma) = I_\sigma$, but $\sigma \in [G]^s$, hence $I_\sigma \in \Ical_s$.
\end{proof}

In particular, in our ongoing forcing, we will use the following corollary.

\begin{corollary}\label{cor:gen-coh-art}
Fix a countable collection of sets~$\Ccal_0, \Ccal_1, \dots \subseteq 2^\omega$ for which
$\ts^s_{d_s+1}$ admits strong c.b-enum avoidance for each $s \in (0, n)$.
For every coloring $f : [\omega]^n \to k$ and every set $C$ which does not compute a c.b-enum of the~$\Ccal$'s,
there exists an infinite set $G$ such that $G \oplus C$ computes no c.b-enum of the~$\Ccal$'s
and a finite set  $I \subseteq \{0,\dots,k-1\}$ such that
$|I| \leq \tilde{d}_n$ and for each $s \in (0,n)$,
$$
(\forall \sigma \in [G]^s)(\exists b)(\forall \tau \in [G \cap (b,+\infty)]^{n-s}) f(\sigma, \tau) \in I
$$
\end{corollary}
\begin{proof}
Apply Theorem~\ref{thm:gen-coh-art} taking $I = \bigcup \Ical$
and Lemma~\ref{lem:1-enum-path-bridge} using strong c.b-enum avoidance of $\ts^s_{d_s+1}$
for~$\vec{\Ccal}$.
\end{proof}

\subsection{Strong c.b-enum of \texorpdfstring{$\Pi^0_1(\ts^n_k)$}{Pi01(TSnk)}}

We now prove step (A1) in the proof of Theorem~\ref{thm:ts-rwkl-pi01-enum-avoidance}.

\begin{theorem}
Fix a sequence of sets~$\Ccal_0, \Ccal_1, \dots \subseteq 2^\omega$.
For every~$n \geq 1$, if~$\ts^s_{d_s+1}$ admits strong c.b-enum avoidance for $\vec{\Ccal}$,
then $\Pi^0_1(\ts^s_{d+1})$ admits c.b-enum avoidance for $\vec{\Ccal}$,
where~$d = 1$ if~$n = 1$ and $d = \sum_{0 < s < n} d_s d_{n-s}$ otherwise.
\end{theorem}
\begin{proof}
Fix a set $C$ computing no c.b-enum of the~$\Ccal$'s,
and let $\Pcal_0$ be a $\Pi^{0,C}_1$ class of colorings of type~$[\omega]^n \to d+1$.
We want to build an infinite $f$-thin set~$G$ for some~$f \in \Pcal_0$
such that~$G \oplus C$ computes no c.b-enum of the~$\Ccal$'s.
Before describing the strategy to build such a set, we need to provide some suitable semantics to the notion
of coloring formula. We shall only consider coloring $d$-formulas over~$\{0, \dots, d\}$.
In the following, we will omit the parameters~$d$ and~$\{0, \dots, d\}$.

\index{assignment (coloring formula)}
\begin{definition}[Assignment]
An \emph{assignment} of a coloring formula $\varphi$
is a function~$\kappa$ such that~$\dom(\kappa) \supseteq \dom(\varphi)$ and for every~$g \in \dom(\kappa)$,
$\kappa(g)$ is a function of type $[\omega]^n \to d+1$.
Given a coloring formula $\varphi = g_0[J_0] \wedge \dots \wedge g_{t-1}[J_{t-1}]$ and 
an assignment~$\kappa$,  a set $F \subseteq \omega$ \emph{satisfies} 
$\varphi$ (written $(F, \kappa) \models \varphi$) if $\kappa(g_j)([F]^n) \subseteq J_j$ for each $j < t$.
\end{definition}

In other words, $(F, \kappa) \models \varphi$
iff there is a valuation~$\pi \models \varphi$ such that $\kappa(g_j)([F]^n) \subseteq \pi(g_j)$ for each~$j < t$.
By Lemma~\ref{lem:artn-pa-avoidance-coloring-type-1}, 
the set~$\vec{\varphi}^0 = (g[J] : J \subseteq [0,d] \wedge |J| = d)$ is a pseudo $(d+1)$-partition
of coloring $d$-formulas over $\{0, \dots, d\}$ with one function symbol~$g$.
We will construct an infinite set~$G$ such that $G \oplus C$ computes no c.b-enum of the $\Ccal$'s.
We furthermore ensure that $(G_i, f) \models \varphi$ for some~$\varphi \in \vec{\varphi}_0$
and some assignment~$\kappa$ such that~$\kappa(g) \in \Pcal_0$.
Therefore, the set~$G$ will be $f$-thin for some~$f \in \Pcal_0$ and will be a solution to~$\Pcal_0$.
The requirements to ensure that $G$ is infinite are
$$
Q_s : (\exists w > s)(w \in G) 
$$

The requirements to ensure that $G \oplus C$ computes no c.b-enum of the $\Ccal$'s are for each~$e,j \in \omega$
$$
R_{e,j} : \Phi^{G \oplus C}_e \mbox{ total}  \imp
 (\exists w)|\Phi^{G \oplus C}_e(w) \neq \Phi_w(w)| > e
	\vee [\Phi^{G \oplus C}_e(w)] \cap \Ccal_j = \emptyset
$$

\begin{definition}Given a Turing functional $\Phi_e$, a finite set $F$, a clopen $V$, a coloring formula $\varphi$,
an assignment~$\kappa$ and a set $X$, we say that $\Phi^{F \oplus C}_e$ \emph{abandons} $V$ on $\varphi$, $\kappa$ and $X$ if
there is a $w \in \omega$ and a finite set $F' \subset X$ such that
$(F', \kappa) \models \varphi$ and 
$$
|\Phi_e^{(F \cup F') \oplus C}(w)| > e \vee [\Phi_e^{(F \cup F') \oplus C}(w)] \cap V = \emptyset
$$
\end{definition}

The following lemma tells us that computing an $e$-enum and not abandoning 
an $e$-disperse sequence of clopen sets is incompatible.

\begin{lemma}\label{lem:pi01ts-disperse-prevents-enum}
Let $V^{(1)}, \dots, V^{(q)}$ be an $e$-disperse sequence of clopen sets.
Suppose $\Phi_e^{F \oplus C}$ does not abandon $V^{(j)}$ on $\varphi$, $\kappa$ and $X$ for every $j = 1, \dots, q$. 
Then for every set $G$ satisfying the Mathias condition $(F, X)$ such that
$(G, \kappa) \models \varphi$, $\Phi_e^{G \oplus C}$
is not total or is not an $e$-enum.
\end{lemma}
\begin{proof*}
Fix such a set $G$ and a $j \in \{1, \dots, q\}$.
Because $\Phi_e^{F \oplus C}$ does not abandon $V^{(j)}$ on $\varphi$, $\kappa$ and $X$, 
for every $w \in \omega$ such that
$\Phi_e^{G \oplus C}(w) \downarrow$, the following holds 
$$
|\Phi_e^{G \oplus C}(w)| \leq e \wedge [\Phi_e^{G \oplus C}(w)] \cap V^{(j)} \neq \emptyset
$$
By convention, if $\rho \in \Phi_e^{G \oplus C}(w)$ then $|\rho| = w$.
Taking $w$ large enough, we have for every $j \in \{1,\dots,q\}$ and every
$\rho \in \Phi_e^{G \oplus C}(w)$
$$
[\rho] \cap V^{(j)} \neq \emptyset \imp [\rho] \subseteq V^{(j)}
$$
For $i < e$, let $\rho_i$ be the $i$th string in $\Phi_e^{G \oplus C}(w)$.
The string $\rho_i$ induces an $e$-cover $(P_i : i < e)$ of the clopen sets defined by
$$
P_i = \{V^{(j)} : [\rho_i] \subseteq V^{(j)}\}
$$
But then for each $i < e$, $[\rho_i] \subseteq \bigcap_{j \in P_i} V^{(j)} \neq \emptyset$
contradicting the assumption that $V^{(1)}, \dots, V^{(q)}$ is $e$-disperse.
\end{proof*}

We are now ready to define the actual notion of forcing
and prove that every sufficiently generic filter yields the desired $p$-tuple of sets.

\begin{definition}[Single condition]\ 
\begin{itemize}
  \item[1.] A \emph{single condition} is a tuple $(F, X, \varphi, \kappa)$ 
	where $(F, X)$ is a Mathias condition, $\varphi$ is a coloring formula
	and~$\kappa$ is an assignment such that for each $s \in (0,n)$,
	$$
	(\forall \sigma \in [F]^s)(\forall \tau \in [F \cup X]^{n-s})(\sigma\tau, \kappa) \models \varphi
	$$

  \item[2.] A single condition $d = (H, Y, \psi, \gamma)$ \emph{extends} $c = (F, X, \varphi, \kappa)$
  if $(H, Y)$ Mathias extends $(F, X)$, $\kappa \subseteq \gamma$ 
	and there exists a coloring formula $\theta$ such that $\psi = \varphi \wedge \theta$.

  \item[3.] A set~$G$ \emph{satisfies} a single condition $(F, X, \varphi, \kappa)$
  if $G$ satisfies the Mathias condition $(F,X)$ and $(G, \kappa) \models \varphi$.
\end{itemize}
\end{definition}

\begin{definition}[Condition]\ 
\begin{itemize}
  \item[1.] A \emph{condition} is a tuple $(k, \vec{F}, X, D, \vec{\varphi}, \Pcal)$
  where $k > 0$, $\vec{F}$ is a $k$-tuple of finite sets $(F_\nu : \nu < k)$,
	$D$ computes no c.b-enum of the~$\Ccal$'s, $X \oplus C \leq_T D$, $\vec{\varphi} = (\varphi_\nu : \nu < k)$
	is a pseudo $k$-partition of coloring formulas,
  $\Pcal$ is a non-empty $\Pi^{0,D}_1$ class of assignments
  and for each $\kappa \in \Pcal$, each $\nu < k$, $(F_\nu, X, \varphi_\nu, \kappa)$ is a single condition.
  
  \item[2.] A condition $d = (m, \vec{H}, Y, E, \vec{\psi}, \Qcal)$ \emph{extends} $c = (k, \vec{F}, X, D, \vec{\varphi}, \Pcal)$
  if $D \leq_T E$ and there is a function $f : m \to k$ with the following property:
  for each $\gamma \in \Qcal$, there is some $\kappa \in \Pcal$
  such that the single condition $(H_\nu, Y, \psi_\nu, \gamma)$ extends 
	$(F_{f(\nu)}, X, \varphi_{f(\nu)}, \kappa)$.
  In this case, the function~$f$ \emph{witnesses} the extension and \emph{part $\nu$ of $d$
  refines part $f(\nu)$ of~$c$}.
  
  \item[3.] A set $G$ \emph{satisfies} some condition $c = (k, \vec{F}, X, D, \vec{\varphi}, \Pcal)$
  on part $\nu$ if there is some $\kappa \in \Pcal$ such that $G$ satisfies the single
  condition $(F_\nu, X, \varphi_\nu, \kappa)$.
	$G$ satisfies $c$ if it satisfies $c$ on some of its parts.
  
  \item[4.] A condition $(k, \vec{F}, X, D, \vec{\varphi}, \Pcal)$ \emph{forces $Q_u$ on part $\nu$}
  if there exists some $w > u$ such that $w \in F_\nu$. 
  
  \item[5.]
  A condition $d$ \emph{forces $R_{e,j}$ on part $\nu$}
  if every set~$G$ satisfying $d$ on part $\nu$ satisfies $R_{e,j}$.

  \item[6.] \emph{Part $\nu$ of $(k, \vec{F}, X, D, \vec{\varphi}, \Pcal)$ is acceptable} if 
  there is an infinite set $Y \subseteq X$ such that $Y \oplus D$ computes no c.b-enum of the $\Ccal$'s
	and there is a $\kappa \in \Pcal$ such that for each $s \in (0,n)$,
	$$
	(\forall \sigma \in [Y]^s)(\exists b)(\forall \tau \in [Y \cap (b,+\infty)]^{n-s})(\sigma\tau, \kappa) \models \varphi_\nu
	$$
\end{itemize}
\end{definition}

\begin{lemma}\label{lem:pi01ts-cb-enum-avoidance-acceptable-part}
Every condition has an acceptable part.
\end{lemma}

\begin{lemma}\label{lem:pi01ts-cb-enum-avoidance-Q}
For every condition $c$ and every $u \in \omega$, there is
a condition $d$ extending $c$ such that $d$ forces $Q_u$
on each of its acceptable parts.
\end{lemma}

\begin{lemma}\label{lem:pi01ts-cb-enum-avoidance-R}
For every condition $c$ and every $e,j \in \omega$
there exists an extension $d$ forcing $R_{e,j}$ on each
of its acceptable parts.
\end{lemma}

The construction of $G$ given the three lemmas above
is strictly the same as in~\cite[Lagniappe]{Hirschfeldt2015Slicing}:
We build an infinite, decreasing sequence of conditions $c_0 \geq c_1 \geq \dots$ 
starting with $c_0 = (d+1, \emptyset, \dots, \emptyset, \omega, C, \vec{\varphi}_0, \Pcal_0)$ 
with the following properties assuming that $c_s = (k_s, \vec{F_s}, X_s, D_s, \vec{\varphi}_s, \Pcal_s)$:
\begin{itemize}
  \item[1.] Each $c_s$ has an acceptable part.
  \item[2.] If part $\nu$ of $c_s$ is acceptable, then $c_s$ forces $R_{e,j}$ on part~$\nu$ if~$s = \tuple{e,j}$.
  \item[3.] If part $\nu$ of $c_s$ is acceptable, then $c_s$ forces $Q_s$ on part $\nu$.
\end{itemize}
If part $\nu$ of $c_{s+1}$ is acceptable and refines part $\mu$ of $c_s$, then part $\mu$ of $c_s$ is also acceptable.
Hence the acceptable parts of the conditions form an infinite finitely branching tree.
By König's lemma, there exists an infinite sequence $\nu_0, \nu_1, \dots$ where part $\nu_{s+1}$ of $c_{s+1}$ refines
part $\nu_s$ of condition $c_s$. One easily checks that $G = \bigcup_s F_{\nu_s, s}$ is the desired set,.

\begin{proof*}[Proof of Lemma~\ref{lem:pi01ts-cb-enum-avoidance-acceptable-part}]
Let $c = (k, \vec{F}, X, D, \vec{\varphi}, \Pcal)$ be a condition. As $\Pcal$ is non-empty,
there exists an assignment~$\kappa \in \Pcal$.
Thanks to Corollary~\ref{cor:gen-coh-art}, define a finite decreasing sequence
$X \supseteq Y_0 \supseteq \dots \supseteq Y_{t-1}$ such that for each $i < t$
\begin{itemize}
	\item[1.] $Y_i \oplus D$ computes no c.b-enum of the $\Ccal$'s
	\item[2.] there is a set~$J_i$ of size~$d$ such that for each $s \in (0,n)$,
	$$
	(\forall \sigma \in [Y_i]^s)(\exists b)(\forall \tau \in [Y_i \cap (b,+\infty)]^{n-s})\kappa(g_i)(\sigma,\tau) \in J_i
	$$
\end{itemize}
Let~$\pi$ be the valuation defined by~$\pi(g_i) = J_i$ for each~$i < t$.
Since~$\vec{\varphi} = (\varphi_\nu : \nu < k)$ is a pseudo $k$-partition, there is some~$\nu < k$
such that~$\pi \models \varphi_\nu$. We claim that~$\nu$ and~$Y_{t-1}$ satisfy the desired properties.
For each~$s \in (0, n)$ and~$i < t$, by definition of~$\pi$, 
$$
	(\forall \sigma \in [Y_{t-1}]^s)(\exists b)(\forall \tau \in [Y_{t-1} \cap (b,+\infty)]^{n-s})\kappa(g_i)(\sigma,\tau) \in \pi(g_i)
$$
Therefore, for each~$s \in (0, n)$,
$$
	(\forall \sigma \in [Y_{t-1}]^s)(\exists b)
	(\forall \tau \in [Y_{t-1} \cap (b,+\infty)]^{n-s})(\forall i < t)\kappa(g_i)(\sigma,\tau) \in \pi(g_i)
$$
Since $(\sigma\tau, \kappa) \models \varphi$ iff 
$(\forall i < t)\kappa(g_i)(\sigma,\tau) \in \pi(g_i)$ for some valuation~$\pi \models \varphi$,
$$
	(\forall \sigma \in [Y_{t-1}]^s)(\exists b)
	(\forall \tau \in [Y_{t-1} \cap (b,+\infty)]^{n-s})(\sigma\tau, \kappa) \models \varphi_\nu
$$
Therefore part $\nu$ of $c$ is acceptable.
\end{proof*}

\begin{proof*}[Proof of Lemma~\ref{lem:pi01ts-cb-enum-avoidance-Q}]
Fix some~$u \in \omega$. It suffices to prove that given a condition $c = (k, \vec{F}, X, D, \vec{\varphi}, \Pcal)$,
if part $\mu$ is acceptable, then there exists an extension $d = (k, \vec{H}, Y, Y \oplus D, \vec{\varphi}, \Qcal)$ which forces
$Q_u$ on part $\mu$ and whose extension is witnessed by the identity map.
By iterating the process, we obtain an extension satisfying the statement of the lemma.

Fix an acceptable part $\mu$. By definition, there exists an assignment $\kappa \in \Pcal$ and
an infinite subset $Y_0 \subseteq X$ such that $Y_0 \oplus D$ computes no c.b-enum of the~$\Ccal$'s
and for each $s \in (0,n)$,
$$
(\forall \sigma \in [Y_0]^s)(\exists b)(\forall \tau \in [Y_0 \cap (b,+\infty)]^{n-s})(\sigma\tau, \kappa) \models \varphi_\mu
$$

By the fact that $(F_\mu, Y_0, \varphi_\mu)$ is a single condition, for each $s \in (0,n)$,
$$
(\forall \sigma \in [F_\mu]^s)(\forall \tau \in [F_\mu \cup Y_0]^{n-s})(\sigma\tau, \kappa) \models \varphi_\mu
$$
therefore by taking $y \in Y_0 \cap (u,+\infty)$ and removing finitely many elements from~$Y_0$,
we obtain a set $Y$ such that for each $s \in (0,n)$,
$$
(\forall \sigma \in [F_\mu \cup \{y\}]^s)
	(\forall \tau \in [F_\mu \cup Y]^{n-s})(\sigma\tau, \kappa) \models \varphi_\mu
$$
Let $H_\nu = F_\nu$ if $\nu \neq \mu$ and $H_\mu = F_\mu \cup \{y\}$ otherwise.
Let $\Qcal$ be the $\Pi^{0, Y \oplus D}_1$ collection of all the assignments $\kappa \in \Pcal$
such that the above formula holds.
The condition $(k, \vec{H}, Y, Y \oplus D, \vec{\varphi}, \Qcal)$ is an extension forcing $Q_u$ on part $\mu$.
\end{proof*}

It remains to prove Lemma~\ref{lem:pi01ts-cb-enum-avoidance-R}.
Given a condition $c$, and any $e,j \in \omega$,
let $U_{e,j}(c)$ be the set of all acceptable parts $\nu$
such that $c$ does not force $R_{e,j}$ on part $\nu$.
If $U_{e,j}(c) = \emptyset$, we are already done as condition $c$ already forces
$R_{e,j}$ on each of its acceptable parts. 
In order to prove Lemma~\ref{lem:pi01ts-cb-enum-avoidance-R},
it suffices to prove and iterate the following lemma.

\begin{lemma}\label{lem:pi01ts-cb-enum-avoidance-R-step}
For every condition $c$ and every $e,j \in \omega$ such that $U_{e,j}(c) \neq \emptyset$, 
there exists an extension $d$ such that $\card{U_{e,j}(d)} < \card{U_{e,j}(c)}$.  
\end{lemma}

The proof of Lemma~\ref{lem:pi01ts-cb-enum-avoidance-R-step} is split into two cases,
according to Lemma~\ref{lem:artn-pa-avoidance-either-disagree-correct-or-pairwise-incompatible}.
In the first case, we can find a
	piece of oracle in a part of $U_{e,j}(d)$, forcing the Turing functional 
	we consider to halt on a ``wrong'' input, i.e., 
	on a clopen set which does not intersect $\Ccal_j$. In the second case, there exist
	many clopen sets which are intersected by the Turing functional whenever it halts.
The first lemma states the existence of a finite extension forcing the Turing functional 
	to halt on a wrong input on a part of $U_{e,j}(d)$ when the first case holds.
The second lemma states the existence of an extension forcing the Turing functionals 
	to diverge or not to be an $e$-enum 
	on each of the parts of $U_{e,j}(d)$ when the second case holds.
Before stating and proving the two lemmas, we need to extend the abandoning terminology to a condition.

\begin{definition} Let $c = (k, \vec{F}, X, D, \vec{\varphi}, \Pcal)$ be a condition and $V$ be a clopen set.
\begin{itemize}
	\item[1.] We say that part $\mu$ of $c$ \emph{abandons} $V$ 
	on some assignment~$\kappa$ if $\Phi_e^{F_\mu \oplus C}$ abandons $V$ on $\varphi_\mu$, $\kappa$ and $X$.

	\item[2.] 
	We say that part $\mu$ of $c$ \emph{abandons} $V$ if for every
	assignment $\kappa \in \Pcal$, part $\mu$ of $c$ abandons $V$ on $\kappa$.
	The condition $c$ \emph{abandons} $V$ if it abandons $V$ on some part $\mu \in U_{e,j}(c)$.
\end{itemize}
\end{definition}

Given a condition  $c = (k, \vec{F}, X, D, \vec{\varphi}, \Pcal)$,
define the following $\Pi^{0,D}_1$ class of assignments for~$\vec{\varphi}$:
$$
\Pcal_V = \{  \kappa \in \Pcal : (\forall \nu \in U_{e,j}(c))
 \Phi_e^{F_\nu \oplus C} \mbox{ does not abandon } V \mbox{ on } \varphi_\nu, \kappa \mbox{ and } X \}
$$

Notice that $c$ abandons $V$ iff $\Pcal_V = \emptyset$,
hence the predicate ``$c$ abandons $V$'' is $\Sigma^{0,D}_1$.
We are now about to prove the first lemma, but need one last definition.
The acceptation of a part $\nu$ of a condition $c = (k, \vec{F}, X, D, \vec{\varphi}, \Pcal)$ intuitively means
that we can find an infinite set $Y \subseteq X$ such that $(k, \vec{F}, Y, Y \oplus D, \vec{\varphi}, \Pcal)$ is a valid extension
and there exists an assignment $\kappa \in \Pcal$
such that for each $s \in (0,n)$,
$$
(\forall \sigma \in [Y]^s)(\exists b)(\forall \tau \in [Y \cap (b,+\infty)]^{n-s})(\sigma\tau, \kappa) \models \varphi_\nu
$$
The condition $(k, \vec{F}, Y, Y \oplus D, \vec{\varphi}, \Pcal)$ has the same number of parts
and its part $\nu$ can take $Y$ as a witness of being acceptable.
This process can be iterated so that we obtain a condition $d = (k, \vec{F}, Z, E, \vec{\varphi}, \Pcal)$
such that for every acceptable part $\nu$ of $d$, there exists
an assignment~$\kappa \in \Pcal$ such that $Z$ is a witness of acceptation of part $\nu$.
Such a condition is said to \emph{witness its acceptable parts}.
Every condition can be extended to a condition witnessing its acceptable parts.
For each $q \in \omega$, define
$$
C_q = \{ \rho \in 2^q : [\rho] \cap \Ccal_j \neq \emptyset \}
$$

\begin{lemma}\label{lem:artn-pa-avoidance-first-case-extension-forcing}
Let $c = (k, \vec{F}, X, D, \vec{\varphi}, \Pcal)$ be a condition witnessing its acceptable parts,
and let $\mu \in U_{e,j}(c)$ such that part $\mu$ of $c$ abandons $C_q$ for some $q \in \omega$.
There exists an extension $d$ with the same parts as $c$,
such that for every set $G$ satisfying $d$ on part $\mu$,
$\Phi_e^{G \oplus C}$ is not an $e$-enum of $\Ccal_j$.
\end{lemma}
\begin{proof*}
By definition of witnessing its acceptable parts, there exists an assignment~$\kappa \in \Pcal$
such that for each $s \in (0,n)$,
$$
(\forall \sigma \in [X]^s)(\exists b)(\forall \tau \in [X \cap (b,+\infty)]^{n-s})(\sigma\tau, \kappa) \models \varphi_\mu
$$
As part $\mu$ of $c$ abandons $C_q$, then $
\Phi_e^{F_\mu \oplus C}$ abandons $C_q$ on $\varphi_\mu$, $\kappa$ and $X$.
Unfolding the definition, there exists a $w \in \omega$ and finite set $F' \subseteq X$, 
such that $(F', \kappa) \models \varphi_\mu$ and 
$$
|\Phi_e^{(F_\mu \cup F') \oplus C}(w)| > e \vee
[\Phi_e^{(F_\mu \cup F') \oplus C}(w)] \cap C_q (\supseteq \Ccal) = \emptyset
$$
Set $H_\nu = F_\mu \cup F'$ if $\mu = \nu$ and $H_\nu = F_\nu$ otherwise.
By removing finitely many elements from $X$, we obtain a set $Y \subseteq X$ such that for each $s \in (0,n)$,
$$
(\forall \sigma \in [F_\mu \cup F']^s)(\forall \tau \in [F_\mu \cup F' \cup Y]^{n-s})(\sigma\tau, \kappa) \models \varphi_\mu
$$
Let $\Qcal$ be the $\Pi^{0,D}_1$ class of all the assignments $\kappa \in \Pcal$
satisfying the above property.
The condition $(k, \vec{H}, Y, D, \vec{\varphi}, \Pcal)$ is a valid extension forcing $R_{e,j}$ on part $\mu$.
\end{proof*}

We now prove the second lemma stating the existence of an extension
forcing $\Phi_e^{G \oplus C}$ to be partial or not to be a $e$-enum
on each of the parts refining a part in $U_{e,j}(c)$.

\begin{lemma}\label{lem:artn-pa-avoidance-second-case-extension-forcing}
Let $V^{(1)}, \dots, V^{(q)}$ be an $e$-disperse sequence of clopen sets
and let $c$ be a condition which does not abandon $V^{(j)}$ for every $j = 1, \dots, q$.
There exists an extension $d$ of~$c$ such that for every set $G$ satisfying $d$,
$\Phi_e^{G \oplus C}$ is either partial or is not an $e$-enum of $\Ccal_j$.
\end{lemma}
\begin{proof*}
Fix some condition~$c = (k, \vec{F}, X, D, \vec{\varphi}, \Pcal)$.
By the hypothesis that $c$ does not abandon $V^{(t)}$ for every $t \in \{1, \dots, q\}$,
the $\Pi^{0,D}_1$ class $\Pcal_{V^{(t)}}$ is non-empty.
Let $\vec{\Kcal}$ be the $e$-supporter of $\{1, \dots, q\}$ constructed in Lemma~\ref{lem:cbenum-avoidance-supporter-from-disperse}
and let $K' = \sum_{\Kcal_i \in \Kcal}|\Kcal_i|$. 
By renaming the function symbols, we can suppose that
$\vec{\varphi}^{V^{(1)}}, \dots, \vec{\varphi}^{V^{(q)}}$ have pairwise disjoint domains.
By Lemma~\ref{lem:partition-cross-from-supporter}, 
$\vec{\psi} = \opcross(\vec{\varphi}^{V^{(1)}}, \dots, \vec{\varphi}^{V^{(q)}}, \vec{\Kcal})$
is a pseudo $K'$-partition of coloring formulas.
Let~$\Qcal$ be the $\Pi^{0,D}_1$ class of all assignments~$\kappa$
such that~$\kappa = \kappa^1 \sqcup \dots \sqcup \kappa^q$
for some~$\kappa^t \in \Pcal_{V^{(t)}}$.
For each $\mu < K'$, let $H_\mu = F_\nu$ if part $\mu$ refines part $\nu$ of~$c$.
Then the condition~$d = (K', \vec{H}, X, D, \vec{\psi}, \Qcal)$
is a valid extension of~$c$. We claim that $d$ forces $\Phi_e^{G \oplus C}$ to be either partial or not
to be an $e$-enum of $\Ccal_j$ for each set~$G$ satisfying~$d$.

Fix such a set $G$ satisfying $c$ on some part $\mu < K'$ of $d$.
If part $\mu$ of $d$ refines a part of $c$ which is not in $U_{e,j}(c)$ then by definition
part $\mu$ of $d$ already forces $\Phi_e^{G \oplus C}$ not to be an $e$-enum of $\Ccal_j$.
So suppose that part $\mu$ of $d$ refines a part $\nu$ of $c$ such that $\nu \in U_{e,j}(c)$.

By definition of satisfaction and the definition of the cross operator, 
there exists a $K \in \Kcal_\nu$ and for each $t \in \{1, \dots, q\}$
an assignment $\kappa^t \in \Pcal_{V^{(t)}}$ such that $G$ satisfies the single condition 
$$
(F_\nu, X, \bigwedge_{t \in K} \varphi^t_\nu, \kappa^t)
$$
In particular, $G$ satisfies the Mathias condition $(F_\nu, X)$ and
$$
(G, \kappa^t) \models \bigwedge_{t \in K} \varphi^t_\nu
$$

By construction of $\vec{\Kcal}$, $\{V^{(t)} : t \in K\}$ is an $e$-disperse
sequence of clopen sets.  By definition of $\Pcal_{V^{(t)}}$,
$\Phi_e^{F_\nu \oplus C}$ does not abandon $V^{(t)}$ on $\varphi^t_\nu$, $\kappa^t$ and $X$.
So in particular $\Phi_e^{F_\nu \oplus C}$ does not abandon $V^{(t)}$ on 
$\bigwedge_{t \in K} \varphi^t_\nu$, $\kappa^1 \sqcup \dots \sqcup \kappa^q$ and~$X$.
Applying Lemma~\ref{lem:pi01ts-disperse-prevents-enum}, we deduce that
$G$ is not total or does not compute an $e$-enum.
\end{proof*}

\begin{proof*}[Proof of Lemma~\ref{lem:pi01ts-cb-enum-avoidance-R-step}]
Fix a condition $c = (k, \vec{F}, X, D, \vec{\varphi}, \Pcal)$.
We can furthermore assume without loss of generality that~$c$ witnesses its acceptable parts.
If $U_{e,j}(c) = \emptyset$, we are done. So suppose $U_{e,j}(c) \neq \emptyset$.
By Lemma~\ref{lem:artn-pa-avoidance-either-disagree-correct-or-pairwise-incompatible}
applied to the $\Sigma^{0,D}_1$ formula ``$c$ abandons $V$'', we have two cases:
\begin{itemize}
	\item[1.] There exists a $q \in \omega$ and a part $\mu \in U_{e,j}(c)$ such that 
	part $\mu$ of $c$ abandons $C_q$. In this case, by Lemma~\ref{lem:artn-pa-avoidance-first-case-extension-forcing}
	there exists an extension $d$ having the same parts as $c$, and such that
	for every set $G$ satisfying $d$ on part $\mu$, $\Phi_e^{G \oplus C}$ is not an $e$-enum of $\Ccal_j$.
	Therefore, $U_{e,j}(d) = U_{e,j}(c) \setminus \{\mu\}$ and we are done.

	\item[2.] There exists an $e$-disperse sequence of clopen sets $V^{(1)}, \dots, V^{(q)}$
	such that for every $t = 1, \dots, q$, $c$ does not abandon $V^{(t)}$.
	By Lemma~\ref{lem:artn-pa-avoidance-second-case-extension-forcing},
	there exists an extension $d$ such that for every set~$G$ satisfying $d$,
	$\Phi_e^{G \oplus C}$ is either partial or is not an $e$-enum of~$\Ccal_j$.
	In this case we have $U_{e,j}(d) = \emptyset$.
\end{itemize}
\end{proof*}

This last lemma finishes the proof.
\end{proof}

\subsection{Strong c.b-enum avoidance of $\ts^1$}

Next, we prove step (A2) of Theorem~\ref{thm:ts-rwkl-pi01-enum-avoidance}.

\begin{theorem}
$\ts^1_{m+1}$ admits strong $m$ c.b-enum avoidance for every~$m \geq 1$. 
\end{theorem}
\begin{proof}
Fix a set~$C$ and $m$ sets~$\Ccal_0, \dots, \Ccal_{m-1} \subseteq 2^\omega$
with no $C$-computable c.b-enum.
Let~$A_0 \cup \dots \cup A_m = \omega$ be an $(m+1)$-cover of~$\omega$,
and assume that there is no infinite set~$H \subseteq \overline{A}_i$ for some~$i \leq m$
such that~$H \oplus C$ computes no c.b-enum of the $\Ccal$'s, otherwise we are done.
We will construct simultaneously $m+1$ sets $G_0, \dots, G_m$ such that
$G_i \oplus C$ computes no c.b-enum of the $\Ccal$'s for some $i \leq m$. We furthermore ensure that for each $i \leq m$,
$G_i \subseteq \overline{A_i}$.
The requirements to ensure that all $G_i$'s are infinite are
$$
Q_s : (\forall i < p)(\exists w > s)(w \in G_i) 
$$

The requirements to ensure that $G_i \oplus C$ computes no c.b-enum of the $\Ccal$'s for some $i \leq m$ are
$$
R_{e_0, \dots, e_m} : \bigwedge_{j < m} R_{e_0, j} \vee \dots \vee \bigwedge_{j < m} R_{e_m,j}
$$
where
$$
R_{e_i, j} : \Phi^{G_i \oplus C}_{e_i} \mbox{ total}  \imp
 (\exists w)|\Phi^{G_i \oplus C}_{e_i}(w) \neq \Phi_w(w)| > e_i
	\vee [\Phi^{G_i \oplus C}_{e_i}(w)] \cap \Ccal_j = \emptyset
$$

The sets~$G_0, \dots, G_m$ are built by a variant of Mathias forcing
whose conditions are tuples~$c = (F_0, \dots, F_m, X)$ where $(F_i, X)$
is a Mathias condition such that~$F_i \subseteq \overline{A_i}$
and~$X \oplus C$ computes no c.b-enum of the~$\Ccal$'s for each~$i \leq m$.
A condition~$d = (H_0, \dots, H_m, Y)$ \emph{extends} a condition~$c = (F_0, \dots, F_m, X)$
if~$(H_i, Y)$ Mathias extends $(F_i, X)$ for each~$i \leq m$.
A set~$G_i$ \emph{satisfies part~$i$ of~$c$} if~$G_i$ satisfies the Mathias condition~$(F_i, X)$.
The first lemma asserts that every sufficiently generic filters yields infinite sets.

\begin{lemma}\label{lem:ts1-cb-enum-avoidance-Q}
For every condition $c = (F_0, \dots, F_m, X)$, every~$i \leq m$ and every $u$, there is
an extension $d = (H_0, \dots, H_m, Y)$ of $c$ such that
$H_i \cap (u, +\infty)$.
\end{lemma}
\begin{proof}
If~$X \cap \overline{A}_i \cap (u, +\infty) = \emptyset$, then $X \subseteq A_i$
contradicting our initial assumption. Therefore, let~$x \in X \cap \overline{A}_i \cap (u, +\infty)$.
The condition~$d = (F_0, \dots, F_{i-1}, F_i \cup \{x\}, F_{i+1}, \dots, F_m, X \setminus (x, +\infty))$
is the desired extension.
\end{proof}

\begin{definition}
Given a Turing functional~$\Phi_e$, a finite set~$F$, a clopen~$V$,
we say that~$\Phi_e^{F \oplus C}$ \emph{abandons} $V$ on a set~$Y$
if there is a $w \in \omega$ and a finite set~$F' \subseteq X$
such that
$$
|\Phi_e^{(F \cup F') \oplus C}(w)| > e \vee [\Phi_e^{(F \cup F') \oplus C}(w)] \cap V = \emptyset
$$
Given a condition~$c = (F_0, \dots, F_m, X)$, a vector of indices~$\vec{e}$ and some~$i_0 < i_1 \leq m$,
we say that~$(i_0, i_1)$ \emph{abandons} $V$ on~$c$ if for every 2-cover~$Z_{i_0} \cup Z_{i_1} = X$,
either $\Phi_{e_i}^{F_i \oplus C}$ abandons~$V$ on~$Z_i$ for some~$i \in \{i_0, i_1\}$.
\end{definition}

Given a condition~$c = (F_0, \dots, F_m, X)$, a vector of indices~$\vec{e}$ and some~$i_0 < i_1 \leq m$,
define the following~$\Pi^{0,X \oplus C}_1$ class:
$$
\Pcal_V = \{Z_{i_0} \oplus Z_{i_1} : Z_{i_0} \cup Z_{i_1} = X \wedge
	(\forall i \in \{i_0, i_1\}) \Phi_{e_i}^{F_i \oplus C} \mbox{ does not abandon } V \emph{ on } Z_i \}
$$

Note that~$(i_0, i_1)$ abandons $c$ on~$V$ iff $\Pcal_V = \emptyset$.
In particular, the formula ``$(i_0, i_1)$ abandons $c$ on~$V$'' is $\Sigma^{0,X \oplus C}_1$.

\begin{lemma}\label{lem:ts1-cb-enum-avoidance-R-iter}
For every condition $c$ and every~$i_0 < i_1 \leq m$, every~$j < m$ and every vector of indices~$\vec{e}$
there exists an extension $d$ of~$c$ forcing $R_{e_{i_0}, j}$ or~$R_{e_{i_1}, j}$.
\end{lemma}
\begin{proof*}
Fix a condition~$c = (F_0, \dots, F_m, X)$ and let~$C_q = \{ \rho \in 2^q : [\rho] \cap \Ccal_j \neq \emptyset \}$.
By Lemma~\ref{lem:artn-pa-avoidance-either-disagree-correct-or-pairwise-incompatible}
applied to the $\Sigma^{0, X \oplus C}_1$ formula ``$(i_0, i_1)$ abandons $c$ on~$V$'',
we have two cases:
\begin{itemize}
	\item Case 1: $(i_0, i_1)$ abandons $c$ on some~$C_q$. 
	Let~$Z_{i_0} \cup Z_{i_1} = X$ be the 2-cover defined for each~$i \in \{i_0, i_1\}$ by~$Z_i = X \cap \overline{A}_i$.
	Let~$i \in \{i_0, i_1\}$ be such that~$\Phi_{e_i}^{F_i \oplus C}$ abandons~$C_q$ on~$X \cap \overline{A}_i$.
	By definition, there is some $w \in \omega$ and a finite set~$F' \subseteq X \cap \overline{A}_i$
	such that
	$$
		|\Phi_{e_i}^{(F_i \cup F') \oplus C}(w)| > e_i \vee [\Phi_{e_i}^{(F_i \cup F') \oplus C}(w)] \cap C_q = \emptyset
	$$
	The extension~$d = (F_0, \dots, F_{i-1}, F_i \cup F', F_{i+1}, \dots, F_m, X \setminus (0, max(F'))$
	forces~$\Rcal_{e_i, j}$.

	\item Case 2: There exists an $(e_{i_0} + e_{i_1})$-disperse sequence of clopen sets $V^{(1)}, \dots, V^{(q)}$
	such that for every $t = 1, \dots, q$, $(i_0, i_1)$ does not abandon~$V^{(t)}$ on~$c$.
	Let~$\vec{\Kcal} = \{ \Kcal_{i_0}, \Kcal_{i_1} \} $ be defined for each~$i \in \{i_0, i_1\}$ by
	$$
	\Kcal_i = \{ K \subseteq \{ 1, \dots, q \} : \{V^{(t)}\}_{t \in K} \mbox{ is an } e_i\mbox{-disperse sequence}\}
	$$
	By Lemma~\ref{lem:cbenum-avoidance-supporter-from-disperse}, $\vec{\Kcal}$ is a 2-supporter of~$\{1, \dots, q\}$.
	By c.b-enum avoidance of~$\rwkl$ for the~$\Ccal$'s, for each~$t \in \{1, \dots, q\}$, 
	there are $2$-covers $Z^{V^{(t)}}_{i_0} \oplus Z^{V^{(t)}}_{i_1} \in \Pcal_{V^{(t)}}$
	such that $\bigoplus_{t = 1}^q Z^{V^{(t)}}_{i_0} \oplus Z^{V^{(t)}}_{i_1} \oplus C$ computes no c.b-enum of the~$\Ccal$'s.
	For each~$t \in \{1, \dots, q\}$, let~$s_t \in \{i_0, i_1\}$ be such that~$Y = \bigcap_{t = 1}^q Z^{V^{(t)}}_{s_t}$ is infinite.
	We claim that the condition~$d = (F_0, \dots, F_m, Y)$ forces
	either~$R_{e_{i_0}, j}$ or~$R_{e_{i_1}, j}$.
	For this, consider the $2$-partition $P_{i_0}, P_{i_1}$ of~$\{1, \dots, q\}$ defined for each~$i \in \{i_0, i_1\}$ by
	$P_i = \{ t \in \{1, \dots, q\} : s_t = i \}$. Since~$\vec{\Kcal}$ is a 2-supporter
	of~$\{1, \dots, q\}$, there is some~$i \in \{i_0, i_1\}$ and some~$K \in \Kcal_i$
	such that~$K \subseteq P_i$. Fix a set~$G_i$ satisfying part~$i$ of~$d$.
	In particular, since~$\Phi_{e_i}^{F_i \oplus C}$ does not abandon~$V^{(t)}$ on~$\bigcap_{t = 1}^q Z^{V^{(t)}}_{s_t}$
	for each~$t \in K$, either~$\Phi_{e_i}^{G_i \oplus C}$ is partial,
	or~$[\Phi_{e_i}^{G_i \oplus C}(w)] \cap V^{(t)} \neq \emptyset$ for every~$t \in K$.
	However, $\{V^{(t)} : t \in K\}$ is an $e_i$-disperse sequence, hence $\bigcap_{t \in K} V^{(t)} = \emptyset$
	so $\Phi_{e_i}^{G_i \oplus C}$ is partial and therefore~$R_{e_i, j}$ is forced.
\end{itemize}
\end{proof*}

\begin{lemma}\label{lem:ts1-cb-enum-avoidance-R}
For every condition $c$ and every vector of indices~$\vec{e}$
there exists an extension $d$ of~$c$ forcing $R_{\vec{e}}$.
\end{lemma}
\begin{proof*}
Fix a condition~$c$, and apply iteratively Lemma~\ref{lem:ts1-cb-enum-avoidance-R-iter}
to obtain an extension~$d$ such that for each~$j < m$, 
$d$ forces~$R_{e_i, j}$ for $m$ different~$i$'s.
By the pigeonhole principle, there exists some~$i \leq m$
such that $d$ forces~$\bigwedge R_{e_i,j}$ for each~$j < m$.
Therefore, $d$ forces~$\Rcal_{\vec{e}}$.
\end{proof*}

Let~$\Fcal = \{c_0, c_1, \dots\}$ be a sufficiently generic filter
containing $(\emptyset, \dots, \emptyset, \omega)$, where $c_s = (F_{0,s}, \dots, \allowbreak F_{m,s}, X_s)$.
The filter~$\Fcal$ yields an $(m+1)$-tuple of reals $G_0, \dots, G_m$
defined by~$G_i = \bigcup_s F_{i,s}$.
By definition of a forcing condition, $G_i \subseteq \overline{A}_i$. 
By Lemma~\ref{lem:ts1-cb-enum-avoidance-Q}, $G_i$ is infinite for each~$i \leq m$,
and by Lemma~\ref{lem:ts1-cb-enum-avoidance-R}, $G_i \oplus C$ computes no c.b-enum of the~$\Ccal$'s
for some~$i \leq m$.
\end{proof}

\subsection{Strong c.b-enum avoidance of~$\ts^n$}

Last, we prove step (A3) of Theorem~\ref{thm:ts-rwkl-pi01-enum-avoidance}.

\begin{theorem}\label{thm:ts-enum-avoidance-step-a3}
Fix a sequence of sets~$\Ccal_0, \Ccal_1, \dots \subseteq 2^\omega$.
Fix some~$n \geq 1$ and suppose that for each~$t \in (0,n)$, $\ts^t_{d_t+1}$ admits strong 1-enum avoidance for $\vec{\Ccal}$
and that $\Pi^0_1(\ts^n_{\tilde{d}_n+1})$ admits 1-enum avoidance for~$\vec{\Ccal}$.
Then $\ts^n_{d_n+1}$ admits strong 1-enum avoidance for $\vec{\Ccal}$.
\end{theorem}
\begin{proof}
The proof is again done in several steps and is very similar to the proof
of preservation of hyperimmunity. Fix a countable collection of sets~$\Ccal_0, \Ccal_1, \dots \subseteq 2^\omega$
and assume that for each~$t \in (0,n)$, $\ts^t_{d_t+1}$ admits strong 1-enum avoidance for $\vec{\Ccal}$
and that~$\Pi^0_1(\ts^n_{\tilde{d}_n+1})$ admits 1-enum avoidance for~$\vec{\Ccal}$.
Fix a coloring $f : [\omega]^n \to d_n+1$ and a set $C$ computing no 1-enum of the $\Ccal$'s.
\bigskip

\begin{itemize}
	\item[(S1)] First, we construct an infinite set $D \subseteq \omega$ such that $D \oplus C$ computes no 1-enum of the $\Ccal$'s
	and a sequence $(I_\sigma : 0 < |\sigma| < n)$ such that for each $t \in (0,n)$ and each $\sigma \in [\omega]^t$
	\begin{itemize}
		\item[(a)] $I_\sigma$ is a subset of $\{0, \dots, d_n\}$ with at most $d_{n-t}$ many elements  
		\item[(b)] $(\exists b)(\forall \tau \in [D \cap (b,+\infty)]^{n-t}) f(\sigma,\tau) \in I_\sigma$
	\end{itemize}

	\item[(S2)]
	Second, we construct an infinite set~$E \subseteq D$ such that $E \oplus C$ computes no 1-enum of the~$\Ccal$'s
	and a sequence $(I_t : 0 < t < n)$ such that for each $t \in (0,n)$
	\begin{itemize}
		\item[(a)] $I_t$ is a subset of~$\{0,\dots, d_n\}$ of size at most $d_t d_{n-t}$
		\item[(b)] $(\forall \sigma \in [E]^t)(\exists b)(\forall \tau \in [E \cap (b,+\infty)]^{n-t}) f(\sigma,\tau) \in I_t$
	\end{itemize}

	\item[(S3)] Third, we construct a sequence~$(\xi_i \in [E]^{<\omega} : i < \omega)$ such that
	\begin{itemize}
		\item[(a)] The set $G = \bigcup_i \xi_i$ is infinite and $G \oplus C$ computes no 1-enum of the~$\Ccal$'s
		\item[(b)] $|f([\xi_i]^n)| \leq \tilde{d}_n$ and~$max(\xi_i) < min(\xi_{i+1})$ for each~$i < \omega$
		\item[(c)] For each~$t \in (0,n)$ and~$\sigma \in [\bigcup_{j < i} \xi_j]^t$,
		$f(\sigma,\tau) \in I_t$ for all~$\tau \in [\bigcup_{j \geq i} \xi_j]^{n-t}$
	\end{itemize}

	\item[(S4)] Finally, we build an infinite set~$H \subseteq G$
	such that~$H \oplus C$ computes no 1-enum of the~$\Ccal$'s
	and~$|f([H]^n)| \leq d_n$.
\end{itemize}
\bigskip

We only show step (S3) since the other steps are exactly the same as for preservation of hyperimmunity.
Given the set~$E$ and the sequence of sets of colors~$(I_t : 0 < t < n)$,
we will construct a sequence~$(\xi_i \in [E]^{<\omega} : i < \omega)$ such that
\begin{itemize}
	\item[(a)] The set $G = \bigcup_i \xi_i$ is infinite and $G \oplus C$ computes no 1-enum of the~$\Ccal$'s
	\item[(b)] $|f([\xi_i]^n)| \leq \tilde{d}_n$ and~$max(\xi_i) < min(\xi_{i+1})$ for each~$i < \omega$
	\item[(c)] For each~$t \in (0,n)$ and~$\sigma \in [\bigcup_{j < i} \xi_j]^t$,
	$f(\sigma,\tau) \in I_t$ for all~$\tau \in [\bigcup_{j \geq i} \xi_j]^{n-t}$
\end{itemize}

We construct our set~$G$ by Mathias forcing~$(\sigma, X)$ where
$X$ is an infinite subset of~$E$ such that $X \oplus C$ computes no 1-enum of the~$\Ccal$'s. 
Using property~(b) of~$E$, we can easily
construct an infinite sequence $(\xi_i \in [E]^{<\omega} : i < \omega)$
satisfying properties~(b) and~(c) of step (S3). The following lemma
shows how to satisfy property~(a).

\begin{lemma}\label{lem:ts-lemma-set-g-enum}
Fix a condition~$(\sigma, X)$ and some~$e,i \in \omega$.
There exists an extension~$(\sigma\xi, Y)$
with~$|f([\xi]^n)| \leq \tilde{d}_n$, forcing~$\Phi_e^{G \oplus C}$ not to be a 1-enum of~$\Ccal_i$.
\end{lemma}
\begin{proof*}
Suppose there exists a string $\rho \in 2^{<\omega}$ such that $[\rho] \cap \Ccal_i = \emptyset$
and a finite set $E \subset X$ such that for every coloring $g : [X]^n \to d_n+1$, there is a set~$\xi \in [X]^{<\omega}$
such that~$|g([\xi]^n)| \leq \tilde{d}_n$ and~$\Phi_e^{\sigma\xi \oplus C}(|\rho|) \downarrow = \rho$.
Then, taking~$f = g$, the condition~$d = (\sigma\xi, X \setminus [0, max(\xi)])$ is a valid extension of~$c$
forcing~$\Phi_e^{G \oplus C}$ not to be a 1-enum of $\Ccal_i$.

So suppose there is no such $\rho \in 2^{<\omega}$.
For each $\rho \in 2^{<\omega}$, let $\Pcal_{\rho}$ denote the collection 
of the colorings $g : [X]^n \to d_n+1$ such that 
for every set~$\xi \in [X]^{<\omega}$ such that~$|g([\xi]^n)| \leq \tilde{d}_n$,
$\Phi_e^{\sigma\xi \oplus C}(\card{\rho}) \uparrow$ or $\Phi_e^{\sigma\xi \oplus C}(\card{\rho}) \neq \rho$.
Note that $\Pcal_{\rho}$ are uniformly $\Pi^{0,X \oplus C}_1$ classes.
Because the previous case does not hold, then by compactness $\Pcal_\rho \neq \emptyset$
for each $\rho$ such that $\Ccal_i \cap [\rho] = \emptyset$.
The set $\{\rho : \Pcal_\rho = \emptyset\}$ is $X \oplus C$-c.e. If for each $u \in \omega$,
there exists a $\rho \in 2^u$ such that $\Pcal_\rho = \emptyset$ then 
$X \oplus C$ computes a 1-enum of $\Ccal_i$, contradicting our hypothesis.
So there must be a $u$ such that $\Pcal_\rho \neq \emptyset$ for each $\rho \in 2^u$.

Thanks to 1-enum avoidance of $\Pi^0_1(\ts^n_{\tilde{d}_n+1})$ for $\vec{\Ccal}$, define a finite decreasing sequence 
$X = X_0 \supseteq \dots \supseteq X_{2^u-1} = \tilde{X}$ such that for each $\rho \in 2^u$
\begin{itemize}
	\item[1.] $|g([X_\rho]^n)| \leq \tilde{d}_n$ for some~$g \in \Pcal_\rho$.
	\item[2.] $X_\rho \oplus C$ computes no 1-enum of any of the $\Ccal$'s.
\end{itemize}
We claim that $(\sigma, \tilde{X})$ is an extension forcing $\Phi_e^{G \oplus C}(u) \uparrow$ or $\Phi_e^{G \oplus C}(u) \not \in 2^u$.
Suppose for the sake of contradiction that there exists a $\rho \in 2^u$ and a set $G$ satisfying $(\sigma, \tilde{X})$
such that $\Phi_e^{G \oplus C}(u) \downarrow = \rho$.
By continuity, there exists a finite $\xi \subseteq G \subseteq \tilde{X}$ such that $\Phi_e^{\sigma\xi \oplus C}(u) \downarrow = \rho$.
In particular, $\xi \subseteq X_\rho$, so $|g([\xi]^n)| \leq \tilde{d}_n$ for some~$g \in \Pcal_\rho$, contradiction.
\end{proof*}

Using Lemma~\ref{lem:ts-lemma-set-g-enum} and property (b) of the set~$E$,
we can construct an infinite descending sequence of conditions~$(\epsilon,E) \geq c_0 \geq \dots$
such that for each~$s \in \omega$
\begin{itemize}
	\item[(i)] $\sigma_{s+1} = \sigma_s\xi_s$ with $|\sigma_s| \geq s$ and ~$f([\xi_s]^n) \leq \tilde{d}_n$
	\item[(ii)] $f(\sigma,\tau) \in I_t$ for each~$t \in (0,n)$, $\sigma \in [\sigma_s]^t$ and~$\tau \in [X]^{n-t}$.
	\item[(iii)] $c_s$ forces~$\Phi_e^{G \oplus C}$ not to be a 1-enum of~$\Ccal_i$ if~$s = \tuple{e,i}$
\end{itemize}
where~$c_s = (\sigma_s, X_s)$. The set~$G = \bigcup_s \sigma_s$
satisfies the desired properties.
This finishes step (S3) and the proof of Theorem~\ref{thm:ts-enum-avoidance-step-a3}.
\end{proof}

\section{The weakness of the free set hierarchy}

Recall that given a coloring $f : [\omega]^n \to \omega$,
a set~$A$ is \emph{$f$-free} if for every~$x_1 < \dots < x_n \in A$,
if $f(x_1, \dots, x_n) \in A$ then $f(x_1, \dots, x_n) \in \set{x_1, \dots, x_n}$.
As Wang did for strong cone avoidance and we did for preservation of hyperimmunity,
we can propagate strong c.b-enum avoidance from the thin set theorem
to the free set theorem. More precisely, we prove the following stronger theorem.

\begin{theorem}\label{thm:fs-strong-enum-avoidance}
Fix a countable collection of sets $\Ccal_0, \Ccal_1, \dots \subseteq 2^\omega$.
If $\ts^s_{d_t+1}$ admits simultaneous strong 1-enum avoidance for~$\vec{\Ccal}$
for each~$s \in (0, n]$, then so does $\fs^n$.
\end{theorem}

\begin{corollary}\label{cor:fs-strong-cb-enum-avoidance}
$\fs$ admits strong $m$ c.b-enum avoidance for every~$m \geq 1$.
\end{corollary}
\begin{proof}
By strong $m$ c.b-enum avoidance of $\ts^t_{d_t+1}$ for every $t \in (0,n]$,
Theorem~\ref{thm:fs-strong-enum-avoidance} and Lemma~\ref{lem:cbenum-1-enum-bridge}.
\end{proof}

In particular, since the rainbow Ramsey theorem is a consequence of the free set theorem
in reverse mathematics, we obtain strong c.b-enum avoidance of the rainbow Ramsey theorem for free.

\begin{corollary}
$\rrt$ admits strong $m$ c.b-enum avoidance for every~$m \geq 1$.
\end{corollary}
\begin{proof}
Wang proved in~\cite{Wang2014Some} that $\rrt \leq_{sc} \fs$.
Apply Lemma~\ref{lem:strong-computably-reducible-strong-avoidance}
and Corollary~\ref{cor:fs-strong-cb-enum-avoidance}.
\end{proof}

The usual lemma about the propagation of strong avoidance 
for $n$-tuples to avoidance for $(n+1)$-tuples holds for the free set theorem as well.

\begin{lemma}[Wang \cite{Wang2014Some}]
For each $n \geq 1$, if $\fs^n$ and~$\coh$ admit strong $\Ccal$ avoidance, 
then $\fs^{n+1}$ admits $\Ccal$ avoidance.
\end{lemma}

The proof of strong 1-enum avoidance of $\fs$ relative to strong 1-enum avoidance
of $\ts$ uses a case analysis only on two kinds of functions: left trapped 
and right trapped functions. Recall that a function $f : [\omega]^n \to \omega$ is \emph{left (resp. right) trapped}
if for every $\sigma \in [\omega]^n$, $f(\sigma) \leq \sigma(n-1)$ (resp. $f(\sigma) > \sigma(n-1)$).
We have already used the notion of trapped functions for the preservation of hyperimmunity
of the free set theorem. We now prove the lemmas in their general form.

\begin{lemma}[Wang in~\cite{Wang2014Some}]\label{lem:avoidance-fs-trapped-to-untrapped}
For each $n \geq 1$, if $\fs^n$ for trapped functions admits (strong) $\Ccal$ avoidance 
for some set~$\Ccal \subseteq \omega^\omega$, then so does $\fs^n$.
\end{lemma}
\begin{proof}
We prove it in the case of strong $\Ccal$ avoidance. The proof of $\Ccal$ avoidance is similar.
Let $f : [\omega]^n \to \omega$ be a coloring and $C$ be set computing no member of $\Ccal$.
For each $\sigma \in [\omega]^n$ and $i \leq n$, let
$$
f_0(\sigma) = min(f(\sigma), max(\sigma))
\hspace{20pt}
f_1(\sigma) = max(f(\sigma), max(\sigma)+1)
$$
By strong $\Ccal$ avoidance of $\fs^n$ for trapped functions,
we can define a finite sequence $\omega \supseteq H_0 \supseteq H_1$ such that
for each $i \leq n$
\begin{itemize}
  \item[1.] $H_i$ is an infinite $f_i$-free set
  \item[2.] $H_i \oplus C$ computes no member of $\Ccal$.
\end{itemize}
We claim that $H_1$ is $f$-free. Let $\sigma \in [H_n]^n$. 
$f(\sigma) = f_i(\sigma)$ for some $i \in \{0,1\}$. As $H_1$ is free for $f_i$, 
$f(\sigma) \not \in H_1 \setminus \sigma$.
\end{proof}

\subsection{Case of right trapped functions}

We now show that the case of right trapped functions
can be reduced to the left trapped functions since they admit
solutions computable by any diagonally non-computable function.

\begin{lemma}\label{lem:enum-avoidance-fsr-r-dnc}
Let $f : [\omega]^n \to \omega$ be a right trapped function.
Every function d.n.c. relative to $f$ computes an infinite set free for $f$.
\end{lemma}
\begin{proof}
By~\cite{Kjos-Hanssen2011Kolmogorov}, every function d.n.c. relative to $f$ computes a function $g$
such that if $|W^f_e| \leq m$ then $g(e, m) \not \in W^f_e$.
Given a finite $f$-free set $F$,
there exists at most ${\card{F} \choose n}$ elements $x$ such that $F \cup \{x\}$
is not $f$-free. We can define an infinite $f$-free set $H$ by stages.
$H_0 = \emptyset$. Given a finite $f$-free set $H_s$ of cardinal $s$,
set $H_{s+1} = H_s \cup \{g(e, {s \choose n})\}$ where $e$ is a Turing index
such that $W^f_e = \{x : F \cup \{x\} \mbox{ is not } f-\mbox{free}\}$.
\end{proof}

\begin{lemma}\label{lem:left-trapped-to-untrapped}
For each $n \geq 1$, if $\fs^n$ for left trapped functions admits (strong) $\Ccal$ avoidance
for some set~$\Ccal \subseteq \omega^\omega$, then so does $\fs^n$.
\end{lemma}
\begin{proof}
Again, we prove it in the case of strong $\Ccal$ avoidance.
By Lemma~\ref{lem:avoidance-fs-trapped-to-untrapped}, it suffices to prove that 
$\fs^n$ for right trapped functions admits strong $\Ccal$ avoidance.
Let~$f : [\omega]^n \to \omega$ be a right trapped function and $C$ be a set computing no member of~$\Ccal$.
By Rice~\cite{RiceThin}, there exists an $f$-computable stable left trapped function $g$ such that every infinite set thin for $g$
computes a function d.n.c.\ relative to $f$. By \cite[Theorem 3.2]{Cholak2001Free}, every infinite set free for $g$
is, up to finite variation, a set thin for $g$. By strong $\Ccal$ avoidance of~$\fs^n$ for left trapped functions,
there is an infinite $g$-free set~$H$ such that~$H \oplus C$ computes no member of~$\Ccal$.
By Lemma~\ref{lem:enum-avoidance-fsr-r-dnc}, $H$ computes an infinite $f$-free set.
\end{proof}

\subsection{Case of left trapped functions}

By Lemma~\ref{lem:left-trapped-to-untrapped} it remains to prove strong 1-enum
avoidance of $\fs^n$ for left trapped functions, assuming
strong 1-enum avoidance of $\ts^t_d$ for each $t \leq n$ and sufficiently large $d$'s.

\begin{theorem}\label{thm:strong-enum-avoidance-fs-left-trapped}
Fix a countable collection of sets $\Ccal_0, \Ccal_1, \dots \subseteq 2^\omega$.
If $\ts^s_{d_t+1}$ admits simultaneous strong 1-enum avoidance for~$\vec{\Ccal}$
for each~$s \in (0, n]$, then so does $\fs^n$ for left trapped functions.
\end{theorem}

The proof will be by induction over $n$. The base case is easy
and follows directly from the strong 1-enum avoidance of $\ts^1_2$ for $\vec{\Ccal}$.

\begin{lemma}\label{lem:rt14-fs1-enum-avoidance}
If $\rt^1_2$ admits strong 1-enum avoidance for $\vec{\Ccal}$ then so does $\fs^1$.
\end{lemma}
\begin{proof}
Cholak et al.\ proved in~\cite{Cholak2001Free} that $\fs^1$ 
for left trapped functions is strongly computably reducible to $\rt^1_4$.
Apply Lemma~\ref{lem:strong-computably-reducible-strong-avoidance} to deduce
strong 1-enum avoidance of $\fs^1$ for left trapped functions for $\vec{\Ccal}$.
Conclude with Lemma~\ref{lem:left-trapped-to-untrapped}.
\end{proof}

The two following lemmas will ensure that the reservoirs of our forcing conditions
will enjoy some good properties which will ensure that the conditions are extendible.

\begin{lemma}\label{lem:fs2-1-trapped-avoid-enum-smaller}
Suppose that $\fs^s$ admits strong $\Ccal$ avoidance for each $s < n$ for some set $\Ccal \subseteq \omega^\omega$.
Fix a set $C$ computing no member of $\Ccal$, a finite set $F$ and an infinite set $X$ computable in $C$.
For every function $f : [X]^n \to \omega$ there exists an infinite
set $Y \subseteq X$ such that $Y \oplus C$ computes no member of $\Ccal$ and 
$(\forall \sigma \in [F]^t)(\forall \tau \in [Y]^{n-t})f(\sigma, \tau) \not \in Y \setminus \tau$
for each $0 < t < n$.
\end{lemma}
\begin{proof}
Fix the finite enumeration $\sigma_1, \dots, \sigma_k$ for all $\sigma \in [F]^t$ for some $0 < t < n$.
Start with $Y_0 = X$.
Suppose that $Y_{m-1} \oplus C$ computes no member of $\Ccal$ and
for all $i < m$, $\forall \tau \in [Y_{m-1}]^{n - |\sigma_i|}f(\sigma_i, \tau) \not \in Y_{m-1} \setminus \tau$.
Define the function $f_{\sigma_m} : [Y_{m-1}]^{n - |\sigma_m|} \to \omega$ 
by $f_{\sigma_m}(\tau) = f(\sigma_m, \tau)$.
By strong $\Ccal$ avoidance of $\fs^{n - |\sigma_m|}$,
there exists an infinite set $Y_m \subseteq Y_{m-1}$ such that $Y_m \oplus C$ computes no member of $\Ccal$
and $(\forall \tau \in [Y_m]^{n - |\sigma_m|})f(\sigma_m, \tau) \not \in Y_m \setminus \tau$.
$Y_k$ is the desired set.
\end{proof}

\begin{lemma}\label{lem:fsr-cohesivity-strong-enum-avoidance}
Suppose that $\ts^s_{d_s+1}$ admits strong 1-enum avoidance for $\vec{\Ccal}$ and for each $0 < s \leq n$
and $\fs^s$ admits strong 1-enum avoidance for $\vec{\Ccal}$ and for each $0 < s < n$.
For every function $f : [\omega]^n \to \omega$
and every set $C$ computing no 1-enum of $\vec{\Ccal}$, there exists an infinite set $X$
such that $X \oplus C$ computes no 1-enum of $\vec{\Ccal}$ and 
for every $\sigma \in [G]^{<\omega}$ such that $0 \leq \card{\sigma} < n$,
$$
(\forall x \in G \setminus \sigma)(\exists b)(\forall \tau \in [G \cap (b, +\infty)]^{n-|\sigma|}) 
  f(\sigma, \tau) \neq x
$$
\end{lemma}
\begin{proof}
Let $X$ be an infinite set satisfying property of Theorem~\ref{thm:generalized-cohesivity-strong-enum-avoidance}
with $t = n$. For each $s < n$ and $i < d_{n-s}$, let $f_{s,i} : [X]^s \to \omega$ be the function
such that $f_{s,i}(\sigma)$ is the $i$th element of the set
$$
\{x : (\forall b)(\exists \tau \in [X \cap (b, +\infty)]^{n-s}) f(\sigma, \tau) = x\}
$$
if it exists, and 0 otherwise.
Define a finite sequence $X \supseteq X_0 \supseteq \dots \supseteq X_{n-1}$ such that for each $s < n$
\begin{itemize}
  \item[1.] $X_s$ is $f_{s,i}$-free for each $i < d_{n-s}$
  \item[2.] $X_s \oplus C$ computes no 1-enum of $\vec{\Ccal}$
\end{itemize}
We claim that $X_{n-1}$ is the desired set. Fix $s < n$ and take any $\sigma \in [X_{n-1}]^{s}$
and any $x \in X_{n-1} \setminus \sigma$. 
If $(\forall b)(\exists \tau \in [G \cap (b, +\infty)]^{n - s})f(\sigma, \tau) = x$,
then by choice of $X$, there exists an $i < d_{n - s}$ such that
$f_{s, i}(\sigma) = x$, contradicting $f_{s,i}$-freeness of $X_{n-1}$.
So $(\exists b)(\forall \tau \in [G \cap (b, +\infty)]^{n - s})f(\sigma, \tau) \neq x$.
\end{proof}

We are now ready to prove Theorem~\ref{thm:strong-enum-avoidance-fs-left-trapped}.

\begin{proof}[Proof of Theorem~\ref{thm:strong-enum-avoidance-fs-left-trapped}]
Fix a countable collection of sets $\Ccal_0, \Ccal_1, \dots \subseteq 2^\omega$
for which $\ts^s_{d_s+1}$ admits strong 1-enum avoidance for each $0 < s \leq n$.
Let $f : [\omega]^n \to \omega$ be a left trapped function and $C$ be a set computing no 1-enum of $\vec{\Ccal}$.
Our forcing conditions are tuples $(k, \vec{F}, X, \vec{g})$ such that
\begin{itemize}
  \item[(a)] $\vec{g}$ is a left trapped $\oplus_k$-function, $\vec{F}$ is a finite $\oplus_k$-set
  \item[(b)] $X$ is an infinite set such that $F_0 < X$ and $X \oplus C$ computes no 1-enum of $\vec{\Ccal}$
  \item[(c)] $(\forall \sigma \in [F_i \cup X]^n) g_i(\sigma) \not \in F_i \setminus \sigma$ for each $i < k$
  \item[(d)] $(\forall \sigma \in [F_i \cup X]^t)(\forall x \in (F_i \cup X) \setminus \sigma)
  (\exists b)(\forall \tau \in [(F_i \cup X) \cap (b, +\infty)]^{n-t})\\
  g_i(\sigma, \tau) \neq x$ for each $i < k$ and $0 \leq t < n$.
  \item[(e)] $(\forall \sigma \in [F_i]^t)(\forall \tau \in [X]^{n-t})
  g_i(\sigma, \tau) \not \in X \setminus \tau$ for each $i < k$ and $0 < t < n$
\end{itemize}
Properties (d) and (e) will be obtained by 
Lemma~\ref{lem:fsr-cohesivity-strong-enum-avoidance} 
and Lemma~\ref{lem:fs2-1-trapped-avoid-enum-smaller} and are present to ensure to have
extensions such that (c) holds.
A set $G$ \emph{satisfies} a condition $(k, \vec{F}, X, \vec{g})$ if it
satisfies the Mathias condition $(F_0, X)$ and $G \setminus (F_0 \setminus F_i)$
if $g_i$-free for each $i < k$.
Our initial condition is $(1, \emptyset, Y, f)$ where $Y$ is obtained by Lemma~\ref{lem:fsr-cohesivity-strong-enum-avoidance}.
A condition $(m, \vec{F}', X', \vec{g}')$ \emph{extends} another condition $(k, \vec{F}, X, \vec{g})$
if $X' \subseteq X$,  $m \geq k$, $(\forall i < k)g_i = g'_i$ and there is a finite $E \subset X$ such that
\begin{itemize}
  \item[(i)] for every $i < k$, $F_i \subseteq F'_i$ and $F'_i \setminus F_i = E$
  \item[(ii)] for every $k \leq i < m$, $F'_i = E$
\end{itemize}

\begin{lemma}\label{lem:fs2-left-trapped-avoid-enum-1}
For every condition $(k, \vec{F}, X, \vec{g})$ there exists an extension
$(k, \vec{H}, \tilde{X}, \vec{g})$ such that $|H_i| > \card{F_i}$ for each $i < k$.
\end{lemma}
\begin{proof*}
Choose an $x \in X$ such that 
$(\forall j < k)(\forall \sigma \in [F_j]^n)g_j(\sigma) \neq x$ and set $H_i = F_i \cup \{x\}$ for each $i < k$.
By property (d) of $(k, \vec{F}, X, \vec{g})$, there exists a $b$
such that $(\forall i < k)(\forall \sigma \in [F_i]^{t})
(\forall \tau \in [X \cap (b, +\infty)]^{n-t})g_i(\sigma, \tau) \neq \{x\} \setminus \sigma$
for each $0 \leq t \leq n$.
By $k$ applications of Lemma~\ref{lem:fs2-1-trapped-avoid-enum-smaller},
there exists a $\tilde{X} \subseteq X \setminus [0, b]$ such that $\tilde{X} \oplus C$ computes no 1-enum of $\vec{\Ccal}$
and property (e) is satisfied for $(k, \vec{H}, \tilde{X}, \vec{g})$.
We claim that $(k, \vec{H}, \tilde{X}, \vec{g})$ is a valid condition.
Properties (a), (b) and (d) trivially hold. It remains to check property (c).
By property (c) of $(k, \vec{F}, X, \vec{g})$, we only need to check that
$(\forall \sigma \in [F_i \cup \tilde{X}]^n)g_i(\sigma) \neq x$ for each $i < k$.
This follows from our choice of~$b$.
\end{proof*}

\begin{lemma}\label{lem:fs2-left-trapped-avoid-enum-2}
For every condition $(k, \vec{F}, X, \vec{g})$ and every $e,i \in \omega$,
there exists an extension $(m, \vec{H}, \tilde{X}, \vec{h})$
forcing $\Phi_e^{G \oplus C}$ not to be a 1-enum of $\Ccal_i$, where $G$ is the forcing variable. 
\end{lemma}
\begin{proof*}
By removing finitely many elements to $X$, we can suppose w.l.o.g. that
$(\forall j < k)(\forall \sigma \in [F_j]^n)g_j(\sigma) \not \in X$.
Suppose there exists a $\sigma \in 2^{<\omega}$ such that $[\sigma] \cap \Ccal_i = \emptyset$
and a finite set $F' \subseteq X$ which is $g_j$-free for each $j < k$
and  $\Phi_e^{(F_0 \cup F') \oplus C}(|\sigma|) \downarrow = \sigma$.
Set $H_j = F_j \cup F'$ for each $j < k$.
By property (d) of $(k, \vec{F}, X, \vec{g})$, there exists a $b$ such that
$(\forall \sigma \in [H_i]^t)(\forall x \in H_i)
(\forall \tau \in [X \cap (b, +\infty)]^{n-t}) 
  g_i(\sigma, \tau) \neq \{x\} \setminus \sigma$ for each $i < k$ and $0 \leq t < n$.
By $k$ applications of Lemma~\ref{lem:fs2-1-trapped-avoid-enum-smaller},
there exists a $\tilde{X} \subseteq X \cap (b, +\infty)$ such that $\tilde{X} \oplus C$ computes no 1-enum of $\vec{\Ccal}$
and property (e) is satisfied for $(k, \vec{H}, \tilde{X}, \vec{g})$.
We claim that $(k, \vec{H}, \tilde{X}, \vec{g})$ is a valid condition.

Properties (a), (b), (d) and (e) trivially hold. It remains to check property (c).
By our choice of $b$, we need only to check that $(\forall \sigma \in [H_i]^n)(\forall x \in H_i)
g_i(\sigma) \neq \{x\} \setminus \sigma$ for each $i < k$.
By property (c) of $(k, \vec{F}, X, \vec{g})$, it suffices to check that
$(\forall \sigma \in [H_i]^n)g_i(\sigma) \not \in F' \setminus \sigma$ for each $i < k$.
By property (e) of $(k, \vec{F}, X, \vec{g})$, it remains the case
$(\forall \sigma \in [F']^n)g_i(\sigma) \not \in F' \setminus \sigma$ for each $i < k$,
which is exactly $\vec{g}$-freeness of $F'$.

Suppose there is no such finite set $F' \subset X$.
For each $\sigma \in 2^{<\omega}$, let $\Fcal_{\sigma}$ denote the collection 
of left trapped $\oplus_k$-functions $\vec{g}$ such that
for each finite set $F' \subset X$ which is $g_j$-free for each $j < k$,
either $\Phi_e^{(F_0 \cup F') \oplus C}(\card{\sigma}) \uparrow$ or $\Phi_e^{(F_0 \cup F') \oplus C}(\card{\sigma}) \neq \sigma$.
Note that $\Fcal_\sigma$ are uniformly $\Pi^{0,X \oplus C}_1$ classes.
Because the above case does not hold, $\vec{g} \in \Fcal_\sigma$ for each $\sigma$ such that $\Ccal_i \cap [\sigma] = \emptyset$.
The set $\{\sigma : \Fcal_\sigma = \emptyset\}$ is $X \oplus C$-c.e. If for each $u \in \omega$
there exists a $\sigma \in 2^u$ such that $\Fcal_\sigma = \emptyset$ then $X \oplus C$ computes a 1-enum of $\Ccal_i$,
contradicting our hypothesis.
So there must be an $u \in \omega$ such that $\Fcal_\sigma \neq \emptyset$ for each $\sigma \in 2^u$.

For each $\sigma \in 2^u$, let $\vec{h}_\sigma \in \Fcal_\sigma$.
Set $H_j = F_j$ for each $j < k$ and $H_j = \emptyset$ for each $k \leq j < (2^u+1)k$.
By $2^u$ applications of Lemma~\ref{lem:fsr-cohesivity-strong-enum-avoidance},
there exists an infinite set $\tilde{X} \subseteq X$
such that $\tilde{X} \oplus C$ computes no 1-enum of $\vec{\Ccal}$ and property (d)
of $((2^u+1)k, \vec{H}, \tilde{X}, \vec{g} \bigoplus_{\sigma \in 2^u} \vec{h}_\sigma)$ holds.
As properties (a-c) and (e) trivially hold, $((2^u+1)k, \vec{H}, \tilde{X}, \vec{g} \bigoplus_{\sigma \in 2^u} \vec{h}_\sigma)$
is a valid condition. Moreover it forces $\Phi_e^{G \oplus C}(u) \uparrow$ or $\Phi_e^{G \oplus C}(u) \downarrow \not \in 2^u$.
\end{proof*}

Let~$\Fcal = \{c_0, c_1, \dots\}$ be a sufficiently generic filter
containing $(1, \emptyset, Y, f)$, where $c_s = (k_s, \vec{F}_s, \vec{g}_s)$.
The filter~$\Fcal$ yields a unique real~$G = \bigcup_s F_{s,0}$.
By definition of a forcing condition, $G$ is an $f$-free set. 
By Lemma~\ref{lem:fs2-left-trapped-avoid-enum-1}, $G$ is infinite,
and by Lemma~\ref{lem:fs2-left-trapped-avoid-enum-2}, $G \oplus C$ computes no 1-enum of $\vec{\Ccal}$.
\end{proof}

\chapter{Strengthenings of Ramsey's theorem}

The last three sections of this chapter are a joint work with Emanuele Frittaion.

Thanks to Seetapun's theorem~\cite{Seetapun1995strength}, we know that Ramsey's theorem for pairs is strictly
weaker than the arithmetic comprehension axiom in reverse mathematics.
Since then, many consequences of Ramsey's theorem have been proven to be strictly weaker than~$\rt^2_2$,
namely, the Erd\H{o}s-Moser theorem, the ascending descending sequence, the thin set theorem for pairs...
However, no natural principle is currently known to be strictly between~$\aca$ and~$\rt^2_2$.
There are two good candidates: a generalization of Ramsey's theorem to colorings
over trees~\cite{Chubb2009Reverse}, and a theorem about partitions of rationals 
due to Erd\H{o}s and Rado~\cite{Erdos1952Combinatorial}.

The tree theorem is a strengthening of Ramsey's theorem in which we do not consider
colorings over tuples of integers, but colorings over tuples of nodes over a binary tree.
Ramsey's theorem can be recovered from the tree theorem by identifying
all nodes at every given level of the tree.

\index{tree theorem}
\index{tt@$\tto^n_k$|see {tree theorem}}
\index{homogeneous!for a coloring (tree)}
\begin{definition}[Tree theorem]
We denote by $[2^{<\omega}]^n$ the collection of \emph{linearly ordered} subsets of~$2^{<\omega}$ of size~$n$.
A subtree $S \subseteq 2^{<\omega}$ is \emph{order isomorphic} to $2^{<\omega}$ (written $S \cong 2^{<\omega}$) 
if there is a bijection $g : 2^{<\omega} \to S$ such that for all $\sigma, \tau \in 2^{<\omega}$,
$\sigma \preceq \tau$ if and only if $g(\sigma) \preceq g(\tau)$.
Given a coloring $f : [2^{<\omega}]^n \to k$, a tree $S$ is $f$-homogeneous if $S \cong 2^{<\omega}$
and $f \uh [S]^n$ is monochromatic.
$\tto^n_k$ is the statement ``Every coloring $f : [2^{<\omega}]^n \to k$ has an $f$-homogeneous tree.''
\end{definition}

Note that if $S \cong 2^{<\omega}$, witnessed by the bijection $g : 2^{<\omega} \to S$,
then $S$ is $g$-computable. Therefore we can consider that $\tto^n$ states the existence of the bijection $g$
instead of the pair~$\tuple{S, g}$.
The tree theorem was first analyzed by McNicholl~\cite{McNicholl1995inclusion}
and by Chubb, Hirst, and McNicholl~\cite{Chubb2009Reverse}. They proved that~$\tto^2_2$
lies between $\aca$ and~$\rt^2_2$ over~$\rca$, and left open whether any of the implications
is strict. Further work was done by Corduan, Groszek, and Mileti~\cite{Corduan2010Reverse}.
Dzhafarov, Hirst and Lakins~\cite{Dzhafarov2010Ramseys} studied stability notions for the
tree theorem and introduced a polarized variant.

Mont\'alban~\cite{Montalban2011Open} asked whether
$\rt^2_2$ implies~$\tto^2_2$ over~$\rca$. We give a negative answer by proving the following stronger theorem.

\begin{theorem}\label{thm:tree-theorem-rt22+wkl-not-tt22}
$\rt^2_2 \wedge \wkl$ does not imply $\tto^2_2$ over~$\rca$.
\end{theorem}

The next three sections are dedicated to a proof of Theorem~\ref{thm:tree-theorem-rt22+wkl-not-tt22}.
We introduce in section~\ref{sect:tree-theorem-partitions-reducibility} the main ideas of the separation,
then we design in section~\ref{sect:fairness-property-for-trees} a weakness property based on the combinatorics of the previous section,
and prove the main preservations in section~\ref{sect:tree-theorem-separating-principles}.
The last sections are devoted to the Erd\H{o}s-Rado theorem.

\section{The tree theorem and strong reducibility}\label{sect:tree-theorem-partitions-reducibility}

In order to get progressively into the framework
used to separate Ramsey's theorem for pairs from the tree theorem for pairs,
we shall first study the singleton version of the considered principles.
Ramsey's theorem for singletons is simply the infinite pigeonhole principle,
stating that for every finite partition of an infinite set,
one of its parts has an infinite subset.
Both $\rt^1_k$ and~$\tto^1_k$ are effective and provable over~$\rca$.
We shall therefore study non-effective instances of~$\rt^1_k$ and~$\tto^1_k$
to see how their combinatorics differ.
The remainder of this section will be dedicated to proving that~$\tto^1_2 \not \leq_{sc} \rt^1_2$.
More precisely, we shall prove the following stronger theorem.

\begin{theorem}\label{thm:tree-theorem-rt12-comp-not-reduc}
There exists a $\Delta^0_2$ $\tto^1_2$-instance $A_0 \cup A_1 = 2^{<\omega}$
such that for every (non-effective) $\rt^1_2$-instance $B_0 \cup B_1 = \omega$,
there is an infinite set homogeneous for the~$B$'s which does not compute
a $\tto^1_2$-solution to the~$A$'s.
\end{theorem}

In section~\ref{sect:tree-theorem-separating-principles}, we will prove a theorem
which implies Theorem~\ref{thm:tree-theorem-rt12-comp-not-reduc}. Therefore we shall focus on the key ideas
of the construction rather than on the technical details.

\smallskip
\emph{Requirements}. Let us first assume that we have constructed our $\tto^1_2$ instance~$A_0 \cup A_1 = 2^{<\omega}$.
Fix some 2-partition~$B_0 \cup B_1 = \omega$. We will construct by forcing an infinite set~$G$
such that both $G \cap B_0$ and~$G \cap B_1$ are infinite.
Let~$\Phi_0, \Phi_1, \dots$ be an enumeration of all partial tree functionals isomorphic to $2^{<\omega}$, that is,
if $\Phi^X(n)$ halts, then $\Phi^X(n)$ outputs $2^n$ pairwise incomparable strings representing the $n$th level
of the tree. We require that the following formula holds
for every pair of indices~$e_0,e_1$.
\[
  \Qcal_{e_0, e_1} : \hspace{20pt} 
		\Rcal_{e_0}^{G \cap B_0} \hspace{20pt} \vee \hspace{20pt}  \Rcal_{e_1}^{G \cap B_1}
\]
where $\Rcal_e^H$ is the statement
\begin{quote}
Either~$\Phi_e^H$ is partial, or~$\Phi_e^H(n)$ is a set~$D$ of $2^n$ incomparable
strings intersecting both~$A_0$ and~$A_1$ for some~$n$.
\end{quote}
If every~$\Qcal$-requirement is satisfied, then by the usual pairing argument, either every~$\Rcal$-requirement
is satisfied for~$G \cap B_0$, or every $\Rcal$-requirement is satisfied for~$G \cap B_1$. 
Call such a set~$H$. Suppose that $H$ computes a tree~$S \cong 2^{<\omega}$ using some procedure~$\Phi_e$.
By the requirement $\Rcal_e^H$, $S$ intersects both~$A_0$ and~$A_1$, and therefore $S$ is not a $\tto^1_2$-solution to the~$A$'s.

\smallskip
\emph{Forcing}. The forcing conditions are Mathias conditions~$(F, X)$ where~$X$ belongs to some fixed \emph{Scott set} $\Scal$,
that is, a Turing ideal satisfying weak König's lemma. By Simpson~\cite[Theorem~VIII.2.17]{Simpson2009Subsystems}, we can
choose~$\Scal$ so that $\Scal = \{ X_i : C = \bigoplus_i X_i \}$ for some low set $C$. This precision will be useful
during the construction of the~$\tto^1_2$-instance. We furthermore assume that~$C$ does not compute a solution to the~$A$'s,
and therefore that there is no $C$-computable infinite set homogeneous for the~$B$'s, otherwise we are done.

The following lemma ensures that we can force both~$G \cap B_0$ and~$G \cap B_1$ to be infinite,
assuming that the $B$'s have no infinite $C$-computable homogeneous set.

\begin{lemma}\label{lem:tree-theorem-rt12-comp-reduc-force-N}
Given a condition~$c = (F, X)$ and some side~$i < 2$, 
there is an extension~$d = (E, Y)$ such that~$|E \cap B_i| > |F \cap B_i|$.
\end{lemma}
\begin{proof}
If $X \cap B_i = \emptyset$ then $X$ is an infinite $C$-computable subset of~$B_{1-i}$, contradicting
our assumption. So there is some~$x \in X \cap B_i$. Take~$d = (F \cup \{x\}, X \setminus [0,x])$
as the desired extension.
\end{proof}

The next step consists in forcing the~$\Qcal$-requirements to be satisfied.
A condition~$c$ \emph{forces} a requirement~$\Qcal_{e_0,e_1}$ if $\Qcal_{e_0,e_1}$
holds for every set~$G$ satisfying~$c$. By choosing our~$\tto^1_2$-instance~$A_0 \cup A_1 = 2^{<\omega}$ carefully,
we claim that we can ensure the following property.

\begin{quote}(P) 
Given a condition~$c = (F, X)$ and some indices~$e_0,e_1$,
there is an extension~$d$ of~$c$ forcing~$\Qcal_{e_0,e_1}$.
\end{quote}

We will first assume that the property (P) holds, and show how to build
our infinite set~$G$ from it. We will construct later a $\tto^1_2$-instance~$A_0 \cup A_1 = 2^{<\omega}$
so that the property (P) is satisfied.

\smallskip
\emph{Construction}.
Thanks to Lemma~\ref{lem:tree-theorem-rt12-comp-reduc-force-N} and the property (P), 
we can define an infinite, decreasing sequence of conditions~$(\emptyset, \omega) \geq c_0 \geq c_1 \dots$
such that for each~$s \in \omega$
\begin{itemize}
	\item[(i)] $|F_s \cap B_0| \geq s$ and~$|F_s \cap B_1| \geq s$
	\item[(ii)] $c_s$ forces~$\Qcal_{e_0,e_1}$ if~$s = \tuple{e_0,e_1}$
\end{itemize}
where~$c_s = (F_s, X_s)$. The set~$G = \bigcup_s F_s$ is such that
both $G \cap B_0$ and~$G \cap B_1$ are infinite by (i),
and either~$G \cap B_0$ or~$G \cap B_1$ does not compute a $\tto^1_2$-solution
to the~$A$'s by~(ii). We now need to satisfy the property (P).

\smallskip
\emph{Satisfying the property (P)}.
Given a condition, the extension stated in the property (P)
cannot be ensured for an arbitrary $\tto^1_2$-instance~$A_0 \cup A_1 = 2^{<\omega}$.
We must design the $A$'s so that the property (P) holds.
To do so, we will apply the ideas developped by Lerman, Solomon and Towsner~\cite{Lerman2013Separating}.
We can see the construction of the set~$G$ as a game. The opponent is the~$\tto^1_2$-instance which
will try everything, not to be diagonalized against. However, the opponent is \emph{$\tto$-fair}, in the sense that
if we have infinitely many occasions to diagonalize against him, then he will let us do it.
More precisely, if given a condition~$c = (F, X)$ and some indices~$e_0,e_1$, we can make both~$\Phi^{G \cap B_0}_{e_0}$
and~$\Phi^{G \cap B_1}_{e_1}$ produce arbitrarily large outputs, then one of those output will intersect both~$A_0$ and~$A_1$.

We now describe how to construct a $\tto$-fair $\tto^1_2$-instance. The construction of the~$A$'s
will be $\Delta^{0,C}_2$, hence~$\Delta^0_2$ since~$C$ is low. The access to the oracle~$C$ enables us to
\emph{code} the conditions $c = (F, X)$ into finite objects, namely, pairs~$(F, i)$ so that~$C = \bigoplus_i X_i$ and $X = X_i$,
and to enumerate them $C$-effectively. More precisely, we can enumerate \emph{preconditions}
since it requires too much computational power to decide whether~$X$ is infinite or not.
Given a precondition~$c = (F, X)$, we can enumerate all possible guesses of~$F \cap B_0$ and~$F \cap B_1$
by considering each 2-partition of~$F$. Last, we can enumerate all pairs of indices~$e_0, e_1$.

The construction of the~$A$'s is done by stages.
At stage~$s$, we have constructed two sets~$A_{0,s} \cup A_{1,s} = 2^{<q}$ for some~$q \in \omega$.
We want to satisfy the property (P) 
given a precondition~$c = (F, X)$, a guess of $F \cap B_0$ and $F \cap B_1$, and a pair of indices~$e_0, e_1$.
If any of $\Phi^{F \cap B_0}_{e_0}(2)$ and $\Phi^{F \cap B_1}_{e_1}(2)$ is not defined, do nothing and go to the next stage.
We can restrict ourselves without loss of generality to preconditions such that both $\Phi^{F \cap B_0}_{e_0}(2)$
and $\Phi^{F \cap B_1}_{e_1}(2)$ are defined. Indeed, if in the property (P),
the condition~$c$ has no such extension, then~$c$ already forces
either $\Phi^{G \cap B_0}_{e_0}$ or $\Phi^{G \cap B_1}_{e_1}$ to be partial and therefore
vacuoulsy forces~$\Qcal_{e_0,e_1}$. The choice of ``2'' as input seems arbitrary. It has not been picked
randomly and this choice will be justified in the next paragraph.

Let~$D_0$ and $D_1$ be the 4-sets of pairwise incomparable strings outputted by 
$\Phi^{F \cap B_0}_{e_0}(2)$ and $\Phi^{F \cap B_1}_{e_1}(2)$, respectively.
Altough the strings are pairwise incomparable \emph{within} $D_0$ or $D_1$,
there may be two comparable strings in $D_0 \cup D_1$. However,
by a simple combinatorial argument, we may always find two strings~$\sigma_0, \tau_0 \in D_0$ and
$\sigma_1, \tau_1 \in D_1$ such that $\sigma_0, \tau_0, \sigma_1$ and $\tau_1$
are pairwise incomparable (see Lemma~\ref{lem:tree-theorem-matrix-combi}). Here, we use the fact that
on input 2, the sets have cardinality 4, which is enough to apply Lemma~\ref{lem:tree-theorem-matrix-combi}.
We are now ready to ask the main question.

\smallskip
``Is it true that for every $2$-partition~$Z_0 \cup Z_1 = X$, there is some side~$i < 2$
and some set~$G \subseteq Z_i$ such that $\Phi^{(F \cap B_i) \cup G}_{e_i}(q)$ halts?''
\smallskip

Note that the question looks $\Pi^{1,X}_2$, but is in fact $\Sigma^{0,X}_1$ by the usual compactness argument.
It is therefore $C'$-decidable since~$X \in \Scal$ and so can be uniformly decided during the construction.
We have two cases.

Case 1: The answer is negative.
In this case, the $\Pi^{0,X}_1$ class $\Ccal$ of all sets $Z_0 \oplus Z_1$
such that $Z_0 \cup Z_1 = X$ and for every~$i < 2$ and every set~$G \subseteq Z_i$,
$\Phi^{(F \cap B_i) \cup G}_{e_i}(q) \uparrow$ is non-empty. In this case,
we do nothing and claim that the property (P) holds for~$c$.
Indeed, since~$\Scal$ is a Scott set containing $X$, there is some~$Z_0 \oplus Z_1 \in \Ccal \cap \Scal$
such that $Z_0 \cup Z_1 = X$. As $X$ is infinite, there is some~$i < 2$ such that~$Z_i$ is infinite.
In this case, $d = (F, Z_i)$ is an extension forcing~$\Phi^{(G \cap B_i)}_{e_i}(q) \uparrow$
and therefore forcing~$\Qcal_{e_0,e_1}$. Note that this extension
cannot be found $C'$-effectively since it requires to decide which of~$Z_0$ and~$Z_1$ is infinite.
However, we do not need to uniformly provide this extension. The property (P) 
simply states the \emph{existence} of such an extension.

Case 2: The answer is positive. Given a string~$\sigma \in 2^{<\omega}$,
let~$S_\sigma = \{\tau \succeq \sigma \}$. 
Since the $\Phi$'s are tree functionals and $\Phi^{(F \cap B_i) \cup G}_{e_i}(2)$
outputs (among others) the strings~$\sigma_i$ and~$\tau_i$, if $\Phi^{(F \cap B_i) \cup G}_{e_i}(q)$
halts, then it outputs a finite set $D$ of size $2^q$ intersecting both~$S_{\sigma_i}$ and~$S_{\tau_i}$.
Therefore, by compactness, there are finite sets $U_0 \subseteq S_{\sigma_0}$,
$V_0 \subseteq S_{\tau_0}$,  $U_1 \subseteq S_{\sigma_1}$ and~$V_1 \subseteq S_{\tau_1}$
such that for every $2$-partition~$Z_0 \cup Z_1 = X$, there is some side~$i < 2$
and some set~$G \subseteq Z_i$ such that $\Phi^{(F \cap B_i) \cup G}_{e_i}(q)$ intersects both $U_i$ and~$V_i$.
In particular, taking $Z_0 = X \cap B_0$ and~$Z_1 = X \cap B_1$,
there is some $i < 2$ and some finite set~$G \subseteq X \cap B_i$ such that 
$\Phi^{(F \cap B_i) \cup G}_{e_i}(q)$ intersects both $U_i$ and~$V_i$.
Notice that all the strings in $U_i$ and $V_i$ have length at least~$q$
and therefore are not yet colored by the~$A$'s.
Put the~$U$'s in~$A_{0,s+1}$ and the~$V$'s in~$A_{1,s+1}$ and complete the coloring
so that~$A_{0,s+1} \cup A_{1, s+1} = 2^{<r}$ for some~$r \geq q$. Then go to the next step.
We claim that the property (P) holds for~$c$.
Indeed, let~$E = F \cup G$ where~$G \subseteq X \cap B_i$ is the finite set witnessing 
$\Phi^{(F \cap B_i) \cup G}_{e_i}(q) \cap U_i \neq \emptyset$ and $\Phi^{(F \cap B_i) \cup G}_{e_i}(q) \cap V_i \neq \emptyset$.
The condition~$d = (E, X \setminus [0, max(E)])$ is an extension forcing~$\Qcal_{e_0,e_1}$
by its $i$th side.
This finishes the construction of the~$\tto^1_2$-instance and the proof of Theorem~\ref{thm:tree-theorem-rt12-comp-not-reduc}.

\section{A fairness property for trees}\label{sect:fairness-property-for-trees}

In this section, we analyse the one-step separation proof
of section~\ref{sect:tree-theorem-partitions-reducibility} in order to extract the core of the argument.
Then, we use the framework of Lerman, Solomon and Towsner
to design the computability-theoretic property
which will enable us to discriminate~$\rt^2_2$ from~$\tto^2_2$.

\smallskip
\emph{The multiple-step case}.
We have seen how to diagonalize against one application of~$\rt^1_2$.
The strength of~$\tto^1_2$ comes from the fact that 
when we build a solution $S$ to some $\tto^1_2$-instance $A_0 \cup A_1 = 2^{<\omega}$,
we must provide finite subtrees $S_n \cong 2^{<n}$ for arbitrarily large~$n$.
However, as soon as we have outputted $S_n$, we \emph{commit}
to provide arbitrarily large extensions to each leaf of~$S_n$. Since the
leaves in $S_n$ are pairwise incomparable, the sets of their extensions
are mutually disjoint. During the construction of the $\tto^1_2$-instance,
we can pick any pair $\sigma, \tau$ of incomparable leaves
in~$S_n$, and put the extensions of~$\sigma$ in $A_0$ and the extensions of~$\tau$ in~$A_1$
since they are disjoint.

In the proof of $\tto^1_2 \not \leq_{sc} \rt^1_2$,
when we create a solution to some~$\rt^1_2$-instance~$B_0 \cup B_1 = \omega$,
we build two candidate solutions~$G \cap B_0$ and~$G \cap B_1$ at the same time.
For each pair of tree functionals $\Phi_{e_0}$ and~$\Phi_{e_1}$, 
we must prevent one of $\Phi^{G \cap B_0}_{e_0}$ and~$\Phi^{G \cap B_1}_{e_1}$
from being a $\tto^1_2$-solution to the~$A$'s. However, the finite
subtrees $S_0$ and~$S_1$ outputted respectively by the left side and the right side
may have comparable leaves. We cannot take any 2 leaves of~$S_0$ and 2 leaves of~$S_1$
to obtain 4 pairwise incomparable strings. Thankfully, if $S_0$ and $S_1$ contain enough leaves (4 is enough),
we can find such strings.

If we try to diagonalize against two applications of~$\rt^1_2$,
below each side $G \cap B_0$ and~$G \cap B_1$ of the first $\rt^1_2$-instance,
we will have again two sides corresponding to the second $\rt^1_2$-instance.
We will then have to diagonalize against four candidate subtrees $S_0$, $S_1$, $S_2$ and~$S_3$.
We need therefore to wait until the subtrees have enough leaves, so that we can find 8 pairwise
incomparable leaves $\sigma_0, \tau_0, \dots, \sigma_3, \tau_3$
such that~$\sigma_i, \tau_i \in S_i$ for each~$i < 4$.

In the general case, we will then have to diagonalize against an arbitrarily large number of subtrees,
and want to ensure that if they contain enough leaves,
we can find two leaves in each, such that they form a set of pairwise incomparable strings.
This leads to the notion of disjoint matrix.

\index{disjoint matrix}
\begin{definition}[Disjoint matrix]
An $m$-by-$n$ \emph{matrix} $M$ is a rectangular array of strings $\sigma_{i,j} \in 2^{<\omega}$
such that $i < m$ and $j < n$. The $i$th \emph{row} $M(i)$ of the matrix $M$
is the $n$-tuple of strings $\sigma_{i,0}, \dots, \sigma_{i,n-1}$.
An $m$-by-$n$ matrix $M$ is \emph{disjoint} if for each row $i < m$,
the strings $\sigma_{i,0}, \dots, \sigma_{i,n-1}$ are pairwise incomparable.
\end{definition}

The following combinatorial lemma gives an explicit bound on the number of leaves 
we require on each subtree to obtain our desired sequence of
pairwise incomparable strings.

\begin{lemma}\label{lem:tree-theorem-matrix-combi}
For every $m$-by-$2m$ disjoint matrix $M$, there are
pairwise incomparable strings $\sigma_0, \tau_0, \dots, \allowbreak \sigma_{m-1}, \tau_{m-1}$
such that $\sigma_i, \tau_i \in M(i)$ for every $i < m$.
\end{lemma}
\begin{proof}
Consider the following greedy algorithm.
At each stage, we maintain a set $P$ of pending rows which is initially the whole matrix $M$.
Pick a string $\rho$ of maximal length among all pending rows.
Let $M(i)$ be a pending row such that $\rho \in M(i)$.
If we have already chosen the value of $\sigma_i$, set $\tau_i = \rho$ and remove 
$M(i)$ from the pending rows. Otherwise, set $\sigma_i = \rho$.
In any case, remove every prefix of $\rho$ from any row of $M$ and go to the next step.

Notice that at any step, we remove at most one string from each row of $M$
since the strings in each row are pairwise incomparable.
Moreover, since we want to construct a sequence of $2m$
pairwise incomparable strings, and at each step we add one string to this sequence,
there are at most $2m$ steps. The algorithm gets stuck at some points only if all pending rows are empty,
which cannot happen since each row contains at least $2m$ strings.

We claim that the chosen strings are pairwise incomparable. Indeed, when at some stage, we pick
a string $\rho$, it is of shorter length than any string we have picked so far,
and cannot be a prefix of any of them since each time we pick a string, we remove its
prefixes from the matrix.
\end{proof}

\emph{Abstracting the requirements}.
The first feature of the framework of Lerman, Solomon and Towsner 
that we already exploited is the ``$\tto$-fairness'' of the $\tto^1_2$-instance
which allows each~$\rt^1_2$-instance to diagonalize him
as soon as the $\rt^1_2$-instance gives him enough occasions to do it.
We will now use the second aspect of this framework which consists in
getting rid of the complexity of the requirements by replacing them
with arbitrary computable predicates (or blackboxes).

Indeed, consider the case of two successive applications of~$\rt^1_2$.
Say that the first instance is $B_0 \cup B_1 = \omega$, and the second $C_0 \cup C_1 = \omega$.
We need to design the $\tto^1_2$-instance $A_0 \cup A_1 = 2^{<\omega}$ so that
there is an infinite set $G \cap B_i$ for some~$i < 2$ and an infinite set~$H \cap B_j$
for some~$j < 2$ such that $(G \cap B_i) \oplus (H \cap B_j)$ does not compute a
solution to the~$A$'s.
While constructing the $A$'s, we enumerate two levels of conditions.
We first enumerate the conditions~$c = (F, X)$ used for constructing the set~$G$,
but we also enumerate the conditions $c_0 = (F_0, X_0)$ and~$c_1 = (F_1, X_1)$
such that~$c_i$ is used to construct a solution~$H$ to the second $\rt^1_2$-instance~$C_0 \cup C_1 = \omega$
below~$G \cap B_i$. The question that the $\tto^1_2$-instance asks during its construction becomes

\smallskip
``For every 2-partition~$Z_0 \cup Z_1 = X$, is there some side~$i < 2$
and some set~$G \subseteq Z_i$ such that for every 2-partition~$W_0 \cup W_1 = X_i$,
there is some side~$j < 2$ and some set~$H \subseteq W_j$ 
such that $\Phi_{e_{i,j}}^{((F \cap B_i) \cup G) \oplus ((F_i \cap C_j) \cup H)}(q)$ halts?''
\smallskip

While staying~$\Sigma^0_1$ (with parameters), the question becomes
arbitrarily complicated to formulate. Moreover, looking at the shape of the question,
we see that the first iteration can box any $\Sigma^0_1$ question asked about the second iteration.
We can therefore abstract the question and make the $\tto$-fairness property independent
of the specificities of the forcing notion used to solve the~$\rt^1_2$-instances.
See~\cite{Patey2015Iterative} for detailed explanations about this abstraction process.

\begin{definition}[Formula, valuation]
An \emph{$m$-by-$n$ formula} is a formula $\varphi$
with distinguished set variables $U_{i,j}$ for each $i < m$ and $j < n$.
Given an $m$-by-$n$ matrix $M = \{\sigma_{i,j} : i < m, j < n\}$, 
an \emph{$M$-valuation} $V$ is a tuple of finite sets 
$A_{i,j} \subseteq \{\tau \in 2^{<\omega} : \tau \succeq \sigma_{i,j}\}$ 
for each $i < m$ and~$j < n$.
The valuation $V$ \emph{satisfies} $\varphi$ if $\varphi(A_{i,j} : i < m, j < n)$ holds.
We write $\varphi(V)$ for $\varphi(A_{i,j} : i < m, j < n)$.
\end{definition}

Given some valuation $V = (A_{i,j} : i < m, j < n)$ and some integer $s$, we write $V > s$
to say that for every $\tau \in A_{i,j}$, $|\tau| > s$. Moreover,
we denote by $V(i)$ the $n$-tuple $A_{i,0}, \dots, A_{i,n-1}$.
Following the terminology of~\cite{Lerman2013Separating}, we define 
the notion of essentiality for a formula (an abstract requirement),
which corresponds to the idea that there is room for diagonalization
since the formula is satisfied for arbitrarily far valuations.

\begin{definition}[Essential formula]
An $m$-by-$n$ formula $\varphi$ is \emph{essential} in an $m$-by-$n$ matrix $M$
if for every $s \in \omega$, there is an $M$-valuation $V > s$ such that $\varphi(V)$ holds.
\end{definition}

The notion of $\tto$-fairness is defined accordingly. If some formula
is essential, that is, gives enough room for diagonalization, then there is
an actual valuation which will diagonalize against the~$\tto^1_2$-instance.

\index{tt-fairness@$\tto$-fairness}
\begin{definition}[$\tto$-fairness]
Fix two sets $A_0, A_1 \subseteq 2^{<\omega}$.
Given an $m$-by-$n$ disjoint matrix $M$, an $M$-valuation~$V$ \emph{diagonalizes} against $A_0, A_1 \subseteq 2^{<\omega}$
if for every $i < m$, there is some~$L, R \in V(i)$ such that $L \subseteq A_0$ and $R \subseteq A_1$.
A set~$X$ is \emph{$n$-$\tto$-fair} for~$A_0, A_1$ if for every $m$ and every $\Sigma^{0,X}_1$ $m$-by-$2^nm$ formula
$\varphi$ essential in some disjoint matrix $M$, there is an $M$-valuation $V$ diagonalizing against $A_0, A_1$ 
such that $\varphi(V)$ holds.
A set $X$ is \emph{$\tto$-fair} for $A_0, A_1$ if it is $n$-$\tto$-fair for~$A_0, A_1$ for some~$n \geq 1$.
\end{definition}

Of course, if $Y \leq_T X$, then every $\Sigma^{0,Y}_1$ formula is $\Sigma^{0,X}_1$.
As an immediate consequence, if $X$ is $n$-$\tto$-fair for some $A_0, A_1$ and $Y \leq_T X$, then $Y$ is $n$-$\tto$-fair for $A_0, A_1$.
Moreover, if $X$ is $n$-$\tto$-fair for $A_0, A_1$ and $p > n$, $X$ is also $p$-$\tto$-fair for $A_0, A_1$ as witnessed
by cropping the rows.

\index{preservation!of $\tto$-fairness}
\begin{definition}[$\tto$-fairness preservation]
Fix a $\Pi^1_2$ statement $\Psf$.
\begin{itemize}
	\item[1.] $\Psf$ admits \emph{$\tto$-fairness (resp. $n$-$\tto$-fairness) preservation} if for all sets $A_0, A_1 \subseteq 2^{<\omega}$,
every set $C$ which is $\tto$-fair (resp. $n$-$\tto$-fair) for~$A_0, A_1$ and every $C$-computable $\Psf$-instance~$X$,
there is a solution $Y$ to $X$ such that $Y \oplus C$ is $\tto$-fair (resp. $n$-$\tto$-fair) for $A_0, A_1$.
	\item[2.] $\Psf$ admits strong \emph{$\tto$-fairness (resp. $n$-$\tto$-fairness) preservation} if for all sets $A_0, A_1 \subseteq 2^{<\omega}$,
every set $C$ which is $\tto$-fair (resp. $n$-$\tto$-fair) for~$A_0, A_1$ and every $\Psf$-instance~$X$,
there is a solution $Y$ to $X$ such that $Y \oplus C$ is $\tto$-fair (resp. $n$-$\tto$-fair) for $A_0, A_1$.
\end{itemize}
\end{definition}

Note that a principle $\Psf$ may admit $\tto$-fairness preservation without preserving $n$-$\tto$-fairness for any fixed~$n$,
as this is the case with~$\rt^2_2$ (see Theorem~\ref{thm:tree-theorem-rt22-tto-fairness-preservation}
and Theorem~\ref{thm:tree-theorem-rt22-not-1-tto-fairness-preservation}). 
On the other hand, if $\Psf$ admits $n$-$\tto$-fairness preservation for every $n$,
then it admits $\tto$-fairness preservation.
The notion of $\tto$-fairness is a weakness property. 
Therefore Lemma~\ref{lem:intro-reduc-preservation-separation} can be applied under the following form.

\begin{lemma}\label{lem:tree-theorem-tto-fairness-provability}
If $\Psf$ admits $\tto$-fairness preservation but not $\Qsf$, then $\Psf$ does not imply $\Qsf$
over $\rca$.
\end{lemma}

Now we have introduced the necessary terminology, we create a non-effective
instance of~$\tto^1_2$ which will serve as a bootstrap for $\tto$-fairness preservation.

\begin{lemma}\label{lem:tree-theorem-partition-emptyset-tto-fair}
There exists a $\Delta^0_2$ partition $A_0 \cup A_1 = 2^{<\omega}$ such that
$\emptyset$ is 1-$\tto$-fair for~$A_0, A_1$.
\end{lemma}
\begin{proof}
The proof is done by a no-injury priority construction.
Let $\varphi_0, \varphi_1, \dots$ be an effective enumeration of all $m$-by-$2m$ $\Sigma^0_1$ formulas
and $M_0, M_1, \dots$ be an enumeration of all $m$-by-$2m$ disjoint matrices for every $m$.
We want to satisfy the following requirements for each pair of integers~$e,k$.

\begin{quote}
$\Rcal_{e,k}$: If $\varphi_e$ is essential in $M_k$, then $\varphi_e(V)$ holds
for some $M_k$-valuation $V$ diagonalizing against $A_0, A_1$.
\end{quote}

The requirements are ordered via the standard pairing function $\tuple{\cdot, \cdot}$.
The sets $A_0$ and $A_1$ are constructed by a $\emptyset'$-computable list of 
finite approximations $A_{i,0} \subseteq A_{i,1} \subseteq \dots$
such that all elements added to~$A_{i,s+1}$ from~$A_{i,s}$
are strictly greater than the maximum of~$A_{i,s}$ for each~$i < 2$. 
We then let $A_i = \bigcup_s A_{i,s}$ which will be a~$\Delta^0_2$ set.
At stage 0, set $A_{0,0} = A_{1,0} = \emptyset$. Suppose that at stage $s$,
we have defined two disjoint finite sets $A_{0,s}$ and $A_{1,s}$ such that
\begin{itemize}
	\item[(i)] $A_{0,s} \cup A_{1,s} = 2^{<b}$ for some integer $b \geq s$
	\item[(ii)] $\Rcal_{e',k'}$ is satisfied for every $\tuple{e',k'} < s$
\end{itemize}
Let $\Rcal_{e,k}$ be the requirement such that $\tuple{e,k} = s$.
Decide $\emptyset'$-computably whether there is some $M_k$-valuation $V > b$
such that $\varphi_e(V)$ holds. If so, effectively fetch such a~$V$ 
and let $d$ be an upper bound on the length of the strings in $V$.
By Lemma~\ref{lem:tree-theorem-matrix-combi}, there are
pairwise incomparable strings $\sigma_0, \tau_0, \dots, \allowbreak \sigma_{m-1}, \tau_{m-1}$
such that $\sigma_i, \tau_i \in M(i)$ for every $i < m$.
For each $i < m$, let $A_{i,l}$ and $A_{i,r}$ be the sets in $V$ corresponding to $\sigma_i$
and $\tau_i$, respectively.
Set $A_{0,s+1} = A_{0,s} \bigcup_{i < m} A_{i,l}$ and $A_{1,s+1} = 2^{<d} \setminus A_{0,s+1}$.
This way, $A_{0,s+1} \cup A_{1,s+1} = 2^{<d}$.
Since the $\sigma$'s and $\tau$'s are pairwise incomparable, the sets $A_{i,l}$ and $A_{i,r}$ are disjoint,
so $\bigcup_{i < m} A_{i,r} \subseteq 2^{<d} \setminus A_{0,s+1}$ and the requirement $\Rcal_{e,i}$ is satisfied.
If no such $M_k$-valuation is found, the requirement $\Rcal_{e,k}$ is vacuously satisfied.
Set $A_{0,s+1} = A_{0,s} \cup 2^b$ and $A_{1,s+1} = A_{1,s}$.
This way, $A_{0,s+1} \cup A_{1,s+1} = 2^{<(b+1)}$.
In any case, go to the next stage. This finishes the construction.
\end{proof}

\begin{theorem}\label{thm:tree-theorem-tt22-not-tto-fairness}
$\tto^2_2$ does not admit $\tto$-fairness preservation.
\end{theorem}
\begin{proof}
Let~$A_0 \cup A_1 = 2^{<\omega}$ be the $\Delta^0_2$ partition constructed in Lemma~\ref{lem:tree-theorem-partition-emptyset-tto-fair}.
By Shoenfield's limit lemma~\cite{Shoenfield1959degrees}, there is a computable function $h : 2^{<\omega} \times \omega \to 2$ such 
that for each~$\sigma \in 2^{<\omega}$, $\lim_s h(\sigma, s)$ exists and $\sigma \in A_{\lim_s h(\sigma, s)}$.
Let~$f : [2^{<\omega}]^2 \to 2$ be the computable coloring defined by $f(\sigma, \tau) = h(\sigma, |\tau|)$
for each~$\sigma \prec \tau \in 2^{<\omega}$. Let~$S \cong 2^{<\omega}$ be a $\tto^2_2$-solution to~$f$
with witness isomorphism~$g : 2^{<\omega} \to S$ and witness color~$c < 2$.

Fix any~$n \geq 1$. 
We claim that $S$ is not $n$-$\tto$-fair for~$A_0, A_1$.
Let~$\varphi(U_j : j < 2^n)$ be the $1$-by-$2^n$ $\Sigma^{0,S}_1$ formula
which holds if for each $j < 2^n$, $U_j$ is a non-empty subset of $S$.
Let $M = (\sigma_j : j < 2^n)$ be the $1$-by-$2^n$ disjoint matrix defined 
for each $j < 2^n$ by
$\sigma_j = g(\tau_j)$ where $\tau_j$ is the $j$th string of length~$n$.
In other words, $\sigma_j$ is the $j$th node at level~$n$ in~$S$.
For every $s$, let $V_s$ be the $M$-valuation defined by $B_j = \{g(\rho)\}$ such that $\rho$ is the least string of length $max(n,s)$ extending~$\tau_j$.  Notice that $V_s > s$ and $\varphi(V_s)$ holds.
Therefore, the formula $\varphi$ is essential in $M$. 
For every $M$-valuation $V = (B_j : j < 2^n)$ such that $\varphi(V)$ holds,
there is no $j < 2^n$ such that $B_j \subseteq A_{1-c}$ since it would contradict
the fact that $B_j$ is a non-empty subset of $A_c$.
Therefore $S$ is not $n$-$\tto$-fair.
\end{proof}

Notice that we actually proved a stronger statement.
Dzhafarov, Hirst and Lakins defined in~\cite{Dzhafarov2010Ramseys} various notions
of stability for the tree theorem for pairs. A coloring~$f : [2^{<\omega}]^2 \to r$
is \emph{1-stable} if for every $\sigma \in 2^{<\omega}$, there is some threshold $t$
and some color~$c < r$ such that $f(\sigma, \tau) = c$ for every $\tau \succ \sigma$
such that $|\tau| \geq t$. In the proof of Theorem~\ref{thm:tree-theorem-tt22-not-tto-fairness}, we showed in fact that
$\tto^2_2$ restricted to 1-stable colorings does not admit $\tto$-fairness preservation.
In the same paper, Dzhafarov et al.\ studied an increasing polarized version of the tree theorem for pairs,
and proved that its 1-stable restriction coincides with the 1-stable tree theorem for pairs over~$\rca$.
Therefore the increasing polarized tree theorem for pairs does not admit $\tto$-fairness preservation.

\section{The strength of the tree theorem for pairs}\label{sect:tree-theorem-separating-principles}

In this section, we prove $\tto$-fairness preservation for various principles in reverse mathematics,
namely, weak K\"onig's lemma, cohesiveness and $\rt^2_2$. We prove independently that
they admit $\tto$-fairness preservation, and then use the compositional
nature of the notion of preservation to deduce that the conjunction of these
principles does not imply $\tto^2_2$ over $\rca$.

We start with weak König's lemma
which is involved in many proofs of Ramsey-type theorems,
and in particular in the proof of Ramsey's theorem for pairs.

\begin{theorem}\label{thm:tree-theorem-wkl-n-tto-fairness}
For every $n \geq 1$, $\wkl$ admits $n$-$\tto$-fairness preservation.
\end{theorem}
\begin{proof}
Let~$C$ be a set $n$-$\tto$-fair for some sets~$A_0, A_1 \subseteq 2^{<\omega}$,
and let~$T \subseteq 2^{<\omega}$ be a $C$-computable infinite binary tree.
We construct an infinite decreasing sequence of computable subtrees $T = T_0 \supseteq T_1 \supseteq \dots$
such that for every path $P$ through $\bigcap_s T_s$, $P \oplus C$ is $n$-$\tto$-fair for $A_0, A_1$.
More precisely, if we interpret $s$ as a tuple~$\tuple{m, \varphi, M}$ where $\varphi(G,U)$
is an $m$-by-$2^nm$ $\Sigma^{0,C}_1$ formula~$\varphi(G,U)$ and $M$ is an $m$-by-$2^nm$ disjoint matrix $M$,
we want to satisfy the following requirement.

\begin{quote}
$\Rcal_s$ : For every path~$P$ through~$T_{s+1}$, either $\varphi(P, U)$ is not essential in $M$,
or~$\varphi(P, V)$ holds for some $M$-valuation $V$ diagonalizing against $A_0, A_1$.
\end{quote}

Given two $M$-valuations $V_0 = (B_{i,j} : i < m, j < 2^nm)$ and $V_1 = (D_{i,j} :  i < m, j < 2^nm)$, 
we write $V_0 \subseteq V_1$
to denote the pointwise subset relation, that is, for every $i < m$ and every $j < 2^nm$, $B_{i,j} \subseteq D_{i,j}$.
At stage~$s = \tuple{m, \varphi, M}$, given some infinite, computable binary tree~$T_s$, 
define the $m$-by-$2^nm$ $\Sigma^{0,C}_1$ formula
\[
\psi(U) = (\exists n)(\forall \tau \in T_s \cap 2^n)(\exists \tilde{V} \subseteq U)\varphi(\tau, \tilde{V})
\]
We have two cases.
In the first case, $\psi(U)$ is not essential in $M$ with some witness~$t$. By compactness,
the following set is an infinite $C$-computable subtree of~$T_s$:
\[
T_{s+1} = \{ \tau \in T_s : (\mbox{for every } M\mbox{-valuation } V > t)\neg \varphi(\tau, V) \}
\]
The tree $T_{s+1}$ has been defined so that $\varphi(P, U)$
is not essential in $M$ for every~$P \in [T_{s+1}]$.

In the second case, $\psi(U)$ is essential in $M$. By $n$-$\tto$-fairness of $C$ for $A_0, A_1$,
there is an $M$-valuation $V$ diagonalizing against $A_0, A_1$ such that $\psi(V)$ holds.
We claim that for every path~$P \in [T_s]$,
$\varphi(P, \tilde{V})$ holds for some $M$-valuation~$\tilde{V}$ diagonalizing against~$A_0, A_1$.
Fix some path~$P \in [T_s]$. Unfolding the definition of~$\psi(V)$, there is some~$n$ such that $\varphi(P \uh n, \tilde{V})$ holds
for some~$M$-valuation~$\tilde{V} \subseteq V$. 
Since $V$ is diagonalizing against~$A_0, A_1$, for every $i < m$, there is some~$L, R \in V(i)$ 
such that $L \subseteq A_0$ and $R \subseteq A_1$.
Let $\tilde{L}, \tilde{R} \in \tilde{V}(i)$ be such that $\tilde{L} \subseteq L$ and $\tilde{R} \subseteq R$.
In particular, $\tilde{L} \subseteq A_0$ and $\tilde{R} \subseteq A_1$ so $\tilde{V}$ diagonalizes against $A_0, A_1$.
Take~$T_{s+1} = T_s$ and go to the next stage.
This finishes the proof of Theorem~\ref{thm:tree-theorem-wkl-n-tto-fairness}.
\end{proof}

As previously noted, preserving $n$-$\tto$-fairness for every~$n$ implies preserving $\tto$-fairness.
However, we really need the fact that $\wkl$ admits $n$-$\tto$-fairness preservation
and not only $\tto$-fairness preservation in the proof of Theorem~\ref{thm:tree-theorem-rt12-strong-tto-fairness}.

\begin{corollary}\label{cor:tree-theorem-wkl-tto-fairness-preservation}
$\wkl$ admits $\tto$-fairness preservation.
\end{corollary}

As usual, we shall decompose Ramsey's theorem for pairs
into a cohesiveness instance and a non-effective partition of the integers.
We shall therefore prove independently $\tto$-fairness preservation of~$\coh$
and strong $\tto$-fairness preservation of~$\rt^1_2$ to deduce that~$\rt^2_2$ admits $\tto$-fairness preservation.

\begin{theorem}\label{thm:tree-theorem-coh-n-tto-fairness}
For every $n \geq 1$, $\coh$ admits $n$-$\tto$-fairness preservation.
\end{theorem}
\begin{proof}
Let~$C$ be a set $n$-$\tto$-fair for some sets~$A_0, A_1 \subseteq 2^{<\omega}$,
and let~$R_0, R_1, \dots$ be a $C$-computable sequence of sets.
We will construct an $\vec{R}$-cohesive set $G$ such that 
$G \oplus C$ is $n$-$\tto$-fair for~$A_0, A_1$.
The construction is done by a Mathias forcing, whose conditions are pairs $(F, X)$
where $X$ is a $C$-computable set. The result is a direct consequence of the following lemma.

\begin{lemma}\label{lem:tree-theorem-coh-preservation-lemma}
For every condition~$(F, X)$, every $m$-by-$2^nm$ $\Sigma^{0,C}_1$ formula~$\varphi(G, U)$
and every $m$-by-$2^nm$ disjoint matrix $M$, there exists an extension~$d = (E, Y)$ such that
$\varphi(G, U)$ is not essential for every set $G$ satisfying $d$, or $\varphi(E, V)$ holds
for some $M$-valuation $V$ diagonalizing against $A_0, A_1$.
\end{lemma}
\begin{proof}
Define the $m$-by-$2^nm$ $\Sigma^{0,C}_1$ formula
$\psi(U) = (\exists G \supseteq F)[G \subseteq F \cup X \wedge \varphi(G, U)]$.
By $n$-$\tto$-fairness of $C$ for $A_0, A_1$, 
either $\psi(U)$ is not essential in $M$, or~$\psi(V)$ holds for some $M$-valuation $V$
diagonalizing against $A_0, A_1$.
In the former case, the condition~$(F,X)$ already satisfies the desired property.
In the latter case, by the finite use property, there exists a finite set~$E$ satisfying~$(F, X)$ such that~$\varphi(E, V)$ holds.
Let $Y = X \setminus [0, max(E)]$. The condition $(E, Y)$ is a valid extension.
\end{proof}

Using Lemma~\ref{lem:tree-theorem-coh-preservation-lemma}, define an infinite descending sequence 
of conditions~$c_0 = (\emptyset, \omega) \geq c_1 \geq \dots$
such that for each~$s \in \omega$
\begin{itemize}
	\item[(i)] $|F_s| \geq s$
	\item[(ii)] $X_{s+1} \subseteq R_s$ or $X_{s+1} \subseteq \overline{R}_s$
	\item[(iii)] $\varphi(G, U)$ is not essential in $M$ for every set $G$ satisfying $c_{s+1}$, 
	or $\varphi(F_{s+1}, V)$ holds
	for some $M$-valuation $V$ diagonalizing against $A_0, A_1$ if~$s = \tuple{\varphi, M}$
\end{itemize}
where~$c_s = (F_s, X_s)$. The set $G = \bigcup_s F_s$ is infinite by (i),
$\vec{R}$-cohesive by (ii) and $G \oplus C$ is $n$-$\tto$-fair for~$A_0, A_1$ by (iii).
This finishes the proof of Theorem~\ref{thm:tree-theorem-coh-n-tto-fairness}.
\end{proof}

\begin{corollary}
$\coh$ admits $\tto$-fairness preservation.
\end{corollary}

The next theorem is the reason why we use the notion of
$\tto$-fairness instead of $n$-$\tto$-fairness
in our separation of~$\rt^2_2$ from~$\tto^2_2$.
Indeed, given an instance of~$\rt^1_2$ and a set~$C$
which is $n$-$\tto$-fair for some sets~$A_0, A_1$, the proof
constructs a solution $H$ such that~$H \oplus C$ is $(n+1)$-$\tto$-fair for~$A_0, A_1$.
We shall see in Corollary~\ref{cor:tree-theorem-rt12-not-strong-1-tto-fairness} that the proof is optimal,
in the sense that $\rt^1_2$ does not admit strong~$n$-$\tto$-fairness preservation.

\begin{theorem}\label{thm:tree-theorem-rt12-strong-tto-fairness}
$\rt^1_2$ admits strong $\tto$-fairness preservation.
\end{theorem}
\begin{proof}
Let~$C$ be a set $n$-$\tto$-fair for some sets~$A_0, A_1 \subseteq 2^{<\omega}$,
and let~$B_0 \cup B_1 = \omega$ be a (non-effective) 2-partition of~$\omega$.
Suppose that there is no infinite set $H \subseteq B_0$ or $H \subseteq B_1$
such that $H \oplus C$ is $n$-$\tto$-fair for~$A_0, A_1$, since otherwise we are done.
We construct a set $G$ such that both $G \cap B_0$ and $G \cap B_1$ are infinite.
We need therefore to satisfy the following requirements for each~$p \in \omega$.
\[
  \Ncal_p : \hspace{20pt} (\exists q_0 > p)[q_0 \in G \cap B_0] 
		\hspace{20pt} \wedge \hspace{20pt} (\exists q_1 > p)[q_1 \in G \cap B_1] 
\]
Furthermore, we want to ensure that one of $(G \cap B_0) \oplus C$ 
and $(G \cap B_1) \oplus C$ is $\tto$-fair for~$A_0, A_1$. To do this, we will satisfy the following requirements
for every integer $m$, every $m$-by-$2^{n+1}m$ $\Sigma^{0,C}_1$ formulas $\varphi_0(H, U)$ and  $\varphi_1(H, U)$
and every $m$-by-$2^{n+1}m$ disjoint matrices $M_0$ and $M_1$.
\[
  \Qcal_{\varphi_0, M_0, \varphi_1, M_1} : \hspace{20pt} 
		\Rcal_{\varphi_0,M_0}^{G \cap B_0} \hspace{20pt} \vee \hspace{20pt}  \Rcal_{\varphi_1,M_1}^{G \cap B_1}
\]
where $\Rcal_{\varphi, M}^H$ holds if $\varphi(H, U)$ is not essential in $M$
or $\varphi(H, V)$ holds for some $M$-valuation $V$ diagonalizing against $A_0, A_1$.
We first justify that if every $\Qcal$-requirement is satisfied, then either $(G \cap B_0) \oplus C$
or $(G \cap B_1) \oplus C$ is $(n+1)$-$\tto$-fair for $A_0, A_1$.
By the usual pairing argument, for every~$m$, there is some side $i < 2$ such that
the following property holds:
\begin{quote}
(P) For every $m$-by-$2^{n+1}m$ $\Sigma^{0,C}_1$ formula $\varphi(G \cap B_i, U)$ 
and every $m$-by-$2^{n+1}m$ disjoint matrix $M$, either $\varphi(G \cap B_i, U)$
is not essential in $M$, or $\varphi(G \cap B_i, V)$ holds for some $M$-valuation~$V$
diagonalizing against $A_0, A_1$.
\end{quote}
By the infinite pigeonhole principle, 
there is some side $i < 2$ such that (P) holds for infinitely many~$m$.
By a cropping argument, if (P) holds for $m$ and $q < m$, then (P) holds for~$q$.
Therefore (P) holds for every $m$ on side~$i$. In other words, $(G \cap B_i) \oplus C$
is $(n+1)$-$\tto$-fair for~$A_0, A_1$.

We construct our set $G$ by forcing. Our conditions are Mathias conditions~$(F, X)$,
such that~$X \oplus C$ is $n$-$\tto$-fair for~$A_0, A_1$.
We now prove the progress lemma, stating that we can force both $G \cap B_0$
and $G \cap B_1$ to be infinite.

\begin{lemma}\label{lem:tree-theorem-rt12-tto-fairness-progress}
For every condition~$c = (F, X)$, every $i < 2$ and every~$p \in \omega$
there is some extension~$d = (E, Y)$ such that~$E \cap B_i \cap (p,+\infty) \neq \emptyset$.
\end{lemma}
\begin{proof}
Fix~$c$, $i$ and~$p$. If $X \cap B_i \cap (p,+\infty) = \emptyset$,
then $X \cap (p,+\infty)$ is an infinite subset of~$B_{1-i}$.
Moreover, $X \cap (p,+\infty)$ is $n$-$\tto$-fair for~$A_0, A_1$, contradicting our hypothesis.
Thus, there is some~$q > p$ such that $q \in X \cap B_i \cap (p,+\infty)$.
Take $d = (F \cup \{q\}, X \setminus [0,q])$ as the desired extension.
\end{proof}

Next, we prove the core lemma stating that we can satisfy each $\Qcal$-requirement.
A condition~$c$ \emph{forces} a requirement~$\Qcal$
if $\Qcal$ is holds for every set~$G$ satisfying~$c$.
This is the place where we really need the fact that~$\wkl$
admits $n$-$\tto$-fairness preservation and not only $\tto$-fairness preservation.

\begin{lemma}\label{lem:tree-theorem-rt12-tto-fairness-forcing}
For every condition~$c = (F, X)$, every integer $m$, every $m$-by-$2^{n+1}m$ $\Sigma^{0,C}_1$ formulas $\varphi_0(H, U)$ and  $\varphi_1(H, U)$ and every $m$-by-$2^{n+1}m$ disjoint matrices $M_0$ and $M_1$,
there is an extension~$d = (E, Y)$ forcing~$\Qcal_{\varphi_0, M_0, \varphi_1, M_1}$.
\end{lemma}
\begin{proof}
Let~$\psi(U_0,U_1)$ be the $2m$-by-$2^{n+1}m$ $\Sigma^{0,X \oplus C}_1$ formula which holds
if for every 2-partition~$Z_0 \cup Z_1 = X$, there is some~$i < 2$,
some finite set $E \subseteq Z_i$
and an $m$-by-$2^{n+1}m$ $M_i$-valuation $V \subseteq U_i$ such that $\varphi_i((F \cap B_i) \cup E, V)$
holds. By $n$-$\tto$-fairness of~$X \oplus C$, we have two cases.

In the first case, $\psi(U_0, U_1)$ is not essential in~$M_0,M_1$, with some witness~$t$.
By compactness, the $\Pi^{0, X \oplus C}_1$ class~$\Ccal$ of sets~$Z_0 \oplus Z_1$
such that $Z_0 \cup Z_1 = \omega$ and for every~$i < 2$
and every finite set $E \subseteq Z_i$, there is no $M_i$-valuation $V > t$
such that $\varphi_i((F \cap B_i) \cup E, V)$ holds is non-empty.
By $n$-$\tto$-fairness preservation of~$\wkl$ (Theorem~\ref{thm:tree-theorem-wkl-n-tto-fairness}), there is a 2-partition
$Z_0 \oplus Z_1 \in \Ccal$ such that~$Z_0 \oplus Z_1 \oplus C$ is $n$-$\tto$-fair for~$A_0, A_1$.
Since $Z_0 \cup Z_1 = X$, there is some~$i < 2$ such that~$Z_i$ is infinite. Take such an~$i$.
The condition~$d = (F, Z_i)$ is an extension forcing $\Qcal_{\varphi_0, M_0, \varphi_1, M_1}$
by the $i$th side.

In the second case, $\psi(V_0, V_1)$ holds for some $(M_0,M_1)$-valuation $(V_0, V_1)$
diagonalizing against $A_0, A_1$. Let~$Z_0 = X \cap B_0$ and~$Z_1 = X \cap B_1$.
By hypothesis, there is some $i < 2$, some finite set~$E \subseteq Z_i = X \cap B_i$
and some $M_i$-valuation~$V \subseteq V_i$ such that $\varphi_i((F \cap B_i) \cup E, V)$ holds.
Since $V \subseteq V_i$, the $M_i$-valuation~$V$ diagonalizes against~$A_0, A_1$.
The condition~$d = (F \cup E, X \setminus [0, max(E)])$ is an extension
forcing $\Qcal_{\varphi_0, M_0, \varphi_1, M_1}$ by the $i$th side.
\end{proof}

Using Lemma~\ref{lem:tree-theorem-rt12-tto-fairness-progress} and Lemma~\ref{lem:tree-theorem-rt12-tto-fairness-forcing}, define an infinite descending sequence 
of conditions~$c_0 = (\emptyset, \omega) \geq c_1 \geq \dots$
such that for each~$s \in \omega$
\begin{itemize}
	\item[(i)] $|F_s \cap B_0| \geq s$ and~$|F_s \cap B_1| \geq s$
	\item[(ii)] $c_{s+1}$ forces~$\Qcal_{\varphi_0, M_0, \varphi_1, M_1}$
	if $s = \tuple{\varphi_0, M_0, \varphi_1, M_1}$
\end{itemize}
where~$c_s = (F_s, X_s)$. Let~$G = \bigcup_s F_s$.
The sets~$G \cap B_0$ and~$G \cap B_1$ are both infinite by (i)
and one of~$G \cap B_0$ and~$G \cap B_1$ is $\tto$-fair for~$A_0, A_1$ by (ii).
This finishes the proof of Theorem~\ref{thm:tree-theorem-rt12-strong-tto-fairness}.
\end{proof}

\begin{theorem}\label{thm:tree-theorem-rt22-tto-fairness-preservation}
$\rt^2_2$ admits $\tto$-fairness preservation.
\end{theorem}
\begin{proof}
Fix any set~$C$ $\tto$-fair for some sets $A_0, A_1 \subseteq 2^{<\omega}$ and any $C$-computable
coloring $f : [\omega]^2 \to 2$.
Consider the uniformly~$C$-computable sequence of sets~$\vec{R}$ defined for each~$x \in \omega$ by
\[
R_x = \{s \in \omega : f(x,s) = 1\}
\]
As~$\coh$ admits $\tto$-fairness preservation, there is
some~$\vec{R}$-cohesive set~$G$ such that $G \oplus C$ is $\tto$-fair for $A_0, A_1$.
The set~$G$ induces a $(G \oplus C)'$-computable coloring~$\tilde{f} : \omega \to 2$ defined by:
\[
(\forall x \in \omega) \tilde{f}(x) = \lim_{s \in G} f(x,s)
\]
As~$\rt^1_2$ admits strong $\tto$-fairness preservation,
there is an infinite $\tilde{f}$-homogeneous set~$H$ such that
$H \oplus G \oplus C$ is $\tto$-fair for $A_0, A_1$.
The set $H \oplus G \oplus C$ computes an infinite $f$-homogeneous set.
\end{proof}

We are now ready to prove our main theorem.

\begin{proof}[Proof of Theorem~\ref{thm:tree-theorem-rt22+wkl-not-tt22}]
By Theorem~\ref{thm:tree-theorem-rt22-tto-fairness-preservation}
and Corollary~\ref{cor:tree-theorem-wkl-tto-fairness-preservation}, $\rt^2_2$ and $\wkl$ admit $\tto$-fairness preservation.
By Theorem~\ref{thm:tree-theorem-tt22-not-tto-fairness}, $\tto^2_2$ does not admit $\tto$-fairness preservation.
We conclude by Lemma~\ref{lem:tree-theorem-tto-fairness-provability}.
\end{proof}

We now prove the optimality of Theorem~\ref{thm:tree-theorem-rt12-strong-tto-fairness} and Theorem~\ref{thm:tree-theorem-rt22-tto-fairness-preservation}
by showing that~$n$-$\tto$-fairness cannot be preserved.

\begin{theorem}\label{thm:tree-theorem-rt22-not-1-tto-fairness-preservation}
$\srt^2_2$ does not admit $n$-$\tto$-fairness preservation for any $n \geq 1$.
\end{theorem}
\begin{proof}
Let~$A_0 \cup A_1 = 2^{<\omega}$ be the $\Delta^0_2$ partition constructed in Lemma~\ref{lem:tree-theorem-partition-emptyset-tto-fair}.
By Shoenfield's limit lemma~\cite{Shoenfield1959degrees}, there is a stable computable function $f : [\omega]^2 \to 2$
such that $x \in A_{\lim_s f(x,s)}$ for every~$x$.
Fix some $n \geq 1$.
For each~$\sigma \in 2^{n+1}$, apply $\srt^2_2$ to the coloring~$f$ restricted to
the set $S_\sigma = \{\tau \succeq \sigma\}$
to obtain an infinite $f$-homogeneous set $H_\sigma$ for some color~$c_\sigma < 2$.
By definition of~$f$, $H_\sigma \subseteq A_{c_\sigma}$.
By the finite pigeonhole principle, there is a color~$c < 2$ and a set $M \subseteq 2^{n+1}$ of size $2^n$
such that $c_\sigma = c$ for every $\sigma \in M$. We can see $M$ as a 1-by-$2^n$ disjoint matrix.
Let $H = \bigoplus_{\sigma \in M} H_\sigma$ and let $\varphi(U_\sigma : \sigma \in M)$ be the $1$-by-$2^n$ $\Sigma^{0,H}_1$ formula
which holds if for every $\sigma \in M$, $U_\sigma$ is a non-empty subset of $H_\sigma$.
Note that $H_\sigma \subseteq A_c$ for every $\sigma \in M$.
The formula $\varphi(U)$ is essential in $M$ but there is no $M$-valuation $V = (V_\sigma : \sigma \in M)$ 
such that $\varphi(V)$ holds and $V_\sigma \subseteq A_{1-c}$ for some $\sigma \in M$.
Therefore $H$ is not $n$-$\tto$-fair for~$A_0, A_1$.
\end{proof}

\begin{corollary}\label{cor:tree-theorem-rt12-not-strong-1-tto-fairness}
$\rt^1_2$ does not admit $n$-$\tto$-fairness preservation for every~$n \geq 1$.
\end{corollary}
\begin{proof}
Fix some~$n \geq 2$. By Theorem~\ref{thm:tree-theorem-rt22-not-1-tto-fairness-preservation},
there is some set $C$ $n$-$\tto$-fair some $A_0, A_1$ and a stable $C$-computable
function $f : [\omega]^2 \to 2$ such that for every infinite $f$-homogeneous set~$H$,
$H \oplus C$ is not $n$-$\tto$-fair for~$A_0, A_1$.
Let~$\tilde{f} : \omega \to 2$ be defined by $\tilde{f}(x) = \lim_s f(x,s)$.
Every infinite $\tilde{f}$-homogeneous set $H$ $C$-computes an infinite 
$f$-homogeneous set $H_1$ such that $H_1 \oplus C$ is not $n$-$\tto$-fair for~$A_0, A_1$.
Therefore $H \oplus C$ is not $n$-$\tto$-fair for~$A_0, A_1$.
\end{proof}

\section{An Erd\H{o}s-Rado theorem}

In this section, we investigate the reverse mathematics of a well-known theorem due to 
Erd\H{o}s and Rado about partitions of rationals~\cite{Frittaion2015Coloring}.
This theorem is a natural strengthening
of Ramsey's theorem for pairs in the following sense. Ramsey's theorem for pairs can be 
stated as $\omega\to(\omega)^2_2$,
where~$\alpha \to (\beta)^2_2$ is the statement ``For every coloring $f\colon[L]^2\to 
2$, where $L$ is a linear order of order type $\alpha$, there is a homogeneous set $H$ 
such that $(H,\leq_L)$ has order type $\beta$''. It turns out that~$\omega$ 
and~$\omega^{*}$
are the only countable order types~$\alpha$ such that~$\alpha \to (\alpha)^2_2$ holds.
In particular, $\eta \to (\eta)^2_2$ does not hold, where $\eta$ is the order type of the 
rationals. A standard counterexample is as follows. Fix a 
one-to-one map 
$j\colon\Qb\to\Nb$. Define $f\colon[\Qb]^2\to 2$ by 
letting
\[   f(x,y)=\begin{cases}
             0   &   \text{if } x<_\Qb y \land j(x)<j(y); \\
             1 &     \text{if } x<_\Qb y\land j(x)>j(y).
            \end{cases} \]
A dense homogeneous set would give an embedding of $\Qb$ into $\omega$ (color $0$) or 
$\omega^*$ (color $1$), which is impossible. Even though Ramsey's theorem for 
rationals fails,  Erd\H{o}s and Rado \cite[Theorem 
4, p. 427]{Erdos1952Combinatorial} proved the following
Ramsey-type theorem (see also Rosenstein \cite[Theorem 11.7, p. 
207]{Rosenstein1982Linear}).

\index{Erd\H{o}s-Rado theorem}
\index{$\alpha \to (\beta)^2_2$}
\begin{theorem}[Erd\H{o}s-Rado theorem]\label{ER52}
The partition relation $\eta\to(\aleph_0,\eta)^2$ holds.
\end{theorem}

The statement $\eta\to(\aleph_0,\eta)^2$ asserts that for every coloring $f\colon[L]^2 
\to 2$, where $L$ is a linear order of order type~$\eta$, there is either an infinite  
0-homogeneous set or a 1-homogeneous set of order type~$\eta$.\smallskip
The partition relation~$\eta\to(\aleph_0,\eta)^2$ can be easily formalized in~$\rca$ by 
fixing a computable presentation $\Qb$ of the rationals. We may safely assume 
that the domain of $\Qb$ is $\Nb$. In order to study  $\erp$ we also consider a 
version of the infinite pigeonhole principle over the rationals, namely the statement 
``For every $n$ and for 
every  $n$-coloring $f\colon\Qb\to n$ there exists a dense homogeneous set'', which we 
denote by $\ers$.  

\index{$\erp$|see {Erd\H{o}s-Rado theorem}}
We start off the analysis of the Erd\H{o}s-Rado theorem by proving that the 
statement $\erp$ lies between $\aca$ and~$\rt^2_2$. On the lower bound hand,
$\erp$ can be seen as an immediate strengthening
of $\rt^2_2$. The upper bound is an effectivization of the original proof
of~$\erp$ by Erd\H{o}s and Rado in~\cite{Erdos1952Combinatorial}. 

\begin{lemma}
$\erp$ implies~$\rt^2_2$ over~$\rca$.
\end{lemma}
\begin{proof}
An instance of $\rt^2_2$ can be regarded as an instance of $\erp$. Moreover, 
provably in $\rca$, a dense set is infinite.
\end{proof}

The rest of this section is devoted to show that $\erp$ is provable in $\aca$. For this 
purpose, we give the following definition.

\index{somewhere dense}
\index{nowhere dense}
\begin{definition}[$\rca$]
By \emph{interval} we mean a set of the form $I=(x,y)_\Qb$ for $x,y\in \Qb$. We say that
$A\subseteq\Qb$ is \emph{somewhere dense} if $A$ is dense in some
interval of $\Qb$,  i.e.,\ there exists an interval $I$ such that for all intervals 
$J\subseteq I$ we have that $A\cap J\neq\emptyset$. We call $A$ \emph{nowhere 
dense}  otherwise.
\end{definition}
Notice that the above notion is the usual notion of topological density with respect to
the order topology of $\Qb$. In general, the nowhere dense sets of a topological 
space form an ideal. This is crucial in the proof by Erd\H{o}s and Rado. For this reason, 
we also use the terminology \emph{positive} and \emph{small} for somewhere dense and 
nowhere dense respectively. In $\rca$ we can show that nowhere dense 
subsets of $\Qb$ are small in the set-theoretic sense, meaning that:
\begin{itemize}
 \item[(1)] If $A\subseteq\Qb$ is small and $B\subseteq A$, then $B$ is small;
 \item[(2)] If $A,B\subseteq\Qb$ are small, then $A\cup B$ is small.
\end{itemize}

With enough induction, it is possible to generalize $(2)$ to finitely 
many sets. 

\begin{lemma}[$\rca+\ist$]\label{lem:erdos-rado-small}
If $A_i$ is a small subset of $\Qb$ for all $i<n$, then $\bigcup_{i<n}A_i$ 
is small. 
\end{lemma}
\begin{proof}
Suppose that $A_i$ is small for every $i<n$. Fix an interval $I$. We aim to show that 
$A^n=\bigcup_{i<n}A_i$ is not dense on $I$.  By $\ist$ (induction on $i$), we prove that
\[ (\forall i\leq n)(\exists J\subseteq I\text{ interval})A^i\cap J=\emptyset, \]
where $A^i=\bigcup_{j<i}A_j$. The conclusion follows from $i=n$. The 
case $i=0$ is vacuously true. Suppose the property holds for $i$ and that $i+1\leq n$. By 
the induction hypothesis there 
exists $J\subseteq I$ such that $A^i\cap J=\emptyset$. As $A_i$ is small, there exists 
$K\subseteq J$ such that $A_i\cap K=\emptyset$. It follows that 
$A^{i+1}\cap K=(A^i\cup A_i)\cap K=\emptyset$, as desired.  
\end{proof}

\begin{theorem}\label{thm:erdos-rado-upper}
$\aca$ implies $\erp$ over~$\rca$.
\end{theorem}
\begin{proof}
Let $f\colon [\Qb]^2\to 2$ be given. For any $x\in \Qb$, let 
$\red(x)=\{y\in\Qb\setminus\{x\}\colon f(x,y)=0\}$. Define 
$\blue(x)$ accordingly. We say that $A\subseteq \Qb$ is
$0$\emph{-admissible} if there exists some $x\in A$ such that $A\cap\red(x)$ is 
positive.

Case I. Every positive subset of $\Qb$ is $0$-admissible. We aim to show that 
there exists an infinite $0$-homogeneous set. We define by arithmetical recursion a 
sequence $(x_n)_{n_\in\Nb}$ as 
follows.
Supppose we have defined $x_i$ for all $i<n$, and assume by arithmetical 
induction that $A_n=\bigcap_{i<n}\red(x_i)$ is positive, and hence $0$-admissible (where $\bigcap_{i < 0} \red(x_i) = \Qb$).  
Search for the $\omega$-least $x_n\in A_n$ such that 
$A_n\cap\red(x_n)=\bigcap_{i<{n+1}}\red(x_i)$ 
is positive. By definition, the set $\{x_n\colon n\in\Nb\}$ is infinite and
$0$-homogeneous.

Case II. There is a positive subset $A$ of $\Qb$ which is not 
$0$-admissible. In this case, we show that there exists a dense 
$1$-homogeneous set. Let $I$ be a witness of $A$ being positive. Fix an 
enumeration $(I_n)_{n\in\Nb}$ of all subintervals of $I$. Notive that by definition $A$ 
intersects every $I_n$.

We define by arithmetical recursion a sequence $(x_n)_{n\in\Nb}$ as follows. Let 
$x_0\in A\cap I_0$. Suppose we have defined $x_i\in A\cap I_i$ for all
$i<n$. By Lemma \ref{lem:erdos-rado-small}, since every $A\cap\red(x_i)$ with $i<n$ is small, 
it follows that $E=\bigcup_{i<n} \big(A\cap\red(x_i)\big)$ is small.
Let $J\subseteq I_n$ be such that $E\cap J=\emptyset$. We may safely assume that no 
$x_i$ with $i<n$ belongs to $J$. Since $A$ is dense in $I$ and $J\subseteq I$, we can 
find $x_n\in A\cap J$. In particular, $x_n\in
\bigcap_{i<n}\blue(x_i)$. Therefore $\{x_n\colon n\in\N\}$ is a dense
$1$-homogeneous set.
\end{proof}

\section{The strength of the Erd\H{o}s-Rado theorem}

In this section, we prove that the Erd\H{o}s-Rado theorem for pairs
does not reduce to Ramsey's theorem for pairs in one step.

\begin{theorem}\label{thm:er22-not-computably-reduces}
$\erp \not \leq_c \rt^2_{<\infty}$.
\end{theorem}

Interestingly, this diagonalization does not seem to be easily
generalizable to a separation over~$\omega$-models.
A reason is that the fairness property  ensured by the~$\erp$-instance
does not seem to be preserved by weak K\"onig's lemma.
This is hitherto the first example of a computable non-reducibility
of a principle~$\Psf$ to~$\rt^2_{<\infty}$ which is not generalizable
to a proof that~$\rt^2_2$ does not imply~$\Psf$ over~$\rca$.
Indeed, a diagonalization against an $\rt^2_4$-instance is similar
to a diagonalization against two $\rt^2_2$-instances. Therefore, diagonalizing against~$\rt^2_{<\infty}$
has some common flavor with a separation over standard models.

The remainder of this section is devoted to a proof of Theorem~\ref{thm:er22-not-computably-reduces}.
The notion of fairness presented below may have some ad-hoc flavor.
It has been obtained by applying the main ideas of the framework
of Lerman, Solomon and Towsner~\cite{Lerman2013Separating,Patey2015Iterative}.
Thanks to an analysis of the combinatorics of Ramsey's theorem for pairs
and the Erd\H{o}s-Rado theorem for pairs, we prove our
computable non-reducibility result by constructing an instance of~$\erp$
ensuring the density of the diagonalizing conditions in the forcing notion of~$\rt^2_2$.
Then we abstract the diagonalization to 
any $\Sigma^0_1$ formula, to get rid of the specificities
of the forcing notion of~$\rt^2_2$ in the notion of $\er$-fairness preservation.

\index{simple partition}
\begin{definition}[Simple partition]
A \emph{simple partition}~$\inter_\Qb(S)$ is a finite sequence of open intervals~$(-\infty, x_0), (x_0, x_1), \dots, (x_{n-1}, +\infty)$
for some set of rationals~$S = \{x_0 <_\Qb \dots <_\Qb x_{n-1}\}$. We set 
$\inter_\Qb(\emptyset) = \{\Qb\}$. 
A simple partition~$I_0, \dots, I_{n-1}$
\emph{refines} another simple partition~$J_0, \dots, J_{m-1}$ if for every~$i < n$,
there is some~$j < m$ such that~$I_i \subseteq J_j$.
Given two simple partitions~$I_0, \dots, I_{n-1}$ and~$J_0, \dots, J_{m-1}$,
the product~$\vec{I} \otimes \vec{J}$ is the simple partition
\[
\{ I \cap J : I \in \vec{I} \wedge J \in \vec{J} \}
\]
\end{definition}
One can easily see that~$\inter_\Qb(S)$ refines~$\inter_\Qb(T)$ if~$T \subseteq S$
and that $\inter_\Qb(S\cup T)=\inter_\Qb(S)\otimes\inter_\Qb(T)$.
Note that every simple partition has a finite description, since the set~$S$
and each rational has a finite description.
Also note that a simple partition is not a true partition of~$\Qb$ since the endpoints do not belong
to any interval. However, we have~$S \cup \bigcup \inter_\Qb(S) = \Qb$.

\begin{definition}[Matrix]
An \emph{$m$-by-$n$ matrix} $M$ is a rectangular array of rationals $x_{i,j} \in \Qb$
such that $x_{i,j} <_\Qb x_{i,k}$ for each~$i < m$ and $j < k < 
n$. 
The $i$th \emph{row} $M(i)$ of the matrix $M$ is the $n$-tuple of rationals $x_{i,0} < 
\dots < x_{i,n-1}$. The simple partition $\inter_\Qb(M)$ is defined by  $\bigotimes_{i < 
m} \inter_\Qb(M(i))$. In particular, $\bigotimes_{i <m} \inter_\Qb(M(i))$ refines the 
simple partition~$\inter_\Qb(M(i))$ for each~$i < m$.
\end{definition}

It is important to notice that an $m$-by-$n$ matrix is formally a 3-tuple~$\tuple{m,n, M}$
and not only the matrix itself~$M$. This distinction becomes important
when dealing with the degenerate cases. An~$m$-by-0 matrix~$M$
and a 0-by-$n$ matrix~$N$ are both empty.  However, they have different sizes.
In particular, we shall define the notion of~$M$-type for a matrix, and this definition
will depend on the number of columns of the matrix~$M$, which is~0 for~$M$, and~$n$ for~$N$.
Notice also that, for a degenerate matrix $M$, the simple partition $\inter_{\Qb}(M)$ 
is the singleton $\{\Qb\}$. 

Given a simple partition~$\vec{I}$, we want to classify the $k$-tuples of rationals
according to which interval of~$\vec{I}$ they belong to.
This leads to the notion of~$(\vec{I},k)$-type.

\begin{definition}[Type]
Given a simple partition~$I_0, \dots, I_{n-1}$ and some~$k \in \omega$,
an \emph{$(\vec{I},k)$-type} is a tuple~$T_0, \dots, T_{k-1}$
such that~$T_i \in \vec{I}$ for each~$i < k$.
Given an $m$-by-$n$ matrix $M$, an~\emph{$M$-type} is an~$(\inter_\Qb(M), n)$-type.
\end{definition}

We now state two simple combinatorial lemmas which will be useful later.
The first trivial lemma simply states that each $m$-tuple of rational (different from
the endpoints of a simple partition) belongs to a type.

\begin{lemma}\label{lem:tuple-has-mtype}
For every simple partition $I_0, \dots, I_{n-1}$ and every~$k$-tuple
of rationals~$x_0, \dots, x_{k-1} \in \bigcup_{i < n} I_i$, there is an $(\vec{I},k)$-type $T_0, \dots, T_{k-1}$
such that~$x_j \in T_j$ for each~$j < k$.
\end{lemma}
\begin{proof}
Fix $k$ rationals $x_0, \dots, x_{k-1}$. For each~$i < k$, 
there is some interval~$T_i \in \vec{I}$ such that~$x_i \in T_i$ since $x_i \in \bigcup_{j < n} 
I_j$. The sequence $T_0, \dots, T_{k-1}$ is the desired~$(\vec{I}, k)$-type.
\end{proof}

The next lemma is a consequence of the pigeonhole principle.

\begin{lemma}\label{lem:mtype-interval-disjoint}
For every $m$-by-$n$ matrix~$M$ and every $M$-type $T_0, \dots, T_{n-1}$,
there is an $m$-tuple of intervals~$J_0, \dots, J_{m-1}$ with~$J_i \in \inter_\Qb(M(i))$ such that
\[
(\bigcup_{j < n} T_j) \cap (\bigcup_{i < m} J_i) = \emptyset
\]
\end{lemma}
\begin{proof}
Let $T_0, \dots, T_{n-1}$ be an~$M$-type.
For every~$i < m$ and~$j < n$, there is some $J \in \inter_\Qb(M(i))$ such that~$T_j \subseteq J$.
Since $|\inter_\Qb(M(i))| = n+1$, there is an interval~$J_i \in \inter_\Qb(M(i))$ such that $(\bigcup_{j < n} T_j) \cap J_i = \emptyset$.
\end{proof}

\begin{definition}[Formula, valuation]
Given an $m$-by-$n$ matrix $M$, an \emph{$M$-formula} 
is a formula $\varphi(\vec{U}, \vec{V})$ with distinguished set variables $U_j$ for each $j < n$
and~$V_{i,I}$ for each~$i < m$ and~$I \in \inter_\Qb(M(i))$.
An \emph{$M$-valuation $(\vec{R}, \vec{S})$} is a tuple of finite sets $R_j \subseteq \Qb$ for each~$j < n$
and~$S_{i,I} \subseteq I$ for each~$i < m$ and~$I \in \inter_\Qb(M(i))$.
The $M$-valuation~$(\vec{R}, \vec{S})$ is of type~$\vec{T}$ for some~$M$-type~$T_0, \dots, T_{n-1}$
if moreover $R_j \subseteq T_j$ for each~$j < n$.
The $M$-valuation $(\vec{R}, \vec{S})$ \emph{satisfies} $\varphi$ if $\varphi(\vec{R}, \vec{S})$ holds.
\end{definition}

Given some valuation $(\vec{R}, \vec{S})$ and some integer $s$, we write $(\vec{R}, \vec{S}) > s$
to say that for every $x \in (\bigcup \vec{R}) \cup (\bigcup \vec{S})$, $x > s$. 
Following the terminology of~\cite{Lerman2013Separating}, we define 
the notion of essentiality for a formula (an abstract requirement),
which corresponds to the idea that there is room for diagonalization
since the formula is satisfied for arbitrarily far valuations.

\begin{definition}[Essential formula]
Given an $m$-by-$n$ matrix~$M$, and $M$-formula $\varphi$
is \emph{essential} if for every $s \in \omega$,
there is an $M$-type $\vec{T}$ and an $M$-valuation $(\vec{R}, \vec{S}) > s$ of type~$\vec{T}$
such that $\varphi(\vec{R}, \vec{S})$ holds.
\end{definition}

The notion of $\er$-fairness is defined accordingly. If some formula
is essential, that is, gives enough room for diagonalization, then there is
an actual valuation which will diagonalize against the~$\erp$-instance.

\index{er-fairness@$\er$-fairness}
\begin{definition}[$\er$-fairness]
Fix two sets $A_0, A_1 \subseteq \Qb$.
Given an $m$-by-$n$ matrix $M$, an $M$-valuation~$(\vec{R}, \vec{S})$ 
\emph{diagonalizes} against $A_0, A_1$
if $\bigcup \vec{R} \subseteq A_1$ and for every $i < m$, there is some~$I \in \inter_\Qb(M(i))$ 
such that $S_{i,I} \subseteq A_0$.
A set~$X$ is \emph{$\er$-fair} for~$A_0, A_1$ if for every $m, n \in \omega$, every $m$-by-$n$ matrix~$M$
and every $\Sigma^{0,X}_1$ essential $M$-formula, there is an $M$-valuation $(\vec{R}, \vec{S})$ diagonalizing against $A_0, A_1$ 
such that $\varphi(\vec{R}, \vec{S})$ holds.
\end{definition}

Of course, if $Y \leq_T X$, then every $\Sigma^{0,Y}_1$ formula is $\Sigma^{0,X}_1$.
As an immediate consequence, if $X$ is $\er$-fair for some $A_0, A_1$ and $Y \leq_T X$, then $Y$ is $\er$-fair for $A_0, A_1$.

Now we have introduced the necessary terminology, we create a non-effective
instance of~$\erps$ which will serve as a bootstrap for $\er$-fairness preservation.

\begin{lemma}\label{lem:partition-emptyset-er-fair}
For every set~$C$, there exists a $\Delta^{0,C}_2$ partition $A_0 \cup A_1 = \Qb$ such that
$C$ is $\er$-fair for~$A_0, A_1$.
\end{lemma}
\begin{proof}
The proof is done by a no-injury priority construction.
Let $M_0, M_1, \dots$ be an enumeration of all $m$-by-$n$ matrices
and $\varphi_0, \varphi_1, \dots$ be an effective enumeration of all $\Sigma^{0,C}_1$ $M_k$-formulas for every $m, n \in \omega$.
We want to satisfy the following requirements for each pair of integers~$e,k$.

\begin{quote}
$\Rcal_{e,k}$: If the~$M_k$-formula $\varphi_e$ is essential, then $\varphi_e(\vec{R}, \vec{S})$ holds
for some $M_k$-valuation $(\vec{R}, \vec{S})$ diagonalizing against $A_0, A_1$.
\end{quote}

The requirements are ordered via the standard pairing function $\tuple{\cdot, \cdot}$.
The sets $A_0$ and $A_1$ are constructed by a $C'$-computable list of 
finite approximations $A_{i,0} \subseteq A_{i,1} \subseteq \dots$
such that all elements added to~$A_{i,s+1}$ from~$A_{i,s}$
are strictly greater than the maximum of~$A_{i,s}$ (in the~$\Nb$ order) for each~$i < 2$. 
We then let $A_i = \bigcup_s A_{i,s}$ which will be a~$\Delta^{0,C}_2$ set.
At stage 0, set $A_{0,0} = A_{1,0} = \emptyset$. Suppose that at stage $s$,
we have defined two disjoint finite sets $A_{0,s}$ and $A_{1,s}$ such that
\begin{itemize}
	\item[(i)] $A_{0,s} \cup A_{1,s} = [0,b]_\Nb$ for some integer $b \geq s$
	\item[(ii)] $\Rcal_{e',k'}$ is satisfied for every $\tuple{e',k'} < s$
\end{itemize}
Let $\Rcal_{e,k}$ be the requirement such that $\tuple{e,k} = s$.
Decide $C'$-computably whether there is some
$M$-type~$\vec{T}$ and some $M_k$-valuation $V=(\vec{R}, \vec{S}) > b$ of type~$\vec{T}$
such that $\varphi_e(V)$ holds. If so, $C$-effectively fetch~$\vec{T} = T_0, \dots, T_{n-1}$
and such a $(\vec{R}, \vec{S}) > b$. Let $d$ be an upper bound (in the~$\Nb$ order) on the rationals in $(\vec{R}, \vec{S})$.
By Lemma~\ref{lem:mtype-interval-disjoint}, for each~$i < m$,
there is some~$J_i \in \inter_\Qb(M(i))$ such that
\[
(\bigcup_{j < n} T_j) \cap (\bigcup_{i < m} J_i) = \emptyset
\]
Set $A_{0,s+1} = A_{0,s} \bigcup_{i < m} J_i \cap (b,d]_\Nb$ and $A_{1,s+1} = [0,d]_\Nb \setminus A_{0,s+1}$.
This way, $A_{0,s+1} \cup A_{1,s+1} = [0, d]_\Nb$.
By the previous equation, $\bigcup_{j < n} T_j \cap (b, d]_\Nb \subseteq [0,d]_\Nb \setminus A_{0,s+1}$ 
and the requirement $\Rcal_{e,k}$ is satisfied.
If no such $M_k$-valuation is found, the requirement $\Rcal_{e,k}$ is vacuously satisfied.
Set $A_{0,s+1} = A_{0,s} \cup \{b\}$ and $A_{1,s+1} = A_{1,s}$.
This way, $A_{0,s+1} \cup A_{1,s+1} = [0, b+1]_\Nb$.
In any case, go to the next stage. This finishes the construction.
\end{proof}

\begin{lemma}\label{lem:er-fair-not-solution}
If~$X$ is $\er$-fair for some sets~$A_0, A_1 \subseteq \Qb$,
then $X$ computes neither an infinite subset of~$A_0$,
nor a dense subset of~$A_1$.
\end{lemma}
\begin{proof}
Since $\er$-fairness is downward-closed under the Turing reducibility,
it suffices to prove that if~$X$ is infinite and $\er$-fair for~$A_0, A_1$,
then it intersects both~$A_0$ and~$A_1$.

We first prove that~$X$ intersects~$A_1$.
Let~$M$ be the 0-by-1 matrix and~$\varphi(U)$ be the $\Sigma^{0,X}_1$
$M$-formula which holds if~$U \cap X \neq \emptyset$.
The only~$M$-type is~$\Qb$ and since $X$ is infinite, $\varphi$ is essential.
By $\er$-fairness of~$X$, there is an $M$-valuation~$R$ diagonalizing against~$A_0, A_1$
such that~$\varphi(R)$ holds. By definition of diagonalization, $R \subseteq A_1$.
Since~$R \cap X \neq \emptyset$, this shows that~$X \cap A_1 \neq \emptyset$.

We now prove that~$X$ interects~$A_0$.
Let~$M$ be the~1-by-0 matrix and~$\varphi(V)$ be the~$\Sigma^{0,X}_1$
$M$-formula which holds if~$V \cap X \neq \emptyset$.
The $M$-formula $\varphi$ is essential since $X$ is infinite.
By $\er$-fairness of~$X$, there is an $M$-valuation~$S$ diagonalizing against~$A_0, A_1$
such that~$\varphi(S)$ holds. By definition of diagonalization, $S \subseteq A_0$.
Since~$R \cap X \neq \emptyset$, this shows that~$X \cap A_0 \neq \emptyset$.
\end{proof}

Note that we did not used the fact that~$X$ is dense to make it intersect $A_0$.
Density will be useful in the proof of Theorem~\ref{thm:rt22-not-er-fairness}.

\begin{definition}
A \emph{Scott set} is a set~$\Scal \subseteq 2^{\omega}$ such that
\begin{itemize}
	\item[(i)] $(\forall X \in \Scal)(\forall Y \leq_T X)[Y \in \Scal]$
	\item[(ii)] $(\forall X, Y \in \Scal)[X \oplus Y \in \Scal]$
	\item[(iii)] Every infinite, binary tree in $\Scal$ has an infinite path in~$\Scal$.
\end{itemize}
\end{definition}

\begin{theorem}\label{thm:rt22-not-er-fairness}
Let~$A_0, A_1 \subseteq \Qb$ and~$\Scal$ be a Scott set whose members are all $\er$-fair for~$A_0, A_1$.
For every set~$C \in \Scal$, every~$C$-computable coloring~$f : [\omega]^2 \to k$,
there is an infinite $f$-homogeneous set~$H$ such that~$H \oplus C$ computes neither
an infinite subset of~$A_0$, nor a dense subset of~$A_1$.
\end{theorem}
\begin{proof}
The proof is done by induction over the number of colors~$k$.
The case~$k = 1$ is ensured by Lemma~\ref{lem:er-fair-not-solution}.
Fix a set~$C \in \Scal$ and let $f : [\omega]^2 \to k$ be a $C$-computable coloring.
If~$f$ has an infinite $f$-thin set $H \in \Scal$, that is, an infinite set over which $f$
avoids at least one color, then $H \oplus C$
computes a coloring~$g : [\omega]^2 \to k-1$ such that every infinite $g$-homogeneneous
set computes relative to~$H \oplus C$ an infinite $f$-homogeneous set. Since~$H \oplus C \in \Scal$,
by induction hypothesis, there is an infinite $g$-homogeneous set~$H_1$
such that~$H_1 \oplus H \oplus C$ computes neither an infinite subset of~$A_0$, nor a dense subset of~$A_1$.
So suppose that $f$ has no infinite $f$-thin set in~$\Scal$.

We construct $k$ infinite sets~$G_0, \dots, G_{k-1}$.
We need therefore to satisfy the following requirements for each~$p \in \omega$.
\[
  \Ncal_p : \hspace{20pt} (\exists q_0 > p)[q_0 \in G_0] 
		\hspace{20pt} \wedge \dots \wedge \hspace{20pt} 
	(\exists q_{k-1} > p)[q_{k-1} \in G_{k-1}] 
\]
Furthermore, we want to ensure that one of the~$G$'s computes
neither an infinite subset of~$A_0$, nor a dense subset of~$A_1$.
To do this, we will satisfy the following requirements
for every $k$-tuple of integers $e_0, \dots, e_{k-1}$.
\[
  \Qcal_{\vec{e}} : \hspace{20pt} 
		\Rcal_{e_0}^{G_0} \hspace{20pt} \vee \dots \vee \hspace{20pt} \Rcal_{e_{k-1}}^{G_{k-1}}
\]
where $\Rcal_e^H$ holds if $W^{H \oplus C}_e$ is neither an infinite subset of~$A_0$, nor a dense subset of~$A_1$.

We construct our sets $G_0, \dots, G_{k-1}$ by forcing. 
Our conditions are variants of Mathias conditions~$(F_0, \dots, F_{k-1}, X)$
such that~$X$ is an infinite set in~$\Scal$ and the following property holds:
\begin{itemize}
	\item[(P)] $(\forall i < k)(\forall x \in X)[F_i \cup \{x\} \mbox{ is } f\mbox{-homogeneous with color } i]$
\end{itemize}
A condition~$d = (E_0, \dots, E_{k-1}, Y)$ \emph{extends} $c = (F_0, \dots, F_{k-1}, X)$
if $(E_i, Y)$ Mathias extends $(F_i, X)$ for every~$i < k$.
We now prove the progress lemma, stating that we can force the~$G$'s to be infinite.
This is where we use the fact that there is no infinite $f$-thin set in~$\Scal$.

\begin{lemma}\label{lem:rt12-er-fairness-progress}
For every condition~$c = (F_0, \dots, F_{k-1}, X)$, every $i < k$ and every~$p \in \omega$
there is some extension~$d = (E_0, \dots, E_{k-1}, Y)$ such that~$E_k \cap (p,+\infty)_\Nb \neq \emptyset$.
\end{lemma}
\begin{proof}
Fix~$c$, $i$ and~$p$. If for every~$x \in X \cap (p,+\infty)_\Nb$
and almost every~$y \in X$, $f(x,y) \neq i$, then $X$ computes an infinite
$f$-thin set, contradicting our hypothesis. Therefore, there is some~$x \in X \cap (p,+\infty)_\Nb$
such that~$f(x, y) = i$ for infinitely many~$y \in X$.
Let~$Y$ be the collection of such $y$'s. The condition~$(F_0, \dots, F_{i-1}, F \cup \{x\}, F_{i+1}, \dots, F_k, Y)$
is the desired extension.
\end{proof}

We now prove the core lemma stating that we can satisfy each $\Qcal$-requirement.
A condition~$c$ \emph{forces} a requirement~$\Qcal$
if $\Qcal$ is holds for every set~$G$ satisfying~$c$.

\begin{lemma}\label{lem:rt12-er-fairness-forcing}
For every condition~$c = (F_0, \dots, F_{k-1}, X)$ and every $k$-tuple of indices~$\vec{e}$,
there is an extension~$d = (E_0, \dots, E_{k-1}, Y)$ forcing~$\Qcal_{\vec{e}}$.
\end{lemma}
\begin{proof}
We can assume that~$W^{F_i \oplus C}_{e_i}$ has already outputted at least $k$ elements
and is either included in~$A_0$ or in~$A_1$ for each~$i < k$.
Indeed, if~$c$ has no such extension, then
$c$ forces $W^{G_i \oplus C}_{e_i}$ to be finite or not to be a valid solution for 
some~$i< k$ 
and therefore forces~$\Qcal_{\vec{e}}$. For each~$i < k$, we associate the label~$\ell_i < 2$ and the number~$p_i$
such that~$W^{F_i \oplus C}_{e_i}$ is the~$(p_i+1)$th set included in~$A_{\ell_i}$.

Let~$n$ be the number of sets~$W^{F_i \oplus C}_{e_i}$ which are included in~$A_0$,
and let~$M$ be the~$(k-n)$-by-$n$ matrix such that the~$j$th row is composed
of the~$n$ first elements already outputted by the set $W^{F_i \oplus C}_{e_i}$
where~$p_i = j$ and~$\ell_i = 1$. In other words, $M(j)$ are the $n$ first elements 
outputted by the $j$th set $W^{F_i \oplus C}_{e_i}$ included in~$A_1$.

Let~$\varphi(\vec{U}, \vec{V})$ be the $\Sigma^{0,X \oplus C}_1$ formula
which holds if for every $k$-partition $Z_0 \cup \dots \cup Z_{k-1} = X$,
there is some~$i < k$ and some finite set~$E \subseteq Z_i$
which is $f$-homogeneous with color~$i$ and such that either~$\ell_i = 0$
and $W^{(F_i \cup E) \oplus C}_{e_i} \cap U_{p_i} \neq \emptyset$,
or $\ell_i = 1$ and~$W^{(F_i \cup E) \oplus C}_{e_i} \cap V_{p_i, I} \neq \emptyset$ for each~$I \in \inter_\Qb(M(p_i))$.
We have two cases.

In the first case, $\varphi(\vec{U}, \vec{V})$ is essential.
Since $X \oplus C$ is $\er$-fair for~$A_0, A_1$, there is 
an $M$-valuation~$(\vec{R}, \vec{S})$ diagonalizing against~$A_0, A_1$
such that~$\varphi(\vec{R}, \vec{S})$ holds.
By compactness and definition of diagonalization against~$A_0, A_1$, 
there is a finite subset~$D \subset X$ such that
for every $k$-partition $D_0 \cup \dots \cup D_{k-1} = D$,
there is some $i < k$ and some finite set~$E \subseteq D_i$
which is $f$-homogeneous with color~$i$ and such that either~$\ell_i = 0$
and $W^{(F_i \cup E) \oplus C}_{e_i} \cap A_1 \neq \emptyset$,
or $\ell_i = 1$ and~$W^{(F_i \cup E) \oplus C}_{e_i} \cap A_0 \neq \emptyset$.

Each~$y \in X \setminus D$ induces a $k$-partition~$D_0 \cup \dots \cup D_{k-1}$ of~$D$
by setting~$D_i = \{ x \in D : f(x, y) = i \}$. Since there are finitely many possible 
$k$-partitions of~$D$, there is a $k$-partition~$D_0 \cup \dots \cup D_{k-1} = D$
and an infinite $X$-computable set~$Y \subseteq X$ such that
\[
(\forall i < k)(\forall x \in D_i)(\forall y \in Y)[f(x,y) = i]
\]
We furthermore assume that~$min(Y)$ is larger than the use of the computations.
Let~$i < k$ and~$E \subseteq D_i$ be the~$f$-homogeneous set with color~$i$
such that either~$\ell_i = 0$ and $W^{(F_i \cup E) \oplus C}_{e_i} \cap A_1 \neq \emptyset$,
or $\ell_i = 1$ and~$W^{(F_i \cup E) \oplus C}_{e_i} \cap A_0 \neq \emptyset$.
The condition~$(F_0, \dots, F_{i-1}, F_i \cup E, F_{i+1}, \dots, F_{k-1}, Y)$
is an extension of~$c$ forcing $\Qcal_{\vec{e}}$ by the $i$th side.

In the second case, there is some threshold~$s \in \omega$
such that for every~$M$-type~$\vec{T}$, there is no $M$-valuation~$(\vec{R}, \vec{S}) > s$
of type~$\vec{T}$ such that~$\varphi(\vec{R}, \vec{S})$ holds.
By compactness, it follows that 
for every $M$-type~$\vec{T}$, the $\Pi^{0,X \oplus C}_1$ class~$\Ccal_{\vec{T}}$
of all $k$-partitions $Z_0 \cup \dots \cup Z_{k-1} = X$
such that for every~$i < k$ and every finite set~$E \subseteq Z_i$
which is $f$-homogeneous with color~$i$, either~$\ell_i = 0$
and $W^{(F_i \cup E) \oplus C}_{e_i} \cap T_{p_i} \cap (s,+\infty)_\Nb = \emptyset$,
or $\ell_i = 1$ and~$W^{(F_i \cup E) \oplus C}_{e_i} \cap I \cap (s,+\infty)_\Nb = \emptyset$
for some~$I \in \inter_\Qb(M(p_i))$ is non-empty.
Since $\Scal$ is a Scott set, for each~$M$-type $\vec{T}$, there is a 
$k$-partition~$\vec{Z}^{\vec{T}}\in \Ccal_{\vec{T}}$
such that~$\bigoplus_{\vec{T}} \vec{Z}^{\vec{T}} \oplus X \oplus C \in \Scal$.

If there is some~$M$-type~$\vec{T}$ and some~$i < k$ such that $\ell_i = 1$ 
and~$Z^{\vec{T}}_i$ is infinite, then the condition~$(F_0, \dots, F_{k-1}, Z^{\vec{T}}_i)$
extends $X$ and forces~$W^{G_i \oplus C}_{e_i}$ not to be dense. 
So suppose that it is not the case. Let~$Y \in \Scal$
be an infinite subset of~$X$ such that for each~$M$-type~$\vec{T}$,
there is some~$i < k$ such that~$Y \subseteq Z_i^{\vec{T}}$. 
Note that by the previous assumption, $\ell_i = 0$ for every such~$i$.
We claim that the condition~$(F_0, \dots, F_{k-1}, Y)$ forces $W^{G_i \oplus C}_{e_i}$ to be finite
for some~$i < k$ such that~$\ell_i = 0$. Suppose for the sake of contradiction that there are some
rationals~$x_0, \dots, x_{n-1} > s$ such that $x_{p_i} \in W^{G_i \oplus C}_{e_i}$ for each~$i < k$ where~$\ell_i = 0$.
Since~$x_0, \dots, x_{n-1} > s$, $x_0, \dots, x_{n-1} \in \inter_\Qb(M)$. 
Therefore, by Lemma~\ref{lem:tuple-has-mtype}, let~$\vec{T}$ be the unique $M$-type
such that~$x_j \in T_j$ for each~$j < n$. By assumption, there is some~$i < k$
such that~$Y \subseteq Z^{\vec{T}}_i$ and~$\ell_i = 0$. By definition of~$Z^{\vec{T}}_i$,
$W^{G_i \oplus C}_{e_i} \cap T_{p_i} \cap (s,+\infty)_\Nb = \emptyset$,
contradicting~$x_{p_i} \in W^{G_i \oplus C}_{e_i}$.
\end{proof}

Using Lemma~\ref{lem:rt12-er-fairness-progress} and Lemma~\ref{lem:rt12-er-fairness-forcing}, define an infinite descending sequence 
of conditions~$c_0 = (\emptyset, \dots, \emptyset, \omega) \geq c_1 \geq \dots$
such that for each~$s \in \omega$
\begin{itemize}
	\item[(i)] $|F_{i,s}| \geq s$ for each~$i < k$
	\item[(ii)] $c_{s+1}$ forces~$\Qcal_{\vec{e}}$ if~$s = \tuple{e_0, \dots, e_{k-1}}$
\end{itemize}
where~$c_s = (F_{0,s}, \dots, F_{k-1,s}, X_s)$. Let~$G_i = \bigcup_s F_{i,s}$ for each~$i < k$.
The~$G$'s are all infinite by (i) and $G_i$ does not compute an $\erps$-solution to the~$A$'s for some~$i < k$ by (ii).
This finishes the proof of Theorem~\ref{thm:rt22-not-er-fairness}.
\end{proof}

We are now ready to prove the main theorem.

\begin{proof}[Proof of Theorem~\ref{thm:er22-not-computably-reduces}]
By the low basis theorem~\cite{Jockusch197201}, there is a low set $P$ of PA degree.
By Scott~\cite{Scott1962Algebras}, every PA degree bounds a Scott set. 
Let~$\Scal$ be a Scott set such that~$X \leq_T P$ for every~$X \in \Scal$.
By Lemma~\ref{lem:partition-emptyset-er-fair}, there is a $\Delta^{0,P}_2$ (hence $\Delta^0_2$)
partition $A_0 \cup A_1 = \Qb$ such that $P$ is $\er$-fair for~$A_0, A_1$. In particular,
every set~$X \in \Scal$ is $\er$-fair for~$A_0, A_1$ since $\er$-fairness is downward-closed under the Turing reducibility.

By Shoenfield's limit lemma~\cite{Shoenfield1959degrees}, there is a computable function $h : [\Qb]^2 \to 2$ such 
that for each~$x \in \Qb$, $\lim_s h(x, s)$ exists and $x \in A_{\lim_s h(x, s)}$.
Note that for every infinite set $D$ 0-homogeneous for~$h$, $D \subseteq A_0$,
and for every dense set $D$ 1-homogeneous for~$h$, $D \subseteq A_1$.

Fix a computable $\rt^2_{<\infty}$-instance~$f : [\omega]^2 \to k$. In particular, $f \in 
\Scal$.
By Theorem~\ref{thm:rt22-not-er-fairness}, there is an infinite $f$-homogeneous set~$H$
such that $H$ computes neither an infinite subset of~$A_0$, nor a dense subset of~$A_1$.
Therefore, $H$ computes no $\erp$-solution to~$h$.
\end{proof}

\section{Discussion and open questions}

Both~$\tto^2_2$ and~$\erp$ lie between the arithmetic comprehension axiom and~$\rt^2_2$,
but more than that, they share a \emph{disjoint extension commitment}.
Let us try to explain this informal notion with a case analysis.

Suppose we want to construct a computable $\rt^1_2$-instance~$f : \Nb \to 2$ which diagonalizes
against two opponents~$W^f_0$ and~$W^f_1$. After some finite amount of time, 
each opponent~$W^f_i$ will have outputted
a finite approximation of a solution to~$f$, that is, a 
finite $f$-homogeneous set~$F_i$.
The two opponents share a common strategy. $W^f_0$ tries to build an infinite $f$-homogeneous set~$H_0$
for color~0, and~$W^f_1$ tries to build an infinite $f$-homogeneous set~$H_1$ for color~1.
It is therefore difficult to defeat both opponents at the same time, since
if from now on we set~$f(x) = 1$, $W^f_1$ will succeed in extending~$F_1$ to an infinite $f$-homogenenous set,
and if we set $f(x) = 0$, $W^f_0$ will succeed with its dual strategy.

Consider now the same situation, where we want to construct a computable~$\tto^1_2$-instance $f : 2^{<\Nb} \to 2$.
After some time, the opponent~$W^f_0$ will have outputted a finite tree~$S_0 \cong 2^{<b}$ which is
$f$-homogeneous for color~$0$, and the opponent~$W^f_1$ will have done the same with a finite tree~$S_1 \cong 2^{<b}$
$f$-homogeneous for color~$1$.
The main difference with the~$\rt^1_2$ case is that each opponent will \emph{commit to extend}
each leaf of his finite tree~$S_i$ into an infinite tree isomorphic to~$2^{<\Nb}$.
In particular, for each tree~$S_i$, the sets~$X_\sigma$ of nodes extending the leaf~$\sigma \in S_i$
are \emph{pairwise disjoint}. Therefore, each opponent commits to extend its partial solution
to disjoint sets. Moreover, by asking~$b$ to be large enough, each opponent will commit
to extend enough pairwise disjoint sets so that we can choose two of them for each opponent
and operate the diagonalization without any conflict.

This combinatorial property works in the same way for~$\ers$-instances.
Indeed, in this case, each opponent will commit to extend its partial solution
to pairwise disjoint intervals due to the density requirement of an~$\ers$-solution.
Since the combinatorial arguments of the Erd\H{o}s-Rado theorem and the tree theorem for pairs
are very similar, one way wonder whether they are equivalent in reverse mathematics.

\begin{question}
How do $\erp$ and~$\tto^2_2$ compare over~$\rca$?
\end{question}

The failure of Seetapun's argument for~$\erp$ comes 
from this disjoint extension commitment feature. 
In particular, there it is hard to find a forcing notion
for~$\erp$ whose conditions are extendible.

\begin{question}
Does any of~$\tto^2_2$ and $\erp$ imply~$\aca$ over~$\rca$?
\end{question}

$\ers$ and~$\tto^1$ have the same state of the art due to their common
combinatorial flavor. However, when looking at their statements for pairs,
$\erp$ and~$\tto^2_2$ have a fundamental difference: $\erp$ has only
a half disjoint extension commitment feature. This weaker property
prevents one from separating~$\rt^2_2$ from~$\erp$ over~$\rca$
by adapting the argument of~$\tto^2_2$.

\begin{question}
Does~$\rt^2_2$ imply $\erp$ over~$\rca$?
\end{question}

\chapter{Controlling iterated jumps}\label{chap-controlling-iterated-jumps}

Effective forcing is a very powerful tool in the computational analysis
of mathematical statements. In this framework, lowness is achieved by
deciding formulas during the forcing argument, while ensuring that the whole construction
remains effective. Thus, the definitional strength of the forcing relation
is very sensitive in effective forcing.

Among the hierarchies of combinatorial principles,
namely, Ramsey's theorem~\cite{Jockusch1972Ramseys,Seetapun1995strength,Cholak2001strength}, 
the rainbow Ramsey theorem~\cite{Csima2009strength,Wang2014Cohesive,Patey2015Somewhere}, 
the free sets and thin set theorems~\cite{Cholak2001Free,Wang2014Some},
only the hierarchy of Ramsey's theorem is known to collapse within the framework of reverse mathematics.
The above-mentioned hierarchies satisfy the lower bounds of Jockusch~\cite{Jockusch1972Ramseys}, 
that is, there exists a computable instance at every level~$n \geq 2$ with no $\Sigma^0_n$ solution.
Thus, a possible strategy for proving that a hierarchy is strict consists
of showing the existence for every computable instance at level~$n$ of a low${}_n$ solution.

The computable analysis of combinatorial principles often uses Mathias forcing, whose forcing
relation is known to be of higher definitional strength than the formula it forces~\cite{Cholak2014Generics}.
Therefore there is a need for new forcing notions with a better-behaving forcing relation.
In this chapter, we present a new forcing argument enabling one to control iterated jumps
of solutions to Ramsey-type theorems.
Therefore it can be seen as a step toward resolving the strictness question
of the combinatorial hierarchies.

\section{Preservation of the arithmetic hierarchy}

The notion of preservation of the arithmetic hierarchy has
been introduced by Wang in~\cite{Wang2014Definability}, in the context of a new analysis
of principles in reverse mathematics in terms of their definitional strength.

\index{preservation!of definitions}
\index{preservation!of the arithmetic hierarchy}
\begin{definition}[Preservation of definitions]\ 
\begin{itemize}
	\item[1.]
A set $Y$ \emph{preserves $\Xi$-definitions} (relative to $X$) for $\Xi$ among $\Delta^0_{n+1}, \Pi^0_n, \Sigma^0_n$ where $n > 0$, if every properly $\Xi$ (relative to $X$) set is properly $\Xi$ relative to $Y$ ($X \oplus Y$). $Y$ \emph{preserves the arithmetic hierarchy} (relative to $X$) if $Y$ preserves $\Xi$-definitions (relative to $X$) for all $\Xi$ among $\Delta^0_{n+1}, \Pi^0_n, \Sigma^0_n$ where $n > 0$.

	\item[2.] A $\Pi^1_2$ statement~$\Psf$ \emph{admits preservation of $\Xi$-definitions} if for each set~$Z$,
every~$Z$-computable $\Psf$-instance admits a solution preserving $\Xi$-definitions relative to $Z$. 
$\Psf$ \emph{admits preservation of the arithmetic hierarchy} if for each set $Z$,
every~$Z$-computable $\Psf$-instance admits a solution preserving the arithmetic hierarchy relative to $Z$.
\end{itemize}
\end{definition}

The preservation of the arithmetic hierarchy seems closely related
to the problem of controlling iterated jumps of solutions.
Indeed, a proof of such a preservation
usually consists of noticing that the forcing relation has the same definitional strength
as the formula it proves and then derive a diagonalization. 
See Lemma 3.14 in~\cite{Wang2014Definability} for a case-in-point.
Wang proved in~\cite{Wang2014Definability} that weak König's lemma,
the rainbow Ramsey theorem for pairs and the atomic model theorem
admit preservation of the arithmetic hierarchy. He conjectured
that this is also the case for cohesiveness and the Erd\H{o}s Moser theorem.
We prove the two conjectures.

\section{An effective forcing for cohesiveness}

In this section, we design a forcing notion enabling us to prove the following theorem.

\begin{theorem}\label{thm:coh-preservation-arithmetic-hierarchy}
$\coh$ admits preservation of the arithmetic hierarchy.
\end{theorem}

Before proving Theorem~\ref{thm:coh-preservation-arithmetic-hierarchy},
we state an immediate corollary.

\begin{corollary}
There exists a cohesive set preserving the arithmetic hierarchy.
\end{corollary}
\begin{proof}
Jockusch~\cite{Jockusch1972Degreesa} proved that every PA degree computes a sequence of sets
containing, among others, all the computable sets.
Wang proved in~\cite{Wang2014Definability} that $\wkl$ preserves the arithmetic hierarchy. 
Therefore there exists a uniform sequence of sets $\vec{R}$ containing all the computable sets 
and preserving the arithmetic hierarchy. By Theorem~\ref{thm:coh-preservation-arithmetic-hierarchy}
relativized to $\vec{R}$, there exists an infinite $\vec{R}$-cohesive set $C$ preserving the arithmetic hierarchy relative to~$\vec{R}$.
In particular $C$ is r-cohesive and preserves the arithmetic hierarchy. By~\cite{Jockusch1993cohesive},
the degrees of r-cohesive and cohesive sets coincide. Therefore $C$ computes a cohesive set which preserves the arithmetic hierarchy.
\end{proof}

Given a uniformly computable sequence of sets~$R_0, R_1, \dots$,
the construction of an $\vec{R}$-cohesive set is usually done
with computable Mathias forcing, that is, using conditions~$(F, X)$
in which~$X$ is computable.
The construction starts with~$(\emptyset, \omega)$ and interleaves two kinds of steps.
Given some condition~$(F,X)$,
\begin{itemize}
	\item[(S1)] the \emph{extension} step consists in taking an element $x$ from $X$ and adding it to~$F$,
	therefore forming the extension $(F \cup \{x\}, X \setminus [0,x])$;
	\item[(S2)] the \emph{cohesiveness} step consists in deciding which one of $X \cap R_i$
	and $X \cap \overline{R}_i$ is infinite, and taking the chosen one as the new reservoir.
\end{itemize} 

As presented in Section~\ref{sect:introduction-forcing-mathias-forcing}, 
computable Mathias forcing has a forcing relation with good definitional properties to decide the first jump, but not iterated jumps.
Indeed, given a computable Mathias condition~$c = (F, X)$ and a $\Sigma^0_1$ formula~$(\exists x)\varphi(G, x)$,
one can $\emptyset'$-effectively decide whether there is an extension~$d$ forcing~$(\exists x)\varphi(G, x)$ by asking the following question:
\begin{quote}
Is there an extension~$d = (E, Y) \leq c$ and some~$n \in \omega$ such that~$\varphi(E, n)$ holds?
\end{quote}
If there is such an extension, then we can choose it to be a \emph{finite extension}, that is, such that~$Y =^* X$.
Therefore, the question is $\Sigma^{0,X}_1$.
Consider now a $\Pi^0_2$ formula~$(\forall x)(\exists y)\varphi(G, x, y)$. The question becomes
\begin{quote}
For every extension~$d \leq c$ and every~$m \in \omega$, is there some extension~$e = (E, Y) \leq d$
and some~$n \in \omega$ such that~$\varphi(E, m, n)$ holds?
\end{quote}
In this case, the extension~$d$ can be arbitrary and therefore the question cannot
be presented in a $\Pi^0_2$ way. In particular, the formula ``$Y$ is an infinite subset of~$X$'' is definitionally complex. 
In general, deciding iterated jumps of a generic set requires to be able to talk about the future
of a given condition, and in particular to describe by simple means the formula ``$d$ is a valid condition''
and the formula ``$d$ is an extension of~$c$''.

Thankfully, in the case of cohesiveness, we do not need the full generality of 
the computable Mathias forcing. As noted in Section~\ref{sect:dominating-erdos-moser-dominating-cohesive}, the reservoirs
have a very special shape. Indeed, after the first application of stage~(S2), 
the set $X$ is, up to finite changes, of the form $\omega \cap R_0$
or $\omega \cap \overline{R_0}$. After the second application of (S2), it is in one of the following forms: $\omega \cap R_0 \cap R_1$,
$\omega \cap R_0 \cap \overline{R}_1$, $\omega \cap \overline{R}_0 \cap R_1$,
$\omega \cap \overline{R}_0 \cap \overline{R}_1$, and so on.
More generally, after $n$ applications of (S2), a condition~$c = (F, X)$ is characterized
by a pair~$(F, \sigma)$ where $\sigma$ is a string of length~$n$ representing the choices made during (S2).

Even within this restricted partial order, the decision of the~$\Pi^0_2$ formula remains too complicated 
sinces it requires to decide whether $R_\sigma$ is infinite.
However, notice that the~$\sigma$'s such that~$R_\sigma$ is infinite are exactly the initial
segments of the $\Pi^{0, \emptyset'}_1$ class~$\Ccal(\vec{R})$ defined in Section~\ref{sect:strength-ramsey-cohesiveness-konig-lemma}.
We can therefore use a compactness argument at the second level to decrease the definitional strength of the forcing relation,
as did Wang~\cite{Wang2014Definability} for weak K\"onig's lemma.

\subsection{The forcing notion}

In order to simplify our presentation, let us first introduce some useful notation.

\begin{notation}
Given two finite sets~$E, F$ and $\sigma \in 2^{<\omega}$,
we write $E \leq_\sigma F$ for the formula $F \subseteq G \subseteq F \cup (R_\sigma \cap (max(F), +\infty))$,
where~$R_\sigma$ is defined in Section~\ref{sect:strength-ramsey-cohesiveness-konig-lemma}. 
A tree~$T$ has \emph{stem} $\sigma \in 2^{<\omega}$
if every node of~$T$ is comparable with~$\sigma$. 
We write~$T^{[\sigma]}$ for~$\{\tau \in T : \tau \preceq \sigma \vee \tau \succeq \sigma \}$.
\end{notation}

We let~$\Tb$ denote the collection of all the infinite $\emptyset'$-primitive recursive trees~$T$
such that~$[T] \subseteq \Ccal(\vec{R})$. By $\emptyset'$-primitive recursive, we mean the class of functions
obtained by adding the characteristic function of $\emptyset'$ to the basic primitive recursive functions,
and closing under the standard primitive recursive operations.
Note that~$\Tb$ is a computable set. We are now ready to defined our partial order.

\begin{definition}
Let $\Pb$ be the partial order whose conditions are tuples $(F, \sigma, T)$
where~$F \subseteq \omega$ is a finite set, $\sigma \in 2^{<\omega}$,
and~$T \in \Tb$ with stem~$\sigma$.
A condition~$d = (E, \tau, S)$ \emph{extends} $c = (F, \sigma, T)$ (written~$d \leq c$)
if~$E \leq_\sigma F$, $\tau \succeq \sigma$ and $S \subseteq T$.
\end{definition}

Given a condition~$c = (F, \sigma, T)$, the string~$\sigma$ imposes a finite restriction on the possible extensions of the set~$F$.
The condition $c$ intuitively denotes the Mathias condition $(F, R_\sigma \cap (max(F), +\infty))$
with some additional constraints on the extensions of~$\sigma$ represented by the tree~$T$.
Accordingly, set~$G$ \emph{satisfies} $(F, \sigma, T)$ if it satisfies
the induced Mathias condition, that is, if~$F \subseteq G \subseteq F \cup (R_\sigma \cap (max(F), +\infty))$. 
We let~$\Ext(c)$ be the collection of all the extensions of~$c$.

Note that although we did not explicitely require $R_\sigma$ to be infinite,
this property holds for every condition~$(F, \sigma, T) \in \Pb$.
Indeed, since $[T] \subseteq \Ccal(\vec{R})$, then $R_\tau$ is infinite
for every extensible node~$\tau \in T$. Since $\sigma$ is a stem of~$T$,
it is extensible and therefore $R_\sigma$ is infinite.

\subsection{Preconditions and forcing $\Sigma^0_1$ ($\Pi^0_1$) formulas}

When forcing complex formulas, we need to be able to consider all possible extensions of some condition~$c$.
Checking that some $d = (E, \tau, S)$ is a valid condition extending $c$ 
requires to decide whether the tree $\emptyset'$-p.r. $S$ is infinite,
which is a $\Pi^0_2$ question.
At some point, we will need to decide a $\Sigma^0_1$ formula without having enough computational power
to check that the tree part is infinite (see clause~(ii) of Definition~\ref{def:forcing-condition}).
As the tree part of a condition is not accurate for such formulas,
we may define the corresponding forcing relation over 
a weaker notion of condition where the tree is not anymore required to be infinite.

\begin{definition}[Precondition]
A \emph{precondition} is a condition~$(F, \sigma, T)$ without the assumption that $T$ is infinite.
\end{definition}

In particular, $R_\sigma$ may be a finite set. The notion of condition extension can be generalized
to the preconditions. The set of all preconditions is computable, contrary to the set~$\Pb$.
Given a precondition~$c$, we denote by~$\Ext_1(c)$ the set of all preconditions~$(E, \tau, S)$
extending~$c$ such that~$\tau = \sigma$ and~$T = S$. Here, $T = S$ in a strong sense,
that is, the Turing indices of~$T$ and~$S$ are the same. This fact is used in clause~a) of Lemma~\ref{lem:extension-complexity}.
We let~$\Ab$ denote the collection of all the finite sets of integers. 
The set $\Ab$ represents the set of finite approximations of the generic set~$G$. 
We also fix a uniformly computable enumeration $\Ab_0 \subseteq \Ab_1 \subseteq \dots$
of finite subsets of~$\Ab$ such that~$\bigcup_s \Ab_s = \Ab$.
We denote by~$\Apx(c)$ the set $\{ E \in \Ab : (E, \sigma, T) \in \Ext_1(c) \}$.
In particular, $\Apx(c)$ is collection of all finite sets~$E$ satisfying~$c$, that is, 
$\Apx(c) = \{ E \in \Ab : E \leq_\sigma F \}$.
Last, we let~$\Apx_s(c) = \Apx(c) \cap \Ab_s$. We start by proving a few trivial statements.

\begin{lemma}\label{lem:basic-statements}
Fix a precondition $c = (F, \sigma, T)$.
\begin{itemize}
	\item[1)] If $c$ is a condition then $\Ext_1(c) \subseteq \Ext(c)$.
	\item[2)] If $c$ is a condition then $\Apx(c) = \{ E : (E, \tau, S) \in \Ext(c) \}$.
	\item[3)] If $d$ is a precondition extending $c$ then 
	$\Apx(d) \subseteq \Apx(c)$ and~$\Apx_s(d) \subseteq \Apx_s(c)$.
\end{itemize}
\end{lemma}
\begin{proof}\ 
\begin{itemize}
	\item[1)]
By definition, if $c$ is a condition, then $T$ is infinite.
If $d \in \Ext_1(c)$ then $d = (E, \sigma, T)$ for some $E \in \Apx(c)$.
As $d$ is a precondition and $T$ is infinite, $d$ is a condition.

	\item[2)] By definition,
	$\Apx(c) = \{ E : (E, \sigma, T) \in \Ext_1(c) \} \subseteq \{ E : (E, \tau, S) \in \Ext(c) \}$.
	In the other direction, fix an extension $(E, \tau, S) \in \Ext(c)$.
	By definition of an extension, $E \leq_\tau F$, so $E \leq_\sigma F$.
	Therefore $(E, \sigma, T) \in \Ext_1(c)$ and by definition of $\Apx(c)$, $E \in \Apx(c)$.

	\item[3)]
	Fix some $(E, \tau, S) \in \Ext_1(d)$. As $d$ extends $c$, $\tau \succeq \sigma$.
	By definition of an extension, $E \leq_\tau F$, so $E \leq_\sigma F$,
	hence $(E, \sigma, T) \in \Ext_1(c)$.
	Therefore $\Apx(d) = \{ E : (E, \tau, S) \in \Ext_1(d) \} \subseteq
\{ E : (E, \sigma, T) \in \Ext_1(c) \} = \Apx(c)$. For any~$s \in \omega$,
	$\Apx_s(d) = \Apx(d) \cap \Ab_s \subseteq \Apx(c) \cap \Ab_s = \Apx_s(c)$.
\end{itemize}
\end{proof}

Note that although the extension relation has been generalized to preconditions, 
$\Ext(c)$ is defined to be the set of all the \emph{conditions} extending $c$. 
In particular, if $c$ is a precondition which is not a condition, $\Ext(c) = \emptyset$,
whereas at least $c \in \Ext_1(c)$.
This is why clause 1 of Lemma~\ref{lem:basic-statements} gives the useful information
that whenever~$c$ is a true condition, so are the members of $\Ext_1(c)$.

\begin{definition}\label{def:forcing-precondition}
Fix a precondition~$c = (F, \sigma, T)$ and a $\Sigma^0_0$ formula $\varphi(G, x)$.
\begin{itemize}
	\item[(i)] $c \Vdash (\exists x)\varphi(G, x)$ iff $\varphi(F, w)$ holds for some~$w \in \omega$
	\item[(ii)] $c \Vdash (\forall x)\varphi(G, x)$ iff $\varphi(E, w)$ holds for every~$w \in \omega$
	and every set~$E \in \Apx(c)$.
\end{itemize}
\end{definition}

As explained, $\sigma$ restricts the possible extensions of the set~$F$ (see clause 3 of Lemma~\ref{lem:basic-statements}), 
so this forcing notion is stable by condition extension. The tree $T$ itself restricts the possible extensions of $\sigma$,
but has no effect of the decision of a $\Sigma^0_1$ formula (Lemma~\ref{lem:tree-no-effect-first-level}).

The following trivial lemma expresses the fact that the tree part of a precondition
has no effect in the forcing relation for a $\Sigma^0_1$ or $\Pi^0_1$ formula.

\begin{lemma}\label{lem:tree-no-effect-first-level}
Fix two preconditions $c = (F, \sigma, T)$ and $d = (F, \sigma, S)$, and some $\Sigma^0_1$ or $\Pi^0_1$
formula~$\varphi(G)$.
$$
c \Vdash \varphi(G) \hspace{10pt} \mbox{ if and only if } \hspace{10pt} d \Vdash \varphi(G)
$$
\end{lemma}
\begin{proof}
Simply notice that the tree part of the condition does not occur in the definition
of the forcing relation, and that $\Apx(c) = \Apx(d)$.
\end{proof}

As one may expect, the forcing relation
for a precondition is closed under extension.

\begin{lemma}\label{lem:forcing-extension-precondition}
Fix a precondition $c$ and a $\Sigma^0_1$ or $\Pi^0_1$ formula $\varphi(G)$.
If $c \Vdash \varphi(G)$ then for every precondition $d \leq c$, $d \Vdash \varphi(G)$.
\end{lemma}
\begin{proof}
Fix a precondition $c = (F, \sigma, T)$ such that $c \Vdash \varphi(G)$ and an extension $d = (E, \tau, S) \leq c$.
\begin{itemize}
  \item
If $\varphi \in \Sigma^0_1$ then $\varphi(G)$ can be expressed
as $(\exists x)\psi(G, x)$ where $\psi \in \Sigma^0_0$.
As $c \Vdash \varphi(G)$, then by clause (i) of Definition~\ref{def:forcing-precondition},
there exists a $w \in \omega$ such that $\psi(F, w)$ holds. 
By definition of $d \leq c$, $E \leq_\sigma F$, so $\psi(E, w)$ holds,
hence $d \Vdash \varphi(G)$.

\item
If $\varphi \in \Pi^0_1$ then $\varphi(G)$ can be expressed
as $(\forall x)\psi(G, x)$ where $\psi \in \Sigma^0_0$.
As $c \Vdash \varphi(G)$, then by clause (ii) of Definition~\ref{def:forcing-precondition},
for every $w \in \omega$ and every $H \in \Apx(c)$, $\varphi(H, w)$ holds.
By clause 3 of Lemma~\ref{lem:basic-statements}, $\Apx(d) \subseteq \Apx(c)$ so $d \Vdash \varphi(G)$.
\end{itemize}
\end{proof}

\subsection{Forcing higher formulas}

We are now able to define the forcing relation for any arithmetic formula.
The forcing relation for arbitrary arithmetic formulas is induced by the forcing 
relation for~$\Sigma^0_1$ formulas. However, the definitional strength of the resulting
relation is too high with respect to the formula it forces. We therefore
design a custom relation with better definitional properties,
and which still preserve the expected properties of a forcing relation,
that is, the density of the set of conditions forcing a formula or its negation,
and the preservation of the forced formulas under condition extension.

\begin{definition}\label{def:forcing-condition}
Let~$c = (F, \sigma, T)$ be a condition and~$\varphi(G)$ be an arithmetic formula.
\begin{itemize}
	\item[(i)] If $\varphi(G) = (\exists x)\psi(G, x)$ where~$\psi \in \Pi^0_{n+1}$ then~$c \Vdash \varphi(G)$
	iff there is a $w < |\sigma|$ such that~$c \Vdash \psi(G, w)$
	\item[(ii)] If $\varphi(G) = (\forall x)\psi(G, x)$ where~$\psi \in \Sigma^0_1$ then $c \Vdash \varphi(G)$
	iff for every~$\tau \in T$, every $E \in \Apx_{|\tau|}(c)$
	and every~$w < |\tau|$, $(E, \tau, T^{[\tau]}) \not \Vdash \neg \psi(G, w)$
	\item[(iii)] If~$\varphi(G) = \neg \psi(G, x)$ where~$\psi \in \Sigma^0_{n+3}$ then $c \Vdash \varphi(G)$
	iff $d \not \Vdash \psi(G)$ for every~$d \leq c$.
\end{itemize}
\end{definition}

Note that in clause (ii) of Definition~\ref{def:forcing-condition},
there may be some $\tau \in T$ such that $T^{[\tau]}$ is finite, hence $(E, \tau, T^{[\tau]})$
is not necessarily a condition. This is where we use the generalization of forcing of $\Sigma^0_1$ and $\Pi^0_1$
formulas to preconditions. We now prove that this relation enjoys the main properties
of a forcing relation.

\begin{lemma}\label{lem:forcing-extension}
Fix a condition $c$ and an arithmetic formula $\varphi(G)$.
If $c \Vdash \varphi(G)$ then for every condition $d \leq c$, $d \Vdash \varphi(G)$.
\end{lemma}
\begin{proof}
We prove by induction over the complexity of the formula $\varphi(G)$
that for every condition $c$, if $c \Vdash \varphi(G)$
then for every condition $d \leq c$, $d \Vdash \varphi(G)$.
Fix a condition $c = (F, \sigma, T)$ such that $c \Vdash \varphi(G)$ and an extension $d = (E, \tau, S)$.
\begin{itemize}
  \item If $\varphi \in \Sigma^0_1 \cup \Pi^0_1$ then it follows from Lemma~\ref{lem:forcing-extension-precondition}.
  
	\item If $\varphi \in \Sigma^0_{n+2}$ then $\varphi(G)$ can be expressed as $(\exists x)\psi(G, x)$ 
  where $\psi \in \Pi^0_{n+1}$.
  By clause~(i) of Definition~\ref{def:forcing-condition}, there exists a $w \in \omega$
  such that $c \Vdash \psi(G, w)$. By induction hypothesis, $d \Vdash \psi(G, w)$
  so by clause~(i) of Definition~\ref{def:forcing-condition}, $d \Vdash \varphi(G)$.

  \item If $\varphi \in \Pi^0_2$ then $\varphi(G)$ can be expressed as $(\forall x)\psi(G, x)$ where $\psi \in \Sigma^0_1$.
	By clause~(ii) of Definition~\ref{def:forcing-condition}, for every $\rho \in T$, every $w < |\rho|$,
	and every $H \in \Apx_{|\rho|}(c)$, $(H, \rho, T^{[\rho]}) \not \Vdash \neg \psi(G, w)$.
  As $S \subseteq T$ and $\Apx(d) \subseteq \Apx(c)$, for every $\rho \in S$, every $w < |\rho|$,
  and every $H \in \Apx_{|\rho|}(d)$, $(H, \rho, T^{[\rho]}) \not \Vdash \neg \psi(G, w)$.
  By Lemma~\ref{lem:tree-no-effect-first-level}, $(H, \rho, S^{[\rho]}) \not \Vdash \neg \psi(G, w)$
  hence by clause~(ii) of Definition~\ref{def:forcing-condition}, $d \Vdash \varphi(G)$.

  \item If $\varphi \in \Pi^0_{n+3}$ then $\varphi(G)$ can be expressed as $\neg \psi(G)$ where $\psi \in \Sigma^0_{n+3}$.
  By clause~(iii) of Definition~\ref{def:forcing-condition}, for every $e \in \Ext(c)$, $e \not \Vdash \psi(G)$.
  As $\Ext(d) \subseteq \Ext(c)$, for every $e \in \Ext(d)$, $e \not \Vdash \psi(G)$,
  so by clause~(iii) of Definition~\ref{def:forcing-condition}, $d \Vdash \varphi(G)$.
\end{itemize}
\end{proof}

\begin{lemma}\label{lem:forcing-dense}
For every arithmetic formula $\varphi$, the following set is dense
$$
\{c \in \Pb : c \Vdash \varphi(G) \mbox{ or } c \Vdash \neg \varphi(G) \}
$$
\end{lemma}
\begin{proof}
We prove by induction over $n > 0$ that if $\varphi$ is a $\Sigma^0_n$ ($\Pi^0_n$)
formula then the following set is dense
$$
\{c \in \Pb : c \Vdash \varphi(G) \mbox{ or } c \Vdash \neg \varphi(G) \}
$$
It suffices to prove it for the case where $\varphi$ is a $\Sigma^0_n$ formula,
as the case where $\varphi$ is a $\Pi^0_n$ formula is symmetric. Fix a condition $c = (F, \sigma, T)$.
\begin{itemize}
	\item In case $n = 1$, the formula $\varphi$ is of the form $(\exists x)\psi(G, x)$ where $\psi \in \Sigma^0_0$.
	Suppose there exist a $w \in \omega$ and a set $E \in \Apx(c)$ such that $\psi(E, w)$ holds.
	The precondition $d = (E, \sigma, T)$ is a condition extending $c$ by clause~1 of Lemma~\ref{lem:basic-statements}
	and by definition of $\Apx(c)$.
  Moreover $d \Vdash (\exists x)\psi(G, x)$ by clause~(i) of Definition~\ref{def:forcing-precondition} hence $d \Vdash \varphi(G)$.
	Suppose now that for every $w \in \omega$ and every $E \in \Apx(c)$, $\psi(E, w)$ does not hold.
	By clause~(ii) of Definition~\ref{def:forcing-precondition}, 
	$c \Vdash (\forall x)\neg \psi(G, x)$, hence $c \Vdash \neg \varphi(G)$.

	\item In case $n = 2$, the formula $\varphi$ is of the form $(\exists x)\psi(G, x)$ where
	$\psi \in \Pi^0_1$. Let
	$$
	S = \{ \tau \in T : (\forall w < |\tau|)(\forall E \in \Apx_{|\tau|}(c))(E, \tau, T^{[\tau]}) \not \Vdash \psi(G, w) \}
	$$
	The set $S$ is obviously $\emptyset'$-p.r. We prove that it is a subtree of $T$.
	Suppose that $\tau \in S$ and $\rho \preceq \tau$. 
	Fix a $w < |\rho|$ and $E \in \Apx_{|\rho|}(c)$. In particular $w < |\tau|$ and $E \in \Apx_{|\tau|}(c)$
	so $(E, \tau, T^{[\tau]}) \not \Vdash \psi(G, w)$. 
	Note that~$(E, \tau,T^{[\tau]})$ is a precondition extending~$(E, \rho, T^{[\rho]})$,
	so by the contrapositive of Lemma~\ref{lem:forcing-extension-precondition}, $(E, \rho, T^{[\rho]}) \not \Vdash \psi(G, w)$.
	Therefore $\rho \in S$. Hence $S$ is a tree, and as $S \subseteq T$, it is a subtree of $T$.
	
	If $S$ is infinite, then $d = (F, \sigma, S)$
	is an extension of $c$ such that for every $\tau \in S$, every $w < |\tau|$
	and every $E \in \Apx_{|\tau|}(c)$, $(E, \tau, T^{[\tau]}) \not \Vdash \psi(G, w)$.
	By Lemma~\ref{lem:tree-no-effect-first-level}, for every $E \in \Apx_{|\tau|}(c)$,
	$(E, \tau, S^{[\tau]}) \not \Vdash \psi(G, w)$ and by clause 3 of Lemma~\ref{lem:basic-statements}, 
	$\Apx_{|\tau|}(d) \subseteq \Apx_{|\tau|}(c)$.
	Therefore, by clause~(ii) of Definition~\ref{def:forcing-condition}, $d \Vdash (\forall x)\neg \psi(G, x)$
	so $d \Vdash \neg \varphi(G)$. If $S$ is finite, then pick some $\tau \in T \setminus S$
	such that $T^{[\tau]}$ is infinite. By choice of $\tau \in T \setminus S$, there exist a $w < |\tau|$
	and an $E \in \Apx_{|\tau|}(c)$ such that $(E, \tau, T^{[\tau]}) \Vdash \psi(G, w)$.
	$d = (E, \tau, T^{[\tau]})$ is a valid condition extending $c$
	and by clause~(i) of Definition~\ref{def:forcing-condition} $d \Vdash \varphi(G)$.

	\item In case $n > 2$, density follows from clause~(iii) of Definition~\ref{def:forcing-condition}.
\end{itemize}
\end{proof}

Any sufficiently generic filter~$\Fcal$ induces a unique generic real~$G$
defined by
$$
G = \bigcup \{ F \in \Ab : (F, \sigma, T) \in \Fcal \}
$$
The following lemma informally asserts that the forcing relation is \emph{sound} and \emph{complete}.
Sound because whenever it forces a property, then this property actually holds over the generic real~$G$.
The forcing is also complete in that every property which holds over~$G$ is forced at some point
whenever the filter is sufficiently generic.

\begin{lemma}\label{lem:coh-holds-filter}
Suppose that $\Fcal$ is a sufficiently generic filter and let~$G$ be the corresponding generic real.
Then for each arithmetic formula $\varphi(G)$,
$\varphi(G)$ holds iff $c \Vdash \varphi(G)$ for some $c \in \Fcal$. 
\end{lemma}
\begin{proof}
We prove by induction over the complexity of the arithmetic formula $\varphi(G)$ that 
$\varphi(G)$ holds iff $c \Vdash \varphi(G)$ for some $c \in \Fcal$.
Note that thanks to Lemma~\ref{lem:forcing-dense}, it suffices to prove that if $c \Vdash \varphi(G)$ for some
$c \in \Fcal$ then $\varphi(G)$ holds. Indeed, conversely if $\varphi(G)$ holds,
then by genericity of $G$ either $c \Vdash \varphi(G)$ or $c \Vdash \neg \varphi(G)$ for some~$c \in \Fcal$,
but if $c \Vdash \neg \varphi(G)$ then $\neg \varphi(G)$ holds, contradicting 
the hypothesis. So $c \Vdash \varphi(G)$.

We proceed by case analysis on the formula~$\varphi$.
Note that in the above argument, the converse of the~$\Sigma$ case is proved assuming
the~$\Pi$ case. However, in our proof, we use the converse of the~$\Sigma^0_{n+3}$
case to prove the $\Pi^0_{n+3}$ case. We need therefore prove to the converse of the~$\Sigma^0_{n+3}$
case without Lemma~\ref{lem:forcing-dense}.
Fix a condition $c = (F, \sigma, T) \in \Fcal$ such that $c \Vdash \varphi(G)$. 
\begin{itemize}
	\item If $\varphi \in \Sigma^0_1$ then $\varphi(G)$ can be expressed as $(\exists x)\psi(G, x)$ where $\psi \in \Sigma^0_0$.
	By clause~(i) of Definition~\ref{def:forcing-precondition}, there exists a $w \in \omega$ such that
	$\psi(F, w)$ holds. As $F \subseteq G$ and~$G \setminus F \subseteq (max(F), +\infty)$, then by continuity $\psi(G, w)$ holds, hence $\varphi(G)$ holds.
	
	\item If $\varphi \in \Pi^0_1$ then $\varphi(G)$ can be expressed as $(\forall x)\psi(G, x)$ where $\psi \in \Sigma^0_0$.
	By clause~(ii) of Definition~\ref{def:forcing-precondition}, for every $w \in \omega$
	and every $E \in \Apx(c)$, $\psi(E, w)$ holds. As $\{E \subset_{fin} G : E \supseteq F \} \subseteq \Apx(c)$,
	then for every $w \in \omega$, $\psi(G, w)$ holds, so $\varphi(G)$ holds.
 
	\item If $\varphi \in \Sigma^0_{n+2}$ then $\varphi(G)$ can be expressed as $(\exists x)\psi(G, x)$ 
  where $\psi \in \Pi^0_{n+1}$.
  By clause~(i) of Definition~\ref{def:forcing-condition}, there exists a $w \in \omega$
  such that $c \Vdash \psi(G, w)$. By induction hypothesis, $\psi(G, w)$ holds, hence $\varphi(G)$ holds.
  
  Conversely, suppose that $\varphi(G)$ holds. Then there exists a $w \in \omega$ such that $\psi(G, w)$ holds,
  so by induction hypothesis $c \Vdash \psi(G, w)$ for some $c \in \Fcal$,
  so by clause~(i) of Definition~\ref{def:forcing-condition}, $c \Vdash \varphi(G)$.

	\item If $\varphi \in \Pi^0_2$ then $\varphi(G)$ can be expressed as $(\forall x)\psi(G, x)$ where $\psi \in \Sigma^0_1$.
	By clause~(ii) of Definition~\ref{def:forcing-condition}, for every $\tau \in T$, every $w < |\tau|$,
	and every $E \in \Apx_{|\tau|}(c)$, $(E, \tau, T^{[\tau]}) \not \Vdash \neg \psi(G, w)$. 
  Suppose by way of contradiction that $\psi(G, w)$ does not hold for some $w \in \omega$.
  Then by induction hypothesis, there exists a $d \in \Fcal$ such that $d \Vdash \neg \psi(G, w)$.
  Let $e = (E, \tau, S) \in \Fcal$ be such that $e \Vdash \neg \psi(G, w)$, $|\tau| > w$
  and $e$ extends both $c$ and $d$. The condition $e$ exists by Lemma~\ref{lem:forcing-extension-precondition}.
  We can furthermore require that $E \in \Apx_{|\tau|}(c)$,
  so $e \not \Vdash \neg \psi(G, w)$ and $e \Vdash \neg \psi(G, w)$. Contradiction.
  Hence for every $w \in \omega$, $\psi(G, w)$ holds, so $\varphi(G)$ holds.

	\item If $\varphi \in \Pi^0_{n+3}$ then $\varphi(G)$ can be expressed as $\neg \psi(G)$ where $\psi \in \Sigma^0_{n+3}$.
  By clause~(iii) of Definition~\ref{def:forcing-condition}, for every $d \in \Ext(c)$, $d \not \Vdash \psi(G)$.
	By Lemma~\ref{lem:forcing-extension}, $d \not \Vdash \psi(G)$ for every~$d \in \Fcal$, and by 
  a previous case, $\psi(G)$ does not hold, so $\varphi(G)$ holds.
\end{itemize}
\end{proof}

We now prove that the forcing relation enjoys the desired definitional properties,
that is, the complexity of the forcing relation is the same as the complexity
of the formula it forces. We start by analysing the complexity of some components
of this notion of forcing.

\begin{lemma}\label{lem:extension-complexity}\ 
\begin{itemize}
	\item[a)] For every precondition $c$, $\Apx(c)$ and $\Ext_1(c)$ are $\Delta^0_1$ uniformly in~$c$.
	\item[b)] For every condition $c$, $\Ext(c)$ is $\Pi^0_2$ uniformly in~$c$.
\end{itemize}
\end{lemma}
\begin{proof}\ 
\begin{itemize}
	\item[a)] Fix a precondition $c = (F, \sigma, T)$.
	A set $E \in \Apx(c)$ iff the following $\Delta^0_1$ predicate holds:
	$$
	(F \subseteq E) \wedge (\forall x \in E \setminus F)[x > max(F) \wedge x \in R_\sigma]
	$$
	Moreover, $(E, \tau, S) \in \Ext_1(c)$ iff the $\Delta^0_1$ predicate 
	$E \in \Apx(c) \wedge \tau = \sigma \wedge S = T$ holds.
	As already mentioned, the equality~$S = T$ is translated into ``the indices
	of~$S$ and~$T$ coincide'' which is a $\Sigma^0_0$ statement.

	\item[b)] Fix a condition $c = (F, \sigma, T)$.
	By clause 2) of Lemma~\ref{lem:basic-statements},
	$(E, \tau, S) \in \Ext(c)$ iff the following $\Pi^0_2$ formula holds
	$$
	\begin{array}{ll}
	E \in \Apx(c) \wedge \sigma \preceq \tau \\
	\wedge (\forall \rho \in S)(\forall \xi)[\xi \preceq \rho \imp \xi \in S] & \mbox{ ($S$ is a tree)}\\
	\wedge (\forall n)(\exists \rho \in 2^n)\rho \in S) & \mbox{ ($S$ is infinite) }\\
	\wedge (\forall \rho \in S)(\sigma \prec \rho \vee \rho \preceq \sigma) & \mbox{ ($S$ has stem $\sigma$)}\\
	\wedge (\forall \rho \in S)(\rho \in T) & \mbox{ ($S$ is a subset of $T$) }\\
	\end{array}
	$$
	
\end{itemize}
\end{proof}

\begin{lemma}\label{lem:complexity-forcing}
Fix an arithmetic formula $\varphi(G)$.
\begin{itemize}
	\item[a)] Given a precondition $c$, if $\varphi(G)$ is a $\Sigma^0_1$ ($\Pi^0_1$) formula 
	then so is the predicate $c \Vdash \varphi(G)$.
	\item[b)] Given a condition $c$, if $\varphi(G)$ is a $\Sigma^0_{n+2}$ ($\Pi^0_{n+2}$) formula 
	then so is the predicate $c \Vdash \varphi(G)$.
\end{itemize}
\end{lemma}
\begin{proof}
We prove our lemma by induction over the complexity of the formula $\varphi(G)$.
Fix a (pre)condition $c = (F, \sigma, T)$.
\begin{itemize}
	\item If $\varphi(G) \in \Sigma^0_1$ then it can be expressed as $(\exists x)\psi(G, x)$ where $\psi \in \Sigma^0_0$.
	By clause~(i) of Definition~\ref{def:forcing-precondition}, $c \Vdash \varphi(G)$ if and only if 
	the formula $(\exists w \in \omega)\psi(F, w)$ holds. This is a $\Sigma^0_1$ predicate.
	
	\item If $\varphi(G) \in \Pi^0_1$ then it can be expressed as $(\forall x)\psi(G, x)$ where $\psi \in \Sigma^0_0$.
	By clause~(ii) of Definition~\ref{def:forcing-precondition}, $c \Vdash \varphi(G)$ if and only if 
	the formula $(\forall w \in \omega)(\forall E \in \Apx(c))\psi(E, w)$ holds. 
	By clause~a) of Lemma~\ref{lem:extension-complexity}, this is a $\Pi^0_1$ predicate.

	\item If $\varphi(G) \in \Sigma^0_{n+2}$ then it can be expressed as $(\exists x)\psi(G, x)$ where $\psi \in \Pi^0_{n+1}$.
	By clause~(i) of Definition~\ref{def:forcing-condition}, $c \Vdash \varphi(G)$ if and only if 
	the formula $(\exists w < |\sigma|)c \Vdash \psi(G, w)$ holds. This is a $\Sigma^0_{n+2}$ predicate
	by induction hypothesis.

	\item If $\varphi(G) \in \Pi^0_2$ then it can be expressed as $(\forall x)\psi(G, x)$ where $\psi \in \Sigma^0_1$.
	By clause~(ii) of Definition~\ref{def:forcing-condition}, $c \Vdash \varphi(G)$ if and only if 
	the formula $(\forall \tau \in T)(\forall w < |\tau|)(\forall E \in \Apx_{|\tau|}(c))
	(E, \tau, T^{[\tau]}) \not \Vdash \neg \psi(G, w)$ holds. 
	By induction hypothesis, $(E, \tau, T^{[\tau]}) \not \Vdash \neg \psi(G, w)$ is a $\Sigma^0_1$ predicate,
	hence by clause~a) of Lemma~\ref{lem:extension-complexity}, $c \Vdash \varphi(G)$ is a $\Pi^0_2$ predicate.

	\item If $\varphi(G) \in \Pi^0_{n+3}$ then it can be expressed as $\neg \psi(G)$ where $\psi \in \Sigma^0_{n+3}$. 
	By clause~(iii) of Definition~\ref{def:forcing-condition}, $c \Vdash \varphi(G)$ if and only if 
	the formula $(\forall d)(d \not \in \Ext(c) \vee d \not \Vdash \psi(G))$ holds.
	By induction hypothesis, $d \not \Vdash \psi(G)$ is a $\Pi^0_{n+3}$ predicate.
	Hence by clause~b) of Lemma~\ref{lem:extension-complexity},
	$c \Vdash \varphi(G)$ is a $\Pi^0_{n+3}$ predicate.
\end{itemize}
\end{proof}

\subsection{Preserving the arithmetic hierarchy}

The following lemma asserts that every sufficiently generic real
for this notion of forcing preserves the arithmetic hierarchy.
The argument deeply relies on the fact that this notion of forcing
admits a forcing relation with good definitional properties.

\begin{lemma}\label{lem:diagonalization}
If $A \not \in \Sigma^0_{n+1}$ and $\varphi(G, x)$ is $\Sigma^0_{n+1}$, 
then the set of $c \in \Pb$ satisfying the following property is dense:
$$
[(\exists w \in A)c \Vdash \neg \varphi(G, w)] \vee [(\exists w \not \in A)c \Vdash \varphi(G, w)]
$$
\end{lemma}
\begin{proof}
Fix a condition $c = (F, \sigma, T)$.
\begin{itemize}
	\item In case $n = 0$, $\varphi(G, w)$ can be expressed as $(\exists x)\psi(G, w, x)$ where $\psi \in \Sigma^0_0$.
	Let $U = \{ w \in \omega : (\exists E \in \Apx(c))(\exists u)\psi(E, w, u) \}$.
	By clause~a) of Lemma~\ref{lem:extension-complexity}, $U \in \Sigma^0_1$, thus $U \neq A$.
	Fix $w \in U \Delta A$. If $w \in U \setminus A$ then by definition of~$U$,
	there exist an $E \in \Apx(c)$ and a $u \in \omega$ such that $\psi(E, w, u)$ holds.
	By definition of $\Apx(c)$ and clause~1) of Lemma~\ref{lem:basic-statements}, $d = (E, \sigma, T)$
	is a condition extending $c$. By clause~(i) of Definition~\ref{def:forcing-precondition},
	$d \Vdash \varphi(G, w)$.
	If $w \in A \setminus U$, then for every $E \in \Apx(c)$ and every $u \in \omega$, 
	$\psi(E, w, u)$ does not hold, so by clause~(ii) of Definition~\ref{def:forcing-precondition},
	$c \Vdash (\forall x)\neg \psi(G, w, x)$, hence $c \Vdash \neg \varphi(G, w)$.

	\item In case $n = 1$, $\varphi(G, w)$ can be expressed as $(\exists x)\psi(G, w, x)$ where $\psi \in \Pi^0_1$.
	Let $U = \{ w \in \omega : (\exists s)(\forall \tau \in 2^s \cap T)
	(\exists u < s)(\exists E \in \Apx_s(c))(E, \tau, T^{[\tau]}) \Vdash \psi(G, w, u) \}$.
	By Lemma~\ref{lem:complexity-forcing} and clause~a) of Lemma~\ref{lem:extension-complexity},
	$U \in \Sigma^0_2$, thus $U \neq A$. Fix $w \in U \Delta A$.
	If $w \in U \setminus A$ then by definition of~$U$, there exist an $s \in \omega$,
	a $\tau \in 2^s \cap T$, a $u < s$ and an $E \in \Apx_s(c)$ such that $T^{[\tau]}$ is infinite
	and $(E, \tau, T^{[\tau]}) \Vdash \psi(G, w, u)$. Thus $d = (E, \tau, T^{[\tau]})$ is a condition
	extending $c$ and by clause~(i) of Definition~\ref{def:forcing-condition}, $d \Vdash \varphi(G, w)$.
	If $w \in A \setminus U$, then let 
	$S = \{ \tau \in T : (\forall u < |\tau|)(\forall E \in \Apx_{|\tau|}(c) (E, \tau, T^{[\tau]}) \not \Vdash \psi(G, w, u) \}$.
	As proven in Lemma~\ref{lem:forcing-dense}, $S$ is a $\emptyset'$-p.r. subtree of~$T$
	and by $w \not \in U$, $S$ is infinite. Thus $d = (F, \sigma, S)$ is a condition extending $c$.
	By clause~3) of Lemma~\ref{lem:basic-statements}, $\Apx(d) \subseteq \Apx(c)$, 
	so for every $\tau \in S$, every $u < |\tau|$, and every $E \in \Apx_{|\tau|}(d)$,
	$(E, \tau, T^{[\tau]}) \not \Vdash \psi(G, w, u)$. By Lemma~\ref{lem:tree-no-effect-first-level},
	$(E, \tau, S^{[\tau]}) \not \Vdash \psi(G, w, u)$, so by clause~(ii) of Definition~\ref{def:forcing-condition},
	$d \Vdash (\forall x)\neg \psi(G, w, u)$ hence $d \Vdash \neg \varphi(G, w)$.

	\item In case $n > 1$, let $U = \{ w \in \omega : (\exists d \in \Ext(c)) d \Vdash \varphi(G, w) \}$.
	By clause~b) of Lemma~\ref{lem:extension-complexity} and Lemma~\ref{lem:complexity-forcing},
	$U \in \Sigma^0_n$, thus $U \neq A$.
	Fix $w \in U \Delta A$. If $w \in U \setminus A$ then by definition of~$U$,
	there exists a condition $d$ extending $c$ such that $d \Vdash \varphi(G, w)$.
	If $w \in A \setminus U$, then for every $d \in \Ext(c) d \not \Vdash \varphi(G, w)$
	so by clause~(iii) of Definition~\ref{def:forcing-condition}, $c \Vdash \neg \varphi(G, w)$. 
\end{itemize}
\end{proof}

We are now ready to prove Theorem~\ref{thm:coh-preservation-arithmetic-hierarchy}.

\begin{proof}[Proof of Theorem~\ref{thm:coh-preservation-arithmetic-hierarchy}]
Let~$C$ be a set and~$R_0, R_1, \dots$ be a uniformly $C$-computable sequence of sets.
Let~$T_0$ be a $C'$-primitive recursive tree such that~$[T_0] \subseteq \Ccal(\vec{R})$.
Let~$\Fcal$ be a sufficiently generic filter containing~$c_0 = (\emptyset, \epsilon, T_0)$.
and let $G$ be the corresponding generic real. By genericity, the set~$G$ is
an infinite $\vec{R}$-cohesive set.
By Lemma~\ref{lem:diagonalization} and Lemma~\ref{lem:complexity-forcing},
$G$ preserves non-$\Sigma^0_{n+1}$ definitions relative to~$C$ for every~$n \in \omega$.
Therefore, by Proposition 2.2 of~\cite{Wang2014Definability}, $G$
preserves the arithmetic hierarchy relative to~$C$.
\end{proof}

\section{An effective forcing for the Erd\H{o}s-Moser theorem}

We now extend the previous result to the Erd\H{o}s-Moser theorem.
This is an extension in that the stable thin set theorem for pairs
does not admit preservation of the arithmetic hierarchy (see Wang~\cite{Wang2014Definability}). 
Therefore, by Theorem~\ref{thm:emo-implies-sts2-or-coh}, if the Erd\H{o}s-Moser admits preservation
of the arithmetic hierarchy, then so does cohesiveness.

\begin{theorem}\label{thm:em-preserves-arithmetic}
$\emo$ admits preservation of the arithmetic hierarchy.
\end{theorem}

Again, the core of the proof consists of finding a good forcing notion
whose generics will preserve the arithmetic hierarchy.
For the sake of simplicity, we will restrict ourselves to stable tournaments
even though it is clear that the forcing notion can be adapted to arbitrary tournaments.
The proof of Theorem~\ref{thm:em-preserves-arithmetic} will be obtained by composing the proof that cohesiveness
and the stable Erd\H{o}s-Moser theorem admit preservation of the arithmetic hierarchy.
We first need to introduce some terminology
and in particular the notion of partition tree.

\subsection{Partition trees}

Recall that given some~$k \in \omega$, $k^{<\omega}$ is the set of all strings $\sigma$
such that~$\sigma(i) < k$ for each~$i < |\sigma|$.
Given a string $\sigma \in k^{<\omega}$,
we denote by $\seto_\nu(\sigma)$ the set $\{ x < |\sigma| : \sigma(x) = \nu \}$
where $\nu < k$. The notion can be extended to sequences $P \in k^{\omega}$
where $\seto_\nu(P) = \{ x \in \omega : P(x) = \nu \}$.

\index{partition tree}
\begin{definition}[Partition tree]
A \emph{$k$-partition tree of $[t, +\infty)$} for some $k, t \in \omega$
is a tuple $(k, t, T)$ such that $T$ is a subtree of $k^{<\omega}$.
A \emph{partition tree} is a $k$-partition tree of $[t, +\infty)$ for some $k, t \in \omega$.
\end{definition}

To simplify our notation, we may use the same letter $T$ to denote both a partition tree $(k, t, T)$
and the actual tree $T \subseteq k^{<\omega}$. We then write $\dom(T)$ for $[t, +\infty)$
and $\parts(T)$ for~$k$. Given a p.r.\ partition tree~$T$,
we write~$\#T$ for its Turing index, and may refer to it as its \emph{code}.

\begin{definition}[Refinement]
Given a function $f : \ell \to k$, a string $\sigma \in \ell^{<\omega}$ \emph{$f$-refines} a string $\tau \in k^{<\omega}$
if $|\sigma| = |\tau|$ and for every $\nu < \ell$, $\seto_\nu(\sigma) \subseteq \seto_{f(\nu)}(\tau)$.
A p.r.\ $\ell$-partition tree $S$ of $[u,+\infty)$ \emph{$f$-refines} a p.r.\ $k$-partition tree $T$ 
of $[t, +\infty)$ (written $S \leq_f T$)
if $\#S \geq \#T$, $\ell \geq k$, $u \geq t$ and for every $\sigma \in S$, $\sigma$ $f$-refines some $\tau \in T$. 
\end{definition}

The collection of partition trees is equipped
with a partial order $\leq$ such that $(\ell, u, S) \leq (k, t, T)$ if there exists a function $f : \ell \to k$ such that $S \leq_f T$.
Given a $k$-partition tree of $[t, +\infty)$ $T$,
we say that part $\nu$ of $T$ is \emph{acceptable} if there exists a path $P$ through $T$
such that $\seto_\nu(P)$ is infinite.
Moreover, we say that part $\nu$ of $T$ is \emph{empty} if 
$(\forall \sigma \in T)[dom(T) \cap \seto_\nu(\sigma) = \emptyset]$.
Note that each partition tree has at least one acceptable part
since for every path~$P$ through~$T$, $\seto_\nu(P)$ is infinite for some~$\nu < k$.
It can also be the case that part~$\nu$ of~$T$ is non-empty,
while for every path $P$ through~$T$, $\seto_\nu(P) \cap \dom(T) = \emptyset$.
However, in this case, we can choose the infinite computable subtree~$S = \{\sigma \in T : \seto_\nu(\sigma) \cap \dom(T) = \emptyset\}$
of~$T$ which has the same collection of infinite paths and such that part~$\nu$ of~$S$ is empty.

Given a $k$-partition tree $T$, a finite set $F \subseteq \omega$ and a part $\nu < k$,
define
$$
T^{[\nu, F]} = \{ \sigma \in T :  F \subseteq \seto_\nu(\sigma) \vee |\sigma| < max(F) \}
$$
The set $T^{[\nu, F]}$ is a (possibly finite) subtree of~$T$ which id-refines $T$
and such that~$F \subseteq \seto_\nu(P)$ for every path~$P$ through $T^{[\nu, F]}$.

We denote by $\TPb$ the set of all ordered pairs $(\nu, T)$
such that $T$ is an infinite, primitive recursive $k$-partition tree of $[t, +\infty)$
for some $t, k \in \omega$ and $\nu < k$. The set $\TPb$ is equipped with a partial ordering $\leq$
such that $(\mu, S) \leq (\nu, T)$ if $S$ $f$-refines $T$ and $f(\mu) = \nu$.
In this case we say that \emph{part $\mu$ of $S$ refines part $\nu$ of~$T$}.
Note that the domain of $\TPb$ and the relation $\leq$ are co-c.e.
We denote by $\TPb[T]$ the set of all $(\nu, S) \leq (\mu, T)$ for some $(\mu, T) \in \TPb$.

\begin{definition}[Promise for a partition tree]\label{def:em-promise}
Fix a p.r. $k$-partition tree of $[t,+\infty)$ $T$.
A class $\Ccal \subseteq \TPb[T]$ is a \emph{promise for $T$} if
\begin{itemize}
	\item[a)] $\Ccal$ is upward-closed under the $\leq$ relation restricted to $\TPb[T]$
	\item[b)] for every infinite p.r. partition tree $S \leq T$,
	$(\mu, S) \in \Ccal$ for some non-empty part $\mu$ of~$S$.
\end{itemize}
\end{definition}

A promise for $T$ can be seen as a two-dimensional tree
with at first level the acyclic digraph of refinement of partition trees.
Given an infinite path in this digraph, the parts of the members of this path
form an infinite, finitely branching tree.

\begin{lemma}\label{lem:em-refinement-complexity}
Let $T$ and $S$ be p.r.\  partition trees such that $S \leq_f T$ for some function $f : \parts(S) \to \parts(T)$
and let $\Ccal$ be a $\emptyset'$-p.r.\ promise for $T$.
\begin{itemize}
	\item[a)] The predicate ``$T$ is an infinite $k$-partition tree of $[t, +\infty)$'' is $\Pi^0_1$ uniformly in $T$, $k$ and $t$.
	\item[b)] The relations ``$S$ $f$-refines $T$''
	and ``part $\nu$ of $S$ $f$-refines part $\mu$ of $T$''
	are $\Pi^0_1$ uniformly in $S$, $T$ and $f$.
	\item[c)] The predicate ``$\Ccal$ is a promise for $T$'' is $\Pi^0_2$ uniformly in an index for~$\Ccal$ and~$T$.
\end{itemize}
\end{lemma}
\begin{proof}\ 
\begin{itemize}
	\item[a)] $T$ is an infinite $k$-partition tree of $[t, +\infty)$ if and only if the $\Pi^0_1$ formula
$[(\forall \sigma \in T)(\forall \tau \preceq \sigma) \tau \in T \cap k^{<\infty}]
 \wedge [(\forall n)(\exists \tau \in k^n) \tau \in T]$  holds.

	\item[b)] Suppose that $T$ is a $k$-partition tree of $[t,+\infty)$ and $S$ is an $\ell$-partition tree of $[u,+\infty)$.
$S$ $f$-refines $T$ if and only if the following $\Pi^0_1$ formula holds:
$$
u \geq t \wedge \ell \geq k \wedge [(\forall \sigma \in S)(\exists \tau \in k^{|\sigma|} \cap T)
(\forall \nu < u)set_\nu(\sigma) \subseteq \seto_{f(\nu)}(\tau)]
$$
Part $\nu$ of $S$ $f$-refines part $\mu$ of $T$ if and only if $\mu = f(\nu)$ and $S$ $f$-refines $T$.

	\item[c)] Given $k, t \in \omega$, let $PartTree(k, t)$ denote 
	the $\Pi^0_1$ set of all the infinite p.r.\ $k$-partition trees of $[t, +\infty)$.
	Given a $k$-partition tree $S$ and a part $\nu$ of $S$,
	let $Empty(S, \nu)$ denote the $\Pi^0_1$ formula ``part $\nu$ of $S$ is empty'',
	that is the formula $(\forall \sigma \in S) \seto_\nu(\sigma) \cap \dom(S) = \emptyset$.

	$\Ccal$ is a promise for $T$ if and only if the following $\Pi^0_2$ formula holds:
	$$
	\begin{array}{l}
	(\forall \ell, u)(\forall S \in PartTree(\ell,u))[S \leq T \imp (\exists \nu < \ell) \neg Empty(S, \nu) \wedge (\nu, S) \in \Ccal)] \\
	\wedge (\forall \ell', u')(\forall V \in PartTree(\ell', u'))(\forall g : \ell \to \ell')[ S \leq_g V \leq T \imp \\
	(\forall \nu < \ell)((\nu, S) \in \Ccal \imp (g(\nu), V) \in \Ccal)]
	\end{array}
	$$
\end{itemize}
\end{proof}

Given a promise $\Ccal$ for $T$ and some infinite p.r. partition tree $S$
refining $T$, we denote by $\Ccal[S]$ the set of all $(\nu, S') \in \Ccal$
below some $(\mu, S) \in \Ccal$, that is, $\Ccal[S] = \Ccal \cap \TPb[S]$.
Note that by clause b) of Lemma~\ref{lem:em-refinement-complexity},
if $\Ccal$ is a $\emptyset'$-p.r. promise for $T$ then $\Ccal[S]$
is a $\emptyset'$-p.r. promise for~$S$.

Establishing a distinction between the acceptable parts and the non-acceptable ones requires
a lot of definitional power. However, we prove that we can always find an
extension where the distinction is $\Delta^0_2$.
We say that an infinite p.r. partition tree $T$
\emph{witnesses its acceptable parts} if its parts are either acceptable or empty.

\begin{lemma}\label{lem:em-promise-extension-witnessing-acceptable}
For every infinite p.r.\ $k$-partition tree $T$ of $[t, +\infty)$,
there exists an infinite p.r.\ $k$-partition tree $S$ of $[u, +\infty)$
refining $T$ with the identity function and such that $S$ witnesses its acceptable parts.
\end{lemma}
\begin{proof}
Given a partition tree $T$, we let $I(T)$ be the set of its empty parts.
Fix an infinite p.r. $k$-partition tree of $[t, +\infty)$ $T$,
It suffices to prove that if $\nu$ is a non-empty and non-acceptable part of~$T$,
then there exists an infinite p.r.\ $k$-partition tree $S$
refining $T$ with the identity function, such that $\nu \in I(S) \setminus I(T)$.
As $I(T) \subseteq I(S)$ and $|I(S)| \leq k$, it suffices to iterate the process
at most $k$ times to obtain a refinement witnessing its acceptable parts.

So fix a non-empty and non-acceptable part $\nu$ of~$T$.
By definition of being non-acceptable, there exists a path $P$ through $T$
and an integer $u > max(t, \seto_\nu(P))$.
Let $S = \{ \sigma \in T : \seto_\nu(\sigma) \cap [u, +\infty) = \emptyset \}$.
The set $S$ is a p.r. $k$-partition tree of $[u, +\infty)$ refining
$T$ with the identity function and such that part $\nu$ of $S$ is empty.
Moreover, $S$ is infinite since~$P \in [S]$. 
\end{proof}

The following lemma strengthens clause b) of Definition~\ref{def:em-promise}.

\begin{lemma}
Let $T$ be a p.r. partition tree and $\Ccal$ be a promise for $T$.
For every infinite p.r. partition tree $S \leq T$, 
$(\mu, S) \in \Ccal$ for some acceptable part $\mu$ of~$S$.
\end{lemma}
\begin{proof}
Fix an infinite p.r. $\ell$-partition tree $S \leq T$.
By Lemma~\ref{lem:em-promise-extension-witnessing-acceptable},
there exists an infinite p.r. $\ell$-partition tree $S' \leq_{id} S$
witnessing its acceptable parts. As $\Ccal$ is a promise for $T$
and $S' \leq T$, there exists a non-empty (hence acceptable) part $\nu$ of $S'$ such that $(\nu, S') \in \Ccal$.
As $\Ccal$ is upward-closed, $(\nu, S) \in \Ccal$.
\end{proof}

\subsection{Forcing conditions}

We now describe the forcing notion for the Erd\H{o}s-Moser theorem.
Recall that an EM condition for an infinite tournament~$R$
is a Mathias condition $(F, X)$ where
$F \cup \{x\}$ is $R$-transitive for each $x \in X$
and $X$ is included in a minimal $R$-interval of $F$.
The main properties of an EM condition are proven in section~\ref{sect:weakness-erdos-moser-theorem}
under Lemma~\ref{lem:emo-cond-beats} and Lemma~\ref{lem:emo-cond-valid}.

\begin{definition}
We denote by~$\Pb$ the forcing notion whose conditions are tuples $(\vec{F}, T, \Ccal)$ where
\begin{itemize}
	\item[(a)] $T$ is an infinite p.r.\ partition tree
	\item[(b)] $\Ccal$ is a $\emptyset'$-p.r. promise for $T$
	\item[(c)] $(F_\nu, \dom(T))$ is an EM condition for $R$ and each $\nu < \parts(T)$
\end{itemize}
A condition $d = (\vec{E}, S, \Dcal)$ \emph{extends} $c = (\vec{F}, T, \Ccal)$
(written $d \leq c$) if there exists a function $f : \ell \to k$ such that
$\Dcal \subseteq \Ccal$ and the followings hold:
\begin{itemize}
	\item[(i)] $(E_\nu, \dom(S))$ EM extends $(F_{f(\nu)}, \dom(T))$ for each $\nu < \parts(S)$ 
	\item[(ii)] $S$ $f$-refines $\bigcap_{\nu < \parts(S)} T^{[f(\nu), E_\nu]}$
\end{itemize}
\end{definition}

We may think of a condition $c = (\vec{F}, T, \Ccal)$ as a collection of EM conditions
$(F_\nu, H_\nu)$ for~$R$, where $H_\nu = \dom(T) \cap \seto_\nu(P)$ for some path $P$ through $T$.
$H_\nu$ must be infinite for at least one of the parts $\nu < \parts(T)$.
At a higher level, $\Dcal$ restricts the possible subtrees $S$
and parts $\mu$ refining some part of $T$ in the condition~$c$.
Given a condition $c = (\vec{F}, T, \Ccal)$, we write $\parts(c)$ for $\parts(T)$.

\begin{lemma}\label{lem:em-forcing-infinite}
For every condition $c = (\vec{F}, T, \Ccal)$ and every $n \in \omega$, there exists an extension $d = (\vec{E}, S, \Dcal)$
such that $|E_\nu| \geq n$ on each acceptable part $\nu$ of~$S$.
\end{lemma}
\begin{proof}
It suffices to prove that for every condition $c = (\vec{F}, T, \Ccal)$
and every acceptable part $\nu$ of $T$, 
there exists an extension $d = (\vec{E}, S, \Dcal)$ such that $S \leq_{id} T$ 
and $|E_\nu| \geq n$. Iterating the process at most $\parts(T)$ times completes the proof.
Fix an acceptable part $\nu$ of $T$ and a path $P$ trough $T$
such that $\seto_\nu(P)$ is infinite. Let $F'$ be an $R$-transitive subset of $\seto_\nu(P) \cap \dom(T)$
of size $n$. Such a set exists by the classical Erd\H{o}s-Moser theorem. Let $\vec{E}$ be defined by $E_\mu = F_\mu$
if $\mu \neq \nu$ and $E_\nu = F_\nu \cup F'$ otherwise. As the tournament~$R$ is stable,
there exists some~$u \geq t$ such that $(E_\nu, [u, +\infty))$ is an EM condition
and therefore EM extends $(F_\nu, \dom(T))$.
Let $S$ be the p.r.\ partition tree $T^{[\nu, E_\nu]}$ of $[u, +\infty)$.
The condition $(\vec{E}, S, \Ccal[S])$ is the desired extension.
\end{proof}

Given a condition~$c \in \Pb$, we denote by $\Ext(c)$ the set of all its extensions. 

\subsection{The forcing relation}

The forcing relation at the first level, namely, for $\Sigma^0_1$ and $\Pi^0_1$ formulas,
is parameterized by some part of the tree of the considered condition. Thanks to the forcing relation
we will define, we can build an infinite decreasing sequence of conditions which decide $\Sigma^0_1$ and~$\Pi^0_1$ formulas
effectively in~$\emptyset'$. The sequence however yields a $\emptyset'$-computably bounded $\emptyset'$-computable tree of (possibly empty)
parts. Therefore, any PA degree relative to~$\emptyset'$ is sufficient to control the first jump of an infinite
transitive subtournament of a stable infinite computable tournament.

We cannot do better since Kreuzer proved in~\cite{Kreuzer2012Primitive} the existence
of an infinite, stable, computable tournament with no low infinite transitive subtournament.
If we ignore the promise part of a condition,
the careful reader will recognize the construction of Cholak, Jockusch and Slaman~\cite{Cholak2001strength} of a low${}_2$ infinite
subset of a $\Delta^0_2$ set or its complement by the first jump control. 
The difference, which at first seems only presentational,
is in fact one of the key features of this notion of forcing. Indeed, forcing iterated jumps require to have
a definitionally weak description of the set of extensions of a condition,
and it requires much less computational power to describe a primitive recursive
tree than an infinite reservoir of a Mathias condition.

\begin{definition}\label{def:em-forcing-precondition}
Fix a condition $c = (\vec{F}, T, \Ccal)$, a $\Sigma^0_0$ formula $\varphi(G, x)$
and a part $\nu < \parts(T)$.
\begin{itemize}
	\item[1.] $c \Vdash_\nu (\exists x)\varphi(G, x)$ iff there exists a $w \in \omega$ such that $\varphi(F_\nu, w)$ holds.
	\item[2.] $c \Vdash_\nu (\forall x)\varphi(G, x)$ iff 
	for every $\sigma \in T$, every $w < |\sigma|$ and every $R$-transitive set $F' \subseteq \dom(T) \cap \seto_\nu(\sigma)$,
	$\varphi(F_\nu \cup F', w)$ holds.
\end{itemize}
\end{definition}

We start by proving some basic properties of the forcing relation over~$\Sigma^0_1$
and $\Pi^0_1$ formulas. As one may expect, the forcing relation at first level
is closed under the refinement relation.

\begin{lemma}\label{lem:em-forcing-extension-level1}
Fix a condition $c = (\vec{F}, T, \Ccal)$ and a $\Sigma^0_1$ ($\Pi^0_1$) formula $\varphi(G)$.
If $c \Vdash_\nu \varphi(G)$ for some $\nu < \parts(T)$, then for every $d = (\vec{E}, S, \Dcal) \leq c$
and every part $\mu$ of $S$ refining part $\nu$ of $T$, $d \Vdash_{\mu} \varphi(G)$.
\end{lemma}
\begin{proof}We have two cases.
\begin{itemize}
	\item If $\varphi \in \Sigma^0_1$ then it can be expressed as $(\exists x)\psi(G, x)$
	where $\psi \in \Sigma^0_0$. By clause~1 of Definition~\ref{def:em-forcing-precondition},
	there exists a $w \in \omega$ such that
	$\psi(F_\nu, w)$ holds. By property (i) of the definition of an extension, $E_\mu \supseteq F_\nu$
	and $(E_\mu \setminus F_\nu) \subset \dom(T)$, therefore $\psi(E_\mu, w)$ holds by continuity,
	so by clause~1 of Definition~\ref{def:em-forcing-precondition}, $d \Vdash_\mu (\exists x)\psi(G, x)$.

	\item If $\varphi \in \Pi^0_1$ then it can be expressed as $(\forall x)\psi(G, x)$
	where $\psi \in \Sigma^0_0$. Fix a $\tau \in S$, a $w < |\tau|$ and an $R$-transitive set
	$F' \subseteq \dom(S) \cap \seto_\mu(\tau)$.
	It suffices to prove that $\varphi(E_\mu \cup F')$ holds to conclude
	that $d \Vdash_\mu (\forall x)\psi(G, x)$ by clause~2 of Definition~\ref{def:em-forcing-precondition}.
 	By property~(ii) of the definition of an extension,
	there exists a $\sigma \in T^{[\nu, E_\mu]}$
	such that $|\sigma| = |\tau|$ and $\seto_\mu(\tau) \subseteq \seto_\nu(\sigma)$. 
	As $\dom(S) \subseteq \dom(T)$, $F' \subseteq \dom(T) \cap \seto_\nu(\sigma)$.
	As $\sigma \in T^{[\nu, E_\mu]}$, $E_\mu \subseteq \seto_\nu(\sigma)$
	and by property~(i) of the definition of an extension, $E_\mu \subseteq \dom(T)$.
	So $E_\mu \cup F' \subseteq \dom(T) \cap \seto_\nu(\sigma)$.
	As $w < |\tau| = |\sigma|$ and~$E_\mu \cup F'$ is an $R$-transitive subset of~$\dom(T) \cap \seto_\nu(\sigma)$,
	then by clause~2 of Definition~\ref{def:em-forcing-precondition} applied to $c \Vdash_\nu (\forall x)\psi(G, x)$, 
	$\varphi(F_\nu \cup (E_\mu \setminus F_\nu) \cup F', w)$ holds, hence $\varphi(E_\mu \cup F')$ holds.
\end{itemize}
\end{proof}

Before defining the forcing relation at higher levels, we prove a density lemma
for $\Sigma^0_1$ and $\Pi^0_1$ formulas. It enables us
in particular to reprove that every degree PA relative to $\emptyset'$
computes the jump of an infinite $R$-transitive set.

\begin{lemma}\label{lem:em-forcing-dense-level1}
For every $\Sigma^0_1$ ($\Pi^0_1$) formula $\varphi$, the following set is dense
$$
\{ c = (\vec{F}, T, \Ccal) \in \Pb : (\forall \nu < \parts(T))[ c \Vdash_\nu \varphi(G) \vee c \Vdash_\nu \neg \varphi(G)] \}
$$
\end{lemma}
\begin{proof}
It suffices to prove the statement for the case where $\varphi$ is a $\Sigma^0_1$ formula,
as the case where $\varphi$ is a $\Pi^0_1$ formula is symmetric. Fix a condition $c = (\vec{F}, T, \Ccal)$
and let $I(c)$ be the set of the parts $\nu < \parts(T)$ such that $c \not \Vdash_\nu \varphi(G)$
and $c \not \Vdash_\nu \neg \varphi(G)$. If $I(c) = \emptyset$ then we are done, so suppose $I(c) \neq \emptyset$
and fix some $\nu \in I(c)$. We will construct an extension $d$ of~$c$ such that $I(d) \subseteq I(c) \setminus \{\nu\}$.
Iterating the operation completes the proof.

The formula $\varphi$ is of the form $(\exists x)\psi(G, x)$ where $\psi \in \Sigma^0_0$.
Define $f : k+1 \to k$ as $f(\mu) = \mu$ if $\mu < k$ and $f(k) = \nu$ otherwise.
Let $S$ be the set of all $\sigma \in (k+1)^{<\omega}$ which $f$-refine some $\tau \in T \cap k^{|\sigma|}$
and such that for every $w < |\sigma|$, every part $\mu \in \{\nu, k\}$
and every finite $R$-transitive set $F' \subseteq \dom(T) \cap \seto_\mu(\sigma)$, 
$\varphi(F_\nu \cup F', w)$ does not hold.

Note that $S$ is a p.r. partition tree of $[t, +\infty)$ refining $T$ with witness function~$f$. 
Suppose that $S$ is infinite. 
Let $\vec{E}$ be defined by $E_\mu = F_\mu$ if $\mu < k$ and $E_k = F_\nu$
and consider the extension $d = (\vec{E}, S, \Ccal[S])$.
We claim that~$\nu, k \not \in I(d)$. Fix a part $\mu \in \{\nu, k\}$ of $S$. By definition of $S$,
for every $\sigma \in S$, every $w < |\sigma|$
and every $R$-transitive set $F' \subseteq \dom(S) \cap \seto_\mu(\sigma)$,
$\varphi(E_\mu \cup F', w)$ does not hold.
Therefore, by clause~2 of Definition~\ref{def:em-forcing-precondition},
$d \Vdash_\mu (\forall x)\neg \psi(G, x)$, hence $d \Vdash_\mu \neg \varphi(G)$.
Note that $I(d) \subseteq I(c) \setminus \{\nu\}$.

Suppose now that $S$ is finite. Fix a threshold $\ell \in \omega$ such that $(\forall \sigma \in S)|\sigma| < \ell$
and a $\tau \in T \cap k^\ell$ such that $T^{[\tau]}$ is infinite.
Consider the 2-partition $E_0 \sqcup E_1$ of $\seto_\nu(\tau) \cap \dom(T)$ defined by
$E_0 = \{ i \geq t : \tau(i) = \nu \wedge (\forall^\infty s) R(i, s) \mbox{ holds}\}$
and $E_1 = \{ i \geq t : \tau(i) = \nu \wedge (\forall^\infty s) R(s, i) \mbox{ holds}\}$.
This is a 2-partition since the tournament~$R$ is stable.
As there exists no $\sigma \in S$ which $f$-refines $\tau$,
there exist a $w < \ell$ and an $R$-transitive set $F' \subseteq E_0$ or $F' \subseteq E_1$
such that $\varphi(F_\nu \cup F', w)$ holds. By choice of the partition,
there exists a $t' > t$ such that $F' \to_R [t', +\infty)$ or $[t', +\infty) \to_R F'$.
By Lemma~\ref{lem:emo-cond-valid}, $(F_\nu \cup F', [t', +\infty))$ is a valid EM extension for $(F_\nu, [t, +\infty))$.
As $T^{[\tau]}$ is infinite, $T^{[\nu, F']}$ is also infinite.
Let $\vec{E}$ be defined by $E_\mu = F_\mu$ if $\mu \neq \nu$ and $E_\mu = F_\nu \cup F'$ otherwise.
Let $S$ be the $k$-partition tree $(k, t', T^{[\nu, F']})$.
The condition $d = (\vec{E}, S, \Ccal[S])$ is a valid extension of $c$.
By clause~1 of Definition~\ref{def:em-forcing-precondition}, $d \Vdash_\mu \varphi(G)$.
Therefore $I(d) \subseteq I(c) \setminus \{\nu\}$.
\end{proof}

As in the previous notion of forcing, the following trivial lemma expresses the fact that the promise part of a condition
has no effect in the forcing relation for a $\Sigma^0_1$ or $\Pi^0_1$ formula.

\begin{lemma}\label{lem:em-promise-no-effect-first-level}
Fix two conditions $c = (\vec{F}, T, \Ccal)$ and $d = (\vec{F}, T, \Dcal)$, and a $\Sigma^0_1$ ($\Pi^0_1$)
formula. For every part $\nu$ of $T$,
$c \Vdash_\nu \varphi(G)$ if and only if $d \Vdash_\nu \varphi(G)$.
\end{lemma}
\begin{proof}
If $\varphi \in \Sigma^0_1$ then $\varphi(G)$ can be expressed as $(\exists x)\psi(G, x)$ 
where $\psi \in \Sigma^0_0$.
By clause~1 of Definition~\ref{def:em-forcing-precondition},
$c \Vdash_\nu \varphi(G)$ iff
there exists a $w \in \omega$ such that $\psi(F_\nu, w)$ holds,
iff $d \Vdash_\nu \varphi(G)$.
Similarily, if $\varphi \in \Pi^0_1$ then $\varphi(G)$ can be expressed as $(\forall x)\psi(G, x)$
where $\psi \in \Sigma^0_0$.
By clause~2 of Definition~\ref{def:em-forcing-precondition},
$c \Vdash_\nu \varphi(G)$ iff
for every $\sigma \in T$, every $w < |\sigma|$ and 
every $R$-transitive set $F' \subseteq \dom(T) \cap \seto_\nu(\sigma)$, $\varphi(F_\nu \cup F', w)$ holds,
iff $d \Vdash_\nu \varphi(G)$.
\end{proof}

We are now ready to define the forcing relation for an arbitrary arithmetic formula.
Again, the natural forcing relation induced by the forcing of~$\Sigma^0_0$ formulas
is too complex, so we design a more effective relation which still
enjoys the main properties of a forcing relation.

\begin{definition}\label{def:em-forcing-condition}
Fix a condition $c = (\vec{F}, T, \Ccal)$ and an arithmetic formula $\varphi(G)$.
\begin{itemize}
	\item[1.] If $\varphi(G) = (\exists x)\psi(G, x)$ where $\psi \in \Pi^0_1$ then
	$c \Vdash \varphi(G)$ iff for every part $\nu < \parts(T)$ such that
	$(\nu, T) \in \Ccal$ there exists a $w < \dom(T)$ such that $c \Vdash_\nu \psi(G, w)$
	\item[2.] If $\varphi(G) = (\forall x)\psi(G, x)$ where $\psi \in \Sigma^0_1$ then
	$c \Vdash \varphi(G)$ iff for every infinite p.r.\ $k'$-partition tree $S$,
	every function~$f : \parts(S) \to \parts(T)$,
	every $w$ and $\vec{E}$ smaller than $\#S$ such that the followings hold
	\begin{itemize}
		\item[i)] $(E_\nu, \dom(S))$ EM extends $(F_{f(\nu)}, \dom(T))$ for each $\nu < \parts(S)$
		\item[ii)] $S$ $f$-refines $\bigcap_{\nu < \parts(S)} T^{[f(\nu), E_\nu]}$
	\end{itemize}
	for every $(\mu, S) \in \Ccal$, $(\vec{E}, S, \Ccal[S]) \not \Vdash_\mu \neg \psi(G, w)$
	\item[3.] If $\varphi(G) = (\exists x)\psi(G, x)$ where $\psi \in \Pi^0_{n+2}$ then
	$c \Vdash \varphi(G)$ iff there exists a $w \in \omega$ such that $c \Vdash \psi(G, w)$
	\item[4.] If $\varphi(G) = \neg \psi(G,x)$ where $\psi \in \Sigma^0_{n+3}$
	then $c \Vdash \varphi(G)$ iff $d \not \Vdash \psi(G)$ for every $d \in \Ext(c)$.
\end{itemize}
\end{definition}

Notice that, unlike the forcing relation for $\Sigma^0_1$ and $\Pi^0_1$ formulas,
the relation over higher formuals does not depend on the part of the relation. 
The careful reader will have recognized the combinatorics of the second jump control 
introduced by Cholak, Jockusch and Slaman in~\cite{Cholak2001strength}.
We now prove the main properties of this forcing relation.

\begin{lemma}\label{lem:em-forcing-extension}
Fix a condition $c$ and a $\Sigma^0_{n+2}$ ($\Pi^0_{n+2}$) formula $\varphi(G)$.
If $c \Vdash \varphi(G)$ then for every $d \leq c$, $d \Vdash \varphi(G)$.
\end{lemma}
\begin{proof}
We prove the statement by induction over the complexity of the formula $\varphi(G)$.
Fix a condition $c = (\vec{F}, T, \Ccal)$ such that $c \Vdash \varphi(G)$ and an extension $d = (\vec{E}, S, \Dcal)$ of~$c$.
\begin{itemize}
	\item If $\varphi \in \Sigma^0_2$ then $\varphi(G)$ can be expressed as $(\exists x)\psi(G, x)$ 
  where $\psi \in \Pi^0_1$. By clause~1 of Definition~\ref{def:em-forcing-condition},
	for every part $\nu$ of $T$ such that $(\nu, T) \in \Ccal$, there exists a $w < \dom(T)$
	such that $c \Vdash_\nu \psi(G, w)$. Fix a part $\mu$ of $S$ such that $(\mu, S) \in \Dcal$.
	As $\Dcal \subseteq \Ccal$, $(\mu, S) \in \Ccal$. By upward-closure of $\Ccal$,
	part $\mu$ of $S$ refines some part $\nu$ of $\Ccal$ such that $(\nu, T) \in \Ccal$.
	Therefore by Lemma~\ref{lem:em-forcing-extension-level1}, $d \Vdash_\mu \psi(G, w)$,
	with $w < \dom(T) \leq \dom(S)$. Applying again clause~1 of Definition~\ref{def:em-forcing-condition},
	we deduce that $d \Vdash (\forall x)\psi(G, x)$, hence $d \Vdash \varphi(G)$.

  \item If $\varphi \in \Pi^0_2$ then $\varphi(G)$ can be expressed as $(\forall x)\psi(G, x)$ where $\psi \in \Sigma^0_1$.
	Suppose by way of contradiction that $d \not \Vdash (\forall x)\psi(G, x)$.
	Let $f : \parts(S) \to \parts(T)$ witness the refinement $S \leq T$.
	By clause~2 of Definition~\ref{def:em-forcing-condition}, there exist an infinite p.r.
	$k'$-partition tree $S'$, a function~$g : \parts(S') \to \parts(S)$, a $w \in \omega$, and $\vec{H}$
	smaller than the code of $S'$ such that
	\begin{itemize}
		\item[i)] $(H_\nu, \dom(S'))$ EM extends $(E_{g(\nu)}, \dom(S))$ for each $\nu < \parts(S')$
		\item[ii)] $S'$ $g$-refines $\bigcap_{\nu < \parts(S')} S^{[g(\nu), H_\nu]}$
		\item[iii)] there exists a $(\mu, S') \in \Dcal$ such that $(\vec{H}, S', \Dcal[S']) \Vdash_\mu \neg \psi(G, w)$.
	\end{itemize}
	To deduce by clause~2 of Definition~\ref{def:em-forcing-condition} that
	$c \not \Vdash (\forall x)\psi(G, x)$ and derive a contradiction, it suffices to prove
	that the same properties hold w.r.t.\ $T$.
	\begin{itemize}
		\item[i)] As by property (i) of the definition of an extension,
		$(E_{g(\nu)}, \dom(S))$ EM extends $(F_{f(g(\nu))}, \dom(T))$
		and $(H_\nu, \dom(S')$ EM extends $(E_{g(\nu)}, \dom(S))$,
		then $(H_\nu, \dom(S'))$ EM extends $(F_{f(g(\nu))}, \dom(T))$.
		\item[ii)] As by property (ii) of the definition of an extension,
		$S$ $f$-refines $\bigcap_{\nu < \parts(S')} T^{[f(g(\nu)), E_{g(\nu)}]}$
		and $S'$ $g$-refines $\bigcap_{\nu < \parts(S')} S^{[g(\nu), H_\nu]}$,
		then $S'$ $(g \circ f)$-refines $\bigcap_{\nu < \parts(S')} T^{[(g \circ f)(\nu), H_\nu]}$.
		\item[iii)] As $\Dcal \subseteq \Ccal$, there exists a part $(\mu, S') \in \Ccal$
		such that $(\vec{H}, S', \Dcal[S']) \Vdash_\mu \neg \psi(G, w)$. 
		By Lemma~\ref{lem:em-promise-no-effect-first-level}, $(\vec{H}, S', \Ccal[S']) \Vdash_\mu \neg \psi(G, w)$.
	\end{itemize}

	\item If $\varphi \in \Sigma^0_{n+3}$ then $\varphi(G)$ can be expressed as $(\exists x)\psi(G, x)$ 
  where $\psi \in \Pi^0_{n+2}$.
  By clause~3 of Definition~\ref{def:em-forcing-condition}, there exists a $w \in \omega$
  such that $c \Vdash \psi(G, w)$. By induction hypothesis, $d \Vdash \psi(G, w)$
  so by clause~3 of Definition~\ref{def:em-forcing-condition}, $d \Vdash \varphi(G)$.

  \item If $\varphi \in \Pi^0_{n+3}$ then $\varphi(G)$ can be expressed as $\neg \psi(G)$ where $\psi \in \Sigma^0_{n+3}$.
  By clause~4 of Definition~\ref{def:em-forcing-condition}, for every $e \in \Ext(c)$, $e \not \Vdash \psi(G)$.
  As $\Ext(d) \subseteq \Ext(c)$, for every $e \in \Ext(d)$, $e \not \Vdash \psi(G)$,
  so by clause~4 of Definition~\ref{def:em-forcing-condition}, $d \Vdash \varphi(G)$.
\end{itemize}
\end{proof}

\begin{lemma}\label{lem:em-forcing-dense}
For every $\Sigma^0_{n+2}$ ($\Pi^0_{n+2}$) formula $\varphi$, the following set is dense
$$
\{c \in \Pb : c \Vdash \varphi(G) \mbox{ or } c \Vdash \neg \varphi(G) \}
$$
\end{lemma}
\begin{proof}
We prove the statement by induction over~$n$.
It suffices to treat the case where $\varphi$ is a $\Sigma^0_{n+2}$ formula,
as the case where $\varphi$ is a $\Pi^0_{n+2}$ formula is symmetric. Fix a condition $c = (\vec{F}, T, \Ccal)$.
\begin{itemize}
	\item In case $n = 0$, the formula $\varphi$ is of the form $(\exists x)\psi(G, x)$ where
	$\psi \in \Pi^0_1$.
	Suppose there exist an infinite p.r. $k'$-partition tree $S$
	for some $k' \in \omega$, a function~$f : \parts(S) \to \parts(T)$ and a $k'$-tuple of finite sets $\vec{E}$ such that
	\begin{itemize}
		\item[i)] $(E_\nu, [\ell, +\infty))$ EM extends $(F_{f(\nu)}, \dom(T))$ for each $\nu < \parts(S)$.
		\item[ii)] $S$ $f$-refines $\bigcap_{\nu < \parts(S)} T^{[f(\nu), E_\nu]}$
		\item[iii)] for each non-empty part $\nu$ of $S$ such that $(\nu, S) \in \Ccal$, 
		$(\vec{E}, S, \Ccal[S]) \Vdash_\nu \psi(G, w)$ for some $w < \#S$
	\end{itemize}
	We can choose $\dom(S)$  so that $(E_\nu, \dom(S))$ EM extends $(F_{f(\nu)}, \dom(T))$ for each $\nu < \parts(S)$.
	Properties i-ii) remain trivially true.
	By Lemma~\ref{lem:em-forcing-extension-level1} and Lemma~\ref{lem:em-promise-no-effect-first-level}, property iii) remains true too.
	Let $\Dcal = \Ccal[S] \setminus \{(\nu, S') \in \Ccal : \mbox{ part } \nu \mbox{ of } S' \mbox{ is empty} \}$. 
	As $\Ccal$ is an $\emptyset'$-p.r. promise for $T$,
	$\Ccal[S]$ is an $\emptyset'$-p.r. promise for $S$.
	As $\Dcal$ is obtained from $\Ccal[S]$ by removing only empty parts, $\Dcal$ is also an $\emptyset'$-p.r. promise for $S$.
	By clause~1 of Definition~\ref{def:em-forcing-condition},
	$d = (\vec{E}, S, \Dcal) \Vdash (\exists x)\psi(G, x)$ hence $d \Vdash \varphi(G)$.

	We may choose a coding of the p.r. trees such that
	the code of $S$ is sufficiently large to witness $\ell$ and $\vec{E}$.
	So suppose now that for every infinite p.r. $k'$-partition tree $S$,
	every function~$f : \parts(S) \to \parts(T)$ and $\vec{E}$ smaller than the code of $S$ such that properties i-ii) hold,
	there exists a non-empty part $\nu$ of $S$ such that $(\nu, S) \in \Ccal$
	and $(\vec{E}, S, \Ccal) \not \Vdash_\nu \psi(G, w)$ for every $w < \ell$.
	Let $\Dcal$ be the collection of all such $(\nu, S)$. The set $\Dcal$ is $\emptyset'$-p.r.\ 
	since by Lemma~\ref{lem:em-complexity-forcing}, both $(\vec{E}, S, \Ccal) \not \Vdash_\nu \psi(G, w)$
	and ``part $\nu$ of~$S$ is non-empty'' are $\Sigma^0_1$.
	By Lemma~\ref{lem:em-forcing-extension-level1} and since we require that~$\#S \geq \#T$
	in the definition of~$S \leq T$, $\Dcal$ is upward-closed under the refinement relation,
	hence is a promise for~$T$. By clause~2
	of Definition~\ref{def:em-forcing-condition}, $d = (\vec{F}, T, \Dcal) \Vdash (\forall x) \neg \psi(G, x)$,
	hence $d \Vdash \neg \varphi(G)$.

	\item In case $n > 0$, density follows from clause~4 of Definition~\ref{def:em-forcing-condition}.
\end{itemize}
\end{proof}

Given any filter~$\Fcal = \{c_0, c_1, \dots \}$ with $c_s = (\vec{F}_s, T_s, \Ccal_s)$, 
the set of the acceptable parts~$\nu$ of~$T_s$ such that~$(\nu, T_s) \in \Ccal_s$ forms
an infinite, directed acyclic graph~$\Gcal(\Fcal)$. Whenever~$\Fcal$ is sufficiently generic,
the graph~$\Gcal(\Fcal)$ has a unique infinite path~$P$.
The path~$P$ induces an infinite set~$G = \bigcup_s F_{P(s), s}$. 
We call~$P$ the \emph{generic path} and $G$ the \emph{generic real}.

\begin{lemma}\label{lem:em-generic-level1}
Suppose that $\Fcal$ is sufficiently generic and let~$P$ and~$G$ be the generic path and the generic real, respectively.
For any $\Sigma^0_1$ ($\Pi^0_1$) formula $\varphi(G)$, 
$\varphi(G)$ holds iff $c_s \Vdash_{P(s)} \varphi(G)$ for some $c_s \in \Fcal$.
\end{lemma}
\begin{proof}
Fix a condition $c_s = (\vec{F}, T, \Ccal) \in \Fcal$ such that $c \Vdash_{P(s)} \varphi(G)$,
and let~$\nu = P(s)$. 
\begin{itemize}
	\item If $\varphi \in \Sigma^0_1$ then $\varphi(G)$ can be expressed as $(\exists x)\psi(G, x)$ where $\psi \in \Sigma^0_0$.
	By clause~1 of Definition~\ref{def:em-forcing-precondition}, there exists a $w \in \omega$ such that
	$\psi(F_\nu, w)$ holds. As $\nu = P(s)$, $F_\nu = F_{P(s)} \subseteq G$
	and~$G \setminus F_\nu \subseteq (max(F_\nu), +\infty)$, so $\psi(G, w)$ holds by continuity, hence $\varphi(G)$ holds.
	
	\item If $\varphi \in \Pi^0_1$ then $\varphi(G)$ can be expressed as $(\forall x)\psi(G, x)$ where $\psi \in \Sigma^0_0$.
	By clause~2 of Definition~\ref{def:em-forcing-precondition}, for every $\sigma \in T$, every $w < |\sigma|$
	and every $R$-transitive set $F' \subseteq \dom(T) \cap \seto_\nu(\sigma)$, $\psi(F_\nu \cup F', w)$ holds. 
	For every $F' \subseteq G \setminus F_\nu$, and $w \in \omega$ there exists a $\sigma \in T$
	such that $w < |\sigma|$ and $F' \subseteq \dom(T) \cap \seto_\nu(\sigma)$. Hence $\psi(F_\nu \cup F', w)$ holds.
	Therefore, for every $w \in \omega$, $\psi(G, w)$ holds, so $\varphi(G)$ holds.
\end{itemize}
The other direction holds by Lemma~\ref{lem:em-forcing-dense-level1}.
\end{proof}

\begin{lemma}
Suppose that $\Fcal$ is sufficiently generic and let~$P$ and~$G$ be the generic path and the generic real, respectively.
For any $\Sigma^0_{n+2}$ ($\Pi^0_{n+2}$) formula $\varphi(G)$, 
$\varphi(G)$ holds iff $c_s \Vdash \varphi(G)$ for some $c_s \in \Fcal$.
\end{lemma}
\begin{proof}
We prove the statement by induction over the complexity of the formula $\varphi(G)$.
As previously noted in Lemma~\ref{lem:coh-holds-filter}, it suffices to prove that if $c_s \Vdash \varphi(G)$ for some
$c_s \in \Fcal$ then $\varphi(G)$ holds. Indeed, conversely if $\varphi(G)$ holds,
then by Lemma~\ref{lem:em-forcing-dense} and by genericity of $\Fcal$ either $c_s \Vdash \varphi(G)$ or $c_s \Vdash \neg \varphi(G)$,
but if $c \Vdash \neg \varphi(G)$ then $\neg \varphi(G)$ holds, contradicting 
the hypothesis. So $c_s \Vdash \varphi(G)$.
Fix a condition $c_s = (\vec{F}, T, \Ccal) \in \Fcal$ such that $c_s \Vdash \varphi(G)$. 
We proceed by case analysis on $\varphi$. 
\begin{itemize}
	\item If $\varphi \in \Sigma^0_2$ then $\varphi(G)$ can be expressed as $(\exists x)\psi(G, x)$ 
  where $\psi \in \Pi^0_1$.
  By clause~1 of Definition~\ref{def:em-forcing-condition}, for every part $\nu$ of $T$
	such that $(\nu, T) \in \Ccal$, there exists a $w < \dom(T)$
  such that $c_s \Vdash_\nu \psi(G, w)$. In particular $(P(s), T) \in \Ccal$, so $c_s \Vdash_{P(s)} \psi(G, w)$.
	By Lemma~\ref{lem:em-generic-level1}, $\psi(G, w)$ holds, hence $\varphi(G)$ holds.

	\item If $\varphi \in \Pi^0_2$ then $\varphi(G)$ can be expressed as $(\forall x)\psi(G, x)$ where $\psi \in \Sigma^0_1$.
	By clause~2 of Definition~\ref{def:em-forcing-condition}, for every infinite $k'$-partition
	tree $S$, every function~$f : \parts(S) \to \parts(T)$, 
	every $w$ and $\vec{E}$ smaller than the code of $S$ such that the followings hold
	\begin{itemize}
		\item[i)] $(E_\nu, \dom(S))$ EM extends $(F_{f(\nu)}, \dom(T))$ for each $\nu < \parts(S)$
		\item[ii)] $S$ $f$-refines $\bigcap_{\nu < \parts(S)} T^{[f(\nu), E_\nu]}$
	\end{itemize}
	for every $(\mu, S) \in \Ccal$, $(\vec{E}, S, \Ccal[S]) \not \Vdash_\mu \neg \psi(G, w)$.
	Suppose by way of contradiction that $\psi(G, w)$ does not hold for some $w \in \omega$.
  Then by Lemma~\ref{lem:em-generic-level1}, there exists a $d_t \in \Fcal$
	such that $d_t \Vdash_{P(t)} \neg \psi(G, w)$.
	Since~$\Fcal$ is a filter, there is a condition $e_r = (\vec{E}, S, \Dcal) \in \Fcal$
	extending both~$c_s$ and~$d_t$.
	Let~$\mu = P(r)$. By choice of~$P$, $(\mu, S) \in \Ccal$, so
	by clause ii), $(\vec{E}, S, \Ccal[S]) \not \Vdash_\mu \psi(G, w)$,
	hence by Lemma~\ref{lem:em-promise-no-effect-first-level}, $e_r \not \Vdash_\mu \neg \psi(G, w)$.
	However, since part~$\mu$ of~$S$ refines part~$P(t)$ of~$d_t$,
	then by Lemma~\ref{lem:em-forcing-extension-level1}, $e_r \Vdash_\mu \neg \psi(G, w)$. Contradiction.
	Hence for every $w \in \omega$, $\psi(G, w)$ holds, so $\varphi(G)$ holds.

	\item If $\varphi \in \Sigma^0_{n+3}$ then $\varphi(G)$ can be expressed as $(\exists x)\psi(G, x)$ 
  where $\psi \in \Pi^0_{n+2}$.
  By clause~3 of Definition~\ref{def:em-forcing-condition}, there exists a $w \in \omega$
  such that $c_s \Vdash \psi(G, w)$. By induction hypothesis, $\psi(G, w)$ holds, hence $\varphi(G)$ holds.

  Conversely, if $\varphi(G)$ holds, then there exists a $w \in \omega$ such that $\psi(G, w)$ holds,
  so by induction hypothesis $c_s \Vdash \psi(G, w)$ for some $c_s \in \Fcal$,
  so by clause~3 of Definition~\ref{def:em-forcing-condition}, $c_s \Vdash \varphi(G)$.

	\item If $\varphi \in \Pi^0_{n+3}$ then $\varphi(G)$ can be expressed as $\neg \psi(G)$ where $\psi \in \Sigma^0_{n+3}$.
  By clause~4 of Definition~\ref{def:em-forcing-condition}, for every $d \in \Ext(c_s)$, $d \not \Vdash \psi(G)$.
	By Lemma~\ref{lem:em-forcing-extension}, $d \not \Vdash \psi(G)$ for every~$d \in \Fcal$
 	and by a previous case, $\psi(G)$ does hold, so $\varphi(G)$ holds.
\end{itemize}
\end{proof}

We now prove that the forcing relation has good definitional properties
as we did with the notion of forcing for cohesiveness.

\begin{lemma}\label{lem:em-extension-complexity}
For every condition $c$, $\Ext(c)$ is $\Pi^0_2$ uniformly in~$c$.
\end{lemma}
\begin{proof}
Recall from Lemma~\ref{lem:em-refinement-complexity} that
given $k, t \in \omega$, $PartTree(k, t)$ denotes 
the $\Pi^0_1$ set of all the infinite p.r.\ $k$-partition trees of $[t, +\infty)$,
and given a $k$-partition tree $S$ and a part $\nu$ of $S$,
the predicate $Empty(S, \nu)$ denotes the $\Pi^0_1$ formula ``part $\nu$ of $S$ is empty'',
that is, the formula $(\forall \sigma \in S)[\seto_\nu(\sigma) \cap \dom(S) = \emptyset]$.
If $T$ is p.r. then so is $T^{[\nu, H]}$ for some finite set $H$.

Fix a condition $c = (\vec{F}, (k, t, T), \Ccal)$.
By definition, $(\vec{H}, (k', t', S), \Dcal) \in \Ext(c)$ iff the following formula holds:
$$
\begin{array}{l@{\hskip 0.5in}r}
(\exists f : k' \to k)\\
(\forall \nu < k')(H_\nu, [t', +\infty)) \mbox{ EM extends } (F_{f(\nu)}, [t, +\infty)) & (\Pi^0_1)\\
\wedge S \in PartTree(k', t') \wedge S \leq_f \bigwedge_{\nu < k'} T^{[f(\nu), H_{\nu}]} & (\Pi^0_1)\\
\wedge \Dcal \mbox{ is a promise for } S \wedge \Dcal \subseteq \Ccal & (\Pi^0_2)\\
\end{array}
$$
By Lemma~\ref{lem:em-refinement-complexity}
and the fact that $\bigwedge_{\nu < k'} T^{[f(\nu), H_{\nu}]}$
is p.r. uniformly in $T$, $f$, $\vec{H}$ and $k'$,
the above formula is $\Pi^0_2$.
\end{proof}

\begin{lemma}\label{lem:em-complexity-forcing}
Fix an arithmetic formula $\varphi(G)$, a condition $c = (\vec{F}, T, \Ccal)$
and a part $\nu$ of $T$.
\begin{itemize}
	\item[a)] If $\varphi(G)$ is a $\Sigma^0_1$ ($\Pi^0_1$) formula 
	then so is the predicate $c \Vdash_\nu \varphi(G)$.
	\item[b)] If $\varphi(G)$ is a $\Sigma^0_{n+2}$ ($\Pi^0_{n+2}$) formula 
	then so is the predicate $c \Vdash \varphi(G)$.
\end{itemize}
\end{lemma}
\begin{proof}
We prove our lemma by induction over the complexity of the formula $\varphi(G)$.
\begin{itemize}
	\item If $\varphi(G) \in \Sigma^0_1$ then it can be expressed as $(\exists x)\psi(G, x)$ where $\psi \in \Sigma^0_0$.
	By clause~1 of Definition~\ref{def:em-forcing-precondition}, $c \Vdash_\nu \varphi(G)$ if and only if 
	the formula $(\exists w \in \omega)\psi(F_\nu, w)$ holds. This is a $\Sigma^0_1$ predicate.
	
	\item If $\varphi(G) \in \Pi^0_1$ then it can be expressed as $(\forall x)\psi(G, x)$ where $\psi \in \Sigma^0_0$.
	By clause~2 of Definition~\ref{def:em-forcing-precondition}, $c \Vdash_\nu \varphi(G)$ if and only if 
	the formula $(\forall \sigma \in T)(\forall w < |\sigma|)(\forall F' \subseteq \dom(T) \cap \seto_\nu(\sigma))
	[F'\ R\mbox{-transitive} \imp \psi(F_\nu \cup F', w)]$ holds. This is a $\Pi^0_1$ predicate.

	\item If $\varphi(G) \in \Sigma^0_2$ then it can be expressed as $(\exists x)\psi(G, x)$ where $\psi \in \Pi^0_1$.
	By clause~1 of Definition~\ref{def:em-forcing-condition}, $c \Vdash \varphi(G)$ if and only if 
	the formula $(\forall \nu < \parts(T))(\exists w < \dom(T))[(\nu, T) \in \Ccal \imp c \Vdash_\nu \psi(G, w)]$ holds.
	This is a $\Sigma^0_2$ predicate by induction hypothesis and the fact that $\Ccal$ is $\emptyset'$-computable.

	\item If $\varphi(G) \in \Pi^0_2$ then it can be expressed as $(\forall x)\psi(G, x)$ where $\psi \in \Sigma^0_1$.
	By clause~2 of Definition~\ref{def:em-forcing-condition}, $c \Vdash \varphi(G)$ if and only if 
	for every infinite $k'$-partition tree $S$, every function~$f : \parts(S) \to \parts(T)$,
	every $w$ and $\vec{E}$ smaller than the code of $S$ such that the followings hold
	\begin{itemize}
		\item[i)] $(E_\nu, \dom(S))$ EM extends $(F_{f(\nu)}, \dom(T))$ for each $\nu < \parts(S)$
		\item[ii)] $S$ $f$-refines $\bigcap_{\nu < \parts(S)} T^{[f(\nu), E_\nu]}$
	\end{itemize}
	for every $(\mu, S) \in \Ccal$, $(\vec{E}, S, \Ccal[S]) \not \Vdash_\mu \neg \psi(G, w)$.
	By Lemma~\ref{lem:em-refinement-complexity}, Properties i-ii) are $\Delta^0_2$.
	Moreover, the predicate $(\mu, S) \in \Ccal$ is $\Delta^0_2$.
	By induction hypothesis, $(\vec{E}, S, \Ccal) \not \Vdash_\mu \neg \psi(G, w)$ is $\Sigma^0_1$.
	Therefore $c \Vdash \varphi(G)$ is a $\Pi^0_2$ predicate.

	\item If $\varphi(G) \in \Sigma^0_{n+3}$ then it can be expressed as $(\exists x)\psi(G, x)$ where $\psi \in \Pi^0_{n+2}$.
	By clause~3 of Definition~\ref{def:em-forcing-condition}, $c \Vdash \varphi(G)$ if and only if 
	the formula $(\exists w \in \omega)c \Vdash \psi(G, w)$ holds. This is a $\Sigma^0_{n+3}$ predicate
	by induction hypothesis.

	\item If $\varphi(G) \in \Pi^0_{n+3}$ then it can be expressed as $\neg \psi(G)$ where $\psi \in \Sigma^0_{n+3}$. 
	By clause~4 of Definition~\ref{def:em-forcing-condition}, $c \Vdash \varphi(G)$ if and only if 
	the formula $(\forall d)(d \not \in \Ext(c) \vee d \not \Vdash \psi(G))$ holds.
	By induction hypothesis, $d \not \Vdash \psi(G)$ is a $\Pi^0_{n+3}$ predicate.
	By Lemma~\ref{lem:em-extension-complexity}, the set $\Ext(c)$ is $\Pi^0_2$-computable uniformly in $c$,
	thus $c \Vdash \varphi(G)$ is a $\Pi^0_{n+3}$ predicate.
\end{itemize}
\end{proof}

\subsection{Preserving the arithmetic hierarchy}

We now prove the core lemmas showing that every sufficiently generic real
preserves the arithmetic hierarchy. The proof is split into two lemmas
since the forcing relation for $\Sigma^0_1$ and~$\Pi^0_1$ formulas
depends on the part of the condition, and therefore has to be treated separately.

\begin{lemma}\label{lem:em-diagonalization-level1}
If $A \not \in \Sigma^0_1$ and $\varphi(G, x)$ is $\Sigma^0_1$,
then the set of $c = (\vec{F}, T, \Ccal) \in \Pb$ satisfying the following property is dense:
$$
(\forall \nu < \parts(T))[(\exists w \in A)c_s \Vdash_\nu \neg \varphi(G, w)] \vee [(\exists w \not \in A)c_s \Vdash_\nu \varphi(G, w)]
$$
\end{lemma}
\begin{proof}
The formula $\varphi(G, w)$ can be expressed as $(\exists x)\psi(G, w, x)$ where $\psi \in \Sigma^0_0$.
Given a condition $c = (\vec{F}, T, \Ccal)$, let $I(c)$ be the set of the parts $\nu$ of $T$
such that for every $w \in A$, $c \not \Vdash_\nu \neg \varphi(G, w)$
and for every $w \in \overline{A}$, $c \not \Vdash_\nu \varphi(G, w)$.
If $I(c) = \emptyset$ then we are done, so suppose $I(c) \neq \emptyset$ and
fix some $\nu \in I(c)$. We will construct an extension $d$ such that 
$I(d) \subseteq I(c) \setminus \{\nu\}$. Iterating the operation completes the proof.

Say that $T$ is a $k$-partition tree of $[t, +\infty)$ for some $k, t \in \omega$.
Define $f : k+1 \to k$ as $f(\mu) = \mu$ if $\mu < k$ and $f(k) = \nu$ otherwise.
Given an integer $w \in \omega$, let $S_w$ be the set of all $\sigma \in (k+1)^{<\omega}$
which $f$-refine some $\tau \in T \cap k^{|\sigma|}$ and such that for every $u < |\sigma|$,
every part $\mu \in \{\nu, k\}$ and every finite $R$-transitive set $F' \subseteq \dom(T) \cap \seto_\mu(\sigma)$,
$\varphi(F_\nu \cup F', w, u)$ does not hold.

The set $S_w$ is a p.r.\ (uniformly in $w$) partition tree of $[t, +\infty)$ refining $T$ with witness function~$f$.
Let $U = \{ w \in \omega : S_w \mbox{ is finite } \}$. $U \in \Sigma^0_1$, thus $U \neq A$.
Fix some $w \in U \Delta A$. Suppose first that $w \in A \setminus U$. By definition of $U$, $S_w$ is infinite. 
Let $\vec{E}$ be defined by $E_\mu = F_\mu$ if $\mu < k$ and $E_k = F_\nu$,
and consider the extension $d = (\vec{E}, S_w, \Ccal[S_w])$. We claim that $I(d) \subseteq I(c) \setminus \{\nu\}$. 
Fix a part $\mu \in \{\nu, k\}$ of $S_w$. By definition of $S_w$,
for every $\sigma \in S_w$, every $u < |\sigma|$ and every $R$-transitive set $F' \subseteq \dom(S_w) \cap \seto_\mu(\sigma)$,
$\varphi(E_\mu \cup F', w, u)$ does not hold. Therefore, by clause~2 of Definition~\ref{def:em-forcing-precondition},
$d \Vdash_\mu (\forall x)\neg \psi(G, w, x)$, hence $d \Vdash_\mu \neg \varphi(G, w)$, and this for some $w \in A$.
Thus $I(d) \subseteq I(c) \setminus \{\nu\}$. 

Suppose now that $w \in U \setminus A$, so $S_w$ is finite.
Fix an $\ell \in \omega$ such that $(\forall \sigma \in S)|\sigma| < \ell$
and a $\tau \in T \cap k^\ell$ such that $T^{[\tau]}$ is infinite.
Consider the 2-partition $E_0 \cup E_1$ of $\seto_\nu(\tau) \cap \dom(T)$
defined by $E_0 = \{i \geq t : \tau(i) = \nu \wedge (\forall^\infty s) R(i, s) \mbox{ holds}\}$
and $E_0 = \{i \geq t : \tau(i) = \nu \wedge (\forall^\infty s) R(s, i) \mbox{ holds}\}$.
As there exists no $\sigma \in S_w$ which $f$-refines $\tau$, there exist a $u < \ell$
and an $R$-transitive set $F' \subseteq E_0$ or $F' \subseteq E_1$ such that $\varphi(F_\nu \cup F', w, u)$ holds.
By choice of the partition, there exists a $t' > t$ such that $F' \to_R [t', +\infty)$
or $[t', +\infty) \to_R F'$. 
By Lemma~\ref{lem:emo-cond-valid}, $(F_\nu \cup F', [t', +\infty))$ is a valid EM extension of $(F_\nu, [t, +\infty))$.
As $T^{[\tau]}$ is infinite, $T^{[\nu, F']}$ is also infinite.
Let $\vec{E}$ be defined by $E_\mu = F_\mu$ if $\mu \neq \nu$ and $E_\mu = F_\nu \cup F'$ otherwise.
Let $S$ be the $k$-partition tree $(k, t', T^{[\nu, F']})$.
The condition $d = (\vec{E}, S, \Ccal[S])$ is a valid extension of $c$.
By clause~1 of Definition~\ref{def:em-forcing-precondition}, $d \Vdash_\mu \varphi(G, w)$ with $w \not \in A$. .
Therefore $I(d) \subseteq I(c) \setminus \{\nu\}$.
\end{proof}

\begin{lemma}\label{lem:em-diagonalization}
If $A \not \in \Sigma^0_{n+2}$ and $\varphi(G, x)$ is $\Sigma^0_{n+2}$,
then the set of $c \in \Pb$ satisfying the following property is dense:
$$
[(\exists w \in A)c \Vdash \neg \varphi(G, w)] \vee [(\exists w \not \in A)c \Vdash \varphi(G, w)]
$$
\end{lemma}
\begin{proof}
Fix a condition $c = (\vec{F}, T, \Ccal)$.
\begin{itemize}
	\item In case $n = 0$, $\varphi(G, w)$ can be expressed as $(\exists x)\psi(G, w, x)$ where $\psi \in \Pi^0_1$.
	Let $U$ be the set of integers $w$ such that there exist an infinite p.r. $k'$-partition tree $S$
	for some $k' \in \omega$, a function~$f : \parts(S) \to \parts(T)$ and a $k'$-tuple of finite sets $\vec{E}$ such that
	\begin{itemize}
		\item[i)] $(E_\nu, [\ell, +\infty))$ EM extends $(F_{f(\nu)}, \dom(T))$ for each $\nu < \parts(S)$.
		\item[ii)] $S$ $f$-refines $\bigcap_{\nu < \parts(S)} T^{[f(\nu), E_\nu]}$
		\item[iii)] for each non-empty part $\nu$ of $S$ such that $(\nu, S) \in \Ccal$, 
		$(\vec{E}, S, \Ccal[S]) \Vdash_\nu \psi(G, w, u)$ for some $u < \#S$
	\end{itemize}
	By Lemma~\ref{lem:em-complexity-forcing} and Lemma~\ref{lem:em-refinement-complexity},
	$U \in \Sigma^0_2$, thus $U \neq A$. Let $w \in U \Delta A$.
	Suppose that $w \in U \setminus A$.
	We can choose $\dom(S)$  so that $(E_\nu, \dom(S))$ EM extends $(F_{f(\nu)}, \dom(T))$ for each $\nu < \parts(S)$.
	By Lemma~\ref{lem:em-forcing-extension-level1} and Lemma~\ref{lem:em-promise-no-effect-first-level}, properties i-ii) remain true.
	Let $\Dcal = \Ccal[S] \setminus \{(\nu, S') \in \Ccal : \mbox{ part } \nu \mbox{ of } S' \mbox{ is empty} \}$. 
	As $\Ccal$ is an $\emptyset'$-p.r. promise for $T$,
	$\Ccal[S]$ is an $\emptyset'$-p.r. promise for $S$.
	As $\Dcal$ is obtained from $\Ccal[S]$ by removing only empty parts, $\Dcal$ is also an $\emptyset'$-p.r. promise for $S$.
	By clause~1 of Definition~\ref{def:em-forcing-condition},
	$d = (\vec{E}, S, \Dcal) \Vdash (\exists x)\psi(G, w, x)$ hence $d \Vdash \varphi(G, w)$ for some $w \not \in A$.

	We may choose a coding of the p.r. trees such that
	the code of $S$ is sufficiently large to witness $u$ and $\vec{E}$.
	So suppose now that $w \in A \setminus U$. Then for every infinite p.r. $k'$-partition tree $S$,
	every $\ell$ and $\vec{E}$ smaller than the code of $S$ such that properties i-ii) hold,
	there exists a non-empty part $\nu$ of $S$ such that $(\nu, S) \in \Ccal$
	and $(\vec{E}, S, \Ccal) \not \Vdash_\nu \psi(G, w, u)$ for every $u < \ell$.
	Let $\Dcal$ be the collection of all such $(\nu, S)$. The set $\Dcal$ is $\emptyset'$-p.r.
	By Lemma~\ref{lem:em-forcing-extension-level1} and since~$\#S \geq \#T$
	whenever~$S \leq_f T$, $\Dcal$ is upward-closed under the refinement relation,
	hence it is a promise for~$T$. By clause~2.
	of Definition~\ref{def:em-forcing-condition}, $d = (\vec{F}, T, \Dcal) \Vdash (\forall x) \neg \psi(G, w, x)$,
	hence $d \Vdash \neg \varphi(G, w)$ for some $w \in A$.

	\item In case $n > 0$, let $U = \{ w \in \omega : (\exists d \in \Ext(c)) d \Vdash \varphi(G, w) \}$.
	By Lemma~\ref{lem:em-extension-complexity} and Lemma~\ref{lem:em-complexity-forcing},
	$U \in \Sigma^0_{n+2}$, thus $U \neq A$.
	Fix some $w \in U \Delta A$. If $w \in U \setminus A$ then by definition of~$U$,
	there exists a condition $d$ extending $c$ such that $d \Vdash \varphi(G, w)$.
	If $w \in A \setminus U$, then for every $d \in \Ext(c)$, $d \not \Vdash \varphi(G, w)$
	so by clause~4 of Definition~\ref{def:em-forcing-condition}, $c \Vdash \neg \varphi(G, w)$. 
\end{itemize}
\end{proof}

We are now ready to prove Theorem~\ref{thm:em-preserves-arithmetic}.
It follows from the preservation of the arithmetic hierarchy
for cohesiveness and the stable Erd\H{o}s-Moser theorem.

\begin{proof}[Proof of Theorem~\ref{thm:em-preserves-arithmetic}]
Since $\rca \vdash \coh \wedge \semo \imp \emo$, then by Theorem~\ref{thm:coh-preservation-arithmetic-hierarchy}
it suffices to prove that $\semo$ admits preservation of the arithmetic hierarchy.
Fix some set~$C$ and a $C$-computable stable infinite tournament~$R$.
Let $\Ccal_0$ be the $C'$-p.r. set of all $(\nu, T) \in \TPb$ such that $(\nu, T) \leq (0, 1^{<\omega})$.
Let~$\Fcal$ be a sufficiently generic filter containing~$c_0 = (\{\emptyset\}, 1^{<\omega}, \Ccal_0)$.
Let~$P$ and $G$ be the corresponding generic path and generic real, respectively.
By definition of a condition, the set~$G$ is $R$-transitive.
By Lemma~\ref{lem:em-forcing-infinite}, $G$ is infinite. By Lemma~\ref{lem:em-diagonalization-level1}
and Lemma~\ref{lem:em-complexity-forcing}, $G$ preserves non-$\Sigma^0_1$ definitions relative to~$C$.
By Lemma~\ref{lem:em-diagonalization} and Lemma~\ref{lem:em-complexity-forcing},
$G$ preserves non-$\Sigma^0_{n+2}$ definitions relative to~$C$ for every~$n \in \omega$.
Therefore, by Proposition 2.2 of~\cite{Wang2014Definability}, $G$
preserves the arithmetic hierarchy relative to~$C$.
\end{proof}

\section{An effective forcing for stable Ramsey's theorem for pairs}

Among the Ramsey-type hierarchies, the $\mathsf{D}$ hierarchy
is conceptually the simplest one. It is therefore natural
to study it in order to understand better the control of iterated jumps
and focus on the core combinatorics without the technicalities specific to
another hierarchy. Recall that for every~$n, k \geq 1$, 
$\mathsf{D}^n_k$ is the statement ``Every~$\Delta^0_n$ $k$-partition of the integers has an infinite subset in of its parts''.
Wang~\cite{Wang2014Definability} studied $\mathsf{D}^2_2$ within his framework of preservation of definitions
and proved that $\mathsf{D}^2_2$ admits preservation of $\Xi$ definitions 
simultaneously for all $\Xi$ in $\{\Sigma^0_{n+2}, \allowbreak \Pi^0_{n+2}, \allowbreak \Delta^0_{n+2} : n \in \omega \}$, 
but not~$\Delta^0_2$ definitions. He used for this a combination of the first jump control
of Cholak, Jockusch and Slaman~\cite{Cholak2001strength} and a relativization
of the preservation of the arithmetic hierarchy by~$\wkl$.

In this section, we design a notion of forcing for~$\mathsf{D}^2_2$
with a forcing relation which has the same definitional complexity
as the formula it forces. It enables us to reprove that $\mathsf{D}^2_2$
admits preservation of $\Xi$ definitions 
simultaneously for all $\Xi$ in $\{\Sigma^0_{n+2}, \allowbreak \Pi^0_{n+2}, \allowbreak \Delta^0_{n+2} : n \in \omega \}$.
The proof is significantly more involved than the previous proofs
of preservation of the arithmetic hierarchy. This is why we will only sketch the proof by explaining
its general structure and state the main lemmas without providing their proofs. We refer the
reader to~\cite{Patey2016Controlling} for the details of the proof.

\subsection{Sides of a sequence of sets}

A main feature in the construction of a solution to an instance $R_0, R_1$ of $\mathsf{D}^2_2$
is the parallel construction of a subset of $R_0$ and a subset of $R_1$.
The intrinsic disjunction in the forcing argument prevents us from applying
the same strategy as for the Erd\H{o}s-Moser theorem and obtain a preservation of the arithmetic hierarchy.
Given some~$\alpha < 2$, we shall refer to $R_\alpha$ or simply $\alpha$ as a \emph{side} of $\vec{R}$.
We also need to define a relative notion of acceptation and emptiness of a part.

\begin{definition}Fix a $k$-partition tree $T$ of $[t, +\infty)$ and a set $X$.
We say that part $\nu$ of $T$ is \emph{$X$-acceptable} if there exists a path $P$ through $T$
such that $\seto_\nu(P) \cap X$ is infinite.
We say that part $\nu$ of $T$ is \emph{$X$-empty} if 
$(\forall \sigma \in T)[\dom(T) \cap \seto_\nu(\sigma) \cap X = \emptyset]$.
\end{definition}

The intended uses of those notions will be $R_\alpha$-acceptation
and $R_\alpha$-emptiness. Every partition tree has an $R_\alpha$-acceptable part for some $\alpha < 2$.
The notion of $X$-emptiness is $\Pi^{0,X}_1$, and therefore $\Pi^0_2$ if $X$ is $\Delta^0_2$,
which raises new problems for obtaining a forcing relation of weak definitional complexity.
We would like to define a stronger notion of ``witnessing its acceptable parts'' and prove
that for every infinite p.r.\ partition tree~$T$, there is a p.r.\ refined tree~$S$
such that for each side $\alpha$ and each part~$\nu$ of~$S$, either~$\nu$ is $R_\alpha$-empty
in~$S$, or~$\nu$ is $R_\alpha$-acceptable. However, the resulting tree~$S$ would be $\emptyset'$-p.r.\
since~$R_\alpha$ is $\emptyset'$-computable. Thankfully, we will be able to circumvent this problem
in Lemma~\ref{lem:cohzp-validity-exists}.

\subsection{Forcing conditions}

Fix a $\Delta^0_2$ 2-partition~$R_0 \cup R_1 = \omega$.
We now describe the notion of forcing to build an infinite subset
of~$R_0$ or of~$R_1$.

\begin{definition}
We denote by~$\Pb$ the forcing notion whose conditions are tuples $(\vec{F}, T, \Ccal)$ where
\begin{itemize}
	\item[(a)] $T$ is an infinite, p.r.\ $k$-partition tree for some $k \in \omega$
	\item[(b)] $\Ccal$ is a $\emptyset'$-p.r.\ promise for $T$
	\item[(c)] $(F^\alpha_\nu, \dom(T))$ is a Mathias condition for each $\nu < k$ and $\alpha < 2$
\end{itemize}
A condition $d = (\vec{E}, S, \Dcal)$ \emph{extends} $c = (\vec{F}, T, \Ccal)$
(written $d \leq c$) if there exists a function $f : \parts(S) \to \parts(T)$ such that
$\Dcal \subseteq \Ccal$ and the followings hold
\begin{itemize}
	\item[(i)] $(E^\alpha_\nu, \dom(S) \cap R_\alpha)$ Mathias extends $(F^\alpha_{f(\nu)}, \dom(T) \cap R_\alpha)$ 
	for each $\nu < \parts(S)$ and $\alpha < 2$ 
	\item[(ii)] $S$ $f$-refines $\bigcap_{\nu < \parts(S), \alpha < 2} T^{[f(\nu), E^\alpha_\nu]}$
\end{itemize}
\end{definition}

In the whole construction, the index $\alpha$ indicates that we are constructing a set which will be almost included in $R_\alpha$.
Given a condition $c = (\vec{F}, T, \Ccal)$, we write again $\parts(c)$ for $\parts(T)$.
The following lemma shows that we can force our constructed set to be infinite
if we choose it among the acceptable parts.

\begin{lemma}\label{lem:cohzp-forcing-infinite}
For every condition $c = (\vec{F}, T, \Ccal)$ and every $n \in \omega$, there exists an extension $d = (\vec{E}, S, \Dcal)$
such that $|E^\alpha_\nu| \geq n$ on each $R_\alpha$-acceptable part $\nu$ of~$S$ for each $\alpha < 2$.
\end{lemma}

Given a condition~$c$, we denote by~$\Ext(c)$ the set of all its extensions.

\subsection{Forcing relation}

We need to define two forcing relations at the first level:
the ``true'' forcing relation, i.e., the one having the good density properties
but whose decision requires too much computational power, and a ``weak'' forcing relation
having better computational properties, but which does not behave well with respect to the forcing.
We start with the definition of the true forcing relation.

\begin{definition}[True forcing relation]\label{def:cohzp-true-forcing-precondition}
Fix a condition $c = (\vec{F}, T, \Ccal)$, a $\Sigma^0_0$ formula $\varphi(G, x)$,
a part $\nu < \parts(T)$, and a side $\alpha < 2$.
\begin{itemize}
	\item[1.] $c \Vvdash^\alpha_\nu (\exists x)\varphi(G, x)$ iff there exists a $w \in \omega$ such that $\varphi(F^\alpha_\nu, w)$ holds.
	\item[2.] $c \Vvdash^\alpha_\nu (\forall x)\varphi(G, x)$ iff 
	for every $\sigma \in T$ such that $T^{[\sigma]}$ is infinite, 
	every $w < |\sigma|$ and every set $F' \subseteq \dom(T) \cap \seto_\nu(\sigma) \cap R_\alpha$,
	$\varphi(F^\alpha_\nu \cup F', w)$ holds.
\end{itemize}
\end{definition}

Given a condition $c$, a side $\alpha < 2$, a part $\nu$ of $c$ and a $\Pi^0_1$ formula $\varphi$,
the relation $c \Vvdash^\alpha_\nu \varphi(G)$ is $\Pi^{0, \emptyset' \oplus R_\alpha}_1$,
hence $\Pi^0_2$ as $R_\alpha$ is $\Delta^0_2$. This relation
enjoys the good properties of a forcing relation, that is, it is downward-closed
under the refinement relation (Lemma~\ref{lem:cohzp-true-forcing-extension-level1}), 
and the set of the conditions forcing either a $\Sigma^0_1$
formula or its negation is dense (Lemma~\ref{lem:cohzp-true-forcing-dense-level1}).

\begin{lemma}\label{lem:cohzp-true-forcing-extension-level1}
Fix a condition $c = (\vec{F}, T, \Ccal)$ and a $\Sigma^0_1$ ($\Pi^0_1$) formula $\varphi(G)$.
If $c \Vvdash^\alpha_\nu \varphi(G)$ for some $\nu < \parts(T)$ and $\alpha < 2$,
then for every $d = (\vec{E}, S, \Dcal) \leq c$ and
every part $\mu$ of $S$ refining part $\nu$ of $T$, $d \Vvdash^\alpha_\mu \varphi(G)$.
\end{lemma}

\begin{lemma}\label{lem:cohzp-true-forcing-dense-level1}
For every $\Sigma^0_1$ ($\Pi^0_1$) formula $\varphi$, the following set is dense in $\Pb$:
$$
\{c \in \Pb : (\forall \nu < \parts(c))(\forall \alpha < 2)
	[c \Vvdash^\alpha_\nu \varphi(G) \mbox{ or } c \Vvdash^\alpha_\nu \neg \varphi(G)] \}
$$
\end{lemma}

We now define the weak forcing relation which is almost the same 
as the true one, expect that the set~$F'$ is not required to be a subset of~$R_\alpha$
in the case of a $\Pi^0_1$ formula.

\begin{definition}[Weak forcing relation]\label{def:cohzp-forcing-precondition}
Fix a condition $c = (\vec{F}, T, \Ccal)$, a $\Sigma^0_0$ formula $\varphi(G, x)$,
a part $\nu < \parts(T)$ and a side $\alpha < 2$.
\begin{itemize}
	\item[1.] $c \Vdash^\alpha_\nu (\exists x)\varphi(G, x)$ iff there exists a $w \in \omega$ such that $\varphi(F^\alpha_\nu, w)$ holds.
	\item[2.] $c \Vdash^\alpha_\nu (\forall x)\varphi(G, x)$ iff 
	for every $\sigma \in T$, every $w < |\sigma|$ and every set $F' \subseteq \dom(T) \cap \seto_\nu(\sigma)$,
	$\varphi(F^\alpha_\nu \cup F', w)$ holds.
\end{itemize}
\end{definition}

As one may expect, the weak forcing relation at the first level
is also closed under the refinement relation.

\begin{lemma}\label{lem:cohzp-forcing-extension-level1}
Fix a condition $c = (\vec{F}, T, \Ccal)$ and a $\Sigma^0_1$ ($\Pi^0_1$) formula $\varphi(G)$.
If $c \Vdash^\alpha_\nu \varphi(G)$ for some $\nu < \parts(T)$ and $\alpha < 2$,
then for every $d = (\vec{E}, S, \Dcal) \leq c$ and
every part $\mu$ of $S$ refining part $\nu$ of $T$, $d \Vdash^\alpha_\mu \varphi(G)$.
\end{lemma}

The following trivial lemma simply reflects the fact that the promise $\Ccal$ is not part 
of the definition of the weak forcing relation
for $\Sigma^0_1$ or~$\Pi^0_1$ formulas, and therefore has no effect on it.

\begin{lemma}\label{lem:cohzp-promise-no-effect-first-level}
Fix two conditions $c = (\vec{F}, T, \Ccal)$ and $d = (\vec{E}, T, \Dcal)$
and a $\Sigma^0_1$ ($\Pi^0_1$) formula. For every part $\nu$ of $T$ such that $F^\alpha_\nu = E^\alpha_\nu$,
$c \Vdash^\alpha_\nu \varphi(G)$ if and only if $d \Vdash^\alpha_\nu \varphi(G)$.
\end{lemma}

We can now define the forcing relation over higher formulas.
It is defined inductively, starting with $\Sigma^0_1$ and~$\Pi^0_1$ formulas.
We extend the weak forcing relation instead of the true one for effectiveness purposes.
We shall see later that the weak forcing relation behaves like the true one
for some parts and some sides of a condition, and therefore that 
it tells us something about the truth of the formula over some carefully defined generic real~$G$.
Note that the forcing relation over higher formulas is still parameterized by the side~$\alpha$ of
the condition.

\begin{definition}\label{def:cohzp-forcing-condition}
Fix a condition $c = (\vec{F}, T, \Ccal)$, a side $\alpha < 2$ and an arithmetic formula $\varphi(G)$.
\begin{itemize}
	\item[1.] If $\varphi(G) = (\exists x)\psi(G, x)$ where $\psi \in \Pi^0_1$ then
	$c \Vdash^\alpha \varphi(G)$ iff for every part $\nu$ of $T$ such that
	$(\nu, T) \in \Ccal$ there exists a $w < \dom(T)$ such that $c \Vdash^\alpha_\nu \psi(G, w)$
	\item[2.] If $\varphi(G) = (\forall x)\psi(G, x)$ where $\psi \in \Sigma^0_1$ then
	$c \Vdash^\alpha \varphi(G)$ iff for every infinite p.r.\ $k'$-partition tree $S$,
	every function $f : \parts(S) \to \parts(T)$,
	every $w$ and $\vec{E}$ smaller than $\#S$ such that the followings hold
	\begin{itemize}
		\item[i)] $E^\beta_\nu = F^\beta_{f(\nu)}$ for each $\nu < \parts(S)$ and $\beta \neq \alpha$
		\item[ii)] $(E^\alpha_\nu, \dom(S) \cap R_\alpha)$ Mathias extends 
			$(F^\alpha_{f(\nu)}, \dom(T) \cap R_\alpha)$ for each $\nu < \parts(S)$
		\item[iii)] $S$ $f$-refines $\bigcap_{\nu < \parts(S)} T^{[f(\nu), E^\alpha_\nu]}$
	\end{itemize}
	for every $(\mu, S) \in \Ccal$, $(\vec{E}, S, \Ccal[S]) \not \Vdash^\alpha_\mu \neg \psi(G, w)$
	\item[3.] If $\varphi(G) = (\exists x)\psi(G, x)$ where $\psi \in \Pi^0_{n+2}$ then
	$c \Vdash^\alpha \varphi(G)$ iff there exists a $w \in \omega$ such that $c \Vdash^\alpha \psi(G, w)$
	\item[4.] If $\varphi(G) = \neg \psi(G)$ where $\psi \in \Sigma^0_{n+3}$
	then $c \Vdash^\alpha \varphi(G)$ iff $d \not \Vdash^\alpha \psi(G)$ for every $d \in \Ext(c)$.
\end{itemize}
\end{definition}

Note that clause 2.ii) of Definition~\ref{def:cohzp-forcing-condition} seems 
to be~$\Pi^0_2$ since~$R_\alpha$ is $\Delta^0_2$. However, in fact, one just needs to ensure
that~$\dom(S) \subseteq \dom(T)$ and~$E^\alpha_\nu \setminus F^\alpha_{f(\nu)} \subseteq \dom(T) \cap R_\alpha$.
This is a $\Delta^0_2$ predicate, and so is its negation, so one can already easily check that
the forcing relation over a $\Pi^0_2$ formula will be also $\Pi^0_2$.
Before proving the usual properties about the forcing relation, we need to discuss the role of the sides
in the forcing relation.
We are now ready to prove that the forcing relation is closed under extension.

\begin{lemma}\label{lem:cohzp-forcing-extension}
Fix a condition $c$, a side $\alpha < 2$ and a $\Sigma^0_{n+2}$ ($\Pi^0_{n+2}$) formula $\varphi(G)$.
If $c \Vdash^\alpha \varphi(G)$ then for every $d \leq c$, $d \Vdash^\alpha \varphi(G)$.
\end{lemma}

Although the weak forcing relation does not satisfy the density property,
the forcing relation over higher formulas does. The reason is that
the extended forcing relation does not involve the weak forcing relation
over~$\Sigma^0_1$ formulas in the clause 2 of Definition~\ref{def:cohzp-forcing-condition}, 
but uses instead the weaker statement ``$c$ does not force the negation of the~$\Sigma^0_1$ formula''.
The link between this statement and the statement ``$c$ has an extension which forces the $\Sigma^0_1$ formula''
is used when proving that $\varphi(G)$ holds iff $c \Vdash \varphi(G)$ for some condition
belonging to a sufficiently generic filter. We now prove the density of the forcing relation for higher formulas.

\begin{lemma}\label{lem:cohzp-forcing-dense}
For every $\Sigma^0_{n+2}$ ($\Pi^0_{n+2}$) formula $\varphi$, 
the following set is dense in $\Pb$:
$$
\{c \in \Pb : (\forall \alpha < 2)[c \Vdash^\alpha \varphi(G) \mbox{ or } c \Vdash^\alpha \neg \varphi(G)] \}
$$
\end{lemma}

We now prove that the weak forcing relation extended to any arithmetic formula
enjoys the desired definability properties. For this, we start with a lemma
showing that the extension relation is $\Pi^0_2$. Therefore, only
the first two levels have to be treated independently, since 
the extension relation does not add some extra complexity to the forcing relation for higher formulas.

\begin{lemma}\label{lem:cohzp-extension-complexity}
For every condition $c$, $\Ext(c)$ is $\Pi^0_2$ uniformly in~$c$.
\end{lemma}

\begin{lemma}\label{lem:cohzp-complexity-forcing}
Fix an arithmetic formula $\varphi(G)$, a condition $c = (\vec{F}, T, \Ccal)$,
a side $\alpha < 2$ and a part $\nu$ of $T$.
\begin{itemize}
	\item[a)] If $\varphi(G)$ is a $\Sigma^0_1$ ($\Pi^0_1$) formula 
	then so is the predicate $c \Vdash^\alpha_\nu \varphi(G)$.
	\item[b)] If $\varphi(G)$ is a $\Sigma^0_{n+2}$ ($\Pi^0_{n+2}$) formula 
	then so is the predicate $c \Vdash^\alpha \varphi(G)$.
\end{itemize}
\end{lemma}

\subsection{Validity}

As we already saw, we have two candidate forcing relations for $\Sigma^0_1$ and $\Pi^0_1$ formulas:
\begin{itemize}
	\item[1.] The ``true'' forcing relation $c \Vvdash^\alpha \varphi(G)$.
	This relation has been shown to have the expected density properties through Lemma~\ref{lem:cohzp-true-forcing-dense-level1}.
	However deciding such a relation requires too much computational power.
	\item[2.] The ``weak'' forcing relation $c \Vdash^\alpha \varphi(G)$.
	Deciding such a relation requires the same definitional power as the formula it forces.
	It provides a sufficient condition for forcing the formula $\varphi(G)$
	as $c \Vdash^\alpha \varphi(G)$ implies $c \Vvdash^\alpha \varphi(G)$,
	but the converse does not hold and we cannot prove the density property in the general case.
\end{itemize}

Thankfully, there exist some sides and parts of any condition
on which those two forcing relations coincide. This leads to the notion of validity.

\begin{definition}[Validity]
Fix an enumeration $\varphi_0(G), \varphi_1(G), \dots$ of all $\Pi^0_1$ formulas.
Fix a condition $c = (\vec{F}, T, \Ccal)$, a side $\alpha < 2$, and a part $\nu$ of $T$.
We say that \emph{side $\alpha$ is $n$-valid in part $\nu$ of $T$} for some~$n \in \omega$
if part $\nu$ of $T$ is $R_\alpha$-acceptable and for every~$i < n$,
$c \Vvdash^\alpha_\nu \varphi_i(G)$ iff $c \Vdash^\alpha_\nu \varphi_i(G)$.
\end{definition}

The following lemma shows that given some~$n \in \omega$,
we can restrict $\Ccal$ so that it ``witnesses its $n$-valid parts''.

\begin{lemma}\label{lem:cohzp-validity-exists}
For every $n \in \omega$, the following set is dense in~$\Pb$:
$$
\{ (\vec{F}, T, \Ccal) \in \Pb : (\forall \nu)(\exists \alpha < 2)
	[(\nu, T) \in \Ccal \imp \mbox{side } \alpha \mbox{ is } n\mbox{-valid in part } \nu \mbox{ of } T] \} 
$$
\end{lemma}

Given any filter~$\Fcal = \{c_0, c_1, \dots \}$ with $c_s = (\vec{F}_s, T_s, \Ccal_s)$
the set of pairs $(\alpha, \nu_s)$ such that $(\nu_s, T_s) \in \Ccal_s$ forms again
an infinite, directed acyclic graph~$\Gcal(\Fcal)$. By Lemma~\ref{lem:cohzp-validity-exists},
whenever~$\Fcal$ is sufficiently generic,
the graph~$\Gcal(\Fcal)$ yields a sequence of parts $P$ such that for every~$s$ if~$c_s$ refines $c_t$, then part~$P(s)$ of~$c_s$
refines part~$P(t)$ of~$c_t$, and such that for every $n$, there is some~$s$ and some side~$\alpha < 2$ such that the side~$\alpha$
is $n$-valid in part $P(s)$ of~$c_s$. The path~$P$ induces an infinite set~$G = \bigcup \{ F^\alpha_{P(s), s} : s \in \omega \}$. 
Since whenever $\alpha$ is $n$-valid in part~$P(s)$ of~$c_s$, then it is $m$-valid in part~$P(s)$ of~$c_s$ for every~$m < n$,
we can fix an $\alpha < 2$ such that for every $n$, there is some~$s$ such that the side~$\alpha$
is $n$-valid in part $P(s)$ of~$c_s$.
We call $\alpha$ the \emph{generic side}, $P$ the \emph{generic path} and $G$ the \emph{generic real}.

By choosing a generic path that goes through valid sides and parts of the conditions, 
we recovered the density property for the weak forcing relation
and can therefore prove that a property holds over the generic real
if and only if it can be forced by some condition belonging to the generic filter.

\begin{lemma}\label{lem:cohzp-generic-level1}
Suppose that $\Fcal$ is sufficiently generic and let $\alpha$,~$P$ and~$G$ be the generic side, the generic path 
and the generic real, respectively.
For every $\Sigma^0_1$ ($\Pi^0_1$) formula $\varphi(G)$,
$\varphi(G)$ holds iff $c_s \Vdash^\alpha_{P(s)} \varphi(G)$ for some $c_s \in \Fcal$.
\end{lemma}

\begin{lemma}
Suppose that $\Fcal$ is sufficiently generic and let $\alpha$ and~$G$ be 
the generic side and the generic real, respectively.
For every $\Sigma^0_{n+2}$ ($\Pi^0_{n+2}$) formula $\varphi(G)$,
$\varphi(G)$ holds iff $c_s \Vdash^\alpha \varphi(G)$ for some $c_s \in \Fcal$.
\end{lemma}

\subsection{Preserving definitions}

The following (and last) lemma shows that every sufficiently generic real
preserves higher definitions. This preservation property cannot be proved
in the case of non-$\Sigma^0_1$ sets since the weak forcing relation
does not have the good density property in general.

\begin{lemma}\label{lem:cohzp-diagonalization}
If $A \not \in \Sigma^0_{n+2}$ and $\varphi(G, x)$ is $\Sigma^0_{n+2}$,
then the set of $c \in \Pb$ satisfying the following property is dense:
$$
(\forall \alpha < 2)[(\exists w \in A)c \Vdash^\alpha \neg \varphi(G, w)] 
\vee [(\exists w \not \in A)c \Vdash^\alpha \varphi(G, w)]
$$
\end{lemma}

We are now ready to reprove Corollary~3.29 from Wang~\cite{Wang2014Definability}.

\begin{theorem}[Wang~\cite{Wang2014Definability}]
$\rt^2_2$ admits preservation of $\Xi$ definitions 
simultaneously for all $\Xi$ in $\{\Sigma^0_{n+2}, \allowbreak \Pi^0_{n+2}, \allowbreak \Delta^0_{n+3} : n \in \omega \}$.
\end{theorem}
\begin{proof}
Since $\rca \vdash \coh \wedge \mathsf{D}^2_2 \imp \rt^2_2$,
and $\coh$ admits preservation of the arithmetic hierarchy, it suffices to prove
that~$\mathsf{D}^2_2$ admits preservation of $\Xi$ definitions 
simultaneously for all $\Xi$ in $\{\Sigma^0_{n+2}, \allowbreak \Pi^0_{n+2}, \allowbreak \Delta^0_{n+3} : n \in \omega \}$.
Fix some set~$C$ and a $\Delta^{0,C}_2$ 2-partition $R_0 \cup R_1 = \omega$.
Let $\Ccal_0$ be the $C'$-p.r.\ set of all $(\nu, T) \in \TPb$ such that $(\nu, T) \leq (0, 1^{<\omega})$.
Let~$\Fcal$ be a sufficiently generic filter containing~$c_0 = (\{\emptyset, \emptyset\}, 1^{<\omega}, \Ccal_0)$.
Let $G$ be the corresponding generic real.
By definition of a condition, the set~$G$ is $\vec{R}$-cohesive.
By Lemma~\ref{lem:cohzp-diagonalization} and Lemma~\ref{lem:cohzp-complexity-forcing},
$G$ preserves non-$\Sigma^0_{n+2}$ definitions relative to~$C$ for every~$n \in \omega$.
Therefore, by Proposition 2.2 of~\cite{Wang2014Definability}, $G$
preserves $\Xi$ definitions relative to~$C$
simultaneously for all $\Xi$ in $\{\Sigma^0_{n+2}, \Pi^0_{n+2}, \Delta^0_{n+3} : n \in \omega \}$.
\end{proof}

\part{Conclusion}

\chapter{Conclusion}

The early study of reverse mathematics has seen the emergence 
of five main subsystems of second-order arithmetic capturing 
the vast majority of mathematical reasoning. These subsystems 
are known as the Big Five~\cite{Montalban2011Open}. Since then, 
the reverse mathematics community tried to understand this phenomenon, 
by searching for its philosophical and mathematical justification on one hand, 
and by studying theorems escaping the Big Five phenomenon on the other hand. 
See Mont\'alban~\cite{Montalban2011Open} for a discussion about the Big Five.

This thesis falls into the second approach, by studying 
the reverse mathematics of Ramsey-type theorems. 
Indeed, Ramsey's theory provides a large class of theorems escaping the Big Five phenomenon.
These theorems are characterized by their lack of \emph{robustness}, i.e., slight variations
of the statements yield different subsystems.
In this thesis, we tried to give a better grasp on the chaotic nature of the Ramsey-type zoo. 
For this, we compared various Ramsey-type theorems using a uniform, minimalistic and arguably natural framework,
namely, Mathias forcing and preservation of hyperimmunity. Its success in providing a precise analysis
of the computational differences between various Ramsey-type statements and its minimalistic nature makes us believe 
that the current framework has reached a good degree of maturation. It follows that whenever this framework
fails to separate two Ramsey-type statements, they are likely to be equivalent. In particular we take it as an argument
in favor of the equivalence between the thin set theorem for pairs and the free set theorem for pairs.

The study of various Ramsey-type theorems yielded the following empirical observation:
few Ramsey-type statements admit a nice computability-theoretic characterization.
The most notable exceptions are cohesiveness~\cite{Jockusch1993cohesive}, 
the Ramsey-type weak weak K\"onig's lemma~\cite{Bienvenu2015logical}
and the rainbow Ramsey theorem for pairs~\cite{MillerAssorted}. These statements admit characterizations
in terms of degrees whose jump is PA over~$\emptyset'$, of DNC degrees, and of DNC degrees over~$\emptyset'$, respectively.
One could argue that the true computability-theoretic nature of the other Ramsey-type statements
has simply not been found yet. Interestingly, all three exceptions admit a universal instance, whereas we proved that
most of the statements coming from the Ramsey-type zoo do not 
(see Chapter~\ref{chap-degrees-bounding-principles}). This is perhaps evidence that 
the remaining Ramsey-type theorems do not admit nice computability-theoretic characterizations.

Last, we pursued the investigations of Liu~\cite{Liu2012RT22,Liu2015Cone} and Flood~\cite{Flood2012Reverse} on
the amount of compactness needed in the proof of Ramsey-type theorems 
by studying variations of the Ramsey-type weak K\"onig's lemma
on one hand, and by extending Liu's notion of c.b-enum avoidance 
to various Ramsey-type theorems such as the free set, thin set
and the rainbow Ramsey theorems on the other hand.

In the remainder of this conclusion, we develop some of the above-mentioned claims, namely,
the current framework used to analyse Ramsey-type theorems is natural
and the computational content of the Ramsey-type statements mainly comes from the sparsity of their solutions. 

\subsection{The naturality of the preservation framework}

The notion of naturality is part of the common language, but it can have many interpretations.
In this section, we adopt the following (biased) definition: a solution is \emph{natural}
if sufficiently smart people working independently on the same problem would come up with the same solution.
Under this definition, an informal proof of naturality consists in giving a simple reasoning with its justification
leading to the design of the solution. In particular, if the reasoning is generated by a short number
of elementary rules, then a large majority of people could arguably eventually come up with the resulting solution.

Although each Ramsey-type theorem comes up with its own notion of forcing reflecting the combinatorics intrinsic to the theorem, 
the framework remains minimalistic, in that these notions of forcing are designed by following a short number of general rules.
Let us take the example of the notion of forcing for Ramsey's theorem for pairs.
Given a coloring $f : [\omega]^2 \to 2$, a naive notion of forcing to construct a solution
to~$f$ would be to use the finite approximation method, that is, use increasing finite approximations
$F_0 \subsetneq F_1 \subsetneq \dots$ of the solution, and let~$H$ be their union.
We would like to prove that each finite approximation is extendible, so that any sufficiently generic filter
yields an infinite solution.

Here, we meet a first issue: which color do we want to make the set~$H$ homogeneous for?
A coloring of pairs does not have in general an infinite set homogeneous for each color.
At this stage, we apply the rule ``If you cannot choose, try all''.
The notion of forcing now becomes over pairs of finite approximations $(F, E)$
where $F$ is homogeneous for color~0 and $E$ for color~1.
We would like to prove that for each such pair~$(F, E)$, either~$F$ or~$E$ is extendible to an infinite homogeneous set.
This is not possible, as witnessed by taking the pair~$(\{x\}, \{y\})$ such that
 $(\forall^{\infty} s)f(x,s) = 1$ and $(\forall^{\infty} s)f(y,s) = 1$.
We need therefore to refine again our notion of forcing.
Without loss of generality, we can consider that a finite approximation~$F$ is extendible 
if it can be ``plugged'' to an infinite set~$H > F$ so that~$F \cup H$ is homogeneous.
If~$F$ and~$H$ are both homogeneous for color~0, then the minimal condition to ensure that~$F \cup H$ is homogeneous for color~0
is $(\forall x \in F)(\forall y \in H)f(x,y) = 0$. The good notion of forcing becomes
$(F, E, X)$, where $F$ is homogeneous for color~0, $E$ is homogeneous for color~1, and $X$ is an infinite set
such that the following holds:
\begin{itemize}
	\item[(i)] for every~$x \in F$ and~$y \in X$, $f(x, y) = 0$
	\item[(ii)] for every~$x \in E$ and~$y \in X$, $f(x, y) = 1$
\end{itemize}
This way, by applying Ramsey's theorem for pairs over the domain~$X$, we obtain an infinite $f$-homogeneous set~$H \subseteq X$,
say for color~0. By (i), $F \cup H$ is $f$-homogeneous for color~0 and therefore~$F$ is extendible.
All the steps leading to the design of this notion of forcing involve simple and natural reasonings.
The condition of extendibility is arguably minimalistic.

Let us now justify the naturality of the preservation of hyperimmunity
by giving some details about the design of the notion.
Suppose we want to separate any statement~$\Psf$ from stable Ramsey's theorem for pairs.
Consider any computable $\srt^2_2$-instance with no computable solution.
Such an instance is nothing, but a stable, computable coloring~$f : [\omega]^2 \to 2$
such that for each~$i < 2$, the set $A_i = \{x : \lim_s f(x, s) \neq i \}$ is immune.
Indeed, any infinite computable set~$H$ will intersect both~$A_0$ and~$A_1$ and will not be $f$-homogeneous.
If we are working over computable reducibility, we need to prove that any computable $\Psf$-instance has a solution~$S$
such that $f$ has no $S$-computable solution, in other words, such that~$A_0$ and~$A_1$ are $S$-hyperimmune.
Over computable entailment, we need to be able to iterate this property, and prove that whenever the~$A$'s are $C$-immune,
then any $C$-computable $\Psf$-instance has a solution~$S$ such that the $A$'s are $S \oplus C$-immune.
This is nothing but the preservation of immunity. One can easily prove that cohesiveness
admits preservation of immunity. 

Some other statements such as the Erd\H{o}s-Moser theorem can prove some weak form of compactness arguments,
namely, the Ramsey-type weak K\"onig's lemma. Because of this, the notion of immunity has to be strengthened
to defeat not only integers, but blocks of integers. This is the notion of hyperimmunity.
In this perspective, preservation of hyperimmunity is the natural notion obtained while trying
to separate a statement containing some weak form of compactness from stable Ramey's theorem for pairs
over computable entailment.

\subsection{The computational content of Ramsey-type theorems}

After an extensive study of various Ramsey-type theorems, it appears
that the reductions from a statement to a Ramsey-type theorem can be grouped into two families.

In the first case, the reduction exploits the ability of the Ramsey-type statement to impose
some strong sparsity constraints. For example, the proof that Ramsey's theorem for triples
implies the arithmetic comprehension axiom uses a coloring whose homogeneous sets
are sparse enough to compute a fast-growing function~\cite{Jockusch1972Ramseys}. Similarly, the proof that 
the thin set theorem for pairs implies the atomic model theorem exploits 
the ability of the thin set theorem for pairs to impose large holes in its solutions,
so that its principal function is not dominated by a $\Delta^0_2$ function (Theorem~\ref{thm:sts2-amt}).

In the second case, the reduction has a combinatorial nature and involves two Ramsey-type statements.
This is then a strong reduction with no computability-theoretic considerations. 
For example, the strong computable equivalence between the free set theorem and the rainbow Ramsey theorem
consists of encoding an instance of one statement into an instance of the other statement,
so that any solution to the latter \emph{is} a solution to the former
(see \cite{Wang2014Some} and Theorem~\ref{thm:rrt2n-fsn}).

It follows that Ramsey-type theorems, and more generally combinatorial theorems
do not seem to carry a strong computational content other than some simple sparsity considerations.
This consideration can be seen as a partial answer to the chaotic behavior of Ramsey-type statements with respect
to the reverse mathematics and computable reducibility.

\subsection{Further developments}

Various notions of preservation have been used among the past few years to separation
statements in reverse mathematics. Seetapun~\cite{Seetapun1995strength} used cone avoidance to separate Ramsey's theorem for pairs
from the arithmetic comprehension axiom, Hirschfeldt and Shore~\cite{Hirschfeldt2007Combinatorial}
used lowness to separate stable chain antichain from the ascending descending sequence principle,
Wang~\cite{Wang2014Definability} used the preservation of non-c.e.\ definitions to separate the Erd\H{o}s-Moser theorem
from the stable thin set theorem for pairs or stable Ramsey's theorem for pairs, Liu~\cite{Liu2015Cone} used c.b-enum avoidance
to separate Ramey's theorem for pairs from weak weak K\"onig's lemma, and the the author 
used preservation of $k$ hyperimmunities to separate the thin set theorem from Ramsey's theorem for pairs, among others.

Those notions of preservation fall into two categories: the effectiveness properties and the genericity properties.
The preservation of an effectiveness property has to be done by an effective construction. For example, lowness is an effectiveness
property, since the whole construction has to be $\Delta^0_2$.
On the other hand, the preservation of a genericity property is proven by designing a notion of forcing such that any sufficiently
generic filter yields a set preserving this property. This is for example the case of cone avoidance and preservation of hyperimmunity.
The genericity properties happen to be easier to prove since they do not involve effectiveness considerations when proving that a set is dense.

The genericity properties often reveal some combinatorial and effective features of the corresponding
notions of forcing. The genericity properties are designed according to characteristics of the notion
of forcing rather than the converse. We detail some aspects taken in account in the design of a genericity property:
\begin{itemize}
	\item[(a)] How many objects are built simultaneously? The notion of forcing for the Erd\H{o}s-Moser theorem
	yields only one set, while the one for Ramsey's theorem for pairs and~$k$ colors yields $k$ sets.
	Whenever multiple sets are built together, the requirements have a disjunctive shape, which prevents
	one from interleaving multiple requirements. This aspect is used in a essential manner
	when separating the thin set for pairs and $k+1$ colors from the thin set for pairs and~$k$ colors.
	\item[(b)] Whenever the decision of a $\Sigma^0_1(G)$ formula is positive,
	how many candidate extensions are there? In computable Mathias forcing, a positive answer to a $\Sigma^0_1$
	question yields a single candidate extension, whereas for weak K\"onig's lemma, this yields a finite, but arbitrarily large
	number of extensions. This aspect is used in the separation of the chain-antichain principle
	from Ramey's theorem for pairs.
	\item[(c)] Whenever the decision of a $\Sigma^0_1(G)$ formula is negative,
	does the current condition already force its negation?
	Again, in computable Mathias forcing, whenever the answer is negative, the current condition already forces the negation
	of the $\Sigma^0_1(G)$ statement, whereas for Ramsey's theorem for pairs, a negative answer simply shows the existence
	of an extension forcing the negation. This aspect is used in the separation of the ascending descending sequence
	from the stable chain antichain principle.
\end{itemize}

One may naturally wonder whether we can establish a formal classification of the notions of forcing
so that the properties of preservation are directly proven over the abstract partial orders.
Then, one would independently prove that the statements studied in reverse mathematics
admit notions of forcing belonging to this classification, and directly deduce
the corresponding separations. Let us take an example.

Consider a notion of forcing~$(P, \leq_P)$ with a computable domain $P = \{c_0, c_1, \dots \}$ together 
with a computable function $f : P \times \Sigma^{0,G}_1 \to \Sigma^0_1$ such that if $c \in P$
and $\varphi(G)$ is some~$\Sigma^{0,G}_1$ formula, $f(c, \varphi)$ is a $\Sigma^0_1$ formula 
which holds iff there is an extension~$d \leq_P c$ forcing~$\varphi(G)$.
We can show that any sufficiently generic filter for this notion of forcing yields a set preserving 
countably many immunities. For this, fix any immune set~$A$ and any Turing index~$e$. We will show that the set 
$D^A_e \subseteq P$ of conditions forcing~$W^G_e$ not to be an infinite subset of~$A$ is $P$-dense.
For any~$n \in \omega$, let~$\varphi_n(G)$ be the statement ``$n \in W_e^G$'',
and consider the set~$B = \{ n : f(\varphi_n) \mbox{ holds} \}$. We have two cases.
In the first case, $B$ is finite. In this case, $c$ forces $n \not \in W_e^G$ for cofinitely many $n$'s
and therefore forces~$W_e^G$ not to be infinite.
In the second case, $B$ is infinite. Since~$f(\varphi_n)$ is $\Sigma^0_1$ uniformly in~$n$,
the set~$B$ is c.e., therefore there is some~$n \in B \setminus A$.
Since~$n \in B$, $f(\varphi_n)$ holds, so there is an extension forcing~$n \in W_e^G$, and therefore
forcing~$W_e^G$ not to be an infinite subset of~$A$.

One easily sees that the computable Mathias forcing satisfies the previous properties
and therefore cohesiveness admits preservation of immunity. On the other hand,
the standard notion of forcing for weak K\"onig's lemma does not satisfy the previous properties.
It would be interesting to extend this observation to various notions of forcing
and various preservation properties.

\chapter{Summary}

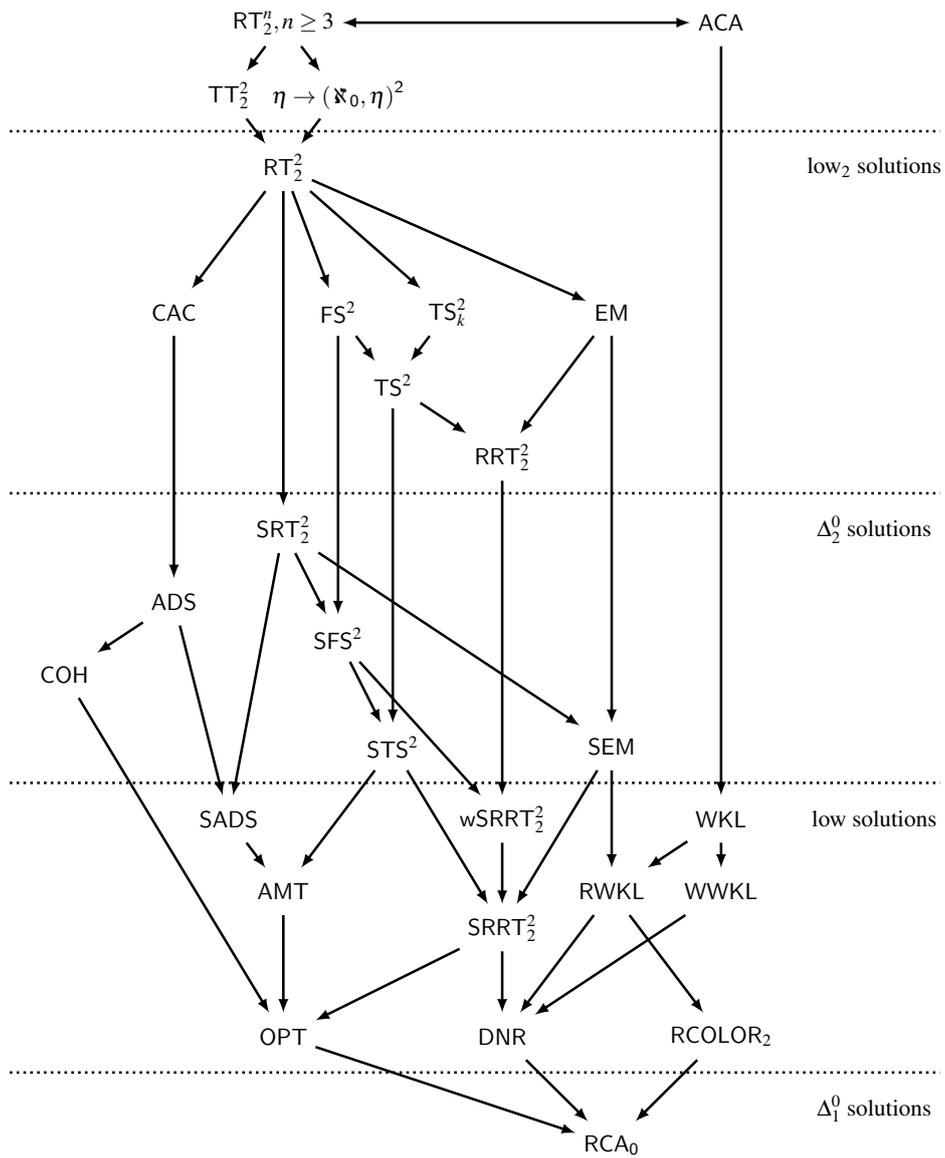
\begin{figure}[htbp]
\begin{center}
\scalebox{.8}{
\begin{tikzpicture}[x=1.8cm, y=1.2cm,
	node/.style={minimum size=2em},
	impl/.style={draw,very thick,-latex},
	equiv/.style={draw, very thick, latex-latex},
	caption/.style={}]

	\node[node] (rt22) at (1,12.5) {$\rt^2_2$};
	\node[node] (tt22) at (0.5, 13.5) {$\tto^2_2$};
	\node[node] (er22) at (1.5, 13.5) {$\erp$};
	\node[node] (rt32) at (1, 14.5)  {$\rt^n_2, n \geq 3$};
	\node[node] (cac) at (0, 10.5) {$\cac$};
	\node[node] (ads) at (0, 6.5) {$\ads$};
	\node[node] (coh) at (-1, 5.5) {$\coh$};
	\node[node] (wkl) at (5, 3.5) {$\wkl$};
	\node[node] (rwkl) at (4, 2.5) {$\rwkl$};
	\node[node] (wwkl) at (5, 2.5) {$\wwkl$};
	\node[node] (dnr) at (3, 0.5) {$\dnr$};
	\node[node] (rcolor2) at (5, 0.5) {$\rcolor_2$};
	\node[node] (aca) at (5, 14.5) {$\aca$};
	\node[node] (rca) at (4, -1) {$\rca$};
	\node[node] (sads) at (0.5, 3.5) {$\sads$};
	\node[node] (amt) at (1,2.5) {$\amt$};
	\node[node] (opt) at (1,0.5) {$\opt$};
	\node[node] (em) at (4, 10.5) {$\emo$};
	\node[node] (rrt22) at (3, 8.5) {$\rrt^2_2$};
	\node[node] (srt22) at (1, 7.5) {$\srt^2_2$};
	\node[node] (sem) at (4, 4.5) {$\semo$};
	\node[node] (wsrrt22) at (3, 3.5) {$\wsrrt^2_2$};
	\node[node] (srrt22) at (3, 2) {$\srrt^2_2$};
	\node[node] (sts2) at (2, 4.5) {$\sts^2$};
	\node[node] (ts2) at (2, 9.5) {$\ts^2$};
	\node[node] (ts2k) at (2.5, 10.5) {$\ts^2_k$};
	\node[node] (fs2) at (1.5, 10.5) {$\fs^2$};
	\node[node] (sfs2) at (1.5, 6) {$\sfs^2$};
	
	\draw[equiv] (aca) -- (rt32);
	\draw[impl] (aca) -- (wkl);
	\draw[impl] (wkl) -- (wwkl);
	\draw[impl] (wkl) -- (rwkl);
	\draw[impl] (wwkl) -- (dnr);
	\draw[impl] (rwkl) -- (dnr);
	\draw[impl] (rwkl) -- (rcolor2);
	\draw[impl] (rcolor2) -- (rca);
	\draw[impl] (dnr) -- (rca);	
	\draw[impl] (rt32) -- (tt22);
	\draw[impl] (rt32) -- (er22);
	\draw[impl] (tt22) -- (rt22);
	\draw[impl] (er22) -- (rt22);
	\draw[impl] (rt22) -- (cac);
	\draw[impl] (rt22) -- (em);
	\draw[impl] (cac) -- (ads);
	\draw[impl] (ads) -- (coh);
	\draw[impl] (ads) -- (sads);
	\draw[impl] (sads) -- (amt);
	\draw[impl] (opt) -- (rca);
	\draw[impl] (amt) -- (opt);
	\draw[impl] (coh) -- (opt);
	\draw[impl] (em) -- (rrt22);
	\draw[impl] (rrt22) -- (wsrrt22);
	\draw[impl] (wsrrt22) -- (srrt22);
	\draw[impl] (srrt22) -- (dnr);
	\draw[impl] (rt22) -- (srt22);
	\draw[impl] (srt22) -- (sads);
	\draw[impl] (sem) -- (srrt22);
	\draw[impl] (srt22) -- (sem);
	\draw[impl] (srrt22) -- (opt);
	\draw[impl] (em) -- (sem);
	\draw[impl] (sem) -- (rwkl);
	\draw[impl] (srt22) -- (sfs2);
	\draw[impl] (sfs2) -- (sts2);
	\draw[impl] (sts2) -- (srrt22);
	\draw[impl] (rt22) -- (ts2k);
	\draw[impl] (rt22) -- (fs2);
	\draw[impl] (ts2k) -- (ts2);
	\draw[impl] (fs2) -- (ts2);
	\draw[impl] (ts2) -- (rrt22);
	\draw[impl] (ts2) -- (sts2);
	\draw[impl] (sts2) -- (amt);
	\draw[impl] (fs2) -- (sfs2);
	\draw[impl] (sfs2) -- (wsrrt22);

	\draw[very thick, dotted] (-1.5, 0) -- (7,0);
	\node[caption] at (6.4,-0.5) {$\Delta^0_1$ solutions};
	
	\draw[very thick, dotted] (-1.5,4) -- (7,4);
	\node[caption] at (6.4,3.5) {low solutions};
	
	\draw[very thick, dotted] (-1.5,8) -- (7,8);
	\node[caption] at (6.4,7.5) {$\Delta^0_2$ solutions};

	\draw[very thick, dotted] (-1.5,13) -- (7,13);
	\node[caption] at (6.4,12.5) {low${}_2$ solutions};
	
\end{tikzpicture}}
\end{center}
\caption{Ramsey-type theorems under computable entailment}
\label{fig:intro-rm-full-ramsey-zoo}
\end{figure}

\def\arraystretch{1.5}
\begin{figure}[htbp]
\begin{center}
\begin{tabular}{l|l|l}
\hline
	& Preservations & Non-preservations\\ \hline
All hyperimmunities & $\wkl$ \cite{Jockusch197201} & Cohen generics \cite{Kurtz1982Randomness} \\
	& & Random reals \cite{Downey2010Algorithmic} \\
	& & $\coh$ \cite{Jockusch1993cohesive} \\ 
	& & $\amt$ \cite{Hirschfeldt2009atomic} \\ 
	& & $\rrt^2_2$ \cite{Csima2009strength} \\ \hline
$\omega$ hyperimmunities & Cohen generics \cite{Patey2015Iterative} & $\sts^2$, $\sfs^2$ (\ref{cor:ts2-non-preservation}) \\
	& Random reals \cite{Patey2015Iterative} & \\
	& $\coh$ (\ref{thm:coh-hyperimmunity-preservation}) & \\ 
	& $\emo$ (\ref{thm:em-hyperimmunity}) & \\ 
	& $\rrt^2_2$ \cite{Patey2015Iterative} & \\
	& $\amt$ \cite{Patey2015Iterative} & \\ \hline
$k$ hyperimmunities & $\ts^2_{k+1}$ (\ref{cor:ts2-preservation}) & $\ts^2_k$ (\ref{cor:ts2-non-preservation}) \\ 
	& $\fs$, $\ts$, $\rrt$ (\ref{thm:fs2-preservation-definitions}) & \\ \hline
\end{tabular}
\caption{Preservations and non-preservations of hyperimmunity}
\end{center}
\end{figure}

\begin{figure}[htbp]
\begin{center}
\begin{tabular}{l|l|l}
\hline
	& Avoidances & Non-avoidances\\ \hline
Path avoidance & Cohen generics (\ref{thm:genericity-path-avoidance}) & $\dnr$ (\ref{cor:dnr-path-avoidance}) \\
	& $\coh$ (\ref{cor:coh-path-avoidance}) & \\
	& $\amt$ (\ref{thm:opt-amt-pizog-path-avoidance}) & \\ \hline

1-enum avoidace & $\rrt^2_2$ (\ref{cor:rrt22-enum-avoidance}) & \\
	& $\dnr$ (\ref{cor:dnr-strong-1-enum-avoidance}) &  \\ \hline
$\omega$ c.b-enum avoidance & $\emo$ (\ref{cor:emo-simultaneous-cbenum-avoidance})
			 & $\sts^2$, $\sfs^2$ (\ref{cor:sts2-not-simultaneous-cb-enum-avoidance}) \\
	& $\rwkl$ (\ref{thm:ts-rwkl-pi01-enum-avoidance}) & $\sads$ \\ \hline
$k$ c.b-enum avoidance & $\ts^2_{k+1}$ (\ref{cor:tsm-cb-enum-avoidance}) 
			& $\ts^2_k$ (\ref{cor:sts2-not-simultaneous-cb-enum-avoidance}) \\ 
	& $\fs$, $\ts$, $\rrt$ (\ref{cor:fs-strong-cb-enum-avoidance}) & Random reals \cite{Liu2015Cone} \\ \hline
\end{tabular}
\caption{Avoidances and non-avoidances of closed sets}
\end{center}
\end{figure}

\chapter{Open questions}\label{chap:open-questions}

The strength of consequences of Ramsey's theorem has been extensively studied in reverse mathematics and under various reducibilities, namely, computable reducibility and Weihrauch reducibility. Our understanding of the combinatorics of Ramsey's theorem and its consequences has been greatly improved over the past decades. In this chapter, we state some questions which naturally arose during this study. The inability to answer those questions reveals some gaps in our understanding of the combinatorics of Ramsey's theorem.

Many of the questions are already stated in various papers.
These questions cover a much thinner branch of reverse mathematics than
the paper of Mont\'alban~\cite{Montalban2011Open} and are of course
influenced by the own interest of the author. We put our focus on a few central open questions,
motivate them and try to give some insights about the reason for their complexity. 
We also detail some promising approaches and ask more technical questions, driven by the resolution of the former ones.
The questions are computability-theoretic oriented, and therefore mainly
involve the relations between statements over computable entailment.

\section{Cohesiveness and partitions}

In order to better understand the combinatorics of Ramsey's theorem for pairs,
Cholak, Jockusch and Slaman~\cite{Cholak2001strength} decomposed it into a cohesive and a stable version.
This approach has been fruitfully reused in the analysis of various consequences
of Ramsey's theorem~\cite{Hirschfeldt2007Combinatorial}.
They claimed that~$\rt^2_2$ is equivalent to~$\srt^2_2+\coh$ over~$\rca$
with an erroneous proof. Mileti~\cite{Mileti2004Partition} and Jockusch \& Lempp [unpublished]
independently fixed the proof. 
Cholak et al.~\cite{Cholak2001strength} proved that~$\coh$ does not imply~$\srt^2_2$
over~$\rca$. The question of the other direction was a long-standing open problem.
Recently, Chong et al.~\cite{Chong2014metamathematics} proved that~$\srt^2_2$ is strictly weaker than~$\rt^2_2$
over~$\rca+\bst$. However they used non-standard models to separate the statements and the question
whether~$\srt^2_2$ and~$\rt^2_2$ coincide over~$\omega$-models remains open.
See the survey of Chong, Li and Yang~\cite{Chong2014Nonstandard} for the approach of non-standard analysis 
applied to reverse mathematics.

\begin{question}\label{quest:srt22-rt22-omega}
Does~$\srt^2_2$ imply~$\rt^2_2$ (or equivalently $\coh$) over $\omega$-models?
\end{question}

Jockusch~\cite{Jockusch1972Ramseys} constructed a computable instance of~$\rt^2_2$
with no $\Delta^0_2$ solution. Cholak et al.~\cite{Cholak2001strength} suggested to build
an $\omega$-model of~$\srt^2_2$ composed only of low sets. However, Downey et al.~\cite{Downey200102}
constructed a $\Delta^0_2$ set with no low subset of either it or its complement.
As often, the failure of an approach should not be seen as a dead-end, but as a starting point.
The construction of Downey et al.\ revealed that~$\srt^2_2$ carries some additional computational
power, whose nature is currently unknown. Indeed, all the natural computability-theoretic consequences
of~$\srt^2_2$ known hitherto admit low solutions. Answering the following question 
would be a significant step towards understanding the strength of~$\srt^2_2$.

\begin{question}
Is there a natural computable $\srt^2_2$-instance with no low solution?
\end{question}

Here, by ``natural'', we mean an instance which carries more informational content than
not having low solutions. Interestingly, Chong et al.~\cite{Chong2014metamathematics} constructed a non-standard
model of~$\rca+\bst\allowbreak+\srt^2_2$ containing only low sets. This shows that the argument
of Downey et al.~\cite{Downey200102} requires more than~$\Sigma^0_2$-bounding,
and suggests that such an instance has to use an elaborate construction.

Hirschfeldt et al.~\cite{Hirschfeldt2008strength} proposed a very promising
approach using an extension of Arslanov's completeness criterion~\cite{Arslanov1981Some,Jockusch1989Recursively}.
A set~$X$ is~\emph{1-CEA over~$Y$} if c.e.\ in and above~$Y$.
A set~$X$ is~\emph{$(n+1)$-CEA over~$Y$} if it is 1-CEA over a set~$Z$ which is $n$-CEA over~$Y$.
In particular the sets which are 1-CEA over~$\emptyset$ are the c.e. sets.
By~\cite{Jockusch1989Recursively}, every set $n$-CEA over~$\emptyset$ computing a set of PA degree is complete.
Hirschfeldt et al.\ asked the following question. 

\begin{question}\label{quest:delta2-low2-delta2}
Does every $\Delta^0_2$ set admit an infinite subset of either it or its complement
which is both low${}_2$ and~$\Delta^0_2$? 
\end{question}

A positive answer to Question~\ref{quest:delta2-low2-delta2} would enable one to build an $\omega$-model $\Mcal$ of~$\srt^2_2$
such that for every set~$X \in \Mcal$, $X$ is low${}_2$ and $X'$ is $n$-CEA over~$\emptyset'$. By Jockusch
\& Stephan~\cite{Jockusch1993cohesive}, if some~$X \in \Mcal$ computes some p-cohesive set, then~$X'$ is of PA degree relative to~$\emptyset'$. By a relativization of Jockusch et al.~\cite{Jockusch1989Recursively}, $X'$ would compute~$\emptyset''$, so would be high,
which is impossible since $X$ is low${}_2$. Therefore, $\Mcal$ would be an $\omega$-model of~$\srt^2_2$ which is not a model of~$\coh$,
answering Question~\ref{quest:srt22-rt22-omega}. Note that the argument also works if we replace~low${}_2$ by low${}_n$, where~$n$ may even
depend on the instance. Hirschfeldt et al.~\cite{Hirschfeldt2008strength} 
proved that every $\Delta^0_2$ set has an infinite incomplete $\Delta^0_2$ subset of either it or its complement. This is the best upper bound currently known. They asked the following question which is the strong negation of Question~\ref{quest:delta2-low2-delta2}.

\begin{question}\label{quest:delta2-jump-pa}
Is there a $\Delta^0_2$ set such that every infinite subset of either it or its complement 
has a jump of PA degree relative to~$\emptyset'$? 
\end{question}

By Arslanov's completeness criterion, a positive answer to Question~\ref{quest:delta2-jump-pa}
would also provide one to the following question.

\begin{question}\label{quest:delta2-high}
Is there a $\Delta^0_2$ set such that every $\Delta^0_2$ infinite subset of either it or its complement 
is high? 
\end{question}

Cohesiveness can be seen as a sequential version of~$\rt^1_2$ with finite errors ($\seqs{\rt^1_2}$).
More formally, given some theorem~$\Psf$, $\seqs{\Psf}$ is the statement ``For every uniform
sequence of $\Psf$-instances~$X_0, X_1, \dots$, there is a set~$Y$ which is, up to finite changes,
a solution to each of the~$X$'s.'' One intuition about the guess that $\srt^2_2$ does not imply~$\coh$
could be that by the equivalence between~$\srt^2_2$ and~$\mathsf{D}^2_2$,
$\srt^2_2$ is nothing but a non-effective instance of~$\rt^1_2$, and that one cannot encode in a single instance
of~$\rt^1_2$ countably many $\rt^1_2$-instances. Note that this argument is of combinatorial
nature, as it does not make any effectiveness assumption on the instance of~$\rt^1_2$.
We express reservations concerning the validity of this argument, as witnessed by considering
a sequential version of the thin set theorem.

Recall that for every~$n \geq 1$ and~$k \geq 2$, $\ts^n_k$ is the statement
``Every coloring $f : [\Nb]^n \to k$ has a thin set''
and $\ts^n$ is the statement ``Every coloring $f : [\Nb]^n \to \Nb$ has a thin set''.
According to the definition, $\seqs{\ts^1_k}$ is the statement ``For every uniform sequence of functions
$f_0, f_1, \dots : \Nb \to k$, there is an infinite set~$H$ which is, up to finite changes, thin
for all the~$f$'s.'' Theorem~3.5 in Jockusch \& Stephan~\cite{Jockusch1993cohesive} can easily be adapted to prove that
the degrees of solutions to primitive recursive instances of~$\seqs{\ts^1}$ (resp.\ $\seqs{\ts^1_k}$)
are exactly those whose jump is of DNC (resp.\ $k$-valued DNC) degree relative to~$\emptyset'$.
By a relativization of Friedberg~\cite{JockuschJr1989Degrees}, the degrees whose jump is PA and those whose
jump bounds a $k$-valued DNC function coincide. Therefore~$\coh$ and~$\seqs{\ts^1_k}$ are computably equivalent.
However, $\seqs{\ts^1}$ is a strictly weaker statement, as for any computable instance of~$\seqs{\ts^1}$,
the measure of oracles computing a solution to it is positive.

Recall our intuition that a single instance of~$\rt^1_2$ cannot encode the information
of countably many instances of~$\rt^1_2$. This intuition is false when considering~$\ts^1$.
Indeed, there is a (non-$\Delta^0_2$) instance of~$\ts^1$ (and \emph{a fortiori} one of~$\rt^1_2$) whose solutions all bound a
function DNC relative to~$\emptyset'$, and therefore computes a solution to any computable instance of~$\seqs{\ts^1}$.
Therefore, before asking whether $\coh$ is a consequence of $\emptyset'$-effective instances of~$\rt^1_2$,
it seems natural to ask whether $\coh$ is a consequence of any coloring over singletons in a combinatorial sense,
that is, with no effectiveness restriction at all.

\begin{question}\label{quest:rt12-combi-jump-pa}
Is there any $\rt^1_2$-instance whose solutions have a jump of PA degree relative to~$\emptyset'$?
\end{question}

Note that a negative answer to Question~\ref{quest:rt12-combi-jump-pa} would have practical reverse mathematical
consequences. There is an ongoing search for natural statements strictly between $\rt^2_2$ and~$\srt^2_2$ over~$\rca$~\cite{Dzhafarov2009polarized}.
Dzhafarov \& Hirst~\cite{Dzhafarov2009polarized} introduced the increasing polarized Ramsey's theorem for pairs ($\ipt^2_2$), 
and proved it to be between~$\rt^2_2$ and~$\srt^2_2$ over~$\rca$. The author~\cite{Patey2015Somewhere} proved that
$\ipt^2_2$ implies the existence of a function DNC relative to~$\emptyset'$, therefore showing that~$\srt^2_2$
does not imply~$\ipt^2_2$ over~$\rca$. The statement $\ipt^2_2$ is equivalent
to $\rwkls{2}$, a relativized variant of the Ramsey-type weak K\"onig's lemma, over~$\rca$,
and therefore is a combinatorial consequence of~$\rt^1_2$. 
An iterable negative answer to Question~\ref{quest:rt12-combi-jump-pa} would prove that~$\ipt^2_2$
does not imply~$\rt^2_2$, hence is strictly between $\rt^2_2$ and~$\srt^2_2$ over~$\rca$.

We have seen that $\seqs{\ts^1}$ is a combinatorial consequence of~$\ts^1$.
This information helps us for tackling the following question.
Indeed, if proven false, it must be answered by effective means and not  by combinatorial ones.

\begin{question}\label{quest:srt-seq-thin-set}
Does~$\srt^2_2$ imply~$\seqs{\ts^1}$ over~$\omega$-models?
\end{question}

\begin{question}\label{quest:rt-jump-DNC}
Is there a $\Delta^0_2$ set such that every infinite subset of either it or its complement 
has a jump of DNC degree relative to~$\emptyset'$? 
\end{question}

Although those questions are interesting in their own right, 
a positive answer to Question~\ref{quest:delta2-low2-delta2} would provide a negative answer to Question~\ref{quest:srt-seq-thin-set},
and a positive answer to Question~\ref{quest:rt-jump-DNC} would provide a positive answer to Question~\ref{quest:delta2-high}.
Indeed, the extended version of Arslanov's completeness criterion states that in fact,
every set $n$-CEA over~$\emptyset$ computing a set of DNC degree is complete. In particular,
if Question~\ref{quest:delta2-low2-delta2} has a positive answer, then there is an $\omega$-model of $\srt^2_2$
which is not a model of the rainbow Ramsey theorem for pairs.

Let us finish this section by a discussion about why those problems are so hard to tackle.
There are different levels of answers, starting from the technical one which is more objective, but probably also less informative,
to the meta discussion might give more insights, but can be more controversial.

From a purely technical point of view, all the forcing notions used so far to produce solutions
to Ramsey-type statements are variants of Mathias forcing. In particular, they restrict the future elements
to a ``reservoir''. Any sufficiently generic filters for those notions
of forcing yield cohesive sets. Therefore, one should not expect to obtain a diagonalization against 
instances of~$\coh$ by exhibiting a particular dense set of conditions. Indeed, one would derive a contradiction
by taking a set sufficiently generic to meet both those diagonalizing sets, and the dense sets producing a cohesive solution.
More generally, as long as we use a forcing notion where we restrict the future elements to a reservoir,
any diagonalization against~$\coh$ has to strongly rely on some effectiveness of the overall construction.
The first and second jump control of Cholak et al.~\cite{Cholak2001strength} form a case in point
of how to restrict the amount of genericity to obtain some stronger properties,
which are provably wrong when taking any sufficiently generic filter.

At a higher level, we have mentioned that $\coh$ is a statement
about the jump of Turing degrees. In other words, by Shoenfield's limit lemma~\cite{Shoenfield1959degrees}, 
$\coh$ is a statement about some limit behavior, and is therefore non-sensitive
to any local modification. However, the computability-theoretic properties
used so far to separate statements below~$\aca$ are mainly acting ``below the first jump'',
in the sense that the diagonalization occurs after a finite amount of time.
With $\coh$, there will be some need for a ``continuous diagonalization'', that is, 
a diagonalization which has to be maintained all along the construction.

\section{A Ramsey-type weak K\"onig's lemma}

Weak K\"onig's lemma plays a central role in reverse mathematics.
It is one of the Big Five and informally captures compactness arguments~\cite{Simpson2009Subsystems}.
$\wkl$ is involved in many constructions of solutions to Ramsey-type statements, e.g.,
cone avoidance~\cite{Seetapun1995strength,Wang2014Some} or control of the jump~\cite{Cholak2001strength,Patey2016Controlling}.
The question of whether $\rt^2_2$ implies~$\wkl$ over~$\rca$
was open for decades, until Liu~\cite{Liu2012RT22} solved it by proving that PA degrees
are not a combinatorial consequence of~$\rt^1_2$.
Recently, Flood~\cite{Flood2012Reverse} clarified the relation between Ramsey-type theorems
and $\wkl$, by introducing a Ramsey-type variant of weak K\"onig's lemma ($\rwkl$),
that we studied in chapter~\ref{chap:ramsey-type-konig-lemma}.
Informally, seeing a set as a 2-coloring of the integers, 
for every $\Pi^0_1$ class of 2-colorings, $\rwkl$ states
the existence of an infinite set homogeneous for one of them.
Very little is currently known about~$\rwkl$.
Despite its complicated formulation, $\rwkl$ is a natural statement
which is worth being studied due to the special status of Ramsey-type theorems in reverse mathematics.

The statements analysed in reverse mathematics are collections of problems (instances)
coming with a class of solutions. Sometimes, it happens that one problem is maximally difficult.
In this case, the strength of the whole statement can be understood by studying this particular (universal) instance.
Weak K\"onig's lemma is known to admit universal instances,
e.g., the tree whose paths are completions of Peano arithmetic,
whereas Mileti~\cite{Mileti2004Partition} proved that $\rt^2_2$ and~$\srt^2_2$ do not admit one.
In chapter~\ref{chap-degrees-bounding-principles}, we studied extensively which theorems
in reverse mathematics admit a universal instance,
and which do not. It happens that most consequences of $\rt^2_2$
in reverse mathematics do not admit a universal instance.
The most notable exceptions are the rainbow Ramsey theorem 
for pairs~\cite{Csima2009strength,Patey2015Somewhere,WangSome,Wang2014Cohesive},
the finite intersection property~\cite{Cholak2015Any,Downey2012finite,Dzhafarov2013strength} and~$\dnr$.
It is natural to wonder, given the fact that~$\rwkl$ is a consequence of both~$\srt^2_2$
and of~$\wkl$, whether~$\rwkl$ admits a universal instance.

\begin{question}\label{quest:rwkl-universal-instances}
Does~$\rwkl$ admit a universal instance?
\end{question}

There is a close link between the $\Psf$-bounding degrees and the existence of a universal $\Psf$-instance.
Indeed, the degrees of the solutions to a universal $\Psf$-instance are $\Psf$-bounding.
Using the contraposition, a statement~$\Psf$ is usually proven not to have a universal instance
by showing that every computable $\Psf$-instance has a solution of degree belonging to a class~$\Ccal$,
and that for every degree~$\dbf \in \Ccal$, there is a computable $\Psf$-instance to which~$\dbf$
bounds no solution~\cite{Mileti2004Partition,Patey2015Degrees}.

Interestingly, Question~\ref{quest:rwkl-universal-instances} has some connections with the~$\srt^2_2$ vs $\coh$ question.
The only construction of solutions to instances of~$\rwkl$ which do not produce
sets of PA degree are variants of Mathias forcing which produce solutions
to~$\seqs{\rwkl}$. In both cases, the solutions are~$\rwkl$-bounding,
and in the latter case, they have a jump of PA degree relative to~$\emptyset'$.
By the previous discussion, one should not expect to prove that~$\rwkl$
admits no universal instance in the usual way, unless one can find a new, suitable, forcing notion.

Some statements do not admit a universal instance because their class of instances is too
restrictive, but they have a natural strengthening which does admit one.
This is for example the case of the rainbow Ramsey theorem for pairs,
which by Miller~\cite{MillerAssorted} admit a universal instance whose solutions are of
DNC degree relative to~$\emptyset'$, but stable variants of the rainbow Ramsey theorem
do not admit one~\cite{Patey2015Degrees}. It is therefore natural to wonder whether there is some 
strengthening of~$\rwkl$ still below~$\rt^2_2$ which admits a universal instance. The Erd\H{o}s-Moser theorem,
defined in the next section, is a good candidate. The following question is a weakening of this interrogation.

\begin{question}
Is there some instance of~$\rt^1_2$ or some computable
instance of~$\rt^2_2$ whose solutions are $\rwkl$-bounding?
\end{question}

Weak K\"onig's lemma is equivalent to several theorems over~$\rca$~\cite{Simpson2009Subsystems}.
However, the Ramsey-type versions of those theorems do not always
give statements equivalent to~$\rwkl$ over~$\rca$. For instance, 
the existence of a separation of two computably inseparable c.e.\ sets
is equivalent to~$\wkl$ over~$\rca$, whereas it is easy to compute
an infinite subset of a separating set.
Among those statements, the Ramsey-type version of the graph coloring problem~\cite{Hirst1990Marriage}
is of particular interest in our study of~$\rwkl$. Indeed,
the $\rwkl$-instances built in~\cite{Bienvenu2015logical,Flood2014Separating,Patey2015Ramsey} are only of two kinds:
trees of positive measure, and trees whose paths code the $k$-colorings
of a computable $k$-colorable graph. The restriction of~$\rwkl$ to the former class of instances
is equivalent to~$\dnr$ over~$\rca$. We shall discuss further the latter one.
Recall that for every~$n \geq 2$, $\rcolor_n$ is the statement ``for every infinite, 
locally $n$-colorable graph $G = (V,E)$, there is an infinite $H \subseteq V$ that is $n$-homogeneous for $G$.''.
The Ramsey-type graph coloring statements have been introduced by Bienvenu, Shafer and the author~\cite{Bienvenu2015logical}.
By Corollary~\ref{cor:RCOLORkisRWKL}, $\rcolor_n$ is equivalent to~$\rwkl$ for~$n \geq 3$,
and by Theorem~\ref{thm-WWKLDoesNotProveRCOLOR2}, 
$\dnr$ (or even weak weak K\"onig's lemma) does not imply~$\rcolor_2$ over~$\omega$-models.

\begin{question}\label{quest:rcolor2-rwkl}
Does~$\rcolor_2$ imply~$\rwkl$ or even~$\dnr$ over~$\rca$?
\end{question}

Separating~$\rcolor_2$ from~$\rwkl$ may require constructing an~$\rwkl$-instance of a new kind,
i.e., not belonging to any of the two classes of instances mentioned before.
The question about the existence of a universal instance also holds for~$\rcolor_2$.

\section{The Erd\H{o}s-Moser theorem}

The Erd\H{o}s-Moser theorem provides, together with the ascending descending principle~($\ads$), 
an alternative decomposition of Ramsey's theorem for pairs. 
Cholak et al.~\cite{Cholak2001strength} proved that every computable $\rt^2_2$-instance has a low${}_2$ solution,
while Mileti~\cite{Mileti2004Partition} and the author (see chapter~\ref{chap-degrees-bounding-principles}) 
showed that various consequences~$\Psf$ of~$\rt^2_2$
do not have~$\Psf$-bounding low${}_2$ degrees, showing therefore that~$\Psf$ does not have
a universal instance. This approach does not apply to~$\emo$ since there is a low${}_2$ $\emo$-bounding degree
(see section~\ref{sect:degrees-bounding-em-low2-degree}).
The Erd\H{o}s-Moser theorem is, together with~$\rwkl$, one of the last Ramsey-type theorems 
for which the existence of a universal instance is unknown.
Since by Theorem~\ref{thm:sem-imp-rwkl}, $\rca \vdash \emo \imp \rwkl$, a positive answer to the following question
would refine our understanding of $\rwkl$-bounding degrees. Note that, like $\rwkl$, the only known
forcing notion for building solutions to $\emo$-instances produces $\emo$-bounding degrees.

\begin{question}
Does $\emo$ admit a universal instance?
\end{question}

Due to the nature of the decomposition of~$\rt^2_2$ into~$\emo$ and~$\ads$,
the Erd\H{o}s-Moser theorem shares many features with~$\rt^2_2$. In particular,
there is a computable $\semo$-instance with no low solution~\cite{Kreuzer2012Primitive}.
The forcing notion used to construct solutions to $\emo$-instances is very similar to the
one used to construct solutions to~$\rt^2_2$-instances. The main difference is that
in the~$\emo$ case, only one object (a transitive subtournament) is constructed,
whereas in the~$\rt^2_2$ case, both a set homogeneous with color~0 and a set homogeneous with color~1
are constructed. As a consequence, the constructions in the $\emo$ case remove the disjunction
appearing in almost every construction of solutions to~$\rt^2_2$,
and therefore simplifies many arguments, while preserving some computational power.
In particular, we proved that~$\coh \leq_c \srt^2_2$ if and only if~$\coh \leq_c \semo$.

Considering that~$\semo$ behaves like~$\srt^2_2$ with respect to~$\coh$,
one would wonder whether $\emo$, like~$\rt^2_2$, implies~$\coh$ over~$\rca$.
The closest result towards an answer is the proof that~$\emo$
implies~$[\sts^2 \vee \coh]$ over~$\rca$, where~$\sts^2$ is the restriction of~$\ts^2$
to stable functions (Theorem~\ref{thm:emo-implies-sts2-or-coh}). 
The Erd\H{o}s-Moser theorem is known not to imply~$\sts^2$ over $\omega$-models (Corollary~\ref{cor:emo-wkl-sts2-sads}),
but the following question remains open.

\begin{question}
Does~$\emo$ imply~$\coh$ over~$\rca$?
\end{question}

A natural first step in the study of the computational strength of a principle
consists in looking how it behaves with respect to ``typical'' sets.
Here, by typical, we mean randomness and genericity. By Liu~\cite{Liu2015Cone},
$\emo$ does not imply the existence of a Martin-Löf random. However, it implies the existence
of a function DNC relative to~$\emptyset'$~\cite{Kang2014Combinatorial,Patey2015Somewhere}, 
which corresponds to the computational power of an infinite subset of a 2-random. 
On the genericity side, $\emo$ implies the existence of a hyperimmune set~\cite{Lerman2013Separating}
and does not imply $\Pi^0_1$-genericity~\cite{Hirschfeldt2009atomic,PateyDominating}.
The relation between~$\emo$ and 1-genericity is currently unclear.
1-genericity received particular attention from the reverse mathematical community recently as it happened
to have close connections with the finite intersection principle~($\fip$).
Dazhafarov \& Mummert~\cite{Dzhafarov2013strength} introduced~$\fip$ in reverse mathematics.
Day, Dzhafarov, and Miller [unpublished] and Hirschfeldt and Greenberg [unpublished] 
independently proved that it is a consequence of the atomic
model theorem (and therefore of~$\srt^2_2$) over~$\rca$.
Downey et al.~\cite{Downey2012finite} first established a link between 1-genericity and~$\fip$
by proving that~$\gens{1}$ implies~$\fip$ over~$\rca$. Later, Cholak et al.~\cite{Cholak2015Any} proved the
other direction.

\begin{question}\label{quest:emo-genericity-rca}
Does~$\emo$ imply~$\gens{1}$ over~$\rca$?
\end{question}

Every 1-generic bounds an $\omega$-model of~$\gens{1}$.
This property is due to a Van Lambalgen-like theorem for genericity~\cite{Yu2006Lowness},
and implies in particular that there is no 1-generic of minimal degree.
Cai et al.~\cite{Cai2015DNR} constructed an $\omega$-model of~$\dnr$
such that the degrees of the second-order part belong to a sequence~$\mathbf{0}, \dbf_1, \dbf_2, \dots$
where~$\dbf_1$ is a minimal degree and~$\dbf_{n+1}$ is a strong minimal cover of~$\dbf_n$.
By the previous remark, such a model cannot contain a 1-generic real. Construcing a similar
model of~$\emo$ would answer negatively Question~\ref{quest:emo-genericity-rca}.
The following question is a first step towards an answer. In particular, 
answering positively would provide a proof that~$\gens{1} \not \leq_c \emo$.

\begin{question}
Does every computable~$\emo$-instance (or even~$\rwkl$-instance) admit 
a solution of minimal degree?
\end{question}

\section{Ramsey-type hierarchies}

Jockusch~\cite{Jockusch1972Ramseys} proved that the hierarchy
of Ramsey's theorem collapses at level 3 over $\omega$-models, that is,
for every~$n, m \geq 3$, $\rt^n_2$ and~$\rt^m_2$ have the same $\omega$-models. 
Simpson~\cite[Theorem III.7.6]{Simpson2009Subsystems} formalized Jockusch's results
within reverse mathematics and proved that $\rt^n_2$
is equivalent to the arithmetic comprehension axiom ($\aca$) over~$\rca$ for every~$n \geq 3$.
Since~$\rt^1_2$ is provable over~$\rca$ and~$\rt^2_2$ is strictly between~$\rt^3_2$ and~$\rca$,
the status of the whole hierarchy of Ramsey's theorem is known.

However, some consequences of Ramsey's theorem form hierarchies
for which the question of the strictness is currently unanswered.
This is for example the case of the thin set theorems, the free set theorem
and the rainbow Ramsey theorem. 
Wang~\cite{Wang2014Some} proved that the full free set hierarchy (hence the thin set hierarchy) lies strictly below~$\aca$
over~$\rca$. We improved this result by proving that 
$\fs$ does not even imply~$\wkl$ (and in fact weak weak König's lemma) over~$\rca$ 
(Corollary~\ref{cor:fs-strong-cb-enum-avoidance}).

\begin{question}
Does the free set theorem (resp.\ the thin set theorem with $\omega$ colors) form a strict hierarchy?
\end{question}

Jockusch~\cite{Jockusch1972Ramseys} proved for every~$n \geq 2$ that every computable~$\rt^n_2$-instance
has a~$\Pi^0_n$ solution and constructed a computable~$\rt^n_2$-instance with no $\Sigma^0_n$ solution.
Cholak et al.~\cite{Cholak2001Free} proved that $\fs^n$ and~$\ts^n$ satisfy the same bounds.
In particular, there is no $\omega$-model of $\fs^n$ or~$\ts^n$ containing only~$\Delta^0_n$ sets.
By Cholak et al.~\cite{Cholak2001strength}, every computable instance of~$\fs^2$ and~$\ts^2$
admit a low${}_2$ solution. Therefore $\fs^2$ (hence $\ts^2$) are strictly weaker than~$\ts^3$ (hence~$\fs^3$) over~$\rca$.
Using this approach, a way to prove the strictness of the hierarchies would be to answer positively the following question.

\begin{question}\label{quest:fs-ts-lown-solution}
For every~$n \geq 3$, does every computable instance of~$\fs^n$ (resp.\ $\ts^n$)
admit a low${}_n$ solution?
\end{question}

The complexity of controlling the $n$th iterate of the Turing jump grows very quickly with~$n$.
While the proof of the existence of low${}_2$ solutions to computable $\rt^2_2$-instances 
using first jump control is relatively simple,
the proof using the second jump control requires already a much more elaborate framework~\cite{Cholak2001strength}.
Before answering Question~\ref{quest:fs-ts-lown-solution}, one may want to get rid of the technicalities
due to the specific combinatorics of the free set and thin set theorems, and focus
on the control of iterated Turing jumps by constructing simpler objects.

Non-effective instances of Ramsey's theorem for singletons is a good starting point,
since the only combinatorial argument involved is the pigeonhole principle.
Moreover, $\rt^1_2$ can be seen as a bootstrap principle, above which the other Ramsey-type
statements are built. For instance, cohesiveness is proven by making $\omega$ applications
of the~$\rt^1_2$ principles,
and $\rt^2_2$ is obtained by making one more application of~$\rt^1_2$ over a non-effective instance.
The proofs of the free set and thin set theorems also make an important use of non-effective instances of~$\rt^1_2$~\cite{Patey2015Combinatorial,Wang2014Some}.

\begin{question}
For every~$n \geq 3$, does every~$\Delta^0_n$ set admit an infinite low${}_n$ subset of either it or its complement?
\end{question}

The solutions to Ramsey-type instances are usually built by forcing.
In order to obtain low${}_n$ solutions, one has in particular 
to $\emptyset^{(n)}$-effectively decide the $\Sigma^0_n$ theory of the generic set.
The forcing relation over a partial order is defined inductively,
and intuitively expresses whether, whatever the choices of conditions extensions we make
in the future, we will still be able to make some progress in satisfying the considered property.

This raises the problem of talking about the ``future'' of a condition~$c$. 
To do that, one needs to be able to describe effectively the conditions extending~$c$.
The problem of the forcing notions used to build solutions to Ramsey-type instances
is that they use variants of Mathias forcing, whose conditions cannot be described effectively.
For instance, let us take the simplest notion of Mathias forcing: pairs~$(F, X)$
where $F$ is a finite set of integers representing the finite approximation of the solution, 
$X$ is a computable infinite ``reservoir'' of integers and~$max(F) < min(X)$. Given a condition~$c = (F, X)$,
the extensions of~$c$ are the pairs~$(H, Y)$ such that~$F \subseteq H$, $Y \subseteq X$
and~$H \setminus F \subseteq X$. Deciding whether a Turing index is the code of 
an infinite computable subset of a fixed computable set requires
a lot of computational power. Cholak et al.~\cite{Cholak2014Generics} studied computable Mathias forcing
and proved that the forcing relation for deciding $\Sigma^0_n$ properties is not~$\Sigma^0_n$ in general.
This is why we need to be more careful in the design of the forcing notions.

In some cases, the reservoir has a particular shape. 
Through their second jump control, Cholak et al.~\cite{Cholak2001strength} first used this idea by noticing
that the only operations over the reservoir were partitioning and finite truncation.
This idea has been reused by Wang~\cite{Wang2014Cohesive} to prove that every computable
$\seqs{\mathsf{D}^2_2}$-instance has a solution of low${}_3$ degree.
We designed in chapter~\ref{chap-controlling-iterated-jumps} new forcing notions
for various Ramsey-type statements, e.g., $\coh$, $\emo$ and $\mathsf{D}^2_2$, in which the forcing relation for
deciding~$\Sigma^0_n$ properties is~$\Sigma^0_n$.

\ifdefined\mypublic
\else
\chapter*{Answered questions}
\addcontentsline{toc}{chapter}{\textcolor{ocre}{Answered questions}}

\begin{question}[Flood and Towsner~\cite{Flood2014Separating}]
For which order $h$ does $\dnr_h$ imply~$\rwkl$ over~$\rca$? 
\end{question}
\begin{answer}
$\dnr_h$ does not imply~$\rwkl$ over~$\rca$ for any computable order~$h$~\cite{Patey2015Ramsey}.
\end{answer}

\begin{question}[Montalban~\cite{Montalban2011Open}]
Is $\tto^2_2$ strictly stronger than~$\rt^2_2$ over~$\rca$?
\end{question}
\begin{answer}
Yes, even stronger than~$\rt^2_2 \wedge \wkl$ by Theorem~\ref{thm:tree-theorem-rt22+wkl-not-tt22}.
\end{answer}

\begin{question}[Mileti~\cite{Mileti2004Partition}]
Do there exist computable $f_1 : [\omega]^2 \to 2$ and $f_2 : [\omega]^2 \to 2$ such there is no computable
$g : [\omega]^2 \to 2$ with the property that every set $H_g$ homogeneous for $g$ computes both a set $H_{f_1}$ homogeneous
for $f_1$ and a set $H_{f_2}$ homogeneous for $f_2$?
\end{question}
\begin{answer}
Yes, by Corollary~\ref{cor:two-ramsey-instances-cannot-be-solved-by-one}.
\end{answer}

\begin{question}[Liu~\cite{Liu2015Cone}]
Is it the case that for every closed set~$\Ccal \subseteq 2^\omega$ with no computable member and every set~$A$,
there is an infinite subset~$H$ of either~$A$ or~$\overline{A}$ such that~$\Ccal$ has no $H$-computable member?
\end{question}
\begin{answer}
No. We can even take~$\Ccal$ to be effectively closed and $A$ to be~$\Delta^0_2$
as a counter-example~\cite{Patey2015Combinatorial}.
\end{answer}

\begin{question}[Flood~\cite{Flood2012Reverse}]
Does $\dnr$ imply~$\rwkl$ over~$\rca$?
\end{question}
\begin{answer}
No. By Theorem~\ref{RCOLOR2hasNRA}, $\rwkl$ has a computable instance such that the measure of oracles computing a solution 
is null. This is a joint work with Laurent Bienvenu and Paul Shafer.
\end{answer}

\begin{question}[Flood~\cite{Flood2012Reverse}]
Does $\srt^2_2$ imply~$\rwkls{2}$ over~$\rca$?
\end{question}
\begin{answer}
No. By~\cite{Bienvenu2015logical}, $\rwkls{2}$ implies~$\dnrs{2}$ over~$\rca+\bst$
and therefore there is no model of~$\rca \wedge \bst \wedge \rwkls{2}$ with only low sets,
whereas~$\rca \wedge \bst \wedge \srt^2_2$ has one.
\end{answer}

\begin{question}[Wang~\cite{Wang2014Definability}]
Does any of~$\coh$ and~$\emo$ admit preservation of the arithmetic hierarchy?
\end{question}
\begin{answer}
Yes, both do~\cite{Patey2016Controlling}.
\end{answer}

\begin{question}[Hirschfeldt~\cite{Hirschfeldt2015Slicing}]
Let $n \geq 3$. Does $\fs^n$ or $\ts^n$ or~$\fs$ or $\ts$ imply~$\wkl$ over~$\rca$?
\end{question}
\begin{answer}
No. By~\cite{Patey2015Combinatorial}, they do not even imply~$\wwkl$ over~$\rca$.
\end{answer}

\begin{question}[Hirschfeldt~\cite{Hirschfeldt2015Slicing}, Montalban~\cite{Montalban2011Open}]
Does $\fs^2$ imply~$\rt^2_2$ (or $\srt^2_2$, $\ads$ or~$\cac$) over~$\rca$?
What about over~$\rca \wedge \wkl$?
\end{question}
\begin{answer}
No. By Corollary~\ref{cor:fs2-not-imply-sts2n}, $\fs \wedge \wkl$ implies none of~$\srt^2_2$, $\ads$ or~$\cac$ over~$\rca$.
\end{answer}

\begin{question}[Montalban~\cite{Montalban2011Open}]
Does~$\srt^2_2$ imply $\ipt^2_2$ over~$\rca$?
\end{question}
\begin{answer}
No. By~\cite{Patey2015Open}, $\ipt^2_2$ and~$\rwkls{2}$ are equivalent over~$\rca$.
\end{answer}

\begin{question}[Murakami~\cite{Murakami2014Ramseyan}]
Does~$\ads$ imply~$\mathsf{psRT}^2_2$ over~$\rca$?
\end{question}
\begin{answer}
Yes. See the Logic Blog 2015.
\end{answer}

\begin{question}[Wang~\cite{Wang2013Omitting}]
Is there a computable sequence of sets~$R_0, R_1, \dots$ with no computable $\vec{R}$-cohesive set
but which admits an one computable from a 2-generic real?
\end{question}
\begin{answer}
Yes. By Lemma~\ref{lem:cohesiveness-unsolvable-but-2-generic}.
\end{answer}

\begin{question}[Wang~\cite{Wang2013Omitting}]
Is there a computable sequence of sets~$R_0, R_1, \dots$ with no computable $\vec{R}$-cohesive set
but which admits an one computable from a 3-random real?
\end{question}
\begin{answer}
No. By Lemma~\ref{lem:randomness-p-cohesive}.
\end{answer}

\begin{question}[Dorais, Dzhafarov, Hirst, Mileti and Shafer~\cite{Dorais2016uniform}]
If $n, j, k \geq 2$ and $j < k$, is it the case that $\rt^n_k \not \leq_W \rt^n_j$?
\end{question}
\begin{answer}
Yes. By Corollary~\ref{cor:ramsey-colors-not-computable-reducible}, we even have~$\rt^n_k \not \leq_c \rt^n_j$.
\end{answer}

\begin{question}[Jockusch and Hirschfeldt~\cite{Hirschfeldtnotions}]
Let $k > j > 2$, where $k$ may also be $<\infty$. Is $\rt^1_k \leq_{sc} \srt^2_j$?
\end{question}
\begin{answer}
No. This is a joint work with Damir Dzhafarov, Reed Solomon and Brown Westrick.
\end{answer}

\begin{question}[Dzhafarov~\cite{DzhafarovStrong}]
For $k > 2$, is it the case that $\coh \leq_{sc} \srt^2_k$?
Is it the case that $\coh \leq_{sc} \srt^2_{<\infty}$?
\end{question}
\begin{answer}
No. This is a joint work with Damir Dzhafarov, Reed Solomon and Brown Westrick.
\end{answer}

\begin{question}[Jockusch and Hirschfeldt~\cite{Hirschfeldtnotions}]
Is it the case that for every function $f : \omega \to 3$ with no computable $f$-thin set and
every function~$g : \omega \to 2$, there is an infinite $g$-homogeneous set that
does not compute any infinite $f$-homogeneous set?
\end{question}
\begin{answer}
No. Build a $\Delta^0_2$ coloring $f : \omega \to 3$ with no computable $f$-thin set, 
but such that any $f$-thin set for color 2 is sparse enough to compute the halting set.
Then, define $g(x) = 1$ iff $f(x) = 2$. Any $g$-homogeneous set for color 1 is $f$-homogeneous for color 2, 
and any $g$-homogeneous set for color 0 or 1 computes the halting set and therefore an $f$-homogeneous set.
\end{answer}

\begin{question}[Montalb\'an, see~\cite{Hirschfeldtnotions}]
Let $2 \leq j < k$.
Is there a coloring $f : \omega \to k$ such that for every $g : \omega \to j$,
there is an infinite $g$-homogeneous set that does not compute any infinite $f$-homogeneous set?
\end{question}
\begin{answer}
Yes. By Corollary~\ref{cor:colors-ramsey-d2-general-colors}, we can even take~$f$ to be low.
\end{answer}

\begin{question}[Kang~\cite{Kang2014Combinatorial}]
Does~$\ts^2$ imply~$\rrt^2_2$ over~$\rca$?
\end{question}
\begin{answer}
Yes. See Theorem~\ref{thm:ts2-dnrzp}.
\end{answer}

\begin{question}[Hirschfeldt, Shore, Slaman~\cite{Hirschfeldt2009atomic}]
Does~$\coh$ (or $\cads$) imply $\amt$ over~$\rca$?
\end{question}
\begin{answer}
No, not even over~$\rca \wedge \wkl$. See~\cite{PateyDominating}.
\end{answer}

\begin{question}[Cholak, Giusto, Hirst and Jockusch~\cite{Cholak2001Free}, Montalb\'an~\cite{Montalban2011Open}]
Does any of $\fs^2$, $\ts^2$, $\fs^2 \wedge \wkl$ or $\fs^2 \wedge \coh$ imply~$\rt^2_2$ over~$\rca$?
\end{question}
\begin{answer}
No. By Corollary~\ref{cor:fs2-not-imply-sts2n}, $\fs \wedge \wkl \wedge \coh$ does not imply~$\srt^2_2$ over~$\rca$.
\end{answer}

\begin{question}[Dzhafarov and Jockusch~\cite{Dzhafarov2009Ramseys}]
Is is the case that for every $\dbf \gg \emptyset'$ and every $C \not \leq_T \emptyset$, 
every computable 2-coloring of pairs admits an infinite homogeneous set $H \not \geq_T C$ with $\deg(H') \leq \dbf$?
\end{question}
\begin{answer}
Yes. If $\deg(C) \not \leq \dbf$, this follows from~\cite{Cholak2001strength}.
If~$C$ is~$\Delta^0_2$, this follows from~\cite{Dzhafarov2009Ramseys}.
If~$C$ is not $\Delta^0_2$ but $\deg(C) \leq \dbf$, then it is similar to the previous case.
We simply use the low basis theorem instead of Theorem 2.1. Indeed, because $C$ is not $\Delta^0_2$,
any low set is necessarily cone avoiding.
Notice that Lemma 4.6 needed the assumption~$C \leq_T \emptyset'$ only because of Theorem~2.1.
Also notice that Lemma 4.7 only needs the assumption that~$\deg(C) \leq \dbf$.
\end{answer}

\begin{question}[Jockusch and Hirschfeldt~\cite{Hirschfeldtnotions}]
Let $n \geq 2$ and $k > j \geq 2$. For which~$m$ do we have~$\rt^n_k \leq_\omega^m \rt^n_j$?
\end{question}
\begin{answer}
See Theorem~\ref{thm:colors-ramsey-general} for a full characterization.
\end{answer}

\begin{question}[Brattka and Rakotoniaina~\cite{Brattka2015Uniform}]
Does $\wkl[\emptyset'] \leq_W \srt^3_2$ hold?
\end{question}
\begin{answer}
No. This does not even hold for computable reducibility.
Indeed, $\srt^3_2 =_c \rt^2_2[\emptyset']$. We now show that $\wkl[\emptyset'] \not \leq_c \rt^2_2[\emptyset']$.
Relativizing Liu's theorem to~$\emptyset'$, we obtain~``For every set~$C$ such that~$C \oplus \emptyset' \not \gg \emptyset'$
and for every set~$A$, there is an infinite subset~$H$ of~$A$ or~$\overline{A}$ such that
$H \oplus C \oplus \emptyset' \not \gg \emptyset'$''.
Fix a $\Delta^0_2$ $\rt^2_2$-instance~$f : [\omega]^2 \to 2$ and let~$R_0, R_1, \dots$
be the sequence of $\Delta^0_2$ sets defined by~$R_x = \{ y : f(x, y) = 1 \}$.
Let~$C$ be a sufficiently generic real for the Mathias forcing whose conditions~$(F, X)$
are such that~$X \oplus \emptyset' \not \gg \emptyset'$. 
By the relativization of Liu's theorem, $C$ is $\vec{R}$-cohesive, and one easily checks that~$C \oplus \emptyset' \not \gg \emptyset'$.
Let~$A = \{ x : \lim_{s \in C} f(x, s) = 1 \}$. Using relativized Liu's theorem,
there is an infinite subset~$H$ of~$A$ or~$\overline{A}$ such that~$H \oplus C \oplus \emptyset' \not \gg \emptyset'$.
In particular, $H \oplus C \oplus \emptyset'$ computes an infinite $f$-homogeneous set.
\end{answer}

\begin{question}[Brattka and Rakotoniaina~\cite{Brattka2015Uniform}]
Does $\ivt \leq_W \rt^2_2$ hold?
\end{question}
\begin{answer}
Yes. Even~$\ivt \leq_W \srt^2_2$ holds.
Fix a computable enumeration~$q_0, q_1, \dots$ of all the rationals.
Let~$f : [0,1] \to \Rb$ be a continuous function such that~$f(0) < 0 < f(1)$.
Define~$C = \{ i \in \omega : f(q_i) = 0 \}$. The set~$C$ is $f$-co-c.e.\ 
and we can fix an $f$-computable sequence of its approximations~$C_0 \supseteq C_1 \supseteq \dots$
such that~$C_i \neq \emptyset$ for every~$i$.
Let~$g : [\omega]^2 \to 2$ be defined by~$g(x, y) = 1$ iff $min(C_x) = min(C_y)$.
The function~$g$ is stable and $f$-computable. By $\srt^2_2$, 
let~$H$ be an infinite $f$-homogeneous set for some color~$c < 2$.
We have two cases.
\begin{itemize}
	\item Case $c = 0$. In this case, we claim that~$C$ is empty. Indeed, if not,
	then $min(C_y) = min(C)$ for almost every~$y$. Taking any two sufficiently large~$y_0, y_1 \in H$,
	$min(C_{y_0}) = min(C) = min(C_{y_1})$, hence~$g(y_0, y_1) = 1$ contradicting~$c = 0$.
	In this case, we know that $f(q) \neq \emptyset$ for every rational~$q$,
	so we can compute some $y$ such that~$f(y) = 0$ using the usual dichotomy method.

	\item Case $c = 1$. Let~$q = min(C_x)$ for some~$x \in H$. We claim that~$f(q) = 0$.
	Indeed, $g(x, y) = 1$ for infinitely many $y$'s. Unfolding the definition of~$g$,
	$min(C_x) = min(C_y)$ for infinitely many~$y$'s, so~$min(C_x) = min(C)$.
\end{itemize}
\end{answer}

\fi


\ifdefined\mypublic
\else
\chapter*{Publications}
\addcontentsline{toc}{chapter}{\textcolor{ocre}{Publications}}
\nocite{Patey2015proof,Patey2014Complexity,Patey2015Iterative,
PateyDominating,Patey2015Combinatorial,Patey2015complexity,
Patey2016Controlling,Patey2015Degrees,Patey2015Open,
Patey2015Ramsey,Patey2015Somewhere,Patey2015strength,
Patey2016weakness,Patey2013note,Frittaion2015Coloring,
Bienvenu2015logical,Bienvenu2016Diagonally,Patey2015Iterativea,
Dzhafarov2016Ramseys,PateyPartial}

\section*{Journal papers}
\printbibliography[heading=bibempty,keyword={patey-journal}]

\section*{Conference papers}
\printbibliography[heading=bibempty,keyword={patey-conference}]

\section*{Submitted papers}
\printbibliography[heading=bibempty,keyword={patey-submitted}]

\section*{Drafts}
\printbibliography[heading=bibempty,keyword={patey-draft}]
\fi

\chapter*{Bibliography}
\addcontentsline{toc}{chapter}{\textcolor{ocre}{Bibliography}}
\printbibliography[heading=bibempty]




\addcontentsline{toc}{chapter}{\textcolor{ocre}{Index}}
\printindex

\appendix
\addcontentsline{toc}{part}{Appendices}

\chapter{A simple separation of~$\ads$ from~$\scac$}

Lerman, Solomon and Towsner~\cite{Lerman2013Separating} developped a framework to separate statements over computable entailment,
and proved in particular that $\ads$ does not imply~$\scac$ over~$\rca$.
The author~\cite{Patey2015Iterative} refined their technique to obtain a simpler framework.
In this appendix, we apply this simplification to the separation of~$\ads$ from~$\scac$.
We also generalize the separation to a pseudo Ramsey theorem.

\section{A pseudo Ramsey theorem}

Pseudo Ramsey's theorem for pairs has been introduced by Murakami, Yamazaki
and Yokoyama in~\cite{Murakami2014Ramseyan} in the context of an analysis of a Ramseyan factorization theorem. 
They proved that it is between the chain antichain principle
and the ascending descending sequence principle over~$\rca$, and asked whether
it was equivalent to one of them. We answer positively. The proof was indenpendently found
by Silvia Steila, see~\cite{SteilaReverse}.

\index{pseudo Ramsey's theorem}
\index{psRT@$\mathsf{psRT}^n_k$|see {pseudo Ramsey's theorem}}
\index{semi-transitive}
\index{pseudo-homogeneous}
\begin{definition}[Pseudo Ramsey's theorem]
A coloring~$f : [\mathbb{N}]^2 \to 2$ is \emph{semi-transitive}
if whenever~$f(x, y) = 1$ and~$f(y, z) = 1$, then~$f(x, z) = 1$ for $x < y < z$.
A set~$H = \{x_0 < x_1 < \dots \}$ is \emph{pseudo-homogeneous} for a coloring~$f : [\mathbb{N}]^n \to k$
if~$f(x_i, \dots, x_{i+n-1}) = f(x_j, \dots, x_{j+n-1})$ for every~$i, j \in \mathbb{N}$.
$\mathsf{psRT}^n_k$ is the statement ``Every coloring~$f : [\mathbb{N}]^n \to k$ has an infinite pseudo-homogeneous set''.
\end{definition}

\begin{theorem}\label{thm:psrt22-is-ads}
$\mathsf{RCA}_0 \vdash \mathsf{psRT}^2_2 \leftrightarrow \ads$
\end{theorem}
\begin{proof}
The direction~$\mathsf{psRT}^2_2 \rightarrow \ads$ is Theorem~24 in~\cite{Murakami2014Ramseyan}.
We prove that~$\mathsf{ADS} \rightarrow \mathsf{psRT}^2_2$. 
Let~$f : [\mathbb{N}]^2 \to 2$ be a coloring. The reduction is in two steps.
We first define a $\Delta^{0,f}_1$ semi-transitive coloring
$g : [\mathbb{N}]^2 \to 2$ such that every infinite set pseudo-homogeneous for~$g$
computes an infinite set pseudo-homogeneous for~$f$.
Then, we define a $\Delta^{0,g}_1$ linear order~$h : [\mathbb{N}]^2 \to 2$
such that every infinite set pseudo-homogeneous for~$h$
computes an infinite set pseudo-homogeneous for~$g$.
We conclude by applying~$\ads$ over~$h$.

\emph{Step 1}:
Define the coloring~$g : [\mathbb{N}]^2 \to 2$ for every~$x < y$ by
$g(x, y) = 1$ if there exists a sequence~$x = x_0 < \dots < x_l = y$
such that $f(x_i, x_{i+1}) = 1$ for every~$i < l$, and~$g(x, y) = 0$ otherwise.
The function~$g$ is a semi-transitive coloring. Indeed,
suppose that~$g(x, y) = 1$ and~$g(y, z) = 1$, witnessed respectively
by the sequences~$x = x_0 < \dots < x_m = y$ and~$y = y_0 < \dots < y_n = z$.
The sequence~$x = x_0 < \dots < x_m = y_0 < \dots < y_n = z$ witnesses~$g(x, z) = 1$.
We claim that every infinite set~$H = \{x_0 < x_1 < \dots \}$ pseudo-homogeneous for~$g$ computes
an infinite set pseudo-homogeneous for~$f$. If $H$ is pseudo-homogeneous with color~0,
then $f(x_i, x_{i+1}) = 0$ for each~$i$, otherwise the sequence~$x_i < x_{i+1}$ would witness
$g(x_i, x_{i+1}) = 1$. Thus $H$ is pseudo-homogeneous for~$f$ with color~0. 
If $H$ is pseudo-homogeneous with color~1,
then define the set~$H_1 \supseteq H$ to be the set of integers in the sequences
witnessing~$g(x_i, x_{i+1}) = 1$ for each~$i$. The set~$H_1$ is~$\Delta^{0,f \oplus H}_1$
and pseudo-homogeneous for~$f$ with color~1.

\emph{Step 2}:
Define the coloring~$h : [\mathbb{N}]^2 \to 2$ for every~$x < y$ by
$h(x, y) = 0$ if there exists a sequence~$x = x_0 < \dots < x_l = y$
such that $g(x_i, x_{i+1}) = 0$ for every~$i < l$, and~$h(x, y) = 1$ otherwise.
For the same reasons as for~$g$, $h(x, z) = 0$ whenever~$h(x, y) = 0$ and~$h(y, z) = 0$ for~$x < y < z$.
We need to prove that if~$h(x, z) = 0$ then either $h(x, y) = 0$ or~$h(y, z) = 0$ for~$x < y < z$.
Let~$x = x_0 < \dots < x_l = z$ be a sequence witnessing $h(x,z) = 0$.
If~$y = x_i$ for some~$i < l$ then the sequence~$x = x_0 < \dots < x_i = y$ witnesses~$h(x,y) = 0$.
If~$y \neq x_i$ for every~$i < l$, then there exists some~$i < l$ such that~$x_i < y < x_{i+1}$. 
By semi-transitivity of~$g$, either $g(x_i, y) = 0$ or~$g(y, x_{i+1}) = 0$. In this case either~$x = x_0 < \dots < x_i < y$ witnesses~$h(x,y) = 0$ or~$y < x_{i+1} < \dots < x_l = z$ witnesses~$h(y,z) = 0$.
Therefore~$h$ is a linear order.
For the same reasons as for~$g$, every infinite set pseudo-homogeneous for~$h$
computes an infinite set pseudo-homogeneous for~$g$.
This last step finishes the proof.
\end{proof}

\section{Dependent hyperimmunity}

In this section, we design the computability-theoretic property
which will enable us to separate $\ads$ from~$\scac$.
The notion of essential formula is inspired by the construction
of Lerman, Solomon and Towsner~\cite{Lerman2013Separating}.

\index{dependent hyperimmunity}
\begin{definition}[Dependent hyperimmunity]
A formula~$\varphi(U, V)$ is \emph{essential}
if for every~$x \in \omega$, there is a finite set~$R > x$ such that for every~$y \in \omega$,
there is a finite set~$S > y$ such that~$\varphi(R, S)$ holds.
A pair of sets~$A_0, A_1 \subseteq \omega$ is \emph{dependently $X$-hyperimmune} 
if for every essential $\Sigma^{0,X}_1$ formula~$\varphi(U, V)$,
$\varphi(R, S)$ holds for some~$R \subseteq \overline{A}_0$ and~$S \subseteq \overline{A}_1$.
\end{definition}

Being essential is a density property.
Intuitively, a formula~$\varphi(U, V)$ is essential if we can find arbitrarily large
values to satisfy it. Therefore, one can ensure by a $\Delta^0_2$ construction of some sets~$A_0, A_1$
based on the finite extension method that every essential formula will have a valuation
included in~$\overline{A}_0$ and~$\overline{A}_1$.

If the pair $A_0, A_1$ is dependently hyperimmune, then $A_0$ and~$A_1$ are both hyperimmune.
A reversal can be obtained for co-c.e.\ sets.

\begin{lemma}\label{lem:dependent-hyperimmunity-for-coce}
Fix two sets~$A_0, A_1$ such that~$A_0$ is $X$-co-c.e.
The pair~$A_0, A_1$ is dependently $X$-hyperimmune iff $A_0$ and $A_1$ are $X$-hyperimmune.
\end{lemma}
\begin{proof}
We first show that if~$A_0$ and~$A_1$ are dependently $X$-hyperimmune
then both~$A_0$ and~$A_1$ are $X$-hyperimmune. Let~$F_0, F_1, \dots$ be a $X$-c.e.\ array.
Let $\varphi(U, V)$ be the $\Sigma^{0,X}_1$ formula which holds if~$U = F_i$ for some~$i \in \omega$.
The formula $\varphi(U, V)$ is essential, therefore there $\varphi(R, S)$ holds for some finite set~$R \subseteq \overline{A}_0$
and $S \subseteq \overline{A}_1$. In particular, $R = F_i$ for some~$i \in \omega$, therefore $F_i \subseteq \overline{A}_0$
and $A_0$ is hyperimmune. Similarly, the $\Sigma^{0,X}_1$ formula $\psi(U, V)$ which holds if $V = F_i$ for some~$i \in \omega$
witnesses that $A_1$ is hyperimmune.

We now prove that if $A_0$ and~$A_1$ are $X$-co-c.e.\ and $X$-hyperimmune, then
the pair $A_0, A_1$ is dependently $X$-hyperimmune.
Let~$\varphi(U, V)$ be an essential $\Sigma^{0,X}_1$ formula.
Define an $X$-c.e.\ sequence of sets $F_0 < F_1 < \dots$ such that for every~$i \in \omega$,
there is some $R < F_i$ such that~$\varphi(R, F_i)$ holds and~$R \subseteq \overline{A}_0$.
First, notice that the sequence is $X$-c.e.\ since $A_0$ is $X$-co-c.e.
Second, we claim that the sequence is infinite. To see this,
define an $X$-c.e.\ array $E_0 < F_1 < \dots$ such that for every~$i \in \omega$, there is some finite set~$S > E_i$
such that $\psi(E_i, S)$ holds. The array is infinite since $\psi(U, V)$ is essential.
Since $A_0$ is $X$-hyperimmune, there are infinitely many $i$'s such that $E_i \subseteq \overline{A}_0$.
Last, by $X$-hyperimmunity of~$A_1$, there is some~$i \in \omega$ such that~$F_i \subseteq \overline{A}_1$.
By definition of~$F_i$, there is some~$R \subseteq \overline{A}_0$ such that~$\varphi(R, F_i)$ holds.
\end{proof}

\index{preservation!of dependent hyperimmunity}
\begin{definition}[Preservation of dependent hyperimmunity]
A $\Pi^1_2$ statement~$\Psf$ admits \emph{preservation of dependent hyperimmunity} if
for every set $Z$, every pair of dependently $Z$-hyperimmune sets~$A_0, A_1 \subseteq \omega$
and every $\Psf$-instance~$X \leq_T Z$, there is a solution~$Y$ to~$X$ such that $A_0, A_1$ are dependently 
$Y \oplus Z$-hyperimmune.
\end{definition}

As usual, note that preservation of dependent hyperimmunity is a preservation property,
therefore if some $\Pi^1_2$ statement~$\Psf$ admits preservation of dependent hyperimmunity while
another statement~$\Qsf$ does not, then by Lemma~\ref{lem:intro-reduc-preservation-separation}, 
$\Psf$ does not computably entail $\Qsf$.
We start by proving that~$\scac$ does not admit preservation of dependent hyperimmunity.

\begin{theorem}\label{thm:stable-semi-transitive-fairness}
There exists a computable, stable semi-transitive coloring~$f : [\omega]^2 \to 2$
such that the pair~$\overline{A}_0, \overline{A}_1$ is 
dep.\ hyperimmune, where
$A_i = \{ x : \lim_s f(x, s) = i \}$.
\end{theorem}
\begin{proof}
Fix an enumeration~$\varphi_0(U,V), \varphi_1(U,V), \dots$ of all $\Sigma^0_1$ formulas.
The construction of the function~$f$ is done by a finite injury priority argument
with a movable marker procedure. We want to satisfy the following scheme of requirements for each~$e$,
where~$A_i = \{ x : \lim_s f(x, s) = i \}$:
\[
\Rcal_e : \varphi_e(U,V) \mbox{ essential} \imp 
	(\exists R \subseteq_{fin} A_0)(\exists S \subseteq_{fin} A_1)\varphi_e(R,S)
\]

The requirements are given the usual priority ordering.
We proceed by stages, maintaining two sets $A_0, A_1$
which represent the limit of the function~$f$.
At stage 0, $A_{0,0} = A_{1,0} = \emptyset$ and $f$ is nowhere defined. 
Moreover, each requirement~$\Rcal_e$ is given a movable marker~$m_e$ initialized to~0.

A strategy for~$\Rcal_e$ \emph{requires attention at stage~$s+1$}
if $\varphi_e(R, S)$ holds for some~$R < S \subseteq (m_e, s]$.
The strategy sets~$A_{0,s+1} = (A_{0,s} \setminus (m_e, min(S)) \cup [min(S), s]$ 
and~$A_{1,s+1} = (A_{1,s} \setminus [min(S), s]) \cup (m_e, min(S))$.
Note that~$(m_e, min(S)) \cap [min(S), s]= \emptyset$ since~$R < S$.
Then it is declared~\emph{satisfied} and does not act until some strategy of higher priority changes its marker. 
Each marker~$m_{e'}$ of strategies of lower priorities is assigned the value~$s+1$.

At stage~$s+1$, assume that~$A_{0,s} \cup A_{1,s} = [0,s)$
and that~$f$ is defined for each pair over~$[0,s)$.
For each~$x \in [0,s)$, set~$f(x,s) = i$ for the unique~$i$ such that~$x \in A_{i,s}$.
If some strategy requires attention at stage~$s+1$, take the least one
and satisfy it.
If no such requirement is found, set~$A_{0,s+1} = A_{0,s} \cup \{s\}$
and~$A_{1,s+1} = A_{1,s}$. Then go to the next stage. This ends the construction.

Each time a strategy acts, it changes the markers of strategies of lower priority, and is declared satisfied.
Once a strategy is satisfied, only a strategy of higher priority can injure it.
Therefore, each strategy acts finitely often and the markers stabilize.
It follows that the $A$'s also stabilize and that~$f$ is a stable function.

\begin{claim}
For every~$x < y < z$, $f(x,y) = 1 \wedge f(y,z) = 1 \imp f(x,z) = 1$
\end{claim}
\begin{proof*}
Suppose that $f(x,y) = 1$ and~$f(y,z) = 1$ but $f(x,z) = 0$. 
By construction of~$f$, $x \in A_{0,z}$, $x \in A_{1,y}$ and~$y \in A_{1,z}$.
Let~$s \leq z$ be the last stage such that~$x \in A_{1,s}$.
Then at stage~$s+1$, some strategy~$\Rcal_e$ receives attention
and moves~$x$ to~$A_{0,s+1}$ and therefore moves~$[x, s]$ to~$A_{0,s+1}$.
In particular~$y \in A_{0,s+1}$ since~$y \in [x,s]$.
Moreover, the strategies of lower priority have had their marker moved to~$s+1$
and therefore will never move any element below~$s$.
Since $f(y, z) = 1$, then $y \in A_{1,z}$.
In particular, some strategy~$\Rcal_i$ of higher priority moved~$y$ to~$A_{1,t+1}$
at stage~$t +1$ for some~$t \in (s, z)$. Since~$\Rcal_i$ has a higher priority, $m_i \leq m_e$, 
and since~$y$ is moved to~$A_{1,t+1}$, then so is~$[m_i, y]$, and in particular
$x \in A_{1,t+1}$ since~$m_i \leq m_e \leq x \leq y$. This contradicts the maximality of~$s$.
\end{proof*}

\begin{claim}
For every~$e \in \omega$, $\Rcal_e$ is satisfied.
\end{claim}
\begin{proof*}
By induction over the priority order. Let~$s_0$ be a stage after which
no strategy of higher priority will ever act. By construction, $m_e$ will not change after stage~$s_0$.
If $\varphi_e(U,V)$ is essential, then $\varphi_e(R,S)$ holds for two sets~ $m_e < R < S$.
Let~$s = 1+max(s_0, S)$. The strategy~$\Rcal_e$ will require attention at some stage before~$s$,
will receive attention, be satisfied and never be injured.
\end{proof*}

This last claim finishes the proof.
\end{proof}

\begin{corollary}\label{cor:scac-not-cac-fairness-preservation}
$\scac$ does not admit preservation of dependent hyperimmunity.
\end{corollary}
\begin{proof}
Let~$f : [\omega]^2 \to 2$ be the coloring of Theorem~\ref{thm:stable-semi-transitive-fairness}.
By construction, the pair $\overline{A}_0$, $\overline{A}_1$ is dependently hyperimmune, 
where~$A_i = \{ x : \lim_s f(x, s) = i \}$.
Let~$H$ be an infinite $f$-homogeneous set. In particular, $H \subseteq A_0$ or~$H \subseteq A_1$.
We claim that the pair~$\overline{A}_0, \overline{A}_1$ is not dependently~$H$-hyperimmune.
The $\Sigma^{0,H}_1$ formula~$\varphi(U, V)$ defined by~$U \neq \emptyset \wedge V \neq \emptyset \wedge U \cup V \subseteq H$
is essential since~$H$ is infinite. However, if there is some~$R \subseteq A_1$ and~$S \subseteq A_0$
such that~$\varphi(R, S)$ holds, then~$H \cap A_0 \neq \emptyset$ and~$H \cap A_1 \neq \emptyset$,
contradicting the choice of~$H$. Therefore $\overline{A}_0, \overline{A}_1$ is not dependently~$H$-hyperimmune.
Hirschfeldt and Shore~\cite{Hirschfeldt2007Combinatorial} proved that
$\scac$ is equivalent to stable semi-transitive Ramsey's theorem for pairs
over~$\rca$. Therefore~$\scac$ does not admit preservation of dependent hyperimmunity.
\end{proof}

We will now prove the positive preservation result.

\begin{theorem}\label{thm:psrt2k-not-cac-fairness-preservation}
For every~$k \geq 2$, $\psrt^2_k$ admits preservation of dep.\ hyperimmunity.
\end{theorem}
\begin{proof}
The proof is done by induction over~$k \geq 2$.
Fix a pair of sets~$A_0, A_1 \subseteq \omega$ dependently $Z$-hyperimmune for some set~$Z$.
Let~$f : [\omega]^2 \to k$ be a $Z$-computable coloring
and suppose that there is no infinite $f$-thin set~$H$ such that
the pair $A_0, A_1$ is dependently $H \oplus Z$-hyperimmune,
as otherwise, we are done by induction hypothesis.
We will build $k$ infinite sets~$G_0, \dots, G_{k-1}$
such that $G_i$ is pseudo-homogeneous for~$f$ with color~$i$ for each~$i < k$
and such that $A_0, A_1$ is dependently $G_i \oplus Z$-hyperimmune for some~$i < k$.
The sets~$G_0, \dots, G_{k-1}$ are built by a variant of Mathias forcing
$(F_0, \dots, F_{k-1}, X)$ such that
\begin{itemize}
	\item[(i)] $F_i \cup \{x\}$ is pseudo-homogeneous for~$f$ with color~$i$ for each~$x \in X$
	\item[(ii)] $X$ is an infinite set such that $A_0, A_1$ is dependently $X \oplus Z$-hyperimmune
\end{itemize}
A condition~$d = (H_0, \dots, H_{k-1}, Y)$ \emph{extends} $c = (F_0, \dots, F_{k-1}, X)$
(written $d \leq c$) if~$(H_i, Y)$ Mathias extends $(F_i, X)$ for each~$i < k$.
A tuple of sets~$G_0, \dots, G_{k-1}$ \emph{satisfies} $c$ if
for every~$n \in \omega$, there is an extension~$d = (H_0, \dots, H_{k-1}, Y)$ of~$c$
such that~$G_i \uh n \subseteq H_i$ for each~$i < k$.
Informally, $G_0, \dots, G_{k-1}$ satisfy~$c$ if the sets are generated by a decreasing
sequence of conditions extending~$c$.
In particular, $G_i$ is pseudo-homogeneous for~$f$ with color~$i$
and satisfies the Mathias condition~$(F_i, X)$. 
The first lemma shows that every sufficiently generic filter yields a $k$-tuple of infinite sets.

\begin{lemma}\label{lem:psrt-cac-fairness-force-Q}
For every condition~$c = (F_0, \dots, F_{k-1}, X)$ and every~$i < k$,
there is an extension~$d = (H_0, \dots, H_{k-1}, Y)$ of~$c$ such that~$|H_i| > |F_i|$.
\end{lemma}

Fix an enumeration~$\varphi_0(G, U,V), \varphi_1(G, U,V), \dots$ of all $\Sigma^{0,Z}_1$ formulas.
We want to satisfy the following requirements for each~$e_0, \dots, e_{k-1} \in \omega$:
$$
\Rcal_{\vec{e}} : \hspace{10pt} \Rcal^{G_0}_{e_0} \hspace{10pt} 
		\vee \hspace{10pt} \dots \hspace{10pt} \vee \hspace{10pt} \Rcal^{G_{k-1}}_{e_{k-1}}
$$
where $\Rcal_e^G$ is the requirement ``$\varphi_e(G, U, V)$ essential $\imp \varphi_e(G, R, S)$ 
for some~$R \subseteq \overline{A}_0$ and~$S \subseteq \overline{A}_1$''.
We say that a condition~$c$ \emph{forces} $\Rcal_{\vec{e}}$ if $\Rcal_{\vec{e}}$  holds
for every $k$-tuple of sets satisfying~$c$. Note that the notion of satisfaction 
has a precise meaning given above.

\begin{lemma}\label{lem:psrt-cac-fairness-force-R}
For every condition~$c$ and every $k$-tuple of indices~$e_0, \dots, e_{k-1} \in \omega$,
there is an extension~$d$ of~$c$ forcing~$\Rcal_{\vec{e}}$.
\end{lemma}
\begin{proof*}
Fix a condition~$c = (F_0, \dots, F_{k-1}, X)$. 
Let~$\psi(U, V)$ be the~$\Sigma^{0,X \oplus Z}_1$ formula which holds if
there is a $k$-tuple of sets~$E_0, \dots, E_{k-1} \subseteq X$ and a $z \in X$ such that for each~$i < k$,
\begin{itemize}
	\item[(i)] $z > max(E_i)$
	\item[(ii)] $F_i \cup E_i \cup \{z\}$ is pseudo-homogeneous for color~$i$.
	\item[(iii)] $\varphi_{e_i}(F_i \cup E_i, U_i, V_i)$ holds for some~$U_i \subseteq U$ and~$V_i \subseteq V$
\end{itemize}
Suppose that~$c$ does not force~$\Rcal_{\vec{e}}$, otherwise we are done.

We claim that~$\psi$ is essential. 
Since~$c$ does not force~$\Rcal_{\vec{e}}$, there is a $k$-tuple of infinite sets~$G_0, \dots, G_{k-1}$
satisfying~$c$ and such that $\varphi_{e_i}(G_i, U, V)$ is essential for each~$i < k$.
Fix some $x \in \omega$. By definition of being essential, there are some finite sets
$R_0, \dots, R_{k-1} > x$ such that for every~$y \in \omega$, there are finite sets~$
S_0, \dots, S_{k-1} > y$ such that~$\varphi_{e_i}(G_i, R_i, S_i)$ holds for each~$i < k$.
Let~$R = \bigcup R_i$ and fix some~$y \in \omega$. There are 
finite sets~$S_0, \dots, S_{k-1} > y$ such that~$\varphi_{e_i}(G_i, R_i, S_i)$ holds for each~$i < k$.
Let~$S = \bigcup S_i$. By continuity, there are finite sets~$E_0, \dots, E_{k-1}$ such that
$G_i \uh max(E_i) = F_i \cup E_i$ and~$\varphi_{e_i}(F_i \cup E_i, R_i, S_i)$ holds for each~$i < k$.
By our precise definition of satisfaction, we can even assume without loss of generality
that $(F_0 \cup E_0, \dots, F_{k-1} \cup E_{k-1}, Y)$ is a valid extension of~$c$ for some infinite set~$Y \subseteq X$.
Let~$z \in Y$. In particular, by the definition of being a condition extending~$c$, 
$z \in X$,  $z > max(E_0, \dots, E_{k-1})$
and $F_i \cup E_i \cup \{z\}$ is pseudo-homogeneous for color~$i$ for each~$i < k$.
Therefore~$\psi(R, S)$ holds, as witnessed by~$E_0, \dots, E_{k-1}$ and~$z$.
Thus $\psi(R, S)$ is essential.

Since $A_0, A_1$ is dependently $X \oplus Z$-hyperimmune, then $\psi(R, S)$ holds for some~$R \subseteq \overline{A}_0$
and some~$S \subseteq \overline{A}_1$. Let~$E_0, \dots, E_{k-1} \subseteq X$
be the $k$-tuple of sets and~$z \in X$ be the integer witnessing $\psi(R,S)$.
Let~$i < k$ be such that the set~$Y = \{w \in X \setminus [0, max(E_i)] : f(z, w) = i\}$ is infinite.
The condition~$d = (F_0, \dots, F_{i-1}, F_i \cup E_i \cup \{z\}, F_{i+1}, \dots, F_{k-1}, Y)$ is
a valid extension of~$c$ forcing~$\Rcal_{\vec{e}}$.
\end{proof*}

Let~$\Fcal = \{c_0, c_1, \dots\}$ be a sufficiently generic filter containing $(\emptyset, \dots, \emptyset, \omega)$,
where~$c_s = (F_{0,s}, \dots, F_{k-1,s}, X_s)$. The filter~$\Fcal$ yields a
$k$-tuple of sets~$G_0, \dots, G_{k-1}$ defined by~$G_i = \bigcup_s F_{i,s}$.
By construction, $G_0, \dots, G_{k-1}$ satisfies every condition in~$\Fcal$.
By Lemma~\ref{lem:psrt-cac-fairness-force-Q}, the set~$G_i$ is infinite for each~$i < k$
and by Lemma~\ref{lem:psrt-cac-fairness-force-R}, the pair $A_0, A_1$ is dependently $G_i \oplus Z$-hyperimmune
for some~$i < k$.
\end{proof}

\begin{corollary}
For every~$k \geq 2$, $\rca \wedge \psrt^2_k \not \vdash \scac$.
\end{corollary}
\begin{proof}
Immediate by Theorem~\ref{thm:psrt2k-not-cac-fairness-preservation},
Corollary~\ref{cor:scac-not-cac-fairness-preservation} and
Lemma~\ref{lem:intro-reduc-preservation-separation}.
\end{proof}

\begin{corollary}\label{cor:ads-not-implies-scac}
$\rca \wedge \ads \not \vdash \scac$.
\end{corollary}
\begin{proof}
Immediate by the previous corollary and Theorem~\ref{thm:psrt22-is-ads}.
\end{proof}

\fi

\end{document}